\documentclass[12pt]{amsart}
\usepackage{latexsym,amsxtra,eucal}

\calclayout

\newcommand{\module}{\+module}
\newcommand{\bimodule}{\+bimodule}
\newcommand{\comodule}{\+comodule}
\newcommand{\bicomodule}{\+bi\-co\-mod\-ule}
\newcommand{\contramodule}{\+con\-tra\-mod\-ule}
\newcommand{\semimodule}{\+semi\-mod\-ule}
\newcommand{\bisemimodule}{\+bi\-semi\-mod\-ule}
\newcommand{\semicontramodule}{\+semi\-con\-tra\-mod\-ule}
\newcommand{\coacyclic}{\+coacyclic}
\newcommand{\contraacyclic}{\+con\-tra\-acyclic}
\newcommand{\contractible}{\+con\-tract\-ible}
\newcommand{\projective}{\+pro\-jec\-tive}
\newcommand{\coprojective}{\+co\-pro\-jec\-tive}
\newcommand{\injective}{\+in\-jec\-tive}
\newcommand{\coinjective}{\+coin\-jec\-tive}
\newcommand{\contraflat}{\+con\-tra\-flat}
\newcommand{\semiflat}{\+semi\-flat}
\newcommand{\semiprojective}{\+semi\-pro\-jec\-tive}
\newcommand{\semiinjective}{\+semi\-in\-jec\-tive}
\newcommand{\semialgebra}{\+semi\-al\-ge\-bra}
\newcommand{\category}{\+cat\-e\-gory}
\newcommand{\subcomodule}{\+sub\-co\-mod\-ule}
\newcommand{\invariant}{\+in\-vari\-ant}
\newcommand{\coinvariant}{\+co\-in\-vari\-ant}
\newcommand{\semiinvariant}{\+semi\-in\-vari\-ant}
\newcommand{\semicontrainvariant}{\+semi\-con\-tra\-in\-vari\-ant}

\newcommand{\Section}[1]{\clearpage\section{#1}\medskip}

\theoremstyle{plain}
\newtheorem*{thm}{Theorem}
\newtheorem*{thm1}{Theorem 1}
\newtheorem*{thm2}{Theorem 2}
\newtheorem*{thm3}{Theorem 3}
\newtheorem*{lem}{Lemma}
\newtheorem*{lem1}{Lemma 1}
\newtheorem*{lem2}{Lemma 2}
\newtheorem*{lem3}{Lemma 3}
\newtheorem*{subl}{Sublemma}
\newtheorem*{subla}{Sublemma A}
\newtheorem*{sublb}{Sublemma B}
\newtheorem*{sublc}{Sublemma C}
\newtheorem*{subld}{Sublemma D}
\newtheorem*{prop}{Proposition}
\newtheorem*{prop1}{Proposition 1}
\newtheorem*{prop2}{Proposition 2}
\newtheorem*{prop3}{Proposition 3}
\newtheorem*{cor}{Corollary}
\newtheorem*{cor1}{Corollary 1}
\newtheorem*{cor2}{Corollary 2}
\newtheorem*{cor3}{Corollary 3}
\newtheorem*{cor4}{Corollary 4}
\theoremstyle{definition}
\newtheorem*{rmk}{Remark}
\newtheorem*{rmk1}{Remark 1}
\newtheorem*{rmk2}{Remark 2}
\newtheorem*{rmk3}{Remark 3}
\newtheorem*{qst}{Question}

\newcommand{\ot}{\otimes}
\newcommand{\oc}{\mathbin{\text{\smaller$\square$}}}
\newcommand{\suboc}{\mathbin{\text{\smaller$\.\scriptstyle\square$}}}
\newcommand{\os}{\lozenge}
\newcommand{\ocn}{\odot}
\newcommand{\Ocn}{\circledcirc}

\newcommand{\rarrow}{\longrightarrow}
\newcommand{\larrow}{\longleftarrow}
\newcommand{\from}{\leftarrow}
\newcommand{\mpsto}{\longmapsto}
\newcommand{\eps}{\varepsilon}
\newcommand{\Ups}{\Upsilon}
\newcommand{\dsb}{\dotsb}
\newcommand{\dsc}{\dotsc}
\newcommand{\ilim}{\varinjlim\nolimits}
\newcommand{\plim}{\varprojlim\nolimits}
\renewcommand{\d}{\partial}
\newcommand{\lan}{\langle}
\newcommand{\ran}{\rangle}
\newcommand{\ract}{\mathbin{\leftarrow}}
\newcommand{\lact}{\mathbin{\rightarrow}}
\newcommand{\vot}{\mathbin{\overset{\rightarrow}
                        {\vphantom{o}\smash\ot}}}
\newcommand{\wot}{\mathbin{\overset{\leftarrow}
                        {\vphantom{o}\smash\ot}}}

\renewcommand{\:}{\colon}
\newcommand{\+}{\nobreakdash-} 
\renewcommand{\.}{\text{$\mskip .5\thinmuskip$}}
\renewcommand{\;}{,\medspace}

\newcommand{\lrarrow}{\.\relbar\joinrel\relbar\joinrel\rightarrow\.}
\newcommand{\llarrow}{\.\leftarrow\joinrel\relbar\joinrel\relbar\.}
\newcommand{\darrow}{\mathrel{\dabar\dabar\dahead}}

\DeclareMathOperator{\id}{id}
\DeclareMathOperator{\Id}{Id}
\DeclareMathOperator{\ev}{ev}
\DeclareMathOperator{\coker}{coker}
\DeclareMathOperator{\im}{im}
\DeclareMathOperator{\cone}{cone}
\DeclareMathOperator{\df}{df}
\DeclareMathOperator{\di}{di}
\DeclareMathOperator{\ord}{ord}
\DeclareMathOperator{\tr}{tr}
\newcommand{\s}{{\mathrm s}}

\DeclareMathOperator{\Hom}{Hom}
\DeclareMathOperator{\Tor}{Tor}
\DeclareMathOperator{\Ext}{Ext}
\DeclareMathOperator{\Cotor}{Cotor}
\DeclareMathOperator{\Cohom}{Cohom}
\DeclareMathOperator{\Coext}{Coext}
\DeclareMathOperator{\SemiTor}{SemiTor}
\DeclareMathOperator{\SemiHom}{SemiHom}
\DeclareMathOperator{\SemiExt}{SemiExt}
\DeclareMathOperator{\ProCotor}{ProCotor}
\DeclareMathOperator{\IndCoext}{IndCoext}
\DeclareMathOperator{\Ctrtor}{Ctrtor}
\DeclareMathOperator{\CtrTor}{CtrTor}

\DeclareMathOperator{\End}{End}
\DeclareMathOperator{\Lie}{Lie}
\DeclareMathOperator{\Sym}{Sym}
\DeclareMathOperator{\Ind}{Ind}
\DeclareMathOperator{\Coind}{Coind}
\DeclareMathOperator{\Spec}{Spec}
\DeclareMathOperator{\Inj}{\mathsf{Ind}}
\DeclareMathOperator{\Pro}{\mathsf{Pro}}
\newcommand{\Diff}{\mathrm{Diff}}
\newcommand{\Cl}{\mathrm{Cl}}
\newcommand{\oCl}{\overline{\Cl\!}\,}
\DeclareMathOperator{\Ad}{Ad}

\DeclareMathOperator{\Br}{Bar}
\DeclareMathOperator{\Cb}{Cob}

\global\let\le\undefined
\global\let\ge\undefined
\DeclareMathSymbol{\le}{\mathrel}{AMSa}{"36}      %\leqslant
\DeclareMathSymbol{\ge}{\mathrel}{AMSa}{"3E}      %\geqslant

\DeclareMathSymbol{\birarrow}{\mathrel}{AMSa}{"13}  %rightrightarrows
\DeclareMathSymbol{\boxtimes}{\mathbin}{AMSa}{"02}
\DeclareMathSymbol{\lozenge}{\mathbin}{AMSa}{"06}
\DeclareMathSymbol{\circledcirc}{\mathbin}{AMSa}{"7D}
\DeclareMathSymbol{\Bbbk}{\mathalpha}{AMSb}{"7C}
\DeclareMathSymbol{\kap}{\mathalpha}{AMSb}{"7B}     %\varkappa
\DeclareMathSymbol{\dabar}{\mathord}{AMSa}{"39}
\DeclareMathSymbol{\dahead}{\mathord}{AMSa}{"4B}

\newcommand{\sD}{{\mathsf D}}
\newcommand{\sH}{{\mathsf H}}
\newcommand{\sE}{{\mathsf E}}
\newcommand{\sM}{{\mathsf M}}
\newcommand{\sN}{{\mathsf N}}
\newcommand{\sK}{{\mathsf K}}
\newcommand{\sS}{{\mathsf S}}
\newcommand{\sF}{{\mathsf F}}
\newcommand{\sO}{{\mathsf O}}
\newcommand{\sA}{{\mathsf A}}
\newcommand{\sB}{{\mathsf B}}
\newcommand{\sC}{{\mathsf C}}
\newcommand{\sT}{{\mathsf T}}
\newcommand{\sP}{{\mathsf P}}
\newcommand{\sJ}{{\mathsf J}}
\newcommand{\DG}{\mathsf{DG}}

\newcommand{\boD}{{\mathbb D}}
\newcommand{\boL}{{\mathbb L}}
\newcommand{\boR}{{\mathbb R}}
\newcommand{\boZ}{{\mathbb Z}}
\newcommand{\boQ}{{\mathbb Q}}
\newcommand{\boT}{{\mathbb T}}
\newcommand{\boG}{{\mathbb G}}
\newcommand{\boH}{{\mathbb H}}

\newcommand{\C}{{\mathcal C}}
\newcommand{\D}{{\mathcal D}}
\newcommand{\M}{{\mathcal M}}
\newcommand{\N}{{\mathcal N}}
\renewcommand{\L}{{\mathcal L}}
\newcommand{\R}{{\mathcal R}}
\newcommand{\K}{{\mathcal K}}
\newcommand{\E}{{\mathcal E}}

\newcommand{\cQ}{{\mathcal Q}}
\newcommand{\cP}{{\mathcal P}}
\newcommand{\cG}{{\mathcal G}}
\newcommand{\cZ}{{\mathcal Z}}
\newcommand{\cU}{{\mathcal U}}
\newcommand{\cV}{{\mathcal V}}
\newcommand{\cJ}{{\mathcal J}}
\newcommand{\cA}{{\mathcal A}}

\renewcommand{\P}{{\mathfrak P}}
\newcommand{\Q}{{\mathfrak Q}}
\newcommand{\gK}{{\mathfrak K}}

\newcommand{\gR}{{\mathfrak R}}
\newcommand{\gI}{{\mathfrak I}}
\newcommand{\gG}{{\mathfrak G}}
\newcommand{\gF}{{\mathfrak F}}
\newcommand{\gE}{{\mathfrak E}}
\newcommand{\gU}{{\mathfrak U}}
\newcommand{\gZ}{{\mathfrak Z}}

\renewcommand{\S}{{\boldsymbol{\mathcal S}}}
\newcommand{\T}{{\boldsymbol{\mathcal T}}}
\newcommand{\bM}{{\boldsymbol{\mathcal M}}}
\newcommand{\bN}{{\boldsymbol{\mathcal N}}}
\newcommand{\bL}{{\boldsymbol{\mathcal L}}}
\newcommand{\bK}{{\boldsymbol{\mathcal K}}}

\newcommand{\bcP}{{\boldsymbol{\mathcal P}}}
\newcommand{\bcQ}{{\boldsymbol{\mathcal Q}}}
\newcommand{\bcJ}{{\boldsymbol{\mathcal J}}}
\newcommand{\bR}{{\boldsymbol{\mathcal R}}}
\newcommand{\bE}{{\boldsymbol{\mathcal E}}}
\newcommand{\bcF}{{\boldsymbol{\mathcal F}}}
\newcommand{\bcU}{{\boldsymbol{\mathcal U}}}
\newcommand{\bcZ}{{\boldsymbol{\mathcal Z}}}
\newcommand{\bC}{{\boldsymbol{\mathcal C}}}
\newcommand{\bcX}{{\boldsymbol{\mathcal X}}}
\newcommand{\bcY}{{\boldsymbol{\mathcal Y}}}

\newcommand{\bP}{{\boldsymbol{\mathfrak P}}}
\newcommand{\bQ}{{\boldsymbol{\mathfrak Q}}}
\newcommand{\bgI}{{\boldsymbol{\mathfrak I}}}
\newcommand{\bgK}{{\boldsymbol{\mathfrak K}}}

\newcommand{\bgR}{{\boldsymbol{\mathfrak R}}}
\newcommand{\bgF}{{\boldsymbol{\mathfrak F}}}
\newcommand{\bgE}{{\boldsymbol{\mathfrak E}}}
\newcommand{\bgU}{{\boldsymbol{\mathfrak U}}}
\newcommand{\bgZ}{{\boldsymbol{\mathfrak Z}}}

\newcommand{\be}{\operatorname{\mathbf e}}
\newcommand{\bm}{\operatorname{\mathbf m}}
\newcommand{\bn}{\operatorname{\mathbf n}}
\newcommand{\bp}{\operatorname{\mathbf p}}

\newcommand{\g}{{\mathfrak g}}
\newcommand{\h}{{\mathfrak h}}
\renewcommand{\b}{{\mathfrak b}}
\renewcommand{\aa}{{\mathfrak a}}
\newcommand{\dd}{{\mathfrak d}}
\newcommand{\gl}{{\mathfrak{gl}}}
\newcommand{\hX}{{\,\widehat{\!X}}}
\newcommand{\tS}{{\S\sptilde}}
\newcommand{\tD}{{\D\sptilde}}
\newcommand{\dbar}{{\.\overline{\!\.\d}}}

\newcommand{\bu}{{\text{\smaller\smaller$\scriptstyle\bullet$}}}
\newcommand{\subbu}{{\text{\smaller\smaller
                                  $\scriptscriptstyle\bullet$}}}
\newcommand{\dual}{\spcheck}
\newcommand{\til}{\sptilde}
\newcommand{\comp}{\sphat\.}

\newcommand{\Hot}{\mathsf{Hot}}
\newcommand{\Acycl}{\mathsf{Acycl}}
\newcommand{\Set}{\mathsf{Set}}
\newcommand{\Rep}{\mathsf{Rep}}

\newcommand{\modl}{{\operatorname{\mathsf{--mod}}}}
\newcommand{\modr}{{\operatorname{\mathsf{mod--}}}}
\newcommand{\bimod}{{\operatorname{\mathsf{--mod--}}}}
\newcommand{\comodl}{{\operatorname{\mathsf{--comod}}}}
\newcommand{\comodr}{{\operatorname{\mathsf{comod--}}}}
\newcommand{\comodrinj}{{\operatorname{\mathsf{comod_{inj}--}}}}
\newcommand{\bcomod}{{\operatorname{\mathsf{--comod--}}}}
\newcommand{\contra}{{\operatorname{\mathsf{--contra}}}}
\newcommand{\contraR}{{\operatorname{\mathsf{contra--}}}}
\newcommand{\simodl}{{\operatorname{\mathsf{--simod}}}}
\newcommand{\simodr}{{\operatorname{\mathsf{simod--}}}}
\newcommand{\bsimod}{{\operatorname{\mathsf{--simod--}}}}
\newcommand{\sicntr}{{\operatorname{\mathsf{--sicntr}}}}
\newcommand{\sicntrR}{{\operatorname{\mathsf{sicntr--}}}}
\newcommand{\vect}{{\operatorname{\mathsf{--vect}}}}
\newcommand{\qcomodl}{{\operatorname{\mathsf{--qcmd}}}}
\newcommand{\qcomodr}{{\operatorname{\mathsf{qcmd--}}}}
\newcommand{\qcontra}{{\operatorname{\mathsf{--qcntr}}}}

\newcommand{\up}{\uparrow}
\newcommand{\down}{\downarrow}
\newcommand{\comodrdown}{{\operatorname{\mathsf{comod}}}^\down
        {\operatorname{\mathsf{--}}}}
\newcommand{\simodrdown}{{\operatorname{\mathsf{simod}}}^\down
        {\operatorname{\mathsf{--}}}}
\newcommand{\simodrup}{{\operatorname{\mathsf{simod}}}^\up
        {\operatorname{\mathsf{--}}}}

\newcommand{\co}{{\mathsf{co}}}
\newcommand{\ctr}{{\mathsf{ctr}}}
\newcommand{\si}{{\mathsf{si}}}
\newcommand{\op}{{\mathsf{op}}}
\newcommand{\rop}{{\mathrm{op}}}
\renewcommand{\ss}{{\mathrm{ss}}}
\newcommand{\gr}{{\mathrm{gr}}}
\newcommand{\sgr}{{\mathsf{gr}}}
\newcommand{\inj}{{\mathsf{inj}}}
\newcommand{\proj}{{\mathsf{proj}}}
\newcommand{\fl}{{\mathsf{fl}}}
\newcommand{\sifl}{{\mathsf{sifl}}}
\newcommand{\sipr}{{\mathsf{sipr}}}
\newcommand{\siin}{{\mathsf{siin}}}
\newcommand{\qq}{{\mathsf{q}}}
\newcommand{\fin}{{\mathsf{fin}}}

\newcommand{\cod}{{\operatorname{\mathsf{co-}}}}
\newcommand{\ctrd}{{\operatorname{\mathsf{ctr-}}}}
\newcommand{\injd}{{\operatorname{\mathsf{inj-}}}}
\newcommand{\projd}{{\operatorname{\mathsf{proj-}}}}
\newcommand{\ctrfld}{{\operatorname{\mathsf{ctrfl-}}}}

\begin{document}

\title{Homological algebra of semimodules \\ and semicontramodules
\\[17.5pt] \smaller\smaller\smaller\normalfont\itshape
Semi-infinite Homological Algebra of Associative \\
Algebraic Structures}
\author{Leonid Positselski}
\dedicatory{\ \\[-2pt] \normalfont
With appendices coauthored by S.~Arkhipov and D.~Rumynin\\[-5pt] \ }

\address{Sector of Algebra and Number Theory, Institute for Information
Transmission Problems, 19 Bolshoj Karetnyj per., Moscow, Russia}
\email{posic@mccme.ru}

\begin{abstract}
 We develop the basic constructions of homological algebra in
the (appropriately defined) unbounded derived categories of modules
over algebras over coalgebras over noncommutative rings
(which we call \emph{semialgebras} over \emph{corings}).
 We define double-sided derived functors SemiTor and SemiExt of
the functors of semitensor product and semihomomorphisms,
and construct an equivalence between the exotic derived categories
of semimodules and semicontramodules.

 Certain (co)flatness and/or (co)projectivity conditions have to be
imposed on the coring and semialgebra to make the module categories
abelian (and the cotensor product associative).
 Besides, for a number of technical reasons we mostly have to assume
that the basic ring has a finite homological dimension (no such
assumptions about the coring and semialgebra are made).

 In the final sections we construct model category structures on
the categories of complexes of semi(contra)modules, and develop
relative nonhomogeneous Koszul duality theory for filtered
semialgebras and quasi-differential corings.

 Our motivating examples come from the semi-infinite cohomology theory.
 Comparison with the semi-infinite (co)homology of Tate Lie algebras
and graded associative algebras is established in appendices; and
the semi-infinite homology of a locally compact topological group
relative to an open profinite subgroup is defined.
\end{abstract}

\maketitle

\clearpage
\hfill \emph{To the memory of my father}
\vspace{15ex}
\tableofcontents
\clearpage

\section*{Introduction}
\medskip

 This monograph grew out of the author's attempts to understand
the definitions of semi-infinite (co)homology of associative algebras
that had been proposed in the literature and particularly in
the works of S.~Arkhipov~\cite{Ar1,Ar2} (see also~\cite{Bezr,Sev}).
 Roughly speaking, the semi-infinite cohomology is defined for a Lie
or associative algebra-like object which is split in two halves;
the semi-infinite cohomology has the features of a homology theory
(left derived functor) along one half of the variables and
a cohomology theory (right derived functor) along the other half.

 In the Lie algebra case, the splitting in two halves only has to be
chosen up to a finite-dimensional space; in particular, the homology
of a finite-dimensional Lie algebra only differs from its cohomology
by a shift of the homological degree and a twist of the module of
coefficients.
 So one can define the semi-infinite homology of a Tate (locally
linearly compact) Lie algebra~\cite{BD1} (see also~\cite{Ar3});
it depends, to be precise, on the choice of a compact open vector
subspace in the Lie algebra, but when the subspace changes it
undergoes only a dimensional shift and a determinantal twist.
 Let us emphasize that what is often called the ``semi-infinite
cohomology'' of Lie algebras should be thought of as their
semi-infinite \emph{homology}, from our point of view.
 What we call the semi-infinite \emph{cohomology} of Tate Lie algebras
is a different and dual functor, defined in this book
(see Appendix~\ref{tate-appendix}).

 In the associative case, people usually considered an algebra $A$ with
two subalgebras $N$ and $B$ such that $N\ot B \simeq A$ and there is
a grading on~$A$ for which $N$ is positively graded and locally
finite-dimensional, while $B$ is nonpositively graded.
 We show that both the grading and the second subalgebra $B$ are
redundant; all one needs is an associative algebra $R$,
a subalgebra $K$ in $R$, and a coalgebra~$\C$ dual to~$K$.
 Certain flatness/projectivity and ``integrability'' conditions have
to be imposed on this data.
 If they are satisfied, the tensor product $\S = \C\ot_K \!\.R$ has
a \emph{semialgebra} structure and all the machinery described below
can be applied.

 Furthermore, we propose the following general setting for
semi-infinite (co)ho\-mol\-ogy of associative algebraic structures.
 Let $\C$ be a coalgebra over a field~$k$.
 Then $\C$\+$\C$\bicomodule s form a tensor category with respect to
the operation of cotensor product over~$\C$; the categories of left and
right $\C$\comodule s are module categories over this tensor category.
 Let $\S$ be a ring object in this tensor category; we call such
an object a \emph{semialgebra} over~$\C$ (due to it being ``an algebra
in half of the variables and a coalgebra in the other half'').
 One can consider module objects over~$\S$ in the module categories
of left and right $\C$\comodule s; these are called left and right
\emph{$\S$\+semimodules}.
 The categories of left and right semimodules are only abelian if $\S$
is an injective right and left $\C$\comodule, respectively;
let us suppose that it is.
 There is a natural operation of \emph{semitensor product} of a right
semimodule and a left semimodule over~$\S$; it can be thought of
as a mixture of the cotensor product in the direction of~$\C$ and
the tensor product in the direction of~$\S$ relative to~$\C$.
 This functor is neither left, nor right exact.
 Its double-sided derived functor $\SemiTor$ is suggested as
the associative version of semi-infinite \emph{homology} theory.

 Before describing the functor $\SemiHom$ (whose derived functor
$\SemiExt$ provides the associative version of semi-infinite
\emph{cohomology}), let us discuss a little bit of abstract nonsense.
 Let $\sE$ be an (associative, but noncommutative) tensor category,
$\sM$ be a left module category over it, $\sN$ be a right module
category, and $\sK$ be a category such that there is a pairing between
the module categories $\sM$ and $\sN$ over~$\sE$ taking values in~$\sK$.
 This means that there are multiplication functors
$\sE\times\sE\to \sE$, \  $\sE\times\sM\to\sM$, \
$\sN\times\sE\to\sN$, and $\sN\times\sM\to\sK$
and associativity constraints for ternary multiplications
$\sE\times\sE\times\sE\to\sE$, \ $\sE\times\sE\times\sM\to\sM$, \
$\sN\times\sE\times\sE\to\sN$, and $\sN\times\sE\times\sM\to\sK$
satisfying the appropriate pentagonal diagram equations.
 Let $A$ be a ring object in~$\sE$.
 Then one can consider the category ${}_A\.\sE_A$ of $A$\+$A$\bimodule s
in~$\sE$, the category ${}_A\.\sM$ of left $A$\module s in~$\sM$, and
the category $\sN_A$ of right $A$\module s in~$\sN$.
 If the categories $\sE$, $\sM$, $\sN$, and $\sK$ are abelian,
there are functors of tensor product over~$A$, making ${}_A\.\sE_A$
into a tensor category, ${}_A\.\sM$ and $\sN_A$ into left and right
module categories over ${}_A\.\sE_A$, and providing a pairing
$\sN_A\times{}_A\.\sM \rarrow\sK$.
 These new tensor structures are associative whenever the original
multiplication functors were right exact.
 
 Suppose that we want to iterate this construction, considering
a coring object~$C$ in ${}_A\.\sE_A$, the categories of
$C$\+$C$\bicomodule s in ${}_A\.\sE_A$ and $C$\comodule s in
${}_A\.\sM$ and $\sN_A$, etc.
 Since the functors of tensor product over~$A$ are not left exact
in general, the cotensor products over~$C$ will be only associative
under certain (co)flatness conditions.
 If one makes the next step and considers a ring object~$S$ in
the category of $C$\+$C$\bicomodule s in ${}_A\.\sE_A$, one discovers
that the functors of tensor products over~$S$ are only partially
defined.
 Considering partially defined tensor structures, one can indeed
build this tower of module-comodule categories and tensor-cotensor
products in them as high as one wishes.
 In this book, we restrict ourselves to 3-story towers of
\emph{semialgebras} over \emph{corings} over (ordinary) rings,
mainly because we don't know how to define unbounded (co)derived
categories of (co)modules for any higher levels (see below).

 Now let us introduce \emph{contramodules}.
 The functor $(V,W)\mpsto\Hom_k(V,W)$ makes the category opposite
to the category of vector spaces into a module category over
the tensor category of vector spaces.
 A contramodule over an algebra~$R$ or a coalgebra~$\C$ is an object
of the category opposite to the category of modules or comodules
in $k\vect^\op$ over the ring object~$R$ or the coring object~$\C$
in $k\vect$.
 One can easily see that an $R$\contramodule{} is just an $R$\module{},
while the vector space of $k$\+linear maps from a $\C$\comodule{} to
a $k$\+vector space provides a typical example of $\C$\contramodule.
 Setting $\sE=\sM=k\vect$ and $\sN=\sK=k\vect^\op$ in the above
construction, one obtains a right module category $\C\contra^\op$
over the tensor category $\C\bcomod\C$ together with a pairing
$\Cohom_\C^\op\:\C\comodl\times\C\contra^\op\rarrow k\vect^\op$.
 Given a semialgebra $\S$ over~$\C$, one can apply the construction
again and obtain the category of \emph{$\S$\semicontramodule s} and
the functor $\SemiHom_\S^\op\:\S\simodl\times\S\sicntr^\op\rarrow
k\vect^\op$ assigning a vector space to an $\S$\semimodule{} and
an $\S$\semicontramodule.
 Though comodules and contramodules are quite different, there is
a strong duality-analogy between them on the one hand, and
an equivalence of their appropriately defined (exotic) unbounded
derived categories on the other hand.

 Let us explain how we define double-sided derived functors.
 While the author knows of no natural way to define a derived functor
of one argument that would not be either a left or a right derived
functor, such a definition of derived functor \emph{of two arguments}
does exist in the balanced case.
 Namely, let $\Theta\: \sH_1\times\sH_2\rarrow\sK$ be a functor and
$\sS_i\subset\sH_i$ be localizing classes of morphisms in categories
$\sH_1$ and~$\sH_2$.
 We would like to define a derived functor $\boD\Theta\:
\sH_1[\sS_1^{-1}]\times \sH_2[\sS_2^{-1}] \rarrow \sK$.
 Let $\sF_1$ be the full subcategory of ``flat objects in $\sH_1$
relative to~$\Theta$'' consisting of all objects
$F\in\sH_1$ such that the morphism $\Theta(F,s)$ is an isomorphism
in~$\sK$ for any morphism $s\in\sS_2$.
 Let $\sF_2$ be the full subcategory in $\sH_2$ defined in
the analogous way.
 Suppose that the natural functors $\sF_i[(\sS_i\cap\sF_i)^{-1}]
\rarrow \sH_i[\sS_i^{-1}]$ are equivalences of categories.
 Then the restriction of the functor~$\Theta$ to the subcategory
$\sF_1\times\sH_2$ of the Carthesian product $\sH_1\times\sH_2$
factorizes through
$\sF_1[(\sS_1\cap\sF_1)^{-1}]\times\sH_2[\sS_2^{-1}]$
and therefore defines a functor on the category
$\sH_1[\sS_1^{-1}] \times \sH_2[\sS_2^{-1}]$.
 The same derived functor can be obtained by restricting
the functor~$\Theta$ to the subcategory $\sH_1\times\sF_2$ of
$\sH_1\times\sH_2$.
 This construction can be even extended to partially defined
functors of two arguments~$\Theta$ (see~\ref{semitor-definition}).

 For this definition of the double-sided derived functor to work
properly, the localizing classes in the homotopy categories have to
be carefully chosen (see~\ref{frobenius-derived-coderived-paradox}).
 That is why our derived functors $\SemiTor$ and $\SemiExt$ are not
defined on the conventional derived categories of semimodules and
semicontramodules, but on their \emph{semiderived categories}.
 The semiderived category of $\S$\+semi(contra)modules is a mixture
of the usual derived category in the module direction (relative
to~$\C$) and the \emph{co}/\emph{contraderived} category in
the $\C$\+co/contramodule direction.
 The coderived category of $\C$\comodule s is equivalent to
the homotopy category of complexes of injective $\C$\comodule s,
and analogously, the contraderived category of
$\C$\contramodule s is equivalent to the homotopy category
of complexes of projective $\C$\contramodule s.
 So the distinction between the derived and co/contraderived
categories is only relevant for unbounded complexes and only in
the case of infinite homological dimension.

 A notable attempt to develop a general theory of semi-infinite
homological algebra was undertaken by A.~Voronov in~\cite{Vor}.
 Let us point out the differences between our approaches.
 First of all, Voronov only considers the semi-infinite homology of
Lie algebras, while we work with associative algebraic structures.
 Secondly, Voronov constructs a double-sided derived functor of
a functor of one argument and the choice of a class of resolutions
becomes an additional ingredient of his construction, while we
define double-sided derived functors of functors of two arguments
and the conditions imposed on resolutions are determined by
the functors themselves.
 Thirdly, Voronov works with graded Lie algebras and his functor
of semivariants is obtained as the image of the invariants with respect
to one half of the Lie algebra in the coinvariants with respect to
the other half, while we consider ungraded Tate Lie algebras with only
one subalgebra chosen, and our functor of \emph{semiinvariants} is
constructed in a more delicate way (see below).
 Finally, no exotic derived categories appear in~\cite{Vor}.

 The coderived category of $\C$\comodule s and the contraderived
category of $\C$\contramodule s turn out to be naturally equivalent.
 This equivalence can be thought of as a covariant analogue of
the contravariant functor $\boR\Hom({-},R)\:\sD(R\modl)\rarrow
\sD(\modr R)$ on the derived category of modules over a ring~$R$.
 Moreover, there is a natural equivalence between the semiderived
categories of $\S$\semimodule s and $\S$\semicontramodule s.
 The functors $\boR\Psi_\S\:\sD^\si(\S\simodl)\rarrow
\sD^\si(\S\sicntr)$ and $\boL\Phi_\S\:\sD^\si(\S\sicntr)\rarrow
\sD^\si(\S\simodl)$ providing this equivalence are defined in terms
of the spaces of homorphisms in the category of $\S$\semimodule s
and the operation of \emph{contratensor product} of
an $\S$\semimodule{} and an $\S$\semicontramodule.
 The latter is a right exact functor which resembles the functor of
tensor product of modules over a ring.
 This equivalence of triangulated categories tranforms the functor
$\SemiExt_\S$ into the functors $\Ext$ in either of the semiderived
categories (and the functor $\SemiTor^\S$ into the left derived
functor $\CtrTor^\S$ of the functor of contratensor product).
 We call this kind of equivalence of triangulated categories
the \emph{comodule-contramodule correspondence} or
the \emph{semimodule-semicontramodule correspondence}.

 The duality-analogy between semimodules and semicontramodules partly
breaks down when one passes from homological algebra to the structure
theory.
 Comodules over a coalgebra over a field are simplistic creatures;
contramodules are quite a bit more complicated, though still much
simpler than modules over a ring, the structure theory of a coalgebra
over a field being much simpler than that of an algebra or a ring.
 We construct some relevant counterexamples.
 There is an analogue of Nakayama's Lemma for contramodules,
a description of contramodules over an infinite direct sum
of coalgebras, etc.
 These results can be extended to contramodules over certain
topological rings (much more general than the topological algebras
dual to coalgebras).
 Contramodules over topological Lie algebras can also be defined;
and an isomorphism of the categories of contramodules over
a topological Lie algebra and its topological enveloping algebra
can be proven under certain assumptions.

 A \emph{coring} $\C$ over a ring~$A$ is a coring object in the tensor
category of bimodules over~$A$.  (In a different terminology, this is
called a \emph{coalgebroid}.)
 A \emph{semialgebra} $\S$ over a coring~$\C$ is a ring object in
the tensor category of bicomodules over~$\C$; for this definition
to make sense, certain (co)flatness conditions have to be imposed
on~$\C$ and~$\S$ to make the cotensor product of bicomodules
well-defined and associative.
 Throughout this monograph (with the exception of
Section~\ref{preliminaries-section} and the appendices) we work with
corings~$\C$ over noncommutative rings~$A$ and semialgebras~$\S$
over~$\C$.
 Mostly we have to assume that the ring $A$ has a finite homological
dimension---for a number of reasons, the most important one being
that otherwise we don't know how to define appropriately the unbounded
(co)derived category of $\C$\comodule s.
 No assumptions about the homological dimension of the coring and
the semialgebra are made.
 Besides, we mostly have to suppose that $\C$ is a flat left and right
$A$\module{} and $\S$ is a coflat left and right $\C$\comodule,
and even certain (co)projectivity conditions have to be imposed in
order to work with contramodules.

 Nonhomogeneous quadratic duality~\cite{Pos,PP} establishes
a correspondence between nonhomogeneous Koszul algebras and Koszul
CDG\+algebras.
 This duality has a relative version with a base ring, assigning,
e.~g., the de Rham DG\+algebra to the filtered algebra of
differential operators (the base ring being the ring of functions,
in this case).
 For a number of reasons, it is advisable to avoid passing to
the dual vector space/module in this construction, working with
CDG\+coalgebras instead of CDG\+algebras; in particular, this allows
to include infinitely-(co)generated Koszul algebras and coalgebras.
 In the relative case, this means considering the graded coring of
polyvector fields, rather than the graded algebra of differential
forms, as the dual object to the differential operators.
 The relevant additional structure on the polyvector fields
(corresponding to the de Rham differential on the differential forms)
is that of a \emph{quasi-differential coring}.
 Another important version of relative nonhomogeneous quadratic
duality uses base coalgebras in place of base rings.
 This situation is simpler in some respects, since one still obtains
CDG\+coalgebras as the dual objects.
 As a generalization of these two dualities, one can consider
nonhomogeneous Koszul semialgebras over corings and assign Koszul
quasi-differential corings over corings to them.
 The Poincare-Birkhoff-Witt theorem for Koszul semialgebras claims
that this correspondence is an equivalence of categories.

 Relative nonhomogeneous Koszul duality theorem provides an equivalence
between the semiderived category of semimodules over a nonhomogeneous
Koszul semialgebra and the coderived category of quasi-differential
comodules over the corresponding quasi-differential coring, and
an analogous equivalence between the semiderived category of
semicontramodules and the contraderived category of quasi-differential
contramodules.
 In particular, for a smooth algebraic variety $M$ and a vector
bundle $E$ over $M$ with a global connection~$\nabla$, there is
an equivalence between the derived category of modules over
the algebra/sheaf of differential operators acting in the sections
of $E$ and the coderived category (and also the contraderived category,
when $M$ is affine) of CDG\+modules over the CDG\+algebra
$\Omega(M,\End(E))$ of differential forms with coefficients in
the vector bundle $\End(E)$, where the differential~$d$ in
$\Omega(M,\End(E))$ is the de Rham differential depending on~$\nabla$
and the curvature element $h\in\Omega^2(M,\End(E))$ is the curvature
of~$\nabla$.

 Natural examples of semialgebras and semimodules come from Lie theory.
 Namely, let $(\g,H)$ be an algebraic Harish-Chandra pair, i.~e.,
$\g$ is a Lie algebra over a field~$k$ and $H$ is a smooth affine
algebraic group corresponding to a finite-dimensional Lie
subalgebra $\h\subset\g$.
 Let $\C(H)$ be the coalgebra of functions on~$H$.
 Then the category $\sO(\g,H)$ of Harish-Chandra modules is isomorphic
to the category of left semimodules over the semialgebra
$\S(\g,H)=U(\g)\ot_{U(\h)}\C(H)$.
 If the group $H$ is unimodular, the semialgebra $\S=\S(\g,H)$ has
an involutive anti-automorphism.
 In general, the opposite semialgebras $\S$ and $\S^\rop$
are Morita-equivalent in some sense; more precisely, there is
a canonical left $\S\ot_k\S$\semimodule{} $\bE=\bE(\g,H)$ such that
the semitensor product with~$\bE$ provides an equivalence between
the categories of right and left $\S$\semimodule s.
 Geometrically, $\bE(\g,H)$ is the bimodule of distributions on
an algebraic group $G$ supported in its subgroup~$H$ and regular
along~$H$.
 So the semitensor product of $\S$\semimodule s can be considered as
a functor on the category $\sO(\g,H)\times\sO(\g,H)$.
 This functor factorizes through the functor of tensor product in
the category $\sO(\g,H)$ and is closely related to the functor of
\emph{$(\g,H)$\+semiinvariants} $\bM\mpsto\bM_{\g,H}$
on the category of $(\g,H)$\module s.
 The semiinvariants are a mixture of invariants over~$H$ and
coinvariants along~$\g/\h$.

 More generally, let $(\g,H)$ be a \emph{Tate Harish-Chandra pair},
that is $\g$ is a Tate Lie algebra and $H$ is an affine proalgebraic
group corresponding to a compact open subalgebra $\h\subset\g$.
 Let $\kap\:(\g',H)\rarrow(\g,H)$ be a morphism of Tate Harish-Chandra
pairs with the same proalgebraic group $H$ such that the Lie algebra
map $\g'\rarrow\g$ is a central extension whose kernel is identified
with~$k$; assume also that $H$ acts trivially in $k\subset\g'$.
 One example of such a central extension of Tate Harish-Chandra pairs
comes from the canonical central extension $\g\til$ of~$\g$; we denote
the corresponding morphism by $\kap_0$.
 There is a semialgebra $\S_\kap(\g,H)=U_\kap(\g)\ot_{U(\h)}\C(H)$ over
the coalgebra~$\C(H)$ such that the category of left semimodules
over $\S_\kap=\S_\kap(\g,H)$ is isomorphic to the category of discrete
$(\g',H)$\module s where the unit central element of~$\g'$ acts by
the identity (Harish-Chandra modules with the central charge~$\kap$).
 Left semicontramodules over the opposite semialgebra $\S_\kap^\rop$
can be described in terms of compatible structures of
$\g'$\contramodule s and $\C(H)$\contramodule s.
 These are called \emph{Harish-Chandra contramodules} with the central
charge~$\kap$; the dual vector spaces to Harish-Chandra modules with
the central charge~$-\kap$ can be found among them.

 The semialgebras $\S_\kap$ and $\S_{-\kap_0-\kap}^\rop$ are
naturally isomorphic, at least, when the pairing $U(\h)\ot_k\C(H)
\rarrow k$ is nondegenerate in $\C(H)$.
 In view of the semimodule-semicontramodule correspondence theorem,
it follows that the semiderived categories of Harish-Chandra modules
with the central charge~$\kap$ and Harish-Chandra contramodules
with the central charge~$\kap+\kap_0$ over $(\g,H)$ are naturally
equivalent.
 So the well-known phenomenon of correspondence between complexes
of modules with complementary central charges over certain
infinite-dimensional Lie algebras can be formulated as an equivalence
of triangulated categories using the notions of contramodules and
semiderived categories.
 Besides, it follows that the category of right semimodules over
$\S_\kap$ is isomorphic to the category of Harish-Chandra modules
with the central charge~$-\kap-\kap_0$.
 When the proalgebraic group $H$ is prounipotent (and $\h$ is
exactly the Lie algebra of~$H$), the object
$\SemiTor^{\S_\kap}(\bN^\bu,\bM^\bu)$ of the derived category of
$k$\+vector spaces is represented by the complex of semi-infinite
forms over~$\g$ with coefficients in~$\bN^\bu\ot_k\bM^\bu$.
 This provides a comparison of our theory of SemiTor with
the semi-infinite homology of Tate Lie algebras.
 Semi-infinite cohomology of Lie algebras, whose coefficients
are contramodules over (the canonical central extensions of)
Tate Lie algebras, is related to SemiExt in the analogous way.

 To a topological group $G$ with an open profinite subgroup $H$ and
a commutative ring~$k$ one can associate a semialgebra $\S_k(G,H)$
over the coring $\C_k(H)$ of $k$\+valued locally constant functions
on $H$ such that the categories of left and right semimodules over
$\S_k(G,H)$ are isomorphic to the category of smooth $G$\module s
over~$k$.
 So the category of semimodules over $\S_k(G,H)$ does not depend on~$H$,
neither does the category of semicontramodules over $\S_k(G,H)$; all
the semialgebras $\S_k(G,H)$ with a fixed $G$ and varying $H$ are
naturally Morita equivalent.
 The semiderived categories of semimodules and semicontramodules
over $\S_k(G,H)$ do depend on $H$ essentially, however, as do
the functors $\SemiTor$ and $\SemiExt$ over $\S_k(G,H)$.
 These double-sided derived functors may be called the semi-infinite
(co)homology of a group with an open profinite subgroup.
 The semi-infinite homology of topological groups is a mixture of
the discrete group homology and the profinite group cohomology.

 Examples of corings~$\C$ over commutative rings~$A$ for which
the left and the right actions of~$A$ in~$\C$ are different come
from the algebraic groupoids theory, and examples of semialgebras
over such corings come from Lie theory of groupoids.
 Namely, let $(M,H)$ be a smooth affine groupoid, i.~e., $M$ and~$H$
are smooth affine algebraic varieties, there are two smooth morphisms
$s_H,\.t_H\: H\birarrow M$ of source and target, and there are unit,
multiplication, and inverse element morphisms satisfying the usual
groupoid axioms.
 Let $A=A(M)$ be the ring of functions on~$M$ and $\C=\C(H)$ be
the ring of functions on~$H$.
 Then $\C$ is a coring over~$A$.
 Moreover, suppose that $(M,H)$ is a closed subgroupoid of
a groupoid $(M,G)$.
 Let $\g$ and $\h$ be the Lie algebroids over the ring~$A$
corresponding to the groupoids $(M,G)$ and $(M,H)$, and let
$U_A(\g)$ and $U_A(\h)$ be their enveloping algebras.
 Then there is a semialgebra $\S=\S_M(G,H) = U_A(\g)\ot_{U_A(\h)}\C(H)$
over the coring~$\C$ and a canonical left $\S\ot_k\S$\semimodule{}
$\bE=\bE_M(G,H)$ providing an equivalence between the categories of
right and left $\S$\semimodule s.
 The semimodule $\bE$ consists of all distributions on~$G$ twisted
with the line bundle $s_G^*(\Omega_M^{-1})\ot t^*_G(\Omega_M^{-1})$,
supported in~$H$ and regular along~$H$ (where $\Omega_M$ denotes
the bundle of top forms on~$M$).
 
 Examples of corings over noncommutative rings come from
Noncommutative Geometry\-~\cite{KR}.
 Noncommutative stacks are represented as quotients of
noncommutative affine schemes corresponding to rings~$A$
by actions of corings~$\C$ over~$A$.
 The cotensor product of $\C$\comodule s can be understood as
the group of global sections of the tensor product of a right
and a left sheaf over a noncommutative stack, while the tensor
product of sheaves itself does not exist.

 Notice that the roles of the ring and coring structures in our
constructions are not symmetric; in particular, we have to consider
conventional derived categories along the algebra variables and
co/contraderived categories along the coalgebra variables.
 The cause of this difference is that the tensor product of modules
commutes with the infinite direct sums, but not with the infinite
products.
 This can be changed by passing to pro-objects; consequently one can
still define versions of derived functors $\Cotor$ and $\Coext$ over
a coring $\C$ without making any homological dimension assumptions
at all by considering pro- and ind-modules (see
Remarks~\ref{semitor-definition} and~\ref{semiext-definition}).
 A problem remains to construct the comodule-contramodule
correspondence without any homological dimension assumptions on
the ring~$A$.
 Here we only manage to weaken the finite homological dimension
assumption to the Gorensteinness assumption.

 Algebras/coalgebras over fields and semialgebras over coalgebras
over fields are briefly discussed in Section~0.
 Semialgebras over corings and the functors of semitensor product
over them are introduced in Section~1, and important constructions
of flat comodules and coflat semimodules are presented there.
 The derived functor $\SemiTor$ is defined in Section~2.
 Contramodules over corings and semicontramodules over semialgebras
are introduced in Section~3, and the derived functor $\SemiExt$
is defined in Section~4.
 Equivalence of exotic derived categories of comodules and
contramodules is proven in Section~5; and the same for semimodules
and semicontramodules is done in Section~6.
 Functors of change of ring and coring for the categories of
comodules and contramodules are introduced in Section~7;
functors of change of coring and semialgebra for the categories
of semimodules and semicontramodules are constructed in Section~8.
 Closed model category structures on the categories of complexes of
semimodules and semicontramodules are defined in Section~9.
 The construction of a semialgebra depending on three embedded rings
and a coring dual to the middle ring is considered in Section~10.
 The Poincare--Birkhoff--Witt theorem and the Koszul duality theorem
for nonhomogeneous Koszul semialgebras are proven in Section~11.
 The basic structure theory of contramodules over a coalgebra
over a field is developed in Appendix~A\hbox{}.
 We compare our theory of $\SemiExt$ and $\SemiTor$ with Arkhipov's
and Sevostyanov's semi-infinite $\Ext$ and~$\Tor$ in Appendix~B\hbox{}.
 Semialgebras corresponding to Harish-Chandra pairs and their
Hopf algebra analogues are discussed in Appendix~C\hbox{}.
 Tate Harish-Chandra pairs are considered in Appendix~D\hbox{}, and
the theorem of comparison with semi-infinite cohomology of Tate
Lie algebras is proven there.
 Semialgebras corresponding to topological groups are discussed in
Appendix~E\hbox{}.
 Pairs of algebraic groupoids are considered in Appendix~F\hbox{}.

 Appendix~C was written in collaboration with Dmitriy Rumynin.
 Appendix~D was written in collaboration with Sergey Arkhipov.

 One terminological note: we will generally use the words
\emph{the homotopy category of} (an additive category $\sA$) and
\emph{the homotopy category of complexes of} (objects from~$\sA$)
as synonymous.
 Analogously, \emph{the homotopy category of complexes} (with
a particular property) \emph{over\/ $\sA$} is a full subcategory
of the homotopy category of~$\sA$.

\subsection*{Acknowledgements}
 I am grateful to B.~Feigin for posing the problem of defining
the semi-infinite cohomology of associative algebras back in
the first half of 1990's.
 Even earlier, I learned about the problem of constructing a derived
equivalence between modules with complementary central charges from
M.~Kapranov's handwritten notes on Koszul duality.
 S.~Arkhipov patiently explained me his ideas about the semi-infinite
cohomology many times over the years, contributing to my efforts to
understand the subject.
 In the Summer of 2000, this work was stimulated by discussions with
S.~Arkhipov and R.~Bezrukavnikov, and my gratitude goes to both
of them.
 I wish to thank J.~Bernstein, B.~Feigin, B.~Keller, V.~Lunts, and
V.~Vologodsky for helpful conversations, and I.~Mirkovic for
stimulating interest.
 Parts of the mathematical content of this monograph were worked out
when the author was visiting Stockholm University, Weizmann Institute,
Independent University of Moscow, Max-Planck-Institut in Bonn, and
Warwick University; I was supported by the European Post-Doctoral
Institute during a part of that time.
 The author was partially supported by grants from CRDF, INTAS, and
P.~Deligne's 2004 Balzan prize while writing the manuscript up.

\setcounter{section}{-1}
\Section{Preliminaries and Summary} \label{preliminaries-section}

 This section contains some known results and some results deemed
to be new, but no proofs.
 Its goal is to prepare the reader for the more technically involved
constructions of the main body of the monograph (where the proofs are
given).
 In particular, we don't have to worry about nonassociativity
of the cotensor product and partial definition of the semitensor
product here, distiguish between the myriad of notions of
absolute/relative coflatness/coprojectivity/injectivity of comodules
and analogously for contramodules, etc., because we only consider
coalgebras over fields.

\subsection{Unbounded Tor and Ext}
 Let $R$ be an algebra over a field~$k$.

\subsubsection{}   \label{prelim-unbounded-tor}
 We would like to extend the familiar definition of the derived
functor of tensor product $\Tor^R\: \sD^-(\modr R) \times
\sD^-(R\modl) \rarrow \sD^-(k\vect)$ on the Carthesian product
of the derived categories of right and left $R$\module s bounded
from above, so as to obtain a functor on the Carthesian
product of unbounded derived categories.

 As always, the tensor product of a complex of right $R$\module s
$N^\bu$ and a complex of left $R$\module s $M^\bu$ is defined as
the total complex of the bicomplex $N^i\ot_R M^j$, constructed by
taking infinite direct sums along the diagonals.
 This provides a functor $\Hot(\modr R)\times\Hot(R\modl)\rarrow
\Hot(k\vect)$ on the Carthesian product of unbounded homotopy
categories of $R$\module s.

 The most straightforward way to define the object
$\Tor^R(N^\bu,M^\bu)$ of $\sD(k\vect)$ is to represent it by
the total complex of the bicomplex
 $$
  \dsb\lrarrow N^\bu\ot_k R\ot_k R\ot_k M^\bu \lrarrow
  N^\bu\ot_k R\ot_k M^\bu \lrarrow N^\bu\ot_k M^\bu,
 $$
constructed by taking infinite direct sums along the diagonals.
 One can check that this bar construction indeed defines a functor
$$
 \Tor^R\: \sD(\modr R)\times\sD(R\modl)\lrarrow\sD(k\vect).
$$

 The unbounded derived functor $\Tor^R$ can be also defined by
restricting the functor of tensor product to appropriate subcategories
of complexes adjusted to the functor of tensor product in the unbounded
homotopy categories of $R$\module s.
 Namely, let us call a complex of left $R$\module s $M^\bu$ \emph{flat}
if the complex of $k$\+vector spaces $M^\bu\ot_R N^\bu$ is acyclic
whenever a complex of right $R$\module s $N^\bu$ is acyclic.
 \emph{Not every complex of flat $R$\module s is a flat complex of
$R$-modules according to this definition.}

 In particular, an acyclic complex of left $R$\module s is flat if
and only if it is \emph{pure}, i.~e., it remains acyclic after taking
the tensor product with any right $R$\module.
 So an acyclic complex of flat $R$\module s is flat if and only if
all of its modules of cocycles are flat.
 On the other hand, any complex of flat $R$\module s bounded from
above is flat.
 If the ring $R$ has a finite weak homological dimension, then
any complex of flat $R$\module s is flat.
 For example, the acyclic complex $M^\bu$ of free modules over the ring
of dual numbers $R=k[\eps]/\eps^2$ whose every term is equal to~$R$
and every differential is the operator of multiplication with~$\eps$
is not flat.
 Indeed, let $N^\bu=(\dsb \to k[\eps/\eps^2]\to k \to 0\to\dsb)$ be
a free resolution of the $R$\module{}~$k$; then the complex
$N^\bu\ot_R M^\bu$ is quasi-isomorphic to $k\ot_R M^\bu$ and has
a one-dimensional cohomology space in every degree, even though
the complex $N^\bu$ is acyclic.

 Any complex of $R$\module s is quasi-isomorphic to a flat complex,
and moreover, the quotient category of the homotopy category
$\Hot_\fl(R\modl)$ of flat complexes of $R$\module s  by the thick
subcategory of acyclic flat complexes $\Acycl(R\modl)\cap
\Hot_\fl(R\modl)$ is equivalent to the derived category $\sD(R\modl)$.
 This result holds for an arbitrary ring~\cite{Spal}, and even
for an arbitrary DG\+ring~\cite{Kel1,BL}.
 The derived functor $\Tor^R$ can be defined by restricting
the functor of tensor product over~$R$ to either of the full
subcategories $\Hot(\modr R)\times \Hot_\fl(R\modl)$ or
$\Hot_\fl(\modr R)\times \Hot(R\modl)$
of the category $\Hot(\modr R)\times\Hot(R\modl)$.

\subsubsection{}    \label{prelim-unbounded-ext}
 The functor $\Hom_R\:\Hot(R\modl)^\op\times\Hot(R\modl)\rarrow
\Hot(k\vect)$ and its derived functor
$\Ext_R\:\sD(R\modl)^\op\times\sD(R\modl)\rarrow \sD(k\vect)$
need no special definition: once the unbounded homotopy and derived
categories are defined, so are the spaces of homomorphisms in them.
 For any (unbounded) complexes of left $R$\module s $L^\bu$ and
$M^\bu$, the total complex of the cobar bicomplex
$$
 \Hom_k(L^\bu,M^\bu)\lrarrow\Hom_k(R\ot_k L^\bu\;M^\bu)
 \lrarrow\Hom_k(R\ot_k R\ot_k L^\bu\;M^\bu)\lrarrow\dsb,
$$
constructed by taking infinite direct products along the diagonals,
represents the object $\Ext_R(L^\bu,M^\bu)$ in $\sD(k\vect)$.

 The unbounded derived functor $\Ext_R$ can be also computed by
restricting the functor $\Hom_R$ to appropriate subcategories in
the Carthesian product of homotopy categories of $R$\module s.
 Let us call a complex of left $R$\module s $L^\bu$ \emph{projective}
if the complex $\Hom_R(L^\bu,M^\bu)$ is acyclic for any acyclic
complex of left $R$\module s $M^\bu$.
 Analogously, a complex of left $R$\module s $M^\bu$ is called
\emph{injective} if the complex $\Hom_R(L^\bu,M^\bu)$ is acyclic
for any acyclic complex of left $R$\module s $L^\bu$.

 Any projective complex of $R$\module s is flat.
 Any complex of projective $R$\module s bounded from above is
projective, and any complex of injective $R$\module s bounded
from below is injective.
 If the ring~$R$ has a finite left homological dimension, then any
complex of projective left $R$\module s is projective and any complex
of injective left $R$\module s is injective.
 
 A complex of $R$\module s is projective if and only if it
belongs to the minimal triangulated subcategory of the homotopy
category of $R$\module s containing the complex
$\dsb\to 0 \to R\to 0 \to\dsb$ and closed under infinite direct sums.
 Analogously, a complex of $R$\module s is injective if and only if 
up to the homotopy equivalence it can be obtained from the complex
$\dsb\to 0\to \Hom_k(R,k)\to 0 \to\dsb$ using the operations of
shift, cone, and infinite direct product.
 The homotopy category $\Hot_\proj(R\modl)$ of projective complexes
of $R$\module s and the homotopy category $\Hot_\inj(R\modl)$ of
injective complexes of $R$\module s are equivalent to the unbounded
derived category $\sD(R\modl)$.
 The results mentioned in this paragraph even hold for an arbitrary
DG\+ring~\cite{Kel1,BL}.
 The functor $\Ext_R$ can be obtained by restricting the functor
$\Hom_R$ to either of the full subcategories $\Hot_\proj(R\modl)^\op
\times\Hot(R\modl)$ or $\Hot(R\modl)^\op\times\Hot_\inj(R\modl)$
of the category $\Hot(R\modl)^\op\times\Hot(R\modl)$.

\subsubsection{}
 The definitions of unbounded Tor and Ext in terms of (co)bar
constructions were known at least since the 1960's.
 The notions of flat, projective, and injective (unbounded) complexes
of $R$\module s were introduced by N.~Spaltenstein~\cite{Spal} (who 
attributes the idea to J.~Bernstein).
 Such complexes were called ``$K$\+flat'', ``$K$\+projective'', and
``$K$\+injective'' in~\cite{Spal}; they are often called
``$H$\+projective'' or ``homotopy projective'' etc.\ nowadays.

\subsection{Coalgebras over fields; Cotor and Coext}
 The notion of a coalgebra over a field is obtained from that of
an algebra by formal dualization.
 Since any coassociative coalgebra is the union of its
finite-dimensional subcoalgebras, the category of coalgebras is
anti-equivalent to the category of profinite-dimensional algebras.
 There are two ways of dualizing the notion of a module over
an algebra: one can consider \emph{comodules} and \emph{contramodules}
over a coalgebra.
 Comodules can be thought of as discrete modules which are unions of
their finite-dimensional submodules, while contramodules are modules
where certain infinite summation operations are defined.
 Dualizing the constructions of the tensor product of modules and
the space of homomorphisms between modules, one obtains the functors
of cotensor product and cohomomorphisms.
 Their derived functors are called $\Cotor$ and $\Coext$.

\subsubsection{}
 A coassociative \emph{coalgebra} with counit over a field $k$ is
a $k$\+vector space $\C$ endowed with a \emph{comultiplication} map
$\mu_\C\: \C\rarrow\C\ot_k\C$ and a \emph{counit} map $\eps_\C\:
\C\rarrow k$ satisfying the equations dual to the associativity and
unity equations on the multiplication and unit maps of an assotiative
algebra with unit.
 More precisely, one should have $(\mu_\C\otimes\id_\C)\circ\mu_\C =
(\id_\C\otimes\mu_\C)\circ\mu_\C$ and $(\eps_\C\otimes\id_\C)\circ\mu_\C
= \id_\C = (\id_\C\otimes\eps_\C)\circ\mu_\C$.

 A \emph{left comodule} $\M$ over a coalgebra $\C$ is a $k$\+vector
space endowed with a \emph{left coaction} map $\nu_\M: \M\rarrow
\C\ot_k\M$ satisfying the equations dual to the associativity and unity
equations on the action map of a module over an associative algebra
with unit.
 More precisely, one should have $(\mu_\C\otimes\id_\M)\circ\nu_\M =
(\id_\C\otimes\nu_\M)\circ\nu_\M$ and $(\eps_\C\otimes\id_\M)
\circ\nu_\M = \id_\M$.
 A \emph{right comodule} $\N$ over a coalgebra $\C$ is a $k$\+vector
space endowed with a \emph{right coaction} map $\nu_\N: \N\rarrow
\N\ot_k\C$ satisfying the coassociativity and counity equations
$(\nu_\N\otimes\id_\C)\circ\nu_\N =(\id_\N\otimes\mu_\C)\circ\nu_\N$
and $(\id_\N\otimes\eps_\C)\circ\nu_\N = \id_\N$.
 For example, the coalgebra $\C$ has natural structures of a left and
a right comodule over itself.

 The categories of left and right $\C$\comodule s are abelian.
 We will denote them by $\C\comodl$ and $\comodr\C$, respectively.
 Both infinite direct sums and infinite products exist in
the category of $\C$\comodule s, but only infinite direct sums
are preserved by the forgetful functor $\C\comodl\rarrow k\vect$
(while the infinite products are not even exact in~$\C\comodl$).
 A \emph{cofree} $\C$\comodule{} is a $\C$\comodule{} of the form
$\C\ot_k V$, where $V$ is a $k$\+vector space.
 The space of comodule homomorphisms into the cofree $\C$\comodule{} is
described by the formula $\Hom_\C(\M\;\C\ot_k V) \simeq \Hom_k(\M,V)$.
 The category of $\C$\comodule s has enough injectives; besides,
a left $\C$\comodule{} is injective if and only if it is
a direct summand of a cofree $\C$\comodule.

 The \emph{cotensor product} $\N\oc_\C \M$ of a right $\C$\comodule{}
$\N$ and a left $\C$\comodule{} $\M$ is defined as the kernel of
the pair of maps
 $$
  (\nu_\N\ot\id_\M\; \id_\N\ot\nu_\M)\:
  \N\ot_k\M \birarrow \N\ot_k\C\ot_k\M.
 $$
 This is the dual construction to the tensor product of a right module
and a left module over an associative algebra.
 There are natural isomorphisms $\N\oc_\C \C\simeq \N$ and
$\C\oc_\C\M\simeq \M$.
 The functor of cotensor product over~$\C$ is left exact.

\subsubsection{}  \label{prelim-unbounded-cotor}
 The cotensor product $\N^\bu\oc_\C \M^\bu$ of a complex of right
$\C$\comodule s $\N^\bu$ and a complex of left $\C$\comodule s $\M^\bu$
is defined as the total complex of the bicomplex $\N^i\oc_\C \M^j$,
constructed by taking infinite direct sums along the diagonals.

 We would like to define the derived functor $\Cotor^\C$ of the functor
of cotensor product in such a way that it could be obtained by
restricting the functor~$\oc_\C$ to appropriate subcategories of
the Carthesian product of homotopy categories $\Hot(\comodr\C)$
and $\Hot(\C\comodl)$.
 In addition, we would like the object $\Cotor^\C(\N^\bu,\M^\bu)$ of
$\sD(k\vect)$ to be represented by the total complex of the cobar
bicomplex
 \begin{equation} \label{cotor-cobar-complex}
  \N^\bu\ot_k \M^\bu \lrarrow \N^\bu\ot_k \C\ot_k \M^\bu 
  \lrarrow\N^\bu\ot_k \C\ot_k \C\ot_k \M^\bu \lrarrow\dsb,
 \end{equation}
constructed by taking infinite direct sums along the diagonals.
 It turns out that a functor $\Cotor^\C$ with these properties does
exist, but it is \emph{not defined on the Carthesian product of
conventional unbounded derived categories\/ $\sD(\comodr\C)$ and\/
$\sD(\C\comodl)$}.

 For example, let $\C$ be the coalgebra dual to the algebra of dual
numbers $\C^*=k[\eps]/\eps^2$, so that $\C$\comodule s are just
$k[\eps]/\eps^2$\module s.
 Let $\M^\bu$ be the acyclic complex of cofree $\C$\comodule s whose
every term is equal to~$\C$ and every differential is the operator of
multiplication with~$\eps$, and let $\N^\bu$ be the complex of
$\C$\comodule s whose only nonzero term is the $\C$\comodule~$k$.
 Then the cobar complex that we want to compute
$\Cotor^\C(\N^\bu,\M^\bu)$ is quasi-isomorphic to the complex
$\N^\bu\oc_\C \M^\bu$ and has a one-dimensional cohomology space
in every degree, even though $\M^\bu$ represents a zero object
in $\sD(\C\comodl)$.
 Therefore, a more refined version of unbounded derived category of
$\C$\comodule s has to be considered.

 A complex of left $\C$\comodule s is called \emph{coacyclic} if it
belongs to the minimal triangulated subcategory of the homotopy
category $\Hot(\C\comodl)$ containing the total complexes of exact
triples ${}'\K^\bu \to \K^\bu \to {}''\K^\bu$ of complexes
of left $\C$\comodule s and closed under infinite direct sums.
 (By the total complex of an exact triple of complexes we mean
the total complex of the corresponding bicomplex with three rows.)
 Any coacyclic complex is acyclic; any acyclic complex bounded
from below is coacyclic.
 The complex $\M^\bu$ from the above example is acyclic, but not
coacyclic.
 (Indeed, the total complex of the cobar
bicomplex~\eqref{cotor-cobar-complex} is acyclic whenever
$\M^\bu$ is coacyclic.)
 The \emph{coderived category} of left $\C$\comodule s
$\sD^\co(\C\comodl)$ is defined as the quotient category of
the homotopy category $\Hot(\C\comodl)$ by the thick subcategory
of coacyclic complexes $\Acycl^\co(\C\comodl)$.

 In the same way one can define the coderived category of any
abelian category with exact functors of infinite direct sum.
 Over a category of finite homological dimension, every acyclic
complex belongs to the minimal triangulated subcategory of
the homotopy category containing the total complexes of exact
triples of complexes, even without the infinite direct sum closure.

 The cotensor product $\N^\bu\oc_\C \M^\bu$ of a complex of right
$\C$\comodule s $\M^\bu$ and a complex of left $\C$\comodule s $\N^\bu$
is acyclic whenever one of the complexes $\M^\bu$ and $\N^\bu$ is
coacyclic and the other one is a complex of injective $\C$\comodule s.
 Besides, the coderived category $\sD^\co(\C\comodl)$ is equivalent to
the homotopy category $\Hot(\C\comodl_\inj)$ of injective
$\C$\comodule s.
 Thus one can define the unbounded derived functor
 $$
  \Cotor^\C\: \sD^\co(\comodr\C)\times\sD^\co(\C\comodl)
  \lrarrow\sD(k\vect)
 $$
by restricting the functor of cotensor product to either of
the full subcategories $\Hot(\comodr\C)\times\Hot(\C\comodl_\inj)$
or $\Hot(\comodrinj\C)\times\Hot(\C\comodl)$ of the category
$\Hot(\comodr\C)\times\Hot(\C\comodl)$.

\subsubsection{}   \label{frobenius-derived-coderived-paradox}
 If one attempts to construct a derived functor of cotensor product
on the Carthesian product of conventional unbounded derived
categories of comodules in the way analogous
to~\ref{prelim-unbounded-tor},
the result may not look like what one expects.

 Consider the example of a finite-dimensional coalgebra $\C$ dual
to a Frobenius algebra $\C^*=F$.
 Let us assume the convention that left $\C$\comodule s are left
$F$\module s and right $\C$\comodule s are right $F$\module s.
 Then the functor $\oc_\C$ is left exact and the functor $\ot_F$
is right exact, but the difference between them is still rather small:
if either a left (co)module $M$, or a right (co)module $N$ is
projective-injective, then there is a natural isomorphism
$N\oc_\C M\simeq(N\oc_\C F)\ot_F M$, and after one chooses
an isomorphism between the left modules $F$ and~$\C$,
the right modules $N$ and $N\oc_\C F$ will only differ by
the Frobenius automorphism of the Frobenius algebra~$F$.
 
 So if one defines ``coflat'' complexes of $\C$\comodule s as
the complexes whose cotensor product with acyclic complexes is
acyclic, then the quotient category of the homotopy category of
``coflat'' complexes by the thick subcategory of acyclic ``coflat''
complexes will be indeed equivalent to the derived category of
comodules, and one will be able to define a ``derived functor
of cotensor product over~$\C$'' in this way, but the resulting
derived functor will coincide, up to the Frobenius twist, with
the functor $\Tor^F$.
 (Indeed, any flat complex of flat modules will be ``coflat''.)
 When the complexes this functor is applied to are concentrated in
degree~$0$, this functor will produce a complex situated in
the negative cohomological degrees, as is characteristic of $\Tor^F$,
and not in the positive ones, as one would expect of $\Cotor^\C$.

 Likewise, if one attempts to construct a derived functor of
tensor product on the Carthesian product of coderived
categories of modules in the way analogous
to~\ref{prelim-unbounded-cotor}, one will find, in the Frobenius
algebra case, that the tensor product of a complex of projective
$F$\module s with a coacyclic complex is acyclic, the homotopy
category of complexes of projective modules is indeed equivalent
to the coderived category of $F$\module s, and one can define
a ``derived functor of tensor product over~$F$'' by restricting to
this subcategory, but the resulting derived functor will coincide,
up to the Frobenius twist, with the functor $\Cotor^\C$.

 Nevertheless, it is well known how to define a derived functor
of cotensor product on the conventional unbounded derived categories
of comodules (see~\ref{derived-first-second-kind},
cf.\ Remark~\ref{semitor-definition}).

\subsubsection{}
 The category $k\vect^\op$ opposite to the category of vector spaces has
a natural structure of a \emph{module category} over the tensor category
$k\vect$ with the action functor $k\vect\times k\vect^\op\rarrow
k\vect^\op$ defined by the rule $(V,W^\op)\mpsto\Hom_k(V,W)^\op$.
 More precisely, there are two module category structures associated
with this functor: the left module category with the associativity
constraint $\Hom_k(U\ot_k V\;W)\simeq \Hom_k(U,\Hom_k(V,W))$ and
the right module category with the associativity constraint
$\Hom_k(U\ot_k V\;W)\simeq \Hom_k(V,\Hom_k(U,W))$.
 The category of \emph{left contramodules} over a coalgebra $\C$
is the opposite category to the category of comodule objects in
the \emph{right} module category $k\vect^\op$ over the coring
object~$\C$ in the tensor category $k\vect$.
 Explicitly, a $\C$\contramodule{} $\P$ is a $k$\+vector space
endowed with a \emph{contraaction} map $\pi_\P\: \Hom_k(\C,\P)
\rarrow \P$ satisfying the \emph{contraassociativity} and
\emph{counity} equations $\pi_\P\circ\Hom(\id_\C,\pi_\P) =
\pi_\P\circ\Hom(\mu_\C,\id_\P)$ and $\pi_\P\circ \Hom(\eps_\C,\id_\P)
= \id_\P$.

 For any right $\C$\comodule{} $\N$ and any $k$\+vector space $V$
the space $\Hom_k(\N,V)$ has a natural structure of left
$\C$\contramodule.
 The category of left $\C$\contramodule s is abelian.
 We will denote it by $\C\contra$.
 Both infinite direct sums and infinite products exist in
the category of contramodules, but only the infinite products
are preserved by the forgetful functor $\C\contra\rarrow k\vect$
(while the infinite direct sums are not even exact in $\C\contra$).
 The category of contramodules has enough projectives.
 Besides, a $\C$\contramodule{} is projective if and only if it
is a direct summand of a \emph{free} $\C$\contramodule{} of
the form $\Hom_k(\C,V)$ for some vector space~$V$.
 The space of contramodule homomorphisms from the free
$\C$\contramodule{} is described by the formula
$\Hom^\C(\Hom_k(\C,V),\.\P) \simeq \Hom_k(V,\.\P)$.

 Let $\M$ be a left $\C$\comodule{} and $\P$ be a left
$\C$\contramodule.
 The \emph{space of cohomomorphisms} $\Cohom_\C(\M,\P)$ is defined
as the cokernel of the pair of maps
\begin{multline*}
 (\Hom(\nu_\M,\id_\P),\Hom(\id_\M,\pi_\P))\: \\
  \Hom_k(\C\ot_k\M\;\P) = \Hom_k(\M,\Hom_k(\C,\P))
  \birarrow \Hom_k(\M,\P).
\end{multline*}
 This is the dual construction to that of the space of homomorphisms
between two modules over a ring.
 There are natural isomorphisms $\Cohom_\C(\C,\P)\simeq\P$ and
$\Cohom_\C(\M,\Hom_k(\N,V)) \simeq \Hom_k(\N\oc_\C\M\;V)$.
 The functor of cohomomorphisms over~$\C$ is right exact.

\subsubsection{}
 The complex of cohomomorphisms $\Cohom_\C(\M^\bu,\P^\bu)$ from
a complex of left $\C$\comodule s $\M^\bu$ to a complex of left
$\C$\contramodule s $\P^\bu$ is defined as the total complex of
the bicomplex $\Cohom_\C(\M^i,\P^j)$, constructed by taking
infinite products along the diagonals.
 Let us define the derived functor $\Coext_\C$ of the functor
of cohomomorphisms.

 A complex of $\C$\contramodule s is called \emph{contraacyclic}
if it belongs to the minimal triangulated subcategory of
the homotopy category $\Hot(\C\contra)$ containing the total
complexes of exact triples ${}'\!\.\gK^\bu\to \gK^\bu \to
{}''\!\.\gK^\bu$ of complexes of $\C$\contramodule s and closed
under infinite products.
 Any contraacyclic complex is acyclic; any acyclic complex bounded
from above is contraacyclic.
 The \emph{contraderived category} of $\C$\contramodule s
$\sD^\ctr(\C\contra)$ is defined as the quotient category of
the homotopy category $\Hot(\C\contra)$ by the thick subcategory
of contraacyclic complexes $\Acycl^\ctr(\C\contra)$.

 The complex of cohomomorphisms $\Cohom_\C(\M^\bu,\P^\bu)$ is acyclic
whenever either $\M$ is a complex of injective $\C$\comodule s
and $\P$ is contraacyclic, or $\M$ is coacyclic and $\P$ is
a complex of projective $\C$\contramodule s.
 Besides, the contraderived category $\sD^\ctr(\C\contra)$ is
equivalent to the homotopy category of projective $\C$\contramodule s
$\Hot(\C\contra_\proj)$.
 Thus one can define the derived functor
 $$
  \Coext_\C\:\sD^\co(\C\comodl)^\op\times\sD^\ctr(\C\contra)
  \lrarrow\sD(k\vect)
 $$
by restricting the functor $\Cohom_\C$ to either of
the subcategories $\Hot(\C\comodl_\inj)^\op\times\Hot(\C\contra)$
or $\Hot(\C\comodl)^\op\times\Hot(\C\contra_\proj)$ of the Carthesian
product $\Hot(\C\comodl)^\op\times\Hot(\C\contra)$.

 The contramodule version of bar construction provides a functorial
complex computing $\Coext_\C$.
 Namely, for any complex of left $\C$\comodule s $\M^\bu$ and complex
of left $\C$\contramodule s $\P^\bu$ the total complex of the bicomplex
$$
 \dsb\lrarrow\Hom_k(\C\ot_k\C\ot_k\M^\bu\;\P^\bu)\lrarrow
 \Hom_k(\C\ot_k\M^\bu\;\P^\bu)\lrarrow\Hom_k(\M^\bu,\P^\bu),
$$
constructed by taking infinite products along the diagonals,
represents the object $\Coext_\C(\M^\bu,\P^\bu)$ in $\sD(k\vect)$.

\subsubsection{}   \label{prelim-co-contra-correspondence}
 The categories of left $\C$\comodule s and left $\C$\contramodule s
are isomorphic if the coalgebra $\C$ is finite-dimensional, but
in general they are quite different.
 However, the coderived category of left $\C$\comodule s is naturally
equivalent to the contraderived category of left $\C$\contramodule s,
$\sD^\co(\C\comodl)\simeq\sD^\ctr(\C\contra)$.

 Indeed, the coderived category $\sD^\co(\C\comodl)$ is equivalent to
the homotopy category $\Hot(\C\comodl_\inj)$ and the contraderived
category $\sD^\ctr(\C\contra)$ is equivalent to the homotopy category
$\Hot(\C\contra_\proj)$.
 Furthermore, the additive category of injective $\C$\comodule s is
the idempotent closure of the category of cofree $\C$\comodule s and
the additive category of projective $\C$\contramodule s is
the idempotent closure of the category of free $\C$\contramodule s.
 One has $\Hom_\C(\C\ot_k U\;\C\ot_k V) = \Hom_k(\C\ot_k U\;V)
=\Hom_k(U,\,\Hom_k(\C,V))=\Hom^\C(\Hom_k(\C,U),\Hom_k(\C,V))$, so
the categories of cofree comodules and free contramodules
are equivalent.
 
 To describe this equivalence of additive categories in a more
invariant way, let us define the operation of contratensor
product of a comodule and a contramodule.

 Let $\N$ be a right $\C$\+comodule{} and $\P$ be a left
$\C$\+contramodule.
 The \emph{contratensor product} $\N\ocn_\C\P$ is defined as
the cokernel of the pair of maps
 $$
  ((\id_N\ot\ev_\C)\circ(\nu_N\ot\id_{\Hom_k(\C,\P)})\;
  \id_N\circ\pi_\P)\:\N\ot_k\Hom_k(\C,\P)\birarrow\N\ot_k\P,
 $$
where $\ev_\C$ denotes the evaluation map $\C\ot_k\Hom_k(\C,\P)\to\P$.
 The contratensor product functor is not a part of any tensor
or module category structure; instead, it is dual to the functors
$\Hom$ in the categories of $\C$\comodule s and $\C$\contramodule s.
 The functor of contratensor product over~$\C$ is right exact.
 There are natural isomorphisms $\N\ocn_\C\Hom_k(\C,V)\simeq\N\ot_k V$
and $\Hom_k(\N\ocn_\C\P\;V) \simeq \Hom^\C(\P,\Hom_k(\N,V))$.

 The desired equivalence between the additive categories of injective
left $\C$\+co\-mod\-ules and projective left $\C$\contramodule s is
provided by the pair of adjoint functors $\Psi_\C(\M)=\Hom_\C(\C,\M)$
and $\Phi_\C(\P)=\C\ocn_\C\P$ between the categories of left
$\C$\comodule s and left $\C$\contramodule s.
 Here the space $\Hom_\C(\C,\M)$ is endowed with a $\C$\contramodule{}
structure as the kernel of a pair of contramodule morphisms
$\Hom_k(\C,\M)\birarrow\Hom_k(\C\;\C\ot_k\M)$ (where the contramodule
structure on $\Hom_k(\C,\M)$ and $\Hom_k(\C\;\C\ot_k\M)$ comes from
the right $\C$\comodule{} structure on~$\C$), while the space
$\C\ocn_\C\P$ is endowed with a left $\C$\comodule{} structure as
the cokernel of a pair of comodule morphisms
$\C\ot_k\Hom_k(\C,\P)\birarrow\C\ot_k\P$.

\subsubsection{}  \label{prelim-co-contra-acyclicity}
 The following examples illustrate the necessity of considering
the exotic derived categories in the above construction of
the derived comodule-contramodule correspondence.
 Let $W$ be a vector space and $\C$ be the symmetric coalgebra
of~$W$.
 One can construct $\C$ as the subcoalgebra of the tensor coalgebra
$\bigoplus_{n=0}^\infty W^{\ot n}$ formed by the symmetric tensors.
 Consider the trivial left $\C$\contramodule{} $k$; it has a left
projective $\C$\contramodule{} resolution of the form
 $$
  \dsb\lrarrow\Hom_k(\C,({\textstyle\bigwedge}^2_k\,W)^*)\lrarrow
  \Hom_k(\C,W^*)\lrarrow\Hom_k(\C,k).
 $$
 Applying the functor $\Phi_\C$ to the above complex of contramodules,
one obtains the complex of injective left $\C$\comodule s
 \begin{equation}  \label{acyclic-comodule}
  \dsb\lrarrow\C\ot_k({\textstyle\bigwedge}^2_k\,W)^*\lrarrow\C\ot_kW^*
  \lrarrow\C.
 \end{equation}
 When $W$ is finite-dimensional, the complex~\eqref{acyclic-comodule}
has its only nonvanishing cohomology in the degree $-\dim W$;
this cohomology a trivial one-dimensional $\C$\comodule{} naturally
isomorphic to $\det(W)^*=(\bigwedge^{\dim W}_kW)^*$ as a vector space.
 When $W$ is infinite-dimensional, the complex~\eqref{acyclic-comodule}
is acyclic; one can think of it as an ``injective resolution
of a one-dimensional $\C$\comodule{} concentrated in
the degree~$-\infty$''.
 So when $\dim W=\infty$ the equivalence of categories
$\sD^\co(\C\comodl)\simeq\sD^\ctr(\C\contra)$ transforms the acyclic
complex of $\C$\comodule s~\eqref{acyclic-comodule} into the trivial
$\C$\contramodule{} $k$ considered as a complex concentrated
in degree~$0$, and back.

 Analogously, consider the trivial left $\C$\comodule{} $k$; it has
a right injective $\C$\comodule{} resolution of the form
 $$
  \C\lrarrow\C\ot_kW\lrarrow\C\ot_k{\textstyle\bigwedge}^2_k\,W
  \lrarrow\dsb
 $$
 Applying the functor $\Psi_\C$ to this complex of comodules,
one obtains the complex of projective left $\C$\contramodule s
 $$
  \Hom_k(\C,k)\lrarrow\Hom_k(\C,W)\lrarrow
  \Hom_k(\C,{\textstyle\bigwedge}^2_k\,W)\lrarrow\dsb
 $$
 When $W$ is finite-dimensional, the latter complex has its only
nonvanishing cohomology in the degree $\dim W$; the cohomology
a trivial one-dimensional $\C$\contramodule{} naturally isomorphic to
$\det(W)$ as a vector space.
 When $W$ is infinite-dimensional, this complex is acyclic;
one can think of it as a ``projective resolution of a one-dimensional
$\C$\contramodule{} concentrated in the degree~$+\infty$''.
 In this case, the equivalence of categories 
$\sD^\co(\C\comodl)\simeq\sD^\ctr(\C\contra)$ transforms
the trivial $\C$\comodule{} $k$ considered as a complex concentrated
in degree~$0$ into this acyclic complex of $\C$\contramodule s,
and back.

 The cohomology computations above are very similar to computing
$\Ext_R(k,R)$ for the algebra $R$ of polynomials in a finite or
infinite number of variables over a field~$k$.

\subsubsection{}
 The functor $\Ext_\C\:\sD^\co(\C\comodl)^\op\times\sD^\co(\C\comodl)
\rarrow\sD(k\vect)$ of homomorphisms in the coderived category
$\sD^\co(\C\comodl)$ can be computed by restricting the functor
$\Hom_\C\:\Hot(\C\comodl)^\op\times\Hot(\C\comodl)\rarrow\Hot(k\vect)$
of homomorphisms in the homotopy category $\Hot(\C\comodl)$ to
the full subcategory $\Hot(\C\comodl)^\op\times\Hot(\C\comodl_\inj)$
of the category $\Hot(\C\comodl)^\op\times\Hot(\C\comodl)$.
 The complex $\Hom_\C(\L^\bu,\M^\bu)$ is acyclic whenever $\L^\bu$
is a coacyclic complex of left $\C$\comodule s and $\M^\bu$ is
a complex of injective left $\C$\comodule s.
 
 Analogously, the functor $\Ext^\C\:\sD^\ctr(\C\contra)^\op\times
\sD^\ctr(\C\contra)\rarrow\sD(k\vect)$ of homomorphisms in
the contraderived category $\sD^\ctr(\C\contra)$
can be computed by restricting the functor
$\Hom^\C\:\Hot(\C\contra)^\op\times\Hot(\C\contra)\rarrow\Hot(k\vect)$
to the full subcategory $\Hot(\C\contra_\proj)^\op\times\Hot(\C\contra)$
of the category $\Hot(\C\contra)^\op\times\Hot(\C\contra)$.
 The complex $\Hom^\C(\P^\bu,\Q^\bu)$ is acyclic whenever $\P^\bu$
is a complex of projective $\C$\contramodule s and $\Q^\bu$ is
a contraacyclic complex of $\C$\contramodule s.

 The contratensor product $\N^\bu\ocn_\C\P^\bu$ of a complex of right
$\C$\comodule s $\N^\bu$ and a complex of left $\C$\contramodule s
$\P^\bu$ is defined as the total complex of the bicomplex
$\N^i\ocn_\C\P^j$, constructed by taking infinite direct sums along
the diagonals.
 The complex $\N^\bu\ocn_\C\P^\bu$ is acyclic whenever $\N^\bu$ is
a coacyclic complex of right $\C$\comodule s and $\P^\bu$ is
a complex of projective left $\C$\contramodule s.
 The left derived functor $\Ctrtor^\C$ of the functor of
contratensor product,
$$
 \Ctrtor^\C\:\sD^\co(\comodr\C)\times\sD^\ctr(\C\contra)
 \lrarrow\sD(k\vect),
$$
is defined by restricting the functor of contratensor product to
the full subcategory $\Hot(\comodr\C)\times\Hot(\C\contra_\proj)$
of the category $\Hot(\comodr\C)\times\Hot(\C\contra)$.
 Notice that the (abelian or homotopy) category of right
$\C$\comodule s does not contain enough objects adjusted to
contratensor product.

 The equivalence of triangulated categories $\sD^\co(\C\comodl)
\simeq \sD^\ctr(\C\contra)$ transforms the functor $\Coext_\C$ into
either of the functors $\Ext_\C$ or $\Ext^\C$ and the functor
$\Cotor^\C$ into the functor $\Ctrtor^\C$.

\subsubsection{}
 A left $\C$\comodule{} $\M$ is called \emph{coflat} if the functor
$\N\mpsto\N\oc_\C\M$ is exact on the category of right $\C$\comodule s.
 A left $\C$\comodule{} $\M$ is called \emph{coprojective} if
the functor $\P\mpsto\Cohom_\C(\M,\P)$ is exact on the category of
left $\C$\contramodule s.
 It is easy to see that an injective comodule is coprojective and
a coprojective comodule is coflat.
 Using the fact that any comodule is a union of its finite-dimensional
subcomodules, one can show that any coflat comodule is injective.
 Thus all the three conditions are equivalent.

 A left $\C$\contramodule{} $\P$ is called \emph{contraflat} if
the functor $\N\mpsto\N\ocn_\C\P$ is exact on the category of
right $\C$\comodule s.
 A left $\C$\contramodule{} $\P$ is called \emph{coinjective} if
the functor $\M\mpsto\Cohom_\C(\M,\P)$ is exact on the category of
left $\C$\comodule s.
 It is easy to see that a projective contramodule is coinjective
and a coinjective contramodule is contraflat.
 We will show in~\ref{cotensor-contratensor-assoc} that any
coinjective contramodule is projective and
in~\ref{contraflat-contramodules} that any contraflat contramodule
is projective.
 Thus all the three conditions are equivalent.

\subsubsection{}    \label{derived-first-second-kind}
 Our definition of the derived functor of cotensor product for
unbounded complexes differs from the most traditional one, which was
introduced (in the greater generality of DG\+coalgebras and
DG\comodule s) by Eilenberg and Moore~\cite{EM1}.
 Husemoller, Moore, and Stasheff~\cite{HMS} call the functor defined by
Eilenberg--Moore \emph{the differential derived functor of cotensor
product of the first kind} and denote it by $\Cotor^{\C,I}$ or simply
$\Cotor^\C$, while the functor $\Cotor^\C$ defined
in~\ref{prelim-unbounded-cotor} is (the nondifferential particular case
of) what they call \emph{the differential derived functor of cotensor
product of the second kind} and denote by $\Cotor^{\C,II}$.

 The functor $\Cotor^{\C,I}$ is computed by the total complex of
the cobar bicomplex~\eqref{cotor-cobar-complex}, constructed by taking
infinite \emph{products} along the diagonals (while the tensor product
complexes $\N^\bu\ot \C\ot\dsb\ot \C\ot \M^\bu$ constituting the cobar
bicomplex are still defined as infinite direct sums).
 It is indeed a functor on the Carthesian product of conventional
unbounded derived categories $\sD(\comodr\C)$ and $\sD(\C\comodl)$.

 The unbounded derived functor $\Tor^R$ defined
in~\ref{prelim-unbounded-tor}
is a derived functor of the first kind in this terminology.
 Roughly, derived functors of the first kind correspond to the
conventional derived categories $\sD$ (which can be therefore called
\emph{derived categories of the first kind}), while derived functors
of the second kind correspond to the coderived and contraderived
categories $\sD^\co$ and $\sD^\ctr$ (which can be called
\emph{derived categories of the second kind}).
 The distinction, which is only relevant for unbounded complexes of
modules (comodules, or contramodules), manifests itself also for
quite finite-dimensional DG\module s (DG\comodule s, or
DG\contramodule s).

 The coderived categories of comodules were introduced by
K.~Lef\`evre-Hasegawa \cite{Lef,Kel2} in the context of Koszul duality;
our definition is equivalent to (the nondifferential case of) his one.
 They first appeared in the author's own research in the very same
context.
 An elaborate discussion of the two kinds of derived categories and
their roles in Koszul duality can be found in~\cite{Pos2}; a proof
of the equivalence of the two definitions is also given there.
 Contramodules were defined by Eilenberg and Moore in~\cite{EM2} and
studied by Barr in~\cite{Bar}.

 All the most important results of this subsection
can be extended straightforwardly to DG\+coalgebras and even
CDG\+coalgebras (see~\ref{cdg-rings}
and~\ref{quasi-differential-corings} for the definition).
 Generally, the constructions of derived categories and functors
of the first kind can be generalized to $A_\infty$\+algebras,
while the constructions of derived categories and functors of
the second kind can be naturally extended to CDG\+coalgebras.

\subsection{Semialgebras over coalgebras over fields}
 The notion of a semialgebra over a coalgebra is dual to that 
of a coring over a noncommutative ring.
 The similarity between the two theories only goes so far,
however.

\subsubsection{}
 Let $\C$ and $\D$ be two coalgebras over a field~$k$.
 A \emph{$\C$\+$\D$\+bicomodule} $\K$ is $k$\+vector space
endowed with a left $\C$\comodule{} and a right $\D$\comodule{}
structures such that the left $\C$\+coaction map
$\nu'_\K\:\K\rarrow\C\ot_k\K$
is a morphism of right $\D$\comodule s, or, equivalently,
the right $\D$\+coaction map $\nu''_\K\:\K\rarrow\K\ot_k\D$
is a morphism of left $\C$\comodule s.
 A bicomodule can be also defined as a vector space endowed with
a \emph{bicoaction} map $\K\rarrow\C\ot_k\K\ot_k\D$ satisfying
the coassociativity and counity equations.
 The category of $\C$\+$\D$\bicomodule s is abelian.
We will denote it by $\C\bcomod\D$.

 Let $\C$, \ $\D$, and $\E$ be three coalgebras, $\N$ be
a $\C$\+$\D$\bicomodule, and $\M$ be a $\D$\+$\E$\bicomodule.
 Then the cotensor product $\N\oc_\D\M$ is endowed with
a $\C$\+$\E$\bicomodule{} structure as the kernel of a pair
of bicomodule morphisms $\N\ot_k\M\birarrow\N\ot_k\D\ot_k\M$.
 The cotensor product of bicomodules is associative: for
any coalgebras $\C$ and~$\D$, any right $\C$\comodule{} $\N$,
left $\D$\comodule{} $\M$, and $\C$\+$\D$\bicomodule{} $\K$
there is a natural isomorphism $\N\oc_\C(\K\oc_\D\M)\simeq
(\N\oc_\C\K)\oc_\D\M$.
 
\subsubsection{}
 In particular, the category of $\C$\+$\C$\bicomodule s is
an associative tensor category with the unit object~$\C$.
 A \emph{semialgebra} $\S$ over~$\C$ is an associative ring object
with unit in this tensor category; in other words, it is
a $\C$\+$\C$\bicomodule{} endowed with a \emph{semimultiplication} map
$\bm_\S\:\S\oc_\C\S\rarrow\S$ and a \emph{semiunit} map 
$\be_\S\:\C\rarrow\S$ which have to be $\C$\+$\C$\bicomodule{}
morphisms satisfying the associativity and unity equations
$\bm_\S\circ(\bm_\S\oc\id_\S)=\bm_\S\circ(\id_\S\oc \bm_\S)$ and
$\bm_\S\circ(\be_\S\oc\id_\S)=\id_\S=\bm_\S\circ(\id_\S\oc \be_\S)$.

 The category of left $\C$\comodule s is a left module category
over the tensor category $\C\bcomod\C$, and the category of right
$\C$\comodule s is a right module category over it.
 A \emph{left semimodule} $\bM$ over~$\S$ is a module object in
this left module category over the ring object $\S$ in this tensor
category; in other words, it is a left $\C$\comodule{} endowed
with a \emph{left semiaction} map $\bn_\bM\:\S\oc_\C\bM\rarrow\bM$,
which has to be a morphism of left $\C$\comodule s satisfying
the associativity and unity equations
$\bn_\bM\circ(\bm_\S\oc\id_\bM)=\bn_\bM\circ(\id_\S\oc \bn_\bM)$
and $\bn_\bM\circ(\be_\S\oc\id_\bM)=\id_\bM$.
 A \emph{right semimodule} $\bN$ over~$\S$ is a right $\C$\comodule{}
endowed with a \emph{right semiaction} morphism of right
$\C$\comodule s $\bn_\bN\:\bN\oc_\C\S \rarrow\bN$ satisfying
the equations $\bn_\bN\circ(\bn_\bN\oc\id_\S)= \bn_\bN\circ
(\id_\bN\oc \bm_\S)$ and $\bn_\bN\circ(\id_\bN\oc \be_\S)=\id_\bN$.

 For any left $\C$\comodule{} $\L$, the cotensor product $\S\oc_\C\L$
has a natural left semimodule structure.
 It is called the left $\S$\semimodule{} \emph{induced} from a left
$\C$\comodule{}~$\L$.
 The space of semimodule homomorphisms from the induced semimodule is
described by the formula $\Hom_\S(\S\oc_\C\L\;\bM)\simeq
\Hom_\C(\L,\bM)$.
 We will denote the category of left $\S$\semimodule s by $\S\simodl$
and the category of right $\S$\semimodule s by $\simodr\S$. 
 The category of left $\S$\semimodule s is abelian provided that $\S$
is an injective right $\C$\comodule.
 Moreover, $\S$ is an injective right $\C$\comodule{} if and only if
the category $\S\simodl$ is abelian and the forgetful functor
$\S\simodl\rarrow\C\comodl$ is exact.

 The operation of cotensor product over~$\C$ provides a pairing functor
$\comodr\C\times\C\comodl\rarrow k\vect$ compatible with the right
module category structure on $\comodr\C$ and the left module category
structure on $\C\comodl$ over the tensor category $\C\bcomod\C$.
 The \emph{semitensor product} $\bN\os_\S\bM$ of a right
$\S$\semimodule{} $\bN$ and a left $\S$\semimodule{} $\bM$ is defined
as the cokernel of the pair of maps
 $$
  (\bn_\bN\oc\id_\bM\;\id_\bN\oc\bn_\bM)\: 
   \bN\oc_\C\S\oc_\C\bM \birarrow \bN\oc_\C\bM.
 $$
 There are natural isomorphisms $\bN\os_\S(\S\oc_\C\L) \simeq
\bN\oc_\C\L$ and $(\R\oc_\C\S)\os_\S\bM \simeq \R\oc_\C\bM$.
 The functor of semitensor product is neither left, nor right exact.

\subsubsection{}
 The semitensor product $\bN^\bu\os_\S\bM^\bu$ of a complex of right
$\S$\semimodule s $\bN^\bu$ and a complex of left $\S$\semimodule s
$\bM^\bu$ is defined as the total complex of the bicomplex
$\bN^i\os_\S\bM^j$, constructed by taking infinite direct sums along
the diagonals.
 Assume that $\S$ is an injective left and right $\C$\comodule.
 We would like to define the double-sided derived functor
$\SemiTor^\S$ of the functor of semitensor product.

 The \emph{semiderived category} of left $\S$\semimodule s
$\sD^\si(\S\simodl)$ is defined as the quotient category of
the homotopy category $\Hot(\S\simodl)$ by the thick subcategory
$\Acycl^{\cod\C}(\S\simodl)$ of complexes of $\S$\semimodule s that are
\emph{coacyclic as complexes of\/ $\C$\comodule s}.
 For example, if the coalgebra $\C$ coincides with the ground field~$k$,
and $\S=R$ is just a $k$\+algebra, then the semiderived category
$\sD^\si(\S\simodl)$ coincides with the derived category $\sD(R\modl)$,
while if the semialgebra $\S$ coincides with the coalgebra~$\C$, 
then the semiderived category $\sD^\si(\S\simodl)$ coincides with
the coderived category $\sD^\co(\C\comodl)$.

 A complex of left $\C$\comodule s~$\bM^\bu$ is called \emph{semiflat}
if the semitensor product $\bN^\bu\os_\S\bM^\bu$ is acyclic for any
$\C$\coacyclic{} complex of right $\S$\semimodule s $\bN^\bu$.
 For example, the complex of $\S$\semimodule s $\S\oc_\C\L^\bu$ induced
from a complex of injective $\C$\comodule s $\L^\bu$ is semiflat.

 The quotient category of the homotopy category $\Hot_\sifl(\S\simodl)$
of semiflat complexes of $\S$\semimodule s by the thick subcategory of
$\C$\coacyclic{} semiflat complexes $\Acycl^{\cod\C}(\S\simodl)\cap
\Hot_\sifl(\S\simodl)$ is equivalent to the semiderived
category of $\S$\semimodule s.
 The derived functor
$$
 \SemiTor^\S\: \sD^\si(\simodr\S)\times\sD^\si(\S\simodl)
 \lrarrow \sD(k\vect)
$$
is defined by restricting the functor of semitensor product over~$\S$
to either of the full subcategories $\Hot(\simodr\S)\times
\Hot_\sifl(\S\simodl)$ or $\Hot_\sifl(\simodr\S)\times\Hot(\S\simodl)$
of the category $\Hot(\simodr\S)\times \Hot(\S\simodl)$.

\subsubsection{}
 Let $\C$ and $\D$ be two coalgebras, $\K$ be a $\C$\+$\D$\bicomodule,
and $\P$ be a left $\C$\contramodule.
 Then the space of cohomomorphisms $\Cohom_\C(\K,\P)$ is endowed with
a left $\D$\contramodule{} structure as the cokernel of a pair of
$\D$\contramodule{} morphisms $\Hom_k(\C\ot_k\K\;\P)=
\Hom_k(\K,\Hom_k(\C,\P)) \birarrow \Hom_k(\K,\P)$.
 For any left $\D$\comodule{} $\M$, left $\C$\contramodule{} $\P$,
and $\C$\+$\D$\bicomodule{} $\K$ there is a natural isomorphism
$\Cohom_\C(\K\oc_\D\M\;\P)\simeq \Cohom_\D(\M,\Cohom_\C(\K,\P))$.

\subsubsection{}
 Therefore, the category opposite to the category of left
$\C$\contramodule s is a right module category over the tensor
category of $\C$\+$\C$\bicomodule s with respect to the action
functor $\Cohom_\C$.
 The category of \emph{left $\S$\+semicontramodules} is the opposite
category to the category of module objects in the right module
category $\C\contra^\op$ over the ring object~$\S$ in the tensor
category $\C\bcomod\C$.
 In other words, a left semicontramodule $\bP$ over~$\S$ is a left
$\C$\contramodule{} endowed with a \emph{left semicontraaction}
map $\bp_\bP\:\bP\rarrow\Cohom_\C(\S,\bP)$, which has to be a morphism
of left $\C$\contramodule s satisfying the associativity and
unity equations $\Cohom(\id_\S,\bp_\bP)\circ\bp_\bP =
\Cohom(\bm_\S,\id_\bP)\circ\bp_\bP$ and
$\Cohom(\be_\S,\id_\bP)\circ\bp_\bP = \id_\bP$.
 
 For example, if the coalgebra $\C$ coincides with the ground
field~$k$, and $\S=R$ is just a $k$\+algebra, then left
$\S$\semicontramodule s are simply left $R$\module s.

 For any right $\S$\semimodule{} $\bN$ and any $k$\+vector space~$V$
the space $\Hom_k(\bN,V)$ has a natural structure of left
$\S$\semicontramodule.
 For any left $\C$\contramodule{} $\Q$, the space of
cohomomorphisms $\Cohom_\C(\S,\Q)$ has a natural structure
of left semicontramodule.
 It is called the left $\S$\semicontramodule{} \emph{coinduced}
from a left $\C$\contramodule{} $\Q$.
 The space of semicontramodule homomorphisms into the coinduced
semicontramodule is described by the formula
$\Hom^\S(\bP,\Cohom_\C(\S,\Q))\simeq\Hom^\C(\bP,\Q)$.
 We will denote the category of left $\S$\semicontramodule s by
$\S\sicntr$.
 The category of left $\S$\semicontramodule s is abelian provided
that $\S$ is an injective left $\C$\comodule.
 Moreover, $\S$ is an injective left $\C$\comodule{} if and only if
the category $\S\sicntr$ is abelian and the forgetful functor
$\S\sicntr\rarrow\C\contra$ is exact.

 The functor $\Cohom_\C^\op\:\C\comodl\times\C\contra^\op\rarrow
k\vect^\op$ is a pairing compatible with the left module
category structure on $\C\comodl$ and the right module
category structure on $\C\contra^\op$ over the tensor category
$\C\bcomod\C$.
 Thus one can define the \emph{space of semihomomorphisms}
$\SemiHom_\S(\bM,\bP)$ from a left $\S$\semimodule{} $\bM$
to a left $\S$\semicontramodule{} $\bP$ as the kernel of
the pair of maps
\begin{multline*}
 (\Cohom(\bn_\bM,\id_\bP),\Cohom(\id_\bM,\bp_\bP))\: \\
 \Cohom_\C(\bM,\bP)\birarrow 
 \Cohom_\C(\S\oc_\C\bM\;\bP) = \Cohom_\C(\bM,\Cohom_\C(\S,\bP)).
\end{multline*}
 There are natural isomorphisms $\SemiHom_\S(\S\oc_\C\L\;\bP)\simeq
\Cohom_\C(\L,\bP)$ and $\SemiHom_\S(\bM,\Cohom_\C(\S,\Q))\simeq
\Cohom_\C(\bM,\Q)$.
 The functor of semihomomorphisms is neither left, nor right exact.

\subsubsection{}
 The complex of semihomomorphisms $\SemiHom_\S(\bM^\bu,\bP^\bu)$ from
a complex of left $\S$\semimodule s $\bM^\bu$ to a complex of
left $\S$\semicontramodule s $\bP^\bu$ is defined as the total complex
of the bicomplex $\SemiHom_\S(\bM^i,\bP^j)$, constructed by taking
infinite products along the diagonals.
 Assume that $\S$ is an injective left and right $\C$\comodule.
 Let us define the double-sided derived functor $\SemiExt_\S$ of
the functor of semihomomorphisms.

 The \emph{semiderived category} $\sD^\si(\S\sicntr)$ of left
$\S$\semicontramodule s is defined as the quotient category of
the homotopy category $\Hot(\S\sicntr)$ by the thick subcategory
$\Acycl^{\ctrd\C}(\S\sicntr)$ of complexes of $\S$\semicontramodule s
that are \emph{contraacyclic as complexes of\/ $\C$\contramodule s}.

 A complex of left $\S$\semimodule s $\bM^\bu$ is called
\emph{semiprojective} if the complex $\SemiHom_\S(\bM^\bu,\bP^\bu)$
is acyclic for any $\C$\contraacyclic{} complex of left
$\S$\semicontramodule s $\bP^\bu$.
 A complex of left $\S$\semicontramodule s $\bP^\bu$ is called
\emph{semiinjective} if the complex $\SemiHom_\S(\bM^\bu,\bP^\bu)$
is acyclic for any $\C$\coacyclic{} complex of left
$\S$\semimodule s $\bM^\bu$.
 For example, the complex of $\S$\semimodule s $\S\oc_\C\L^\bu$ induced
from a complex of injective $\C$\comodule s $\L^\bu$ is semiprojective.
 Any semiprojective complex of semimodules is semiflat.
 The complex of $\S$\semicontramodule s $\Cohom_\C(\S,\Q^\bu)$
coinduced from a complex of projective $\C$\contramodule s $\Q^\bu$
is semiinjective.

 The quotient category of the homotopy category $\Hot_\sipr(\S\simodl)$
of semiprojective complexes of $\S$\semimodule s by the thick
subcategory of $\C$\coacyclic{} semiprojective complexes
$\Acycl^{\cod\C}(\S\simodl)\cap \Hot_\sipr(\S\simodl)$
is equivalent to the semiderived category of $\S$\semimodule s.
 Analogously, the quotient category of the homotopy category
$\Hot_\siin(\S\sicntr)$ of semiinjective complexes of
$\S$\semicontramodule s by the thick subcategory of
$\C$\contraacyclic{} semiinjective complexes
$\Acycl^{\ctrd\C}(\S\sicntr)\cap\Hot_\siin(\S\sicntr)$ is equivalent to
the semiderived category of $\S$\semicontramodule s.
 The derived functor
$$
 \SemiExt_\S\: \sD^\si(\S\simodl)^\op\times\sD^\si(\S\sicntr)
 \lrarrow \sD(k\vect)
$$
is defined by restricting the functor of semihomomorphisms
to either of the full subcategories $\Hot_\sipr(\S\simodl)^\op\times
\Hot(\S\sicntr)$ or $\Hot(\S\simodl)^\op\times\Hot_\siin(\S\sicntr)$
of the category $\Hot(\S\simodl)^\op\times\Hot(\S\sicntr)$.

\subsubsection{}
 Assume that $\S$ is an injective left and right $\C$\comodule.
 One can check that the adjoint functors $\Psi_\C\:\C\comodl\rarrow
\C\contra$ and $\Phi_\C\:\C\contra\rarrow\C\comodl$ transform
left $\C$\comodule s with an $\S$\semimodule{} structure into
left $\C$\contramodule s with an $\S$\semicontramodule{} structure
and vice versa.
 Therefore, there is a pair of adjoint functors $\Psi_\S\:\S\simodl
\rarrow\S\sicntr$ and $\Phi_\S\:\S\sicntr\rarrow\S\simodl$ agreeing
with the functors $\Psi_\C$ and~$\Phi_\C$ and providing
an equivalence between the exact categories of $\C$\+injective left
$\S$\semimodule s and $\C$\+projective left $\S$\semicontramodule s.
 
 To construct this pair of adjoint functors in a natural way,
let us define the operation of contratensor product of
a semimodule and a semicontramodule.

 Let $\bN$ be a right $\S$\semimodule{} and $\bP$ be a left
$\S$\semicontramodule.
 The \emph{contratensor product} $\bN\Ocn_\S\bP$ is defined
as the cokernel of the pair of maps
 $$
  (\bn_\bN\ocn\id_\bP\;\eta_\S\circ(\id_{\bN\.\suboc_\C\S}\ocn\bp_\bP))
  \:(\bN\oc_\C\S)\ocn_\C\bP \birarrow \bN\ocn_\C\bP
 $$
where the natural ``evaluation'' map
$\eta_\K\:(\N\oc_\C\K)\ocn_\D\Cohom_\C(\K,\P)\rarrow\N\ocn_\C\P$
exists for any right $\C$\comodule{} $\N$, left
$\C$\contramodule{} $\P$, and $\C$\+$\D$\bicomodule{} $\K$
and is dual to the map
\begin{setlength}{\multlinegap}{0pt}
\begin{multline*}
  \Hom_k(\eta_\K,V) = \Cohom_\C(\K,{-})\: \\
  \Hom^\C(\P,\Hom_k(\N,V))
  \lrarrow \Hom^\D(\Cohom_\C(\K,\P),\Cohom_\C(\K,\Hom_k(\N,V)))
\end{multline*}
for any $k$\+vector space~$V$.
 There are natural isomorphisms $(\R\oc_\C\S)\Ocn_\S\bP\simeq
\R\ocn_\C\bP$ and $\Hom_k(\bN\Ocn_\S\bP\;V)\simeq
\Hom^\S(\bP,\Hom_k(\bN,V))$.
 The functor of contratensor product over~$\S$ is right exact
whenever $\S$ is an injective left $\C$\comodule.

 The adjoint functors $\Psi_\S$ and $\Phi_\S$ can be defined by
the formulas $\Psi_\S(\bM)=\Hom_\S(\S,\bM)$ and $\Phi_\S(\bP)=
\S\Ocn_\S\bP$.
 Here the space $\Hom_\S(\S,\bM)$ is endowed with a left
$\S$\semicontramodule{} structure as a subsemicontramodule of
the semicontramodule $\Hom_k(\S,\bM)$, while the space
$\S\Ocn_\S\bP$ is endowed with a left $\S$\semimodule{} structure
as a quotient semimodule of the semimodule $\S\ot_k\bP$.
\end{setlength}

 The quotient category of the homotopy category of $\C$\injective{}
$\S$\semimodule s $\Hot(\S\simodl_{\injd\C})$ by the thick subcategory
of $\C$\contractible{} complexes of $\C$\injective{} $\S$\semimodule s
is equvalent to the semiderived category of $\S$\semimodule s.
 Analogously, the quotient category of the homotopy category 
$\Hot(\S\sicntr_{\projd\C})$ of $\C$\+projective{}
$\S$\semicontramodule s by the thick subcategory of $\C$\contractible{}
complexes of $\C$\projective{} $\S$\semicontramodule s is equivalent to
the semiderived category of $\S$\semicontramodule s.
 Thus the semiderived categories of left $\S$\semimodule s and
left $\S$\semicontramodule s are equivalent, $\sD^\si(\S\simodl)
\simeq\sD^\si(\S\sicntr)$.

 When $\S$ is not an injective left or right $\C$\comodule,
the exact categories of $\C$\+injective{} $\S$\semimodule s and
$\C$\projective{} $\S$\semicontramodule s are still equivalent,
even though the functors $\Psi_\S$ and $\Phi_\S$ are not defined
on the whole categories of all comodules and contramodules.

\subsubsection{}
  The functor $\Ext_\S:\sD^\si(\S\simodl)^\op\times \sD^\si(\S\simodl)
\rarrow\sD(k\vect)$ of homomorphisms in the semiderived category
$\sD^\si(\S\simodl)$ can be computed by restricting the functor
$\Hom_\S\:\Hot(\S\simodl)^\op\times\Hot(\S\simodl)\rarrow\Hot(k\vect)$
of homomorphisms in the homotopy category $\Hot(\S\simodl)$ to
an appropriate subcategory of the Carthesian product
$\Hot(\S\simodl)^\op\times\Hot(\S\simodl)$.
 Namely, a complex of left $\S$\semimodule s $\bL^\bu$ is called
\emph{projective relative to\/~$\C$} ($\S/\C$\+projective)
if the complex $\Hom_\S(\bL^\bu,\bM^\bu)$ is acyclic for any
$\C$\contractible{} complex of $\C$\injective{}
left $\S$\semimodule s~$\bM^\bu$.
 For example, the complex of $\S$\semimodule s $\S\oc_\C\L^\bu$
induced from a complex of $\C$\comodule s~$\L^\bu$ is projective
relative to~$\C$.
 The quotient category of the homotopy category
$\Hot_{\projd\S/\C}(\S\simodl)$ of $\S/\C$\projective{} complexes
of $\S$\semimodule s by the thick subcategory
$\Acycl^{\cod\C}(\S\simodl)\cap\Hot_{\projd\S/\C}(\S\simodl)$
of $\C$\coacyclic{} $\S/\C$\projective{} complexes is equivalent
to the semiderived category of $\S$\semimodule s.
 The functor $\Ext_\S$ can be obtained by restricting the
functor $\Hom_\S$ to the full subcategory
$\Hot_{\projd\S/\C}(\S\simodl)^\op\times\Hot(\S\simodl_{\injd\C})$
of the category $\Hot(\S\simodl)^\op\times\Hot(\S\simodl)$.

 Analogously, the functor $\Ext^\S\:\sD^\si(\S\sicntr)^\op\times
\sD^\si(\S\sicntr)\rarrow\sD(k\vect)$ of homomorphisms in
the semiderived category $\sD^\si(\S\sicntr)$ can be computed by
restricting the functor $\Hom^\S\:\Hot(\S\sicntr)^\op\times
\Hot(\S\sicntr)\rarrow\Hot(k\vect)$ to an appropriate subcategory of
the Carthesian product $\Hot(\S\sicntr)^\op\times\Hot(\S\sicntr)$.
 A complex of $\S$\semicontramodule s $\bQ^\bu$ is called
\emph{injective relative to\/~$\C$} ($\S/\C$\+injective)
if the complex $\Hom^\S(\bP^\bu,\bQ^\bu)$ is acyclic for any
$\C$\contractible{} complex of $\C$\projective{}
$\S$\semicontramodule s~$\bP^\bu$.
 For example, the complex of $\S$\semicontramodule s
$\Cohom_\C(\S,\Q^\bu)$ coinduced from a complex of
$\C$\contramodule s~$\Q^\bu$ is %injective relative to~$\C$.
$\S/\C$\injective.
 The quotient category of the homotopy category
$\Hot_{\injd\S/\C}(\S\sicntr)$ of $\S/\C$\injective{} complexes
of $\S$\semicontramodule s by the thick subcategory
$\Acycl^{\ctrd\C}(\S\sicntr)\cap\Hot_{\injd\S/\C}(\S\sicntr)$
of $\C$\contraacyclic{} $\S/\C$\injective{} complexes is equivalent
to the semiderived category of $\S$\semicontramodule s.
 The functor $\Ext^\S$ can be obtained by restricting the
functor $\Hom^\S$ to the full subcategory
$\Hot(\S\sicntr_{\projd\C})^\op\times\Hot_{\injd\S/\C}(\S\sicntr)$
of the category $\Hot(\S\sicntr)^\op\times\Hot(\S\sicntr)$.

 The contratensor product $\bN^\bu\Ocn_\S\bP^\bu$ of a complex
of right $\S$\semimodule s $\bN^\bu$ and a complex of left
$\S$\semicontramodule s $\bP^\bu$ is defined as the total complex
of the bicomplex $\bN^i\Ocn_\S\bP^j$, constructed by taking
infinite direct sums along the diagonals.
 Let us define the left derived functor $\CtrTor^\S$ of the functor
of contratensor product over~$\S$.
 A complex of right $\S$\semimodule s $\bN^\bu$ is called
\emph{contraflat relative to\/~$\C$} ($\S/\C$\+contraflat) if
the complex $\bN^\bu\Ocn_\S\bP^\bu$ is acyclic for any
$\C$\contractible{} complex of $\C$\projective{}
$\S$\semicontramodule s~$\bP^\bu$.
  For example, the complex of $\S$\semimodule s $\R^\bu\oc_\C\S$
induced from a complex of right $\C$\comodule s~$\R^\bu$ is
contraflat relative to~$\C$.
 A complex of right $\S$\semimodule s $\bN^\bu$ is contraflat
relative to~$\C$ if and only if the complex of left
$\S$\semimodule s $\Hom_k(\bN^\bu,k)$ is injective relative to~$\C$.
 The quotient category of the homotopy category
$\Hot_{\ctrfld\S/\C}(\simodr\C)$ of $\S/\C$\+contraflat complexes
of right $\S$\semimodule s by the thick subcategory
$\Acycl^{\cod\C}(\simodr\S)\cap\Hot_{\ctrfld\S/\C}(\simodr\C)$
of $\C$\coacyclic{} $\S/\C$\+contraflat complexes is equivalent
to the semiderived category of right $\S$\semimodule s.
 The left derived functor
$$
 \CtrTor^\S\:\sD^\si(\simodr\S)\times\sD^\si(\S\sicntr)
 \lrarrow\sD(k\vect)
$$
is defined by restricting the functor of contratensor product
to the full subcategory $\Hot_{\ctrfld\S/\C}(\simodr\S)\times
\Hot(\S\sicntr_{\projd\C})$ of the category
$\Hot(\simodr\S)\times \Hot(\S\sicntr)$.

 The equivalence of triangulated categories $\sD^\si(\S\simodl)\simeq
\sD^\si(\S\sicntr)$ transforms the double-sided derived functor
$\SemiExt_\S$ into the functor $\Ext$ in either of the semiderived
categories and the double-sided derived functor $\SemiTor^\S$
into the left derived functor $\CtrTor^\S$.

\subsubsection{}
 Any semiprojective complex of $\S$\semimodule s is $\S/\C$\projective.
 An $\S/\C$\projective{} complex of $\C$\injective{} $\S$\semimodule s
is semiprojective.
 The homotopy category of semiprojective complexes of $\C$\injective{}
$\S$\semimodule s is equivalent to the semiderived category of
$\S$\semimodule s.

 Analogously, any semiinjective complex of $\S$\semicontramodule s is
$\S/\C$\injective.
 An $\S/\C$\injective{} complex of $\C$\projective{}
$\S$\semicontramodule s is semiinjective.
 The homotopy category of semiinjective complexes of $\C$\injective{}
$\S$\semicontramodule s is equivalent to the semiderived category of
$\S$\semicontramodule s.

 Our definitions of $\S/\C$\projective{} and $\S/\C$\injective{}
complexes differ from the traditional ones;
cf.~\ref{semijective-complexes} and Remark~\ref{semi-model-struct}.1.

\subsection{Nonhomogeneous Koszul duality over a base ring}
\label{koszul-over-ring}
 This subsection is intended to supply preliminary material for
Section~\ref{nonhom-koszul-section} and Appendix~\ref{tate-appendix}.

\subsubsection{}
 A graded ring $S=S_0\oplus S_1\oplus S_2\oplus\dsb$ is called
\emph{quadratic} if it is generated by $S_1$ over $S_0$ with
relations of degree~$2$ only. 
 In other words, this means that if one considers the graded ring
freely generated by the $S_0$\+$S_0$\bimodule{} $S_1$
(the ``tensor ring'' of the $S_0$\+$S_0$\bimodule{} $S_1$), i.~e.,
the graded ring $\boT_{S_0,S_1}$ with components $S_1^{\ot_{S_0}\.n}=
S_1\ot_{S_0} S_1\ot_{S_0}\dsb\ot_{S_0}S_1$, then the ring $S$ should be
isomorphic to the quotient ring of $\boT_{S_0,S_1}$ by the ideal
generated by a certain subbimodule $I_S$ in $S_1\ot_{S_0}S_1$.

 A quadratic ring $S$ is called \emph{$2$\+left finitely projective}
if both left $S_0$\module s $S_1$ and $S_2$ are projective and
finitely generated.
 A quadratic ring is called \emph{$3$\+left finitely projective} if
the same applies to $S_1$, \ $S_2$, and $S_3$.
 Further conditions of this kind are not very sensible to consider
for general quadratic rings.
 Analogously one defines \emph{$2$\+right finitely projective} and
\emph{$3$\+right finitely projective} quadratic rings.

 There is an anti-equivalence between the category of $2$\+left
finitely projective quadratic rings and the category of $2$\+right
finitely projective quadratic rings, called the \emph{quadratic duality}.
 The duality functors are defined by the formulas $R_0=S_0$, \ 
$R_1=\Hom_{S_0}(S_1,S_0)$, \ $R_2=\Hom_{S_0}(I_S,S_0)$, \
$I_R\simeq\Hom_{S_0}(S_2,S_0)$, and conversely, $S_1=\Hom_{R_0^\rop}(R_1,R_0)$,
\ $S_2=\Hom_{R_0^\rop}(I_R,R_0)$, \ $I_S\simeq\Hom_{R_0^\rop}(R_2,R_0)$.
 Here we use the natural isomorphism
$$
 \Hom_{S_0}(N,S_0)\ot_{S_0}\Hom_{S_0}(M,S_0)\simeq
 \Hom_{S_0}(M\ot_{S_0} N\;S_0)
$$ for $S_0$\+$S_0$\bimodule s $M$ and $N$
that are projective and finitely generated left $S_0$\module s,
and the analogous isomorphism
$$
 \Hom_{R_0^\rop}(N,R_0)\ot_{R_0}\Hom_{R_0^\rop}(M,R_0)\simeq
 \Hom_{R_0^\rop}(M\ot_{R_0} N\;R_0)
$$
for $R_0$\+$R_0$\bimodule s $M$ and $N$ that are projective and finitely
generated right $R_0$\module s.

 The duality functor sends $3$\+left finitely projective quadratic rings
to $3$\+right finitely projective quadratic rings and vice versa.
 Indeed, set $J_S=I_S\ot_{S_0}S_1\cap S_1\ot_{S_0}I_S\subset
S_1\ot_{S_0}S_1\ot_{S_0}S_1$; then
 $$
  0\lrarrow I_S\lrarrow I_A\ot_{S_0}S_1\oplus S_1\ot_{S_0}I_A
  \lrarrow S_1\ot_{S_0} S_1\ot_{S_0} S_1\lrarrow S_3\lrarrow0
 $$
is an exact sequence of finitely generated projective left
$S_0$\module s, and $R_3\simeq\Hom_{S_0}(J_S,S_0)$, since
the sequence remains exact after applying $\Hom_{S_0}({-},S_0)$.

\subsubsection{}
 A graded ring $S=S_0\oplus S_1\oplus S_2\oplus\dsb$ is called
\emph{left flat Koszul} if it is flat as a left $S_0$\module{}
and one has $\Tor^S_{i,j}(S_0,S_0)=0$ for $i\ne j$.
 Here $S_0$ is endowed with the right and left $S$\module{}
structures via the augmentation map $S\rarrow S_0$ and the second
grading $j$ on the $\Tor$ is induced by the grading of~$S$.
 \emph{Right flat Koszul} graded rings are defined in the analogous
way.
 A left/right flat Koszul ring is called \emph{left}/\emph{right}
(\emph{finitely}) \emph{projective Koszul}, if it is a projective
(with finitely generated grading components) left/right $S_0$\module.

 Notice that when $S$ is a flat left $S_0$\module, the reduced
relative bar construction
$$
 \dsb\lrarrow S\ot_{S_0}S/S_0\ot S/S_0\lrarrow S\ot S/S_0\lrarrow S
$$
is a flat resolution of the left $S$\module{} $S_0$, so one can
use it to compute $\Tor^S(S_0,S_0)$.
 When $S$ is a projective left $S$\module, the same resolution
can be used to compute $\Ext_S(S_0,S_0)$.
 Assume that the grading components of $S$ are finitely generated
projective left $S_0$\module s; then it follows that $S$ is
left finitely projective Koszul if and only if
$\Ext_S^{i,j}(S_0,S_0)=0$ for $i\ne j$ and $\Ext_S^{i,i}(S_0,S_0)$ are
projective right $S_0$\module s.

 Assume that a graded ring $S$ is a flat left $S_0$\module.
 Then $S$ is left flat Koszul if and only if it is quadratic and
for each degree~$n$ the lattice of subbimodules in $S_1^{\ot_{S_0}\.n}$
generated by the $n-1$ subbimodules $S_1^{\ot_{S_0}\.i-1}\ot_{S_0}I_S
\ot_{S_0}S_1^{\ot_{S_0}\.n-i-1}$ is distributive.
 This means that for any three subbimodules $X$, $Y$, $Z$ that can be
obtained from the generating subbimodules by applying the operations
of sum and intersection one should have $(X+Y)\cap Z = X\cap Z +
Y\cap Z$.
 Furthermore, if $S$ is a left finitely projective Koszul ring,
then the ring $R$ quadratic dual to $S$ is right finitely projective 
Koszul, and vise versa; besides, in this case the graded ring
$\Ext_S(S_0,S_0)$ is isomorphic to $R^\rop$ and the graded ring
$\Ext_{R^\rop}(R_0,R_0)$ is isomorphic to~$S$.

\subsubsection{}    \label{cdg-rings}
 Let $S$ be a $3$\+left finitely projective quadratic ring.
 Suppose that we are given a ring $S\til$ endowed with an increasing
filtration $F_0S\til\subset F_1S\til\subset F_2S\til\subset\dsb$
such that $S=\bigcup_n F_nS\til$ and the associated graded ring
$\gr_FS\til$ is identified with $S$.
 Such a ring $S\til$ will be called a \emph{$3$\+left finitely
projective nonhomogeneous quadratic ring}.
 If the graded ring $S$ is left finitely projective Koszul,
the filtered ring $S\til$ is called a \emph{left finitely
projective nonhomogeneous Koszul ring}.

 Let $R$ be the $3$\+right finitely projective quadratic ring
dual to~$S$.
 We would like to describe the additional structure on the ring
$R$ corresponding to the data of a filtered ring $S\til$ endowed
with an isomorphism $\gr_FS\til\simeq S$.

 A \emph{CDG\+ring} (\emph{curved differential graded ring}) is
a graded ring $R=\bigoplus_n R^n$ endowed with an odd derivation $d$
of degree~$1$ and a ``curvature element'' $h\in R^2$ such that
$d^2(x)=[h,x]$ for all $x\in R$ and $d(h)=0$.
 A morphism of CDG\+rings ${}'\!R\rarrow{}''\!R$ is a pair $(f,a)$,
where $f\:{}'\!R\rarrow{}''\!R$ is a morphism of graded rings and
$a$ is an element in ${}''\!R^1$ such that $f(d'(x))=d''(f(x))+
[a,f(x)]$ (the supercommutator) for all $x\in{}'\!R$ and $f(h')=
h''+d''(a)+a^2$.
 Composition of morphisms is defined by the rule $(g,b)(f,a)=
(gf\;b+g(a))$.
 Identity morphisms are the morphisms $(\id,0)$.

 So the \emph{category of CDG\+rings} is defined.
 Notice that the natural functor from the category of DG\+rings
to the category of CDG\+rings is faithful, but not fully faithful.
 In other words, two DG\+rings may be isomorphic in the category of
CDG\+rings without being isomorphic as DG\+rings.
 Furthermore, two CDG\+rings of the form $(R\; d+[a,\cdot]\; h+d(a)+a^2)$
and $(R,d,h)$ are always naturally isomorphic, the isomorphism being
given by the pair $(\id,a)$.

 There is a fully faithful contravariant functor from the category of
$3$\+left finitely projective nonhomogeneous quadratic rings $S\til$
with a fixed ring $F_0S\til$ to the category of CDG\+rings $(R,d,h)$
with the same component $R^0=F_0S\til$ such that the underlying
graded ring $R$ of the CDG\+ring $(R,d,h)$ corresponding to~$S\til$
is the $3$\+right finitely projective quadratic ring dual to
the ring $S=\gr_FS\til$ (in the grading $R_i=R^i$).
 
 This functor is constructed as follows.
 For each $3$\+left finitely projective nonhomogeneous quadratic ring
$S\til$ choose a complementary left $S_0=F_1S\til$\+submodule $V$ to
the submodule $F_0S\til$ in the left $S_0$\module{} $F_1S\til$.
 This can be done, because the quotient module $S_1=F_1S\til/F_0S\til$
is projective.
 Since $V$ maps isomorphically to $S_1=F_0S\til/F_1S\til$, it is
endowed with a structure of an $S_0$\+$S_0$\bimodule.
 The embedding $V\rarrow F_0S\til$ is only a morphism of left
$S_0$\module s, however; the right actions of $S_0$ in $V$ and
$F_1S\til$ are compatible modulo $F_0S\til$.
 Put $q(v,s)=m(v,s)-vs$ for $v\in V$, $s\in S_0$, where $m(v,s)$
is the product in $S\til$ and $vs$ denotes the right action of
$S_0$ in $V$.
 This defines a map $q\:V\ot_\boZ S_0\rarrow S_0$.
 
 Let $I\til$ be the full preimage of the subbimodule
$I_S\subset S_1\ot_{S_0}S_1$ under the surjective map
$S_1\ot_\boZ S_1\rarrow S_1\ot_{S_0}S_1$.
 Using the identification of $V$ with $S_1$, we will consider $I\til$
as the full preimage of $F_1S\til$ under the multiplication map
$m\: V\ot_\boZ V\rarrow S\til$.
 Let us split the map $m\: I\til\rarrow F_1S\til$ into two components
$(g,-h)$ according to the direct sum decomposition $F_1S\til\simeq
V\oplus S_0$, so that $g\:I\til\rarrow V$ and $h\:I\til\rarrow S_0$.

 The differentials $d_0\: R^0\rarrow R^1$ and $d_1\: R^1\rarrow R^2$
are defined in terms of the maps $q$ and~$g$ by the formulas
$$
 \langle v, d_0(s) \rangle = q(v,s), \qquad
 \langle i, d_1(r) \rangle = \langle g(\tilde{\imath}), r \rangle
 - q(\tilde{\imath}_1, \langle \tilde{\imath}_2,r\rangle),
$$
where $\langle\,\.,\,\rangle$ denotes the pairing of $V$ with $R^1$
and of $I_S$ with $R^2$, and $\tilde{\imath}$ is any preimage of $i$
in $I\til$, written also as $\tilde{\imath}=
\tilde{\imath}_1\ot\tilde{\imath}_2$. 
 The map $h$ factorizes through $I_S$, providing the curvature
element in $R^2=\Hom_{S_0}(I_S,S_0)$.

 Finally, to a morphism of nonhomogeneous quadratic rings
$f\:S''{}\til\rarrow S'{}\til$ with chosen complementary submodules
$V''\subset F''_1S''{}\til$ and $V'\subset F'_1S'{}\til$ one
assigns a morphism of dual CDG\+rings $(g,a)\:({}'\!R,d',h')\rarrow
({}''\!R,d'',h'')$ defined as follows.
 The morphism of quadratic rings $g\:{}'\!R\rarrow{}''\!R$ is
the quadratic dual map to the associated graded morphism
$\gr\,f\:S''\rarrow S'$, while the element $a\in{}''\!R^1=
\Hom_{S_0}(S_1'',S_0)$ is equal to minus the composition
$V''\rarrow F_1''S''{}\til\rarrow F_1'S'{}\til\rarrow S_0$ of
the embedding $V''\rarrow F_1''S''{}\til$, the map $f$, and
the projection $F_1'S'{}\til\rarrow S_0$ along $V'$.
 In particular, for a given nonhomogeneous quadratic ring $S\til$
changing the splitting of $F_1S\til$ by the rule
$V''=\{\.v'-a(v')\mid v'\in V'\.\}$ leads to a natural morphism of
CDG\+rings $(\id,a)\:(R,d',h')\rarrow(R,d'',h'')$.

 One has to make quite some computations in order to check that
everything is well-defined and compatible in this construction.
 In particular, the $3$\+left projectivity is actually used in
the form of the duality between $J_S$ (where some self-consistency
equations on the defining relations of $S\til$ live) and $R^3$
(where the equations $d(e)=0$ for $e\in I_R$, \ $d^2(r)=[h,r]$, 
and $d(h)=0$ have to be verified).

 The nonhomogeneous quadratic duality functor restricted to
the categories of left finitely projective nonhomogeneous Koszul
rings and right finitely projective Koszul CDG\+rings becomes
an equivalence of categories.
 In other words, any CDG\+ring whose underlying graded
ring is right finitely projective Koszul corresponds to a left
finitely projective nonhomogeneous Koszul ring.
 This is the statement of the Poincare--Birkhoff--Witt theorem
for finitely projective nonhomogeneous Koszul rings.

\subsubsection{}  \label{quasi-differential-rings}
 A \emph{quasi-differential ring} $R\til$ is a graded ring
$R\til=\bigoplus_n R^n{}\til$ endowed with an odd derivation~$\d$
of degree~$-1$ with zero square such that the cohomology of~$\d$
vanish (equivalently, the unit element of $R\til$ lies in
the image of~$\d$).
 A \emph{quasi-differential structure} on a graded ring~$R$ is
the data of a quasi-differential ring $(R\til,\d)$ together with
an isomorphism of graded rings $\ker\d\simeq R$.

 The category of quasi-differential rings is equivalent to the category
of CDG\+rings.
  This equivalence assigns to a CDG\+ring $(R,d,h)$
the quasi-differential ring $R\til=R[\delta]$ with an added
generator~$\delta$ of degree~$1$, the relations $[\delta,x]=d(x)$
(the supercommutator) for $x\in R$ and $\delta^2=h$, and the derivation
$\d=\d/\d\delta$ (the partial derivative in~$\delta$, meaning the unique
odd derivation $\d$ of $R\til$ for which $\d(R)=0$ and $\d(\delta)=1$).
 Conversely, to construct a CDG\+ring structure on the kernel $R$
of the derivation $\d$ of a quasi-differential ring $R\til$,
it suffices to choose an element $\delta\in R^1{}\til$ such that
$\d(\delta)=1$ and set $d(x)=[\delta,x]$, \ $h=\delta^2$.
 Choosing two different elements $\delta$ leads to two naturally
isomorphic CDG\+rings.

 A \emph{left CDG\module} $M$ over a CDG\+ring $(R,d,h)$ is
a graded left $R$\module{} endowed with a $d$\+derivation $d_M$
(that is a homogeneous map $M\rarrow M$ of degree~$1$ for which
$d_M(rx)=d(r)x+(-1)^{|r|}rd(x)$ for $r\in R$, \ $x\in M$, where
$|r|$ denotes the degree of a homogeneous element $r$) such
that $d_M^2(x)=hx$.
 A \emph{quasi-differential left module} over a quasi-differential
ring $R\til$ is just a graded left $R\til$\module{} (without any
differential).
 The category of left CDG\module s over a CDG\+ring $(R,d,h)$ is 
isomorphic to the category of quasi-differential left modules over
the quasi-differential ring $R\til$ corresponding to $(R,d,h)$;
this isomorphism of categories assigns to a graded $R\til$\module{}
structure on a graded left $R$\module{} $M$ the derivation
$d_M(x)=\delta x$ on~$M$.

 Analogously, a \emph{right CDG\module} $N$ over $(R,d,h)$ is
a graded right $R$\module{} endowed with a $d$\+derivation $d_N$
(that is a homogeneous map $N\rarrow N$ of degree~$1$ for which
$d_N(xr)=d_N(x)r+(-1)^{|x|}xd(r)$ for $x\in N$, \ $r\in R$) such
that $d_N^2(x)=-xh$.
 A \emph{quasi-differential right module} over a quasi-differential
ring $R\til$ is just a graded left $R\til$\module.
 The category of right CDG\module s over $(R,d,h)$ is isomorphic
to the category of quasi-differential $R\til$\module s when $R\til$
corresponds to $(R,d,h)$; this isomorphism of categories assigns
to a graded $R\til$\module{} structure on a graded right $R$\module{}
$N$ the derivation $d_N(x)=(-1)^{|x|+1}x\delta$ on~$N$.

\subsubsection{}
 CDG\module s over a CDG\+ring form a \emph{DG\category}, i.~e.,
a category where for any two given objects there is a complex of
morphisms between them.
 We will consider the cases of left and right CDG\module s
separately.

 Let $L$ and $M$ be two left CDG\module s over a CDG\+ring $(R,d,h)$.
 The complex $\Hom_R^\bu(L,M)$ is defined as follows.
 The component $\Hom_R^n(L,M)$ consists of all homogeneous maps
$L\rarrow M$ of degree $n$ supercommuting with the $R$\module{}
structures in $L$ and~$M$.
 This means that for $f\in\Hom_R^n(L,M)$ and $r\in R$, \ $x\in L$
one should have $f(rx)={-1}^{n|r|}rf(x)$.
 The differential is defined by the formula $(df)(x)=d_Mf(x)-
(-1)^{|f|}fd_L(x)$.
 One has $d^2(f)=0$, because $f(hx)=hf(x)$.

 Let $K$ and $N$ be two right CDG\module s over $(R,d,h)$.
 The component $\Hom_R^n(K,N)$ of the complex $\Hom_R^\bu(K,N)$
consists of all homogeneous maps $K\rarrow N$ of degree $n$
commuting with the $R$\module{} structures in $L$ and $M$
(without any signs).
 The differential is defined by the formula $(df)(x)=d_Nf(x)-
(-1)^{|f|}fd_K(x)$.

 One can see that shifts and cones exist in the DG\+categories of
(left or right) CDG\module s, and moreover, a CDG\module{} structure
can be twisted with any cochain in the complex of endomorphisms
satisfying the Maurer--Cartan equation~\cite{BK}.
 It follows that the homotopy categories of CDG\module s, defined
as the categories of zero cohomology of the DG\+categories of
CDG\module s, are triangulated.

 Furthermore, one can speak about the total CDG\module s of complexes
of CDG\module s, constructed by taking infinite direct sums or infinite
products along the diagonals.
 In particular, there are total CDG\module s of exact triples of
CDG\module s.
 This allows one to define the \emph{coderived} and \emph{contraderived
categories of CDG\module s} over $(R,d,h)$ as the quotient categories
of the homotopy categories of CDG\module s by the minimal triangulated
subcategories containing the total CDG\module s of exact triples of
CDG\module s and closed under infinite direct sums and infinite
products, respectively.

 Notice that one \emph{cannot} define the conventional derived
category of CDG\module s, as CDG\module s don't have any
cohomology groups.

\subsubsection{}
 Let $S\til$ be a left finitely projective nonhomogeneous Koszul
ring and $(R,d,h)$ be the dual CDG\+ring.
 Assume that the ring $S_0$ has a finite right homological dimension.
 Then the Koszul duality theorem claims that the derived category
of right $S\til$\module s is equivalent to the coderived category
of right CDG\module s $N$ over $(R,d,h)$ such that every element
of $N$ is annihilated by $R^n$ for $n\gg0$.
 Assuming that $S_0$ has a finite left homological dimension, 
the derived category of left $S\til$\module s is also described
as being equivalent to the contraderived category of left
CDG\module s over $(R,d,h)$ in which certain infinite summation
operations are defined.

 One can drop the homological dimension assumptions, replacing
the derived categories of $S\til$\module s in the formulations
of these results with certain semiderived categories relative
to~$S_0$ (see Theorem~\ref{koszul-duality-theorem} and
Remark~\ref{quasi-differential-contramodules}).
 And the conventional derived category of right $S\til$\module s
without the homological dimension assumption on~$S_0$ is equivalent
to the quotient category of the coderived category of locally
nilpotent (in the above sense) right CDG\module s over $(R,d,h)$
by its minimal triangulated subcategory closed under infinite
direct sums and containing all the CDG\module s $N$, where $R^n$
act by zero for all $n>0$ and which are acyclic with respect
to~$d_N$ (one has $d_N^2=0$, since $Nh=0$).
 The latter result has an obvious analogue in the case of left
CDG\+modules with infinite summation operations.

\subsubsection{}
 The following example is thematic.
 Let $M$ be a smooth affine algebraic variety and $E$ be
a vector bundle over~$M$.
 Let $\Diff_{M,E}$ denote the ring of differential operators acting
in the sections of~$E$.
 The natural filtration of $\Diff_{M,E}$ by the order of differential
operators makes it a left (and right) finitely projective
nonhomogeneous Koszul ring.
 To construct the dual CDG\+ring, choose a global connection
$\nabla_E$ in~$E$.
 Let $\Omega(M,\End(E))$ be the graded algebra of differential forms
with coefficients in the vector bundle $\End(E)$ of endomorphisms
of~$E$, endowed with the de Rham differential $d_\nabla$ depending on
the connection $\nabla_{\End(E)}$ corresponding to~$\nabla_E$ and
the element $h_\nabla\in\Omega^2(M,\End(E))$ equal to the curvature
of~$\nabla_E$.
 The Koszul duality theorem provides an equivalence between the derived
category of right $\Diff_{M,E}$\module s and the coderived category of
right CDG\module s over $\Omega(M,\End(E))$.
 The proof of this result given in~\ref{koszul-duality-theorem}
generalizes easily to nonaffine varieties (the approach with
quasi-differential structures allows to get rid of the choice of
a global connection).

 These results are even valid in prime characteristic, describing
the derived category of modules over the ring/sheaf of crystalline
differential operators (those generated by endomorphisms and vector
fields with commutation relations analogous to the zero
characteristic case).
 Furthermore, it is not difficult to see that the quotient category
of the homotopy category of finitely generated right CDG\module s
over $\Omega(M,\End(E))$ by its minimal thick subcategory containing
the total CDG\module s of exact triples of finitely generated
CDG\module s is a full subcategory of the coderived category of
CDG\module s.
 This full subcategory is equivalent to the bounded derived category
of finitely generated (coherent) right $\Diff_{M,E}$\module s.
 All of this is applicable to any smooth varieties, not necessarily
affine.

 For a smooth affine variety~$M$, the derived category of left
$\Diff_{M,E}$\module s is equivalent to the contraderived
category of left CDG\module s over $\Omega(M,\End(E))$.

\subsubsection{}
 Koszul algebras were introduced by S.~Priddy; the standard
contemporary sources are~\cite{BGS,PP}.
 Nonhomogeneous quadratic duality (the equivalence of categories of
nonhomogeneous Koszul algebras and Koszul CDG\+algebras) was developed
in~\cite{Pos,PP}.
 Homogeneous Koszul duality (the equivalence of derived categories of
graded modules over dual Koszul algebras) was established in~\cite{BGS}.
 Koszul duality in the context of $A_\infty$\+algebras and
DG\+coalgebras was worked out in~\cite{Lef}.
 All of these papers only consider duality over the ground field
(or, in the case of~\cite{BGS}, a semisimple algebra) rather than
over an arbitrary ring, as above.

 Notable attempts to define a version of derived category of
DG\module s over the de Rham complex so that the derived category of
modules over the differential operators would be equivalent to it
were undertaken in~\cite{Kap} and~\cite[subsection~7.2]{BD2}.
 They were not entirely successful, in the present author's view, in
that in~\cite{Kap} the analytic topology and analytic functions were
used in the definition of an essentially purely algebraic category,
while in~\cite{BD2} the right hand side of the purpoted equivalence of
categories is to a certain extent defined in terms of the left hand
side.
 The latter problem is also present in Lef\`evre-Hasegawa's Koszul
duality~\cite{Lef,Kel2}.

\Section{Semialgebras and Semitensor Product}

 Throught Sections 1--11,\, $k$~is a commutative ring.
 All our rings, bimodules, abelian groups~\dots\ will be $k$\module s;
all additive categories will be $k$\+linear.

\subsection{Corings and comodules}
 Let $A$ be an associative $k$\+algebra (with unit).

\subsubsection{}
 A \emph{coring} $\C$ over~$A$ is a coring object in the tensor category
of $A$\+$A$\bimodule s; in other words, it is a $k$\module{} endowed
with an $A$\+$A$\bimodule{} structure and two $A$\+$A$\bimodule{} maps
of \emph{comultiplication} $\C\rarrow\C\ot_A\C$ and \emph{counit}
$\C\rarrow A$ satisfying the coassociativity and counity equations:
two compositions of the comultiplication map $\C\rarrow\C\ot_A\C$
with the maps $\C\ot_A\C\birarrow\C\ot_A\C\ot_A\C$ induced by
the comultiplication map  should coincide with each other and
two compositions $\C\rarrow\C\ot_A\C\birarrow\C$ of the comultiplication
map with the maps $\C\ot_A\C\birarrow\C$ induced by the counit map
should coincide with the identity map of~$\C$.
 
 A \emph{left comodule} $\M$ over a coring~$\C$ is a comodule
object in the left module category of left $A$\module s over
the coring object~$\C$ in the tensor category of $A$\+$A$\bimodule s;
in other words, it is a left $A$\module{} endowed with a left
$A$\module{} map of \emph{left coaction} $\M\rarrow\C\ot_A\M$
satisfying the coassociativity and counity equations: two compositions
of the coaction map $\M\rarrow\C\ot_A\M$ with the maps
$\C\ot_A\M\birarrow\C\ot_A\C\ot_A\M$ induced by the comultiplication
and coaction maps should coincide with each other and the composition
$\M\rarrow\C\ot_A\M\rarrow\M$ of the coaction map with the map 
$\C\ot_A\M\rarrow\M$ induced by the counit map should coincide with
the identity map of~$\M$.
 A \emph{right comodule} $\N$ over~$\C$ is a comodule object in
the right module category of right $A$\module s over the coring
object~$\C$ in the tensor category of $A$\+$A$\bimodule s; in other
words, it is a right $A$\module{} endowed with a right $A$\module{}
map of \emph{right coaction} $\N\rarrow\N\ot_A\C$ satisfying
the coassociativity and counity equations for the compositions
$\N\rarrow\N\ot_A\C\birarrow\N\ot_A\C\ot_A\C$ and
$\N\rarrow\N\ot_A\C\rarrow\N$.

\subsubsection{}   \label{induced-coinduced}
 If $V$ is a left $A$\module, then the left $\C$\comodule{}
$\C\ot_A V$ is called the left $\C$\comodule{} \emph{coinduced}
from an $A$\+module~$V$.
 The $k$\module{} of comodule homomorphisms from an arbitrary
$\C$\comodule{} into the coinduced $\C$\comodule{} is described
by the formula $\Hom_\C(\M\;\C\ot_A V)\simeq\Hom_A(\M,V)$.
 This is an instance of the following general fact, which we prefer to
formulate in the tensor (monoidal) category language, though it can be
also formulated in the monad language.

\begin{lem}
 Let\/ $\sE$ be a (not necessarily additive) associative tensor
category with a unit object, $\sM$ be a left module category over it,
$R$ be a ring object with unit in\/~$\sE$, and ${}_R\.\sM$ be
the category of $R$\module{} objects in\/~$\sM$.
 Then the induction functor\/ $\sM\rarrow {}_R\.\sM$ defined by
the rule $V\mpsto R\ot V$ is left adjoint to the forgetful
functor ${}_R\.\sM\rarrow\sM$.
\end{lem}

\begin{proof}
 For any object $V$ and any $R$\module{} $M$ in $\sM$, the map
$\Hom_\sM(V,M)\rarrow\Hom_\sM(R\ot V\;M)$ is a split equalizer
(see~\cite{McL}) of the pair of maps $\Hom_\sM(R\ot V\;M)\birarrow
\Hom_\sM(R\ot R\ot V\;M)$ in the category of sets, with
the splitting maps $\Hom_\sM(V,M)\larrow\Hom_\sM(R\ot V\;M)\larrow
\Hom_\sM(R\ot R\ot V\; M)$ induced by the unit morphism of~$R$
(applied at the rightmost factor~$R$).
\end{proof}

 We will denote the category of left $\C$\comodule s by
$\C\comodl$ and the category of right $\C$\comodule s by $\comodr\C$.
 The category of left $\C$\comodule s is abelian whenever $\C$ is
a flat right $A$\module.
 Moreover, the right $A$\module{} $\C$ is flat if and only if
the category $\C\comodl$ is abelian and the forgetful functor
$\C\comodl\rarrow A\modl$ is exact.
 This is an instance of a general fact applicable to any monad
over an abelian category.
 The ``only if'' assertion is straightforwardly checked, 
while the ``if'' part is deduced from the observations
that the coinduction functor $V\mpsto\C\ot_A V$ is right adjoint
to the forgetful functor and a right adjoint functor is left exact.

 At the same time, for any coring~$\C$ there are four natural exact
categories of left comodules: the exact category of $A$\+projective
$\C$\comodule s, the exact category of $A$\+flat $\C$\+comodules,
the exact category of arbitrary $\C$\comodule s with $A$\+split exact
triples, and the exact category of arbitrary left $\C$\comodule s
with \emph{$A$\+pure} exact triples, i.~e., the exact triples which
as triples of left $A$\module s remain exact after the tensor product
with any right $A$\module.
 Besides, any morphism of $\C$\comodule s has a cokernel and
the forgetful functor $\C\comodl\rarrow A\modl$ preserves cokernels.
 When a morphism of $\C$\comodule s has the property that its kernel
in the category of $A$\module s is preserved by the functors of
tensor product with~$\C$ and $\C\ot_A\C$ over~$A$, this kernel has
a natural $\C$\comodule{} structure, which makes it the kernel of
that morphism in the category of $\C$\comodule s.

 Infinite direct sums always exist in the category of $\C$\comodule s
and the forgetful functor $\C\comodl\rarrow A\modl$ preserves them.
 The coinduction functor $A\modl\rarrow\C\comodl$ preserves both
infinite direct sums and infinite products.
 To construct products of $\C$\comodule s, one can present them as
kernels of morphisms of coinduced comodules, so the category
of $\C$\comodule s has infinite products if it has kernels.

 If $\C$ is a projective right $A$\module, or $\C$ is a flat right
$A$\module{} and $A$ is a left Noetherian ring, then any left
$\C$\comodule{} is a union of its subcomodules that are finitely
generated as $A$\module s~\cite{BW}.

\subsubsection{}  \label{flat-comodule-surjection}
 Assume that the coring $\C$ is a flat left and right $A$\module{} and
the ring~$A$ has a finite weak homological dimension (Tor-dimension).

\begin{lem}
 There exists a (not always additive) functor assigning to
any\/ $\C$\comodule{} a surjective map onto it from
an $A$\+flat\/ $\C$\comodule.
\end{lem}

\begin{proof}
 Let $G(M)\rarrow M$ be a surjective map onto an $A$\module{} $M$
from a flat $A$\module{} $G(M)$ functorially depending on~$M$.
 For example, one can take $G(M)$ to be the direct sum of copies of
the $A$\module{} $A$ over all elements of~$M$.
 Let $\M$ be a left $\C$\comodule.
 Consider the coaction map $\M\rarrow\C\ot_A\M$; it is
an injective morphism of left $\C$\comodule s; let $\K(\M)$ denote
its cokernel.
 Let $\cQ(\M)$ be the kernel of the composition $\C\ot_A G(\M)\rarrow
\C\ot_A\M \rarrow \K(\M)$.
 Then the composition of maps $\cQ(\M)\rarrow \C\ot_A G(\M)\rarrow
\C\ot_A \M$ factorizes through the injection $\M\rarrow\C\ot_A\M$,
so there is a natural sujective morphism of $\C$\comodule s
$\cQ(\M)\rarrow\M$.
 Let us show that the flat dimension $\df_A\cQ(\M)$ of the $A$\module{}
$\cQ(\M)$ is smaller than that of~$\M$.
 Indeed, the $A$\module{} $\C\ot_A G(\M)$ is flat, hence $\df_A \cQ(\M)
= \df_A \K(\M) - 1 \le \df_A (\C\ot_A\M) - 1 \le \df_A\M - 1$, because
the $A$\module{} $\K(\M)$ is a direct summand ot the $A$\module{}
$\C\ot_A\M$ and a flat resolution of the $A$\module{} $\C\ot_A\M$
can be constructed by taking the tensor product of a flat resolution
of the $A$\module~$\M$ with the $A$\+$A$\bimodule~$\C$.
 It remains to iterate the functor $\M\mpsto\cQ(M)$ sufficiently
many times.
 Notice that the comodule $\cQ(\M)$ is an extension of~$\M$ by
a coinduced comodule $\C\ot_A \ker(G(\M)\to\M)$.
\end{proof}

\subsection{Cotensor product}
 
\subsubsection{}   \label{induced-tensor-cotensor}
 The \emph{cotensor product} $\N\oc_\C\M$ of a right $\C$\comodule{}
$\N$ and a left $\C$\comodule{} $\M$ is a $k$\module{} defined as
the kernel of the pair of maps $\N\ot_A\M\birarrow\N\ot_A\C\ot_A\M$
one of which is induced by the $\C$\+coaction in~$\N$ and the other
by the $\C$\+coaction in~$\M$.
 The functor of cotensor product is neither left, nor right exact in
general; it is left exact if the ring~$A$ is absolutely flat.
 For any right $A$\module{} $V$ and any left $\C$\comodule{} $\M$
there is a natural isomorphism $(V\ot_A\C)\oc_\C\M\simeq V\ot_A\M$.
 This is an instance of the following general fact.

\begin{lem}
 Let\/ $\sE$ be a tensor category, $\sM$ be a left module category over
it, $\sN$ be a right module category, $\sK$ be an additive category,
and\/ ${\ot}\:\sN\times\sM\rarrow\sK$ be a pairing functor compatible
with the module category structures on\/ $\sM$ and\/~$\sN$.
 Let $R$ be a ring object with unit in\/~$\sE$, \ $M$ be
an $R$\module{} object in\/~$\sM$, and\/ $V$ be an object of\/~$\sN$.
 Then the morphism $V\ot R\ot M\rarrow V\ot M$ induced by the action
of~$R$ in~$M$ is a cokernel of the pair of morphisms
$V\ot R\ot R\ot M \birarrow V\ot R\ot M$, one of which is induced by
the multiplication in~$R$ and the other by the $R$\+action in~$M$.
\end{lem}

\begin{proof}
 The whole bar complex $\dsb\rarrow V\ot R\ot R\ot M\rarrow
V\ot R\ot M\rarrow V\ot M\rarrow 0$ is contractible 
with contracting homotopy $\dsb\larrow V\ot R\ot R\ot M\larrow
V\ot R\ot M\larrow V\ot M$
induced by the unit morphism of~$R$
(applied at the leftmost factor~$R$).
\end{proof}

\subsubsection{}   \label{absolute-relative-coflat}
 Assume that $\C$ is a flat right $A$\module.
 A right comodule~$\N$ over~$\C$ is called \emph{coflat} if the functor
of cotensor product with~$\N$ is exact on the category of left
$\C$\comodule s.
 It is easy to see that any coflat $\C$\comodule{} is a flat $A$\module.
 The $\C$\comodule{} coinduced from a flat $A$\module{} is coflat.
 A left comodule~$\M$ over~$\C$ is called \emph{coflat relative to~$A$}
($\C/A$\+coflat) if its cotensor product with any exact triple of
$A$\+flat right $\C$\comodule s is an exact triple.
 Any coinduced $\C$\comodule{} is $\C/A$\+coflat.

 The definition of a relatively coflat $\C$\comodule{} does not really
depend on the flatness assumption on~$\C$, but appears to be useful
when this assumption holds.

\begin{lem}
 The classes of coflat right\/ $\C$\comodule s and\/ $\C/A$\+coflat
left\/ $\C$\comodule s are closed under extensions.
 The quotient comodule of a\/ $\C/A$\+coflat left\/ $\C$\comodule{}
by a\/ $\C/A$\+coflat subcomodule is\/ $\C/A$\+coflat; an $A$\+flat
quotient comodule of a coflat right\/ $\C$\comodule{} by a coflat
subcomodule is coflat.
 The cotensor product of an exact triple of coflat right\/
$\C$\comodule s with any left\/ $\C$\comodule{} is an exact triple
and the cotensor product of an $A$\+flat right\/ $\C$\comodule{}
with an exact triple of\/ $\C/A$\+coflat left\/ $\C$\comodule s is
an exact triple.
\end{lem}

\begin{proof}
 All of these results follow from the standard properties of the right
derived functor of the left exact functor of cotensor product on
the Carthesian product of the exact category of $A$\+flat right
$\C$\comodule s and the abelian category of left $\C$\comodule s.
 One can simply define the $k$\module s $\Cotor^\C_i(\N,\M)$, \
$i=0$,~$-1$,~\dots\ as the cohomology of the cobar complex
$\N\ot_A\M\rarrow\N\ot_A\C\ot_A\M\rarrow \N\ot_A\C\ot_A\C\ot_A\M
\rarrow\dsb$ for any $A$\+flat right $\C$\comodule{} $\N$ and any
left $\C$\comodule{} $\M$.
 Then $\Cotor^\C_0(\N,\M)\simeq\N\oc_\C\M$, and there are long exact
sequences of $\Cotor^\C_*$ associated with exact triples of
$\C$\comodule s in either argument, since in both cases the cobar
complexes form an exact triple.
 Now an $A$\+flat right $\C$\comodule{} $\N$ is coflat if and only if
$\Cotor^\C_i(\N,\M)=0$ for any left $\C$\comodule{} $\M$ and all $i<0$.
 Indeed, the ``if'' assertion follows from the homological exact
sequence, and ``only if'' holds since the cobar complex is
the cotensor product of the comodule $\N$ with the cobar resolution of
the comodule $\M$, which is exact except in degree~$0$.
 Analogously, a left $\C$\comodule{} $\M$ is $\C/A$\+coflat if and only
if $\Cotor^\C_i(\N,\M)=0$ for any $A$\+flat right $\C$\comodule{} $\N$
and all $i<0$, since the cobar resolution of the comodule~$\N$ is
a complex of $A$\+flat right $\C$\comodule s, exact except in
degree~$0$ and split over~$A$.
 The rest is obvious.
\end{proof}

\begin{rmk}
 A much more general construction of the double-sided derived
functor $\Cotor^\C_*(\N,\M)$ defined for arbitrary $\C$\comodule s
$\M$ and~$\N$ will be given, in the assumptions
of~\ref{flat-comodule-surjection}, in Section~\ref{semitor-section}.
 Using this construction, one can prove somewhat stronger results.
 In particular, $\Cotor^\C_i(\M,\N)=0$ for any $\C/A$\+coflat
left $\C$\comodule{} $\M$, any right $\C$\comodule{} $\N$, and
all $i<0$, since the $k$\module s $\Cotor^\C_i(\M,\N)$
can be computed using a left resolution of~$\N$ consisting
of $A$\+flat right $\C$\comodule s (see~\ref{relatively-semiflat}).
 Therefore, any $A$\+flat $\C/A$\+coflat $\C$\+comodule{} is coflat.
 It follows that the construction of
Lemma~\ref{flat-comodule-surjection} assigns to any $\C/A$\+coflat
$\C$\comodule{} a surjective map onto it from a coflat $\C$\comodule{}
with a $\C/A$\+coflat kernel.
\end{rmk}

\subsubsection{}     \label{tensor-cotensor-assoc}
 Now let $\C$ be an arbitrary coring.
 Let us call a left $\C$\comodule{} $\M$ \emph{quasicoflat}
if the functor of cotensor product with~$\M$ is right exact
on the category of right $\C$\comodule s, i.~e., this
functor preserves cokernels.
 Any coinduced $\C$\comodule{} is quasicoflat.
 Any quasicoflat $\C$\comodule{} is $\C/A$\+coflat.

\begin{prop}
 Let $\N$ be a right\/ $\C$\comodule, $\K$ be a left\/ $\C$\comodule{}
endowed with a right action of a $k$\+algebra $B$ by comodule
endomorphisms, and $M$ be a left $B$\module.
 Then there is a natural $k$\module{} map\/ $(\N\oc_\C\K)\ot_B M
\rarrow\N\oc_\C (\K\ot_B M)$, which is an isomorphism, at least, in
the following cases:
\begin{enumerate}
 \item $M$ is a flat left $B$\module;
 \item $\N$ is a quasicoflat right\/ $\C$\comodule;
 \item $\C$ is a flat right $A$\module, $\N$ is a flat right $A$\module,
       $\K$ is a\/ $\C/A$\+coflat left\/ $\C$\comodule,
       $\K$ is a flat right $B$\module, and the ring $B$ has
       a finite weak homological dimension;
 \item $\K$ as a left\/ $\C$\comodule{} with a right $B$\module{}
       structure is coinduced from an $A$\+$B$\bimodule.
\end{enumerate}
 Besides, in the case\/~\textup{(c)} the cotensor product\/
$\N\oc_\C\K$ is a flat right $B$\module.
\end{prop}

\begin{proof}
 The map $(\N\oc_\C\K)\ot_B M\rarrow \N\ot_A\K\ot_B M$ obtained by
taking the tensor product of the map $\N\oc_\C\K\rarrow\N\ot_A\K$
with the $B$\module{} $M$ has equal compositions with two maps
$\N\ot_A\K\ot_B M\birarrow\N\ot_A\C\ot_A\K\ot_B M$, hence there
is a natural map $(\N\oc_\C\K)\ot_B M\rarrow \N\oc_\C (\K\ot_B M)$.
 The case~(a) is obvious.
 In the case~(b), it suffices to present $M$ as the cokernel of
a map of flat $B$\module s.
 To prove (c) and~(d), consider the cobar complex
\begin{equation}  \label{tensor-cotensor}
 \N\oc_\C\K\lrarrow \N\ot_A\K\lrarrow \N\ot_A\C\ot_A\K
 \lrarrow \N\ot_A\C\ot_A\C\ot_A\K\lrarrow\dsb
\end{equation}
 In the case~(c) this complex is exact, since it is the cotensor
product of a $\C/A$\+coflat $\C$\comodule{} $\K$ with an
$A$\+split exact complex of $A$\+flat $\C$\comodule s
$\N\rarrow\N\ot_A\C\rarrow\N\ot_A\C\ot_A\C\rarrow\dsb$\,
  Since all the terms of the complex~\eqref{tensor-cotensor}, except
possibly the leftmost one, are flat right $B$\module s and the weak
homological dimension of the ring~$B$ is finite, the leftmost term
$\K\oc_\C\M$ is also a flat $B$\module{} and the tensor product
of this complex with the left $B$\module~$M$ is exact.
 In the case~(d), the complex~\eqref{tensor-cotensor} is exact and
split as a complex of right $B$\module s.
\end{proof}

\subsubsection{}    \label{bicomodule-cotensor}
 Let $\C$ be a coring over a $k$\+algebra $A$ and $\D$ be a coring
over a $k$\+algebra $B$.
 A \emph{$\C$\+$\D$\+bicomodule} $\K$ is an $A$\+$B$\bimodule{} in
the category of $k$\module s endowed with a left $\C$\comodule{}
and a right $\D$\comodule{} structures such that the right
$\D$\+coaction map $\K\rarrow\K\ot_B\D$ is a morphism of left
$\C$\comodule s and the left $\C$\+coaction map $\K\rarrow\C\ot_A\K$
is a morphism of right $B$\module s, or equivalently, the right
$\D$\+coaction map is a morphism of left $A$\module s and the left
$\C$\+coaction map is a morphism of right $\D$\comodule s.
 Equivalently, a $\C$\+$\D$\bicomodule{} is a $k$\module{} endowed
with an $A$\+$B$\bimodule{} structure and an $A$\+$B$\bimodule{} map
of \emph{bicoaction} $\K\rarrow\C\ot_A\K\ot_B\D$ satisfying
the coassociativity and counity equations.
 We will denote the category of $\C$\+$\D$\bicomodule s by
$\C\bcomod\D$.

 Assume that $\C$ is a flat right $A$\module{} and $\D$ is a flat
left $B$\module.
 Then the category of $\C$\+$\D$\bicomodule s is abelian and
the forgetful functor $\C\bcomod\D\rarrow k\modl$ is exact.
 Let $\E$ be a coring over a $k$\+algebra~$F$.
 Let $\N$ be a $\C$\+$\E$\bicomodule{} and $\M$ be
a $\E$\+$\D$\bicomodule.
 Then the cotensor product $\N\oc_\E\M$ can be endowed with
a $\C$\+$\D$\bicomodule{} structure as the kernel of a pair
of bicomodule morphisms $\N\ot_F\M\birarrow\N\ot_F\E\ot_F\M$.

 More generally, let $\C$, $\D$, and $\E$ be arbitrary corings.
 Assume that the functor of tensor product with $\C$ over $A$
and with $\D$ over $B$ preserves the kernel of the pair of maps
$\N\ot_F\M\birarrow\N\ot_F\E\ot_F\M$, that is the natural map
$\C\ot_A(\N\oc_\E\M)\ot_B\D\rarrow(\C\ot_A\N)\oc_\E(\M\ot_B\D)$
is an isomorphism.
 Then one can define a bicoaction map $\N\oc_\E\M\rarrow
\C\ot_A(\N\oc_\E\M)\ot_B\D$ taking the cotensor product over~$\E$
of the left $\C$\+coaction map $\N\rarrow\C\ot_A\N$ and
the right $\D$\+coaction map $\M\rarrow\M\ot_B\D$.
 One can easily see that this bicoaction is counital and coassociative,
at least, if the natural maps $\C\ot_A\C\ot_A(\N\oc_\E\M)\rarrow
(\C\ot_A\C\ot_A\N)\oc_\E\M$ and $(\N\oc_\E\M)\ot_B\D\ot_B\D\rarrow
\N\oc_\E(\M\ot_B\D\ot_B\D)$ are also isomorphisms.

 In particular, if $\C$ is a flat right $A$\module{} and either 
$\D$ is a flat left $B$\module, or $\N$ is a quasicoflat right
$\E$\comodule, or $\N$ is a flat right $F$\module, $\E$ is a flat right
$F$\module, $\M$ is an $\E/F$\+coflat left $\E$\comodule, $\M$ is
a flat right $B$\module, and $B$ has a finite weak homological
dimension, or $\M$ as a left $\E$\comodule{} with a right $B$\module{}
structure is coinduced from an $F$\+$B$\bimodule, then the cotensor
product $\N\oc_\E\M$ has a natural $\C$\+$\D$\bicomodule{} structure.

\subsubsection{}     \label{cotensor-associative}
 Let $\C$ be a coring over a $k$\+algebra $A$ and $\D$ be a coring
over a $k$\+algebra $B$.

\begin{prop}
 Let\/ $\N$ be a right\/ $\C$\comodule, $\K$ be
a\/ $\C$\+$\D$\+bicomodule, and\/ $\M$ be a left\/ $\D$\comodule.
 Then the iterated cotensor products\/ $(\N\oc_\C\K)\oc_\D\M$ and\/
$\N\oc_\C(\K\oc_\D\M)$ are naturally isomorphic, at least,
in the following cases:
\begin{enumerate}
 \item $\C$ is a flat right $A$\module, $\N$ is a flat right
       $A$\module, $\D$ is a flat left $B$\module, and\/ $\M$ is
       a flat left $B$\module;
 \item $\C$ is a flat right $A$\module{} and\/ $\N$ is a coflat
       right\/ $\C$\comodule;
 \item $\C$ is a flat right $A$\module, $\N$ is a flat right
       $A$\module, $\K$ is a\/ $\C/A$\+coflat left\/ $\C$\comodule,
       $\K$ is a flat right $B$\module,
       and the ring $B$ has a finite weak homological dimension;
 \item $\C$ is a flat right $A$\module, $\N$ is a flat right
       $A$\module, and\/ $\K$ as a left\/ $\C$\+co\-mod\-ule with
       a right $B$\module{} structure is coinduced from
       an $A$\+$B$\bimodule;
 \item $\M$ is a quasicoflat left\/ $\C$\comodule{} and\/ $\K$ as
       a left\/ $\C$\comodule{} with a right $B$\module{} structure
       is coinduced from an $A$\+$B$\bimodule;
 \item $\K$ as a left\/ $\C$\comodule{} with a right $B$\module{}
       structure is coinduced from an $A$\+$B$\bimodule{} and\/ $\K$
       as a right\/ $\D$\comodule{} with a left $A$\module{} structure
       is coinduced from an $A$\+$B$\bimodule.
\end{enumerate}
 More precisely, in all cases in this list the natural maps from both
iterated cotenzor products to the $k$\module\/ $\N\ot_A\K\ot_B\M$ are
injective, their images coincide and are equal to the intersection
of two submodules\/ $(\N\ot_A\K)\oc_\D\M$ and\/ $\N\oc_\C(\K\ot_B \M)$
in the $k$\module\/ $\N\ot_A\K\ot_B\M$.
\end{prop}

\begin{proof}
 One can easily see that whenever both maps
$(\N\oc_\C\K)\ot_B\M\rarrow \N\oc_\C(\K\ot_B\M)$
and $(\N\oc_\C\K)\ot_B\D\ot_B\M\rarrow \N\oc_\C(\K\ot_B\D\ot_B\M)$
are isomorphisms, the natural map $(\N\oc_\C\K)\oc_\D\M\rarrow
\N\ot_A\K\ot_B\M$ is injective and its image coincides with
the desired intersection of two submodules in $\N\ot_A\K\ot_B\M$.
 Thus it remains to apply Proposition~\ref{tensor-cotensor-assoc}.
\end{proof}

 When associativity of cotensor product of four or more (bi)comodules
is an issue, it becomes important to know that the pentagonal diagrams
of associativity isomorphisms are commutative.
 Since each of the five iterated cotensor products of four factors 
of the form $\N\oc_\C\K\oc_\E\L\oc_\D\M$ is endowed with a natural
map into the tensor product $\N\ot_A\K\ot_F\L\ot_B\M$ and
the associativity isomorphisms are, presumably, compatible with
these maps, it suffices to check that at least one of these five maps
is injective in order to show that the pentagonal diagram commutes.
 In particular, if the above Proposition provides all the five
associativity isomorphisms constituting the pentagonal diagram
and either $\M$ is a flat left $B$\module, or $\N$ is a flat
right $A$\module, or both $\K$ and $\L$ as left (right) comodules
with right (left) module structures are coinduced from bimodules,
then the pentagonal diagram is commutative.

 We will say that a multiple cotensor product of several bicomodules
$\N\oc_\C\dsb\oc_\D\M$ is associative if for any way of putting
parentheses in this product all the intermediate cotensor products
can be endowed with bicomodule structures via the construction
of~\ref{bicomodule-cotensor}, all possible associativity isomorphisms
between intermediate cotensor products exist in the sense of the last
assertion of Proposition and preserve bicomodule structures, and all
the pentagonal diagrams commute.
 This definition allows to consider associativity of cotensor products
as a \emph{property} rather than an additional structure.
 In particular, associativity isomorphisms and bicomodule structures on
associative multiple cotensor products are preserved by the morphisms
between them induced by any bicomodule morphisms of the factors.

\subsection{Semialgebras and semimodules}

\subsubsection{}
 Assume that the coring $\C$ over~$A$ is a flat right $A$\module.

 It follows from Proposition~\ref{cotensor-associative}(b) that
the category of $\C$\+$\C$\bicomodule s which are coflat right
$\C$\comodule s is an associative tensor category with a unit
object~$\C$, the category of left $\C$\comodule s is a left module
category over it, and the category of coflat right $\C$\comodule s
is a right module category over this tensor category.
 Furthermore, it follows from Proposition~\ref{cotensor-associative}(c)
that whenever the ring $A$ has a finite weak homological dimension,
the $\C$\+$\C$\bicomodule s that are flat right $A$\module s
and $\C/A$\+coflat left $\C$\comodule s also form a tensor category,
left $\C$\comodule s form a left module category over it, and
$A$\+flat right $\C$\comodule s form a right module category over
this tensor category.
 Finally, it follows from Proposition~\ref{cotensor-associative}(a)
that whenever the ring $A$ is absolutely flat, the categories of left
and right $\C$\comodule s are left and right module categories over
the tensor category of $\C$\+$\C$\bicomodule s.
 In each case, the cotensor product operation provides a pairing
between these left and right module categories compatible with
their module category structures and taking values in the category
of $k$\module s.

 A \emph{semialgebra} over~$\C$ is a ring object with unit in
one of the tensor categories of $\C$\+$\C$\bicomodule s of
the kind described above.
 In other words, a semialgebra~$\S$ over~$\C$
is a $\C$\+$\C$\bicomodule{} satisfying appropriate (co)flatness
conditions guaranteeing associativity of cotensor products 
$\S\oc_\C\dsb\oc_\C\S$ of any number of copies of~$\S$ and
endowed with two bicomodule morphisms of \emph{semimultiplication}
$\S\oc_\C\S\rarrow\S$ and \emph{semiunit} $\C\rarrow\S$ satisfying
the associativity and unity equations.
 Namely, two compositions $\S\oc_\C\S\oc_\C\S\birarrow\S\oc_\C\S
\rarrow\S$ of the morphisms $\S\oc_\C\S\oc_\C\S\birarrow\S\oc_\C\S$
induced by the semimultiplication morphism with the semimultiplication
morphism $\S\oc_\C\S\rarrow\S$ should coincide with each other
and two compositions $\S\birarrow\S\oc_\C\S\rarrow\S$ of the morphisms
$\S\birarrow\S\oc_\C\S$ induced by the semiunit morphism with
the semimultiplication morphism should coincide with 
the identity morphism of~$\S$.

 A \emph{left semimodule} over~$\S$ is a module object in one
of the left module categories of $\C$\comodule s of the above
kind over the ring object~$\S$ in the corresponding tensor category
of $\C$\+$\C$\bicomodule s.
 In other words, a left $\S$\semimodule{} $\bM$ is a left
$\C$\comodule{} endowed with a left $\C$\comodule{} morphism of
\emph{left semiaction} $\S\oc_\C\bM\rarrow\bM$ satisfying
the associativity and unity equations.
 Namely, two compositions $\S\oc_\C\S\oc_\C\bM\birarrow\S\oc_\C\bM
\rarrow\bM$ of the morphisms $\S\oc_\C\S\oc_\C\bM\birarrow\S\oc_\C\bM$
induced by the semimultiplication and the semiaction morphisms with
the semiaction morphism $\S\oc_\C\bM\rarrow\bM$ should coincide with
each other and the composition $\bM\rarrow\S\oc_\C\bM\rarrow\bM$
of the morphism $\bM\rarrow\S\oc_\C\bM$ induced by the semiunit
morphism with the semiaction morphism should coincide with
the identity morphism of~$\bM$.
 For this definition to make sense, (co)flatness conditions imposed
on $\S$ and/or $\bM$ must guarantee associativity of multiple
cotensor products of the form $\S\oc_\C\dsb\oc_\C\S\oc_\C\bM$.
 \emph{Right semimodules} over~$\S$ are defined in the analogous way.

 If $\L$ is a left $\C$\comodule{} for which the multiple cotensor
products $\S\oc_\C\dsb\oc_\C\S\oc_\C\L$ are associative, then there
is a natural left $\S$\semimodule{} structure on the cotensor product
$\S\oc_\C\L$.
 The left semimodule $\S\oc_\C\L$ is called the left $\S$\semimodule{}
\emph{induced} from a $\C$\comodule~$\L$.
 According to Lemma~\ref{induced-coinduced}, the $k$\module{} of
semimodule homomorphisms from the induced $\S$\semimodule{} to
an arbitrary $\S$\semimodule{} is described by the formula
$\Hom_\S(\S\oc_\C\L\;\bM)\simeq\Hom_\C(\L,\bM)$.

 We will denote the category of left $\S$\semimodule s by $\S\simodl$
and the category of right $\S$\semimodule s by $\simodr\S$. 
 This notation presumes that one can speak of (left or right)
$\S$\semimodule s with no flatness conditions imposed on them.
 If $\S$ is a coflat right $\C$\comodule, the category of left
semimodules over~$\S$ is abelian and the forgetful functor
$\S\simodl\rarrow\C\comodl$ is exact.

 Assume that either $\S$ is a coflat right $\C$\comodule, or $\S$ is
a flat right $A$\module{} and a $\C/A$\+coflat left $\C$\comodule{}
and $A$ has a finite weak homological dimension, or $A$ is absolutely
flat.
 Then both infinite direct sums and infinite products exist in
the category of left $\S$\semimodule s, and both are preserved by
the forgetful functor $\S\simodl\rarrow\C\comodl$, even though only
infinite direct sums are preserved by the full forgetful functor
$\S\simodl\rarrow A\modl$.

 If $\S$ is a flat right $A$\module{} and a $\C/A$\+coflat left
$\C$\comodule{} and $A$ has a finite weak homological dimension,
then the category of $A$\+flat right $\S$\semimodule s is exact.
 Of course, if $\S$ is a coflat right $\C$\comodule, then
the category of $A$\+flat left $\S$\semimodule s is exact.
 In both cases there are exact categories of $\C$\+coflat right
$\S$\semimodule s and $\C/A$\+coflat left $\S$\semimodule s.
 If $A$ is absolutely flat, there are exact categories of $\C$\+coflat
left and right $\S$\semimodule s.
 Infinite direct sums exist in all of these exact categories,
and the forgetful functors preserve them.

\subsubsection{}
 Assume that the coring $\C$ is a flat left and right
$A$\module, the semialgebra $\S$ is a flat left $A$\module{}
and a coflat right $\C$\comodule, and the ring~$A$ has a finite weak
homological dimension.

\begin{lem}   \label{flat-semimodule-surjection}
 There exists a (not always additive) functor assigning to any
left\/ $\S$\semimodule{} a surjective map onto it from
an $A$\+flat left\/ $\S$\semimodule.
\end{lem}

\begin{proof}
 Let $\cP(\M)\rarrow\M$ denote the functorial surjective morphism
onto a $\C$\comodule{} $\M$ from an $A$\+flat $\C$\comodule{}
$\cP(\M)$ constructed in Lemma~\ref{flat-comodule-surjection}.
 Then for any left $\S$\semimodule{} $\bM$ the composition of maps
$\S\oc_\C\cP(\bM)\rarrow \S\oc_\C\bM\rarrow\bM$ provides the desired
surjective morphism of $\S$\semimodule s.
 According to the last assertion of
Proposition~\ref{tensor-cotensor-assoc}
(with the left and right sides switched),
the $A$\module{} $\bcP(\bM)=\S\oc_\C\cP(\bM)$ is flat.
\end{proof}

\begin{rmk}
 In the above assumptions, the same construction provides also
a (not always additive) functor assigning to any $\C/A$\+coflat right
$\S$\semimodule{} a surjective map onto it from a semiflat right
$\S$\semimodule{} (see~\ref{semiflat-subsect}) with
a $\C/A$\+coflat kernel.
 This follows from Lemma~\ref{absolute-relative-coflat} and
Remark~\ref{absolute-relative-coflat}, since the cotensor product with
$\S$ over~$\C$ preserves the kernel of the morphism $\cP(\bN)\rarrow
\bN$ and the kernel of the map $\bN\oc_\C\S\rarrow\bN$ is isomorphic
to a direct summand of $\bN\oc_\C\S$ as a right $\C$\comodule.
\end{rmk}

\subsubsection{}   \label{coflat-semimodule-injection}
 Assume that the coring $\C$ is a flat right $A$\module,
the semialgebra $\S$ is a $\C/A$\+coflat left $\C$\comodule{}
and a coflat right $\C$\comodule, and the ring $A$ has
a finite weak homological dimension.
 
\begin{lem}
 There exists an exact functor assigning to any $A$\+flat
right\/ $\S$\semimodule{} an injective morphism from it into a coflat
right\/ $\S$\semimodule{} with an $A$\+flat quotient semimodule.
 Besides, there exists an exact functor assigning to any left\/
$\S$\semimodule{} an injective morphism from it into a\/ $\C/A$\+coflat
left\/ $\S$\semimodule.
\end{lem}

\begin{proof}
 For any $A$\+flat right $\C$\comodule{} $\N$, set $\cG(\N)=\N\ot_A\C$.
 Then the coaction map $\N\rarrow\cG(\N)$ is an injective morphism of
$\C$\comodule s, the comodule $\cG(\N)$ is coflat, and the quotient
comodule $\cG(\N)/\N$ is $A$\+flat.
 Now let $\bN$ be an $A$\+flat right $\S$\semimodule.
 The semiaction map $\bN\oc_\C\S\rarrow\bN$ is a surjective morphism
of $A$\+flat $\S$\semimodule s; let $\bK(\bN)$ denote its kernel.
 The map $\bN\oc_\C\S\rarrow\cG(\bN)\oc_\C\S$ is an injective morphism
of $A$\+flat $\S$\semimodule s with an $A$\+flat quotient semimodule
$(\cG(\bN)/\bN)\oc_\C\S$.
 Let $\bcQ(\bN)$ be the cokernel of the composition
$\bK(\bN)\rarrow\bN\oc_\C\S\rarrow\cG(\bN)\oc_\C\S$.
 Then the composition of maps $\bN\oc_\C\S\rarrow\cG(\bN)\oc_\C\S
\rarrow\bcQ(\bN)$ factorizes through the surjection
$\bN\oc_\C\S\rarrow\bN$, so there is a natural injective morphism
of $\S$\semimodule s $\bN\rarrow\bcQ(\bN)$.
 The quotient semimodule $\bcQ(\bN)/\bN$ is isomorphic
to $(\cG(\bN)/\bN)\oc_\C\S$, hence both $\bcQ(\bN)/\bN$ and $\bcQ(\bN)$
are flat $A$\module s.

 Notice that the semimodule morphism $\bN\rarrow\bcQ(\bN)$
can be lifted to a comodule morphism $\bN\rarrow\cG(\bN)\oc_\C\S$.
 Indeed, the map $\bN\rarrow\bcQ(\bN)$ can be presented as
the composition $\bN\rarrow\bN\oc_\C\S\rarrow\cG(\bN)\oc_\C\S
\rarrow\bcQ(\bN)$, where the map $\bN\rarrow\bN\oc_\C\S$ is
induced by the semiunit morphism $\C\rarrow\S$ of the semialgebra~$\S$.
 
 Iterating this construction, we obtain an inductive system of
$\C$\comodule{} morphisms $\bN\rarrow\cG(\bN)\oc_\C\S\rarrow \bcQ(\bN)
\rarrow \cG(\bcQ(\bN))\oc_\C\S\rarrow\bcQ(\bcQ(\bN))\rarrow\dsb$,
where the maps $\bN\rarrow\bcQ(\bN)\rarrow\bcQ(\bcQ(\bN))
\rarrow\dsb$ are injective morphisms of $\S$\semimodule s with
$A$\+flat cokernels, while the $\C$\comodule s $\cG(\bN)\oc_\C\S$,
\ $\cG(\bcQ(\bN))\oc_\C\S$,~\dots\ are coflat.
 Denote by $\bcJ(\bN)$ the inductive limit of this system; then
$\bN\rarrow\bcJ(\bN)$ is an injective morphism of $\S$\semimodule s
with an $A$\+flat cokernel and the $\C$\comodule{} $\bcJ(\bN)$
is coflat (since the functor of cotensor product preserves filtered
inductive limits).

 A functorial injection $\bM\rarrow\bcJ(\bM)$ of any left
$\S$\semimodule{} $\bM$ into a $\C/A$\+coflat left $\S$\semimodule{}
$\bcJ(\bM)$ is provided by the same construction (with the left and
right sides switched).
 The only changes are that $A$\module s are no longer flat, for any
left $\C$\comodule{} $\M$ the $\C$\comodule{} $\cG(\M)=\C\ot_A\M$ is
$\C/A$\+coflat, and therefore the $\S$\semimodule{} $\S\oc_\C\cG(M)$
is $\C/A$\+coflat.

 Both functors~$\bcJ$ are exact, since the kernels of surjective
maps, the cokernels of injective maps, and the filtered inductive
limits preserve exact triples.
\end{proof}

\subsection{Semitensor product}

\subsubsection{}  \label{semitensor-subsect}
 Assume that the coring $\C$ is a flat right $A$\module,
the semialgebra $\S$ is a flat right $A$\module{} and a $\C/A$\+coflat
left $\C$\comodule, and the ring $A$ has a finite weak homological
dimension.
 Let $\bM$ be a left $\S$\semimodule{} and $\bN$ be an $A$\+flat
right $\S$\semimodule.
 The \emph{semitensor product} $\bN\os_\S\bM$ is a $k$\module{}
defined as the cokernel of the pair of maps $\bN\oc_\C\S\oc_\C\bM
\birarrow\bN\oc_\C\bM$ one of which is induced by the $\S$\+semiaction
in $\bN$ and another by the $\S$\+semiaction in~$\bM$.
 Even under the strongest of our (co)flatness conditions on $\C$
and~$\S$, the flatness of either~$\bN$ or $\bM$ is still needed
to guarantee that the triple cotensor product $\bN\oc_\C\S\oc_\C\bM$
is associative.

 For any $A$\+flat right $\S$\semimodule{} $\bN$ and any left
$\C$\comodule{} $\L$ there is a natural isomorphism
$\bN\os_\S(\S\oc_\C\L)\simeq \bN\oc_\C\L$.
 Analogously, for any $A$\+flat right $\C$\comodule{} $\R$ and
any left $\S$\semimodule{} $\bM$ there is a natural isomorphism
$(\R\oc_\C\S)\os_\S\bM\simeq \R\oc_\C\bM$.
 These assertions follow from Lemma~\ref{induced-tensor-cotensor}.

\subsubsection{}  \label{semiflat-subsect}
 If the coring $\C$ is a flat right $A$\module{} and the semialgebra
$\S$ is a coflat right $\C$\comodule, one can define the semitensor
product of a $\C$\+coflat right $\S$\semimodule{} and an arbitrary
left $\S$\semimodule.
 In these assumptions, a $\C$\+coflat right $\S$\semimodule{} $\bN$
is called \emph{semiflat} if the functor of semitensor product
with~$\bN$ is exact on the abelian category of left $\S$\semimodule s.
 The $\S$\semimodule{} induced from a coflat $\C$\comodule{} is
semiflat.

 If $\C$ is a flat right $A$\module, $\S$ is a coflat right
$\C$\comodule{} and $\C/A$\+coflat left $\C$\comodule, and the ring $A$ 
has a finite weak homological dimension, one can define semiflat
$\S$\semimodule s as $A$\+flat right $\S$\semimodule s such that
the functors of semitensor product with them are exact.
 Then one can prove that any semiflat $\S$\semimodule{} is a coflat
$\C$\comodule.

 When the ring $A$ is absolutely flat, the semitensor product of
arbitrary two $\S$\semimodule s is defined without any conditions
on the coring~$\C$ and the semialgebra~$\S$.

\subsubsection{}
 Let $\S$ be a semialgebra over a coring $\C$ over a $k$\+algebra $A$
and $\T$ be a semialgebra over a coring~$\D$ over a $k$\+algebra~$B$.
 Let $\bK$ denote a $\C$\+$\D$\bicomodule.
 One can speak about \emph{$\S$\+$\T$\+bisemimodule} structures on
$\bK$ if the (co)flatness conditions imposed on~$\S$, $\T$, and $\bK$
guarantee associativity of multiple cotensor products of the form
$\S\oc_\C\dsb\oc_\C\S\oc_\C\bK\oc_\D \T\oc_\D\dsb\oc_\D\T$.
 Assuming that this is so, $\bK$ is called an $\S$\+$\T$\bisemimodule{}
if it is endowed with a left $\S$\semimodule{} and a right
$\T$\semimodule{} structures such that the right $\T$\+semiaction map
$\bK\oc_\D\T\rarrow\bK$ is a morphism of left $\S$\semimodule s and
the left $\S$\+semiaction map $\S\oc_\C\bK\rarrow\bK$ is a morphism
of right $\D$\comodule s, or equivalently, the right $\T$\+semiaction
map is a morphism of left $\C$\comodule s and the left $\S$\+semiaction
map is a morphism of right $\T$\semimodule s.
 Equivalently, the $\C$\+$\D$\bicomodule{} $\bK$ is called
an $\S$\+$\T$\bisemimodule{} if it is endowed with
a $\C$\+$\D$\bicomodule{} morphism of \emph{bisemiaction}
$\S\oc_\C\bK\oc_\D\T\rarrow\bK$ satisfying the associativity and
unity equations.

 In particular, one can speak about $\S$\+$\T$\bisemimodule s $\bK$
without imposing any (co)flatness conditions on $\bK$ if $\C$ is a flat
right $A$\module{} and either $\S$ is a coflat right $\C$\comodule,
or $\S$ is a flat right $A$\module{} and a $\C/A$\+coflat left
$\C$\comodule{} and $A$ has a finite weak homological dimension,
while $\D$ is a flat left $B$\module{} and either $\T$ is a coflat
left $\D$\comodule, or $\T$ is a flat left $B$\module{} and
a $\D/B$\+coflat right $\D$\comodule{} and $B$ has a finite weak
homological dimension.
 We will denote the category of $\S$\+$\T$\bisemimodule s by
$\S\bsimod\T$.
 Besides, one can consider $B$\+flat $\S$\+$\T$\bisemimodule s
if $\C$ is a flat right $A$\module{} and $\S$ is a coflat right
$\C$\comodule, while $\D$ is a flat right $B$\module, $\T$ is
a flat right $B$\module{} and a $\D/B$\+coflat right $\D$\comodule,
and $B$ has a finite weak homological dimension; and one can
consider $\D$\+coflat $\S$\+$\T$\bisemimodule s if $\C$ is a flat
right $A$\module{} and $\S$ is a coflat right $\C$\comodule, while
$\D$ is a flat right $B$\module{} and $\T$ is a coflat right
$\D$\comodule.

\subsubsection{}    \label{semitensor-associative}
 Let $\bR$ be a semialgebra over a coring~$\E$ over a $k$\+algebra~$F$.
 Let $\bN$ be an $\S$\+$\bR$\bisemimodule{} and $\bM$ be an
$\bR$\+$\T$\bisemimodule.
 We would like to define an $\S$\+$\T$\bisemimodule{} structure on
the semitensor product $\bN\os_\bR\bM$.

 Assume that multiple cotensor products of the form
$\S\oc_\C\dsb\oc_\C\S\oc_\C\bN\oc_\E\bR\oc_\E\bM\oc_\D\T\oc_\D
\dsb\oc_\D\T$ are associative.
 Then, in particular, the semitensor products
$(\S^{\suboc\.n}\oc_\C\bN) \os_\bR (\bM\oc_\D\T^{\suboc\.m})$
can be defined.
 Assume in addition that multiple cotensor products of the form
$\S\oc_\C\dsb\oc_\C\S\oc_\C\bN\oc_\E\bM\oc_\D\T\oc_\D \dsb\oc_\D\T$
are associative.
 Then the semitensor products $(\S^{\suboc\.n}\oc_\C\bN) \os_\bR
(\bM\oc_\D\T^{\suboc\.m})$ have natural $\C$\+$\D$\bicomodule{}
structures as cokernels of $\C$\+$\D$\bicomodule{} morphisms.
 Assume that multiple cotensor products of the form
$\S\oc_\C\dsb\oc_\C\S\oc_\C(\bN\os_\bR\bM)\oc_\D\T\oc_\D \dsb\oc_\D\T$
are also associative.
 Finally, assume that the semitensor product with $\S^{\suboc\.n}$
over~$\C$ and with $\T^{\suboc\.m}$ over~$\D$ preserves the cokernel
of the pair of morphisms $\bN\oc_\E\bR\oc_\E\bM\birarrow\bN\oc_\E\bM$
for $n+m=2$, that is the bicomodule morphisms
$(\S^{\suboc\.n}\oc_\C\bN)\os_\bR(\bM\oc_\D\T^{\suboc\.m})\rarrow
\S^{\suboc\.n}\oc_\C(\bN\os_\bR\bM)\oc_\D\T^{\suboc\.m}$
are isomorphisms.
 Then one can define an associative and unital bisemiaction morphism
$\S\oc_\C(\bN\os_\bR\bM)\oc_\D\T\rarrow \bN\os_\bR\bM$ taking
the semitensor product over~$\bR$ of the morphism of $\S$\+semiaction
in~$\bN$ and the morphism of $\T$\+semiaction in~$\bM$.

 For example, if $\C$ is a flat right $A$\module, $\S$ is a coflat
right $\C$\comodule, $\D$ is a flat left $B$\module, $\T$ is a coflat
right $\D$\comodule, $\E$ is a flat right $F$\module, $\bR$ is
a flat right $F$\module{} and a $\E/F$\+coflat left $\E$\comodule,
and $F$ has a finite weak homological dimension, then the semitensor
product of any $F$\+flat $\S$\+$\bR$\bisemimodule{} $\bN$ and any
$\bR$\+$\T$\bisemimodule{} $\bM$ has a natural $\S$\+$\T$\bisemimodule{}
structure.
 Since the category of $\S$\+$\T$\bisemimodule s is abelian in this
case, the bisemimodule $\bN\os_\bR\bM$ can be simply defined as
the cokernel of the pair of bisemimodule morphisms
$\bN\oc_\E\bR\oc_\E\bM\birarrow\bN\oc_\E\bM$.

\begin{prop}
 Let\/ $\bN$ be a right\/ $\S$\semimodule, $\bK$ be
an $\S$\+$\T$\bisemimodule, and\/ $\bM$ be a left\/ $\T$\semimodule.
 Then the iterated semitensor products $(\bN\os_\S\bK)\os_\T\bM$
and\/ $\bN\os_\S(\bK\os_\T\bM)$ are well-defined and naturally
isomorphic, at least, in the following cases:
\begin{enumerate}
 \item $\C$ is a flat right $A$\module, $\S$ is a coflat right\/
       $\C$\comodule, $\bN$ is a coflat right\/ $\C$\comodule,
       $\D$ is a flat left $B$\module, $\T$ is a coflat left\/
       $\D$\comodule, and\/ $\bM$ is a coflat left\/ $\D$\comodule;
 \item $\C$ is a flat right $A$\module, $\S$ is a coflat right\/
       $\C$\comodule, $\bN$ is a semiflat right\/ $\S$\semimodule,
       and either
       \begin{itemize}
          \item $\D$ is a flat right $B$\module, $\T$ is a coflat
                right\/ $\D$\comodule, and\/ $\bK$ is a coflat right\/
                $\D$\comodule, or
          \item $\D$ is a flat right $B$\module, $\T$ is a flat right
                $B$\module{} and a\/ $\D/B$\+coflat left\/
                $\D$\comodule, the ring $B$ has a finite weak
                homological dimension, and\/ $\bK$ is a flat right
                $B$\module, or
          \item $\D$ is a flat left $B$\module, $\T$ is a flat left
                $B$\module{} and a\/ $\D/B$\+coflat right\/
                $\D$\comodule, the ring $B$ has a finite weak
                homological dimension, and\/ $\bM$ is a flat left
                $B$\module, or
          \item the ring $B$ is absolutely flat;
       \end{itemize}
 \item $\C$ is a flat right $A$\module, $\S$ is a coflat right\/
       $\C$\comodule, $\bN$ is a coflat right\/ $\C$\comodule,
        and either
       \begin{itemize}
          \item $\D$ is a flat right $B$\module, $\T$ is a coflat
                right\/ $\D$\comodule, and\/ $\bK$ as a left\/
                $\S$\semimodule{} with a right\/ $\D$\comodule{}
                structure is induced from a\/ $\D$\+coflat\/
                $\C$\+$\D$\bicomodule, or
          \item $\D$ is a flat right $B$\module, $\T$ is a flat right
                $B$\module{} and a\/ $\D/B$\+coflat left\/
                $\D$\comodule, the ring $B$ has a finite weak
                homological dimension, and\/ $\bK$ as a left\/
                $\S$\semimodule{} with a right\/ $\D$\comodule{}
                structure is induced from a $B$\+flat\/
                $\C$\+$\D$\+bicomodule, or
          \item $\D$ is a flat left $B$\module, $\T$ is a flat left
                $B$\module{} and a\/ $\D/B$\+coflat right\/
                $\D$\comodule, the ring $B$ has a finite weak
                homological dimension, $\bK$ as a left\/
                $\S$\semimodule{} with a right\/ $\D$\comodule{}
                structure is induced from a\/ $\C$\+$\D$\+bicomodule,
                and\/ $\bM$ is a flat left $B$\module, or
          \item the ring $B$ is absolutely flat and\/ $\bK$ as a left\/
                $\S$\semimodule{} with a right\/ $\D$\comodule{}
                structure is induced from a $B$\+flat\/
                $\C$\+$\D$\+bicomodule.
       \end{itemize}
\end{enumerate}
More precisely, in all cases in this list the natural maps into
both iterated semitensor products from the $k$\module\/
$(\bN\oc_\C\bK)\oc_\D\bM\simeq\bN\oc_\C(\bK\oc_\D\bM)$ are surjective,
their kernels coincide and are equal to the sum of the kernels
of two maps from this module onto its quotient modules\/
$(\bN\oc_\C\bK)\os_\T\bM$ and\/ $\bN\os_\S(\bK\oc_\D\bM)$.
\end{prop}

\begin{proof}
 It follows from Proposition~\ref{cotensor-associative} that
all multiple cotensor products of the form $\bN\oc_\C\S\oc_\C\dsb
\oc_\C\S\oc_\C\bK\oc_\D\T\oc_\D\dsb\oc_\D\T\oc_\D\bM$ are
associative.
 Multiple cotensor products $\bN\oc_\C\S\oc_\C\dsb\oc_\C\S\oc_\C
(\bK\os_\T\bM)$ and $(\bN\os_\S\bK)\oc_\D\T\oc_\D\dsb\oc_\D\T\oc_\D\bM$
are also associative by the same Proposition (here one has to notice
that the semitensor product $\bN\os_\S\bK$ is a coflat right
$\D$\comodule{} whenever $\bK$ is a coflat right $\D$\comodule{} and
$\bN$ is a semiflat right $\S$\semimodule).
 The map $\bN\oc_\C\bK\oc_\D\bM\rarrow(\bN\os_\S\bK)\oc_\D\bM$
factorizes through the surjection $\bN\oc_\C\bK\oc_\D\bM\rarrow
\bN\os_\S(\bK\oc_\D\bM)$, hence there is a natural map
$\bN\os_\S(\bK\oc_\D\bM)\rarrow(\bN\os_\S\bK)\oc_\D\bM$.
 One can easily see that whenever this map and the analogous maps
for $\T$, \ $\T\oc_\D\T$, and $\T\oc_\D\bM$ in place of $\bM$ are
isomorphisms, the iterated semitensor product $(\bN\os_\S\bK)\os_\T\bM$
is defined, the natural map $\bN\oc_\C\bK\oc_\D\bM\rarrow
(\bN\os_\S\bK)\os_\T\bM$ is surjective, and its kernel is equal to
the desired sum of two kernels of maps from $\bN\oc_\C\bK\oc_\D\bM$
onto its quotient modules.
 Thus it remains to prove that the map $\bN\os_\S(\bK\oc_\D\bM)
\rarrow(\bN\os_\S\bK)\oc_\D\bM$ is an isomorphism, i.~e., the exact
sequence of right $\D$\comodule s $\bN\oc_\C\S\oc_\C\bK\rarrow
\bN\oc_\C\bK\rarrow\bN\os_\S\bK\rarrow0$ remains exact after taking
the cotensor product with~$\bM$ over~$\D$.
 This is obvious if $\bM$ is a quasicoflat $\D$\comodule.
 If $\bN$ is a semiflat $\S$\semimodule, it suffices to present $\bM$
as a kernel of a morphism of (quasi)coflat $\D$\comodule s.
 Finally, if $\bK$ as a left $\S$\semimodule{} with a right
$\D$\comodule{} structure is induced from a $\C$\+$\D$\+bicomodule,
then our exact sequence of right $\D$\comodule s splits.
\end{proof}

\Section{Derived Functor SemiTor}  \label{semitor-section}

\subsection{Coderived categories}     \label{coderived-categories}
 A complex $C^\bu$ over an exact category~$\sA$ is called exact
if it is composed of exact triples $Z^i\to C^i\to Z^{i+1}$ in~$\sA$.
 A complex over~$\sA$ is called acyclic if it is homotopy equivalent
to an exact complex (or equivalently, if it is a direct summand of
an exact complex).
 Acyclic complexes form a thick subcategory $\Acycl(\sA)$ of
the homotopy category $\Hot(\sA$) of complexes over~$\sA$.
 All acyclic complexes over~$\sA$ are exact if and only if
$\sA$ contains images of idempotent endomorphisms~\cite{Nee1}.
 The quotient category $\sD(\sA)=\Hot(\sA)/\Acycl(\sA)$ is called
the derived category of~$\sA$.

 Let $\sA$ be an exact category where all infinite direct sums
exist and the functors of infinite direct sum are exact.
 By the total complex of an exact triple ${}'\!K^\bu\to K^\bu\to
{}''\!K^\bu$ of complexes over~$\sA$ we mean the total complex of
the corresponding bicomplex with three rows.
 A complex $C^\bu$ over~$\sA$ is called \emph{coacyclic} if it
belongs to the minimal triangulated subcategory $\Acycl^\co(\sA)$
of the homotopy category $\Hot(\sA)$ containing all the total
complexes of exact triples of complexes over~$\sA$ and closed under
infinite direct sums.
 Any coacyclic complex is acyclic.
 Acyclic complexes are not always coacyclic
(see~\ref{prelim-unbounded-cotor}).
 It follows from the next Lemma that any acyclic complex bounded
from below is coacyclic.

\begin{lem}
 Let\/ $0\to M^{0,\bu}\to M^{1,\bu}\to\dsb$ be an exact sequence,
bounded from below, of arbitrary complexes over~$\sA$.
 Then the total complex $T^\bu$ of the bicomplex $M^{\bu,\bu}$
constructed by taking infinite direct sums along the diagonals
is coacyclic.
\end{lem}

\begin{proof}
 An exact sequence of complexes $0\to M^{0,\bu}\to M^{1,\bu}\to\dsb$
can be presented as the inductive limit of finite exact sequences
of complexes $0\to M^{0,\bu}\to\dsb\to M^{n,\bu}\to Z^{n+1,\bu}\to0$.
 The total complex $T_n^\bu$ of the latter finite exact sequence
is homotopy equivalent to a complex obtained from total complexes
of the exact triples $Z^{n,\bu}\to M^{n,\bu}\to Z^{n+1,\bu}$
using the operations of shift and cone.
 Hence the complexes $T_n^\bu$ are coacyclic.
 The complex $T^\bu$ is their inductive limit; moreover, the inductive
system of~$T_n^\bu$ is obtained by applying the functor of total
complex to a locally stabilizing inductive system of bicomplexes.
 Therefore, the construction of homotopy inductive limit provides
an exact triple of complexes
$\bigoplus_n T_n^\bu\rarrow \bigoplus_n T_n^\bu \rarrow T^\bu$.
 Since the total complex of this exact triple is coacyclic and
the direct sum of coacyclic complexes is coacyclic,
the complex $T^\bu$ is coacyclic.
 (In fact, this exact triple of complexes is split in every degree,
so its total complex is even contractible.)
\end{proof}

 The category of coacyclic complexes $\Acycl^\co(\sA)$ is
a thick subcategory of the homotopy category $\Hot(\sA)$,
since it is a triangulated subcategory with infinite
direct sums~\cite{Nee1,Nee2}.
 The \emph{coderived category} $\sD^\co(\sA)$ of an exact
category~$\sA$ is defined as the quotient category
$\Hot(\sA)/\Acycl^\co(\sA)$.

\begin{rmk}
 If an exact category~$\sA$ has a finite homological dimension,
then the minimal triangulated subcategory of the homotopy
category $\Hot(\sA)$ containing the total complexes of exact
triples of complexes over~$\sA$ coincides with the subcategory
of acyclic complexes.
 Indeed, let $C^\bu$ be an exact complex over~$\sA$ and $n$ be
a number greater than the homological dimension of~$\sA$.
 Let $Z^i$ be the objects of cycles of the complex $C^\bu$.
 Then for any integer~$j$ the Yoneda extension class represented by
the extension $Z^{2jn}\to C^{2jn}\to\dsb\to C^{2jn+n-1}\to Z^{2jn+n}$
is trivial, and therefore, this extension can be
connected with the split extension by a pair of extension morphisms
$(Z^{2jn}\to C^{2jn}\to\dsb\to C^{2jn+n-1}\to Z^{2jn})\rarrow
(Z^{2jn}\to {}'C^{2jn}\to\dsb\to{}'C^{2jn+n-1} \to Z^{2jn+n})
\larrow (Z^{2jn}\to Z^{2jn}\to0\to\dsb\to0\to Z^{2jn+n}\to Z^{2jn+n})$.
 Let ${}'C^\bu$ be the complex obtained by replacing all the even
segments $C^{2jn}\to\dsb\to C^{2jn+n-1}$ of the complex $C^\bu$ with
the segments ${}'C^{2jn}\to\dsb\to{}'C^{2jn+n-1}$ while leaving
the odd segments $C^{2jn+n}\to\dsb\to C^{2(j+1)n-1}$ in place, and
let ${}''C^\bu$ be the complex obtained by replacing the same even
segments of the complex $C^\bu$ with the segments $Z^{2jn}\to 0\to
\dsb\to 0\to Z^{2jn+n}$ while leaving the odd segments in place.
 Then the complex ${}''C^\bu$ and the cones of both morphisms
$C^\bu\rarrow {}'C^\bu$ and ${}''C^\bu\rarrow{}'C^\bu$ are homotopy
equivalent to complexes obtained from total complexes of exact
triples of complexes with zero differentials using the operation of
cone repeatedly.
\end{rmk}

\subsection{Coflat complexes}    \label{coflat-complexes}
 Let $\C$ be a coring over a $k$\+algebra~$A$.
 The cotensor product $\N^\bu\oc_\C\M^\bu$ of a complex of right
$\C$\comodule s $\N^\bu$ and a complex of left $\C$\comodule s $\M^\bu$
is defined as the total complex of the bicomplex $\N^i\oc_\C\M^j$,
constructed by taking infinite direct sums along the diagonals.

 Assume that $\C$ is a flat right $A$\module.
 Then the category of left $\C$\comodule s is an abelian category
with exact functors of infinite direct sums, so the coderived
category $\sD^\co(\C\comodl)$ is defined.
 When speaking about \emph{coacyclic complexes} of $\C$\comodule s,
we will always mean coacyclic complexes with respect to
the abelian category of $\C$\comodule s, unless another
exact category of $\C$\comodule s is explicitly mentioned.

 A complex of right $\C$\comodule s $\N^\bu$ is called \emph{coflat}
if the complex $\N^\bu\oc_\C\M^\bu$ is acyclic whenever a complex
of left $\C$\comodule s $\M^\bu$ is coacyclic.

\begin{lem}
 Any complex of coflat\/ $\C$\comodule s is coflat. 
\end{lem}
 
\begin{proof}
 Let $\N^\bu$ be a complex of coflat $\C$\comodule s.
 Since the functor of cotensor product with $\N^\bu$ preserves
shifts, cones, and infinite direct sums, it suffices to show
the complex $\N^\bu\oc_\C\M^\bu$ is acyclic whenever $\M^\bu$
is the total complex of an exact triple of complexes of left
$\C$\comodule s ${}'\K^\bu\to\K^\bu\to{}''\K^\bu$.
 In this case, the triple of complexes $\N^\bu\oc_\C{}'\K^\bu\rarrow
\N^\bu\oc_\C\K^\bu\rarrow\N^\bu\oc_\C{}''\K^\bu$ is also exact,
because $\N^\bu$ is a complex of coflat $\C$\comodule s, and
the complex $\N^\bu\oc_\C\M^\bu$ is the total complex of this exact
triple.
\end{proof}

 If the ring $A$ has a finite weak homological dimension, then any
coflat complex of $\C$\comodule s is a flat complex of $A$\module s
in the sense of~\ref{prelim-unbounded-tor}.
 (Indeed, if $V^\bu$ is a complex of right $A$\module s such that
the tensor product of $V^\bu$ with any coacyclic complex of left
$A$\module s is acyclic, then the tensor product of $V^\bu$ with any
acyclic complex $U^\bu$ of left $A$\module s is also acyclic, since
one can construct a morphism into $U^\bu$ from an acyclic complex
of flat $A$\module s with a coacyclic cone.)
 The complex of $\C$\comodule s $V^\bu\ot_A\C$ coinduced from
a flat complex of $A$\module s $V^\bu$ is coflat.

\begin{rmk}
 The coderived category $\sD^\co(\C\comodl)$ can be only thought of
as the ``right'' version of exotic unbounded derived category of
$\C$\comodule s (e.~g., for the purposes of defining the derived
functors $\Cotor^\C$ and $\Coext_\C$, constructing the equivalence
of derived categories of $\C$\comodule s and $\C$\contramodule s,
etc.)\ when the ring $A$ has a finite (weak or left) homological
dimension.
 Indeed, what is needed is a definition of ``relative coderived
category'' of $\C$\comodule s such that for $\C=A$ it would
coincide with the derived category of $A$\module s, while when
$\C$ is a coalgebra over a field it would be the coderived category
of $\C$\comodule s defined above.
 (The same applies to the semiderived category $\sD^\si(\S\simodl)$
of $\S$\semimodule s---it only appears to be the ``right''
definition when the ring $A$ has a finite homological dimension.)
\end{rmk}

\subsection{Semiderived categories}
 Let $\S$ be a semialgebra over a coring~$\C$.
 Assume that $\C$ is a flat right $A$\module{} and $\S$ is a coflat
right $\C$\comodule, so that the category of left $\S$\semimodule s
is abelian.
 The \emph{semiderived category} of left $\S$\semimodule s
$\sD^\si(\S\simodl)$ is defined as the quotient category of
the homotopy category $\Hot(\S\simodl)$ by the thick subcategory
$\Acycl^{\cod\C}(\S\simodl)$ of complexes of $\S$\semimodule s
that are \emph{coacyclic as complexes of\/ $\C$\comodule s}. 

\begin{rmk}
 There is no claim that the semiderived category \emph{exists}
in the sense that morphisms between a given pair of objects form
a set rather than a class.
 Rather, we think of our localizations of categories as of ``very
large'' categories with classes of morphisms instead of sets.
 We will explain in~\ref{co-contra-ctrtor-definition}
and~\ref{semi-ctrtor-definition} how to compute the modules of
homomorphisms in semiderived categories in terms of resolutions;
then it will follow that the semiderived category does exist,
under certain assumptions.
\end{rmk}

\subsection{Semiflat complexes}   \label{semiflat-complexes}
 Let $\S$ be a semialgebra.
 The semitensor product $\bN^\bu\os_\S\bM^\bu$ of a complex of right
$\S$\semimodule s $\bN^\bu$ and a complex of left $\S$\semimodule s
$\bM^\bu$ is defined as the total complex of the bicomplex
$\bN^i\os_\S\bM^j$, constructed by taking infinite direct sums
along the diagonals.
 Of course, appropriate (co)flatness conditions must be imposed on
$\S$, $\bN^\bu$, and $\bM^\bu$ for this definition to make sense.

 Assume that the coring $\C$ is a flat right $A$\module,
the semialgebra $\S$ is a coflat right $\C$\comodule{} and a
$\C/A$\+coflat left $\C$\comodule, and the ring $A$ has a finite
weak homological dimension.
 A complex of $A$\+flat right $\S$\semimodule s $\bN^\bu$
is called \emph{semiflat} if the complex $\bN^\bu\os_\S\bM^\bu$
is acyclic whenever a complex of left $\S$\semimodule s $\bM^\bu$
is $\C$\coacyclic.
 Any semiflat complex of $\S$\semimodule s is a coflat complex of
$\C$\comodule s. 
 The complex of $\S$\semimodule s $\R^\bu\oc_\C\S$ induced from
a coflat complex of $A$\+flat $\C$\comodule s $\R^\bu$ is
semiflat.

 If it is only known that $\C$ is a flat right $A$\module{} and
$\S$ is a coflat right $\C$\comodule, one can define semiflat
complexes of $\C$\+coflat right $\S$\semimodule s.
 Then the complex of $\S$\semimodule s induced from a complex of
coflat $\C$\comodule s is semiflat; it is also a complex of
semiflat semimodules.

 Notice that \emph{not every complex of semiflat semimodules is
semiflat} (see~\ref{prelim-unbounded-tor}).
 In particular, it follows from Theorem~\ref{semitor-main-theorem}
and Lemma~\ref{semitor-definition} below that (in the assumptions
of~\ref{semitor-main-theorem}) a $\C$\coacyclic{} complex of
$A$\+flat right $\S$\semimodule s $\bN^\bu$ is semiflat if and only if
its semitensor product with any complex of left $\S$\semimodule s
$\bM^\bu$ (or just with any left $\S$\semimodule{} $\bM$) is acyclic.
 Thus a $\C$\coacyclic{} complex of semiflat $\S$\semimodule s
is semiflat if and only if all of its semimodules of cocycles
are semiflat.

 On the other hand, any complex of semiflat semimodules bounded
from above is semiflat.
 Moreover, if $\dsb\to\bN^{-1,\bu}\to\bN^{0,\bu}\to0$ is a complex,
bounded from above, of semiflat complexes of $\S$\semimodule s,
then the total complex $\bE^\bu$ of the bicomplex $\bN^{\bu,\bu}$
constructed by taking infinite direct sums along the diagonals
is semiflat.
 Indeed, the category of semiflat complexes is closed under shifts,
cones, and infinite direct sums, so one can apply the following
Lemma.

\begin{lem}
 Let\/ $\dsb\to N^{-1,\bu}\to N^{0,\bu}\to0$ be a complex, bounded
from above, of arbitrary complexes over an additive category\/~$\sA$
where infinite direct sums exist.
 Then the total complex $E^\bu$ of the bicomplex $N^{\bu,\bu}$
up to the homotopy equivalence can be obtained from the complexes
$N^{-i,\bu}$ using the operations of shift, cone, and infinite
direct sum.
\end{lem}

\begin{proof}
 Let $E_n^\bu$ be the total complex of the finite complex of
complexes $0\to N^{-n,\bu}\to\dsb\to N^{0,\bu}\to0$.
 Then the complex $E^\bu$ is the inductive limit of the complexes
$E_n^\bu$, and in addition, the embeddings of complexes
$E_n^\bu\rarrow E_{n+1}^\bu$ split in every degree.
 Thus the triple of complexes $\bigoplus_n E_n^\bu\rarrow
\bigoplus_n E_n^\bu\rarrow E^\bu$ is split exact in every degree
and the complex $E^\bu$ is homotopy equivalent to the cone of
the morphism $\bigoplus_n E_n^\bu\rarrow \bigoplus_n E_n^\bu$
(the homotopy inductive limit of the complexes~$E_n^\bu$).
\end{proof}

\subsection{Main theorem for comodules}   \label{cotor-main-theorem}
 Assume that the coring $\C$ is a flat left and right $A$\module{}
and the ring $A$ has a finite weak homological dimension.

\begin{thm}
 The functor mapping the quotient category of the homotopy category
of complexes of coflat\/ $\C$\comodule s (coflat complexes of\/
$\C$\comodule s) by its intersection with the thick subcategory of
coacyclic complexes of\/ $\C$\comodule s into the coderived category
of\/ $\C$\comodule s is an equivalence of triangulated categories.
\end{thm}

\begin{proof}
 We will show that any complex of $\C$\comodule s $\K^\bu$ can be
connected with a complex of coflat $\C$\comodule s in a functorial
way by a chain of two morphisms $\K^\bu\larrow\boR_2(\K^\bu)\rarrow
\boR_2\boL_1(\K^\bu)$ with coacyclic cones.
 Moreover, if the complex $\K^\bu$ is a complex of coflat
$\C$\comodule s (coflat complex of $\C$\comodule s), then
the intermediate complex $\boR_2(\K^\bu)$ in this chain is also
a complex of coflat $\C$\comodule s (coflat complex of
$\C$\comodule s).
 Then we will apply the following Lemma.

\begin{lem}
 Let\/ $\sC$ be a category and\/ $\sF$ be its full subcategory.
 Let\/ $\sS$ be a class of morphisms in\/~$\sC$ containing the third
morphism of any triple of morphisms $s$, $t$, and $st$ when it
contains two of them.
 Suppose that for any object $X$ in\/ $\sC$ there is a chain of
morphisms $X\from F_1(X)\to\dsb\from F_{n-1}(X)\to F_n(X)$ belonging
to\/ $\sS$ and functorially depending on~$X$ such that the object
$F_n(X)$ belongs to\/ $\sF$ for any $X\in\sC$ and all the objects
$F_i(X)$ belong to\/ $\sF$ for any $X\in\sF$.
 Then the functor\/ $\sF[(\sS\cap\sF)^{-1}]\rarrow\sC[\sS^{-1}]$
induced by the embedding\/ $\sF\rarrow \sC$ is an equivalence
of categories.
\end{lem}

\begin{proof}
 It is obvious that the functor between the localized categories
is surjective on the isomorphism classes of objects; let us show that
it is bijective on morphisms.
 It follows from the condition on the class $\sS$ that the functors
$F_i$ preserve it.
 Let $U$ and $V$ be two objects of~$\sF$ and $\phi\:U\rarrow V$
be a morphism between them in the category $\sC[\sS^{-1}]$.
 Applying the functor $F_n\:\sC\rarrow\sF$, we obtain a morphism
$F_n(\phi)\:F_n(U)\rarrow F_n(V)$ in the category
$\sF[(\sS\cap\sF)^{-1}]$.
 The square diagram of morphisms in the category $\sC[\sS^{-1}]$ formed
by the morphism $\phi$, the isomorphism between $U$ and $F_n(U)$,
the morphism $F_n(\phi)$, and the isomorphism between $V$ and
$F_n(V)$ is commutative, since it is composed from commutative
squares of morphisms in the category $\sC$.
 Since the other three morphisms in this commutative square lift to
$\sF[(\sS\cap\sF)^{-1}]$, the morphism $\phi$ belongs to the image
of the functor $\sF[(\sS\cap\sF)^{-1}]\rarrow\sC[\sS^{-1}]$.
 Now suppose that two morphisms $\phi$ and $\psi\:U\rarrow V$
in the category $\sF[(\sS\cap\sF)^{-1}]$ map to the same morphism
in $\sC[\sS^{-1}]$.
 Applying the functor $F_n$, we see that the morphisms $F_n(\phi)$
and $F_n(\psi)$ are equal in $\sF[(\sS\cap\sF)^{-1}]$.
 So we have two commutative squares in the category
$\sF[(\sS\cap\sF)^{-1}]$ with the same vertices $U$, $V$, \ 
$F_n(U)$, and $F_n(V)$, the same morphism $F_n(U)\rarrow F_n(V)$,
the same isomorphisms $U\simeq F(U)$ and $V\simeq F(V)$, and
two morphisms $\phi$ and $\psi\:U\rarrow V$.
 It follows that the latter two morphisms are equal.
\end{proof}

 Let $\K^\bu$ be a complex of $\C$\comodule s.
 Let $\cP(\M)\rarrow\M$ denote the functorial surjective morphism onto
an arbitrary $\C$\comodule{} $\M$ from an $A$\+flat $\C$\comodule{}
$\cP(\M)$ constructed in Lemma~\ref{flat-comodule-surjection}.

 The functor~$\cP$ is not always additive, but as any functor
from an additive category to an abelian one it is the direct sum of
a constant functor $\M\mpsto\cP(0)$ and a functor $\cP^+(\M) =
\ker(\cP(\M)\to\cP(0)) = \coker(\cP(0)\to\cP(\M))$ sending
zero objects to zero objects and zero morphisms to zero morphisms.
 For any $\C$\comodule{} $\M$, the comodule $\cP^+(\M)$ is
$A$\+flat and the morphism $\cP^+(\M)\rarrow\M$ is surjective.

 Set $\cP_0(\K^\bu) = \cP^+(\K^\bu)$, \ $\cP_1(\K^\bu) =
\cP^+(\ker(\cP^0(\K^\bu)\to\K^\bu))$, etc.
 For $d$ large enough, the kernel $\cZ(\K^\bu)$ of the morphism
$\cP_{d-1}(\K^\bu)\rarrow\cP_{d-2}(\K^\bu)$ will be a complex
of $A$\+flat $\C$\comodule s.
 Let $\boL_1(\K^\bu)$ be the total complex of the bicomplex
$$
 \cZ(\K^\bu)\lrarrow\cP_{d-1}(\K^\bu)
 \lrarrow\dsb\lrarrow\cP_1(\K^\bu)\lrarrow\cP_0(\K^\bu).
$$
 Then $\boL_1(\K^\bu)$ is a complex of $A$\+flat $\C$\comodule s
and the cone of the morphism $\boL_1(\K^\bu)\rarrow\K^\bu$ is
the total complex of a finite exact sequence of complexes of
$\C$\comodule s, and therefore, a coacyclic complex.

 Now let $\L^\bu$ be a complex of $A$\+flat left $\C$\comodule s.
 Consider the cobar construction
$$
 \C\ot_A\L^\bu \lrarrow \C\ot_A\C\ot_A\L^\bu \lrarrow
 \C\ot_A\C\ot_A\C\ot_A\L^\bu\lrarrow\dsb
$$
 Let $\boR_2(\L^\bu)$ be the total complex of this bicomplex,
constructed by taking infinite direct sums along the diagonals.
 Then $\boR_2(\L^\bu)$ is a complex of coflat $\C$\comodule s.
 The functor $\boR_2$ can be extended to arbitrary complexes of
$\C$\comodule s; for any complex $\K^\bu$, the cone of
the morphism $\K^\bu\rarrow\boR_2(\K^\bu)$ is coacyclic
by Lemma~\ref{coderived-categories}.

 Finally, if $\K^\bu$ is a coflat complex of $\C$\comodule s,
then $\boR_2(\K^\bu)$ is also a coflat complex of $\C$\comodule s,
since the cotensor product of $\boR_2(\K^\bu)$ with a complex of right
$\C$\comodule s $\N^\bu$ coincides with the cotensor product of
$\K^\bu$ with the total cobar complex $\boR_2(\N^\bu)$, and the latter
is coacyclic whenever $\N^\bu$ is coacyclic.

 We have constructed the chain of morphisms $\K^\bu\larrow
\boR_2(\K^\bu)\rarrow \boR_2\boL_1(\K^\bu)$ with the desired
properties.
 The only remaining problem is that the functor $\boL_1$ is not
additive and therefore not defined on the homotopy category of
complexes of $\C$\comodule s, but only on the (abelian) category
of complexes and their morphisms.
 So we have to apply Lemma to the category $\sC$ of complexes of
$\C$\comodule s, the full subcategory $\sF$ of complexes of coflat
$\C$\comodule s (coflat complexes of $\C$\comodule s) in it, and
the class $\sS$ of morphisms with coacyclic cones.

 The corresponding localizations will coincide with the desired
quotient categories of homotopy categories due to the following
general fact~\cite[III.4.2-3]{GM}.
 For any DG\+category $\DG$ where shifts and cones exist
the localization of the category of closed morphisms in $\DG$ with
respect to the class of homotopy equivalences coincides with
the homotopy category of $\DG$ (i.~e., closed morphisms homotopic
in~$\DG$ become equal after inverting homotopy equivalences).
 In particular, this is true for any category of complexes over
an additive category that is closed under shifts and cones.
\end{proof}

\begin{rmk}
 Another proof of Theorem (for complexes of coflat comodules
or coflat complexes of $A$\+flat comodules) can be found
in~\ref{semitor-main-theorem}.
 After Theorem has been proven, it turns out that the functors $\boL_1$
and $\boR_2$ can be also applied in the reverse order: for any complex
of $\C$\comodule s $\L^\bu$, the complex $\boR_2(\L^\bu)$ is a complex
of $\C/A$\+coflat $\C$\comodule s, and for any complex of
$\C/A$\+coflat $\C$\comodule s $\K^\bu$, the complex $\boL_1(\K^\bu)$
is a complex of coflat $\C$\comodule s
(by Remark~\ref{absolute-relative-coflat}, which depends on Theorem).
\end{rmk}

\subsection{Main theorem for semimodules}  \label{semitor-main-theorem}
 Assume that the coring $\C$ is a flat left and right $A$\module,
the semialgebra $\S$ is a coflat left and right $\C$\comodule,
and the ring $A$ has a finite weak homological dimension.

\begin{thm}
 The functor mapping the quotient category of the homotopy category of
semiflat complexes of $A$\+flat (\.$\C$\+coflat, semiflat)\/
$\S$\semimodule s by its intersection with the thick subcategory of\/
$\C$\coacyclic{} complexes of\/ $\S$\semimodule s into the semiderived
category of\/ $\S$\semimodule s is an equivalence of triangulated
categories.
\end{thm}

\begin{proof}
 We will show that in the chain of functors mapping the quotient
category of (the homotopy category of) semiflat complexes of
$\C$\+coflat (semiflat) $\S$\semimodule s by $\C$\coacyclic{} semiflat
complexes of $\C$\+coflat $\S$\semimodule s into the quotient category
of complexes of $\C$\+coflat $\S$\semimodule s by $\C$\coacyclic{}
complexes of $\C$\+coflat $\S$\semimodule s into the quotient category
of complexes of $A$\+flat $\S$\semimodule s by $\C$\coacyclic{}
complexes of $A$\+flat $\S$\semimodule s into the semiderived category
of $\S$\semimodule s all the three functors are equivalences of
categories.
 Analogously, in the chain of functors mapping the quotient category
of (the homotopy category of) semiflat complexes of $A$\+flat
$\S$\semimodule s by $\C$\coacyclic{} semiflat complexes of $A$\+flat
$\S$\semimodule s into the quotient category of $\C$\+coflat complexes
of $A$\+flat $\S$\semimodule s by $\C$\coacyclic{} $\C$\+coflat
complexes of $A$\+flat $\S$\semimodule s into the quotient category
of complexes of $A$\+flat $\S$\semimodule s by $\C$\coacyclic{}
complexes of $A$\+flat $\S$\semimodule s into the semiderived category
of $\S$\semimodule s all the three functors are equivalences
of categories.

 In order to prove this, we will construct for any complex of
$\S$\semimodule s $\bK^\bu$ a morphism $\boL_1(\bK^\bu)\rarrow\bK^\bu$
into $\bK^\bu$ from a complex of $A$\+flat $\S$\semimodule s
$\boL_1(\bK^\bu)$, for any complex of $A$\+flat $\S$\semimodule s
$\bL^\bu$ a morphism $\bL^\bu\rarrow\boR_2(\bL^\bu)$ from $\bL^\bu$
into a complex of $\C$\+coflat $\S$\semimodule s $\boR_2(\bL^\bu)$,
and for any $\C$\+coflat complex of $A$\+flat $\S$\semimodule s
(complex of $\C$\+coflat $\S$\semimodule s) $\bM^\bu$ a morphism
$\boL_3(\bM^\bu)\rarrow\bM^\bu$ into $\bM^\bu$ from a semiflat complex
of $A$\+flat (semiflat) $\S$\semimodule s $\boL_3(\bM^\bu)$  such that
in each case the cone of this morphism will be a $\C$\+coacyclic
complex of $\S$\semimodule s.
 Then we will apply the following Lemma.

\begin{lem}
 Let\/ $\sH$ be a category and\/ $\sF$ be its full subcategory.
 Let\/ $\sS$ be a localizing (i.~e., satisfying the Ore
conditions) class of morphisms in\/~$\sH$.
 Assume that for any object $X$ of\/~$\sH$ there exists an object
$U$ of\/~$\sF$ together with a morphism $U\rarrow X$ belonging
to\/~$\sS$ (or for any object $X$ of\/ $\sH$ there exists an object
$V$ of\/~$\sF$ together with a morphism $X\rarrow V$ belonging
to\/~$\sS$).
 Then the functor\/ $\sF[(\sS\cap\sF)^{-1}]\rarrow\sH[\sS^{-1}]$
induced by the embedding\/ $\sF\rarrow \sH$ is an equivalence
of categories.
\end{lem}

\begin{proof}
 It is obvious that the functor between the localized categories
is surjective on the isomorphism classes of objects; let us show that
it is bijective on morphisms.
 Any morphism in the category $\sH[\sS^{-1}]$ between two objects $U$
and $V$ from~$\sF$ can be represented by a fraction $U\from X\to V$,
where $X$ is an object of~$\sH$ and the morphism $X\to U$
belongs to~$\sS$.
 By our assumption, there is an object $W$ from~$\sF$ together with
a morphism $W\to X$ from~$\sS$.
 Then the fractions $U\from X\to V$ and $U\from W\to V$ represent
the same morphism in $\sH[\sS^{-1}]$, while the second fraction
represents also a certain morphism in $\sF[(\sS\cap\sF)^{-1}]$.
 Furthermore, any two morphisms from an object~$U$ to an object~$V$
in the category $\sF[(\sS\cap\sF)^{-1}]$ can be represented by two
fractions of the form $U\from U'\birarrow V$, with the same morphism
$U\to U'$ from $\sS\cap\sF$ and two different morphisms $U'\birarrow V$
(since the class of morphisms $\sS\cap\sF$ in the category $\sF$
satisfies the right Ore conditions).
 If the images of these morphisms in the category $\sH[\sS^{-1}]$
are equal, then there is a morphism $X\to U'$ from~$\sS$ with
an object $X$ from~$\sH$ such that two compositions
$X\to U'\birarrow V$ coincide.
 Again there is an object $W$ from~$\sF$ together with a morphism
$W\to X$ belonging to~$\sS$.
 Since the two compositions $W\to U'\birarrow V$ coincide in~$\sF$,
the morphisms represented by the two fractions $U\from U'\birarrow V$
are equal in $\sF[(\sS\cap\sF)^{-1}]$.
\end{proof}

 Let $\bK^\bu$ be a complex of $\S$\semimodule s.
 Let $\bcP(\bM)\rarrow\bM$ denote the functorial surjective morphism
onto an arbitrary $\S$\semimodule{} $\bM$ from an $A$\+flat
$\S$\semimodule{} $\bcP(\bM)$ constructed in
Lemma~\ref{flat-semimodule-surjection}.
 As explained in the proof of Theorem~\ref{cotor-main-theorem},
the functor~$\bcP$ is the direct sum of a constant functor
$\bM\mpsto\bcP(0)$ and a functor $\bcP^+$ sending zero morphisms to
zero morphisms.
 For any $\S$\semimodule{} $\bM$, the semimodule $\bcP^+(\bM)$ is
$A$\+flat and the morphism $\bcP^+(\bM)\rarrow\bM$ is surjective.

 Set $\bcP_0(\bK^\bu) = \bcP^+(\bK^\bu)$, \ $\bcP_1(\bK^\bu) =
\bcP^+(\ker(\bcP^0(\bK^\bu)\to\bK^\bu))$, etc.
 For $d$ large enough, the kernel $\bcZ(\bK^\bu)$ of the morphism
$\bcP_{d-1}(\bK^\bu)\rarrow\bcP_{d-2}(\bK^\bu)$ will be a complex
of $A$\+flat $\S$\semimodule s.
 Let $\boL_1(\bK^\bu)$ be the total complex of the bicomplex
$$
 \bcZ(\bK^\bu)\lrarrow\bcP_{d-1}(\bK^\bu)
 \lrarrow\dsb\lrarrow\bcP_1(\bK^\bu)\lrarrow\bcP_0(\bK^\bu).
$$
 Then $\boL_1(\bK^\bu)$ is a complex of $A$\+flat $\S$\semimodule s
and the cone of the morphism $\boL_1(\bK^\bu)\rarrow\bK^\bu$ is
the total complex of a finite exact sequence of complexes of
$\S$\semimodule s, and therefore, a $\C$\coacyclic{} complex
(and even an $\S$\coacyclic{} complex).

 Now let $\bL^\bu$ be a complex of $A$\+flat $\S$\semimodule s.
 Let $\bM\rarrow\bcJ(\bM)$ denote the functorial injective morphism
from an arbitrary $A$\+flat $\S$\semimodule{} $\bM$ into
a $\C$\+coflat $\S$\semimodule{} $\bcJ(\bM)$ with an $A$\+flat cokernel
$\bcJ(\bM)/\bM$ constructed in Lemma~\ref{coflat-semimodule-injection}.
 Set $\bcJ^0(\bL^\bu)=\bcJ(\bL^\bu)$, \ $\bcJ^1(\bL^\bu) =
\bcJ(\coker(\bL^\bu\to\bcJ^0(\bL^\bu)))$, etc.
 Let $\boR_2(\bL^\bu)$ be the total complex of the bicomplex
$$
 \bcJ^0(\bL^\bu)\lrarrow\bcJ^1(\bL^\bu)\lrarrow\bcJ^2(\bL^\bu)
 \lrarrow\dsb,
$$
constructed by taking infinite direct sums along the diagonals.
 Then $\boR_2(\bL^\bu)$ is a complex of $\C$\+coflat $\S$\semimodule s
and the cone of the morphism $\bL^\bu\rarrow\boR_2(\bL^\bu)$ is
a $\C$\coacyclic{} (and even $\S$\coacyclic) complex by
Lemma~\ref{coderived-categories}.

 Finally, let $\bM^\bu$ be a $\C$\+coflat complex of $A$\+flat
left $\S$\semimodule s.
 Then the complex $\S\oc_\C\bM^\bu$ is a semiflat complex of
$A$\+flat $\S$\semimodule s.
 Moreover, if $\bM^\bu$ is a complex of $\C$\+coflat
$\S$\+semimodule s, then $\S\oc_\C\bM^\bu$ is a semiflat complex
of semiflat $\S$\semimodule s.
 Consider the bar construction
$$
 \dsb\lrarrow \S\oc_\C\S\oc_\C\S\oc_\C\bM^\bu\lrarrow
 \S\oc_\C\S\oc_\C\bM^\bu\lrarrow\S\oc_\C\bM^\bu.
$$
 Let $\boL_3(\bM^\bu)$ be the total complex of this bicomplex,
constructed by taking infinite direct sums along the diagonals.
 Then the complex $\boL_3(\bM^\bu)$ is semiflat
by~\ref{semiflat-complexes} and the cone of the morphism
$\boL_3(\bM^\bu)\rarrow\bM^\bu$ is not only $\C$\coacyclic,
but even $\C$\contractible{} (the contracting homotopy being
induced by the semiunit morphism $\C\rarrow\S$).
\end{proof}

\begin{rmk}
 It is clear that the constructions of complexes $\boR_2(\bL^\bu)$
and $\boL_3(\bM^\bu)$ can be applied to arbitrary complexes of
$\S$\semimodule s, with no (co)flatness conditions imposed on them.
 For example, an alternative way of proving Theorem is to show that
the functors mapping the quotient category of semiflat complexes
of $\C$\+coflat (semiflat) $\S$\semimodule s by $\C$\coacyclic{}
semiflat complexes into the quotient category of complexes of
$\C/A$\+coflat $\S$\semimodule s by $\C$\coacyclic{} complexes
into the semiderived category of $\S$\semimodule s are both
equivalences of categories.
 Indeed, for any complex of $\S$\semimodule s $\bL^\bu$ the complex
$\boR_2(\bL^\bu)$ is a complex of $\C/A$\+coflat $\S$\semimodule s
by Lemma~\ref{coflat-semimodule-injection} and for any complex of
$\C/A$\+coflat $\S$\semimodule s $\bK^\bu$ the complex
$\boL_1(\bK^\bu)$ is a complex of $\C$\+coflat $\S$\semimodule s
by Remark~\ref{flat-semimodule-surjection} (hence the complex
$\boL_3\boL_1(\bK^\bu)$ is a semiflat complex of semiflat
$\S$\semimodule s).
 Yet another useful approach to proving Theorem was presented
in~\ref{cotor-main-theorem}: any complex of $\S$\semimodule s $\bK^\bu$
can be connected with a semiflat complex of semiflat $\S$\semimodule s
in a functorial way by a chain of three morphisms $\bK^\bu\larrow
\boL_3(\bK^\bu)\rarrow \boL_3\boR_2(\bK^\bu) \larrow 
\boL_3\boR_2\boL_1(\bK^\bu)$ with $\C$\coacyclic{} cones, and when
$\bK^\bu$ is a semiflat complex of ($A$\+flat, $\C$\+coflat,
or semiflat) $\S$\semimodule s, all the complexes in this chain are
also semiflat complexes of ($A$\+flat, $\C$\+coflat, or semiflat)
$\S$\semimodule s.
\end{rmk}

\begin{qst}
 Is the quotient category of $\C$\+coflat complexes of
$\S$\semimodule s by the thick subcategory of $\C$\+coacyclic
$\C$\+coflat complexes equivalent to the semiderived category
of $\S$\semimodule s?
\end{qst}

\subsection{Derived functor SemiTor}  \label{semitor-definition}
 The following Lemma provides a general approach to double-sided
derived functors of (partially defined) functors of two arguments.

\begin{lem}
 Let\/ $\sH_1$ and\/ $\sH_2$ be two categories, $\sH$ be
a (not necessarily full) subcategory in\/ $\sH_1\times\sH_2$,
and\/ $\sS_1$ and\/ $\sS_2$ be localizing classes of morphisms
in\/ $\sH_1$ and\/ $\sH_2$.
 Let\/ $\sK$ be a category and\/ $\Theta\:\sH\rarrow\sK$ be a functor.
 Let\/ $\sF_1$ and\/ $\sF_2$ be subcategories in\/ $\sH_1$ and\/
$\sH_2$.
 Assume that both functors\/ $\sF_i[(\sS_i\cap\sF_i)^{-1}] \rarrow
\sH_i[\sS_i^{-1}]$ induced by the embeddings\/ $\sF_i\rarrow\sH_i$
are equivalences of categories and the subcategory\/ $\sH$ contains
both subcategories\/ $\sF_1\times \sH_2$ and\/ $\sH_1\times \sF_2$.
 Furthermore, assume that the morphisms\/ $\Theta(U,t)$ and\/
$\Theta(s,V)$ are isomorphisms in the category\/~$\sK$ for any
objects $U\in \sF_1$, $V\in\sF_2$ and any morphisms $s\in\sS_1$,
$t\in\sS_2$.
 Then the restrictions of the functor\/~$\Theta$ to the subcategories\/
$\sF_1\times \sH_2$ and\/ $\sH_1\times \sF_2$ factorize through
their localizations by their intersections with\/ $\sS_1\times\sS_2$,
so one can define derived functors\/ $\boD_1\Theta$, $\boD_2\Theta\:
\sH_1[\sS_1^{-1}]\times\sH_2[\sS_2^{-1}]\rarrow\sK$ by restricting
the functor $\Theta$ to these subcategories.
 Moreover, the derived functors\/ $\boD_1\Theta$ and\/ $\boD_2\Theta$
are naturally isomorphic to each other and therefore do not depend
on the choice of subcategories\/ $\sF_1$ and\/~$\sF_2$, provided
that both subcategories exist.
\end{lem}

\begin{proof}
 Let us show that for any morphism $s\in \sS_1\cap\sF_1$ and any
object $X\in\sH_2$ the morphism $\Theta(s,X)$ is an isomorphism
in~$\sK$.
 By assumptions of Lemma, the image of~$X$ in $\sH_2[\sS_2^{-1}]$
is isomorphic to the image of a certain object $V\in \sF_2$.
 First suppose that there exists a fraction $X\from Y\to V$
of morphisms from~$\sS_2$ connecting $X$ and $V$.
 Then both morphisms of morphisms $\Theta(s,Y)\rarrow\Theta(s,X)$
and $\Theta(s,Y)\rarrow\Theta(s,V)$ are isomorphisms of morphisms,
since the source and the target of~$s$ belong to $\sF_1$.
 Now the morphism $\Theta(s,X)$ is an isomorphism, because
the morphism $\Theta(s,V)$ is an isomorphism.
 In the general case, there exist a fraction $X\from Y\to V$
connecting $X$ and~$V$ and two morphisms $Y'\to Y$ and $V\to V'$
such that the morphism $Y\to X$ and two compositions
$Y'\to Y\to V$ and $Y\to V\to V'$ belong to~$\sS_2$.
 Then the compositions of morphisms of morphisms
$\Theta(s,Y')\rarrow\Theta(s,Y)\rarrow\Theta(s,V)$ and
$\Theta(s,Y)\rarrow\Theta(s,V)\rarrow\Theta(s,V')$ are
isomorphisms of morphisms, so the morphism of morphisms
$\Theta(s,Y)\rarrow\Theta(s,V)$ is both left and right invertible,
and therefore, is an isomorphism of morphisms.
 Since the morphism of morphisms $\Theta(s,Y)\rarrow\Theta(s,X)$
is also an isomorphism of morphisms and the morphism $\Theta(s,V)$
is an isomorphism, one can conclude that the morphism $\Theta(s,X)$
is also an isomorphism.

 Thus the derived functor $\boD_1\Theta$ is defined; it remains to
construct an isomorphism between $\boD_1\Theta$ and $\boD_2\Theta$.
 But the compositions of the functors $\boD_1\Theta$ and $\boD_2\Theta$
with the functor
$\sF_1[(\sS_1\cap\sF_1)^{-1}]\times \sF_2[(\sS_2\cap\sF_2)^{-1}]
\rarrow \sH_1[\sS_1^{-1}]\times \sH_2[\sS_2^{-1}]$
coincide by definition, and the latter functor is an equivalence
of categories.
\end{proof}

 Assume that the coring $\C$ is a flat left and right $A$\module,
the semialgebra $\S$ is a coflat left and right $\C$\comodule,
and the ring $A$ has a finite weak homological dimension.

 The double-sided derived functor $\SemiTor^\S$ on the Carthesian
product of the semiderived categories of right and left
$\S$\semimodule s is defined as follows.
 Consider the partially defined functor of semitensor product of
complexes of $\S$\semimodule s $\os_\S\:\Hot(\simodr\S)\times
\Hot(\S\simodl)\darrow \Hot(k\modl)$.
 This functor is defined on the full subcategory of the Carthesian
product of homotopy categories that consists of pairs
of complexes $(\N^\bu,\M^\bu)$ such that either $\N^\bu$ or $\M^\bu$
is a complex of $A$\+flat $\S$\semimodule s.
 Compose it with the functor of localization $\Hot(k\modl)\rarrow
\sD(k\modl)$ and restrict to the Carthesian product of the homotopy
category of semiflat complexes of $A$\+flat right $\S$\semimodule s and
the homotopy category of complexes of left $\S$\semimodule s.

 By the definition, the functor so obtained factorizes through
the semiderived category of left $\S$\semimodule s in the second
argument, and it follows from Theorem~\ref{semitor-main-theorem}
and the above Lemma that it factorizes through the quotient category
of the homotopy category of semiflat complexes of $A$\+flat right
$\S$\semimodule s by its intersection with the thick subcategory of
$\C$\coacyclic{} complexes in the first argument.

 Explicitly, let $\bN^\bu$ be a $\C$\coacyclic{} semiflat complex
of $A$\+flat right $\S$\semimodule s and $\bM^\bu$ be a complex
of left $\S$\semimodule s.
 Using the constructions from the proof of
Theorem~\ref{semitor-main-theorem}, connect $\bM^\bu$ with
a semiflat complex of $A$\+flat left $\S$\semimodule s $\bL^\bu$
by a chain of morphisms with $\C$\coacyclic{} cones.
 Then the complexes $\bN^\bu\os_\S\bM^\bu$ and $\bN^\bu\os_\S\bL^\bu$
are connected by a chain of quasi-isomorphisms, and since the complex
$\bN^\bu\os_\S\bL^\bu$ is acyclic, the complex $\bN^\bu\os_\S\bM^\bu$
is acyclic, too.

 Thus we have constructed the double-sided derived functor
$$
 \SemiTor^\S\:\sD^\si(\simodr\S)\times\sD^\si(\S\simodl)
 \lrarrow\sD(k\modl).
$$
 According to Lemma, the same derived functor can be obtained by
restricting the functor of semitensor product to the Carthesian product
of the homotopy category of complexes of left $\S$\semimodule s
and the homotopy category of semiflat complexes of $A$\+flat right
$\S$\semimodule s, or indeed, to the Carthesian product of the homotopy
categories of semiflat complexes of $A$\+flat right and
left $\S$\semimodule s.
 One can also use semiflat complexes of $\C$\+coflat $\S$\semimodule s
or semiflat complexes of semiflat $\S$\semimodule s instead of
semiflat complexes of $A$\+flat $\S$\semimodule s.

 In particular, when the coring $\C$ is a flat left and right
$A$\module{} and the ring $A$ has a finite weak homological dimension,
one defines the double-sided derived functor
$$
 \Cotor^\C\:\sD^\co(\comodr\C)\times\sD^\co(\C\comodl)
 \lrarrow \sD(k\modl)
$$
by composing the functor of cotensor product ${\oc_\C}\:\Hot(\comodr\C)
\times\Hot(\C\comodl) \rarrow \Hot(k\modl)$ with the functor of
localization $\Hot(k\modl)\rarrow \sD(k\modl)$ and restricting it to
the Carthesian product of the homotopy category of complexes of coflat
right $\C$\comodule s and the homotopy category of arbitrary complexes
of left $\C$\comodule s.
 The same derived functor is obtained by restricting the functor of
cotensor product to the Carthesian product of the homotopy category
of arbitrary complexes of right $\C$\comodule s and the homotopy
category of complexes of coflat left $\C$\comodule s, or indeed, to
the Carthesian product of the homotopy categories of coflat right
$\C$\comodule s and coflat left $\C$\comodule s.
 One can also use coflat complexes of $\C$\comodule s or coflat
complexes of $A$\+flat $\C$\comodule s instead of complexes of
coflat $\C$\comodule s.

\begin{rmk}
 One can define a version of derived functor $\Cotor$ without making
any homological dimension assumptions by considering pro-objects
in the spirit of~\cite{GK1,GK2}.
 Let $k\modl^\omega$ denote the category of pro-objects over
the category $k\modl$ that can be represented by countable filtered
projective systems of $k$\module s; this is an abelian tensor
category with exact functors of countable filtered projective limits
and a right exact functor of tensor product commuting with countable
filtered projective limits.
 Let $\cA$ be a ring object in $k\modl^\omega$; then one can
consider right and left $\cA$\module{} objects and
$\cA$\+$\cA$\bimodule{} objects in $k\modl^\omega$, which we will
simply call right and left $\cA$\module s and $\cA$\+$\cA$\bimodule s.
 Furthermore, let $\bC$ be a coring object in the tensor category
of $\cA$\+$\cA$\bimodule s; we will consider $\bC$\comodule{}
objects in the categories of right and left $\cA$\module s and
call them right and left $\bC$\comodule s.
 Define the functor of cotensor product over $\bC$ taking values in
the category $k\modl^\omega$ in the usual way and extend it to
the Carthesian product of the homotopy categories of complexes of
right and left $\bC$\comodule s by taking infinite products along
the diagonals in the bicomplex of cotensor products.
 The categories of right and left $\cA$\module s are abelian.
 Assume that $\bC$ is a flat left and right $\cA$\module; then
the categories of right and left $\bC$\comodule s are also abelian.
 Define the semiderived categories of right and left $\bC$\comodule s
as the quotient categories of the homotopy categories by the thick
subcategories of $\cA$\contraacyclic{} complexes (the contraacyclic
complexes being defined in terms of countable products).
 Then one can use Lemma to define the double-sided derived functor
$\ProCotor^\bC$ of cotensor product on the Carthesian product of
the semiderived categories of right and left $\bC$\comodule s in terms
of coflat complexes of $\bC$\comodule s.
 In order to obtain for any complex of $\bC$\comodule s $\bM^\bu$
a coflat complex of $\bC$\comodule s connected with $\bM^\bu$
by a functorial chain of two morphisms with $\cA$\contraacyclic{} 
cones one needs to construct a surjective morphism onto any
$\bC$\comodule{} $\bM$ from an $\cA$\+flat $\bC$\comodule{}
$\bcF(\bM)$.
 This construction is dual to that of
Lemma~\ref{coflat-semimodule-injection} and uses the surjective map
onto any $\cA$\module{} $\M$ from an $\cA$\+flat $\cA$\module{}
$\cG(\M)=\cA\ot_k^\omega M'$, where $M'$ is a pro-$k$-module
represented by a countable filtered projective system of flat
$k$\module s mapping onto the pro-$k$-module $\M$ and $\ot_k^\omega$
denotes the functor of tensor product in $k\modl^\omega$.
 The $\cA$\+flat $\bC$\comodule{} $\bcF(\bM)$ is obtained as
the projective limit in $k\modl^\omega$ of the projective system of
$\bC$\comodule s $\bM\larrow\bcQ(\bM)\larrow\bcQ(\bcQ(\bM))\larrow
\dsb$\,
 Given a complex of $\cA$\+flat $\bC$\comodule s $\bM^\bu$, a coflat
complex of $\bC$\comodule s endowed with a morphism from the complex
$\bM^\bu$ with an $\cA$\contractible{} cone is obtained as the total
complex of the cobar complex of $\bM^\bu$, constructed by taking
infinite products along the diagonals.
 One can also consider the category of arbitrary pro-$k$-modules
in place of $k\modl^\omega$.
 Notice that for a conventional coalgebra $\bC$ over a field $\cA=k$
and complexes of $\bC$\comodule s $\bN^\bu$ and $\bM^\bu$ in
the category of $k$\+vector spaces that are both bounded from above
or from below the object of the derived category of $k$\+vector
spaces obtained by applying the derived functor of projective limit
to the object $\ProCotor^\bC(\bN^\bu,\bM^\bu)$ of the derived
category $\sD(k\vect^\omega)$ coincides with $\Cotor^{\bC,I}
(\bN^\bu,\bM^\bu)$ (see~\ref{derived-first-second-kind}).
\end{rmk}

\subsection{Relatively semiflat complexes}  \label{relatively-semiflat}
 We keep the assumptions and notation of~\ref{cotor-main-theorem},
\ref{semitor-main-theorem}, and~\ref{semitor-definition}.

 One can compute the derived functor $\Cotor^\C$ using resolutions of
a different kind.
 Namely, the cotensor product $\N^\bu\oc_\C\M^\bu$ of a complex of
$A$\+flat right $\C$\comodule s $\N^\bu$ and a complex of
$\C/A$\+coflat $\C$\comodule s $\M^\bu$ represents an object naturally
isomorphic to $\Cotor^\C(\M^\bu,\N^\bu)$ in the derived category
of $k$\module s.
 Indeed, the complex $\boR_2(\N^\bu)$ is a complex of coflat
$\C$\comodule s and the cone of the morphism $\N^\bu\rarrow
\boR_2(\N^\bu)$ is coacyclic with respect to the exact category of
$A$\+flat right $\C$\comodule s, hence the morphism $\N^\bu\oc_\C\M^\bu
\rarrow \boR_2(\N^\bu)\oc_\C\M^\bu$ is a quasi-isomorphism.
 One can prove that the cotensor product of a complex coacyclic
with respect to the exact category of $A$\+flat $\C$\comodule s and
a complex of $\C/A$\+coflat $\C$\comodule s is acyclic in the way
completely analogous to the proof of Lemma~\ref{coflat-complexes}.

 One can also compute the derived functor $\SemiTor^\S$ using
resolutions of different kinds.
 Namely, a complex of left $\S$\semimodule s is called \emph{semiflat
relative to~$A$} if its semitensor product with any complex of
$A$\+flat right $\S$\semimodule s that as a complex of $\C$\comodule s
is coacyclic with respect to exact category of $A$\+flat right
$\C$\comodule s is acyclic
(cf.\ Theorem~\ref{co-contra-push-well-defined}(a)).
 For example, the complex of $\S$\semimodule s induced from a complex
of $\C/A$\+coflat $\C$\comodule s is semiflat relative to~$A$, hence
the complex $\boL_3\boR_2(\bK^\bu)$ is semiflat relative to~$A$ for
any complex of left $\S$\semimodule s $\bK^\bu$.
 The semitensor product $\bN^\bu\os_\S\bM^\bu$ of a complex of
$A$\+flat right $\S$\semimodule s $\bN^\bu$ and a complex of left
$\S$\semimodule s $\bM^\bu$ semiflat relative to~$A$ represents
an object naturally isomorphic to $\SemiTor^\S(\bN^\bu,\bM^\bu)$ in
the derived category of $k$\module s.
 Indeed, $\boL_3\boR_2(\bN^\bu)$ is a semiflat complex of right
$\S$\semimodule s connected with $\bN^\bu$ by a chain of morphisms
$\bN^\bu\rarrow\boR_2(\bN^\bu)\larrow\boL_3\boR_2(\bN^\bu)$ whose
cones are coacyclic with respect to the exact category of $A$\+flat
$\C$\comodule s and contractible over~$\C$, respectively.
 Hence there is a chain of two quasi-isomorphisms connecting
$\bN^\bu\os_\S\bM^\bu$ with $\boL_3\boR_2(\bN^\bu)\os_\S\bM^\bu$.

 Analogously, a complex of left $\S$\semimodule s is called
\emph{semiflat relative to\/~$\C$} if its semitensor product with any
$\C$\contractible{} complex of $\C$\+coflat right $\S$\semimodule s
is acyclic.
 For example, the complex of $\S$\semimodule s induced from any
complex of $\C$\comodule s is semiflat relative to~$\C$, hence
the complex $\boL_3(\bK^\bu)$ is semiflat relative to~$\C$ for any
complex of left $\S$\semimodule s $\bK^\bu$.
 The semitensor product $\bN^\bu\os_\S\bM^\bu$ of a complex of
$\C$\+coflat right $\S$\semimodule s $\bN^\bu$ and a complex of
left $\S$\semimodule s $\bM^\bu$ semiflat relative to~$\C$ represents
an object naturally isomorphic to $\SemiTor^\S(\bN^\bu,\bM^\bu)$
in the derived category of $k$\module s.
 Indeed, $\boL_3(\bN^\bu)$ is a semiflat complex of right
$\S$\semimodule s and the cone of the morphism $\boL_3(\bN^\bu)
\rarrow\bN^\bu$ is a $\C$\contractible{} complex of $\C$\+coflat
right $\S$\semimodule s.
 It follows that the semitensor product of a complex of left
$\S$\semimodule s semiflat relative to~$\C$ with a $\C$\coacyclic{}
complex of $\C$\+coflat right $\S$\semimodule s is acyclic.

 At last, a complex of $A$\+flat right $\S$\semimodule s is called
\emph{semiflat relative to\/ $\C$ relative to~$A$} ($\S/\C/A$\+semiflat)
if its semitensor product with any $\C$\contractible{} complex of
$\C/A$\+coflat left $\S$\semimodule s is acyclic.
 For example, the complex of $\S$\semimodule s induced from a complex
of $A$\+flat $\C$\comodule s is $\S/\C/A$\+semiflat, hence the complex
$\boL_3\boL_1(\bK^\bu)$ is $\S/\C/A$\+semiflat for any complex of
right $\S$\semimodule s $\bK^\bu$.
 The semitensor product $\bN^\bu\os_\S\bM^\bu$ of
an $\S/\C/A$\+semiflat complex of $A$\+flat right $\S$\semimodule s
$\bN^\bu$ and a complex of $\C/A$\+coflat left $\S$\semimodule s
$\bM^\bu$ represents an object naturally isomorphic to
$\SemiTor^\S(\bN^\bu,\bM^\bu)$ in the derived category of $k$\module s.
 Indeed, $\boL_3(\bM^\bu)$ is a complex of left $\S$\semimodule s
semiflat relative to~$A$ and the cone of the morphism $\boL_3(\bM^\bu)
\rarrow\bM^\bu$ is a $\C$\contractible{} complex of $\C/A$\+coflat
right $\S$\semimodule s.
 It follows that the semitensor product of an $\S/\C/A$\semiflat{}
complex of $A$\+flat right $\S$\semimodule s with a $\C$\coacyclic{}
complex of $\C/A$\+coflat left $\S$\semimodule s is acyclic.

 The functors mapping the quotient categories of the homotopy
categories of complexes of $\S$\semimodule s semiflat relative to~$A$,
complexes of $\S$\semimodule s semiflat relative to~$\C$, and
$\S/\C/A$\+semiflat complexes of $A$\+flat $\S$\semimodule s by
their intersections with the thick subcategory of $\C$\coacyclic{}
complexes into the semiderived category of $\S$\semimodule s are
equivalences of triangulated categories.
 The same applies to complexes of $A$\+flat, $\C$\+coflat, or
$\C/A$\+coflat $\S$\semimodule s.
 These results follow easily from either of
Lemmas~\ref{cotor-main-theorem} or~\ref{semitor-main-theorem}.
 So one can define the derived functor $\SemiTor^\S$ by restricting
the functor of semitensor product to these categories of complexes
of $\S$\semimodule s as explained above.

\begin{rmk}
 Assuming that $\C$ is a flat right $A$\module, $\S$ is a coflat right
and a $\C/A$\+coflat left $\C$\comodule, and $A$ has a finite weak
homological dimension, one can define the double-sided derived functor
$\SemiTor^\S$ on the Carthesian product of the semiderived category
of $A$\+flat right $\S$\semimodule s and the semiderived category of
left $\S$\semimodule s.
 The former is defined as the quotient category of the homotopy
category of complexes of $A$\+flat right $\S$\semimodule s by
the thick subcategory of complexes that as complexes of $\C$\comodule s
are coacyclic with respect to the exact category of $A$\+flat right
$\C$\comodule s.
 The derived functor is constructed by restricting the functor of
semitensor product to the Carthesian product of the homotopy category
of complexes of $A$\+flat right $\S$\semimodule s and the homotopy
category of complexes of left $\S$\semimodule s semiflat relative
to~$A$, or the Carthesian product of the homotopy category of semiflat
complexes of $A$\+flat right $\S$\semimodule s and the homotopy
category of complexes of left $\S$\semimodule s.
 Assuming that $\C$ is a flat left and right $A$\module, $\S$ is
a flat left $A$\module{} and a coflat right $\C$\comodule, and $A$ has
a finite weak homological dimension, one can define the left derived
functor $\SemiTor^\S$ on the Carthesian product of the semiderived
category of $\C/A$\+coflat right $\S$\semimodule s and the semiderived
category of left $\S$\semimodule s.
 The former is defined as the quotient category of the homotopy
category of complexes of $\C/A$\+flat right $\S$\semimodule s by
the thick subcategory of complexes that as complexes of $\C$\comodule s
are coacyclic with respect to the exact category of $\C/A$\+coflat
right $\C$\comodule s (cf.\ Remark~\ref{co-contra-push-well-defined}).
 The derived functor is constructed by restricting the functor of
semitensor product to the Carthesian product of the homotopy category
of complexes of $\C/A$\+coflat right $\S$\semimodule s and the homotopy
category of $\S/\C/A$\semiflat{} complexes of $A$\+flat left
$\S$\semimodule s, or the Carthesian product of the homotopy category
of semiflat complexes of $\C$\+coflat right $\S$\semimodule s and
the homotopy category of complexes of left $\S$\semimodule s.
 Both of these definitions of derived functors are particular cases of
Lemma~\ref{semitor-definition}.
\end{rmk}

\subsection{Remarks on derived semitensor product of bisemimodules}
\label{remarks-derived-semitensor-bi}
 We would like to define the double-sided derived functor of
semitensor product of bisemimodules and in such a way that derived
semitensor products of several factors would be associative.
 It appears that there are two approaches to this problem, even
in the case of modules over rings.
 First suppose that we only wish to have associative derived
semitensor products of three factors.
 Let $\S$ be a semialgebra over a coring~$\C$ and $\T$ be
a semialgebra over a coring~$\D$, both satisfying the conditions
of~\ref{semitor-main-theorem}.

 The semiderived category of $\S$\+$\T$\bisemimodule s
$\sD^\si(\S\bsimod\T)$ is defined as the quotient category of
the homotopy category $\Hot(\S\bsimod\T)$ by the thick subcategory
of complexes of bisemimodules that as complexes of
$\C$\+$\D$\bicomodule s are coacyclic with respect to
the abelian category of $\C$\+$\D$\bicomodule s.
 We would like to define derived functors of semitensor product
\begin{alignat*}{3}
 &\os_\S^\boD \: & &\sD^\si(\simodr\S)\times\sD^\si(\S\bsimod\T)
    & &\lrarrow \sD^\si(\simodr\T) \\
 &\os_\T^\boD \: & &\sD^\si(\S\bsimod\T)\times\sD^\si(\T\simodl)
    & &\lrarrow \sD^\si(\S\simodl)
\end{alignat*}
and prove the associativity isomorphism
 $$
  \SemiTor^\T(\bN^\bu\os_\S^\boD\bK^\bu\;\bM^\bu)\simeq
  \SemiTor^\S(\bN^\bu\;\bK^\bu\os_\T^\boD\bM^\bu).
 $$

 Let us call a complex of $\C$\+coflat right $\S$\semimodule s
\emph{quite semiflat} if it belongs to the minimal triangulated
subcategory of the homotopy category of $\S$\semimodule s
containing the complexes induced from complexes of coflat right
$\C$\comodule s and closed under infinite direct sums.
 One can show (see Remark~\ref{co-contra-push-well-defined} and
the proof of Theorem~\ref{semi-pull-push-adjusted}.2) that
the quotient category of the category of quite semiflat complexes of
$\C$\+coflat $\S$\semimodule s by its minimal triangulated subcategory
containing the complexes of $\S$\semimodule s induced from complexes
of $\C$\comodule s coacyclic with respect to the exact category of
$\C$\+coflat $\C$\comodule s and closed under infinite direct sums
is equivalent to the semiderived category of $\S$\semimodule s.
 In other words, any $\C$\coacyclic{} quite semiflat complex of
$\C$\+coflat $\S$\semimodule s can be obtained from the complexes of
$\S$\semimodule s induced from the total complexes of exact triples
of complexes of coflat $\C$\comodule s using the operations of
cone and infinite direct sum.

 It follows (by Lemmas~\ref{coflat-complexes}
and~\ref{absolute-relative-coflat}) that the restriction of the functor
of semitensor product $\Hot(\simodr\S)\times\Hot(\S\bsimod\T)\darrow
\sD^\si(\simodr\T)$ to the Carthesian product of the homotopy category
of quite semiflat complexes of $\C$\+coflat right $\S$\semimodule s and
the homotopy category of complexes of $\S$\+$\T$\bisemimodule s
factorizes through the Carthesian product of semiderived categories
of right $\S$\semimodule s and $\S$\+$\T$\bisemimodule s.
 So the desired derived functors are defined; and the associativity
isomorphism follows from Proposition~\ref{semitensor-associative}.
 Notice that this definition of a double-sided derived functor is
\emph{not} a particular case of the construction of
Lemma~\ref{semitor-definition}.

\begin{qst}
 Can one use arbitrary semiflat complexes of $\C$\+coflat
$\S$\semimodule s or, at least, semiflat complexes of semiflat
$\S$\semimodule s instead of quite semiflat complexes in this
construction?
 In other words, assume that $\bN^\bu$ is a $\C$\coacyclic{} semiflat
complex of semiflat right $\S$\semimodule s and $\bK$ is
an $\S$\+$\T$\bisemimodule.
 Is the complex $\bN^\bu\os_\S\bK$ necessarily $\D$\coacyclic?
(Cf.~\ref{remarks-derived-semihom-bi}.)
\end{qst}

 Now suppose that we want to have derived semitensor products of
any number of factors.
 Let $\S$ be a semialgebra over a coring $\C$ over a $k$\+algebra $A$, 
\ $\T$ be a semialgebra over a coring $\D$ over a $k$\+algebra $B$,
and $\bR$ be a semialgebra over a coring $\E$ over a $k$\+algebra $F$,
all three satisfying the conditions of~\ref{semitor-main-theorem}.
 We would like to define the derived functor of semitensor product
$$
  \os_\bR^\boD\: \sD^\si(\S\bsimod\bR)\times\sD^\si(\bR\bsimod\T)
  \lrarrow \sD^\si(\S\bsimod\T).
$$
 This can be done, assuming that the $k$\+algebras $A$, $B$, and $F$
are flat $k$\module s.

 Let us call a complex of $F$\+flat $\S$\+$\bR$\bisemimodule s
\emph{strongly\/ $\bR$\+semiflat} if its semitensor product over~$\bR$
with any $\E$\+$\D$\coacyclic{} complex of $\bR$\+$\T$\bisemimodule s
is a $\C$\+$\D$\+coacyclic{} complex of $\S$\+$\T$\bisemimodule s for
any semialgebra~$\T$.
 Using bimodule versions of the constructions of Lemmas
\ref{flat-semimodule-surjection} and~\ref{coflat-semimodule-injection},
one can prove that the quotient category of the homotopy category of
strongly $\bR$\+semiflat complexes of $F$\+flat
$\S$\+$\bR$\bisemimodule s by its intersection with the thick
subcategory of $\C$\+$\E$\coacyclic{} bisemimodules is equivalent
to the semiderived category of $\S$\+$\bR$\bisemimodule s, and
the analogous result holds for the homotopy category of strongly
$\S$\+semiflat and strongly $\bR$\+semiflat complexes of
$A$\+flat and $F$\+flat $\S$\+$\bR$\bisemimodule s.
 One just uses the functor $G(M)=\bigoplus_{m\in M}A\ot_k F$
in the construction of Lemma~\ref{flat-comodule-surjection}, considers
the bicoaction and bisemiaction morphisms in place of the coaction
and semiaction morphisms, etc.
 (As we only want our $A$\+$F$\bimodule s to be $A$\+flat and $F$\+flat,
no assumption about the homological dimension of $A\ot_k F$ is needed.)
 So Lemma~\ref{semitor-definition} is applicable to the functor of
semitensor product $\Hot(\S\bsimod\bR)\times\Hot(\bR\bsimod\T)\darrow
\sD^\si(\S\bsimod\T)$ and we obtain the desired double-sided
derived functor.
 There is an associativity isomorphism
$(\bN^\bu\os_\S^\boD\bK^\bu)\os_\T^\boD\bM^\bu \simeq
\bN^\bu\os_\S^\boD(\bK^\bu\os_\T^\boD\bM^\bu)$.

 In the case of derived cotensor product of bicomodules, one does
not need to introduce quite coflat or strongly coflat complexes.
 It suffices to consider complexes of $\C$\+coflat $\C$\comodule s or
complexes of ($\C$\+coflat and) $\E$\+coflat $\C$\+$\E$\bicomodule s.
 One can define double-sided derived functors
\begin{alignat*}{3}
 &\oc_\C^\boD \: & &\sD^\co(\comodr\C)\times\sD^\co(\C\bcomod\D)
    & &\lrarrow \sD^\co(\comodr\D) \\
 &\oc_\D^\boD \: & &\sD^\co(\C\bcomod\D)\times\sD^\co(\D\comodl)
    & &\lrarrow \sD^\co(\C\comodl)
\end{alignat*}
and prove the associativity isomorphism
 $$
  \Cotor^\D(\N^\bu\oc_\C^\boD\K^\bu\;\M^\bu)\simeq
  \Cotor^\C(\N^\bu\;\K^\bu\oc_\D^\boD\M^\bu)
 $$
by replacing the complex of right $\C$\comodule s $\N^\bu$ with
a complex of coflat right $\C$\comodule s and the complex of left
$\D$\comodule s $\M^\bu$ by a complex of coflat left $\D$\comodule s
representing the same object in the coderived category of comodules.
 The derived functors $\oc_\C^\boD$ and $\oc_\D^\boD$ are well-defined,
since any coacyclic complex of coflat comodules is coacyclic with
respect to the exact category of coflat comodules
(see~\ref{co-contra-push-well-defined}).
 If the $k$\module s $A$ and $F$ are flat, one can prove that
the quotient category of the homotopy category of $\E$\+coflat
$\C$\+$\E$\bicomodule s by its intersection with the thick
subcategory of coacyclic complexes of $\C$\+$\E$\bicomodule s is
equivalent to the coderived category of bicomodules, and the same
applies to the homotopy category of $\C$\+coflat and $\E$\+coflat
$\C$\+$\E$\bicomodule s.
 Then one can apply Lemma~\ref{semitor-definition} in order to define
the double-sided derived functor
$$
  \oc_\E^\boD\: \sD^\co(\C\bcomod\E)\times\sD^\co(\E\bcomod\D)
  \lrarrow \sD^\co(\C\bcomod\D)
$$
and there is an associativity isomorphism
$(\N^\bu\oc_\C^\boD\K^\bu)\oc_\D^\boD\M^\bu \simeq
\N^\bu\oc_\C^\boD(\K^\bu\oc_\D^\boD\M^\bu)$.

\Section{Semicontramodules and Semihomomorphisms}

 Throughout Sections 3--11,\, $k\dual$~is an injective cogenerator of
the category of $k$\module s.
 One can always take $k\dual=\Hom_\boZ(k,\boQ/\boZ)$.

\subsection{Contramodules}
 For two $k$\+algebras $A$ and $B$, we will denote by $A\bimod B$
the category of $k$\module s with an $A$\+$B$\bimodule{} structure.

\subsubsection{}
 The identity $\Hom_A(K\ot_A M\;P)\simeq\Hom_A(M,\Hom_A(K,P)$ for
left $A$\module s $M$, $P$ and an $A$\+$A$\bimodule{} $K$ means that
the category opposite to the category of left $A$\module s is a right
module category over the tensor category of $A$\+$A$\bimodule s with
the functor of right action $(N,P^{\op})\mpsto\Hom(N,P)^\op$.
 Therefore, one can consider module objects in this module category
over ring objects in $A\bimod A$ an comodule objects in this module
category over coring objects in $A\bimod A$.

 Clearly, a ring object $B$ in $A\bimod A$ is just a $k$\+algebra
endowed with a $k$\+algebra morphism $A\rarrow B$.
 A $B$\module{} in $A\modl^\op$ is an $A$\module{} $P$ endowed with
a map $P\rarrow\Hom_A(B,P)$; so one can easily see that $B$\module s
in $A\modl^\op$ are just (objects of the category opposite to
the category of) usual left $B$\module s.

 Let $\C$ be a coring over $A$.
 The category of \emph{left contramodules} over $\C$ is the opposite
category to the category of comodule objects in the right module
category $A\modl^\op$ over the coring object $\C$ in the tensor
category $A\bimod A$.
 In other words, a left $\C$\contramodule{} $\P$ is a left $A$\module{}
endowed with a \emph{left contraaction} map $\Hom_A(\C,\P)\rarrow\P$,
which should be a morphism of left $A$\module s satisfying
the following \emph{contraassociativity} and \emph{counity} equations.
 First, two maps from the module $\Hom_A(\C\ot_A\C\;\P)=
\Hom_A(\C,\Hom_A(\C,\P)$ to the module $\Hom_A(\C,\P)$, one of which 
is induced by the comultiplication map of~$\C$ and the other by
the contraaction map, should have equal compositions with
the contraaction map $\Hom_A(\C,\P)\rarrow\P$, and second,
the composition $\P=\Hom_A(A,\P)\rarrow\Hom_A(\C,\P) \rarrow\P$ of
the map induced by the counit map of~$\C$ with the contraaction map
should be equal to the identity map of~$\P$.
 A \emph{right contramodule} $\gR$ over~$\C$ is a right $A$\module{}
endowed with a \emph{right contraaction} map $\Hom_{A^\rop}(\C,\gR)
\rarrow \gR$, which should be a map of right $A$\module s satisfying
the analogous equations.

\subsubsection{}  \label{contramodule-std-example}
 The standard example of a $\C$\contramodule: for any right
$\C$\comodule{} $\N$ endowed with a left action of a $k$\+algebra $B$
by $\C$\comodule{} endomorphisms and any left $B$\module{} $V$,
the left $A$\module{} $\Hom_B(\N,V)$ has a natural left
$\C$\contramodule{} structure.
 The left $\C$\contramodule{} $\Hom_A(\C,V)$ is called
the $\C$\contramodule{} \emph{induced} from a left $A$\module{} $V$.
 According to Lemma~\ref{induced-coinduced}, the $k$\module{} of
contramodule homomorphisms from the induced $\C$\contramodule{} to
an arbitrary $\C$\contramodule{} is described by the formula
$\Hom^\C(\Hom_A(\C,V),\.\P)\simeq\Hom_A(V,\P)$.

 We will denote the category of left $\C$\contramodule s by $\C\contra$
and the category of right $\C$\contramodule s by $\contraR\C$.
 The category of left $\C$\contramodule s is abelian whenever $\C$ is
a projective left $A$\module.
 Moreover, the left $A$\module{} $\C$ is projective if and only if
the category $\C\contra$ is abelian and the forgetful functor
$\C\contra\rarrow A\modl$ is exact.
 This can be proven by the same adjoint functor argument as
the analogous result for $\C$\comodule s.

 For any coring~$\C$, there are two natural exact categories of left
contramodules: the exact category of $A$\+injective $\C$\contramodule s
and the exact category of arbitrary $\C$\contramodule s with $A$\+split
exact triples.
 Besides, any morphism of $\C$\contramodule s has a kernel and
the forgetful functor $\C\contra\rarrow A\modl$ preserves kernels.
 When a morphism of $\C$\contramodule s has the property that its
cokernel in the category of $A$\module s is preserved by the functors
of homomorphisms from $\C$ and $\C\ot_A\C$ over~$A$, this cokernel has
a natural $\C$\contramodule{} structure, which makes it the cokernel
of that morphism in the category of $\C$\contramodule s.

 Infinite products always exist in the category of $\C$\contramodule s
and the forgetful functor $\C\contra\rarrow A\modl$ preserves them.
 The induction functor $A\modl\rarrow\C\contra$ preserves both infinite
direct sums and infinite products.
 To construct direct sums of $\C$\contramodule s, one can present them
as cokernels of morphisms of induced contramodules, and all cokernels
exist in the category of $\C$\contramodule s~\cite{Bar}, so
the category of $\C$\contramodule s has infinite direct sums.

\begin{qst}
 If $\C$ is a flat right $A$\module, then subcomodules of finite
direct sums of copies of~$\C$ constitute a set of generators of
the category of left $\C$\comodule s~\cite{BW}.
 Does the category of $\C$\contramodule s have a set of
cogenerators?
\end{qst}

\subsubsection{}   \label{proj-inj-co-contra-module}
 Assume that the coring $\C$ is a projective left and a flat right
$A$\module{} and the ring $A$ has a finite left homological dimension
(homological dimension of the category of left $A$\module s).

\begin{lem}
 \textup{(a)} There exists a (not always additive) functor assigning
to any left\/ $\C$\comodule{} a surjective map onto it from
an $A$\projective\/ $\C$\comodule. 
 Moreover, the kernel of this map is an iterated extension of
coinduced\/ $\C$\comodule s. \par
 \textup{(b)} There exists a (not always additive) functor assigning
to any left\/ $\C$\contramodule{} an injective map from it into
an $A$\injective\/ $\C$\contramodule.
 Moreover, the cokernel of this map is an iterated extension of
induced\/ $\C$\contramodule s.
\end{lem}

\begin{proof}
 The proof of part~(a) is completely analogous to the proof of
Lemma~\ref{flat-comodule-surjection} and part~(b) is proven in
the following way.
 Let $P\rarrow G(P)$ be an injective map from an $A$\module{} $P$
into an injective $A$\module{} $G(P)$ functorially depending on~$P$.
 For example, one can take $G(P)$ to be the direct product of copies
of the $A$\module{} $\Hom_A(A,k\dual)$ numbered by all $k$\module{}
homomorphisms $P\rarrow k\dual$.
 Let $\P$ be a left $\C$\contramodule.
 Consider the contraaction map $\Hom_A(\C,\P)\rarrow\P$; it is
a surjective morphism of $\C$\contramodule s; let $\gK(\P)$ denote
its kernel.
 Let $\Q(\P)$ be the cokernel of the composition $\gK(\P)\rarrow
\Hom_A(\C,\P)\rarrow\Hom_A(\C,G(\P))$.
 Then the composition of maps $\Hom_A(\C,\P)\rarrow\Hom_A(\C,G(\P))
\rarrow\Q(\P)$ factorizes through the surjection $\Hom_A(\C,\P)
\rarrow\P$, so there is a natural injective morphism of
$\C$\contramodule s $\P\rarrow\Q(\P)$.
 Let us show that the injective dimension $\di_A\Q(\P)$ of
the $A$\module{} $\Q(\P)$ is smaller than that of~$\P$.
 Indeed, the $A$\module{} $\Hom_A(\C,G(\P))$ is injective, hence
$\di_A\Q(\P)=\di_A\gK(\P)-1 \le \di_A\Hom_A(\C,\P)-1 \le \di_A(\P)-1$,
because the $A$\module{} $\gK(\P)$ is a direct summand of
the $A$\module{} $\Hom_A(\C,\P)$ and an injective resolution of
the $A$\module{} $\Hom_A(\C,\P)$ can be constructed by applying
the functor $\Hom_A(\C,{-})$ to an injective resolution of~$\P$.
 Notice that the cokernel of the map $\P\rarrow\Q(\P)$ is an induced
$\C$\contramodule{} $\Hom_A(\C,G(\P)/\P)$.
 It remains to iterate the functor $\P\mpsto\Q(\P)$ sufficiently
many times.
\end{proof}

\subsection{Cohomomorphisms}

\subsubsection{}
 The $k$\module{} of \emph{cohomomorphisms} $\Cohom_\C(\M,\P)$ from
a left $\C$\comodule{} $\M$ to a left $\C$\contramodule{} $\P$ is
defined as the cokernel of the pair of maps $\Hom_A(\C\ot_A\M\;\P) =
\Hom_A(\M,\Hom_A(\C,\P)) \birarrow \Hom_A(\M,\P)$ one of which is
induced by the $\C$\+coaction in~$\M$ and the other by
the $\C$\+contraaction in~$\P$.
 The functor of cohomomorphisms is neither left nor right exact in
general; it is right exact if the ring $A$ is semisimple.
 For any left $A$\module{} $U$ and any left $\C$\contramodule{} $\P$
there is a natural isomorphism $\Cohom_\C(\C\ot_A U\;\P)\simeq
\Hom_A(U,\P)$, and for any left $\C$\comodule{} $\M$ and any left
$A$\module{} $V$ there is a natural isomorphism
$\Cohom_\C(\M,\Hom_A(\C,V))\simeq \Hom_\C(\M,V)$.
 These assertions follow from Lemma~\ref{induced-tensor-cotensor}.
 Explicitly, the first isomorphism can be obtained by applying
the functor $\Hom_A(U,{-})$ to the split exact sequence of $A$\module s
$\Hom_A(\C\ot_A\C\;\P)\rarrow\Hom_A(\C,\P)\rarrow\P$ and the second one
can be obtained by applying the functor $\Hom_A({-},V)$ to the split
exact sequence of $A$\module s $\M\rarrow\C\ot_A\M\rarrow
\C\ot_A\C\ot_A\M$.

\subsubsection{}   \label{absolute-relative-coproj-coinj}
 Assuming that $\C$ is a projective left $A$\module, a left comodule
$\M$ over $\C$ is called \emph{coprojective} if the functor of
cohomomorphisms from~$\M$ is exact on the category of left
$\C$\contramodule s.
 It is easy to see that any coprojective $\C$\comodule{} is
a projective $A$\module.
 The $\C$\comodule{} coinduced from a projective $A$\module{}
is coprojective.
 Assuming that $\C$ is a flat right $A$\module, a left contramodule $\P$
over $\C$ is called \emph{coinjective} if the functor of cohomomorphisms
into~$\P$ is exact on the category of left $\C$\comodule s.
 Any coinjective $\C$\contramodule{} is an injective $A$\module.
 The $\C$\contramodule{} induced from an injective $A$\module{}
is coinjective.

 A left comodule $\M$ over $\C$ is called \emph{coprojective relative
to~$A$} ($\C/A$\+coprojective) if the functor of cohomomorphisms
from $\M$ maps exact triples of $A$\injective{} $\C$\contramodule s
to exact triples.
 A left contramodule $\P$ over $\C$ is called \emph{coinjective
relative to~$A$} ($\C/A$\+coinjective) if the functor of
cohomomorphisms into~$\P$ maps exact triples of $A$\projective{}
$\C$\comodule s to exact triples.
 Any coinduced $\C$\comodule{} is $\C/A$\coprojective{} and
any induced $\C$\contramodule{} is $\C/A$\coinjective.

 For any right $\C$\comodule{} $\N$ and any left $\C$\comodule{} $\M$
there is a natural isomorphism $\Hom_k(\N\oc_\C\M\;k\dual) \simeq
\Cohom_\C(\M,\Hom_k(\N,k\dual))$.
 Therefore, any coprojective $\C$\comodule{} $\M$ is coflat and any
$\C/A$\coprojective{} $\C$\comodule{} $\M$ is $\C/A$\+coflat.
 Besides, a right $\C$\comodule{} $\N$ is coflat if and only if
the left $\C$\contramodule{} $\Hom_k(\N,k\dual)$ is coinjective;
if a right $\C$\comodule{} $\N$ is $\C/A$\+coflat, then the left
$\C$\contramodule{} $\Hom_k(\N,k\dual)$ is $\C/A$\coinjective{}
(and the converse can be deduced from
Lemma~\ref{proj-inj-co-contra-module}(a) and the proof of Lemma
below in the assumptions of~\ref{proj-inj-co-contra-module}).

 It appears that the notion of a relatively coprojective left
$\C$\comodule{} is useful when $\C$ is a flat right $A$\module,
and the notion of a relatively coinjective left $\C$\contramodule{}
is useful when $\C$ is a projective left $A$\module.

\begin{lem}
 \textup{(a)} Assume that\/ $\C$ is a flat right $A$\module.
 Then the class of\/ $\C/A$\coprojective{} left\/ $\C$\comodule s is
closed under extensions and cokernels of injective morphisms.
 The functor of cohomomorphisms into an $A$\injective{} left\/
$\C$\contramodule{} maps exact triples of\/ $\C/A$\coprojective{}
left\/ $\C$\comodule s to exact triples. \par
 \textup{(b)} Assume that\/ $\C$ is a projective left $A$\module.
 Then the class of\/ $\C/A$\coinjective{} left\/ $\C$\contramodule s
is closed under extensions and kernels of surjective morphisms.
 The functor of cohomomorphisms from an $A$\projective{}
left\/ $\C$\comodule{} maps exact triples of\/ $\C/A$\coinjective{}
left\/ $\C$\contramodule s to exact triples.
\end{lem}

\begin{proof}
 Part~(a): these results follow from the standard properties of
the left derived functor of the right exact functor of cohomomorphisms
on the Carthesian product of the abelian category of left
$\C$\comodule s and the exact category of $A$\injective{} left
$\C$\contramodule s.
 One can define the $k$\module s $\Coext_\C^i(\M,\P)$, \
$i=0$,~$-1$,~\dots\ as the homology of the bar complex $\dsb\rarrow
\Hom_A(\C\ot_A\C\ot_A\M\;\P)\rarrow\Hom_A(\C\ot_A\M\;\P)\rarrow
\Hom_A(\M,\P)$ for any left $\C$\comodule{} $\M$ and any
$A$\injective{} left $\C$\contramodule{} $\P$.
 Then $\Coext_\C^0(\M,\P)\simeq\Cohom_\C(\M,\P)$ and there are long
exact sequences of $\Coext_\C^*$ associated with exact triples of
comodules and contramodules.
 Now a left $\C$\comodule{} $\M$ is $\C/A$\coprojective{} if and only
if $\Coext_\C^i(\M,\P)=0$ for any $A$\injective{} left
$\C$\contramodule{} $\P$ and all $i<0$.
 Indeed, the ``if'' assertion follows from the homological exact
sequence, and ``only if'' holds since the bar complex is isomorphic
to the complex of cohomomorphisms from the $\C$\comodule{} $\M$ into
the bar resolution $\dsb\rarrow\Hom_A(\C,\Hom_A(\C,\P))\rarrow
\Hom_A(\C,\P)$ of the $\C$\contramodule{} $\P$, which is
a complex of $A$\injective{} $\C$\contramodule s, exact except
in degree~$0$ and split over~$A$.
 The proof of part~(b) is completely analogous; it uses the left
derived functor of the functor of cohomomorphisms on the Carthesian
product of the exact category of $A$\projective{} left $\C$\comodule s
and the abelian category of left $\C$\contramodule s.
\end{proof}

\begin{rmk}
 It follows from Lemma~\ref{cotensor-contratensor-assoc} that any
extension of an $A$\projective{} $\C$\+co\-mod\-ule by a coprojective
$\C$\comodule{} splits, and any extension of a coinjective
$\C$\contramodule{} by an $A$\injective{} $\C$\contramodule{} splits.
 The analogues of the results of Remark~\ref{absolute-relative-coflat}
also hold for (relatively) coprojective comodules and coinjective
contramodules in the assumptions of~\ref{proj-inj-co-contra-module};
see the proof of Lemma~\ref{rel-inj-proj-co-contra-mod}.2 for details.
\end{rmk}

\begin{qst}
 Are all relatively coflat $\C$\comodule s relatively coprojective?
 Are all $A$\projective{} coflat $\C$\comodule s coprojective?
\end{qst}

\subsubsection{}   \label{co-tensor-co-hom-assoc}
 Let $\C$ be an arbitrary coring.
 Let us call a left $\C$\comodule{} $\M$ \emph{quasicoprojective}
if the functor of cohomomorphisms from $\M$ is left exact on
the category of left $\C$\contramodule s, i.~e., this functor
preserves kernels.
 Any coinduced $\C$\comodule{} is quasicoprojective.
 Any quasicoprojective comodule is quasicoflat.
 Let us call a left $\C$\contramodule{} $\P$ \emph{quasicoinjective}
if the functor of cohomomorphisms into $\P$ is left exact on
the category of left $\C$\comodule s, i.~e., this functor maps
cokernels to kernels.
 Any induced $\C$\contramodule{} is quasicoinjective.
 (Cf.\ Lemma~\ref{cotensor-contratensor-assoc}.)

\begin{prop1}
 Let $\M$ be a left\/ $\C$\comodule, $\K$ be a right\/ $\C$\comodule{}
endowed with a left action of a $k$\+algebra $B$ by comodule
endomorphisms, and $P$ be a left $B$\module.
 Then there is a natural $k$\module{} map\/ $\Cohom_\C(\M,\Hom_B(\K,P))
\rarrow \Hom_B(\K\oc_\C\M\;P)$, which is an isomorphism, at least, in
the following cases:
\begin{enumerate}
 \item $P$ is an injective left $B$\module;
 \item $\M$ is a quasicoprojective left\/ $\C$\comodule;
 \item $\C$ is a projective left $A$\module, $\M$ is a projective
       left $A$\module, $\K$ is a\/ $\C/A$\+coflat right\/
       $\C$\comodule, $\K$ is a projective left $B$\module, and
       the ring $B$ has a finite left homological dimension;
 \item $\K$ as a right\/ $\C$\comodule{} with a left $B$\module{}
       structure is coinduced from a $B$\+$A$\bimodule.
\end{enumerate}
 Besides, in the case\/~\textup{(c)} the left $B$\module\/
$\K\oc_\C\M$ is projective.
\end{prop1}

\begin{proof}
 The map $\Hom_B(\K\ot_A\M\;P)\rarrow\Hom_B(\K\oc_\C\M\;P)$ annihilates
the difference of two maps $\Hom_B(\K\ot_A\C\ot_A\M\;P)\birarrow
\Hom_B(\K\ot_A\M\;P)$ and this pair of maps can be identified with
the pair of maps $\Hom_A(\C\ot_A\M\;\Hom_B(\K,P))\birarrow
\Hom_A(\M,\Hom_B(\K,P))$ whose cokernel is, by the definition,
the cohomomorphism module $\Cohom_\C(\M,\Hom_B(\K,P))$.
 Hence there is a natural map $\Cohom_\C(\M,\allowbreak\Hom_B(\K,P))
\rarrow\Hom_B(\K\oc_\C\M\;P)$.
 The case~(a) is obvious.
 In the case~(b), it suffices to present $P$ as the kernel of a map
of injective $B$\module s.
 The rest of the proof is completely analogous to the proof of
Proposition~\ref{tensor-cotensor-assoc} (with flat modules replaced
by projective ones and the left and right sides switched).
\end{proof}

\begin{prop2}
 Let $\P$ be a left\/ $\C$\contramodule, $\K$ be a left\/
$\C$\comodule{} endowed with a right action of a $k$\+algebra $B$ by
comodule endomorphisms, and $M$ be a left $B$\module.
 Then there is a natural $k$\module{} map\/ $\Cohom_\C(\K\ot_B M\;\P)
\rarrow \Hom_B(M,\Cohom_\C(\K,\P))$, which is an isomorphism, at least,
in the following cases:
\begin{enumerate}
 \item $M$ is a projective left $B$\module;
 \item $\P$ is a quasicoinjective left\/ $\C$\contramodule;
 \item $\C$ is a flat right $A$\module, $\P$ is an injective
       left $A$\module, $\K$ is a $\C/A$\coprojective{} left\/
       $\C$\comodule, $\K$ is a flat right $B$\module, and
       the ring $B$ has a finite left homological dimension;
 \item $\K$ as a left\/ $\C$\comodule{} with a right $B$\module{}
       structure is coinduced from an $A$\+$B$\bimodule.
\end{enumerate}
 Besides, in the case\/~\textup{(c)} the left $B$\module\/
$\Cohom_\C(\K,\P)$ is injective.
\end{prop2}

\begin{proof}
 The map $\Hom_B(M,\Hom_A(\K,\P))\rarrow\Hom_B(M,\Cohom_\C(\K,\P))$
annihilates the difference of two maps $\Hom_B(M\;\Hom_A(\C\ot_A\K\;\P))
\birarrow\Hom_B(M,\Hom_A(\K,\P))$ and this pair of maps can be
identified with the pair of maps $\Hom_A(\C\ot_A\K\ot_B M\;\allowbreak
\P)\birarrow \Hom_A(\K\ot_B M\;\P)$ whose cokernel is, by
the definition, the cohomomorphism module $\Cohom_\C(\K\ot_B M\;\P)$.
 Hence there is a natural map $\Cohom_\C(\K\ot_B M\;\P)
\rarrow \Hom_B(M,\Cohom_\C(\K,\P))$.
 The case~(a) is obvious.
 In the case~(b), it suffices to present $M$ as the cokernel of
a map of projective $B$\module s.
 To prove (c) and~(d), consider the bar complex
\begin{multline}  \label{hom-cohom}
 \dsb\lrarrow\Hom_A(\C\ot_A\C\ot_A\K\;\P)
 \lrarrow\Hom_A(\C\ot_A\K\;\P) \\
 \lrarrow\Hom_A(\K,\P)\lrarrow\Cohom_\C(\K,\P).
\end{multline}
 In the case~(c) this complex is exact, since it is the complex
of cohomomorphisms from a $\C/A$\coprojective{} $\C$\comodule{} $\K$
into an $A$\+split exact complex of $A$\injective{} $\C$\contramodule s
$\dsb\rarrow\Hom_A(\C\ot_A\C\;\P)\rarrow\Hom_A(\C,\P)\rarrow\P$.
 Since all the terms of the complex~\eqref{hom-cohom}, except possibly
the rightmost one, are injective left $B$\module s and the left
homological dimension of the ring~$B$ is finite, the rightmost term
$\Cohom_\C(\K,\P)$ is also an injective $B$\module, the complex of
left $B$\module s~\eqref{hom-cohom} is contractible, and the complex
of $B$\module{} homomorphisms from the left $B$\module{} $M$
into~\eqref{hom-cohom} is exact.
 In the case~(d), the complex~\eqref{hom-cohom} is also a split exact
complex of left $B$\module s.
\end{proof}

\subsubsection{}   \label{bicomodule-cohom}
 Let $\C$ be a coring over a $k$\+algebra $A$ and $\D$ be a coring
over a $k$\+algebra~$B$.
 Assume that $\D$ is a projective left $B$\module.
 Let $\K$ be a $\C$\+$\D$\bicomodule{} and $\P$ be a left
$\C$\contramodule.
 Then the module of cohomomorphisms $\Cohom_\C(\K,\P)$ is endowed\- 
with a left $\D$\contramodule{} structure as the cokernel of a pair of
contramodule morphisms $\Hom_A(\C\ot_A\K\;\P)\birarrow\Hom_A(\K,\P)$.
%(where $\Hom_A(\K,\P)$ and $\Hom_A(\C\ot_A\K\;\P)=
%\Hom_A(\K,\Hom_A(\C,\P))$ are left $\D$\contramodule s
%by~\ref{contramodule-std-example}).

 More generally, let $\C$ and $\D$ be arbitrary corings.
 Assume that the functor of homomorphisms from $\D$ over $B$ preserves
the cokernel of the pair of maps $\Hom_A(\C\ot_A\K\;\P)\birarrow
\Hom_A(\K,\P)$, that is the natural map $\Cohom_\C(\K\ot_B\D\;\P)
\rarrow\Hom_B(\D,\Cohom_\C(\K,\P))$ is an isomorphism.
 Then one can define a left contraaction map
$\Hom_B(\D,\Cohom_\C(\K,\P))\rarrow\Cohom_\C(\K,\P)$ taking
the cohomomorphisms over~$\C$ from the right $\D$\+coaction map
$\K\rarrow\K\ot_B\D$ into the contramodule~$\P$.
 This contraaction is counital and
contraassociative, at least, if the natural map
$\Cohom_\C(\K\ot_B\D\ot_B\D\;\P)\allowbreak\rarrow
\Hom_B(\D\ot_B\D\;\Cohom_\C(\K,\P))$ is also an isomorphism.

 In particular, if one of the conditions of
Proposition~\ref{co-tensor-co-hom-assoc}.2 is satisfied (for $M=\D$),
then the left $B$\module{} $\Cohom_\C(\K,\P)$ has a natural
$\D$\contramodule{} structure.

\subsubsection{}   \label{cohom-associative}
 Let $\C$ be a coring over a $k$\+algebra $A$ and $\D$ be a coring
over a $k$\+algebra~$B$.

\begin{prop}
 Let\/ $\M$ be a left\/ $\D$\comodule, $\K$ be a $\C$\+$\D$\bicomodule,
and\/ $\P$ be a left $\C$\contramodule.
 Then the iterated cohomomorphism modules\/ $\Cohom_\C(\K\oc_\D\M\;\P)$
and\/ $\Cohom_\D(\M,\Cohom_\C(\K,\P))$ are naturally isomorphic,
at least, in the following cases:
\begin{enumerate}
 \item $\D$ is a projective left $B$\module, $\M$ is a projective
       left $B$\module, $\C$ is a flat right, and\/ $\P$ is
       an injective left $A$\module;
 \item $\D$ is a projective left $B$\module{} and\/ $\M$ is
       a coprojective left\/ $\D$\comodule;
 \item $\C$ is a flat right $A$\module{} and\/ $\P$ is
       a coinjective left\/ $\C$\contramodule;
 \item $\D$ is a projective left $B$\module, $\M$ is a projective
       left $B$\module, $\K$ is a $\D/B$\+coflat right\/
       $\D$\comodule, $\K$ is a projective left $A$\module, and
       the ring $A$ has a finite left homological dimension;
 \item $\C$ is a flat right $A$\module, $\P$ is an injective left
       $A$\module, $\K$ is a $\C/A$\coprojective{} left\/
       $\C$\comodule, $\K$ is a flat right $B$\module, and
       the ring $B$ has a finite left homological dimension;
 \item $\D$ is a projective left $B$\module, $\M$ is a projective
       left $B$\module, and\/ $\K$ as a right\/ $\D$\comodule{} with
       a left $A$\module{} structure is coinduced from
       an $A$\+$B$\bimodule;
 \item $\C$ is a flat right $A$\module, $\P$ is an injective left
       $A$\module, and\/ $\K$ as a left\/ $\C$\comodule{} with
       a right $B$\module{} structure is coinduced from
       an $A$\+$B$\bimodule;
 \item $\M$ is a quasicoprojective left\/ $\D$\comodule{} and\/ $\K$ as
       a left\/ $\C$\comodule{} with a right $B$\module{} structure
       is coinduced from an $A$\+$B$\bimodule;
 \item $\P$ is a quasicoinjective left $\C$\contramodule{} and\/ $\K$
       as a right\/ $\D$\comodule{} with a left $A$\module{} structure
       is coinduced from an $A$\+$B$\bimodule;
 \item $\K$ as a left\/ $\C$\comodule{} with a right $B$\module{}
       structure is coinduced from an $A$\+$B$\bimodule{} and\/ $\K$
       as a right\/ $\D$\comodule{} with a left $A$\module{} structure
       is coinduced from an $A$\+$B$\bimodule.
\end{enumerate}
 More precisely, in all cases in this list the natural maps from
the $k$\+module\/ $\Hom_A(\K\ot_B\M\;\P)=\Hom_B(\M,\Hom_A(\K,\P))$
into both iterated cohomomorphism modules under consideration
are surjective, their kernels coincide and are equal to the sum of
the kernels of two maps from this module onto its quotient modules\/
$\Cohom_\C(\K\ot_B\M\;\P)$ and\/ $\Cohom_\D(\M,\Hom_A(\K,\P))$. 
\end{prop}

\begin{proof}
 One can easily see that whenever both maps
$\Cohom_\D(\M,\Hom_A(\K,\P)) \rarrow \Hom_A(\K\oc_\D\M\;\P)$
and $\Cohom_\D(\M,\.\Hom_A(\K,\.\Hom_A(\C,\P)))\rarrow
\Hom_A(\K\oc_\D\M\;\Hom_A(\C,\P))$ are isomorphisms, the natural map
$\Hom_A(\K\ot_B\M\;\P)\rarrow\Cohom_\C(\K\oc_\D\M\;\P)$ is surjective
and its kernel coincides with the desired sum of two kernels of maps
from $\Hom_A(\K\ot_B\M\;\P)$ onto its quotient modules.
 Analogously, whenever both maps $\Cohom_\C(\K\ot_B\M\;\P)\rarrow
\Hom_B(\M,\Cohom_\C(\K,\P))$ and $\Cohom_\C(\K\ot_B\D\ot_B\M\;\P)
\rarrow\Hom_B(\D\ot_B\M,\Cohom_\C(\K,\P))$ are isomorphisms,
the natural map $\Hom_B(\M,\Hom_A(\K,\P))\rarrow
\Cohom_\D(\M,\Cohom_\C(\K,\P))$ is surjective and it kernel
coincides with the desired sum of two kernels in
$\Hom_B(\M,\Hom_A(\K,\P))$.
 Thus it remains to apply Propositions~\ref{co-tensor-co-hom-assoc}.1
and~\ref{co-tensor-co-hom-assoc}.2.
\end{proof}

 Commutativity of pentagonal diagrams of associativity isomorphisms
between iterated cohomomorphism modules can be established in the way
analogous to the case of iterated cotensor products.
 Namely, each of the five iterated cohomomorphism modules
$\Cohom_\C((\K\oc_\E\L)\oc_\D\M\;\P)$, \
$\Cohom_\C(\K\oc_\E(\L\oc_\D\M)\;\P)$, \
$\Cohom_\E(\L\oc_\D\M\;\Cohom_\C(\K,\P))$, \ 
$\Cohom_\D(\M,\.\Cohom_\E(\L,\.\Cohom_\C(\K,\P)))$, and
$\Cohom_\D(\M\;\Cohom_\C(\K\oc_\E\L\;\P))$
is endowed with a natural map into it from the homomorphism module
$\Hom_A(\K\ot_F\L\ot_B\M\;\P)$, and since the associativity
isomorphisms are, presumably, compatible with these maps, it suffices
to check that at least one of these five maps is surjective in order
to show that the pentagonal diagram commutes.
 In particular, if the above Proposition together with
Proposition~\ref{cotensor-associative} provide all the five
isomorphisms constituting the pentagonal diagram and either $\M$
is a projective left $B$\module, or $\P$ is an injective left
$A$\module, or both $\K$ and $\L$ as left (right) comodules
with right (left) module structures are coinduced from bimodules,
then the pentagonal diagram is commutative.

 We will say that multiple cohomomorphisms between several bicomodules
and a contramodule $\Cohom_\C(\K\oc_\E\dsb\oc_\D\M\;\P)$
are associative if the multiple cotensor product $\K\oc_\E\dsb\oc_\D\M$
is associative and for any possible way of representing this multiple
cohomomorphism module in terms of iterated cotensor product and
cohomomorphism operations all the intermediate cohomomorphism modules
can be endowed with contramodule structures via the construction
of~\ref{bicomodule-cohom}, all possible associativity isomorphisms
between iterated cohomomorphism modules exist in the sense of the last
assertion of Proposition and preserve contramodule structures, and all
the pentagonal diagrams commute.
 Associativity isomorphisms and contramodule structures on associative
multiple cohomomorphisms are preserved by the morphisms between them
induced by any bicomodule and contramodule morphisms of the factors.

\subsection{Semicontramodules} 

\subsubsection{}   \label{semicontramodules-definition}
 Depending on the (co)flatness, (co)projectivity, and/or
(co)injectivity conditions imposed, there are several ways to make
the category opposite to a category of left $\C$\contramodule s
into a right module category over a tensor category of
$\C$\+$\C$\bicomodule s with respect to the functor $\Cohom_\C$.
 Moreover, a category of left $\C$\comodule s typically can be made
into a left module category over the same tensor category, so that
the functor $\Cohom_\C$ would provide also a pairing between these
left and right module categories taking values in the category
$k\modl^\op$.

 It follows from Proposition~\ref{cohom-associative}(b) that whenever
$\C$ is a projective left $A$\module, the category opposite to
the category of left $\C$\contramodule s is a right module category
over the tensor category of $\C$\+$\C$\bicomodule s that are
coprojective left $\C$\comodule s; the category of coprojective
left $\C$\comodule s is a left module category over this tensor
category.
 If follows from Proposition~\ref{cohom-associative}(c) that whenever
$\C$ is a flat right $A$\module, the category opposite to the category
of coinjective left $\C$\contramodule s is a right module category
over the tensor category of $\C$\+$\C$\bicomodule s that are
coflat right $\C$\comodule s; the category of left $\C$\comodule s
is a left module category over this tensor category.
 It follows from Proposition~\ref{cohom-associative}(d) that whenever
$\C$ is a projective left $A$\module{} and the ring $A$ has a finite
left homological dimension, the category opposite to the category
of left $\C$\contramodule s is a right module category over the tensor
category of $\C$\+$\C$\bicomodule s that are projective left
$A$\module s and $\C/A$\+coflat right $\C$\comodule s; the category
of $A$\+projective left $\C$\comodule s is a left module category
over this tensor category.
 It follows from Proposition~\ref{cohom-associative}(e) that
whenever $\C$ is a flat right $A$\module{} and the ring $A$ has
a finite left homological dimension, the category opposite to
the category of $A$\injective{} left $\C$\contramodule s is a right
module category over the tensor category of $\C$\+$\C$\bicomodule s
that are flat right $A$\module s and $\C/A$\coprojective{} left
$\C$\comodule s; the category of left $\C$\comodule s is a left module
category over this tensor category.
 Finally, it follows from Proposition~\ref{cohom-associative}(a) that
whenever the ring $A$ is semisimple, the category opposite to
the category of left $\C$\contramodule s is a right module category
over the tensor category of $\C$\+$\C$\bicomodule s; the category of
left $\C$\comodule s is a left module category over this tensor
category.
 In each case, there is a pairing between these left and right module
categories compatible with their module category structures and taking
values in the category opposite to the category of $k$\module s.

 A \emph{left semicontramodule} over a semialgebra $\S$ is
an object of the category opposite to the category of module objects
in one of the right module categories of the above kind (opposite to
a category of left $\C$\contramodule s) over the ring object~$\S$
in the corresponding tensor category of $\C$\+$\C$\bicomodule s.
 In other words, a left $\S$\semicontramodule{} $\bP$ is a left
$\C$\contramodule{} endowed with a left $\C$\contramodule{} morphism
of \emph{left semicontraaction} $\bP\rarrow\Cohom_\C(\S,\bP)$
satisfying the associativity and unity equations.
 Namely, two compositions $\bP\rarrow\Cohom_\C(\S,\bP)\birarrow
\Cohom_\C(\S\oc_\C\S\;\bP)$ of the semicontraaction morphism
$\bP\rarrow\Cohom_\C(\S,\bP)$ with the morphisms
$\Cohom_\C(\S,\bP)\birarrow\Cohom_\C(\S\oc_\C\S\;\bP)=
\Cohom_\C(\S,\Cohom_\C(\S,\bP))$ induced by the semimultiplication
morphism of~$\S$ and the semicontraaction morphism should coincide
with each other and the composition $\bP\rarrow\Cohom_\C(\S,\bP)
\rarrow\bP$ of the semicontraaction morphism with the morphism
induced by the semiunit morphism of $\S$ should coincide with
the identity morphism of~$\bP$.
 For this definition to make sense, (co)flatness, (co)projectivity,
and/or (co)injectivity conditions imposed on $\S$ and/or $\bP$
must guarantee associativity of multiple cohomomorphism modules of
the form $\Cohom_\C(\S\oc_\C\dsb\oc_\C\S\;\bP)$.
 \emph{Right semicontramodules} over~$\S$ are defined in
the analogous way.

 If $\Q$ is a left $\C$\contramodule{} for which multiple
cohomomorphisms $\Cohom_\C(\S\oc_\C\dsb\oc_\C\S\;\Q)$ are associative,
then there is a natural left $\S$\semicontramodule{} structure on
the cohomomorphism module $\Cohom_\C(\S,\Q)$.
 The semicontramodule $\Cohom_\C(\S,\Q)$ is called
the $\S$\semicontramodule{} \emph{coinduced} from
a $\C$\contramodule{}~$\Q$.
 According to Lemma~\ref{induced-coinduced}, the $k$\module{}
of semicontramodule homomorphisms from an arbitrary
$\S$\semicontramodule{} into the coinduced $\S$\semicontramodule{}
is described by the formula $\Hom^\S(\bP,\Cohom_\C(\S,\Q))\simeq
\Hom^\C(\bP,\Q)$.

 We will denote the category of left $\S$\semicontramodule s by
$\S\sicntr$ and the category of right $\S$\semicontramodule s by
$\sicntrR\S$.
 This notation presumes that one can speak of (left or right)
$\S$\semicontramodule s with no (co)injectivity conditions imposed
on them.
 If $\C$ is a projective left $A$\module{} and $\S$ is a coprojective
left $\C$\comodule, then the category of left semicontramodules
over~$\S$ is abelian and the forgetful functor
$\S\sicntr\rarrow\C\contra$ is exact.

 If $\C$ is a projective left $A$\module{} and either $\S$ is
a coprojective left $\C$\comodule, or $\S$ is a projective left
$A$\module{} and a $\C/A$\+coflat right $\C$\comodule{} and $A$ has
a finite left homological dimension, or $A$ is semisimple, then both
infinite direct sums and infinite products exist in the category of
left $\S$\semicontramodule s and both are preserved by the forgetful
functor $\S\sicntr\rarrow\C\contra$, even though only infinite products
are preserved by the full forgetful functor $\S\sicntr\rarrow A\modl$.

 If $\C$ is a flat right $A$\module, $\S$ is a flat right $A$\module{}
and a $\C/A$\coprojective{} left $\C$\comodule, and $A$ has a finite
left homological dimension, then the category of $A$\injective{}
left $\S$\semicontramodule s is exact.
 If $\C$ is a projective left $A$\module, $\S$ is a projective left
$A$\module{} and a $\C/A$\+coflat right $\C$\comodule, and $A$ has
a finite left homological dimension, then the category of
$\C/A$\coinjective{} left $\S$\semicontramodule s is exact.
 If $\C$ is a flat right $A$\module{} and $\S$ is a coflat right
$\C$\comodule, then the category of $\C$\coinjective{} left
$\S$\semicontramodule s is exact.
 If $A$ is semisimple, the category of $\C$\coinjective{}
$\S$\semimodule s is exact. 
 Infinite products exist in all of these exact categories, and
the forgetful functors preserve them.

\begin{qst}
 When $\C$ is a flat right $A$\module{} and $\S$ is a coflat right
$\C$\comodule, a right adjoint functor to the forgetful functor
$\S\simodl\rarrow\C\comodl$ exists according to the abstract adjoint
functor existence theorem~\cite{McL}.
 Indeed, the forgetful functor preserves colimits and
the category of left $\S$\semimodule s has a set of generators
(since the category of left $\C$\comodule s does; see
Question~\ref{contramodule-std-example}).
 Does a left adjoint functor to the forgetful functor $\S\sicntr
\rarrow\C\contra$ exist?
 Can one describe these functors more explicitly?
\end{qst}

\subsubsection{}   \label{proj-inj-semi-mod-contra}
 Assume that the coring $\C$ is a projective left and a flat right
$A$\module{} and the ring $A$ has a finite left homological dimension.

\begin{lem}
 \textup{(a)} If the semialgebra\/ $\S$ is a coflat right\/
$\C$\comodule{} and a projective left $A$\module, then there exists
a (not always additive) functor assigning to any left\/
$\S$\semimodule{} a surjective map onto it from an $A$\projective\/
$\S$\semimodule. \par
 \textup{(b)} If the semialgebra\/ $\S$ is a coprojective left\/
$\C$\comodule{} and a flat right $A$\module, then there exists
a (not always additive) functor assigning to any left\/
$\S$\semicontramodule{} an injective map from it into
an $A$\injective\/ $\S$\semicontramodule.
\end{lem}

\begin{proof}
 The proof of part~(a) is completely analogous to the proof
of Lemma~\ref{flat-semimodule-surjection} (with the last assertion of
Proposition~\ref{co-tensor-co-hom-assoc}.1 used as needed);
and part~(b) is proven in the following way.
 Let $\P\rarrow\gI(\P)$ denote the functorial injective morphism
from a $\C$\contramodule{} $\P$ into an $A$\injective{}
$\C$\contramodule{} $\gI(\P)$ constructed in
Lemma~\ref{proj-inj-co-contra-module}.
 Then for any $\S$\semicontramodule{} $\bP$ the composition of maps
$\bP\rarrow\Cohom_\C(\S,\bP)\rarrow\Cohom_\C(\S,\gI(\bP))$
provides the desired injective morphism of $\S$\semicontramodule s.
 According to the last assertion of
Proposition~\ref{co-tensor-co-hom-assoc}.2, the $A$\module{}
$\bgI(\bP)=\Cohom_\C(\S,\gI(\bP))$ is injective.
\end{proof}

\begin{rmk}
 The analogues of the result of Remark~\ref{flat-semimodule-surjection}
hold for $\C/A$\coprojective/semi\-pro\-jec\-tive{} $\S$\semimodule s
and $\C/A$\coinjective/semi\-in\-jec\-tive{} $\S$\semicontramodule s;
see the proof of Lemma~\ref{semi-model-struct}.1 for details.
\end{rmk}

\subsubsection{}    \label{coproj-coinj-semi-mod-contra}
 Let $\S$ be a semialgebra over a coring $\C$ over a $k$\+algebra $A$.

\begin{lem}
 \textup{(a)} Assume that\/ $\C$ is a projective left $A$\module,
$\S$ is a coprojective left\/ $\C$\comodule{} and a\/ $\C/A$\+coflat
right\/ $\C$\comodule, and the ring $A$ has a finite left homological
dimension.
 Then there exist
\begin{itemize}
 \item an exact functor assigning to any $A$\projective{} left\/
       $\S$\semimodule{} an $A$\+split injective morphism from it
       into a\/ $\C$\coprojective\/ $\S$\semimodule, and
 \item an exact functor assigning to any left\/ $\S$\semicontramodule{}
       a surjective morphism onto it from a\/ $\C/A$\coinjective\/
       $\S$\semicontramodule.
\end{itemize} \par
 \textup{(b)} Assume that\/ $\C$ is a flat right $A$\module,
$\S$ is a coflat right\/ $\C$\comodule{} and a\/ $\C/A$\coprojective{}
left\/ $\C$\comodule, and the ring $A$ has a finite left homological
dimension.
 Then there exist
\begin{itemize}
 \item an exact functor assigning to any $A$\injective\/ left\/
       $\S$\semicontramodule{} an $A$\+split surjective morphism onto
       it from a\/ $\C$\coinjective\/ $\S$\semicontramodule, and
 \item an exact functor assigning to any left\/ $\S$\semimodule{}
       an injective morphism from it into a\/ $\C/A$\coprojective\/
       $\S$\semimodule.
\end{itemize} \par
 \textup{(c)} When both the assumptions of\/ \textup{(a)} and\/
\textup{(b)} are satisfied, the two functors acting in categories of
semimodules (can be made to) agree and the two functors acting in
categories of semicontramodules (can be made to) agree.
\end{lem}

\begin{proof}
 The proof of the first assertion of part~(a) and the second assertion
of part~(b) is based on the construction completely analogous to that
of the proof of Lemma~\ref{coflat-semimodule-injection}, with (co)flat
(co)modules replaced by (co)projective ones, and the left and right
sides switches as needed.
 The only difference is that the inductive limit of a sequence of
coprojective comodules does not have to be coprojective, because
even the inductive limit of a sequence of projective modules does
not have to be projective.
 This obstacle is dealt with in the following way.

\begin{subla}
 Assume that\/ $\C$ is a projective left $A$\module.
 Let\/ $\cU_1\rarrow\cU_2\rarrow\cU_3\rarrow\cU_4\rarrow\dsb$ be
an inductive system of left\/ $\C$\comodule s, where the comodules\/
$\cU_{2i}$ are coprojective, while the morphisms of comodules\/
$\cU_{2i-1}\rarrow\cU_{2i+1}$ are injective and split over~$A$.
 Then the inductive limit\/ $\ilim \cU_j$ is a coprojective\/
$\C$\comodule.
\end{subla}

\begin{proof}
 Let us first show that for any $\C$\contramodule{} $\P$ there is
an isomorphism $\Cohom_\C(\ilim\cU_j,\P)=
\plim\Cohom_\C(\cU_j,\P)$.
 Denote by $G_j^\bu$ the bar complex
$$
 \dsb\lrarrow\Hom_A(\C\ot_A\C\ot_A\cU_j\;\P)\lrarrow
 \Hom_A(\C\ot_A\cU_j\;\P)\lrarrow\Hom_A(\cU_j,\P);
$$
we will denote the terms of this complex by upper indices, so that
$G_j^n=0$ for $n>0$ and $H^0(G_j^\bu)=\Cohom_\C(\cU_j,\P)$.
 Clearly, we have $H^0(\plim G_j^\bu)=\Cohom_\C(\ilim\cU_j,\P)$.
 Since the comodules $\cU_{2i}$ are coprojective, 
$H^n(G_{2i}^\bu)=0$ for $n\ne0$, as the complex $G_{2i}^\bu$
can be obtained by applying the functor $\Cohom_\C(\cU_{2i},{-})$
to the complex of $\C$\contramodule s $\dsb\rarrow\Hom_A(\C\ot_A\C\;\P)
\rarrow\Hom_A(\C,\P)$, which is exact except at degree~0.
 Since the maps of $A$\module s $\cU_{2i-1}\rarrow\cU_{2i+1}$ are
split injective, the morphisms of complexes $G_{2i+1}^\bu\rarrow
G_{2i-1}^\bu$ are surjective.
 Therefore, $\plim^1 G_j^\bu=\plim^1 G_{2i-1}^\bu=0$, hence there is
a ``universal coefficients'' sequence~\cite{Wei}
 $$
  0\lrarrow\plim^1 H^{n-1}(G_j^\bu)\lrarrow H^n(\plim G_j^\bu)
  \lrarrow\plim H^n(G_j^\bu)\lrarrow0.
 $$ 
 In particular, for $n=0$ we obtain the desired isomorphism
$H^0(\plim G_j^\bu)=\plim H^0(G_j^\bu)$, because
$\plim^1 H^{-1}(G_j^\bu)=\plim^1 H^{-1}(G_{2i}^\bu)=0$.

 Now for any exact triple of $\C$\contramodule s $\P'\to\P \to\P''$
we have an exact triple of projective systems $\Cohom_\C(\cU_{2i},\P')
\rarrow\Cohom_\C(\cU_{2i},\P)\rarrow \Cohom_\C(\cU_{2i},\P'')$ and
$\plim^1\Cohom_\C(\cU_{2i},\P')=\plim^1\Cohom_\C(\cU_{2i-1},\P')=0$,
hence the triple remains exact after passing to the projective limit.
\end{proof}

\begin{sublb}
 Assume that\/ $\C$ is a flat right $A$\module.
 Let\/ $\cU_1\rarrow\cU_2\rarrow\cU_3\rarrow\cU_4\rarrow\dsb$ be
an inductive system of left\/ $\C$\comodule s, where the comodules\/
$\cU_{2i}$ are\/ $\C/A$\coprojective, while the morphisms of
comodules\/ $\cU_{2i-1}\rarrow\cU_{2i+1}$ are injective.
 Then the inductive limit\/ $\ilim \cU_j$ is a\/ $\C/A$\coprojective\/
$\C$\comodule.
\end{sublb}

\begin{proof}
 Analogous to the proof of Sublemma A, the only changes being that
$\P$, $\P'$, $\P''$ are now $A$\injective{} $\C$\contramodule s and
the complex $\dsb\rarrow\Hom_A(\C\ot_A\C\;\P)\rarrow\Hom_A(\C,\P)
\rarrow\P$ is an $A$\+split exact sequence of $A$\injective{}
$\C$\contramodule s.
\end{proof}

 Proof of the first assertion of part~(b): for any $A$\injective{}
$\C$\contramodule{} $\P$, set $\gG(\P)=\Hom_A(\C,\P)$.
 Then the contraaction map $\gG(\P)\rarrow\P$ is a surjective morphism
of $\C$\contramodule s, the contramodule $\gG(\P)$ is coinjective, and
the kernel of this morphism is $A$\injective.
 Now let $\bP$ be an $A$\injective{} left $\S$\semicontramodule.
 The semicontraaction map $\bP\rarrow\Cohom_\C(\S,\bP)$ is an injective
morphism of $A$\injective{} $\S$\semicontramodule s; let $\bgK(\bP)$
denote its cokernel.
 The map $\Cohom_\C(\S,\gG(\bP))\rarrow\Cohom_\C(\S,\bP))$
is a surjective morphism of $\S$\semicontramodule s with
an $A$\injective{} kernel $\Cohom_\C(\S\;\ker(\gG(\bP)\to\bP))$.
 Let $\bQ(\bP)$ be the kernel of the composition
$\Cohom_\C(\S,\gG(\bP))\rarrow\Cohom_\C(\S,\bP)\rarrow\bgK(\bP)$.
 Then the composition of maps $\bQ(\bP)\rarrow\Cohom_\C(\S,\gG(\bP))
\rarrow\Cohom_\C(\S,\bP)$ factorizes through the injection
$\bP\rarrow\Cohom_\C(\S,\bP)$, so there is a natural surjective
morphism of $\S$\semicontramodule s $\bQ(\bP)\rarrow\bP$.
 The kernel of the map $\bQ(\bP)\rarrow\bP$ is isomorphic to the kernel
of the map $\Cohom_\C(\S,\gG(\bP))\rarrow\Cohom_\C(\S,\bP)$, hence
both $\ker(\bQ(\bP)\to\bP)$ and~$\bQ(\bP)$ are injective $A$\module s.

 Notice that the semicontramodule morphism $\bQ(\bP)\rarrow\bP$ can be
extended to a contramodule morphism $\Cohom_\C(\S,\gG(\bP))\rarrow\bP$.
 Indeed, the map $\bQ(\bP)\rarrow\bP$ can be presented as
the composition $\bQ(\bP)\rarrow\Cohom_\C(\S,\gG(\bP))\rarrow
\Cohom_\C(\S,\bP)\rarrow\bP$, where the map $\Cohom_\C(\S,\bP)
\rarrow\bP$ is induced by the semiunit morphism $\C\rarrow\S$ of
the semialgebra~$\S$.

 Iterating this construction, we obtain a projective system of
$\C$\contramodule{} morphisms $\bP\larrow\Cohom_\C(\S,\gG(\bP))
\larrow\bQ(\bP)\larrow\Cohom_\C(\S,\gG(\bQ(\bP)))\larrow
\bQ(\bQ(\bP))\larrow\dsb$, where the maps $\bP\larrow\bQ(\bP)
\larrow\bQ(\bQ(\bP))\larrow\dsb$ are $A$\+split surjective morphisms
of $A$\injective{} $\S$\contramodule s, while the $\C$\contramodule s
$\Cohom_\C(\S,\gG(\bP))$, \ $\Cohom_\C(\S,\gG(\bQ(\bP)))$,~\dots\
are coinjective.
 Denote by $\bgF(\bP)$ the projective limit of this system;
then $\bgF(\bP)\rarrow\bP$ is an $A$\+split surjective morphism
of $\S$\semicontramodule s, while coinjectivity of
the $\C$\contramodule{} $\bgF(\bP)$ follows from the next Sublemma.

\begin{sublc}
 Assume that\/ $\C$ is a flat right $A$\module.
 Let\/ $\gU_1\larrow\gU_2\larrow\gU_3\larrow\gU_4\larrow\dsb$ be
a projective system of left\/ $\C$\contramodule s, where
the contramodules\/ $\gU_{2i}$ are coinjective, while the morphisms of
contramodules\/ $\gU_{2i+1}\rarrow\gU_{2i-1}$ are surjective and
split over~$A$.
 Then the projective limit\/ $\plim\gU_j$ is a coinjective\/
$\C$\contramodule.
\end{sublc}

\begin{proof}
 Completely analogous to the proof of Sublemma~A.
 One considers the projective system of bar-complexes
$\dsb\rarrow\Hom_A(\C\ot_A\C\ot_A\M\;\gU_j)\rarrow
\Hom_A(\C\ot_A\M\;\gU_j)\rarrow\Hom_A(\M,\gU_j)$, etc.
\end{proof}
 
 The proof of the second assertion of part~(a) is based on the same
construction; the only changes are that $A$\module s are no longer
injective, for any left $\C$\contramodule{} $\bP$
the $\C$\contramodule{} $\gG(\bP)=\Hom_A(\C,\bP)$ is
$\C/A$\coinjective, and therefore the $\S$\semicontramodule{}
$\Cohom_\C(\S,\gG(\bP))$ is $\C/A$\coinjective.
 The projective limit $\bgF(\bP)$ is $\C/A$\coinjective{} according
to the following Sublemma.

\begin{subld}
 Assume that\/ $\C$ is a projective left $A$\module.
 Let\/ $\gU_1\larrow\gU_2\larrow\gU_3\larrow\gU_4\larrow\dsb$ be
a projective system of left\/ $\C$\contramodule s, where
the contramodules\/ $\gU_{2i}$ are $\C/A$\coinjective, while
the morphisms of contramodules\/ $\gU_{2i+1}\rarrow\gU_{2i-1}$
are surjective.
 Then the projective limit\/ $\plim\gU_j$ is a\/ $\C/A$\coinjective\/
$\C$\contramodule. \qed
\end{subld}
%\begin{proof}
% Analogous to the proofs of Sublemmas~B and~C.
%\end{proof}

 Both functors $\bgF$ are exact, since the cokernels of injective
maps, the kernels of surjective maps, and the projective limits
of Mittag-Leffler sequences of $k$\module s preserve exact triples.
 Part~(c) is clear from the constructions. 
\end{proof}

\subsection{Semihomomorphisms}

\subsubsection{}
 Assume that the coring $\C$ is a projective left $A$\module,
the semialgebra $\S$ is a projective left $A$\module{} and
a $\C/A$\+coflat right $A$\module, and the ring $A$ has a finite
left homological dimension.
 Let $\bM$ be an $A$\projective{} left $\S$\semimodule{} and
$\bP$ be a left $\S$\semicontramodule.
 The $k$\module{} of \emph{semihomomorphisms} $\SemiHom_\S(\bM,\bP)$
is defined as the kernel of the pair of maps $\Cohom_\C(\bM,\bP)
\birarrow\Cohom_\C(\S\oc_\C\bM\;\bP)=\Cohom_\C(\bM,\Cohom_\C(\S,\bP))$
one of which is induced by the $\S$\+semiaction in~$\bM$ and
the other by the $\S$\+semicontraaction in~$\bP$.

 For any $A$\projective{} left $\C$\comodule{} $\L$ and any
left $\S$\semicontramodule{} $\bP$ there is a natural isomorphism
$\SemiHom_\S(\S\oc_\C\L\;\bP)\simeq\Cohom_\C(\L,\bP)$.
 Analogously, for any $A$\projective{} left $\S$\semimodule{} $\bM$
and any left $\C$\contramodule{} $\Q$ there is a natural isomorphism
$\SemiHom_\S(\bM,\Cohom_\C(\S,\Q))\simeq\Cohom_\C(\bM,\Q)$.
These assertions follow from Lemma~\ref{induced-tensor-cotensor}.

\subsubsection{}
 Assume that the coring $\C$ is a flat right $A$\module,
the semialgebra $\S$ is a flat right $A$\module{} and
a $\C/A$\coprojective{} left $A$\module, and the ring $A$ has a finite
left homological dimension.
 Let $\bM$ be a left $\S$\semimodule{} and $\bP$ be an $A$\injective{}
left $\S$\semicontramodule.
 As above, the $k$\module{} of \emph{semihomomorphisms}
$\SemiHom_\S(\bM,\bP)$ is defined as the kernel of the pair
of maps $\Cohom_\C(\bM,\bP)\birarrow\Cohom_\C(\S\oc_\C\bM\;\bP)=
\Cohom_\C(\bM,\Cohom_\C(\S,\bP))$ one of which is induced by
the $\S$\+semiaction in~$\bM$ and the other by
the $\S$\+semicontraaction in~$\bP$. 

 For any left $\C$\comodule{} $\L$ and any $A$\injective{} left
$\S$\semicontramodule{} $\bP$ there is a natural isomorphism
$\SemiHom_\S(\S\oc_\C\L\;\bP)\simeq\Cohom_\C(\L,\bP)$.
 Analogously, for any left $\S$\semimodule{} $\bM$ and any
$A$\injective{} left $\C$\contramodule{} $\Q$ there is a natural
isomorphism $\SemiHom_\S(\bM,\Cohom_\C(\S,\Q))\simeq\Cohom_\C(\bM,\Q)$.

 Notice that even under the strongest of our assumptions on $A$, $\C$
and~$\S$, the $A$\+projectivity of $\bM$ or the $A$\+injectivity
of~$\bP$ is still needed to guarantee that the triple cohomomorphisms
$\Cohom_\C(\S\oc_\C\bM\;\bP)$ are associative.

\subsubsection{}
 If the coring $\C$ is a projective left $A$\module{} and
the semialgebra $\S$ is a coprojective left $\C$\comodule, one
can define the module of semihomomorphisms from a $\C$\coprojective{}
left $\S$\semimodule{} into an arbitrary left $\S$\semicontramodule.
 In these assumptions, a $\C$\coprojective{} left $\S$\semimodule{}
$\bM$ is called \emph{semiprojective} if the functor of
semihomomorphisms from~$\bM$ is exact on the abelian category of left
$\S$\semicontramodule s.
 The $\S$\semimodule{} induced from a coprojective $\C$\comodule{}
is semiprojective.
 Any semiprojective $\S$\semimodule{} is semiflat.

 If the coring $\C$ is a flat right $A$\module{} and
the semialgebra $\S$ is a coflat right $\C$\comodule, one can define
the module of semihomomorphisms from an arbitrary left
$\S$\semimodule{} into a $\C$\coinjective{} left $\S$\semicontramodule.
 In these assumptions, a $\C$\coinjective{} left $\S$\semicontramodule{}
$\bP$ is called \emph{semiinjective} if the functor of
semihomomorphisms into~$\bP$ is exact on the abelian category of
left $\S$\semimodule s.
 The $\S$\semicontramodule{} coinduced from a coinjective
$\C$\contramodule{} is semiinjective.

 When the ring $A$ is semisimple, the module of semihomomorphisms
from an arbitrary $\S$\semimodule{} into an arbitrary
$\S$\semicontramodule{} is defined without any conditions on
the coring~$\C$ and the semialgebra~$\S$.

\subsubsection{}    \label{semihom-associative}
 Let $\S$ be a semialgebra over a coring~$\C$ over a $k$\+algebra $A$
and $\T$ be a semialgebra over a coring~$\D$ over a $k$\+algebra $B$.
 Let $\bK$ be an $\S$\+$\T$\bisemimodule{} and $\bP$ be 
a left $\S$\semicontramodule.
 We would like to define a left $\T$\semicontramodule{} structure on
the module of semihomomorphisms $\SemiHom_\S(\bK,\bP)$.

 Assume that multiple cohomomorphisms of the form
$\Cohom_\C(\S\oc_\C\bK\oc_\D\T\oc_\D\dsb\oc_\D\T\;\bP)$
are associative.
 Then, in particular, the $k$\module s of semihomomorphisms
$\SemiHom_\S(\bK\oc_\D\T\oc_\D\dsb\oc_\D\T\;\bP)$ can be defined.
 Assume in addition that multiple cohomomorphisms of the form
$\Cohom_\C(\bK\oc_\D\T\oc_\D\dsb\oc_\D\T\;\bP)$
are associative.
 Then the semihomomorphism modules
$\SemiHom_\S(\bK\oc_\D\T\oc_\D\dsb\oc_\D\T\;\bP)$
have natural left $\D$\contramodule{} structures as kernels of
$\D$\contramodule{} morphisms.
 Assume that multiple cohomomorphisms of the form
$\Cohom_\D(\T\oc_\D\dsb\oc_\D\T\;\SemiHom_\S(\bK,\bP))$
are also associative.
 Finally, assume that the cohomomorphisms from $\T^{\suboc\.m}$
preserve the kernel of the pair of morphisms
$\Cohom_\C(\bK,\bP)\birarrow\Cohom_\C(\S\oc_\C\bK,\bP)$ for $m=1$
and~$2$, that is the contramodule morphisms $\Cohom_\D(\T^{\suboc\.m},
\SemiHom_\S(\bK,\bP))\lrarrow\SemiHom_\S(\bK\oc_\D\T^{\suboc\.m}\;\bP)$
are isomorphisms.
 Then one can define an associative and unital semicontraaction
morphism\- $\SemiHom_\S(\bK,\bP)\rarrow
\Cohom_\D(\T,\SemiHom_\S(\bK,\bP))$ taking
the semihomomorphisms over~$\S$ from the right $\T$\+semiaction
morphism $\bK\oc_\D\T\rarrow\bK$ into the semicontramodule~$\bP$.

 For example, if $\D$ is a projective left $B$\module, $\T$ is
a coprojective left $\D$\comodule, $A$ has a finite left homological
dimension, and either $\C$ is a projective left $A$\module, $\S$ is
a projective left $A$\module{} and a $\C/A$\+coflat right $\C$\comodule,
and $\bK$ is a projective left $A$\module, or $\C$ is a flat right
$A$\module, $\S$ is a flat right $A$\module{} and
a $\C/A$\coprojective{} left $\C$\comodule, and $\bP$ is an injective
left $A$\module, then the module of semihomomorphisms
$\SemiHom_\S(\bK,\bP)$ has a natural left $\T$\semicontramodule{}
structure.
 Since the category of left $\T$\semicontramodule s is abelian in this
case, the $\T$\semicontramodule{} $\SemiHom_\S(\bK,\bP)$ can be simply
defined as the kernel of the pair of semicontramodule morphisms
$\Cohom_\C(\bK,\bP)\birarrow\Cohom_\C(\S\oc_\C\bK,\bP)$.

\begin{prop}
 Let\/ $\bM$ be a left $\T$\semimodule, $\bK$ be an\/
$\S$\+$\T$\bisemimodule, and\/ $\P$ be a left $\S$\semicontramodule.
 Then the iterated semihomomorphism modules\/ $\SemiHom_\S(\bK\os_\T
\bM\;\bP)$ and\/ $\SemiHom_\T(\bM,\SemiHom_\S(\bK,\bP))$ are 
well-defined and naturally isomorphic, at least, in the following cases:
\begin{enumerate}
 \item $\D$ is a projective left $B$\module, $\T$ is a coprojective
       left\/ $\D$\comodule, $\bM$ is a coprojective left\/
       $\D$\comodule, $\C$ is a flat right $A$\module, $\S$ is a coflat
       right\/ $\C$\comodule, and\/ $\bP$ is a coinjective left\/
       $\C$\contramodule;
 \item $\D$ is a projective left $B$\module, $\T$ is a coprojective
       left\/ $\D$\comodule, $\bM$ is a semiprojective left\/
       $\T$\semimodule, and either
       \begin{itemize}
         \item $\C$ is a projective left $A$\module, $\S$ is 
               a coprojective left\/ $\C$\comodule, and\/ $\bK$ is
               a coprojective left\/ $\C$\comodule, or
         \item $\C$ is a projective left $A$\module, $\S$ is
               a projective left $A$\module{} and a\/ $\C/A$\+coflat
               right\/ $\C$\comodule, the ring $A$ has a finite left
               homological dimension, and\/ $\bK$ is a projective left
               $A$\module, or
         \item $\C$ is a flat right $A$\module, $\S$ is a flat
               right $A$\module{} and a\/ $\C/A$\coprojective{} left\/
               $\C$\comodule, and\/ $\bP$ is an injective left
               $A$\module, or
         \item the ring $A$ is semisimple;
       \end{itemize}
 \item $\C$ is a flat right $A$\module, $\S$ is a coflat
       right\/ $\C$\comodule, $\bP$ is a semiinjective left\/
       $\S$\semicontramodule, and either
       \begin{itemize}
         \item $\D$ is a flat right $B$\module, $\T$ is a coflat
               right\/ $\D$\comodule, and\/ $\bK$ is a coflat right\/
               $\D$\comodule, or
         \item $\D$ is a flat right $B$\module, $\T$ is a flat right
               $B$\module{} and a\/ $\D/B$\coprojective{} left\/
               $\D$\comodule, the ring $B$ has a finite left
               homological dimension, and\/ $\bK$ is a flat right
               $B$\module, or
         \item $\D$ is a projective left $B$\module, $\T$ is
               a projective left $B$\module{} and a\/ $\D/B$\+coflat
               right\/ $\D$\comodule, the ring $B$ has a finite left
               homological dimension, and\/ $\bM$ is a projective left
               $B$\module, or
         \item the ring $B$ is semisimple;
       \end{itemize}
 \item $\D$ is a projective left $B$\module, $\T$ is a coprojective
       left\/ $\D$\comodule, $\bM$ is a coprojective left\/
       $\D$\comodule, and either
       \begin{itemize}
          \item $\C$ is a projective left $A$\module, $\S$ is
                a coprojective left\/ $\C$\comodule, and\/ $\bK$ as
                a right\/ $\T$\semimodule{} with a left\/
                $\C$\comodule{} structure is induced from a\/
                $\C$\coprojective\/ $\C$\+$\D$\bicomodule, or
          \item $\C$ is a projective left $A$\module, $\S$ is
                a projective left $A$\module{} and a\/ $\C/A$\+coflat
                right\/ $\C$\comodule, the ring $A$ has a finite left
                homological dimension, and\/ $\bK$ as a right\/
                $\T$\semimodule{} with a left\/ $\C$\comodule{}
                structure is induced from an $A$\projective\/
                $\C$\+$\D$\bicomodule, or
          \item $\C$ is a flat right $A$\module, $\S$ is a flat right
                $A$\module{} and a\/ $\C/A$\coprojective{} left\/
                $\C$\comodule, the ring $A$ has a finite left
                homological dimension, $\bK$ as a right\/
                $\T$\semimodule{} with a left\/ $\C$\comodule{}
                structure is induced from a\/ $\C$\+$\D$\bicomodule,
                and\/ $\bP$ is an injective left $A$\module, or
          \item the ring $A$ is semisimple and\/ $\bK$ as a right\/
                $\T$\semimodule{} with a left\/ $\C$\comodule{}
                structure is induced from a\/ $\C$\+$\D$\bicomodule;
       \end{itemize}
 \item $\C$ is a flat right $A$\module, $\S$ is a coflat right\/
       $\C$\comodule, $\bP$ is a coinjective left\/ $\C$\contramodule,
       and either
       \begin{itemize}
          \item $\D$ is a flat right $B$\module, $\T$ is a coflat
                right $\D$\comodule, and\/ $\bK$ as a left\/
                $\S$\semimodule{} with a right\/ $\D$\comodule{}
                structure is induced from a\/ $\D$\+coflat\/
                $\C$\+$\D$\bicomodule, or
          \item $\D$ is a flat right $B$\module, $\T$ is a flat
                right $B$\module{} and a\/ $\D/B$\coprojective{}
                left\/ $\D$\comodule, the ring $B$ has a finite left
                homological dimension, and\/ $\bK$ as a left\/
                $\S$\semimodule{} with a right\/ $\D$\comodule{}
                structure is induced from a $B$\+flat\/
                $\C$\+$\D$\bicomodule, or
          \item $\D$ is a projective left $B$\module, $\T$ is
                a projective left $B$\module{} and a\/ $\D/B$\+coflat
                right\/ $\D$\comodule, the ring $B$ has a finite left
                homological dimension, $\bK$ as a left\/
                $\S$\semimodule{} with a right\/ $\D$\comodule{}
                structure is induced from a\/ $\C$\+$\D$\bicomodule,
                and\/ $\bM$ is a projective left $B$\module, or
          \item the ring $B$ is semisimple and\/ $\bK$ as a left\/
                $\S$\semimodule{} with a right\/ $\D$\comodule{}
                structure is induced from a\/ $\C$\+$\D$\bicomodule.
       \end{itemize}
\end{enumerate}
 More precisely, in all cases in this list the natural maps from both
iterated semihomomorphism modules under consideration into
the iterated cohomomorphism module\/ $\Cohom_\C(\bK\oc_\D\bM\;\bP)
\simeq\Cohom_\D(\bM,\Cohom_\C(\bK,\bP))$ are injective, their images
coincide and are equal to the intersection of two submodules\/
$\SemiHom_\S(\bK\oc_\D\bM\;\bP)$ and\/
$\SemiHom_\T(\bM,\Cohom_\C(\bK,\bP))$ in this $k$\module.
\end{prop}

\begin{proof}
Analogous to the proof of Proposition~\ref{semitensor-associative}
(see also the proof of Proposition~\ref{cohom-associative}).
\end{proof}

\Section{Derived Functor SemiExt}

\subsection{Contraderived categories}  \label{contraderived-categories}
 Let $\sA$ be an exact category in which all infinite products exist
and the functors of infinite product are exact.
 A complex $C^\bu$ over~$\sA$ is called \emph{contraacyclic} if
it belongs to the minimal triangulated subcategory $\Acycl^\ctr(\sA)$
of the homotopy category $\Hot(\sA)$ containing all the total
complexes of exact triples ${}'\!K^\bu\to K^\bu\to {}''\!K^\bu$
of complexes over~$\sA$ and closed under infinite products.
 Any contraacyclic complex is acyclic.
 It follows from the next Lemma that any acyclic complex bounded
from above is contraacyclic.

\begin{lem}
 Let $\dsb\to P^{-1,\bu} \to P^{0,\bu}\to0$ be an exact sequence,
bounded from above, of arbitrary complexes over~$\sA$.
 Then the total complex $T^\bu$ of the bicomplex $P^{\bu,\bu}$
constructed by taking infinite products along the diagonals is
contraacyclic.
\end{lem}

\begin{proof}
 See the proof of Lemma~\ref{coderived-categories}.
\end{proof}

 The category of contraacyclic complexes $\Acycl^\ctr(\sA)$ is a thick
subcategory of the homotopy category $\Hot(\sA)$, since it is 
a triangulated subcategory with infinite products.
 The \emph{contraderived category} $\sD^\ctr(\sA)$ of an exact
category~$\sA$ is defined as the quotient category
$\Hot(\sA)/\Acycl^\ctr(\sA)$.

\begin{rmk}
 One can check that for any exact category $\sA$ and any thick
subcategory $\sT$ in $\Hot(\sA)$ contained in the thick subcategory
of acyclic complexes, containing all bounded acyclic complexes, and
containing with every exact complex its subcomplexes and quotient
complexes of canonical filtration, the groups of homomorphisms
$\Hom_{\Hot(\sA)/\sT}(X,Y[i])$ between complexes with a single nonzero
term coincide with the Yoneda extension groups $\Ext_\sA^i(X,Y)$.
 Moreover, the natural functors $\Hot^{+/-/b}(\sA)/
(\sT\cap\Hot^{+/-/b}(\sA))\rarrow\Hot(\sA)/\sT$ between
the ``$\sT$-derived categories'' with various bounding conditions
are all fully faithful.
% Besides, when $\sA$ is an abelian category there is a (degenerate)
%t-structure on $\Hot(\sA)/\sT$ formed by the full subcategories of
%complexes concentrated in nonpositive and nonnegative degrees.
 In particular, these assertions hold if $\sT\subset\Hot(\sA)$ consists
of acyclic complexes and contains either all acyclic complexes
bounded from above or all acyclic complexes bounded from below.
\end{rmk}

\subsection{Coprojective and coinjective complexes}
\label{co-proj-inj-complexes}
 Let $\C$ be a coring over a $k$\+al\-ge\-bra~$A$.
 The complex of cohomomorphisms $\Cohom_\C(\M^\bu,\P^\bu)$ from
a complex of left $\C$\comodule s $\M^\bu$ into a complex of
left $\C$\contramodule s $\P^\bu$ is defined as the total complex
of the bicomplex $\Cohom_\C(\M^i,\P^j)$, constructed by taking
infinite products along the diagonals.

 If $\C$ is a projective left $A$\module, the category of left
$\C$\contramodule s is an abelian category with exact functors
of infinite products, so the contraderived category
$\sD^\ctr(\C\contra)$ is defined.
 When speaking about \emph{contraacyclic complexes} of
$\C$\contramodule s, we will always mean contraacyclic complexes
with respect to the abelian category of $\C$\contramodule s,
unless another exact category of $\C$\contramodule s is explicitly
mentioned.

 Assuming that $\C$ is a projective left $A$\module, a complex
of left $\C$\comodule s $\M^\bu$ is called \emph{coprojective} if
the complex $\Cohom_\C(\M^\bu,\P^\bu)$ is acyclic whenever
a complex of left $\C$\contramodule s $\P^\bu$ is contraacyclic.
 Assuming that $\C$ is a flat right $A$\module, a complex of left
$\C$\contramodule s $\P^\bu$ is called \emph{coinjective} if
the complex $\Cohom_\C(\M^\bu,\P^\bu)$ is acyclic whenever
a complex of left $\C$\comodule s $\M^\bu$ is coacyclic.

\begin{lem}
 \textup{(a)} Any complex of coprojective\/ $\C$\comodule s is
coprojective. \par
 \textup{(b)} Any complex of coinjective\/ $\C$\contramodule s
is coinjective.
\end{lem}

\begin{proof}
 Argue as in the proof of Lemma~\ref{coflat-complexes}, using the fact
that the functor of cohomomorphisms of complexes maps infinite direct
sums in the first argument into infinite products and preserves
infinite products in the second argument.
\end{proof}

 If the ring $A$ has a finite left homological dimension, then any
coprojective complex of left $\C$\comodule s is a projective complex
of $A$\module s in the sense of~\ref{prelim-unbounded-ext} and any
coinjective complex of left $\C$\contramodule s is an injective
complex of $A$\module s.
 The complex of $\C$\comodule s $\C\ot_A U^\bu$ coinduced from
a projective complex of $A$\module s $U^\bu$ is coprojective and
the complex of $\C$\contramodule s $\Hom_A(\C,V^\bu)$ induced
from an injective complex of $A$\module s is coinjective.

\subsection{Semiderived categories}
 Let $\S$ be a semialgebra over a coring $\C$.
 Assume that $\C$ is a projective left $A$\module{} and the semialgebra
$\S$ is a coprojective left $\C$\comodule, so that the category of
left $\S$\semicontramodule s is abelian.
 The \emph{semiderived category} of left $\S$\semicontramodule s 
$\sD^\si(\S\sicntr)$ is defined as the quotient category of the homotopy
category $\Hot(\S\sicntr)$ by the thick subcategory
$\Acycl^{\ctrd\C}(\S\sicntr)$ of complexes of $\S$\semicontramodule s
that are \emph{contraacyclic as complexes of\/ $\C$\contramodule s}.

\subsection{Semiprojective and semiinjective complexes}
\label{semi-inj-proj-complexes}
 Let $\S$ be a semialgebra.
 The complex of semihomomorphisms $\SemiHom_\S(\bM^\bu,\bP^\bu)$
from a complex of left $\S$\semimodule s $\bM^\bu$ to a complex of
left $\S$\semicontramodule s $\bP^\bu$ is defined as the total
complex of the bicomplex $\SemiHom_\S(\bM^i,\bP^j)$, constructed
by taking infinite products along the diagonals.
 Of course, appropriate conditions must be imposed on $\S$, $\bM^\bu$,
and $\bP^\bu$ for this definition to make sense.

 Assume that the coring $\C$ is a projective left $A$\module{} and
a flat right $A$\module, the semialgebra $\S$ is a coprojective
left $\S$\semimodule{} and a coflat right $\S$\semimodule, and
the ring $A$ has a finite left homological dimension.

 A complex of $A$\projective{} left $\S$\semimodule s $\bM^\bu$
is called \emph{semiprojective} if the complex
$\SemiHom_\S(\bM^\bu,\bP^\bu)$ is acyclic whenever a complex of left
$\S$\semicontramodule s $\bP^\bu$ is $\C$\contraacyclic.
 Any semiprojective complex of $\S$\semimodule s is a coprojective
complex of $\C$\comodule s.
 The complex of $\S$\semimodule s $\S\oc_\C\L^\bu$ induced from
a coprojective complex of $A$\+flat $\C$\comodule s is semiprojective.
 Any semiprojective complex of $\S$\semimodule s is semiflat.
 Analogously, a complex of $A$\injective{} left $\S$\semicontramodule s
$\bP^\bu$ is called \emph{semiinjective} if the complex
$\SemiHom_\S(\bM^\bu,\bP^\bu)$ is acyclic whenever a complex of left
$\S$\semimodule s $\bM^\bu$ is $\C$\coacyclic.
 Any semiinjective complex of $\S$\semicontramodule s is a coinjective
complex of $\C$\contramodule s.
 The complex of $\S$\semicontramodule s $\Cohom_\C(\S,\Q^\bu)$
coinduced from a coinjective complex of $A$\injective{}
$\C$\contramodule s is semiinjective.

 Notice that not every complex of semiprojective semimodules
is semiprojective and not every complex of semiinjective
semicontramodules is semiinjective.
 On the other hand, any complex of semiprojective semimodules bounded
from above is semiprojective.
 Moreover, if $\dsb\to\bM^{-1,\bu}\to\bM^{0,\bu}\to0$ is a complex,
bounded from above, of semiprojective complexes of $\S$\semimodule s,
then the total complex $\bE^\bu$ of the bicomplex $\bM^{\bu,\bu}$
constructed by taking infinite direct sums along the diagonals
is semiprojective.
 Indeed, the category of semiprojective complexes is closed under
shifts, cones, and infinite direct sums, so one can apply
Lemma~\ref{semiflat-complexes}.
 Analogously, any complex of semiinjective semicontramodules bounded
from below is semiinjective.
 Moreover, if $0\to\bP^{0,\bu}\to\bP^{1,\bu}\to\dsb$ is a complex,
bounded from below, of semiinjective complexes of 
$\S$\semicontramodule s, then the total complex $\bgE^\bu$ of
the bicomplex $\bP^{\bu,\bu}$ constructed by taking infinite products
along the diagonals is semiinjective.
 Indeed, the category of semiinjective complexes is closed under shifts,
cones, and infinite products, so one can apply the following Lemma.

\begin{lem}
 Let\/ $0\to P^{0,\bu}\to P^{1,\bu}\to\dsb$ be a complex, bounded from
below, of arbitrary complexes over an additive category\/~$\sA$ where
infinite products exist.
 Then the total complex $E^\bu$ of the bicomplex $P^{\bu,\bu}$ up to
the homotopy equivalence can be obtained from the complexes
$P^{i,\bu}$ using the operations of shift, cone, and infinite product.
\end{lem}

\begin{proof}
 See the proof of Lemma~\ref{semiflat-complexes}.
\end{proof}

\subsection{Main theorem for comodules and contramodules}
\label{coext-main-theorem}
 Assume that the coring $\C$ is a projective left and a flat right
$A$\module{} and the ring $A$ has a finite left homological dimension.

\begin{thm}
 \textup{(a)} The functor mapping the quotient category of
the homotopy category of complexes of coprojective left\/
$\C$\comodule s (coprojective complexes of left\/ $\C$\comodule s)
by its intersection with the thick subcategory of coacyclic complexes
of\/ $\C$\comodule s into the coderived category of left\/
$\C$\comodule s is an equivalence of triangulated categories. \par
 \textup{(b)} The functor mapping the quotient category of
the homotopy category of complexes of coinjective left\/
$\C$\contramodule s (coinjective complexes of left\/
$\C$\contramodule s) by its intersection with the thick subcategory
of contraacyclic complexes of\/ $\C$\contramodule s into
the contraderived category of left\/ $\C$\contramodule s
is an equivalence of triangulated categories.
\end{thm}

\begin{proof}
 The proof of part~(a) is completely analogous to the proof of
Theorem~\ref{cotor-main-theorem}.
 It is based on the same constructions of resolutions
$\boL_1$ and $\boR_2$, and uses the result of
Lemma~\ref{proj-inj-co-contra-module}(a) instead of
Lemma~\ref{flat-comodule-surjection}.

 To prove part~(b), we will show that any complex of left
$\C$\contramodule s $\gK^\bu$ can be connected with a complex
of coinjective $\C$\contramodule s in a functorial way by
a chain of two morphisms $\gK^\bu\rarrow\boL_2(\gK^\bu)\larrow
\boL_2\boR_1(\gK^\bu)$ with contraacyclic cones.
 Moreover, if the complex $\gK^\bu$ is a complex of coinjective
$\C$\contramodule s (coinjective complex of $\C$\contramodule s),
then the intermediate complex $\boL_2(\gK^\bu)$ is also
a complex of coinjective $\C$\contramodule s (coinjective complex
of $\C$\contramodule s).
 Then we will apply Lemma~\ref{cotor-main-theorem} in the way
explained in the end of the proof of Theorem~\ref{cotor-main-theorem}.

 Let $\gK^\bu$ be a complex of left $\C$\contramodule s.
 Let $\P\rarrow\gI(\P)$ denote the functorial injective morphism from
an arbitrary left $\C$\contramodule{} $\P$ into an $A$\injective{}
$\C$\contramodule{} $\gI(\P)$ constructed in
Lemma~\ref{proj-inj-co-contra-module}(b).
 The functor $\gI$ is the direct sum of a constant functor
$\P\mpsto\gI(0)$ and a functor $\gI^+$ sending zero morphisms
to zero morphisms.
 For any $\C$\contramodule{} $\P$, the contramodule $\gI^+(\P)$
is $A$\injective{} and the morphism $\P\rarrow\gI^+(\P)$ is injective.
 Set $\gI^0(\gK^\bu)=\gI^+(\gK^\bu)$, \
$\gI^1(\gK^\bu)=\gI^+(\coker(\gK^\bu\to\gI^0(\gK^\bu)))$, etc.
 For $d$ large enough, the cokernel $\gZ(\gK^\bu)$ of the morphism
$\gI^{d-2}(\gK^\bu)\rarrow\gI^{d-1}(\gK^\bu)$ will be a complex of
$A$\injective{} $\C$\contramodule s.
 Let $\boR_1(\gK^\bu)$ be the total complex of the bicomplex
$$
 \gI^0(\gK^\bu)\lrarrow\gI^1(\gK^\bu)\lrarrow\dsb\lrarrow
 \gI^{d-1}(\gK^\bu)\lrarrow\gZ(\gK^\bu).
$$
 Then $\boR_1(\gK^\bu)$ is a complex of $A$\injective{}
$\C$\contramodule s and the cone of the morphism
$\gK^\bu\rarrow\boR_1(\gK^\bu)$ is the total complex of a finite
exact sequence of complexes of $\C$\contramodule s, and therefore,
a contraacyclic complex.

 Now let $\gR^\bu$ be a complex of $A$\injective{} left
$\C$\contramodule s.
 Consider the bar construction
$$
 \dsb\lrarrow\Hom_A(\C,\Hom_A(\C,\gR^\bu))\lrarrow
 \Hom_A(\C,\gR^\bu).
$$
 Let $\boL_2(\gR^\bu)$ be the total complex of this bicomplex,
constructed by taking infinite products along the diagonals.
 Then $\boL_2(\gR^\bu)$ is a complex of coinjective
$\C$\contramodule s.
 The functor $\boL_2$ can be extended to arbitrary complexes of
$\C$\contramodule s; for any complex $\gK^\bu$, the cone of
the morphism $\boL_2(\gK^\bu)\rarrow\gK^\bu$ is
contraacyclic by Lemma~\ref{contraderived-categories}.

 Finally, if $\gK^\bu$ is a coinjective complex of $\C$\contramodule s,
then $\boL_2(\gK^\bu)$ is also a coinjective complex of
$\C$\contramodule s, since the complex of cohomomorphisms from
a complex of left $\C$\comodule s $\M^\bu$ into $\boL_2(\gK^\bu)$
coincides with the complex of cohomomorphisms into $\gK^\bu$ from
the total cobar complex $\boR_2(\M^\bu)$, and the latter is
coacyclic whenever $\M^\bu$ is coacyclic.
\end{proof}

\begin{rmk}
 Another proof of Theorem (for complexes of coprojective comodules
and complexes of coinjective contramodules) can be deduced from
the results of Section~\ref{co-contra-correspondence-section}.
 In addition, it will follow that any coacyclic complex of coprojective
left $\C$\comodule s is contractible and any contraacyclic complex of
coinjective left $\C$\contramodule s is contractible
(see Remark~\ref{co-contra-ctrtor-definition}).
\end{rmk}

\subsection{Main theorem for semimodules and semicontramodules}
\label{semiext-main-theorem}
 Assume that the coring $\C$ is a projective left and a flat right
$A$\module, the semialgebra $\S$ is a coprojective left and
a coflat right $\C$\comodule, and the ring $A$ has a finite left
homological dimension.

\begin{thm}
 \textup{(a)} The functor mapping the quotient category of
the homotopy category of semiprojective complexes of $A$\projective{}
(\.$\C$\coprojective, semiprojective) left\/ $\S$\semimodule s by
its intersection with the thick subcategory of\/ $\C$\coacyclic{}
complexes of\/ $\S$\semimodule s into the semiderived category
of left\/ $\S$\semimodule s is an equivalence of triangulated
categories. \par
 \textup{(b)} The functor mapping the quotient category of
the homotopy category of semiinjective complexes of $A$\injective{}
(\.$\C$\coinjective, semiinjective) left\/ $\S$\semicontramodule s by
its intersection with the thick subcategory of\/ $\C$\contraacyclic{}
complexes of\/ $\S$\semicontramodule s into the semiderived category
of left\/ $\S$\semicontramodule s is an equivalence of triangulated
categories.
\end{thm}

\begin{proof}
 There are two approaches: one can argue as in \ref{cotor-main-theorem}
or as in~\ref{semitor-main-theorem}.
 Either way, the proof is based on the constructions of intermediate
resolutions $\boL_i$ and~$\boR_j$.
 For part~(a), it is the same constructions that were presented in
the proof of Theorem~\ref{semitor-main-theorem}.
 One just has to use the results of
Lemmas~\ref{proj-inj-semi-mod-contra}(a)
and~\ref{coproj-coinj-semi-mod-contra}(a)
instead of Lemmas~\ref{flat-semimodule-surjection}
and~\ref{coflat-semimodule-injection}.
 Let us introduce the analogous constructions for part~(b).

 Let $\bgK^\bu$ be a complex of left $\S$\semicontramodule s.
 Let $\bP\rarrow\bgI(\bP)$ denote the functorial injective
morphism from an arbitrary left $\S$\semicontramodule{} $\bP$
into an $A$\injective{} $\S$\semicontramodule{} $\bgI(\bP)$
constructed in Lemma~\ref{proj-inj-semi-mod-contra}(b).
 The functor $\bgI$ is the direct sum of a constant functor
$\bP\mpsto\bgI(0)$ and a functor $\bgI^+$ sending zero
morphisms to zero morphisms.
 For any $\S$\semicontramodule{} $\bP$, the semicontramodule
$\bgI^+(\bP)$ is $A$\injective{} and the morphism
$\bP\rarrow\bgI^+(\bP)$ is injective.
 Set $\bgI^0(\bgK^\bu)=\bgI^+(\bgK^\bu)$, \
$\bgI^1(\bgK^\bu)=\bgI^+(\coker(\bgK^\bu\to\bgI^0(\bgK^\bu)))$, etc.
 For $d$ large enough, the cokernel $\bgZ(\bgK^\bu)$ of the morphism
$\bgI^{d-2}(\bgK^\bu)\rarrow\bgI^{d-1}(\bgK^\bu)$ will be a complex
of $A$\injective{} $\S$\semicontramodule s.
 Let $\boR_1(\bgK^\bu)$ be the total complex of the bicomplex
$$
 \bgI^0(\bgK^\bu)\lrarrow\bgI^1(\bgK^\bu)\lrarrow\dsb\lrarrow
 \bgI^{d-1}(\bgK^\bu)\lrarrow\bgZ(\bgK^\bu).
$$
 Then $\boR_1(\bgK^\bu)$ is a complex of $A$\injective{}
$\S$\semicontramodule s and the cone of the morphism
$\bgK^\bu\rarrow\boR_1(\bgK^\bu)$ is the total complex of
a finite exact sequence of complexes of $\S$\semicontramodule s,
and therefore, a $\C$\contraacyclic{} complex
(and even an $\S$\contraacyclic{} complex).

 Now let $\bgR^\bu$ be a complex of $A$\injective{} left
$\S$\semicontramodule s.
 Let $\bgF(\bP)\rarrow\bP$ denote the functorial surjective morphism
onto an arbitrary $A$\injective{} $\S$\semicontramodule{} $\bP$ from
a $\C$\coinjective{} $\S$\semicontramodule{} $\bgF(\bP)$ with
an $A$\injective{} kernel $\ker(\bgF(\bP)\to\bP)$ constructed in 
Lemma~\ref{coproj-coinj-semi-mod-contra}(b).
 Set $\bgF_0(\bgR^\bu)=\bgF(\bgR^\bu)$, \ 
$\bgF_1(\bgR^\bu)=\bgF(\ker(\bgF_0(\bgR^\bu)\to\bgR^\bu))$, etc.
 Let $\boL_2(\bgR^\bu)$ be the total complex of the bicomplex
$$
 \dsb\lrarrow\bgF_2(\bgR^\bu)\lrarrow\bgF_1(\bgR^\bu)
 \lrarrow\bgF_0(\bgR^\bu),
$$
constructed by taking infinite products along the diagonals.
 Then $\boL_2(\bgR^\bu)$ is a complex of $\C$\coinjective{}
$\S$\semicontramodule s.
 Since the surjection $\bgF(\bP)\rarrow\bP$ can be defined for
arbitrary left $\S$\semicontramodule s, the functor $\boL_2$ can
be extended to arbitrary complexes of $\S$\semicontramodule s.
 For any complex $\bgK^\bu$, the cone of the morphism
$\boL_2(\bgK^\bu)\rarrow\bgK^\bu$ is a $\C$\contraacyclic{} complex
(and even an $\S$\contraacyclic{} complex) by
Lemma~\ref{contraderived-categories}.

 Finally, let $\bP^\bu$ be a $\C$\coinjective{} complex of
$A$\injective{} left $\S$\semicontramodule s.
 Then the complex $\Cohom_\C(\S,\bP^\bu)$ is a semiinjective complex
of $A$\injective{} left $\S$\semicontramodule s.
 Moreover, if $\bP^\bu$ is a complex of $\C$\coinjective{}
$\S$\semicontramodule s, then $\Cohom_\C(\S,\bP^\bu)$ is
a semiinjective complex of semiinjective $\S$\semicontramodule s.
 Consider the cobar construction
$$
 \Cohom_\C(\S,\bP^\bu)\lrarrow\Cohom_\C(\S,\Cohom_\C(\S,\bP^\bu))
 \lrarrow\dsb
$$
 Let $\boR_3(\bP^\bu)$ be the total complex of this bicomplex,
constructed by taking infinite products along the diagonals.
 Then complex $\boR_3(\bP^\bu)$ is semiinjective
by Lemma~\ref{semi-inj-proj-complexes}.
 The functor $\boR_3$ can be extended to arbitrary complexes of
$\S$\semicontramodule s; for any complex $\bgK^\bu$, the cone of
the morphism $\bgK^\bu\rarrow\boR_3(\bgK^\bu)$ is not only
$\C$\contraacyclic, but even $\C$\contractible{} (the contracting
homotopy being induced by the semiunit morphism $\C\rarrow\S$.)

 It follows that the natural functors between the quotient categories
of the homotopy categories of semiinjective complexes of semiinjective
$\S$\semicontramodule s, semiinjective complexes of $\C$\coinjective{}
$\S$\semicontramodule s, complexes of $\C$\coinjective{}
$\S$\semicontramodule s, semiinjective complexes of $A$\injective{}
$\S$\semicontramodule s, $\C$\coinjective{} complexes of
$A$\injective{} $\S$\semicontramodule s, complexes of $A$\injective{}
$\S$\semicontramodule s by their intersections with the thick
subcategory of $\C$\contraacyclic complexes and the semiderived
category of left $\S$\semicontramodule s are all equivalences of
triangulated categories.
 Moreover, any complex of left $\S$\semicontramodule s $\bgK^\bu$
can be connected with a semiinjective complex of semiinjective
$\S$\semicontramodule s in a functorial way by a chain of three
morphisms $\bgK^\bu\rarrow\boR_3(\bgK^\bu)\larrow\boR_3\boL_2(\bgK^\bu)
\rarrow\boR_3\boL_2\boR_1(\bgK^\bu)$ with $\C$\contraacyclic{} cones,
and when $\bgK^\bu$ is a semiinjective complex of ($A$\injective,
$\C$\coinjective, or semiinjective) $\S$\semicontramodule s, all
complexes in this chain are also semiinjective complexes of
($A$\injective, $\C$\coinjective, or semiinjective)
$\S$\semicontramodule s.
\end{proof}

\begin{rmk}
 One can show using the methods developed in
Section~\ref{semi-correspondence-section} that any
$\C$\coacyclic{} semiprojective complex of $\C$\coprojective{} left
$\S$\semimodule s is contractible, and analogously, any
$\C$\contraacyclic{} semiinjective complex of $\C$\coinjective{}
left $\S$\semicontramodule s is contractible
(see Remark~\ref{birelatively-adjusted}).
\end{rmk}

\subsection{Derived functor SemiExt}
\label{semiext-definition}
 Assume that the coring $\C$ is a projective left and a flat right
$A$\module, the semialgebra $\S$ is a coprojective left and
a coflat right $\C$\comodule, and the ring $A$ has a finite left
homological dimension.

 The double-sided derived functor
$$
 \SemiExt_\S\:\sD^\si(\S\simodl)\times\sD^\si(\S\sicntr)
 \lrarrow\sD(k\modl)
$$
is defined as follows.
 Consider the partially defined functor of semihomomorphisms of
complexes $\SemiHom_\S\:\Hot(\S\simodl)^\op\times\Hot(\S\sicntr)
\darrow \Hot(k\modl)$.
 This functor is defined on the full subcategory of the Carthesian
product of homotopy categories that consists of pairs of complexes
$(\M^\bu,\P^\bu)$ such that either $\M^\bu$ is a complex of
$A$\projective{} $\S$\semimodule s, or $\P^\bu$ is a complex of 
$A$\injective{} $\S$\semicontramodule s.
 Compose it with the functor of localization $\Hot(k\modl)\rarrow
\sD(k\modl)$ and restrict either to the Carthesian product of
the homotopy category of semiprojective complexes of $A$\projective{}
$\S$\semimodule s and the homotopy category of $\S$\semicontramodule s,
or to the Carthesian product of the homotopy category of
$\S$\semimodule s and the homotopy category of semiinjective complexes
of $A$\injective{} $\S$\semicontramodule s.

 By Theorem~\ref{semiext-main-theorem} and
Lemma~\ref{semitor-definition}, both functors so obtained
factorize through the Carthesian product of semiderived categories
of left semimodules and left semicontramodules and the derived
functors so defined are naturally isomorphic.
 The same derived functor is obtained by restricting the functor
of semihomomorphisms to the Carthesian product of the homotopy
categories of semiprojective complexes of $A$\projective{}
$\S$\semimodule s and semiinjective complexes of $A$\injective{}
$\S$\semicontramodule s.
 One can also use semiprojective complexes of $\C$\coprojective{}
$\S$\semimodule s or semiinjective complexes of $\C$\coinjective{}
$\S$\semicontramodule s, etc.

{\hbadness=3000
 In particular, when the coring $\C$ is a projective left and a flat
right $A$\module{} and the ring $A$ has a finite left homological
dimension, one defines the double-sided derived functor
$$
 \Coext_\C\:\sD^\co(\C\comodl)^\op\times\sD^\ctr(\C\contra)
 \lrarrow\sD(k\modl)
$$
by composing the functor of cohomomorphisms $\Cohom_\C\:
\Hot(\C\comodl)^\op\times\Hot(\C\contra)\rarrow \Hot(k\modl)$
with the functor of localization $\Hot(k\modl)\rarrow
\sD(k\modl)$ and restricting it to either the Carthesian product
of the homotopy category of complexes of coprojective
$\C$\comodule s and the homotopy category of arbitrary complexes
of $\C$\contramodule s, or the Carthesian product of the homotopy
category of arbitrary complexes of $\C$\comodule s and the homotopy
category of complexes of coinjective $\C$\contramodule s.
 The same derived functor is obtained by restricting the functor
of cohomomorphisms to the Carthesian product of the homotopy
categories of coprojective $\C$\comodule s and coinjective
$\C$\contramodule s.
 One can also use coprojective complexes of $\C$\comodule s or
coinjective complexes of $\C$\contramodule s. \par}

\begin{qst}
 Assuming only that $\C$ is a flat left and right $A$\module, one
can define the double-sided derived functor $\Cotor^\C$ on the Carthesian
product of coderived categories of the exact categories of right and
left $\C$\comodule s of flat dimension over~$A$ not exceeding~$d$,
for any given~$d$, using Lemma~\ref{semitor-definition} and
the corresponding version of Lemma~\ref{flat-comodule-surjection}.
 Analogously, assuming that $\C$ is a projective left and a flat right
$A$\module, one can define the double-sided derived functor $\Coext_\C$
on the Carthesian product of the coderived category of left
$\C$\comodule s of projective dimension over~$A$ not exceeding~$d$ and
the contraderived category of $\C$\contramodule s of injective
dimension over~$A$ not exceeding~$d$.
 One can even do with the homological dimension assumption on only one
of the arguments of $\Cotor^\C$ and $\Coext_\C$, using the corresponding
versions of the results of Theorem~\ref{co-contra-push-well-defined}.
 Can one define, at least, a derived functor $\SemiTor^\S$ for complexes
of $A$\+flat $\S$\semimodule s and a derived functor $\SemiExt_\S$ for 
complexes of $A$\projective{} $\S$\semimodule s and $A$\injective{}
$\S$\semicontramodule s without the homological dimension assumptions
on~$A$?
 The only problem one encounters attempting to do so comes from
the homological dimension conditions in
Propositions~\ref{tensor-cotensor-assoc}(c) and
\ref{co-tensor-co-hom-assoc}.1\+2(c) and consequently in
Lemmas~\ref{coflat-semimodule-injection}
and~\ref{coproj-coinj-semi-mod-contra}; when $\S$ satisfies
the conditions of Proposition~\ref{cotensor-associative}(f)
there is no problem.
\end{qst}

\begin{rmk}
 In the way completely analogous to Remark~\ref{semitor-definition},
without any homological dimension assumptions one can define
the double-sided derived functor $\IndCoext_\bC$ for complexes of left
$\bC$\comodule s in $k\modl^\omega$ and complexes of left
$\bC$\contramodule s in the category $k\modl_\omega$ of ind-objects
over $k\modl$ representable by countable filtered inductive systems of
$k$\module s.
 Here the category opposite to $k\modl_\omega$ is considered as
a module category over the tensor category $k\modl^\omega$ and $\bC$
is a coring over a ring $\cA$ in $k\modl^\omega$.
 Appropriate coflatness and ``contraprojectivity'' conditions have to
be imposed on~$\bC$.
 The countability assumption can be dropped.
\end{rmk}

\subsection{Relatively semiprojective and semiinjective complexes}
\label{relatively-semiproj-semiinj}
 We keep the assumptions and notation of~\ref{coext-main-theorem},
\ref{semiext-main-theorem}, and~\ref{semiext-definition}.

 One can compute the derived functor $\Coext_\C$ using resolutions of
other kinds.
 Namely, the complex of cohomomorphisms $\Cohom_\C(\M^\bu,\P^\bu)$ from
a complex of $\C/A$\coprojective{} left $\C$\comodule s $\M^\bu$ into
a complex of $A$\injective{} left $\C$\contramodule s $\P^\bu$
represents an object naturally isomorphic to $\Coext_\C(\M^\bu,\P^\bu)$
in the derived category of $k$\module s.
 Indeed, the complex $\boL_2(\P^\bu)$ is a complex of coinjective
$\C$\contramodule s and the cone of the morphism $\boL_2(\P^\bu)\rarrow
\P^\bu$ is contraacyclic with respect to the exact category of
$A$\injective{} $\C$\contramodule s, hence the morphism
$\Cohom_\C(\M^\bu,\boL_2(\P^\bu))\rarrow\Cohom_\C(\M^\bu,\P^\bu)$
is an isomorphism.
 Analogously, the complex of cohomomorphisms $\Cohom_\C(\M^\bu,\P^\bu)$
from a complex of $A$\projective{} left $\C$\comodule s $\M^\bu$ into
a complex of $\C/A$\coinjective{} left $\C$\contramodule s $\P^\bu$
represents an object naturally isomorphic to $\Coext_\C(\M^\bu,\P^\bu)$
in the derived category of $k$\module s. 

 One can also compute the derived functor $\SemiExt_\C$ using
resolutions of other kinds.
 Namely, a complex of left $\S$\semimodule s is called
\emph{semiprojective relative to~$A$} if the complex of
semihomomorphisms from it into any complex of $A$\injective{} left
$\S$\semicontramodule s that as a complex of $\C$\contramodule s is
contraacyclic with respect to the exact category of $A$\injective{}
$\C$\contramodule s is acyclic
(cf.\ Theorem~\ref{co-contra-push-well-defined}(c)).
 The complex of semihomomorphisms $\SemiHom_\C(\bM^\bu,\bP^\bu)$
from a complex of left $\S$\semimodule s $\bM^\bu$ semiprojective
relative to~$A$ into a complex of $A$\injective{} left
$\S$\semicontramodule s $\bP^\bu$ represents an object naturally
isomorphic to $\SemiExt_\S(\bM^\bu,\bP^\bu)$ in the derived category
of $k$\module s.
 Indeed, $\boR_3\boL_2(\bP^\bu)$ is a semiinjective complex of
$\S$\semicontramodule s connected with $\bP^\bu$ by a chain of
morphisms $\bP^\bu\larrow\boL_2(\bP^\bu)\rarrow\boR_3\boL_2(\bP^\bu)$
whose cones are contraacyclic with respect to the exact category of
$A$\injective{} $\C$\contramodule s and contractible over~$\C$,
respectively.
 Analogously, a complex of left $\S$\semicontramodule s is called
\emph{semiinjective relative to~$A$} if the complex of
semihomomorphisms into it from any complex of $A$\projective{} left
$\S$\semimodule s that as a complex of $\C$\comodule s is coacyclic
with respect to the exact category of $A$\projective{} $\C$\comodule s
is acyclic (cf.\ Theorem~\ref{co-contra-push-well-defined}(b)).
 The complex of semihomomorphisms $\SemiHom_\C(\bM^\bu,\bP^\bu)$
from a complex of $A$\projective{} left $\S$\semimodule s to a complex
of left $\S$\semicontramodule s semiinjective relative to~$A$
represents an object naturally isomorphic to $\SemiExt_\S(\bM^\bu,
\bP^\bu)$ in the derived category of $k$\module s.
 For example, the complex of $\S$\semimodule s induced from a complex
of $\C/A$\coprojective{} left $\C$\comodule s is semiprojective
relative to~$A$ and the complex of $\S$\semicontramodule s coinduced
from a complex of $\C/A$\coinjective{} left $\C$\contramodule s
is semiinjective relative to~$A$.

 A complex of left $\S$\semimodule s is called \emph{semiprojective
relative to\/~$\C$} if the complex of semihomomorphisms from it into
any $\C$\contractible{} complex of $\C$\coinjective{} left
$\S$\semicontramodule s is acyclic.
 The complex of semihomomorphisms $\SemiHom_\C(\bM^\bu,\bP^\bu)$
from a complex of left $\S$\semimodule s $\bM^\bu$ semiprojective
relative to~$\C$ into a complex of $\C$\coinjective{} left
$\S$\semicontramodule s $\bP^\bu$ represents an object naturally
isomorphic to $\SemiExt_\S(\bM^\bu,\bP^\bu)$ in the derived category
of $k$\module s.
 Indeed, $\boR_3(\bP^\bu)$ is a semiinjective complex of
$\S$\semicontramodule s and the cone of the morphism
$\bP^\bu\rarrow\boR_3(\bP^\bu)$ is a $\C$\contractible{} complex
of $\C$\coinjective{} $\S$\semicontramodule s.
 Analogously, a complex of left $\S$\semicontramodule s is called
\emph{semiinjective relative to\/~$\C$} if the complex of
semihomomorphisms into it from any $\C$\contractible{} complex of
$\C$\coprojective{} left $\S$\semimodule s is acyclic.
 The complex of semihomomorphisms $\SemiHom_\C(\bM^\bu,\bP^\bu)$
from a complex of $\C$\coprojective{} left $\S$\semimodule s $\bM^\bu$
into a complex of left $\S$\semicontramodule s $\bP^\bu$ semiinjective
relative to~$\C$ represents an object naturally
isomorphic to $\SemiExt_\S(\bM^\bu,\bP^\bu)$ in the derived category
of $k$\module s.
 It follows that the complex of semihomomorphisms from a complex of
left $\S$\semimodule s semiprojective relative to~$\C$ into
a $\C$\contraacyclic{} complex of $\C$\coinjective{} left
$\S$\semicontramodule s is acyclic, and the complex of
semihomomorphisms into a complex of left $\S$\semicontramodule s
semiinjective relative to~$\C$  from a $\C$\coacyclic{} complex of
$\C$\coprojective{} left $\S$\semimodule s is acyclic.
 For example, the complex of $\S$\semimodule s induced from a complex
of left $\C$\comodule s is semiprojective relative to~$\C$ and
the complex of $\S$\semicontramodule s coinduced from a complex of
left $\C$\contramodule s is semiinjective relative to~$\C$.

 At last, a complex of $A$\projective{} left $\S$\semimodule s
is called \emph{semiprojective relative to\/ $\C$ relative to~$A$}
($\S/\C/A$\semiprojective) if the complex of semihomomorphisms from it
into any $\C$\contractible{} complex of $\C/A$\coinjective{}
left $\S$\semicontramodule s is acyclic.
 The complex of semihomomorphisms $\SemiHom_\C(\bM^\bu,\bP^\bu)$
from an $\S/\C/A$\semiprojective{} complex of $A$\projective{} left
$\S$\semimodule s $\bM^\bu$ into a complex of $\C/A$\coinjective{} left
$\S$\semicontramodule s $\bP^\bu$ represents an object naturally
isomorphic to $\SemiExt_\S(\bM^\bu,\bP^\bu)$ in the derived category
of $k$\module s.
 Indeed, $\boR_3(\bP^\bu)$ is a complex of left $\S$\semicontramodule s
semiinjective relative to~$A$ and the cone of the morphism $\bP^\bu
\rarrow\boR_3(\bP^\bu)$ is a $\C$\contractible{} complex of
$\C/A$\coinjective{} $\S$\semicontramodule s.
 Analogously, a complex of $A$\injective{} left $\S$\semicontramodule s
is called \emph{semiinjective relative to\/ $\C$ relative to~$A$}
($\S/\C/A$\semiinjective) if the complex of semihomomorphisms into it
from any $\C$\contractible{} complex of $\C/A$\coprojective{}
left $\S$\semimodule s is acyclic.
 The complex of semihomomorphisms $\SemiHom_\C(\bM^\bu,\bP^\bu)$
from a complex of $\C/A$\coprojective{} left $\S$\semimodule s
$\bM^\bu$ into an $\S/\C/A$\semiinjective{} complex of $A$\injective{}
left $\S$\semicontramodule s $\bP^\bu$ represents an object naturally
isomorphic to $\SemiExt_\S(\bM^\bu,\bP^\bu)$ in the derived category
of $k$\module s.
 It follows that the complex of semihomomorphisms from an
$\S/\C/A$\semiprojective{} complex of $A$\projective{} left
$\S$\semimodule s into a $\C$\contraacyclic{} complex of
$\C/A$\coinjective{} left $\S$\semicontramodule s is acyclic, and
the complex of semihomomorphisms into an $\S/\C/A$\semiinjective{}
complex of $A$\injective{} left $\S$\semicontramodule s from a
$\C$\coacyclic{} complex of $\C/A$\coprojective{} left $\S$\semimodule s
is acyclic.
 For example, the complex of $\S$\semimodule s induced from a complex
of $A$\projective{} left $\C$\comodule s is $\S/\C/A$\semiprojective{}
and the complex of $\S$\semicontramodule s coinduced from a complex of
$A$\injective{} left $\C$\contramodule s is $\S/\C/A$\semiinjective.

 The functors mapping the quotient categories of the homotopy
categories of complexes of $\S$\semimodule s semiprojective relative
to~$A$, complexes of $\S$\semimodule s semiprojective relative to~$\C$,
and $\S/\C/A$\semiprojective{} complexes of $\S$\semimodule s by their
intersections with the thick subcategory of $\C$\coacyclic{} complexes
into the semiderived category of left $\S$\semimodule s are equivalences
of triangulated categories.
 Analogously, the functors mapping the quotient categories of
the homotopy categories of complexes of $\S$\semicontramodule s
semiinjective relative to~$A$, complexes of $\S$\semicontramodule s
semiinjective relative to~$\C$, and $\S/\C/A$\semiinjective{} complexes
of $\S$\semicontramodule s by their intersections with the thick
subcategory of $\C$\contraacyclic{} complexes into the semiderived
category of left $\S$\semicontramodule s are equivalences of triangulated
categories.
 The same applies to complexes of $A$\projective, $\C$\coprojective,
and $\C/A$\coprojective{} $\S$\semimodule s and complexes of
$A$\injective, $\C$\coinjective, and $\C/A$\coinjective{}
$\S$\semicontramodule s.
 These results follow easily from either of
Lemmas~\ref{cotor-main-theorem} or~\ref{semitor-main-theorem}.
 So one can define the derived functor $\SemiExt_\S$ by restricting
the functor of semihomomorphisms to these categories of complexes
as explained above.

\begin{rmk}
 One can define the double-sided or right derived functor
$\SemiExt_\S$ in the assumptions analogous to those of
Remark~\ref{relatively-semiflat} in the completely analogous ways.
\end{rmk}

\subsection{Remarks on derived semihomomorphisms from bisemimodules}
\label{remarks-derived-semihom-bi}
 Let $\S$ be a semialgebra over a coring~$\C$ and $\T$ be a semialgebra
over a coring~$\D$, both satisfying the conditions
of~\ref{semiext-main-theorem}.
 One can define the double-sided derived functor 
 $$
  \boD\SemiHom_\S\:\sD^\si(\S\bsimod\T)^\op\times\sD^\si(\S\sicntr)
  \lrarrow\sD^\si(\T\sicntr)
 $$
by restricting the functor of semihomomorphisms $\SemiHom_\S\:
\sD^\si(\S\bsimod\T)^\op\times\sD^\si(\S\sicntr)\darrow
\sD^\si(\T\sicntr)$ to the Carthesian product of the homotopy
category of complexes of $\S$\+$\T$\bisemimodule s and the homotopy
category of semiinjective complexes of $\C$\coinjective{} left
$\S$\semicontramodule s (using the result of
Remark~\ref{birelatively-adjusted}).
 There is an associativity isomorphism
$$
 \SemiExt_\S(\bK^\bu\os_\T^\boD\bM^\bu\;\bP^\bu)\simeq
 \SemiExt_\T(\bM^\bu,\boD\SemiHom_\S(\bK^\bu,\bP^\bu)).
$$

 Let $\bR$ be a semialgebra over a coring~$\E$  satisfying
the conditions of~\ref{semiext-main-theorem}. 
 If the $k$\+algebra $A$ is a flat $k$\module{} and the $k$\+algebras
$B$ and $F$ are projective $k$\module s, then the derived functor
$\boD\SemiHom$ can be defined using Lemma~\ref{semitor-definition}
in terms of \emph{strongly\/ $\S$\semiprojective} complexes of
$A$\projective{} $\S$\+$\T$\bisemimodule s and semiinjective complexes
of $\C$\coinjective{} left $\S$\semicontramodule s (or \emph{strongly
semiinjective} complexes of $A$\injective{} left
$\S$\semicontramodule s).
 Here a complex of $A$\projective{} $\S$\+$\T$\bisemimodule s
$\bK^\bu$ is called strongly $\S$\semiprojective{} if for any
$\C$\contraacyclic{} complex of left $\S$\semicontramodule s
$\bP^\bu$ the complex of left $\T$\semicontramodule s
$\SemiHom_\S(\bK^\bu,\bP^\bu)$ is $\D$\contraacyclic; strongly
semiinjective complexes are defined in the analogous way.
 In this case, there is an associativity isomorphism
$$
 \boD\SemiHom_\S(\bK^\bu\os_\T^\boD\bM^\bu\;\bP^\bu)\simeq
 \boD\SemiHom_\T(\bM^\bu,\boD\SemiHom_\S(\bK^\bu,\bP^\bu))
$$
for any complex of $\T$\+$\bR$\bisemimodule s $\bM^\bu$, any complex
of $\S$\+$\T$\bisemimodule s $\bK^\bu$, and any complex of left
$\S$\semicontramodule s $\bP^\bu$.

 In particular, without any conditions on the $k$\module{} $A$ for
any complex of right $\S$\semimodule s $\bN^\bu$ and any complex of
left $\S$\semimodule s $\bM^\bu$ there is a natural isomorphism
$\Hom_k(\SemiTor^\S(\bN^\bu,\bM^\bu),\.k\dual)\simeq
\SemiExt_\S(\bM^\bu,\Hom_k(\bN^\bu,k\dual))$.

\Section{Comodule-Contramodule Correspondence}
\label{co-contra-correspondence-section}

\subsection{Contratensor product and comodule/contramodule
homomorphisms}
 Let $\C$ be a coring over a $k$\+algebra~$A$.

\subsubsection{}
 The \emph{contratensor product} $\N\ocn_\C\P$ of a right
$\C$\comodule{} $\N$ and a left $\C$\contramodule{} $\P$ is
a $k$\module{} defined as the cokernel of the pair of maps
$\N\ot_A\Hom_A(\C,\P)\birarrow\N\ot_A\P$ one of which is induced by
the $\C$\+contraaction in $\P$, while the other is the composition
of the map induced by the $\C$\+coaction in~$\N$ and the map induced
by the evaluation map $\C\ot_A\Hom_A(\C,\P)\rarrow\P$.

 The contratensor product operation is dual to homomorphisms
in the category of contramodules: for any right $\C$\comodule{} $\N$
with a left action of a $k$\+algebra $B$ by $\C$\comodule{}
endomorphisms, any left $\C$\contramodule{} $\P$, and any left
$B$\module{} $U$ there is a natural isomorphism
$\Hom_B(\N\ocn_\C\P\;U)\simeq\Hom^\C(\P,\Hom_B(\N,U))$.
 Indeed, both $k$\module s are isomorphic to the kernel of
the same pair of maps $\Hom_A(\P,\Hom_B(\N,U))\birarrow
\Hom_A(\Hom_A(\C,\P),\Hom_B(\N,U))$.
 Taking $B=k$, one can conclude that for any right $\C$\comodule{} $\N$
and any left $A$\module{} $V$ there is a natural isomorphism
$\N\ocn_\C\Hom_A(\C,V)\simeq\N\ot_A V$.

 When $\C$ is a projective left $A$\module, the functor of
contratensor product over $\C$ is right exact in both its
arguments.

\subsubsection{}   \label{psi-phi-adjoint}
 Let $\D$ be a coring over a $k$\+algebra~$B$.
 For any $\C$\+$\D$\bicomodule{} $\K$ and any left $\C$\comodule{}
$\M$, the $k$\module{} $\Hom_\C(\K,\M)$ has a natural left
$\D$\contramodule{} structure as the kernel of a pair of
$\D$\contramodule{} morphisms $\Hom_A(\K,\M)\birarrow
\Hom_A(\K\;\M\ot_B\D)$.
 Analogously, for any $\D$\+$\C$\bicomodule{} $\K$ and any left
$\C$\contramodule{} $\P$, the $k$\module{} $\K\ocn_\C\P$ has
a natural left $\D$\comodule{} structure as the cokernel of a pair
of $\D$\comodule{} morphisms $\K\ot_A\Hom_A(\C,\P)\birarrow\K\ot_A\P$.

 For any left $\D$\comodule{} $\M$, any $\D$\+$\C$\bicomodule{} $\K$,
and any left $\C$\contramodule{} $\P$ there is a natural isomorphism
$\Hom_\D(\K\ocn_\C\P\;\M)\simeq \Hom^\C(\P,\Hom_\D(\K,\M))$.
 Indeed, a $B$\module{} map $\K\ot_A\P\rarrow\M$ factorizes through
$\K\ocn_\C\M$ if and only if the corresponding $A$\module{} map
$\P\rarrow\Hom_B(\K,\M)$ is a $\C$\contramodule{} morphism, and
a $B$\module{} map $\K\ot_A\P\rarrow\M$ is a $\D$\comodule{} morphism
if and only if the corresponding $A$\module{} map
$\P\rarrow\Hom_B(\K,\M)$ factorizes through $\Hom_\D(\K,\M)$.

 In particular, there is a pair of adjoint functors
$\Psi_\C\:\C\comodl\rarrow\C\contra$ and $\Phi_\C\:\C\contra\rarrow
\C\comodl$ between the categories of left $\C$\comodule s and left
$\C$\contramodule s defined by the rules $\Psi_\C(\M)=\Hom_\C(\C,\M)$
and $\Phi_\C(\P)=\C\ocn_\C\P$.

\subsubsection{}   \label{quite-injective-projective}
 A left $\C$\comodule{} $\M$ is called \emph{quite injective relative
to~$A$} (quite $\C/A$\+injective) if the functor of $\C$\comodule{}
homomorphisms into~$\M$ maps $A$\+split exact triples of left
$\C$\comodule s to exact triples.
 It is easy to see that a $\C$\comodule{} is quite $\C/A$\injective{}
if and only if it is a direct summand of a coinduced $\C$\comodule.
 Analogously, a left $\C$\contramodule{} $\P$ is called \emph{quite
projective relative to~$A$} (quite $\C/A$\+projective) if the functor
of $\C$\contramodule{} homomorphisms from~$\P$ maps $A$\+split
exact triples of left $\C$\contramodule s to exact triples.
 A $\C$\contramodule{} is quite $\C/A$\projective{} if and only if
it is a direct summand of an induced $\C$\contramodule.

 The restrictions of the functors $\Psi_\C$ and $\Phi_\C$ on
the subcategories of quite $\C/A$\injective{} left $\C$\comodule s
and quite $\C/A$\projective{} left $\C$\contramodule s are mutually
inverse equivalences between these subcategories.
 Indeed, one has $\Hom_\C(\C\;\C\ot_A V) = \Hom_A(\C,V)$
and $\C\ocn\Hom_A(\C,V) = \C\ot_A V$.

\subsubsection{}
 A left $\C$\comodule{} $\M$ is called \emph{injective relative to~$A$}
($\C/A$\+injective) if the functor of homomorphisms into~$\M$ maps
exact triples of $A$\projective{} left $\C$\comodule s to exact triples.
 A left $\C$\contramodule{} $\P$ is called \emph{projective relative
to~$A$} ($\C/A$\+projective) if the functor of homomorphisms from~$\P$
maps exact triples of $A$\injective{} left $\C$\contramodule s
to exact triples.
 (Cf.\ Lemma~\ref{rel-inj-proj-co-contra-mod}.2.)
 
\begin{rmk}
 What we call quite relatively injective comodules are usually
called relatively injective comodules~\cite{BW}.
 We chose this nontraditional terminology for coherence with our
definitions of relative coflatness, etc., and also because what we call
relatively injective comodules is a more important notion from our
point of view.
\end{rmk}

\begin{qst}
 One can compute modules Ext in the exact category of left
$\C$\comodule s with $A$\+split exact triples in terms of
the cobar resolution.
 When $\C$ is a projective left $A$\module, this resolution can be
also used to compute modules Ext in the exact category of
$A$\projective{} left $\C$\comodule s, which therefore turn out
to be the same.
 How can one compute modules Ext in the exact category of
$A$\projective{} $\C$\comodule s without making any projectivity
assumptions on~$\C$?
\end{qst}

\subsubsection{}    \label{injective-projective-co-contra-mod}
 When $\C$ is a flat right $A$\module, the coinduction functor
$A\modl\rarrow\C\comodl$ preserves injective objects.
 It follows easily that any left $\C$\comodule{} is a subcomodule of
an injective $\C$\comodule; a $\C$\comodule{} is injective if
and only if it is a direct summand of a $\C$\comodule{} coinduced from
an injective $A$\module.
 Analogously, when $\C$ is a projective left $A$\module, the induction
functor $A\modl\rarrow\C\contra$ preserves projective objects.
 Hence any left $\C$\contramodule{} is a quotient contramodule of
a projective $\C$\contramodule; a $\C$\contramodule{} is
projective if and only if it is a direct summand of
a $\C$\contramodule{} induced from a projective $A$\module.

\subsubsection{}
 When $\C$ is a flat left $A$\module, a left $\C$\contramodule{}
$\P$ is called \emph{contraflat} if the functor of contratensor
product with $\P$ is exact on the category of right $\C$\comodule s.
 The $\C$\contramodule{} induced from a flat $A$\module{} is
contraflat.
 Any projective $\C$\contramodule{} is contraflat.

 A left $\C$\contramodule{} $\P$ is called \emph{quite\/
$\C/A$\+contraflat} if the functor of contratensor product with~$\P$
maps those exact triples of right $\C$\comodule s which as exact
triples of $A$\module s remain exact after the tensor product with
any left $A$\module{} to exact triples.
 Any quite $\C/A$\projective{} $\C$\contramodule{} is quite
$\C/A$\contraflat.
 A left $\C$\contramodule{} $\P$ is called \emph{$\C/A$\+contraflat}
if the functor of contratensor product with~$\P$ maps exact triples
of $A$\+flat right $\C$\comodule s to exact triples.
 Using the dualization functor $\Hom_k({-},k\dual)$, one can easily
check that any $\C/A$\projective{} $\C$\comodule{} is $\C/A$\contraflat.

\subsection{Associativity isomorphisms}
\label{cotensor-contratensor-assoc}
 Let $\C$ be a coring over a $k$\+algebra $A$ and $\D$ be a coring
over a $k$\+algebra $B$.
 The following three Propositions will be mostly applied to
the case of $\K=\D=\C$ in the sequel.

\begin{prop1}
 Let\/ $\N$ be a right\/ $\D$\comodule, $\K$ be a\/
$\D$\+$\C$\bicomodule, and\/ $\P$ be a left\/ $\C$\contramodule.
 Then there is a natural map\/ $(\N\oc_\D\K)\ocn_\C\P\rarrow
\N\oc_\D(\K\ocn_\C\P)$ whenever the cotensor product\/ $\N\oc_\D\K$
is endowed with a right\/ $\C$\comodule{} structure such that the map\/
$\N\oc_\D\K\rarrow\N\ot_B\K$ is a\/ $\C$\comodule{} morphism.
 This natural map is an isomorphism, at least, in the following cases:
\begin{enumerate}
 \item $\C$ is a flat left $A$\module{} and\/ $\P$ is a contraflat{}
       left\/ $\C$\contramodule,
 \item $\P$ is a quite\/ $\C/A$\contraflat{} left\/ $\C$\contramodule{}
       and\/ $\K$ as a left\/ $\D$\comodule{} with a right $A$\module{}
       structure is conduced from a $B$\+$A$\bimodule;
 \item $\P$ is a\/ $\C/A$\contraflat{} left\/ $\C$\contramodule, $\D$
       is a flat right $B$\module, $\N$ is a flat right $B$\module,
       and\/ $\K$ as a left\/ $\D$\comodule{} with a right $A$\module{}
       structure is coinduced from an $A$\+flat $B$\+$A$\bimodule;
 \item $\P$ is a\/ $\C/A$\contraflat{} left\/ $\C$\contramodule, $\D$
       is a flat right $B$\module, $\N$ is a flat right $B$\module,
       $\K$ is a flat right $A$\module, $\K$ is a\/ $\D/B$\+coflat
       left\/ $\D$\comodule, and the ring $A$ has a finite weak
       homological dimension;
 \item $\N$ is a quasicoflat right\/ $\D$\comodule.
\end{enumerate} 
\end{prop1}

\begin{proof}
 The map $(\N\oc_\D\K)\ocn_\C\P\rarrow\N\ot_B\K\ocn_\C\P$ has equal
compositions with two maps $\N\ot_B\K\ocn_\C\P\birarrow
\N\ot_B\D\ot_B\K\ocn_\C\P$, so there is a natural map
$(\N\oc_\D\K)\ocn_\C\P\rarrow\N\oc_\D(\K\ocn_\C\P)$.
 Besides, the composition $(\N\oc_\D\K)\ot_A\P\rarrow
\N\oc_\D(\K\ot_A\P)\rarrow\N\oc_\D(\K\ocn_\C\P)$ annihilates
the difference between two maps $(\N\oc_\D\K)\ot_A\Hom_A(\C,\P)
\birarrow(\N\oc_\D\K)\ot_A\P$, which leads to the same natural map
$(\N\oc_\D\K)\ocn_\C\P\rarrow\N\oc_\D(\K\ocn_\C\P)$.
 To prove cases~(a-d), one shows that the sequence
$0\rarrow\N\oc_\D\K\rarrow\N\ot_B\K\rarrow\N\ot_B\D\ot_B\K$
remains exact after taking the contratensor product with~$\P$.
 Indeed, the case~(a) is obvious, in the cases~(b-c) this exact
sequence of right $A$\module s splits, and in the cases~(c-d)
this sequence of right $A$\module s is exact with respect to
the exact category of flat $A$\module s (see the proof of
Proposition~\ref{tensor-cotensor-assoc}).
 To prove~(e), one notices that the sequence
$\K\ot_A\Hom_A(\C,\P)\rarrow\K\ot_A\P\rarrow\K\ocn_\C\P\rarrow0$
remains exact after taking the cotensor product with~$\N$ and uses
Proposition~\ref{tensor-cotensor-assoc}(b).
\end{proof}

\begin{prop2}
 Let\/ $\L$ be a left\/ $\D$\comodule, $\K$ be a\/
$\C$\+$\D$\bicomodule, and\/ $\M$ be a left\/ $\C$\comodule.
 Then there is a natural map\/ $\Cohom_\D(\L,\Hom_\C(\K,\M))\rarrow
\Hom_\C(\K\oc_\D\L\;\M)$ whenever the cotensor product\/ $\K\oc_\D\L$
is endowed with a left\/ $\C$\comodule{} structure such that the map\/
$\K\oc_\D\L\rarrow\K\ot_B\L$ is a $\C$\comodule{} morphism.
 This natural map is an isomorphism, at least, in the following cases:
\begin{enumerate}
 \item $\C$ is a flat right $A$\module{} and\/ $\M$ is an injective
       left\/ $\C$\comodule;
 \item $\M$ is a quite\/ $\C/A$\injective{} left\/ $\C$\comodule{}
       and\/ $\K$ as a right\/ $\D$\comodule{} with a left $A$\module{}
       structure is coinduced from an $A$\+$B$\bimodule;
 \item $\M$ is a\/ $\C/A$\injective{} left\/ $\C$\comodule, $\D$ is
       a projective left $B$\module, $\L$ is a projective left 
       $B$\module, and\/ $\K$ as a right $\D$\comodule{} with a left
       $A$\module{} structure is coinduced from an $A$\projective{}
       $A$\+$B$\bimodule;
 \item $\M$ is a\/ $\C/A$\injective{} left\/ $\C$\comodule, $\D$ is
       a projective left $B$\module, $\L$ is a projective left 
       $B$\module, $\K$ is a projective left $A$\module, $\K$ is
       a\/ $\D/B$\+coflat right\/ $\D$\comodule, and the ring $A$ has
       a finite left homological dimension;
 \item $\L$ is a quasicoprojective left\/ $\D$\comodule.
\end{enumerate}
\end{prop2}

\begin{proof}
 Analogous to the proof of Proposition~1 and Proposition~3 below
(see also the proof of Proposition~\ref{co-tensor-co-hom-assoc}.1).
 In particular, to prove~(e) one notices that the sequence
$0\rarrow\Hom_\C(\K,\M)\rarrow\Hom_A(\K,\M)\rarrow
\Hom_A(\K\;\C\ot_A\M)$ remains exact after taking
the cohomomorphisms from~$\L$.
\end{proof}

\begin{prop3}
 Let\/ $\P$ be a left\/ $\C$\contramodule, $\K$ be
a\/ $\D$\+$\C$\bicomodule, and\/ $\Q$ be a left\/ $\D$\contramodule.
 Then there is a natural map\/ $\Cohom_\D(\K\ocn_\C\P\;\Q)\rarrow
\Hom^\C(\P,\Cohom_\D(\K,\Q))$ whenever the cohomomorphism module\/
$\Cohom_\D(\K,\Q)$ is endowed with a left\/ $\C$\contramodule{}
structure such that the map\/ $\Hom_B(\K,\Q)\rarrow\Cohom_\D(\K,\Q)$
is a\/ $\C$\contramodule{} morphism.
 This natural map is an isomorphism, at least, in the following cases:
\begin{enumerate}
 \item $\C$ is a projective left $A$\module{} and\/ $\P$ is
       a projective left\/ $\C$\contramodule;
 \item $\P$ is a quite\/ $\C/A$\projective{} left\/ $\C$\contramodule{}
       and\/ $\K$ as a left\/ $\D$\comodule{} with a right $A$\module{}
       structure is coinduced from a $B$\+$A$\bimodule,
 \item $\P$ is a\/ $\C/A$\projective{} left\/ $\C$\contramodule, $\D$ is
       a flat right $B$\module, $\Q$ is an injective left $B$\module,
       and\/ $\K$ as a left\/ $\D$\comodule{} with a right $A$\module{}
       structure is coinduced from an $A$\+flat $B$\+$A$\bimodule;
 \item $\P$ is a\/ $\C/A$\projective{} left\/ $\C$\contramodule, $\D$ is
       a flat right $B$\module, $\Q$ is an injective left $B$\module,
       $\K$ is a flat right $A$\module, $\K$ is a\/ $\D/B$\projective{}
       left\/ $\D$\comodule, and the ring $A$ has a finite left
       homological dimension;
 \item $\Q$ is a quasicoinjective left\/ $\D$\contramodule.
\end{enumerate}
\end{prop3}

\begin{proof}
 The map $\Hom^\C(\P,\Hom_B(\K,\Q))\rarrow\Hom^\C(\P,\Cohom_\D(K,\Q))$
annihilates the difference of two maps
$\Hom^\C(\P,\Hom_B(\D\ot_B\K,\Q))\rarrow\Hom^\C(\P,\Hom_B(\K,\Q))$
and this pair of maps can be identified with the pair of maps
$\Hom_B(\D\ot_B\K\ocn_\C\P\;\allowbreak\Q)\birarrow
\Hom_B(\K\ocn_\C\P\;\Q)$ whose cokernel is, by the definition,
the cohomomorphism module $\Cohom_\D(\K\ocn_\C\P\;\Q)$.
 Hence there is a natural map $\Cohom_\D(\K\ocn_\C\P\;\Q)\rarrow
\Hom^\C(\P,\Cohom_\D(\K,\Q))$.
 Besides, the composition $\Cohom_\D(\K\ocn_\C\P\;\Q)\rarrow
\Cohom_\D(\K\ot_A\P\;\Q)\rarrow\Hom_A(\P,\Cohom_\D(\K,\Q))$ has
equal compositions with two maps $\Hom_A(\P,\Cohom_\D(\K,\Q))
\birarrow\Hom_A(\Hom_A(\C,\P),\allowbreak\Cohom_\D(\K,\Q))$,
which leads to the same natural map $\Cohom_\D(\K\ocn_\C\P\;\Q)\rarrow
\Hom^\C(\P,\Cohom_\D(\K,\Q))$.
 To prove cases~(a-d), one shows that the sequence
$\Hom_B(\D\ot_B\K\;\Q)\rarrow\Hom_B(\K,\Q)\rarrow\Cohom_\D(\K,\Q)
\rarrow0$ remains exact after applying the functor $\Hom^\C(\P,{-})$.
 Indeed, the case~(a) is obvious, in the cases~(b-d) this sequence
of left $A$\module s splits, and in the cases~(c-d) it is also
an exact sequence of injective $A$\module s (see the proof of
Proposition~\ref{co-tensor-co-hom-assoc}.2).
 To prove~(e), one notices that the sequence 
$\K\ot_A\Hom_A(\C,\P)\rarrow\K\ot_A\P\rarrow\K\ocn_\C\P\rarrow0$
remains exact after taking the cohomomorphisms into~$\Q$ and
uses Proposition~\ref{co-tensor-co-hom-assoc}.2(b).
\end{proof}

 In the case of $\K=\D=\C$, the natural maps defined in
Propositions~2--3 have the following property of compatibility with
the adjoint functors $\Psi_\C$ and $\Phi_\C$: for any left
$\C$\comodule{} $\M$ and any left $\C$\contramodule{} $\P$ the maps
$\Cohom_\C(\Phi_\C(\P),\Psi_\C(\M))\rarrow\Hom_\C(\Phi_\C(\P),\M)$ and
$\Cohom_\C(\Phi_\C(\P),\Psi_\C(\M))\rarrow\Hom^\C(\P,\Psi_\C(\M))$ form
a commutative diagram with the isomorphism $\Hom_\C(\Phi_\C(\P),\M)
\simeq\Hom^\C(\P,\Psi_\C(\M))$.

\smallskip
 The following important Lemma is deduced as a corollary
of Propositions~2--3.

\begin{lem}
 \textup{(a)} A $\C$\comodule{} is quasicoprojective if and only if
it is quite\/ $\C/A$\injective.
 If\/ $\C$ is a projective left $A$\module, then a left\/
$\C$\comodule{} is coprojective if and only if it is a direct summand
of a comodule coinduced from a projective $A$\module. \par
 \textup{(b)} A $\C$\contramodule{} is quasicoinjective if and only if
it is quite\/ $\C/A$\projective.
 If\/ $\C$ is a flat right $A$\module, then a left\/
$\C$\contramodule{} is coinjective if and only if it is a direct
summand of a contramodule induced from an injective $A$\module.
\end{lem}

\begin{proof}
 Part~(a): let $\M$ be a quasicoprojective left $\C$\comodule.
 Denote by $l$ the coaction map $\M\rarrow\C\ot_A\M$.
 It is an $A$\+split injective morphism of quasicoprojective
$\C$\comodule s.
 According to Proposition~2(e), we have an isomorphism of morphisms
$\Hom_\C(l,\M)\simeq\Cohom_\C(l,\Hom_\C(\C,\M))$.
 But the map $\Cohom_\C(l,\P)$ is surjective for any left
$\C$\contramodule{} $\P$.
 Therefore, the map $\Hom_\C(l,\M)$ is also surjective, hence
the morphism~$l$ splits and the comodule $\M$ is quite
$\C/A$\injective.
 Now suppose that $\M$ is coprojective; then we already know that
$\M$ is quite $\C/A$\injective.
 Set $\P=\Psi_\C(\M)$.
 It follows from Proposition~3(b) that there is an isomorphism of
functors $\Hom^\C(\P,{-})\simeq\Cohom_\C(\M,{-})$ on the category
of left $\C$\contramodule s.
 Therefore, the $\C$\contramodule{} $\P$ is projective, hence it is
a direct summand of a $\C$\contramodule{} induced from a projective
$A$\module{} and $\M$ is a direct summand of the $\C$\comodule{}
coinduced from the same projective $A$\module.
 The proof of part~(b) is completely analogous; it uses
Propositions~3(e) and~2(b).
\end{proof}

\begin{qst}
 Are there any analogues of the results of Lemma for (quasi)coflat
comodules and (quite relatively) contraflat contramodules?
\end{qst}

\subsection{Relatively injective comodules and relatively projective
contramodules}   \label{rel-inj-proj-co-contra-mod}
 Assume that $\C$ is a projective left $A$\module.
 For any right $\C$\comodule{} $\N$ and any left $\C$\contramodule{}
$\P$ denote by $\Ctrtor^\C_i(\N,\P)$ the sequence of left derived
functors in the second argument of the right exact functor of
contratensor product $\N\ocn_\C\P$.
 By the definition, the $k$\module s $\Ctrtor^\C_i(\N,\P)$ are
computed using a left projective resolution of the
$\C$\contramodule{} $\P$.
 Since projective contramodules are contraflat, the functor
$\Ctrtor^\C_*(\N,\P)$ assigns long exact sequences to exact
triples in either of its arguments.

\begin{qst}
 Can one compute the derived functor $\Ctrtor$ using contraflat
resolutions of the second argument?
 In other words, is it true that $\Ctrtor^\C_i(\N,\P)=0$ for any
right $\C$\comodule{} $\N$, any contraflat left $\C$\contramodule{}
$\P$, and all $i>0$?
 Also, is it true that $\Ctrtor^\C_i(\N,\P)=0$ for any $A$\+flat
right $\C$\comodule{} $\N$, any (quite) $\C/A$\contraflat{} left
$\C$\contramodule{} $\P$, and all $i>0$?
 A related question: is $\Ctrtor^\C_{>0}(\N,\P)$ an effaceable
functor of its first argument?
\end{qst}

 Now assume that $\C$ is a projective left and a flat right
$A$\module{} and the ring $A$ has a finite left homological dimension.

\begin{lem1}
 \textup{(a)} A left\/ $\C$\comodule\/ $\M$ is\/ $\C/A$\injective{}
if and only if for any $A$\projective{} left\/ $\C$\comodule\/ $\L$
the $k$\module s\/ $\Ext_\C^i(\L,\M)$ of Yoneda extensions in
the abelian category of left\/ $\C$\comodule s vanish for all $i>0$.
 In particular, the functor of\/ $\C$\comodule{} homomorphisms from
an $A$\projective{} left\/ $\C$\comodule\/ $\L$ maps exact triples of\/
$\C/A$\injective{} left\/ $\C$\comodule s to exact triples.
 Besides, the class of\/ $\C/A$\injective{} left\/ $\C$\comodule s
is closed under extensions and cokernels of injective morphisms. \par
 \textup{(b)} A left\/ $\C$\contramodule\/ $\P$ is\/ $\C/A$\projective{}
if and only if for any $A$\injective{} left\/ $\C$\contramodule\/ $\Q$
the $k$\module s $\Ext^{\C,i}(\P,\Q)$ of Yoneda extensions in
the abelian category of left\/ $\C$\contramodule s vanish for all $i>0$.
 In particular, the functor of\/ $\C$\contramodule{} homomorphisms into
an $A$\injective{} left\/ $\C$\contramodule\/ $\Q$ maps exact triples
of\/ $\C/A$\projective{} left\/ $\C$\contramodule s to exact triples.
 Besides, the class of\/ $\C/A$\projective{} left\/ $\C$\contramodule s
is closed under extensions and kernels of surjective morphisms. \par
 \textup{(c)} For any\/ $\C/A$\projective{} left\/ $\C$\contramodule\/
$\P$ and any $A$\+flat right\/ $\C$\comodule\/ $\N$ the $k$\module s\/
$\Ctrtor^\C_i(\N,\P)$ vanish for all $i>0$.
 In particular, the functor of contratensor product with an $A$\+flat
right\/ $\C$\comodule{} maps exact triples of\/ $\C/A$\projective{}
left\/ $\C$\contramodule s to exact triples.
\end{lem1}

\begin{proof}
 Part~(a): the ``if'' part of the first assertion is obvious;
let us prove the ``only if'' part.
 An arbitrary element of $\Ext_\C^i(\L,\M)$ can be represented by 
a morphism of degree~$i$ from an exact complex $\dsb\to \L_i\to\dsb
\to \L_0\to \L\to0$ to the comodule~$\M$.
 According to Lemma~\ref{proj-inj-co-contra-module}(a), any left
$\C$\comodule{} is a surjective image of an $A$\projective{}
$\C$\comodule.
 Therefore, one can assume that the comodules $\L_i$ are
$A$\projective.
 Now if $\L$ is also $A$\projective, then our exact complex of
$\C$\comodule s is composed of exact triples of $A$\projective{}
$\C$\comodule s, so if $\M$ is $\C/A$\injective, then the complex
of homomorphisms into $\M$ from this complex of $\C$\comodule s
is acyclic.
 The remaining two assertions follow from the first one.
 The proof of part~(b) is completely analogous.
 To prove~(c), notice the isomorphism $\Hom_k(\Ctrtor^\C_i(\N,\P),
\.k\dual)\simeq\Ext^{\C,i}(\P,\Hom_k(\N,k\dual))$.
\end{proof}

\begin{rmk}
 Analogues of the third assertion of Lemma~1(a) and the third
assertion of Lemma~1(b) are \emph{not} true for quite relatively
injective comodules and quite relatively projective contramodules
(see Remark~\ref{co-contra-faithful-change-derived}; cf.\
Remark~\ref{co-contra-model-struct}).
\end{rmk}

\begin{thm}
 For any\/ $\C/A$\injective{} left\/ $\C$\comodule\/ $\M$ the left\/
$\C$\contramodule\/ $\Psi_\C(\M)$ is\/ $\C/A$\projective{} and
for any\/ $\C/A$\projective{} left\/ $\C$\contramodule\/ $\P$
the left\/ $\C$\+co\-mod\-ule\/ $\Phi_\C(\M)$ is\/ $\C/A$\injective.
 The restrictions of the functors\/ $\Psi_\C$ and\/ $\Phi_\C$ to
the full subcategories of\/ $\C/A$\injective\/ $\C$\comodule s
and\/ $\C/A$\projective\/ $\C$\contramodule s are mutually inverse
equivalences between these subcategories.
\end{thm}

\begin{proof}
 Let us first show that the injective dimension of a $\C/A$\injective{}
left $\C$\+co\-mod\-ule $\M$ in the abelian category of $\C$\comodule s
does not exceed the left homological dimension~$d$ of the ring~$A$.
 Indeed, it follows from Lemma~\ref{proj-inj-co-contra-module}(a) that
any left $\C$\comodule{} $\L$ has a finite resolution
$0\to\L_d\to\dsb\to\L_0\to\L\to0$ with $A$\projective{}
$\C$\comodule s $\L_j$; and since $\Ext^i(\L_j,\M)=0$ for all~$j$
and all $i>0$, the complex $\Hom_\C(\L_\bu,\M)$ computes
$\Ext^*_\C(\L,\M)$.
 So the $\C$\comodule{} $\M$ has a finite injective resolution, and
consequently it has a finite resolution $0\to \M\to \K^0\to\dsb\to
\K^d\to0$ consisting of quite $\C/A$\injective{} $\C$\comodule s
$\K^j$.
 According to Lemma~1(a), this exact sequence is composed of exact
triples of $\C/A$\injective{} $\C$\comodule s, which the functor
$\Psi_\C$ maps to exact triples; so the sequence $0\rarrow\Psi_\C(\M)
\rarrow\Psi_\C(\K^0)\rarrow\dsb\rarrow\Psi_\C(\K^d)\rarrow0$
is also exact.
 Since the $\C$\contramodule s $\Psi_\C(\K^j)$ are quite
$\C/A$\projective, it follows from Lemma~1(b) that
the $\C$\contramodule{} $\Psi_\C(\M)$ is $\C/A$\projective{} and
the latter exact sequence is composed of exact triples of
$\C/A$\projective{} $\C$\contramodule s.
 Thus it follows from Lemma~1(c) that the sequence $0\rarrow
\Phi_\C\Psi_\C(\M)\rarrow\Phi_\C\Psi_\C(\K^0)\rarrow\dsb\rarrow
\Phi_\C\Psi_\C(\K^d)\rarrow0$ is also exact.
 Now since the adjunction maps $\Phi_\C\Psi_\C(\K^j)\rarrow\K^j$
are isomorphisms, the adjunction map $\Phi_\C\Psi_\C(\M)\rarrow\M$
is also an isomorphism.
 The remaining assertions are proven in the completely analogous way.
\end{proof}

\begin{lem2}
 \textup{(a)} In the above assumptions, a left\/ $\C$\comodule{} is\/
$\C/A$\coprojective{} if and only if it is\/ $\C/A$\injective. \par
 \textup{(b)} In the above assumptions, a left\/ $\C$\contramodule{}
is\/ $\C/A$\coinjective{} if and only if it is\/ $\C/A$\projective.
\end{lem2}

\begin{proof}
 Part~(a) in the ``if'' direction: it follows from
Proposition~\ref{cotensor-contratensor-assoc}.3(c) that whenever
a left $\C$\contramodule{} $\P$ is $\C/A$\projective,
the $\C$\comodule{} $\Phi_\C(\P)$ is $\C/A$\coprojective.
 Now if a left $\C$\comodule{} $\M$ is $\C/A$\injective, then
the $\C$\contramodule{} $\P=\Psi_\C(\M)$ is $\C/A$\projective{}
and $\M=\Phi_\C(\P)$ by the above Theorem.
 Part~(a) in the ``only if'' direction: in view of Lemma~1(a),
the construction of Lemmas~\ref{flat-comodule-surjection}
and~\ref{proj-inj-co-contra-module}(a) represents any left
$\C$\comodule{} $\M$ as the quotient comodule of an $A$\projective{}
$\C$\comodule{} $\cP(\M)$ by a $\C/A$\injective{} $\C$\comodule.
 We will show that whenever $\M$ is a $\C/A$\coprojective{}
$\C$\comodule, $\cP(\M)$ is a coprojective $\C$\comodule; then if
will follow that $\M$ is a $\C/A$\injective{} $\C$\comodule{} by
Lemma~\ref{cotensor-contratensor-assoc}(a) and Lemma~1(a).
 Indeed, an extension of $\C/A$\coprojective{} left $\C$\comodule s is
$\C/A$\coprojective{} by Lemma~\ref{absolute-relative-coproj-coinj}(a);
let us check that an $A$\projective{} $\C/A$\coprojective{}
$\C$\comodule{} is coprojective.
 For any left $\C$\comodule{} $\M$ and any left $\C$\contramodule{}
$\P$ denote by $\Coext^i_\C(\M,\P)$ the cohomology of the object
$\Coext_\C(\M,\P)$ of the derived category $\sD(k\modl)$ that was
constructed in~\ref{semiext-definition}.
 This definition agrees with the definition of $\Cotor_\C^*(\M,\P)$
for an $A$\projective{} $\C$\comodule{} $\M$ or an $A$\injective{}
$\C$\contramodule{} $\P$ given in the proof of
Lemma~\ref{absolute-relative-coproj-coinj}.
 The functor $\Cotor_\C^*(\M,\P)$ assigns long exact sequences to
exact triples in either of its arguments.
 For any $A$\projective{} left $\C$\comodule{} $\M$ and any left
$\C$\contramodule{} $\P$ one has $\Coext^i_\C(\M,\P)=0$ for all $i>0$
and $\Coext^0_\C(\M,\P)\simeq\Cohom_\C(\M,\P)$.
 Therefore, an $A$\projective{} left $\C$\comodule{} $\M$ is
coprojective if and only if $\Cohom^i_\C(\M,\P)=0$ for any left
$\C$\contramodule{} $\P$ and all $i\ne0$.
 For any $\C/A$\coprojective{} left $\C$\comodule{} $\M$ and any left
$\C$\comodule{} $\P$ one has $\Coext^i_\C(\M,\P)=0$ for all $i<0$,
since one can compute $\Coext_\C(\M,\P)$ using a finite $A$\injective{}
right resolution of~$\P$ by the result
of~\ref{relatively-semiproj-semiinj}.
 Thus an $A$\projective{} $\C/A$\coprojective{} left $\C$\comodule{}
is coprojective.
 The proof of part~(b) is completely analogous; it uses 
Proposition~\ref{cotensor-contratensor-assoc}.2(c) and
Lemma~\ref{proj-inj-co-contra-module}(b).
\end{proof}

\begin{qst}
 It follows from Proposition~1(c) that if $\C$ is a flat right
$A$\module, then whenever a left $\C$\contramodule{} $\P$ is
$\C/A$\contraflat{} the $\C$\comodule{} $\Phi_\C(\P)$ is
$\C/A$\+coflat.
 Does the converse hold?
\end{qst}

\subsection{Comodule-contramodule correspondence}
\label{comodule-contramodule-subsect}
 Assume that the coring $\C$ is a projective left and a flat right
$A$\module{} and the ring $A$ has a finite left homological dimension.

 The categories of $\C/A$\injective{} left $\C$\comodule s and
$\C/A$\projective{} left $\C$\contramodule s have natural exact
category structures as full subcategories, closed under extensions,
of the abelian categories of left $\C$\comodule s and left
$\C$\contramodule s.

\begin{thm}
 \textup{(a)} The functor mapping the quotient category of the homotopy
category of complexes of\/ $\C/A$\injective{} left\/ $\C$\comodule s by
its minimal triagulated subcategory containing the total complexes of
exact triples of complexes of\/ $\C/A$\injective\/ $\C$\comodule s into
the coderived category of left\/ $\C$\comodule s is an equivalence of
triangulated categories. \par 
 \textup{(b)} The functor mapping the quotient category of the homotopy
category of complexes of\/ $\C/A$\projective{} left\/
$\C$\contramodule s by its minimal triangulated subcategory containing
the total complexes of\/ $\C/A$\projective\/ $\C$\contramodule s into
the contraderived category of left\/ $\C$\contramodule s is
an equivalence of triangulated categories.
\end{thm}

\begin{proof}
 Part~(a): let $\M^\bu$ be a complex of left $\C$\comodule s.
 Then the total complex of the cobar bicomplex
$\C\ot_A\M^\bu\rarrow\C\ot_A\C\ot_A\M^\bu\rarrow\dsb$
is a complex of (quite) $\C/A$\injective{} $\C$\comodule s,
the complex $\M^\bu$ maps into this total complex, and the cone
of this map is coacyclic.
 Hence it follows from Lemma~\ref{semitor-main-theorem} that
the coderived category of left $\C$\comodule s is equivalent to
the quotient category of the homotopy category of complexes of
$\C/A$\injective{} $\C$\comodule s by its intersection with the thick
subcategory of coacyclic complexes of $\C$\comodule s.
 It remains to show that this intersection of subcategories coincides
with the minimal triangulated subcategory containing the total
complexes of exact triples of complexes of $\C/A$\injective{}
$\C$\comodule s.

\begin{lem}
 \textup{(a)} For any exact category\/~$\sA$ where infinite direct sums
exist and preserve exact triples, the complex of homomorphisms from
a coacyclic complex over\/~$\sA$ into a complex of injective objects
with respect to\/~$\sA$ is acyclic. \par
 \textup{(b)} For any exact category\/~$\sA$ where infinite products
exist and preserve exact triples, the complex of homomorphisms from
a complex of projective objects with respect to\/~$\sA$ into
a contraacyclic complex over\/~$\sA$ is acyclic.
\end{lem}

\begin{proof}
 Analogous to the proofs of Lemmas~\ref{coflat-complexes}
and~\ref{co-proj-inj-complexes}.
 Part~(a): let $M^\bu$ be a complex of injective objects with
respect to~$\sA$.
 Since the functor of homomorphisms into~$M^\bu$ maps distinguished
triangles in the homotopy category to distinguished triangles and
infinite direct sums to infinite products, it suffices to check that
the complex $\Hom_\sA(L^\bu,M^\bu)$ is acyclic whenever $L^\bu$ is
the total complex of an exact triple
${}'\!K^\bu\to K^\bu\to {}''\!K^\bu$ of complexes over~$\sA$.
 But the complex $\Hom_\sA(L^\bu,M^\bu)$ is the total complex of
an exact triple of complexes of abelian groups
$\Hom_\sA(\.{}''\!K^\bu,M^\bu)\rarrow\Hom_\sA(K^\bu,M^\bu)\rarrow
\Hom_\sA(\.{}'\!K^\bu,M^\bu)$ in this case.
 The proof of part~(b) is dual.
\end{proof}

 We will show that (i)~the minimal triangulated subcategory containing
the total complexes of exact triples of complexes of $\C/A$\injective{}
$\C$\comodule s and (ii)~the homotopy category of complexes of
injective $\C$\comodule s form a semiorthogonal decomposition of
the homotopy category of complexes of $\C/A$\injective{} left
$\C$\comodule s.
 This means, in addition to the subcategory~(i) being left orthogonal
to the subcategory~(ii), that for any complex $\K^\bu$ of
$\C/A$\injective{} $\C$\comodule s there exists a (unique and
functorial) distinguished triangle $\L^\bu\to\K^\bu\to\M^\bu\to
\L^\bu[1]$ in the homotopy category of $\C$\comodule s, where $\L^\bu$
belongs to the subcategory~(i) and $\M^\bu$ belongs to
the subcategory~(ii).
 It will follow that the subcategory~(i) is the maximal subcategory of
the homotopy category of complexes of $\C/A$\injective{}
$\C$\comodule s left orthogonal to the subcategory~(ii), hence
the subcategory~(i) contains the intersection of the homotopy category
of complexes of $\C/A$\injective{} $\C$\comodule s with the thick
subcategory of coacyclic complexes of $\C$\comodule s.

 Indeed, let $\K^\bu$ be a complex of $\C/A$\injective{} left
$\C$\comodule s.
 Choose for every~$n$ an injection~$j^n$ of the $\C$\comodule{} $\K^n$
into an injective $\C$\comodule{} $\cJ^n$.
 Consider the complex $\E^\bu=\E(\K^\bu)$ whose terms are
the $\C$\comodule s $\E^n=\cJ^n\oplus\cJ^{n+1}$ and
the differential $d_\E^n\:\E^n\rarrow\E^{n+1}$ maps $\cJ^{n+1}$ into
itself by the identity map and vanishes in the restriction to
$\cJ^n$ and in the projection to $\cJ^{n+2}$.
 There is a natural injective morphism of complexes
$\K^\bu\rarrow\E^\bu$ formed by the $\C$\comodule{} maps
$\K^n\rarrow\E^n$ whose components are $j^n\:\K^n\rarrow\cJ^n$
and $j^{n+1}d_\K^n\:\K^n\rarrow\cJ^{n+1}$.
 Set ${}^0\E^\bu=\E(\K^\bu)$, \ ${}^1\E^\bu=\E({}^0\E^\bu/\K^\bu)$,
etc.
 As it was shown in the proof of
Theorem~\ref{rel-inj-proj-co-contra-mod}, the injective dimension of
a $\C/A$\injective{} left $\C$\comodule{} does not exceed the left
homological dimension~$d$ of the ring~$A$.
 Therefore, the complex $\cZ^\bu=\coker({}^{d-2}\E^\bu\to
{}^{d-1}\E^\bu)$ is a complex of injective $\C$\comodule s.
 Now it is clear that the total complex $\M^\bu$ of the bicomplex
${}^0\E^\bu\rarrow{}^1\E^\bu\rarrow\dsb\rarrow{}^{d-1}\E^\bu\rarrow
\cZ^\bu$ is a complex of injective $\C$\comodule s and the cone
$\L^\bu$ of the morphism $\K^\bu\rarrow\M^\bu$ belongs to the minimal
triangulated subcategory containing the total complexes of exact
triples of complexes of $\C/A$\injective{} $\C$\comodule s by
Lemma~\ref{rel-inj-proj-co-contra-mod}.1(a).

 Part~(a) is proven; the proof of part~(b) is completely analogous
and uses Lemma~\ref{rel-inj-proj-co-contra-mod}.1(b).
\end{proof}

\begin{rmk}
 Let $\sA$ be an exact category where infinite direct sums exist
and preserve exact triples, every object admits an admissible
monomorphism into an object injective relative to~$\sA$, and
the class of such injective objects is closed under infinite
direct sums.
 Then the thick subcategory of coacyclic complexes with respect
to~$\sA$ and the triangulated subcategory of complexes of injective
objects form a semiorthogonal decomposition of the homotopy
category $\Hot(\sA)$, so the coderived category $\sD^\co(\sA)$ is
equivalent to the homotopy category of complexes of injectives in~$\sA$.
 Indeed, orthogonality is already proven in Lemma, so it remains
to construct a morphism from any complex $C^\bu$ over $\sA$
into a complex of injectives $M^\bu$ with a coacyclic cone.
 To do so, one proceeds as in the proof of Theorem, constructing
a morphism from $C^\bu$ into a complex of injectives ${}^0\.\!E^\bu$
that is an admissible monomorphism in every degree, taking
the quotient complex, constructing an analogous morphism from
it into a complex of injectives ${}^1\.\!E^\bu$, etc.
 Finally, one constructs the total complex $M^\bu$ of the bicomplex
${}^\bu\.\!E^\bu$ by taking infinite direct sums along the diagonals;
then $M^\bu$ is a complex of injectives and the cone of the morphism
$C^\bu\rarrow M^\bu$ is coacyclic by Lemma~\ref{coderived-categories}. 
 Consequently, the homotopy category of acyclic complexes of
injectives in~$\sA$ is equivalent to the quotient category
$\Acycl(\sA)/\Acycl^{\co}(\sA)$ and to the kernel of the localization
functor $\sD^{\co}(\sA)\rarrow\sD(\sA)$ (cf.~\cite{Kra}).
 When $\sA$ has a finite homological dimension, the condition that
the class of injectives is closed under infinite direct sums is
not needed in this argument.
 This is a somewhat trivial situation, though; see
Remark~\ref{coderived-categories}.
 Moreover, let $\sA$ be an exact category where infinite direct
sums exist and preserve exact triples and every object admits
an admissible monomorphism into an injective.
 Let $\sF\subset\sA$ be a class of objects closed under
cokernels of admissible monomorphisms, containing the injectives, and
consisting of objects of finite injective dimension.
 Then every complex over $\sF$ that is coacyclic as a complex over
$\sA$ belongs to the minimal triangulated subcategory of $\Hot(\sA)$
containing the total complexes of exact triples of complexes
over~$\sF$.
 When there is a class $\sF\subset\sA$ closed under infinite direct
sums, consisting of objects of finite injective dimension, and such
that every object of $\sA$ admits an admissible monomorphism into
an object of~$\sF$, the coderived category $\sD^\co(\sA)$ is equivalent
to the homotopy category of complexes of injectives in~$\sA$
(cf.~\cite{Jorg}, where in the dual situation the role of
the class~$\sF$ is played by flat modules).
 To show this, one has to repeat twice the above construction of
a resolution ${}^\bu\.\!E^\bu$, taking infinite direct sums along
the diagonals for the first time and finite directs sums along
the diagonals of the canonical truncation for the second time.
 When $\sA$ is the abelian category of $\C$\comodule s, one can take
$\sF$ to be the class of $\C/A$\injective{} $\C$\comodule s
or quite $\C/A$\injective{} $\C$\comodule s.
 The related results for comodules and contramodules are obtained in
Theorem~\ref{co-contra-ctrtor-definition} and
Remark~\ref{co-contra-ctrtor-definition}.
\end{rmk}

\begin{cor}
 The restrictions of the functors\/ $\Psi_\C$ and\/ $\Phi_\C$ (applied
to complexes term-wise) to the homotopy category of complexes of\/
$\C/A$\injective\/ $\C$\comodule s and the homotopy category of
complexes of\/ $\C/A$\projective\/ $\C$\contramodule s define
mutually inverse equivalences\/ $\boR\Psi_\C$ and\/ $\boL\Phi_\C$
between the coderived category of left\/ $\C$\comodule s and
the contraderived category of left\/ $\C$\contramodule s.
\end{cor}

\begin{proof}
 By Theorem~\ref{rel-inj-proj-co-contra-mod}, the functors $\Psi_\C$
and $\Phi_\C$ induce mutually inverse equivalences between
the homotopy categories of $\C/A$\injective{} left $\C$\comodule s
and $\C/A$\projective{} left $\C$\contramodule s.
 According to Lemma~\ref{rel-inj-proj-co-contra-mod}.1(a) and~(c),
the total complexes of exact triples of complexes of $\C/A$\injective{}
$\C$\comodule s correspond to the total complexes of exact triples
of complexes of $\C/A$\projective{} $\C$\contramodule s under
this equivalence.
 So it remains to apply the above Theorem.
\end{proof}

\begin{qst}
 Can one obtain a version of the derived comodule-contramodule
correspondence (an equivalence between appropriately defined exotic
derived categories of left $\C$\comodule s and left $\C$\contramodule s)
not depending on any assumptions about the homological dimension of
the ring~$A$?
 It is not difficult to see that one can weaken the assumption that
$A$ has a finite left homological dimension to the assumption that
$A$ is left Gorenstein, i.~e., the classes of left $A$\module s of
finite projective and injective dimensions coincide.
 In this case, the coderived category of left $\C$\comodule s and
the contraderived category of left $\C$\contramodule s are naturally
equivalent whenever the coring $\C$ is a projective left and a flat
right $A$\module.
 Indeed, arguing as in Theorem~\ref{co-contra-ctrtor-definition} below,
one can show that the coderived category of left $\C$\comodule s is
equivalent to the quotient category of the homotopy category of
complexes of $\C$\comodule s coinduced from left $A$\module s of
finite projective (injective) dimension by its minimal triangulated
subcategory containing the total complexes of exact triples of
complexes of $\C$\comodule s that at every term are exact triples of
$\C$\comodule s coinduced from exact triples of $A$\module s of
finite projective (injective) dimension.
 Analogously, the contraderived category of left $\C$\contramodule s
is equivalent to the quotient category of the homotopy category of
complexes of $\C$\contramodule s induced from left $A$\module s of
finite projective (injective) dimension by its minimal triangulated
subcategory containing the total complexes of exact triples of
complexes of $\C$\contramodule s that at every term are exact triples
of $\C$\contramodule s induced from exact triples of $A$\module s
of finite projective (injective) dimension.
 The key step is to notice that the class of left $A$\module s of
finite projective (injective) dimension is closed under infinite
direct sums and products.
\end{qst}

\subsection{Derived functor Ctrtor}
\label{co-contra-ctrtor-definition}
 The following analogue of Theorem~\ref{comodule-contramodule-subsect}
holds under slightly weaker conditions.

\begin{thm}
 \textup{(a)} Assume that the coring\/ $\C$ is a flat right
$A$\module{} and the ring~$A$ has a finite left homological
dimension.
 Then the functor mapping the homotopy category of complexes of
injective left\/ $\C$\comodule s into the coderived category of left\/
$\C$\comodule s is an equivalence of triangulated categories.
 In addition, the functor mapping the quotient category of the homotopy
category of complexes of quite\/ $\C/A$\injective{} left\/
$\C$\comodule s by the minimal triangulated subcategory containing
the total complexes of exact triples of complexes of coinduced\/
$\C$\comodule s that at every term are exact triples of\/
$\C$\comodule s coinduced from exact triples of $A$\module s into
the coderived category of left\/ $\C$\comodule s is an equivalence of
triangulated categories. \par
 \textup{(b)} Assume that the coring\/ $\C$ is a projective left
$A$\module{} and the ring~$A$ has a finite left homological
dimension.
 Then the functor mapping the homotopy category of complexes of
projective left\/ $\C$\contramodule s into the contraderived category
of left\/ $\C$\contramodule s is an equivalence of triangulated
categories.
 In addition, the functor mapping the quotient category of the homotopy
category of complexes of quite\/ $\C/A$\projective{} left\/
$\C$\contramodule s by the minimal triangulated subcategory containing
the total complexes of exact triples of complexes of induced\/
$\C$\contramodule s that at every term are exact triples of\/
$\C$\contramodule s induced from exact triples of $A$\module s into
the contraderived category of left\/ $\C$\contramodule s is
an equivalence of triangulated categories.
\end{thm}

\begin{proof}
 Part~(a): when $\C$ is also a projective left $A$\module, the first
assertion follows from the proof of
Theorem~\ref{comodule-contramodule-subsect}.
 To prove both assertions in the general case, we will show that
(i)~the minimal triangulated subcategory containing the total complexes
of exact triples of complexes of coinduced $\C$\comodule s that at
every term are exact triples of $\C$\comodule s coinduced from exact
triples of $A$\module s and (ii)~the homotopy category of complexes
of injective $\C$\comodule s form a semiorthogonal decomposition of
the homotopy category of complexes of quite $\C/A$\injective{} left
$\C$\comodule s.
 Then we will argue as in the proof of
Theorem~\ref{comodule-contramodule-subsect}.

 Any complex of quite $\C/A$\injective{} $\C$\comodule s is homotopy
equivalent to a complex of coinduced $\C$\comodule s.
 Let $\K^\bu$ be a complex of coinduced left $\C$\comodule s; then
$\K^n\simeq\C\ot_A V^n$ for certain $A$\module s $V^n$.
 Let $V^i\rarrow I^i$ be injective maps of the $A$\module s $V^n$
into injective $A$\module s $I^n$; set $\cJ^n=\C\ot_A I^n$.
 Then $\cJ^n$ are injective $\C$\comodule s endowed with injective
$\C$\comodule{} morphisms $\K^n\rarrow\cJ^n$.
 As in the proof of Theorem~\ref{comodule-contramodule-subsect}, we
construct the complex of injective $\C$\comodule s $\E^\bu$ with
$\E^n=\cJ^n\oplus\cJ^{n+1}$ and the injective morphism of complexes
$\K^\bu\rarrow\E^\bu$.
 Let us show that there exists an automorphism of the $\C$\comodule{}
$\E^n$ such that its composition with the injection $\K^n\rarrow\E^n$
is the injection whose components are $j_n\:\K^n\rarrow\cJ^n$ and
the zero map $\K^n\rarrow\cJ^{n+1}$.
 Since the comodule $\cJ^{n+1}$ is injective, the component
$\K^n\rarrow\cJ^{n+1}$ of the morphism $\K^n\rarrow\E^n$
can be extended from the comodule $\K^n$ to comodule $\cJ^n$
containing it.
 Denote the morphism so obtained by $h^n\:\cJ^n\rarrow\cJ^{n+1}$;
then the automorphism of the comodule $\E^n$ whose components are
$-h^n$, the identity automorphisms of $\cJ^n$ and $\cJ^{n+1}$,
and zero has the desired property.
 Now it is clear that the triple $\K^\bu\rarrow\E^\bu\rarrow
\E^\bu/\K^\bu$ is an exact triple of complexes of coinduced
$\C$\comodule s which at every term is an exact triple of
$\C$\comodule s coinduced from an exact triple of $A$\module s.
 Moreover, $\E^n/\K^n\simeq\C\ot_A W^n$, where the injective dimension
$\di_A W^n$ is equal to $\di_A V^n-1$.
 So we can iterate this (nonfunctorial) construction, setting
${}^0\E^\bu=\E(\K^\bu)=\E^\bu$, \ ${}^1\E^\bu = \E({}^0\E^\bu/\K^\bu)$,
etc., and $\cZ^\bu = \coker({}^{d-2}\E^\bu\to{}^{d-1}\E^\bu)$.
 Then the total complex $\M^\bu$ of the bicomplex ${}^0\E^\bu\rarrow
{}^1\E^\bu\rarrow\dsb\rarrow{}^{d-1}\E^\bu\rarrow\cZ^\bu$ is a complex
of injective $\C$\comodule s and the cone $\L^\bu$ of the morphism
$\K^\bu\rarrow\M^\bu$ belongs to the minimal triangulated subcategory
containing the total complexes of exact triples of complexes of
coinduced $\C$\comodule s that at every term is an exact triple of
$\C$\comodule s coinduced from an exact triple of $A$\module s.
\end{proof}

\begin{rmk}
 The above Theorem provides an alternative way of proving
Corollary~\ref{comodule-contramodule-subsect}.
 Besides, it follows from~\ref{quite-injective-projective}, 
Lemma~\ref{cotensor-contratensor-assoc}, and the above Theorem that
in the assumptions of~\ref{comodule-contramodule-subsect} the functor
mapping the homotopy category of complexes of coprojective left
$\C$\comodule s into the coderived category of left $\C$\comodule s and
the functor mapping the homotopy category of complexes of coinjective
left $\C$\contramodule s into the contraderived category of left
$\C$\contramodule s are equivalences of triangulated categories.
 This is a stronger result than Theorem~\ref{coext-main-theorem}.
\end{rmk}

 The contratensor product $\N^\bu\ocn_\C\P^\bu$ of a complex of right
$\C$\comodule s $\N^\bu$ and a complex of left $\C$\contramodule s
$\P^\bu$ is defined as the total complex of the bicomplex
$\N^i\ocn_\C\P^j$, constructed by taking infinite direct sums along
the diagonals.
 Assume that the coring $\C$ is a projective left $A$\module{} and
the ring $A$ has a finite left homological dimension.
 One can prove in the way completely analogous to the proof of
Lemma~\ref{coflat-complexes} that the contratensor product of
a coacyclic complex of right $\C$\comodule s and a complex of
contraflat (and in particular, projective) left $\C$\contramodule s
is acyclic.
 The left derived functor of contratensor product
$$
 \Ctrtor^\C\:\sD^\co(\comodr\C)\times\sD^\ctr(\C\contra)\lrarrow
 \sD(k\modl)
$$
is defined by restricting the functor of contratensor product to
the Carthesian product of the homotopy category of right
$\C$\comodule s and the homotopy category of complexes of projective
left $\C$\contramodule s.

 The same derived functor can be obtained by restricting the functor
of contratensor product to the Carthesian product of the homotopy
category of complexes of $A$\+flat right $\C$\comodule s and
the homotopy category of complexes of quite $\C/A$\projective{} left
$\C$\contramodule s.
 Indeed, it follows from part~(b) of Theorem that the contratensor
product of a complex of $A$\+flat right $\C$\comodule s and
a contraacyclic complex of quite $\C/A$\projective{} left
$\C$\contramodule s is acyclic.
 Now if $\N^\bu$ is a complex of $A$\+flat right $\C$\comodule s,
$\P^\bu$ is a complex of quite $\C/A$\projective{} left
$\C$\contramodule s, and ${}'\P^\bu\rarrow\P^\bu$ is a morphism from
a complex of projective $\C$\contramodule s ${}'\P^\bu$ into $\P^\bu$
with a contraacyclic cone, then the map $\N^\bu\ocn_\C\P^\bu\rarrow
\N^\bu\ocn_\C{}'\P^\bu$ is a quasi-isomorphism.
 In particular, if the complex $\N^\bu$ is coacyclic, then the complex
$\N^\bu\ocn_\C\P^\bu$ is acyclic, since the complex
$\N^\bu\ocn_\C{}'\P^\bu$ is.
 When $\C$ is also a flat right $A$\module, one can use complexes
of $\C/A$\projective{} $\C$\contramodule s instead of complexes of
quite $\C/A$\projective{} $\C$\contramodule s, because the contratensor
product of a complex of $A$\+flat right $\C$\comodule s and
a contraacyclic complex of $\C/A$\projective{} left $\C$\contramodule s
is acyclic by Theorem~\ref{comodule-contramodule-subsect}(b)
and Lemma~\ref{rel-inj-proj-co-contra-mod}.1(c).
 Notice that this definition of the derived functor $\Ctrtor^\C$ is
\emph{not} a particular case of Lemma~\ref{semitor-definition}
(instead, it is a particular case of
Lemma~\ref{semi-ctrtor-definition}.2 below).

 Analogously, assume that the coring $\C$ is a flat right $A$\module{}
and the ring~$A$ has a finite left homological dimension.
 According to Lemma~\ref{comodule-contramodule-subsect}, the complex of
homomorphisms from a coacyclic complex of left $\C$\comodule s into
a complex of injective left $\C$\comodule s is acyclic.
 Therefore, the natural map $\Hom_{\Hot(\C\comodl)}(\L^\bu,\M^\bu)
\rarrow\Hom_{\sD^\co(\C\comodl)}(\L^\bu,\M^\bu)$ is an isomorphism
whenever $\M^\bu$ is a complex of injective $\C$\comodule s.
 So the functor of homomorphisms in the coderived category of left
$\C$\comodule s can be lifted to a functor
$$
 \Ext_\C\:\sD^\co(\C\comodl)^\op\times\sD^\co(\C\comodl)\lrarrow
 \sD(k\modl),
$$
which is defined by restricting the functor of homomorphisms of
complexes of $\C$\comodule s to the Carthesian product of the homotopy
category of left $\C$\comodule s and the homotopy category of
complexes of injective left $\C$\comodule s.

 The same functor $\Ext_\C$ can be obtained by restricting the functor
of homomorphisms to the Carthesian product of the homotopy category of
complexes of $A$\projective{} left $\C$\comodule s and the homotopy
category of complexes of quite $\C/A$\injective{} left $\C$\comodule s.
 Indeed, it follows from part~(a) of Theorem that the complex of
homomorphisms from a complex of $A$\projective{} left $\C$\comodule s
into a coacyclic complex of quite $\C/A$\injective{} left
$\C$\comodule s is acyclic.
 Now if $\L^\bu$ is a complex of $A$\projective{} left $\C$\comodule s,
$\M^\bu$ is a complex of quite $\C/A$\injective{} left $\C$\comodule s,
and $\M^\bu\rarrow{}'\M^\bu$ is a morphism from~$\M^\bu$ into a complex
of injective $\C$\comodule s ${}'\M^\bu$ with a coacyclic cone, then
the map $\Hom_\C(\L^\bu,\M^\bu)\rarrow\Hom_\C(\L^\bu,{}'\M^\bu)$ is
a quasi-isomorphism.
 When $\C$ is also a projective left $A$\module, one can use complexes
of $\C/A$\injective{} $\C$\comodule s instead of complexes of quite
$\C/A$\injective{} $\C$\comodule s.

 Finally, assume that the coring $\C$ is a projective left $A$\module{}
and the ring $A$ has a finite left homological dimension.
 By Lemma~\ref{comodule-contramodule-subsect}, the natural map
$\Hom_{\Hot(\C\contra)}(\P^\bu,\Q^\bu)\rarrow
\Hom_{\sD^\ctr(\C\contra)}(\P^\bu,\Q^\bu)$ is an isomorphism whenever
$\P^\bu$ is a complex of projective $\C$\contramodule s.
 So the functor of homomorphisms in the contraderived category of left
$\C$\contramodule s can be lifted to a functor
$$
 \Ext^\C\:\sD^\ctr(\C\contra)^\op\times\sD^\ctr(\C\contra)\lrarrow
 \sD(k\modl),
$$
which is defined by restricting the functor of homomorphisms of
complexes of $\C$\contramodule s to the Carthesian product of
the homotopy category of complexes of projective left
$\C$\contramodule s and the homotopy category of
left $\C$\contramodule s.

 The same functor $\Ext^\C$ can be obtained by restricting the functor
of homomorphisms to the Carthesian product of the homotopy category of
complexes of quite $\C/A$\projective{} left $\C$\contramodule s and
the homotopy category of complexes of $A$\injective{} left
$\C$\contramodule s.
 When $\C$ is also a flat right $A$\module, one can use complexes
of $\C/A$\projective{} $\C$\contramodule s instead of complexes of
quite $\C/A$\projective{} $\C$\contramodule s.

\subsection{Coext and Ext, Cotor and Ctrtor}
 Assume that the coring $\C$ is a projective left and a flat right
$A$\module{} and the ring $A$ has a finite left homological dimension.

{\hbadness=10000
\begin{cor}
 \textup{(a)} There are natural isomorphisms of functors\/
$\Coext_\C(\M^\bu,\P^\bu)\simeq\Ext_\C(\M^\bu,\boL\Phi_\C(\P^\bu))
\simeq\Ext^\C(\boR\Psi_\C(\M^\bu),\P^\bu)$ on the Carthesian product
of the category opposite to the coderived category of left\/
$\C$\comodule s and the contraderived category of left\/
$\C$\contramodule s. \par
 \textup{(b)} There is a natural isomorphism of functors\/
$\Cotor^\C(\N^\bu,\M^\bu)\simeq\Ctrtor^\C(\N^\bu,\boR\Psi_\C(\M^\bu))$
on the Carthesian product of the coderived category of right\/
$\C$\comodule s and the coderived category of left\/ $\C$\comodule s.
\end{cor}}

\begin{proof}
 Clearly, it suffices to construct natural isomorphisms
$\Coext_\C(\L^\bu,\boR\Psi_\C(\M^\bu))\simeq\Ext_\C(\L^\bu,\M^\bu)$, \ 
$\Coext_\C(\boL\Phi_\C(\P^\bu),\Q^\bu)\simeq\Ext^\C(\P^\bu,\Q^\bu)$,
and $\Cotor^\C(\N^\bu,\boL\Phi_\C(\P^\bu))\simeq
\Ctrtor^\C(\N^\bu,\P^\bu)$.
 In the first case, represent the image of $\M^\bu$ in
$\sD^\co(\C\comodl)$ by a complex of injective $\C$\comodule s,
notice that the functor $\Psi_\C$ maps injective comodules to
coinjective contramodules, and use
Proposition~\ref{cotensor-contratensor-assoc}.2(a).
 Alternatively, represent the image of $\M^\bu$ in
$\sD^\co(\C\comodl)$ by a complex of $\C/A$\injective{}
$\C$\comodule s and the image of $\L^\bu$ in $\sD^\co(\C\comodl)$
by a complex of coprojective{} $\C$\comodule s, and use
Proposition~\ref{cotensor-contratensor-assoc}.2(e); or
represent the image of $\M^\bu$ in $\sD^\co(\C\comodl)$ by a complex
of $\C/A$\injective{} $\C$\comodule s and the image of $\L^\bu$ in
$\sD^\co(\C\comodl)$ by a complex of $A$\projective{} $\C$\comodule s,
and use Proposition~\ref{cotensor-contratensor-assoc}.2(c),
Lemma~\ref{rel-inj-proj-co-contra-mod}.2(b), and the result
of~\ref{relatively-semiproj-semiinj}.
 In the second case, represent the image of $\P^\bu$ in
$\sD^\ctr(\C\contra)$ by a complex of projective $\C$\contramodule s,
notice that the functor $\Phi_\C$ maps projective contramodules to
coprojective comodules, and use 
Proposition~\ref{cotensor-contratensor-assoc}.3(a).
 Alternatively, represent the image of $\P^\bu$ in
$\sD^\ctr(\C\contra)$ by a complex of $\C/A$\projective{}
$\C$\contramodule s and the image of $\Q^\bu$ in $\sD^\ctr(\C\contra)$
by a complex of coinjective $\C$\contramodule s, and use
Proposition~\ref{cotensor-contratensor-assoc}.3(e); or
represent the image of $\P^\bu$ in $\sD^\ctr(\C\contra)$ by a complex
of $\C/A$\projective{} $\C$\contramodule s and the image of $\Q^\bu$ in
$\sD^\ctr(\C\contra)$ by a complex of $A$\injective{}
$\C$\contramodule s, and use
Proposition~\ref{cotensor-contratensor-assoc}.3(c),
Lemma~\ref{rel-inj-proj-co-contra-mod}.2(a), and the result
of~\ref{relatively-semiproj-semiinj}.
 In the third case, represent the image of $\P^\bu$ in
$\sD^\ctr(\C\contra)$ by a complex of projective $\C$\contramodule s,
notice that the functor $\Phi_\C$ maps projective contramodules to
coprojective comodules, and use 
Proposition~\ref{cotensor-contratensor-assoc}.1(a).
 Alternatively, represent the image of $\P^\bu$ in
$\sD^\ctr(\C\contra)$ by a complex of $\C/A$\projective{}
$\C$\contramodule s and the image of $\N^\bu$ in $\sD^\co(\comodr\C)$
by a complex of coflat $\C$\comodule s, and use
Proposition~\ref{cotensor-contratensor-assoc}.1(e); or
represent the image of $\P^\bu$ in $\sD^\ctr(\C\contra)$ by a complex
of $\C/A$\projective{} $\C$\contramodule s and the image of $\N^\bu$
in $\sD^\co(\comodr\C)$ by a complex of $A$\+flat $\C$\comodule s,
and use Proposition~\ref{cotensor-contratensor-assoc}.1(c),
Lemma~\ref{rel-inj-proj-co-contra-mod}.2(a), and the result
of~\ref{relatively-semiflat}.

 Finally, to show that the three pairwise isomorphisms between
the functors $\Coext_\C(\M^\bu,\P^\bu)$, \
$\Ext_\C(\M^\bu,\boL\Phi_\C(\P^\bu))$, and
$\Ext^\C(\boR\Psi_\C(\M^\bu),\P^\bu)$ form a commutative diagram,
one can represent the image of $\M^\bu$ in $\sD^\co(\C\comodl)$
by a complex of coprojective $\C$\comodule s and the image of
$\P^\bu$ in $\sD^\ctr(\C\contra)$ by a complex of coinjective
$\C$\contramodule s (having in mind
Lemma~\ref{cotensor-contratensor-assoc}), and use a result
of~\ref{cotensor-contratensor-assoc}.
\end{proof}

\Section{Semimodule-Semicontramodule Correspondence}
\label{semi-correspondence-section}

\subsection{Contratensor product and semimodule/semicontramodule
homomorphisms}
 Let $\S$ be a semialgebra over a coring $\C$.

\subsubsection{}
 We would like to define the operation of contratensor product of
a right $\S$\semimodule{} and a left $\S$\semicontramodule.
 Depending on the (co)flatness and/or (co)projectivity conditions
on $\C$ and~$\S$, one can speak of $\S$\semimodule s and
$\S$\semicontramodule s with various (co)flatness and
(co)injectivity conditions imposed on them.
 In particular, when $\C$ is a projective left $A$\module{} and
either $\S$ is a coprojective left $\C$\comodule, or $\S$ is
a projective left $A$\module{} and a $\C/A$\+coflat right
$\C$\comodule{} and $A$ has a finite left homological dimension,
or $A$ is semisimple, one can consider right $\S$\semimodule s and
left $\S$\semicontramodule s with no (co)flatness or (co)injectivity
conditions imposed.
 When $\C$ is a flat right $A$\module, $\S$ is a flat right
$A$\module{} and a $\C/A$\projective{} left $\C$\comodule, and
$A$ has a finite left homological dimension, one can consider
$A$\+flat right $\S$\semimodule s and $A$\injective{} left
$\S$\semicontramodule s.
 When $\C$ is a flat right $A$\module{} and $\S$ is a coflat right
$\C$\comodule, one can consider $\C$\+coflat right $\S$\semimodule s
and $\C$\coinjective{} left $\S$\semicontramodule s.

 The \emph{contratensor product} $\bN\Ocn_\S\bP$ of a right
$\S$\semimodule{} $\bN$ and a left $\S$\semicontramodule{} $\bP$ is
a $k$\module{} defined as the cokernel of the following pair of maps
$(\bN\oc_\C\S)\ocn_\C\bP\birarrow\bN\ocn_\C\bP$.
 The first map is induced by the right $\S$\+semiaction morphism
$\bN\oc_\C\S\rarrow\bN$.
 The second map is the composition of the map
induced by the left $\S$\+semicontraaction morphism $\bP\rarrow
\Cohom_\C(\S,\bP)$ and the natural ``evaluation'' map
$\eta_\S\:(\bN\oc_\C\S)\ocn_\C\Cohom_\C(\S,\bP)\rarrow\bN\ocn_\C\bP$.

 The latter is defined in the following generality.
 Let $\C$ be a coring over a $k$\+algebra $A$ and $\D$ be a coring
over a $k$\+algebra~$B$.
 Let $\K$ be a $\C$\+$\D$\bicomodule, $\N$ be a right $\C$\comodule, 
and $\P$ be a left $\C$\contramodule.
 Suppose that the cotensor product $\N\oc_\C\K$ is endowed with
a right $\D$\comodule{} structure via the construction
of~\ref{bicomodule-cotensor} and the cohomomorphism module
$\Cohom_\C(\K,\P)$ is endowed with
a left $\D$\contramodule{} structure via the construction
of~\ref{bicomodule-cohom}.
 Then the composition of maps $(\N\oc_\C\K)\ot_B\Hom_A(\K,\P)
\rarrow\N\ot_A\K\ot_B\Hom_A(\K,\P)\rarrow\N\ot_A\P\rarrow\N\ocn_\C\P$
factorizes through the surjection $(\N\oc_\C\K)\ot_B\Hom_A(\K,\P)
\rarrow(\N\oc_\C\K)\ocn_\D\Cohom_\C(\K,\P)$, so there is a natural
map $\eta_\K\:(\N\oc_\C\K)\ocn_\D\Cohom_\C(\K,\P)\rarrow\N\ocn_\C\P$.

 Indeed, the kernel of this surjection is equal to the sum of
the difference of two maps $(\N\oc_\C\K)\ot_B\Hom_A(\K\ot_B\D\;\P)
\birarrow(\N\oc_\C\K)\ot_B\Hom_A(\K,\P)$ and the difference of two maps
$(\N\oc_\C\K)\ot_B\Hom_A(\C\ot_A\K\;\P)\birarrow(\N\oc_\C\K)\ot_B
\Hom_A(\K,\P)$.
 The difference of the first pair of maps vanishes already in
the composition with the map $(\N\oc_\C\K)\ot_B\Hom_A(\K,\P)\rarrow
\N\ot_A\P$, while the second pair of maps can be presented as
the composition of the map $(\N\oc_\C\K)\ot_B\Hom_A(\C\ot_A\K\;\P)
\rarrow\N\ot_A\Hom_A(\C,\P)$ and the pair of maps $\N\ot_A\Hom_A(\C,\P)
\birarrow\N\ot_A\P$ whose cokernel is, by the definition,
$\N\ocn_\C\P$.
 The ``evaluation'' map $\eta_\K$ is dual to the map
\begin{setlength}{\multlinegap}{0pt}
\begin{multline*}
  \Hom_k(\eta_\K,k\dual) = \Cohom_\C(\K,{-})\: \\
  \Hom^\C(\P,\Hom_k(\N,k\dual))
  \lrarrow \Hom^\D(\Cohom_\C(\K,\P),\Cohom_\C(\K,\Hom_k(\N,k\dual))).
\end{multline*}

\subsubsection{}
 The operation of contratensor product over~$\S$ is dual to
homomorphisms in the category of left $\S$\semicontramodule s:
for any right $\S$\semimodule{} $\bN$ and any left
$\S$\semicontramodule{} $\bP$ there is a natural isomorphism
$\Hom_k(\bN\Ocn_\S\bP\;k\dual)\simeq\Hom^\S(\bP,\Hom_k(\bN,k\dual))$.
 Indeed, both $k$\module s are isomorphic to the kernel of the same
pair of maps $\Hom^\C(\bP,\Hom_k(\bN,k\dual))\birarrow
\Hom^\C(\bP,\Cohom_\C(\S,\Hom_k(\bN,k\dual)))$.
 It follows that for any right $\C$\comodule{} $\R$ for which
the induced right $\S$\semimodule{} $\R\oc_\C\S$ is defined and any
left $\S$\semicontramodule{} $\bP$ the composition of the map
$(\R\oc_\C\S)\ocn_\C\bP\rarrow(\R\oc_\C\S)\ocn_\C\Cohom_\C(\S,\bP)$
induced by the $\S$\+semicontraaction in~$\bP$ with
the ``evaluation'' map $(\R\oc_\C\S)\ocn_\C\Cohom_\C(\S,\bP)\rarrow
\R\ocn_\C\bP$ induces a natural isomorphism 
$(\R\oc_\C\S)\Ocn_\S\bP\simeq\R\ocn_\C\bP$.
\end{setlength}

 When $\C$ is a projective left $A$\module{} and $\S$ is a coprojective
left $\C$\comodule, the functor of contratensor product over~$\S$ is
right exact in both its arguments.

\subsubsection{}    \label{bisemimodule-contratensor}
 Let $\S$ be a semialgebra over a coring~$\C$ and $\T$ be
a semialgebra over a coring~$\D$.
 We would like to define a $\T$\semimodule{} structure on
the contratensor product of a $\T$\+$\S$\bisemimodule{} and
an $\S$\semicontramodule, and an $\S$\semicontramodule{} structure
on semimodule{} homomorphisms from a $\T$\+$\S$\bisemimodule{}
to a $\T$\semimodule.

 Let $\N$ be a right $\D$\comodule, $\bK$ be $\D$\+$\C$\bicomodule{}
with a right $\S$\semimodule{} structure such that the multiple
cotensor products $\N\oc_\D\bK\oc_\C\S\oc_\C\dsb\oc_\C\S$ are
associative and the semiaction map $\bK\oc_\C\S \rarrow\bK$ is a left
$\D$\comodule{} morphism, and $\bP$ be a left $\S$\semicontramodule.
 Then the contratensor product $\bK\Ocn_\S\bP$ has a natural left
$\D$\comodule{} structure as the cokernel of a pair of $\D$\comodule{}
morphisms $(\bK\oc_\C\S)\ocn_\C\bP\birarrow\bK\ocn_\C\bP$.
 The composition of maps $(\N\oc_\D\bK)\ocn_\C\bP\rarrow
\N\oc_\D(\bK\ocn_\C\bP)\rarrow\N\oc_\D(\bK\Ocn_\S\bP)$ factorizes
through the surjection $(\N\oc_\D\bK)\ocn_\C\bP\rarrow
(\N\oc_\D\bK)\Ocn_\S\bP$, so there is a natural map
$(\N\oc_\D\bK)\Ocn_\S\bP\rarrow\N\oc_\D(\bK\Ocn_\S\bP)$.

 Indeed, the composition of the pair of maps $(\N\oc_\D\bK\oc_\C\S)
\ocn_\C\bP\birarrow(\N\oc_\D\bK)\ocn_\C\bP$ whose cokernel is,
by the definition, $(\N\oc_\D\bK)\Ocn_\S\bP$, with the map
$(\N\oc_\D\bK)\ocn_\C\bP\rarrow\N\oc_\D(\bK\ocn_\C\bP)$ is equal to
the composition of the map $(\N\oc_\D\bK\oc_\C\S)\ocn_\C\bP\rarrow
\N\oc_\D((\bK\oc_\C\S)\ocn_\C\bP)$ with the pair of maps
$\N\oc_\D((\bK\oc_\C\S)\ocn_\C\bP)\birarrow\N\oc_\D(\bK\ocn_\C\bP)$.

 Now let $\bK$ be a $\T$\+$\S$\bisemimodule{} and $\bP$ be a left
$\S$\semicontramodule.
 Assume that the multiple cotensor products $\T\oc_\D\dsb\oc_\D\T
\oc_\D(\bK\Ocn_\S\bP)$ are associative and the $\D$\comodule{}
morphisms $(\T^{\suboc\.m}\oc_\D\bK)\Ocn_\S\bP\rarrow
\T^{\suboc\.m}\oc_\D(\bK\Ocn_\S\bP)$ are isomorphisms for $m\le2$.
 Then one can define an associative and unital semiaction morphism
$\T\oc_\D(\bK\Ocn_\S\bP)\rarrow\bK\Ocn_\S\bP$ taking the contratensor
product over~$\S$ of the semiaction morphism $\T\oc_\D\bK
\rarrow\bK$ with the semicontramodule~$\bP$.

 Analogously, let $\L$ be a left $\C$\comodule, $\bK$ be
a $\D$\+$\C$\bicomodule{} with a left $\T$\semimodule{} structure
such that the multiple cotensor products $\T\oc_\D\dsb\oc_\D\T
\oc_\D\bK\oc_\C\L$ are associative and the semiaction map
$\T\oc_\D\bK\rarrow\K$ is a right $\C$\comodule{} morphism, and
$\bM$ be a left $\T$\semimodule.
 Then the module of homomorphisms $\Hom_\T(\bK,\bM)$ has a natural
left $\C$\contramodule{} structure as the kernel of a pair of
$\C$\contramodule{} morphisms $\Hom_\D(\bK,\bM)\birarrow
\Hom_\D(\T\oc_\D\bK\;\bM)$.
 The composition of maps $\Cohom_\C(\L,\Hom_\T(\bK,\bM))\rarrow
\Cohom_\C(\L,\Hom_\D(\bK,\bM))\rarrow\Hom_\D(\bK\oc_\C\L\;\bM)$
factorizes through the injection $\Hom_\T(\bK\oc_\C\L\;\bM)\rarrow
\Hom_\D(\bK\oc_\C\L\;\bM)$, so there is a natural map
$\Cohom_\C(\L,\Hom_\T(\bK,\bM))\rarrow\Hom_\T(\bK\oc_\C\L\;\bM)$.

 Now let $\bK$ be a $\T$\+$\S$\bisemimodule{} and $\bM$ be a left
$\T$\semimodule.
 Assume that the multiple cohomomorphisms $\Cohom_\C(\S\oc_\C\dsb
\oc_\C\S\;\Hom_\T(\bK,\bM))$ are associative and
the $\C$\contramodule{} morphisms $\Cohom_\C(\S^{\suboc\.n},
\Hom_\T(\bK,\bM))\rarrow\Hom_\T(\bK\oc_\C\S^{\suboc\.n}\;\bM)$
are isomorphisms for $n\le2$.
 Then one can define an associative and unital semicontraaction
morphism $\Hom_\T(\bK,\bM)\rarrow\Cohom_\C(\S,\Hom_\T(\bK,\bM))$
taking the $\T$\semimodule{} homomorphisms from the semiaction
morphism $\bK\oc_\C\S\rarrow\bK$ into the semimodule~$\bM$.

\subsubsection{}    \label{semimod-semicontra-adjunction}
 Let $\bM$ be a left $\T$\semimodule, $\bK$ be
a $\T$\+$\S$\bisemimodule, and $\bP$ be a left $\S$\semicontramodule.
 Assume that a left $\T$\semimodule{} structure on $\bK\Ocn_\S\bP$
and a left $\S$\semicontramodule{} structure on $\Hom_\T(\bK,\bM)$
are defined via the constuctions of~\ref{bisemimodule-contratensor}.
 Then there is a natural adjunction isomorphism
$\Hom_\T(\bK\Ocn_\S\bP\;\bM)\simeq\Hom^\S(\bP,\Hom_\T(\bK,\bM))$.

 Indeed, the module $\Hom_\T(\bK\Ocn_\S\bP\;\bM)$ is the kernel of
the pair of maps $\Hom_\D(\bK\Ocn_\S\bP\;\bM)\birarrow
\Hom_\D(\T\oc_\D\bK\Ocn_\S\bP\;\bM)$ and there is an injection
$\Hom_\D(\T\oc_\D\bK\Ocn_\S\bP\;\bM)\rarrow\Hom_\D((\T\oc_\D\bK)
\ocn_\C\bP\;\bM)$.
 The module $\Hom_\D(\bK\Ocn_\S\bP\;\bM)$ is the kernel of
the pair of maps $\Hom_\D(\bK\ocn_\C\bP\;\bM)\birarrow
\Hom_\D((\bK\oc_\C\S)\ocn_\C\bP\;\bM)$.
 There is a pair of natural maps $\Hom_\D(\bK\ocn_\C\bP\;\M)
\birarrow\Hom_\D((\T\oc_\D\bK)\ocn_\C\bP\;\bM)$ (one of which
goes through $\Hom_\D(\T\oc_\D(\bK\ocn_\C\bP)\;\bM)$) extending
the pair of maps $\Hom_\D(\bK\Ocn_\S\bP\;\bM)\birarrow
\Hom_\D(\T\oc_\D\bK\Ocn_\S\bP\;\bM)$.
 Therefore, the module $\Hom_\T(\bK\Ocn_\S\bP\;\bM)$ is isomorphic
to the intersection of the kernels of two pairs of maps
$\Hom_\D(\bK\ocn_\C\bP\;\bM)\birarrow\Hom_\D((\bK\oc_\C\S)
\ocn_\C\bP\;\bM)$ and $\Hom_\D(\bK\ocn_\C\bP\;\M)
\birarrow\Hom_\D((\T\oc_\D\bK)\ocn_\C\bP\;\bM)$.
 Analogously, the module $\Hom^\S(\bP,\Hom_\T(\bK,\bM))$ is embedded
into $\Hom^\C(\bP,\Hom_\D(\bK,\bM))$ by the composition of maps
$\Hom^\S(\bP,\Hom_\T(\bK,\bM))\rarrow\Hom^\C(\bP,\Hom_\T(\bK,\bM))
\rarrow\Hom^\C(\bP,\Hom_\D(\bK,\bM))$ and its image coincides with
the intersection of the kernels of two pairs of maps
$\Hom^\C(\bP,\Hom_\D(\bK,\bM))\birarrow\Hom^\C(\bP\;
\Hom_\D(\T\oc_\D\bK\;\bM))$ and $\Hom^\C(\bP,\Hom_\D(\bK,\bM))
\birarrow\Hom^\C(\bP\;\Hom_\D(\bK\oc_\C\S\;\bM))$.
 These are the same two pairs of maps.

 In order to obtain adjoint functors and equivalences between
specific categories of left semimodules and left semicontramodules,
we will have to prove associativity isomorphisms needed for
the constructions of~\ref{bisemimodule-contratensor} to work.

\subsection{Associativity isomorphisms}
\label{semitensor-contratensor-assoc}
 Let $\S$ be a semialgebra over a coring~$\C$ over a $k$\+algebra~$A$
and $\T$ be a semialgebra over a coring~$\D$ over a $k$\+algebra~$B$.
 The following three Propositions will be mostly applied to
the cases of $\bK=\T=\S$ or $\T=\D=\C$, \ $\bK=\S$ in the sequel.

\begin{prop1}
 Let\/ $\bN$ be a right\/ $\T$\semimodule, $\bK$ be a\/
$\T$\+$\S$\bisemimodule, and\/ $\bP$ be a left\/ $\S$\semicontramodule.
 Then there is a natural map\/ $(\bN\os_\T\bK)\Ocn_\S\bP\rarrow
\bN\os_\S(\bK\Ocn_\S\bP)$ whenever both modules are defined via
the constructions of~\ref{semitensor-associative}
and~\ref{bisemimodule-contratensor}.
 This map is an isomorphism, at least, in the following cases:
\begin{enumerate}
 \item $\D$ is a flat left $B$\module, $\C$ is a projective left
       $A$\module,  $\bP$ is a contraflat left\/ $\C$\contramodule,
       and either
       \begin{itemize}
        \item $\T$ is a coflat left\/ $\D$\comodule, $\S$ is
              a coprojective left\/ $\C$\comodule, and\/ $\bK$ as
              a right\/ $\S$\semimodule{} with a left\/ $\D$\comodule{}
              structure is induced from a\/ $\D$\+coflat\/
              $\D$\+$\C$\bicomodule, or
        \item $\T$ is a flat left $B$\module{} and a\/ $\D/B$\+coflat
              right\/ $\D$\comodule, $\S$ is a projective left
              $A$\module{} and a\/ $\C/A$\+coflat right\/
              $\C$\comodule, the ring $A$ (resp., $B$) has a finite
              left (resp., weak) homological dimension, $\bK$ as
              a right\/ $\S$\semimodule{} with a left\/ $\D$\comodule{}
              structure is induced from a $B$\+flat and\/
              $\C/A$\+coflat\/ $\D$\+$\C$\bicomodule, and\/ $\bK$ as
              a left\/ $\T$\semimodule{} with a right\/ $\C$\comodule{}
              structure is induced from a $B$\+flat\/
              $\D$\+$\C$\bicomodule, or
        \item the ring~$A$ is semisimple, the ring~$B$ is absolutely
              flat, $\bK$ as a right\/ $\S$\semimodule{} with a left\/
              $\D$\comodule{} structure is induced from a\/
              $\D$\+$\C$\bicomodule, and\/ $\bK$ as a left\/
              $\T$\semimodule{} with a right\/ $\C$\comodule{}
              structure is induced from a\/ $\D$\+$\C$\bicomodule;
       \end{itemize}
 \item $\bN$ is a flat right $B$\module, $\D$ is a flat right
       $B$\module, $\T$ is a flat right $B$\module{} and a\/
       $\D/B$\+coflat left $\D$\comodule, $\C$ is a flat right
       $A$\module, $\S$ is a flat right $A$\module{} and a\/
       $\C/A$\coprojective{} left $\C$\comodule, the ring $A$
       (resp., $B$) has a finite left (resp., weak) homological
       dimension, $\bK$ as a right\/ $\S$\semimodule{} with a left\/
       $\D$\comodule{} structure is induced from an $A$\+flat and
       $\D/B$\+coflat\/ $\D$\+$\C$\bicomodule, $\bK$ as a left\/
       $\T$\semimodule{} with a right\/ $\C$\comodule{} structure is
       induced from an $A$\+flat and\/ $\D/B$\+coflat\/
       $\D$\+$\C$\bicomodule, and\/ $\bP$ is an $A$\injective{} and\/
       $\C/A$\contraflat{} left\/ $\C$\contramodule;
 \item $\bN$ is a flat right $B$\module, $\D$ is a flat right
       $B$\module, $\T$ is a flat right $B$\module{} and a\/
       $\D/B$\+coflat left $\D$\comodule, the ring $B$ has a finite
       weak homological dimension, $\bK$ as a right\/ $\S$\semimodule{}
       with a left\/ $\D$\comodule{} structure is induced from
       an $A$\+flat\/ $\D$\+$\C$\bicomodule, $\C$ is a projective left
       $A$\module, $\bP$ is a\/ $\C/A$\contraflat{} left\/
       $\C$\contramodule,  and either
       \begin{itemize}
        \item $\S$ is a coprojective left\/ $\C$\comodule{} and
              the ring $A$ has a finite weak homological dimension, or
        \item $\S$ is a projective left $A$\module{} and a\/
              $\C/A$\+coflat right\/ $\C$\comodule, the ring $A$ has
              a finite left homological dimension, and\/ $\bK$ as
              a left\/ $\T$\semimodule{} with a right\/ $\C$\comodule{}
              structure is induced from a\/ $\D$\+$\C$\bicomodule;
       \end{itemize}
 \item $\D$ is a flat right $B$\module, $\T$ is a coflat right\/
       $\D$\comodule, $\bN$ is a coflat right\/ $\D$\comodule,
       and either
       \begin{itemize}
        \item $\C$ is a projective left $A$\module{} and\/ $\S$ is
              a coprojective left\/ $\C$\comodule, or
        \item $\C$ is a projective left $A$\module, $\S$ is
              a projective left $A$\module{} and a\/ $\C/A$\+coflat
              right\/ $\C$\comodule, the right $A$ has a finite left
              homological dimension, and\/ $\bK$ as a left\/
              $\T$\semimodule{} with a right\/ $\C$\comodule{}
              structure is induced from a\/ $\D$\+$\C$\bicomodule, or
        \item $\C$ is a flat right $A$\module, $\S$ is a flat right
              $A$\module{} and a\/ $\C/A$\coprojective{} left\/
              $\C$\comodule, the ring $A$ has a finite left homological
              dimension, $\bK$ as a left\/ $\T$\semimodule{} with
              a right\/ $\C$\comodule{} structure is induced from
              an $A$\+flat\/ $\D$\+$\C$\bicomodule, and\/ $\bP$ is
              an injective left $A$\module, or
        \item $\C$ is a flat right $A$\module, $\S$ is a coflat right\/
              $\C$\comodule, $\bK$ as a left\/ $\T$\semimodule{} with
              a right\/ $\C$\comodule{} structure is induced from
              a\/ $\C$\+coflat\/ $\D$\+$\C$\bicomodule, and\/ $\bP$ is
              a coinjective{} left\/ $\C$\contramodule.
       \end{itemize}
\end{enumerate}
\end{prop1}

\begin{proof}
 If $\bN'''\to\bN''\to\bN'\to0$ is a sequence of right
$\S$\semimodule{} morphisms which is exact in the category of
$A$\module s and remains exact after taking the cotensor product
with $\S$ over~$\C$, then for any left $\S$\semicontramodule{} $\bP$
there is an exact sequence $\bN'''\Ocn_\S\bP\rarrow\bN''\Ocn_\S\bP
\rarrow\bN'\Ocn_\S\bP\rarrow0$.
 Hence whenever a right $\S$\semimodule{} structure on $\bN\os_\T\bK$
is defined via the construction of~\ref{semitensor-associative},
the $k$\module{} $(\bN\os_\T\bK)\Ocn_\S\bP$ is the cokernel of the pair
of maps $(\bN\oc_\D\T\oc_\D\bK)\Ocn_\S\bP\birarrow(\bN\oc_\D\bK)
\Ocn_\S\bP$.
 By the definition, the semitensor product $\bN\os_\T(\bK\Ocn_\S\bP)$
is the cokernel of the pair of maps $\bN\oc_\D\T\oc_\D(\bK\Ocn_\S\bP)
\birarrow\bN\oc_\D(\bK\Ocn_\S\bP)$.
 There are natural maps $(\bN\oc_\D\bK)\Ocn_\S\bP\rarrow
\bN\oc_\D(\bK\Ocn_\S\bP)$ and $(\bN\oc_\D\T\oc_\D\bK)\Ocn_\S\bP
\rarrow\bN\oc_\D\T\oc_\D(\bK\Ocn_\S\bP)$ constructed
in~\ref{bisemimodule-contratensor}.
 Whenever the left $\T$\semimodule{} structure on $\bK\Ocn_\S\bP$
is defined via the construction of~\ref{bisemimodule-contratensor},
the corresponding (two) square diagrams commute.
 So there is a natural map $(\bN\os_\T\bK)\Ocn_\S\bP\rarrow
\bN\os_\T(\bK\Ocn_\S\bP)$, which is an isomorphism provided that
the map $(\bN\oc_\D\bK)\Ocn_\S\bP\rarrow\bN\oc_\D(\bK\Ocn_\S\bP)$
and the analogous map for $\bN\oc_\D\T$ in place of $\bN$ are
isomorphisms; and the left $\T$\semimodule{} structure on
$\bK\Ocn_\S\bP$ is defined provided that the analogous map for $\T$
in place of $\bN$ is an isomorphism.
 It is straightforward to check that in each case~(a-d) a right
$\S$\semimodule{} structure on $\bN\os_\T\bK$ is defined via
the construction of~\ref{semitensor-associative} (that is where
the conditions that $\bK$ as a left $\T$\semimodule{} with a right
$\C$\comodule{} structure is induced from a $\D$\+$\C$\bimodule{}
are used).
 It is also easy to verify the (co)flatness conditions on
$\bK\Ocn_\S\bP$ that are needed to guarantee that the semitensor
product $\bN\os_\T(\bK\Ocn_\S\bP)$ is defined in the case~(a).
 Thus it remains to show that the map $(\bN\oc_\D\bK)\Ocn_\S\bP
\rarrow\bN\oc_\D(\bK\Ocn_\S\bP)$ is an isomorphism.

 In the case~(d), the map $(\bN\oc_\D\bK)\ocn_\C\bP\rarrow
\bN\oc_\D(\bK\ocn_\C\bP)$ and the analogous map for $\bK\oc_\C\S$
in place of $\bK$ are isomorphisms by
Proposition~\ref{cotensor-contratensor-assoc}.1(e) and
the module $\bN\oc_\D(\bK\Ocn_\S\bP)$ is the cokernel of the pair
of maps $\bN\oc_\D((\bK\oc_\C\S)\ocn_\C\bP)\birarrow
\bN\oc_\D(\bK\ocn_\C\bP)$, so it is clear from the construction
of the map $(\bN\oc_\D\bK)\Ocn_\S\bP \rarrow\bN\oc_\D(\bK\Ocn_\S\bP)$
that it is an isomorphism.
 In the cases~(a-c), one has $\bK\simeq\K\oc_\C\S$ and
the multiple cotensor products $\bN\oc_\D\K\oc_\C\S\oc_\C\dsb\oc_\C\S$
are associative.
 So the map $(\bN\oc_\D\bK)\Ocn_\S\bP\rarrow\bN\oc_\D(\bK\Ocn_\S\bP)$
is naturally isomorphic to the map $(\bN\oc_\D\K)\ocn_\C\bP
\rarrow\bN\oc_\D(\K\ocn_\C\bP)$.
 The latter is an isomorphism by
Proposition~\ref{cotensor-contratensor-assoc}.1(a) in the case~(a)
and by Proposition~\ref{cotensor-contratensor-assoc}.1(d) in
the cases~(b-c).
\end{proof}

\begin{prop2}
 Let\/ $\bL$ be a left\/ $\C$\semimodule, $\bK$ be
a\/ $\T$\+$\S$\bisemimodule, and\/ $\bM$ be a left\/ $\T$\semimodule.
 Then there is a natural map\/ $\SemiHom_\S(\bL,\Hom_\T(\bK,\bM))
\rarrow\Hom_\T(\bK\os_\S\bL\;\bM)$ whenever both modules are
defined via the constructions of~\ref{semitensor-associative}
and~\ref{bisemimodule-contratensor}.
 This map is an isomorphism, at least, in the following cases:
\begin{enumerate}
 \item $\C$ is a flat right $A$\module, $\D$ is a flat right
       $B$\module, $\bM$ is an injective left\/ $\D$\comodule,
       and either
       \begin{itemize}
        \item $\S$ is a coflat right\/ $\C$\comodule, $\T$ is a coflat
              right\/ $\D$\comodule, and\/ $\bK$ as a left\/
              $\T$\semimodule{} with a right\/ $\C$\comodule{}
              structure is induced from a\/ $\C$\+coflat\/
              $\D$\+$\C$\bicomodule, or
        \item $\S$ is a flat right $A$\module{} and a\/
              $\C/A$\coprojective{} left\/ $\C$\comodule, $\T$ is
              a flat right $B$\module{} and a\/ $\D/B$\+coflat left\/
              $\D$\comodule, the ring $A$ (resp., $B$) has a finite
              left (resp., weak) homological dimension, $\bK$ as
              a left\/ $\T$\semimodule{} with a right\/ $\C$\comodule{}
              structure is induced from an $A$\+flat and\/
              $\D/B$\+coflat\/ $\D$\+$\C$\bicomodule, and\/ $\bK$ as
              a right\/ $\S$\semimodule{} with a left\/ $\D$\comodule{}
              structure is induced from an $A$\+flat\/
              $\D$\+$\C$\bicomodule, or              
        \item the ring $A$ is semisimple, the ring $B$ is absolutely
              flat, $\bK$ as a left\/ $\T$\semimodule{} with a right\/
              $\C$\comodule{} structure is induced from a\/
              $\D$\+$\C$\bicomodule, and\/ $\bK$ as a right\/
              $\S$\semimodule{} with a left\/ $\D$\comodule{}
              structure is induced from a\/ $\D$\+$\C$\bicomodule;
       \end{itemize}
 \item $\bL$ is a projective left $A$\module, $\C$ is a projective
       left $A$\module, $\S$ is a projective left $A$\module{} and
       a\/ $\C/A$\+coflat right\/ $\C$\comodule, $\D$ is a flat
       left $B$\module, $\T$ is a flat left $B$\module{} and
       a\/ $\D/B$\+coflat right\/ $\D$\comodule, the rings $A$ and $B$
       have finite left homological dimensions, $\bK$ as a left\/
       $\T$\semimodule{} with a right\/ $\C$\comodule{} structure is
       induced from a $B$\projective{} and\/ $\C/A$\+coflat\/
       $\D$\+$\C$\bicomodule, $\bK$ as a right\/ $\S$\semimodule{} with
       a left\/ $\D$\comodule{} structure is induced from a $B$\+flat
       and\/ $\C/A$\+coflat\/ $\D$\+$\C$\bicomodule, and\/ $\bM$ is
       a $B$\+flat and\/ $\D/B$\injective{} left\/ $\D$\comodule;
 \item $\bL$ is a projective left $A$\module, $\C$ is a projective
       left $A$\module, $\S$ is a projective left $A$\module{} and
       a\/ $\C/A$\+coflat right\/ $\C$\comodule, the rings $A$ and $B$
       have finite left homological dimensions, $\bK$ as a left\/
       $\T$\semimodule{} with a right\/ $\C$\comodule{} structure is
       induced from a $B$\projective{} $\D$\+$\C$\bicomodule, $\D$
       is a flat right $B$\module, $\bM$ is a\/ $\D/B$\injective{}
       left\/ $\D$\comodule, and either
       \begin{itemize}
        \item $\T$ is a coflat right\/ $\D$\comodule, or
        \item $\T$ is a flat right $B$\module{} and
              a\/ $\D/B$\+coflat\/ left\/ $\D$\comodule, and\/ $\bK$ as
              a right\/ $\S$\semimodule{} with a left\/ $\D$\comodule{}
              structure is induced from a\/ $\D$\+$\C$\bicomodule;
       \end{itemize}
 \item $\C$ is a projective left $A$\module, $\S$ is a coprojective
       left\/ $\C$\comodule, $\bL$ is a coprojective left\/
       $\C$\comodule, and either
       \begin{itemize}
        \item $\D$ is a flat right $B$\module{} and\/ $\T$ is a coflat
              right\/ $\D$\comodule, or
        \item $\D$ is a flat right $B$\module, $\T$ is a flat right
              $B$\module{} and a\/ $\D/B$\+coflat left\/ $\D$\comodule,
              and\/ $\bK$ as a right\/ $\S$\semimodule{} with a left\/
              $\D$\comodule{} structure is induced from a\/
              $\D$\+$\C$\bicomodule, or
        \item $\D$ is a flat left $B$\module, $\T$ is a flat left
              $B$\module{} and a\/ $\D/B$\+coflat right $\D$\comodule,
              $\bK$ as a right\/ $\S$\semimodule{} with a left\/
              $\D$\comodule{} structure is induced from a $B$\+flat\/
              $\D$\+$\C$\bicomodule, and\/ $\bM$ is a flat left 
              $B$\module, or
        \item $\D$ is a flat left $B$\module, $\T$ is a coflat left
              $\D$\comodule, $\bK$ as a right\/ $\S$\semimodule{} with
              a left\/ $\D$\comodule{} structure is induced from
              a $\D$\+coflat\/ $\D$\+$\C$\bicomodule, and\/ $\bM$ is
              a coflat left\/ $\D$\comodule.
       \end{itemize}
\end{enumerate}
\end{prop2}

\begin{proof}
 Any sequence $\bL''\to\bL''\to\bL'\to0$ of $\T$\semimodule{}
morphisms which is exact in the category of $B$\module s and remains
exact after taking the cotensor product with~$\T$ over~$\D$ is exact
in the category of $\T$\semimodule s, i.~e., for any $\T$\semimodule{}
$\bM$ there is an exact sequence $0\rarrow\Hom_\T(\bL',\bM)\rarrow
\Hom_\T(\bL'',\bM)\rarrow\Hom_\T(\bL''',\bM)$.
 Hence whenever a left $\T$\semimodule{} structure is defined on
$\bK\os_\S\bL$ via the construction of~\ref{semitensor-associative},
the $k$\module{} $\Hom_\T(\bK\os_\S\bL\;\bM)$ is the kernel of
the pair of maps $\Hom\T(\bK\oc_\C\bL\;\bM)\birarrow\Hom_\T(\bK\oc_\C
\S\oc_\C\bL\;\bM)$.
 By the definition, the $k$\module{}
$\SemiHom_\S(\bL,\Hom_\T(\bK,\bM))$ is the kernel of the pair of maps
$\Cohom_\C(\bL,\Hom_\T(\bK,\bM))\birarrow\Cohom_\C(\S\oc_\C\bL\;
\Hom_\T(\bK,\bM))=\Cohom_\C(\bL,\.\Cohom_\C(\S,\.\Hom_\T(\bK,\bM)))$.
 There are natural maps $\Cohom_\C(\bL,\Hom_\T\allowbreak(\bK,\bM))
\rarrow\Hom_\T(\bK\oc_\C\bL\;\bM)$ and $\Cohom_\C(\S\oc_\C\bL\;
\Hom_\T(\bK,\bM))\rarrow\Hom_\T(\bK\oc_\C\S\oc_\C\bL\;\bM)$
constructed in~\ref{bisemimodule-contratensor}.
 Whenever the left $\S$\semicontramodule{} structure on
$\Hom_\T(\bK,\bM)$ is defined via the construction
of~\ref{bisemimodule-contratensor}, the corresponding (two) square
diagrams commute.
 So there is a natural map $\SemiHom_\S(\bL,\Hom_\T(\bK,\bM))
\rarrow\Hom_\T(\bK\os_\S\bL\;\bM)$, which is an isomorphism provided
that the map $\Cohom_\C(\bL,\Hom_\T(\bK,\bM))\rarrow
\Hom_\T(\bK\oc_\C\bL\;\bM)$ and the analogous map for $\S\oc_\C\bL$
in place of~$\bL$ are isomorphisms; and the left
$\S$\semicontramodule{} structure on $\Hom_\T(\bK,\bM)$ is defined
provided that the analogous map for $\S$ in place of~$\bL$ is
an isomorphism.
 It is straightforward to check that in each case~(a-d) a left
$\T$\semimodule{} structure on $\bK\os_\S\bL$ is defined via
the construction of~\ref{semitensor-associative}.
 It is also easy to veryfy (using
Proposition~\ref{cotensor-contratensor-assoc}.2(a))
the (co)injectivity conditions on $\Hom_\T(\bK,\bM)$ that are needed
to guarantee that the semihomomorphism module
$\SemiHom_\S(\bL,\Hom_\T(\bK,\bM))$ is defined in the case~(a).
 Thus it remains to show that the map $\Cohom_\C(\bL,\Hom_\T(\bK,\bM))
\rarrow\Hom_\T(\bK\oc_\C\bL\;\bM)$ is an isomorphism.

 In the case~(d), the map $\Cohom_\C(\bL,\Hom_\D(\bK,\bM))
\rarrow\Hom_\D(\bK\oc_\C\bL\;\bM)$ and the analogous map for
$\T\oc_\D\bK$ in place of $\bK$ are isomorphisms by
Proposition~\ref{cotensor-contratensor-assoc}.2(e) and the module
$\Cohom_\C(\bL,\Hom_\T(\bK,\bM))$ is the kernel of the pair of
maps $\Cohom_\C(\bL,\Hom_\T(\bK,\bM))\birarrow
\Cohom_\C(\bL\;\Hom_\T(\T\oc_\D\bK\;\bM))$, so it is clear from
the construction of the map $\Cohom_\C(\bL,\Hom_\T(\bK,\bM))
\rarrow\Hom_\T(\bK\oc_\C\bL\;\bM)$ that it is an isomorphism.
 In the cases~(a-c), one has $\bK=\T\oc_\D\K$ and multiple
cotensor products $\T\oc_\D\dsb\oc_\D\T\oc_\D\K\oc_\C\L$ are
associative.
 So the map $\Cohom_\C(\bL,\Hom_\T(\bK,\bM))\rarrow
\Hom_\T(\bK\oc_\C\bL\;\bM)$ is naturally isomorphic to the map
$\Cohom_\C(\bL,\Hom_\D(\K,\bM))\rarrow\Hom_\D(\K\oc_\C\bL\;\bM)$.
 The latter is an isomorphism by
Proposition~\ref{cotensor-contratensor-assoc}.2(a) in the case~(a)
and by ~\ref{cotensor-contratensor-assoc}.2(d) in the cases~(b-c).
\end{proof}

\begin{prop3}
 Let\/ $\bP$ be a left\/ $\S$\semicontramodule, $\bK$ be a\/
$\T$\+$\S$\bisemimodule, and\/ $\bQ$ be a left $\T$\semicontramodule.
 Then there is a natural map\/ $\SemiHom_\T(\bK\Ocn_\S\bP\;\bQ)\rarrow
\Hom^\S(\bP,\SemiHom_\T(\bK,\bQ))$ whenever both modules are defined
via the constructions of~\ref{semihom-associative}
and~\ref{bisemimodule-contratensor}.
 This map is an isomorphism, at least, in the following cases:
\begin{enumerate}
 \item $\D$ is a projective left $B$\module, $\C$ is a projective left
       $A$\module,  $\bP$ is a projective left\/ $\C$\contramodule,
       and either
       \begin{itemize}
        \item $\T$ is a coprojective left\/ $\D$\comodule, $\S$ is
              a coprojective left\/ $\C$\comodule, and\/ $\bK$ as
              a right\/ $\S$\semimodule{} with a left\/ $\D$\comodule{}
              structure is induced from a\/ $\D$\coprojective\/
              $\D$\+$\C$\bicomodule, or
        \item $\T$ is a projective left $B$\module{} and a\/
              $\D/B$\+coflat right\/ $\D$\comodule, $\S$ is
              a projective left $A$\module{} and a\/ $\C/A$\+coflat
              right\/ $\C$\comodule, the rings $A$ and $B$ have finite
              left homological dimensions, $\bK$ as a right\/
              $\S$\semimodule{} with a left\/ $\D$\comodule{}
              structure is induced from a $B$\projective{} and\/
              $\C/A$\+coflat\/ $\D$\+$\C$\bicomodule, and\/ $\bK$ as
              a left\/ $\T$\semimodule{} with a right\/ $\C$\comodule{}
              structure is induced from a $B$\projective\/
              $\D$\+$\C$\bicomodule, or
        \item the rings $A$ and $B$ are semisimple,
              $\bK$ as a right\/ $\S$\semimodule{} with a left\/
              $\D$\comodule{} structure is induced from a\/
              $\D$\+$\C$\bicomodule, and\/ $\bK$ as a left\/
              $\T$\semimodule{} with a right\/ $\C$\comodule{}
              structure is induced from a\/ $\D$\+$\C$\bicomodule;
       \end{itemize}
 \item $\D$ is a flat right $B$\module, $\T$ is a flat right
       $B$\module{} and a\/ $\D/B$\+coflat left $\D$\comodule,
       $\bQ$ is an injective left $B$\module, $\C$ is a flat right
       $A$\module, $\S$ is a flat right $A$\module{} and a\/
       $\C/A$\coprojective{} left $\C$\comodule, the rings $A$ and $B$
       have finite left homological dimensions, $\bK$ as a right\/
       $\S$\semimodule{} with a left\/ $\D$\comodule{} structure is
       induced from an $A$\+flat and $\D/B$\coprojective\/
       $\D$\+$\C$\bicomodule, $\bK$ as a left\/ $\T$\semimodule{}
       with a right\/ $\C$\comodule{} structure is
       induced from an $A$\+flat and\/ $\D/B$\coprojective\/
       $\D$\+$\C$\bicomodule, and\/ $\bP$ is a coinjective left\/
       $\C$\contramodule;
 \item $\D$ is a flat right $B$\module, $\T$ is a flat right
       $B$\module{} and a\/ $\D/B$\+coflat left $\D$\comodule,
       $\bQ$ is an injective left $B$\module, the rings $A$ and $B$
       have finite left homological dimensions, $\bK$ as a right\/
       $\S$\semimodule{} with a left\/ $\D$\comodule{} structure is
       induced from an $A$\+flat\/ $\D$\+$\C$\bicomodule, $\C$ is
       a projective left $A$\module, $\bP$ is a\/ $\C/A$\projective{}
       left\/ $\C$\contramodule, and either
       \begin{itemize}
        \item $\S$ is a coprojective left\/ $\C$\comodule, or
        \item $\S$ is a projective left $A$\module{} and a\/
              $\C/A$\+coflat right\/ $\C$\comodule, and\/ $\bK$ as
              a left\/ $\T$\semimodule{} with a right\/ $\C$\comodule{}
              structure is induced from a\/ $\D$\+$\C$\bicomodule;
       \end{itemize}
 \item $\D$ is a flat right $B$\module, $\T$ is a coflat right\/
       $\D$\comodule, $\bQ$ is a coinjective left\/ $\D$\comodule, and
       one of the conditions of the list of Proposition~1(d) holds.
\end{enumerate}
\end{prop3}

\begin{proof}
 Let $\Q$ be a left $\D$\contramodule, $\bK$ be
a $\D$\+$\C$\bicomodule{} with a right $\S$\semimodule{} structure
such that multiple cohomomorphisms $\Cohom_\D(\bK\oc_\C\S\oc_\C\dsb
\oc_\C\S\;\Q)$ are associative and the semiaction map $\bK\oc_\C\S
\rarrow\bK$ is a left $\D$\comodule{} morphism, and $\bP$ be a left
$\S$\semicontramodule.
 Then there is a natural left $\S$\semicontramodule{} structure on
the module $\Cohom_\D(\bK,\Q)$.
 The composition of maps $\Cohom_\D(\bK\Ocn_\S\bP\;\Q)\rarrow
\Cohom_\D(\bK\ocn_\C\bP\;\Q)\rarrow\Hom^\C(\bP,\Cohom_\D(\bK,\Q))$
factorizes through the injection $\Hom^\S(\bP,\Cohom_\D(\bK,\Q))
\rarrow\Hom^\C(\bP,\Cohom_\D(\bK,\Q))$, so there is a natural map
$\Cohom_\D(\bK\Ocn_\S\bP\;\Q)\rarrow\Hom^\S(\bP,\Cohom_\D(\bK,\Q))$.
 The rest of the proof is analogous to the proofs of
Propositions~1 and~2.
\end{proof}

 Assume that $\C$ is a projective left $A$\module, $\S$ is
a coprojective left $\C$\comodule, $\D$ is a flat right $B$\module,
and $\T$ is a coflat right $\D$\comodule.
 Then it follows from~\ref{semimod-semicontra-adjunction} together
with Propositions~1(d) and~2(d) that for any left $\T$\semimodule{}
$\bP$, and $\T$\+$\S$\bisemimodule{} $\bK$, and any left
$\S$\semicontramodule{} $\bP$ there is a natural isomorphism
$\Hom_\T(\bK\Ocn_\S\bP\;\bM)\simeq\Hom^\S(\bP,\Hom_\T(\bK,\bM))$.

 In particular, when $\C$ is a projective left and a flat right
$A$\module{} and $\S$ is a coprojective left and a coflat right
$\C$\comodule, there is a pair of adjoint functors $\Psi_\S\:
\S\simodl\rarrow\S\sicntr$ and $\Phi_\S\:\S\sicntr\rarrow\S\simodl$
compatible with the functors $\Psi_\C\:\C\comodl\rarrow\C\contra$
and $\Phi_\C\:\C\contra\rarrow\C\comodl$.
 In other words, the $\S$\semimodule{} $\Psi_\S(\bM)$ as
a $\C$\comodule{} is naturally isomorphic to $\Psi_\C(\bM)$ and
the $\S$\semicontramodule{} $\Phi_\S(\bP)$ as a $\C$\contramodule{}
is naturally isomorphic to $\Phi_\C(\bP)$.

 Assume that $\C$ is a projective left $A$\module{} and either
$\S$ is a coprojective left $\C$\comodule, or $\S$ is a projective
left $A$\module{} and a $\C/A$\+coflat right $\C$\comodule{} and $A$
has a finite left homological dimension.
 Then it follows from Propositions~1(a) and~2(b,d) that the categories
of $\C$\coprojective{} left $\S$\semimodule s and $\C$\projective{}
left $\S$\semicontramodule s are naturally equivalent.

 Assume that $\C$ is a flat right $A$\module{} and either $\S$ is
a coflat right $\C$\comodule, or $\S$ is a flat right $A$\module{}
and a $\C/A$\coprojective{} left $\C$\comodule{} and $A$ has a finite
left homological dimension.
 Then if follows from Propositions~1(b,d) and~2(a) that the categories
of $\C$\injective{} left $\S$\semimodule s and $\C$\coinjective{}
left $\S$\semicontramodule s are naturally equivalent.

 Assume that $\C$ is a projective left $A$\module{} and a flat right
$A$\module, $A$ has a finite left homological dimension, and either
$\S$ is a coprojective left $\C$\comodule{} and a flat right
$A$\module, or $\S$ is a projective left $A$\module{} and a coflat
right $\C$\comodule.
 Then it follows from Propositions~1(c,d) and~2(c,d) that
the categories of $\C/A$\injective{} left $\S$\semimodule s and
$\C/A$\projective{} left $\S$\semicontramodule s are naturally
equivalent.

 Finally, assume that the ring $A$ is semisimple.
 Then it follows from Propositions~1(a) and~2(a) that the categories
of $\C$\injective{} left $\S$\semimodule s and $\C$\projective{}
left $\S$\semicontramodule s are naturally equivalent.

\smallskip
 In each of these cases, the natural maps defined in Propositions~2--3
in the case of $\bK=\T=\S$ have the following property of compatibility
with the adjoint functors between categories of $\S$\semimodule s and
$\S$\semicontramodule s.
 For any left $\S$\semimodule{} $\bM$ and any
left $\S$\semicontramodule{} $\bP$ such that the $\S$\semimodule{}
$\Phi_\S(\bP)=\S\Ocn_\S\bP$, the $\S$\semicontramodule{}
$\Psi_\S(\bM)=\Hom_\S(\S,\bM)$, and the $k$\module{} of
semihomomorphisms $\SemiHom_\S(\Phi_\S(\bP),\Psi_\S(\bM))$ are
defined via the constructions of~\ref{bisemimodule-contratensor}
and~\ref{semihom-associative}, the maps
$\SemiHom_\S(\Phi_\S(\bP),\Psi_\S(\bM))\rarrow
\Hom_\S(\Phi_\S(\bP),\bM)$ and $\SemiHom_\S(\Phi_\S(\bP),\Psi_\S(\bM))
\rarrow\Hom^\S(\bP,\Psi_\S(\bM))$ form a commutative diagram with
the adjunction isomorphism $\Hom_\S(\Phi_\S(\bP),\bM)\simeq
\Hom^\S(\bP,\Psi_\S(\bM))$.

\subsection{Semimodule-semicontramodule correspondence}
\label{semimodule-semicontramodule-subsect}
 Assume that the coring $\C$ is a projective left and a flat right
$A$\module, the semialgebra $\S$ is a coprojective left
and a coflat right $\C$\comodule, and the ring $A$ has a finite
left homological dimension.

\begin{thm}
 \textup{(a)} The functor mapping the quotient category of complexes
of\/ $\C/A$\injective{} left\/ $\S$\semimodule s by the thick
subcategory of\/ $\C$\coacyclic{} complexes of\/ $\C/A$\injective\/
$\S$\semimodule s into the semiderived category of left\/
$\S$\semimodule s is an equivalence of triangulated categories. \par
 \textup{(b)} The functor mapping the quotient category of
complexes of\/ $\C/A$\projective{} left\/ $\S$\semicontramodule s by
the thick subcategory of\/ $\C$\contraacyclic{} complexes of\/
$\C/A$\projective\/ $\S$\semicontramodule s into the semiderived
category of left\/ $\S$\semicontramodule s is an equivalence of
triangulated categories.
\end{thm}

\begin{proof}
 Part~(b) follows from Lemma~\ref{rel-inj-proj-co-contra-mod}.2(b)
and Lemma~\ref{semitor-main-theorem} applied to the construction of
the morphism of complexes $\boL_2(\bP^\bu)\rarrow\bP^\bu$ from
the proof of Theorem~\ref{semiext-main-theorem}(b).
 As an alternative to using Lemma~\ref{rel-inj-proj-co-contra-mod}.2,
one can show that $\boL_2(\bP^\bu)$ is a complex of $\C/A$\projective{}
$\S$\semicontramodule s in the following way.
 Use Lemma~\ref{proj-inj-semi-mod-contra}(b) to construct a finite
right $A$\injective{} resolution of every term of the complex of left
$\S$\semicontramodule s $\bP^\bu$, then apply the functor $\boL_2$,
which maps exact triples of complexes to exact triples, and use
Lemmas~\ref{coproj-coinj-semi-mod-contra}(c),
\ref{cotensor-contratensor-assoc}(b),
and~\ref{rel-inj-proj-co-contra-mod}.1(b).
 The proof of part~(a) is completely analogous.
\end{proof}

\begin{rmk}
 The analogue of Theorem for complexes of quite $\C/A$\injective{}
$\S$\semimodule s and quite $\C/A$\projective{} $\S$\semicontramodule s
is true.
 Moreover, for any complex of left $\S$\semimodule s $\bM^\bu$ there
exists a morphism from $\bM^\bu$ into a complex of $\C$\injective{}
$\S$\semimodule s with a $\C$\coacyclic{} cone, and for any complex
of left $\S$\semicontramodule s $\bP^\bu$ there exists a morphism
into $\bP^\bu$ from a complex of $\C$\projective{}
$\S$\semicontramodule s with a $\C$\contraacyclic{} cone.
 Indeed, consider the complex of $\C/A$\injective{} $\S$\semimodule s
$\Phi_\S\boL_2(\bP^\bu)$ and apply to it the construction of
the morphism of complexes $\boL_1(\bK^\bu)\rarrow\bK^\bu$
from the proof of Theorems~\ref{semitor-main-theorem}
and~\ref{semiext-main-theorem}(a).
 For any complex of $\C/A$\injective{} $\S$\semimodule s $\bK^\bu$,
the complex $\boL_1(\bK^\bu)$ is a complex of coprojective
$\S$\semimodule s by Remark~\ref{absolute-relative-coproj-coinj} and
Lemma~\ref{rel-inj-proj-co-contra-mod}.2(a) (or simply because
the class of $\C/A$\injective{} left $\C$\comodule s is closed under
extensions and any $A$\projective{} $\C/A$\injective{} left
$\C$\comodule{} is coprojective, which is easy to check).
 So the complex of $\C$\coprojective{} $\S$\semimodule s
$\boL_1(\Phi_\S\boL_2(\bP^\bu))$ maps into $\Phi_\S\boL_2(\bP^\bu)$
with a $\C$\coacyclic{} cone, hence the complex of $\C$\projective{}
$\S$\semicontramodule s $\Psi_\S\boL_1(\Phi_\S\boL_2(\bP^\bu))$ maps
into $\boL_2\bP^\bu$ and $\bP^\bu$ with $\C$\contraacyclic{} cones.
\end{rmk}

\begin{cor}
 The restrictions of the functors\/ $\Psi_\S$ and\/ $\Phi_\S$
(applied to complexes term-wise) to the homotopy category of complexes
of\/ $\C/A$\injective\/ $\S$\semimodule s and\/ $\C/A$\projective\/
$\S$\semicontramodule s define mutually inverse equivalences\/
$\boR\Psi_\S$ and\/ $\boL\Phi_\S$ between the semiderived category
of left\/ $\S$\semimodule s and the semiderived category of left\/
$\S$\semicontramodule s.
\end{cor}

\begin{proof}
 By Corollary~\ref{comodule-contramodule-subsect}, the restrictions
of functors $\Psi_\S$ and $\Phi_\S$ induce mutually inverse
equivalences between the quotient category of the homotopy category
of $\C/A$\injective{} $\S$\semimodule s by its intersection with
the thick subcategory of $\C$\coacyclic{} complexes and the quotient
category of the homotopy category of $\C/A$\projective{}
$\S$\semicontramodule s by its intersection with the thick subcategory
of $\C$\contraacyclic{} complexes.
 Thus it remains to take in account the above Theorem.
\end{proof}

\subsection{Birelatively contraflat, projective, and injective
complexes}   \label{birelatively-adjusted}
 We keep the assumptions of~\ref{semimodule-semicontramodule-subsect}.

 A complex of left $\S$\semimodule s $\bM^\bu$ is called
\emph{projective relative to\/ $\C$ relative to\/~$A$} 
($\S/\C/A$\projective) if the complex of homomorphisms over~$\S$ from
$\bM^\bu$ into any $\C$\coacyclic{} complex of $\C/A$\injective{}
$\S$\semimodule s is acyclic.
 A complex of left $\S$\semicontramodule s $\bP^\bu$ is called
\emph{injective relative to\/ $\C$ relative to\/~$A$} 
($\S/\C/A$\injective) if the complex of homomorphisms over~$\S$ into
$\bP^\bu$ from any $\C$\contraacyclic{} complex of $\C/A$\projective{}
$\S$\semicontramodule s is acyclic.

 The contratensor product $\bN^\bu\Ocn_\S\bP^\bu$ of a complex
$\bN^\bu$ of right $\S$\semimodule s and a complex $\bP^\bu$ of
left $\S$\semicontramodule s is defined as the total complex of
the bicomplex $\bN^i\Ocn_\S\bP^j$, constructed by taking infinite
direct sums along the diagonals.
 A complex of right $\S$\semimodule s $\bN^\bu$ is called
\emph{contraflat relative to\/ $\C$ relative to\/~$A$}
($\S/\C/A$\contraflat) if the contratensor product over~$\S$
of the complex $\bN^\bu$ any any $\C$\contraacyclic{} complex of
$\C/A$\projective{} left $\S$\semicontramodule s is acyclic.

 It follows from Theorem~\ref{comodule-contramodule-subsect} and
Lemma~\ref{rel-inj-proj-co-contra-mod}.1 that the complex of left
$\S$\semimodule s induced from a complex of $A$\projective{}
$\C$\comodule s is $\S/\C/A$\projective, the complex of left
$\S$\semicontramodule s coinduced from a complex of $A$\injective{}
$\C$\contramodule s is $\S/\C/A$\injective, and the complex of
right $\S$\semimodule s induced from a complex of $A$\+flat
$\C$\comodule s is $\S/\C/A$\contraflat.

\begin{lem}
 \textup{(a)} Any\/ $\S/\C/A$\semiflat{} complex of $A$\+flat right\/
$\S$\semimodule s (in the sense of~\ref{relatively-semiflat})
is\/ $\S/\C/A$\contraflat. \par
 \textup{(b)} A complex of $A$\projective{} left $\S$\semimodule s
is\/ $\S/\C/A$\projective{} if and only if it is\/
$\S/\C/A$\semiprojective{} (in the sense
of~\ref{relatively-semiproj-semiinj}). \par
 \textup{(c)} A complex of $A$\injective{} left\/
$\S$\semicontramodule s is\/ $\S/\C/A$\injective{} if and only if
it is\/ $\S/\C/A$\semiinjective{} (in the sense
of~\ref{relatively-semiproj-semiinj}).
\end{lem}

\begin{proof}
 The functors $\Psi_\S$ and $\Phi_\S$ define an equivalence between
the category of $\C$\coacyclic{} complexes of $\C/A$\injective{}
left $\S$\semimodule s and the category of $\C$\contraacyclic{}
complexes of $\C/A$\projective{} left $\S$\semicontramodule s.
 Therefore, part~(a) follows from
Proposition~\ref{semitensor-contratensor-assoc}.1(c) (applied to
$\bK=\T=\S$) and Lemma~\ref{rel-inj-proj-co-contra-mod}.2(a), part~(b)
follows from Proposition~\ref{semitensor-contratensor-assoc}.2(c)
and Lemma~\ref{rel-inj-proj-co-contra-mod}.2(b), and part~(c) follows
from Proposition~\ref{semitensor-contratensor-assoc}.3(c) and
Lemma~\ref{rel-inj-proj-co-contra-mod}.2(a).
\end{proof}

 In view of the relevant results of~\ref{relatively-semiproj-semiinj},
it is also clear that a complex of $A$\projective{} left $\S$\semimodule s
is $\S/\C/A$\projective{} if the complex of $\S$\semimodule{}
homomorphisms from it into any $\C$\contractible{} complex of
$\C/A$\injective{} $\S$\semimodule s is acyclic.
 Analogously, a complex of $A$\injective{} left $\S$\semicontramodule s
is $\S/\C/A$\injective{} if the complex of $\S$\semicontramodule{}
homomorphisms into it from any $\C$\contractible{} complex of
$\C/A$\projective{} $\S$\semicontramodule s is acyclic.

\begin{qst}
 One can show using the construction of the morphism of complexes
of left $\S$\semimodule s $\bL^\bu\rarrow\boR_2(\bL^\bu)$ and
Lemma~\ref{absolute-relative-coflat} that any $\S/\C/A$\contraflat{}
complex of (appropriately defined) $\S/\C/A$\semiflat{} right
$\S$\semimodule s is $\S/\C/A$\semiflat.
 One can also show using the functor $\SemiTor^\S$ that any $A$\+flat
$\S/\C/A$\contraflat{} right $\S$\semimodule{} (defined in terms of
exact triples of $\C/A$\projective{} or $\C/A$\contraflat{}
left $\S$\semicontramodule s) is $\S/\C/A$\semiflat; the converse
is clear (cf.~\ref{semi-model-struct}).
 Are all $\S/\C/A$\contraflat{} (in either definition) right
$\S$\semimodule s $A$\+flat?
 Are all $\S/\C/A$\contraflat{} complexes of $A$\+flat right
$\S$\semimodule s $\S/\C/A$\semiflat?
\end{qst}

 The functor mapping the quotient category of $\S/\C/A$\contraflat{}
complexes of right $\S$\semimodule s by its intersection with the thick
subcategory of $\C$\coacyclic{} complexes into the semiderived category
of right $\S$\semimodule s is an equivalence of triangulated categories,
since the complex $\boL_3\boL_1(\bK^\bu)$ is  $\S/\C/A$\contraflat{}
for any complex of right $\S$\semimodule s $\bK^\bu$.
 The analogous results for $\S/\C/A$\projective{} complexes of left
$\S$\semimodule s and $\S/\C/A$\injective{} complexes of left
$\S$\semicontramodule s follow from the corresponding results
of~\ref{relatively-semiproj-semiinj}.

\begin{rmk}
 It follows from the above Lemma and
Lemma~\ref{cotensor-contratensor-assoc} that any $\C$\coacyclic{}
semiprojective complex of $\C$\coprojective{} left $\S$\semimodule s
is contractible.
 Indeed, such a complex is simultaneously an $\S/\C/A$\projective{}
complex and a $\C$\coacyclic{} complex of $\C/A$\injective{}
$\S$\semimodule s.
 Analogously, any $\C$\contraacyclic{} semiinjective complex of
$\C$\coinjective{} left $\S$\semicontramodule s is contractible.
 Hence the homotopy category of semiprojective complexes of
$\C$\coprojective{} $\S$\semimodule s is equivalent to
the semiderived category of left $\S$\semimodule s and the homotopy
category of semiinjective complexes of $\C$\coinjective{}
$\S$\semicontramodule s is equivalent to the semiderived category
of left $\S$\semicontramodule s.
 Furthermore, it follows that the homotopy category of semiprojective
complexes of $\C$\coprojective{} $\S$\semimodule s is the minimal
triangulated subcategory containing the complexes of left
$\S$\semimodule s induced from complexes of $\C$\coprojective{}
$\C$\comodule s and closed under infinite direct sums.
 Analogously, the homotopy category of semiinjective complexes of
$\C$\coinjective{} $\S$\semicontramodule s is the minimal triangulated
subcategory containing the complexes of left $\S$\semicontramodule s
coinduced from complexes of $\C$\coinjective{} $\C$\contramodule s and
closed under infinite products.
 (Cf.~\ref{remarks-derived-semitensor-bi}.)
\end{rmk}

\subsection{Derived functor CtrTor}   \label{semi-ctrtor-definition}
 The following Lemmas provide a general approach to one-sided derived
functors of any number of arguments.
 They are essentially due to P.~Deligne~\cite{Del}.

\begin{lem1}
 Let\/ $\sH$ be a category and\/ $\sS$ be a localizing class of
morphisms in\/ $\sH$.
 Let\/ $\sP$ and\/ $\sJ$ be full subcategories of\/ $\sH$ such that
either
\begin{enumerate}
 \item the map\/ $\Hom_\sH(Q,j)$ is bijective for any object $Q\in\sP$
       and any morphism $j\in\sS\cap\sJ$, and for any object $Y\in\sH$
       there is an object $J\in\sJ$ together with a morphism
       $Y\rarrow J$ belonging to\/ $\sS$, or
 \item the map\/ $\Hom_\sH(q,J)$ is bijective for any morphism
       $q\in\sS\cap\sP$ and any object $J\in\sJ$, and for any object
       $X\in\sH$ there is an object $Q\in\sP$ together with a morphism
       $Q\rarrow X$ belonging to\/ $\sS$.
\end{enumerate}
 Then for any objects\/ $P\in\sP$ and\/ $I\in\sJ$ the natural map\/
$\Hom_\sH(P,I)\rarrow\Hom_{\sH[\sS^{-1}]}(P,I)$ is bijective.
\end{lem1}

\begin{proof}
 Part~(b): any element of $\Hom_{\sH[\sS^{-1}]}(P,I)$ can be
represented by a fraction of morphisms $P\larrow X\rarrow I$ in~$\sH$,
where the morphism $X\rarrow P$ belongs to~$\sS$.
 Choose an object $Q\in\sP$ together with a morphism $Q\rarrow X$
belonging to~$\sS$.
 Then the composition $Q\rarrow X\rarrow P$ belongs to $\sS\cap\sP$,
hence the map $\Hom_\sH(P,I)\rarrow\Hom_\sH(Q,I)$ is bijective and
there exists a morphism $P\rarrow I$ that forms a commutative triangle
with the morphisms $Q\rarrow X\rarrow P$ and $Q\rarrow X\rarrow I$.
 Obviously, this morphism $P\rarrow I$ represents the same morphism in
$\sH[\sS^{-1}]$ that the fraction $P\larrow X\rarrow I$.
 Now suppose that there are two morphisms $P\birarrow I$ in~$\sH$
whose images in $\sH[\sS^{-1}]$ coincide.
 Then there exists a morphism $X\rarrow P$ belonging to~$\sH$
which has equal compositions with the morphisms $P\birarrow I$.
 Choose an object $Q\in\sP$ together with a morphism $Q\rarrow X$
belonging to~$\sH$ again.
 The composition $Q\rarrow X\rarrow P$ has equal compositions with
the morphisms $P\birarrow I$, and since the map $\Hom_\sH(P,I)
\rarrow\Hom_\sH(Q,I)$ is bijective, our two morphisms $P\birarrow I$
are equal.
 Proof of part~(a) is dual.
\end{proof}

\begin{lem2}
 Let\/ $\sH_i$, \ $i=1,\dsc,n$ be several categories, $\sS_i$
be localizing classes of morphisms in\/ $\sH_i$, and\/ $\sF_1$ be
full subcategories of\/ $\sH_i$.
 Assume that for any object $X$ in\/ $\sH_i$ there is an object $U$
in\/ $\sF_i$ together with a morphism $U\rarrow X$ from $\sS_i$.
 Let\/ $\sK$ be a category and\/ $\Theta\:\sH_1\times\dsb
\times\sH_n\rarrow\sK$ be a functor such that the morphism\/
$\Theta(U_1,\dsc,U_{i-1},t,U_{i+1},\dsc,U_n)$ is an isomorphism for
any objects $U_j\in\sF_j$ and any morphism $t\in\sS_i\cap\sF_i$.
 Then the left derived functor\/ $\boL\Theta\:\sH_1[\sS_1^{-1}]\times
\dsb\times\sH_n[\sS_n^{-1}]\rarrow\sK$ obtained by restricting $\Theta$
to\/ $\sF_1\times\dsb\times\sF_n$ is a universal final object in
the category of all functors\/ $\Xi\:\sH_1\times\dsb\times\sH_n\rarrow
\sK$ factorizable through\/ $\sH_1[\sS_1^{-1}]\times\dsb\times
\sH_n[\sS_n^{-1}]$ and endowed with a morphism of functors\/
$\Xi\rarrow\Theta$.
\end{lem2}

\begin{proof}
 It suffices to consider a single category $\sH=\sH_1\times\dsb\times
\sH_n$ with the class of morphisms $\sS=\sS_1\times\dsb\times\sS_n$,
the full subcategory $\sF=\sF_1\times\dsb\times\sF_n$, and the functor
of one argument $\Theta\:\sH\rarrow\sK$.
 The functor $\sF[(\sS\cap\sF)^{-1}]\rarrow\sH[\sS^{-1}]$ is
an equivalence of categories by Lemma~\ref{semitor-main-theorem},
so the derived functor $\boL\Theta$ can be defined.
 For any object $X\in\sH$, choose an object $U_X\in\sF$ together with
a morphism $U_X\rarrow X$ from $\sS$; then we have the induced
morphism $\boL\Theta(X)=\Theta(U_X)\rarrow \Theta(X)$.
 For any morphism $X\rarrow Y$ in $\sH$ there exists an object $V$
in $\sF$ together with a morphism $V\rarrow U_X$ belonging to~$\sS$
and a morphism $V\rarrow U_Y$ in $\sH$ forming a commutative diagram
with the morphisms $U_X\rarrow X\rarrow Y$ and $V_X\rarrow Y$.
 So we have constructed a morphism of functors
$\boL\Theta\rarrow\Theta$.
 Now if a functor $\Xi\:\sH\rarrow\sK$ factorizable through
$\sH[\sS^{-1}]$ is endowed with a morphism of functors
$\Xi\rarrow\Theta$, then the desired morphism of functors
$\Xi\rarrow\boL\Theta$ can be obtained by restricting the morphism of
functors $\Xi\rarrow\Theta$ to the subcategory $\sF\subset\sH$.
\end{proof}

 Notice the difference between the construction of a double-sided
derived functor of two arguments in Lemma~\ref{semitor-definition}
and the construction of a left derived functor of any number of
arguments in Lemma~2.
 While in the former construction only \emph{one} of the two arguments
is resolved, and the conditions imposed on the resolutions guarantee
that the two derived functors obtained in this way coincide, in
the latter construction \emph{all} of the arguments are resolved
at once and it would not suffice to resolve only some of them.
 In other words, the construction of Lemma~\ref{semitor-definition}
only works to define \emph{balanced} double-sided derived functors,
while construction of Lemma~2 allows to define \emph{nonbalanced}
one-sided derived functors.

\smallskip
 Assume that the semialgebra $\S$ satisfies the conditions
of~\ref{semimodule-semicontramodule-subsect}.

 According to Lemma~1(a) and (the proof of)
Theorem~\ref{semimodule-semicontramodule-subsect}(a), the natural map
$\Hom_{\Hot(\S\simodl)}(\bL^\bu,\bM^\bu)\rarrow
\Hom_{\sD^\si(\S\simodl)}(\bL^\bu,\bM^\bu)$ is an isomorphism
whenever $\bL^\bu$ is a complex of $\S/\C/A$\projective{}
$\S$\semimodule s and $\bM^\bu$ is a complex of $\C/A$\injective{}
$\S$\semimodule s.
 So the functor of homomorphisms in the semiderived category of left
$\S$\semimodule s can be lifted to a functor
$$
 \Ext_\S\: \sD^\si(\S\simodl)^\op\times\sD^\si(\S\simodl)\lrarrow
 \sD(k\modl),
$$
which is defined by restricting the functor of homomorphisms of
complexes of left $\S$\semimodule s to the Carthesian product of
the homotopy category of $\S/\C/A$\projective{} complexes of
$\S$\semimodule s and the homotopy category of complexes of
$\C/A$\injective{} $\S$\semimodule s.
 By Lemma~2, this construction of the right derived functor $\Ext_\S$
does not depend on the choice of subcategories of adjusted complexes.

 Analogously, according to Lemma~1(b) and (the proof of)
Theorem~\ref{semimodule-semicontramodule-subsect}(b), the natural map
$\Hom_{\Hot(\S\sicntr)}(\bP^\bu,\bQ^\bu)\rarrow
\Hom_{\sD^\si(\S\sicntr)}(\bP^\bu,\bQ^\bu)$ is an isomorphism
whenever $\bP^\bu$ is a complex of $\C/A$\injective{}
$\S$\semicontramodule s and $\bQ^\bu$ is a complex of
$\S/\C/A$\injective{} $\S$\semicontramodule s.
 So the functor of homomorphisms in the semiderived category of left
$\S$\semicontramodule s can be lifted to a functor
$$
 \Ext^\S\: \sD^\si(\S\sicntr)^\op\times\sD^\si(\S\sicntr)\lrarrow
 \sD(k\modl),
$$
which is defined by restricting the functor of homomorphisms of
complexes of left $\S$\semicontramodule s to the Carthesian product
of the homotopy category of complexese of $\C/A$\projective{}
$\S$\semicontramodule s and the homotopy category of
$\S/\C/A$\injective{} complexes of $\S$\semicontramodule s.

 Finally, the left derived functor of contratensor product
$$
 \CtrTor^\S\:\sD^\si(\simodr\S)\times\sD^\si(\S\sicntr)\lrarrow
 \sD(k\modl)
$$
is defined by restricting the functor of contratensor product
over~$\S$ to the Carthesian product of the homotopy category of
$\S/\C/A$\contraflat{} complexes of right $\S$\semimodule s and
the homotopy category of complexes of $\C/A$\projective{} left
$\S$\semicontramodule s.
 By the definition, this restriction factorizes through
the semiderived category of left $\S$\semicontramodule s in the second
argument; let us show that it also factorizes through the semiderived
category of right $\S$\semimodule s in the first argument.
 The complex of left $\S$\semicontramodule s $\Hom_k(\bN^\bu,k\dual)$
is $\S/\C/A$\injective{} whenever a complex of right $\S$\semimodule s
$\bN^\bu$ is $\S/\C/A$\contraflat; and the complex
$\Hom_k(\bN^\bu,k\dual)$ is $\C$\contraacyclic{} whenever the complex
$\bN^\bu$ is $\C$\coacyclic.
 Hence if $\bN^\bu$ is a $\C$\coacyclic{} $\S/\C/A$\contraflat{}
complex of right $\S$\semimodule s and $\bP^\bu$ is a complex of
$\C/A$\projective{} left $\S$\semicontramodule s, then the complex
$\Hom^\S(\bP^\bu,\Hom_k(\bN^\bu,k\dual))$ is acyclic, so the complex 
$\bN^\bu\Ocn_\S\bP^\bu$ is also acyclic.
 By Lemma~2, this construction of the left derived functor $\CtrTor^\S$
does not depend on the choice of subcategories of adjusted complexes.

 Notice that the constructions of derived functors $\boR\Psi_\S$ and
$\boL\Phi_\S$ in Corollary~\ref{semimodule-semicontramodule-subsect}
are also particular cases of Lemma~2.

\begin{rmk}
 To define/compute the composition multiplication
$\Ext_\S(\bL^\bu,\bM^\bu)\ot_k^\boL\Ext_\S(\bK^\bu,\bL^\bu)\rarrow
\Ext_\S(\bK^\bu,\bM^\bu)$ it suffices to represent the images of
$\bK^\bu$, $\bL^\bu$, and $\bM^\bu$ in the semiderived category of
left $\S$\semimodule s by semiprojective complexes of
$\C$\coprojective{} $\S$\semimodule s.
 The same applies to the functor $\Ext^\S$ and semiinjective complexes
of $\C$\coinjective{} $\S$\semicontramodule s.
 Besides, one can compute the functors $\Ext_\S$, $\Ext^\S$, and
$\CtrTor^\S$ using resolutions of other kinds.
 In particular, one can use complexes of $\C$\injective{}
$\S$\semimodule s and complexes of $\C$\projective{}
$\S$\semicontramodule s (see
Remark~\ref{semimodule-semicontramodule-subsect}) together with
(appropriately defined) $\S/\C$\projective{} complexes of left
$\S$\semimodule s, $\S/\C$\injective{} complexes of left
$\S$\semicontramodule s, and $\S/\C$\contraflat{} complexes of
right $\S$\semimodule s.
 One can also compute the functor $\Ext_\S$ in terms of injective
complexes of $\S$\semimodule s (defined as complexes right orthogonal
to $\C$\coacyclic{} complexes in $\Hot(\S\simodl)$) and the functor
$\Ext^\S$ in terms of projective complexes of $\S$\semicontramodule s.
 These can be obtained by applying the functor $\Phi_\S$ to
semiinjective complexes of $\C$\coinjective{} $\S$\semicontramodule s
and the functor $\Psi_\S$ to semiprojective complexes of
$\C$\coprojective{} $\S$\semimodule s, and using
Propositions~\ref{semitensor-contratensor-assoc}.2(a)
and~\ref{semitensor-contratensor-assoc}.3(a).
 Injective complexes of $\S$\semimodule s can be also constructed using
the functor right adjoint to the forgetful functor $\S\simodl\rarrow
\C\comodl$ (see Question~\ref{semicontramodules-definition}) and
infinite products of complexes of $\S$\semimodule s; this approach
works assuming only that $\C$ is a flat right $A$\module,
$\S$ is a coflat right $\C$\comodule, and $A$ has a finite left
homological dimension.
\end{rmk}

\subsection{SemiExt and Ext, SemiTor and CtrTor}
\label{semiext-and-ext}
 We keep the assumptions of~\ref{semimodule-semicontramodule-subsect}.

{\hbadness=5000
\begin{cor}
 \textup{(a)} There are natural isomorphisms of functors\/
$\SemiExt_\S(\bM^\bu,\bP^\bu)\simeq\Ext_\S(\bM^\bu,\boL\Phi_\S
(\bP^\bu))\simeq\Ext^\S(\boR\Psi_\S(\bM^\bu),\bP^\bu)$ on
the Carthesian product of the category opposite to the semiderived
category of left\/ $\S$\semimodule s and the semiderived category of
left\/ $\S$\semicontramodule s. \par
 \textup{(b)} There is a natural isomorphism of functors\/
$\SemiTor^\S(\bN^\bu,\bM^\bu)\simeq\CtrTor^\S(\bN^\bu,\boR\Psi_\S
(\bM^\bu))$ on the Carthesian product of the semiderived category
of right\/ $\S$\semimodule s and the semiderived category of
left\/ $\S$\semimodule s.
\end{cor}}

{\hbadness=2000
\begin{proof}
 It suffices to construct natural isomorphisms $\SemiExt_\S
(\bL^\bu,\boR\Psi_\S(\bM^\bu))\simeq\Ext_\S(\bL^\bu,\bM^\bu)$, \ 
$\SemiExt_\S(\boL\Phi_\S(\bP^\bu),\bQ^\bu)\simeq\Ext^\S
(\bP^\bu,\bQ^\bu)$, and $\SemiTor^\S(\bN^\bu,\allowbreak\boL\Phi_\S
(\bP^\bu))\simeq\CtrTor^\S(\bN^\bu,\bP^\bu)$.
 In the first case, represent the image of $\bM^\bu$ in
$\sD^\si(\S\simodl)$ by a complex of $\C/A$\injective{}
$\S$\semimodule s and the image of $\bL^\bu$ in $\sD^\si(\S\simodl)$
by a semiprojective complex of $\C$\coprojective{}
$\S$\semimodule s, and use
Proposition~\ref{semitensor-contratensor-assoc}.2(d) and
Lemma~\ref{birelatively-adjusted}(b).
 Alternatively, represent the image of $\bM^\bu$ in
$\sD^\si(\S\simodl)$ by a complex of $\C/A$\injective{}
$\S$\semimodule s and the image of $\bL^\bu$ in $\sD^\si(\S\simodl)$
by an $\S/\C/A$\semiprojective{} complex of $A$\projective{}
$\S$\semimodule s (see~\ref{relatively-semiproj-semiinj}), and use
Proposition~\ref{semitensor-contratensor-assoc}.2(c),
Lemma~\ref{birelatively-adjusted}(b), and
Lemma~\ref{rel-inj-proj-co-contra-mod}.2(b).
 In the second case, represent the image of $\bP^\bu$ in
$\sD^\si(\S\sicntr)$ by a complex of $\C/A$\projective{}
$\S$\semicontramodule s and the image of $\bQ^\bu$ in
$\sD^\si(\S\sicntr)$ by a semiinjective complex of $\C$\coinjective{}
$\S$\semimodule s, and use
Proposition~\ref{semitensor-contratensor-assoc}.3(d) and
Lemma~\ref{birelatively-adjusted}(c).
 Alternatively, represent the image of $\bP^\bu$ in
$\sD^\si(\S\sicntr)$ by a complex of $\C/A$\projective{}
$\S$\semicontramodule s and the image of $\bQ^\bu$ in
$\sD^\si(\S\sicntr)$ by an $\S/\C/A$\semiinjective{} complex of
$A$\injective{} $\S$\semicontramodule s
(see~\ref{relatively-semiproj-semiinj}), and use
Proposition~\ref{semitensor-contratensor-assoc}.3(c),
Lemma~\ref{birelatively-adjusted}(c), and
Lemma~\ref{rel-inj-proj-co-contra-mod}.2(a).
 In the third case, represent the image of $\bP^\bu$ in
$\sD^\si(\S\sicntr)$ by a complex of $\C/A$\projective{}
$\S$\semicontramodule s and the image of $\bN^\bu$ in
$\sD^\si(\simodr\S)$ by a semiflat complex of $\C$\+coflat{}
$\S$\semimodule s, and use
Proposition~\ref{semitensor-contratensor-assoc}.1(d) and
Lemma~\ref{birelatively-adjusted}(a).
 Alternatively, represent the image of $\bP^\bu$ in
$\sD^\si(\S\sicntr)$ by a complex of $\C/A$\projective{}
$\S$\semicontramodule s and the image of $\bN^\bu$ in
$\sD^\si(\simodr\S)$ by an $\S/\C/A$\semiflat{} complex of
$A$\+flat{} $\S$\semimodule s (see~\ref{relatively-semiflat}),
and use Proposition~\ref{semitensor-contratensor-assoc}.1(c),
Lemma~\ref{birelatively-adjusted}(a), and
Lemma~\ref{rel-inj-proj-co-contra-mod}.2(a).

 Finally, to show that the three pairwise isomorphisms between
the functors $\SemiExt_\S(\bM^\bu,\bP^\bu)$, \
$\Ext_\S(\bM^\bu,\boL\Phi_\S(\bP^\bu))$, and
$\Ext^\S(\boR\Psi_\S(\bM^\bu),\bP^\bu)$ form a commutative diagram,
one can represent the image of $\bM^\bu$ in $\sD^\si(\S\simodl)$
by a semiprojective complex of $\C$\coprojective{} $\S$\semimodule s
and the image of $\bP^\bu$ in $\sD^\si(\S\sicntr)$ by
a semiinjective complex of $\C$\coinjective{} $\S$\semicontramodule s
(having in mind Lemmas~\ref{birelatively-adjusted}
and~\ref{cotensor-contratensor-assoc}), and use a result
of~\ref{semitensor-contratensor-assoc}.
\end{proof}}

\Section{Functoriality in the Coring}

\subsection{Compatible morphisms}
\label{co-contra-compatible-morphisms}
 Let $\C$ be a coring over a $k$\+algebra $A$ and $\D$ be a coring
over a $k$\+algebra $B$.

\subsubsection{}
 We will say that a map $\C\rarrow\D$ is compatible with a $k$\+algebra
morphism $A\rarrow B$ if the biaction maps $A\ot_k\C\ot_k A\rarrow\C$
and $B\ot_k\D\ot_k B\rarrow\D$ form a commutative diagram with
the maps $\C\rarrow\D$ and $A\ot_k\C\ot_k A\rarrow B\ot_k\D\ot_k B$
(in other words, the map $\C\rarrow\D$ is an $A$\+$A$\bimodule{}
morphism) and the comultiplication maps $\C\rarrow\C\ot_A\C$ and
$\D\rarrow\D\ot_B\D$, as well as the counit maps $\C\rarrow A$ and
$\D\rarrow B$, form commutative diagrams with the maps $A\rarrow B$, \ 
$\C\rarrow\D$, and $\C\ot_A\C\rarrow\D\ot_B\D$.
 
 Let $\C\rarrow\D$ be a map of corings compatible with a $k$\+algebra
map $A\rarrow B$.
 Let $\M$ be a left comodule over $\C$ and $\N$ be a left comodule
over $B$.
 We will say that a map $\M\rarrow\N$ is compatible with the maps
$A\rarrow B$ and $\C\rarrow\D$ if the action maps $A\ot_k\M\rarrow\M$
and $B\ot_k\N\rarrow\N$ form a commutative diagram with the maps
$\M\rarrow\N$ and $A\ot_k\M\rarrow B\ot_k\N$ (that is the map
$\M\rarrow\N$ is an $A$\module{} morphism) and the coaction maps
$\M\rarrow\C\ot_A\M$ and $\N\rarrow\D\ot_B\N$ form a commutative
diagram with the maps $\M\rarrow\N$ and $\C\ot_A\M\rarrow\D\ot_B\M$.
 Analogously, let $\P$ be a left contramodule over $\C$ and $\Q$ be
a left contramodule over $\D$.
 We will say that a map $\Q\rarrow\P$ is compatible with the maps
$A\rarrow B$ and $\C\rarrow\D$ if the action maps $\P\rarrow
\Hom_k(A,\P)$ and $\Q\rarrow\Hom_k(B,\Q)$ form a commutative diagram
with the maps $\Q\rarrow\P$ and $\Hom_k(B,\Q)\rarrow\Hom_k(A,\P)$
(that is the map $\Q\rarrow\P$ is an $A$\module{} morphism) and
the contraaction maps $\Hom_A(\C,\P)\rarrow\P$ and $\Hom_B(\D,\Q)
\rarrow\Q$ form a commutative diagram with the maps $\Q\rarrow\P$ and
$\Hom_B(\D,\Q)\rarrow\Hom_A(\C,\P)$.

 Let $\M'\rarrow\N'$ be a map from a right $\C$\comodule{} $\M'$
to a right $\D$\comodule{} $\N'$ compatible with the maps $A\rarrow B$
and $\C\rarrow\D$, and $\M''\rarrow\N''$ be a map from a left
$\C$\comodule{} $\M''$ to a left $\D$\comodule{} $\N''$ compatible
with the maps $A\rarrow B$ and $\C\rarrow\D$.
 Then there is a natural map $\M'\oc_\C\M''\rarrow\N'\oc_\D\N''$.
 Analogously, let $\M\rarrow\N$ be a map from a left $\C$\comodule{}
$\M$ to a left $\D$\comodule{} $\N$ compatible with the maps
$A\rarrow B$ and $\C\rarrow\D$, and $\Q\rarrow\P$ be a map from a left
$\D$\contramodule{} $\Q$ to a left $\C$\contramodule{} $\P$ compatible
with the maps $A\rarrow B$ and $\C\rarrow\D$.
 Then there is a natural map $\Cohom_\D(\N,\Q)\rarrow\Cohom_\C(\M,\P)$.

\subsubsection{}   \label{co-contra-pull-push}
 Let $\C\rarrow\D$ be a map of corings compatible with a $k$\+algebra
map $A\rarrow B$.
 Then there is a functor from the category of left $\C$\comodule s
to the category of left $\D$\comodule s assigning to a $\C$\comodule{}
$\M$ the $\D$\comodule{} ${}_B\M=B\ot_A\M$ with the $\D$\+coaction map
defined as the composition $B\ot_A\M\rarrow B\ot_A\C\ot_A\M\rarrow
\D\ot_A\M=\D\ot_B(B\ot_A\M)$ of the map induced by the $\C$\+coaction
in $\M$ and the map induced by the map $\C\rarrow\D$ and the left
$B$\+action in~$\D$.
 The functor $\M\mpsto\M_B$ from the category of right $\C$\comodule s
to the category of right $\D$\comodule s is defined in the analogous
way.
 Furthermore, there is a functor from the category of left
$\C$\contramodule s to the category of left $\D$\contramodule s
assigning to a $\C$\contramodule{} $\P$ the $\D$\contramodule{}
${}^B\P=\Hom_A(B,\P)$ with the contraaction map defined as
the composition $\Hom_B(\D,\Hom_A(B,\P))=\Hom_A(\D,\P)\rarrow
\Hom_A(\C\ot_A B\;\P)=\Hom_A(B,\Hom_A(\C,\P))\rarrow\Hom_A(B,\P)$
of the map induced by the map $\C\rarrow\D$ and the right $B$\+action
in $\D$ with the map induced by the $\C$\+contraaction{} in~$\P$.

 If $\C$ is a flat right $A$\module, then the functor $\M\mpsto{}_B\M$
has a right adjoint functor assigning to a left $\D$\comodule{} $\N$
the left $\D$\comodule{} ${}_\C\N=\C_B\oc_\D\N$, where $\C_B=\C\ot_A B$
is a $\C$\+$\D$\bicomodule{} with the right $\D$\comodule{} structure
provided by the above construction.
 These functors are adjoint since both $k$\module s
$\Hom_\D({}_B\M,\N)$ and $\Hom_\C(\M,{}_\C\N)$ are isomorphic to
the $k$\module{} of all maps of comodules $\M\rarrow\N$ compatible
with the maps $A\rarrow B$ and $\C\rarrow\D$.
 Without any assumptions on the coring~$\C$, the functor $\N\mpsto
{}_\C\N$ is defined on the full subcategory of left $\D$\comodule s
such that the cotensor product $\C_B\oc_\D\N$ can be endowed with
a left $\C$\comodule{} structure via the construction
of~\ref{bicomodule-cotensor}; this includes, in particular,
quasicoflat $\D$\comodule s.
 Analogously, if $\C$ is a flat left $A$\module, then the functor
$\M\mpsto\M_B$ has a right adjoint functor assigning to a right
$\D$\comodule{} $\N$ the right $\C$\comodule{} $\N_\C=\N\oc_\D{}_B\C$,
where ${}_B\C=B\ot_A\C$ is a $\D$\+$\C$\bicomodule{} with the left
$\D$\comodule{} structure provided by the above construction.

 Furthermore, if $\C$ is a projective left $A$\module, then the functor
$\P\mpsto{}^B\P$ has a left adjoint functor assigning to a left
$\D$\contramodule{} $\Q$ the left $\C$\contramodule{}
${}^\C\Q=\Cohom_\D({}_B\C,\Q)$.
 These functors are adjoint since both $k$\module s
$\Hom^\D(\Q,{}^B\P)$ and $\Hom^\C({}^\C\Q,\P)$ are isomorphic to
the $k$\module{} of all maps of contramodules $\Q\rarrow\P$
compatible with the maps $A\rarrow B$ and $\C\rarrow\D$.
 Without any assumptions on the coring~$\C$, the functor
$\Q\mpsto{}^\C\Q$ is defined on the full subcategory of left
$\D$\contramodule s such that the cohomomorphism module 
$\Cohom_\D({}_B\C,\Q)$ can be endowed with a left $\C$\contramodule{}
structure via the construction of~\ref{bicomodule-cohom}; this
includes, in particular, quasicoinjective $\D$\contramodule s.

 If $\C$ is a projective left $A$\module, then for any right
$\C$\comodule{} $\M$ and any left $\D$\contramodule{} $\Q$ there is
a natural isomorphism $\M_B\ocn_\D\Q\simeq\M\ocn_\C{}^\C\Q$.
 Indeed, both $k$\module s are isomorphic to the cokernel of the pair
of maps $\M\ot_A\Hom_B(\D,\Q)\birarrow\M\ot_A\Q$, one of which is
induced by the $\D$\+contraaction in $\Q$ and the other is
the composition of the map induced by the $\C$\+coaction in $\M$ and
the map induced by the evaluation map $\C_B\ot_B\Hom_B(\D,\Q)
\rarrow\Q$.
 This is obvious for $\M_B\ocn_\D\Q$, and in order to show this for
$\M\ocn_\C{}^\C\Q$ it suffices to represent ${}^\C\Q$ as
the cokernel of the pair of $\C$\contramodule{} morphisms
$\Hom_B({}_B\C,\Hom_B(\D,\Q))\birarrow\Hom_B({}_B\C,\Q)$.
 Without any assumptions on the coring~$\C$, there is a natural
isomorphism $\M_B\ocn_\D\Q\simeq\M\ocn_\C{}^\C\Q$ for any right
$\C$\comodule{} $\M$ and any left $\D$\contramodule{} $\Q$ for which
the $\C$\contramodule{} ${}^\C\Q=\Cohom_\D({}_B\C,\Q)$ is defined
via the construction of~\ref{bicomodule-cohom}.
%quasicoinjective left $\D$\contramodule{} $\Q$.

\subsubsection{}    \label{pull-push-cotensor}
 Let $\C\rarrow\D$ be a map of corings compatible with a $k$\+algebra
map $A\rarrow B$.

\begin{prop}
 \textup{(a)} For any left\/ $\C$\comodule\/ $\M$ and any right\/
$\D$\comodule\/ $\N$ for which the right\/ $\C$\comodule\/ $\N_\C$
is defined there is a natural map\/ $\N_\C\oc_\C\M\rarrow
\N\oc_\D{}_B\M$, which is an isomorphism, at least, when\/ $\C$ and\/
$\M$ are flat left $A$\module s or\/ $\N$ is a quasicoflat right\/
$\D$\comodule. \par
 \textup{(b)} For any left\/ $\C$\contramodule\/ $\P$ and any left\/
$\D$\comodule\/ $\N$ for which the left\/ $\C$\comodule{} ${}_\C\N$
is defined there is a natural map\/ $\Cohom_\D(\N,{}^B\P)\rarrow
\Cohom_\C({}_\C\N,\P)$, which is an isomorphism, at least, when
either\/ $\C$ is a flat right $A$\module{} and\/ $\P$ is an injective
left $A$\module, or\/ $\N$ is a quasicoprojective left\/
$\D$\comodule. \par
 \textup{(c)} For any left $\C$\comodule\/ $\M$ and any left
$\D$\contramodule\/ $\Q$ for which the left\/ $\C$\contramodule{}
${}^\C\Q$ is defined there is a natural map\/ $\Cohom_\D({}_B\M,\Q)
\rarrow\Cohom_\C(\M,{}^\C\Q)$, which is an isomorphism, at least,
when\/ $\C$ and $\M$ are projective left $A$\module s or\/ $\Q$ is
a quasicoinjective left\/ $\D$\contramodule.
\end{prop}

\begin{proof}
 Part~(a): for any left $\C$\comodule{} $\M$ and any right
$\D$\comodule{} $\N$ there are maps of comodules $\M\rarrow{}_B\M$
and $\N_\C\rarrow\N$ compatible with the maps $A\rarrow B$ and
$\C\rarrow\D$.
 So there is the induced map $\N_\C\oc_\C\M\rarrow\N\oc_\D{}_B\M$.
 On the other hand, for any left $\C$\comodule{} $\M$ there is
a natural isomorphism of left $\D$\comodule s ${}_B\M\simeq{}
_B\C\oc_\C\M$, hence $\N_\C\oc_\C\M=(\N\oc_\D{}_B\C)\oc_\C\M$ and
$\N\oc_\D{}_B\M=\N\oc_\D({}_B\C\oc_\C\M)$.
 Let us check that the maps $\N_\C\oc_\C\M\rarrow\N\oc_\D{}_B\M$, \ 
$\N_\C\oc_\C\M\rarrow\N\ot_B{}_B\C\ot_A\M$ and $\N\oc_\D{}_B\M\rarrow
\N\ot_B{}_B\C\ot_A\M$ form a commutative diagram.
 Indeed, the map $(\N\oc_\D{}_B\C)\ot_A\M\rarrow\N\ot_B{}_B\C\ot_A\M$
is equal to the composition of the map $(\N\oc_\D{}_B\C)\ot_A\M
\rarrow(\N\oc_\D{}_B\C)\ot_A\C\ot_A\M$ induced by the $\C$\+coaction
in $\N\oc_\D{}_B\C$ with the map $(\N\oc_\D{}_B\C)\ot_A\C\ot_A\M
\rarrow\N\ot_B{}_B\C\ot_A\M$ induced by the maps $\N\oc_\D{}_B\C
\rarrow\N$ and $\C\rarrow{}_B\C$; while the composition of maps
$(\N\oc_\D{}_B\C)\ot_A\M\rarrow\N\ot_B({}_B\C\oc_\C\M)\rarrow
\N\ot_B{}_B\C\ot_A\M$ is equal to the composition of the map
$(\N\oc_\D{}_B\C)\ot_A\M \rarrow(\N\oc_\D{}_B\C)\ot_A\C\ot_A\M$
induced by the $\C$\+coaction in $\M$ with the same map
$(\N\oc_\D{}_B\C)\ot_A\C\ot_A\M\rarrow\N\ot_B{}_B\C\ot_A\M$.
 It remains to apply Proposition~\ref{cotensor-associative}(d)
and~(e) with the left and right sides switched.
 The proofs of parts~(b) and~(c) are completely analogous;
the proof of~(b) uses Proposition~\ref{cohom-associative}(g,h)
and the proof of~(c) uses Proposition~\ref{cohom-associative}(f,i).
\end{proof}

\subsubsection{}    \label{pull-psi-phi}
 Let $\C\rarrow\D$ be a map of corings compatible with a $k$\+algebra
map $A\rarrow B$.
 Assume that $\C$ is a projective left and a flat right $A$\module.

 Then for any left $\D$\contramodule{} $\Q$ there is a natural morphism
of $\C$\comodule s $\Phi_\C({}^\C\Q)\rarrow{}_\C(\Phi_\D\Q)$, which is
an isomorphism, at least, when $\D$ is a flat right $B$\module{} and
$\Q$ is a $\D/B$\contraflat{} left $\D$\contramodule.
 Indeed, $\Phi_\C({}^\C\Q)=\C\ocn_\C{}^\C\Q\simeq\C_B\ocn_\D\Q$ as
a left $\C$\comodule{} and ${}_\C(\Phi_\D\Q)=\C_B\oc_\D(\D\ocn_\D\Q)$,
so it remains to apply
Proposition~\ref{cotensor-contratensor-assoc}.1(c).
 Analogously, for any left $\D$\comodule{} $\N$ there is a natural
morphism of $\C$\contramodule s ${}^\C(\Psi_\D\N)\rarrow
\Psi_\C({}_\C\N)$, which is an isomorphism, at least, when $\D$ is
a projective left $B$\module{} and $\N$ is a $\D/B$\injective{}
left $\D$\comodule.
 Indeed, $\Psi_\C({}_\C\N)=\Hom_\C(\C,{}_\C\N)\simeq\Hom_\D({}_B\C,\N)$
as a left $\C$\contramodule{} and ${}^\C(\Psi_\D\N)=
\Cohom_\D({}_B\C,\Hom_\D(\D,\N))$, so it remains to apply
Proposition~\ref{cotensor-contratensor-assoc}.2(c).

 Without any assumptions on the corings $\C$ and $\D$, there is
a natural isomorphism $\Phi_\C({}^\C\Q)\simeq{}_\C(\Phi_\D\Q)$ for
any quite $\D/B$\projective{} $\D$\contramodule{} $\Q$ and
a natural isomorphism ${}^\C(\Psi_\D\N)\simeq\Psi_\C({}_\C\N)$ for
any quite $\D/B$\injective{} $\D$\comodule{} $\N$.

 The natural morphisms $\Phi_\C({}^\C\Q)\rarrow{}_\C(\Phi_\D\Q)$ and
${}^\C(\Psi_\D\N)\rarrow\Psi_\C({}_\C\N)$ have the following
compatibility property.
 For any left $\D$\comodule{} $\N$ and left $\D$\contramodule{} $\Q$
for which the $\C$\comodule{} ${}_\C\N$ and the $\C$\contramodule{}
${}^\C\Q$ are defined via the constructions
of~\ref{bicomodule-cotensor} and~\ref{bicomodule-cohom}, for any pair
of morphisms $\Phi_\D\Q\rarrow\N$ and $\Q\rarrow\Psi_\D\N$
corresponding to each other under the adjunction of functors $\Psi_\D$
and $\Phi_\D$, the compositions $\Phi_\C({}^\C\Q)\rarrow{}_\C(\Phi_\D\Q)
\rarrow{}_\C\N$ and ${}^\C\Q\rarrow{}^\C(\Psi_\D\N)\rarrow
\Psi_\C({}_\C\N)$ correspond to each other under the adjunction
of functors $\Psi_\C$ and~$\Phi_\C$.

\subsection{Properties of the pull-back and push-forward functors}

\subsubsection{}   \label{co-contra-pull-well-defined}
 Let $\C\rarrow\D$ be a map of corings compatible with a $k$\+algebra
map $A\rarrow B$.

\begin{thm}
 \textup{(a)} Assume that\/ $\C$ is a flat right $A$\module.
 Then the functor\/ $\N\mpsto{}_\C\N$ maps\/ $\D/B$\+coflat
(\.$\D/B$\coprojective) left\/ $\D$\comodule s to\/ $\C/A$\+coflat
(\.$\C/A$\coprojective) left\/ $\C$\comodule s.
 Assume additionally that\/ $\D$ is a flat right $B$\module.
 Then the same functor applied to complexes maps coacyclic complexes
of\/ $\D/B$\+coflat\/ $\D$\comodule s to coacyclic complexes of\/
$\C$\comodule s. \par
 \textup{(b)} Assume that\/ $\C$ is a projective left $A$\module.
 Then the functor\/ $\Q\mpsto{}^\C\Q$ maps\/ $\D/B$\coinjective{}
left\/ $\D$\contramodule s to $\C/A$\coinjective{} left\/
$\C$\contramodule s.
 Assume additionally that\/ $\D$ is a projective left $B$\module.
 Then the same functor applied to complexes maps contraacyclic
complexes of\/ $\D/B$\coinjective{} $\D$\contramodule s to
contraacyclic complexes of\/ $\C$\contramodule s.
\end{thm}

\begin{proof}
 Part~(a): the first assertion follows from parts~(a)
(with the left and right sides switched) and~(b) of
Proposition~\ref{pull-push-cotensor}.
 To prove the second assertion, denote by $\K^\bu$ the cobar resolution
$\C_B\ot_B\D\rarrow\C_B\ot_B\D\ot_B\D\rarrow\dsb$ of the right
$\D$\comodule{} $\C_B$.
 Then $\K^\bu$ is a complex of $\D$\+coflat $\C$\+$\D$\bicomodule s and
the cone of the morphism $\C_B\rarrow\K^\bu$ is coacyclic with respect
to the exact category of $B$\+flat $\C$\+$\D$\bicomodule s.
 Thus if $\N^\bu$ is a coacyclic complex of left $\D$\comodule s, then
the complex of left $\C$\comodule s $\K^\bu\oc_\D\N^\bu$ is coacyclic
and if $\N^\bu$ is a complex of $\D/B$\+coflat left $\D$\comodule s,
then the cone of the morphism $\C_B\oc_\D\N^\bu\rarrow
\K^\bu\oc_\D\N^\bu$ is coacyclic.
 The proof of part~(b) is completely analogous.
\end{proof}

\subsubsection{}   \label{co-contra-push-well-defined}
 It is obvious that the functor $\M\mpsto{}_B\M$ maps complexes of
$A$\+flat $\C$\comodule s to complexes of $B$\+flat $\D$\comodule s.
 It will follow from the next Theorem that it maps coacyclic complexes
of $A$\+flat $\C$\comodule s to coacyclic complexes of $\D$\comodule s.

\begin{thm}
 \textup{(a)} Assume that the coring\/ $\C$ is a flat left and right
$A$\module{} and the ring $A$ has a finite weak homological dimension.
 Then any complex of $A$\+flat\/ $\C$\comodule s that is
coacyclic as a complex of\/ $\C$\comodule s is coacyclic with respect
to the exact category of $A$\+flat\/ $\C$\comodule s. \par
 \textup{(b)} Assume that the coring\/ $\C$ is a projective left and
a flat right $A$\module{} and the ring $A$ has a finite left
homological dimension.
 Then any complex of $A$\projective{} left\/ $\C$\comodule s
that is coacyclic as a complex of\/ $\C$\comodule s is coacyclic with
respect to the exact category of $A$\projective{} left\/
$\C$\comodule s. \par
 \textup{(c)} In the assumptions of part~(b), any complex of
$A$\injective{} left $\C$\contramodule s that is contraacyclic as
a complex of\/ $\C$\contramodule s is contraacyclic with respect to
the exact category of $A$\injective{} left\/ $\C$\contramodule s.
\end{thm}

\begin{proof}
 The proof is not difficult when $k$ is a field, as in this case
the functors of Lemmas~\ref{flat-comodule-surjection}
and~\ref{proj-inj-co-contra-module} can be made additive and exact.
 Then it follows that for any coacyclic complex of $\C$\comodule s
$\M^\bu$ the complex $\boL_1(\M^\bu)$ is coacyclic with respect to
the exact category of $A$\+flat $\C$\comodule s, while it is clear
that for any complex of $A$\+flat $\C$\comodule s $\M^\bu$ the cone
of the morphism $\boL_1(\M^\bu)\rarrow\M^\bu$ is coacyclic with
respect to the exact category of $A$\+flat $\C$\comodule s.
 Besides, parts~(b) and~(c) can be derived from the result of
Remark~\ref{co-contra-ctrtor-definition} using the cobar and bar
constructions for $\C$\comodule s and $\C$\contramodule s.
 Finally, part~(a) can be deduced from part~(b) using
Lemma~\ref{proj-inj-co-contra-module}(a), but this argument
requires stronger assumptions on~$\C$ and~$A$.

 Here is a direct proof of part~(a).
 Let us call a complex of $\C$\comodule s $m$\+flat if its terms
considered as $A$\module s have weak homological dimensions not
exceeding~$m$, and let us call an $m$\+flat complex of $\C$\comodule s
$m$\coacyclic{} if it is coacyclic with respect to the exact category
of $\C$\comodule s whose weak homological dimension over~$A$
does not exceed~$m$.
 We will show that for any $m$\coacyclic{} complex of $\C$\comodule s
$\M^\bu$ there exists an $(m-1)$\coacyclic{} complex of $\C$\comodule s
$\L^\bu$ together with a surjective morphism of complexes
$\L^\bu\rarrow\M^\bu$ whose kernel $\K^\bu$ is also $(m-1)$\coacyclic.
 It will follow that any $(m-1)$\+flat $m$\+coacyclic complex of
$\C$\comodule s $\M$ is $(m-1)$\coacyclic, since the total complex of
the exact triple $\K^\bu\to\L^\bu\to\M^\bu$ is $(m-1)$\coacyclic, as is
the cone of the morphism $\K^\bu\rarrow\L^\bu$.
 By induction we will deduce that any $0$\+flat $d$\coacyclic{}
complex of $\C$\comodule s is $0$\coacyclic, where $d$ denotes
the weak homological dimension of the ring~$A$; that is
a reformulation of the assertion~(a).
 
 Let $\M^\bu$ be the total complex of an exact triple of $m$\+flat
complexes of $\C$\comodule s ${}'\M^\bu\to{}''\M^\bu\to{}'''\M^\bu$.
 Let us choose for each degree~$n$ projective $A$\module s ${}'G^n$
and ${}'''G^n$ endowed with surjective $A$\module{} maps
${}'G^n\rarrow{}'\M^n$ and ${}'''G^n\rarrow{}'''\M^n$.
 The latter map can be lifted to an $A$\module{} map
${}'''G^n\rarrow{}''\M^n$, leading to a surjective map from the exact
triple of $A$\module s ${}'G^n\to{}'G^n\oplus{}'''G^n\to{}'''G^n$ to
the exact triple of $\C$\comodule s ${}'\M^n\to{}''\M^n\to{}'''\M^n$.
 Applying the construction of Lemma~\ref{flat-comodule-surjection},
one can obtain a surjective map from an exact triple of $A$\+flat
$\C$\comodule s ${}'\cP^n\to{}''\cP^n\to{}'''\cP^n$ to the exact triple
of $\C$\comodule s ${}'\M^n\to{}''\M^n\to{}'''\M^n$.
 Now consider three complexes of $\C$\comodule s ${}'\L^\bu$,
${}''\L^\bu$, and ${}'''\L^\bu$ whose terms are ${}^{(i)}\L^n=
{}^{(i)}\cP^{n-1}\oplus{}^{(i)}\cP^n$ and the differential
$d_{{}^{(i)}\L}^n\:{}^{(i)}\L^n\rarrow{}^{(i)}\L^{n+1}$ maps
${}^{(i)}\cP^n$ into itself by the identity map and vanishes in
the restriction to ${}^{(i)}\cP^{n-1}$ and in the projection to
${}^{(i)}\cP^{n+1}$.
 There are natural surjective morphisms of complexes ${}^{(s)}\L^\bu
\rarrow{}^{(s)}\M^\bu$ constructed as in the proof of 
Theorem~\ref{comodule-contramodule-subsect}.
 Taken together, they form a surjective map from the exact triple of
complexes ${}'\L^\bu\to{}''\L^\bu\to{}'''\L^\bu$ onto the exact triple
of complexes ${}'\M^\bu\to{}''\M^\bu\to{}'''\M^\bu$.
 Let ${}'\K^\bu\to{}''\K^\bu\to{}'''\K^\bu$ be the kernel of this map
of exact triples of complexes; then the complexes ${}^{(s)}\L^\bu$ are
$0$\+flat, while the complexes ${}^{(s)}\K^\bu$ are $(m-1)$\+flat.
 Therefore, the total complex $\L^\bu$ of the exact triple
${}'\L^\bu\to{}''\L^\bu\to{}'''\L^\bu$ is $0$\coacyclic, while
the total complex $\K^\bu$ of the exact triple
${}'\K^\bu\to{}''\K^\bu\to{}'''\K^\bu$ is $(m-1)$\coacyclic.
 There is a surjective morphism of complexes $\L^\bu\rarrow\M^\bu$
with the kernel $\K^\bu$.

 Now let ${}'\K^\bu\to{}'\L^\bu\to{}'\M^\bu$ and ${}''\K^\bu\to
{}''\L^\bu\to{}''\M^\bu$ be exact triples of complexes of
$\C$\comodule s where the complexes ${}'\K^\bu$, ${}'\L^\bu$,
${}''\K^\bu$, and ${}''\L^\bu$ are $(m-1)$\coacyclic, and suppose
that there is a morphism of complexes ${}'\M^\bu\rarrow{}''\M^\bu$.
 Let us construct for the complex $\M^\bu=\cone({}'\M^\bu\to
{}''\M^\bu)$ an exact triple of complexes $\K^\bu\to\L^\bu\to\M^\bu$
with $(m-1)$\coacyclic{} complexes $\K^\bu$ and $\L^\bu$.
 Denote by ${}'''\L^\bu$ the complex ${}'\L^\bu\oplus{}''\L^\bu$;
there is the embedding of a direct summand ${}'\L^\bu\rarrow
{}'''\L^\bu$ and the surjective morphism of complexes ${}'''\L^\bu
\rarrow{}''\M^\bu$ whose components are the composition
${}'\L^\bu\rarrow{}'\M^\bu\rarrow{}''\M^\bu$ and the surjective
morphism ${}''\L^\bu\rarrow{}''\M^\bu$.
 These two morphisms form a commutative square with the morphisms
${}'\L^\bu\rarrow{}'\M^\bu$ and ${}'\M^\bu\rarrow{}''\M^\bu$.
 The kernel ${}'''\K^\bu$ of the morphism ${}'''\L^\bu\rarrow
{}''\M^\bu$ is the middle term of an exact triple of complexes
${}''\K^\bu\rarrow{}'''\K^\bu\rarrow{}'\L^\bu$.
 Since the complexes ${}''\K^\bu$ and ${}'\L^\bu$ are
$(m-1)$\coacyclic, the complex ${}'''\K^\bu$ is also $(m-1)$\coacyclic.
 Set $\L^\bu=\cone({}'\L^\bu\rarrow{}'''\L^\bu)$ and
$\K^\bu=\cone({}'\K^\bu\rarrow{}'''\K^\bu)$; then there is an exact
triple of complexes $\K^\bu\to\L^\bu\to\M^\bu$ with the desired
properties.

 Obviously, if certain complexes of $\C$\comodule s $\M^\bu_\alpha$
can be presented as quotient complexes of $(m-1)$\coacyclic{}
complexes by $(m-1)$\coacyclic{} subcomplexes, then their direct sum
$\bigoplus\M^\bu_\alpha$ can be also presented in this way.

 Finally, let $\M^\bu\rarrow{}'\M^\bu$ be a homotopy equivalence of
$m$\+flat complexes of $\C$\comodule s, and suppose that there is
an exact triple ${}'\K^\bu\to{}'\L^\bu\to{}'\M^\bu$ with
$(m-1)$\coacyclic{} complexes ${}'\K^\bu$ and ${}'\L^\bu$.
 Let us construct an exact triple of complexes $\K^\bu\to\L^\bu\to
\M^\bu$ with $(m-1)$\coacyclic{} complexes $\K^\bu$ and $\L^\bu$.
 Consider the cone of the morphism $\M^\bu\rarrow{}'\M^\bu$; it is
contractible, and therefore isomorphic to the cone of the identity
endomorphism of a complex of $\C$\comodule s $\N^\bu$ with zero
differential.
 The complex $\N^\bu$ is $m$\+flat, so it can be presented as
the quotient complex of a complex of $A$\+flat $\C$\comodule s
$\cP^\bu$ by its $(m-1)$\+flat subcomplex $\cQ^\bu$.
 Hence the complex $\cone(\M^\bu\to{}'\M^\bu)$ is isomorphic to
the quotient complex of a $0$\+flat contractible complex
$\cone(\id_{\cP^\subbu\!\.})$ by an $(m-1)$\+flat contractible
subcomplex $\cone(\id_{\cQ^\subbu\!\.})$.

 As we have proven, for the cocone ${}''\M^\bu$ of the morphism
${}'\M^\bu\rarrow\cone(\M^\bu\to{}'\M^\bu)$ there exists an exact
triple ${}''\K^\bu\to{}''\L^\bu\to{}''\M^\bu$ with $(m-1)$\coacyclic{}
complexes ${}''\K^\bu$ and ${}''\L^\bu$.
 The complex ${}''\M^\bu$ is isomorphic to the direct sum of
the complex $\M^\bu$ and the cocone of the identity endomorphism of
the complex ${}'\M^\bu$.
 (Indeed, there is a term-wise split exact triple of complexes
$\cone(\id_{\.{}'\M^\subbu\!\.})[-1]\rarrow{}''\M^\bu\rarrow\M^\bu$
and the complex $\cone(\id_{\.{}'\M^\subbu\!\.})[-1]$ is contractible.)
 The latter cocone can be presented as the quotient complex of
an $(m-1)$\+flat contractible complex ${}'\cP^\bu$ by
an $(m-1)$\+flat contractible subcomplex ${}'\cQ^\bu$, e.~g.,
by taking ${}'\cP^\bu=\cone(\id_{\.{}'\L^\subbu\!\.})[-1]$ and
${}'\cQ^\bu=\cone(\id_{\.{}'\K^\subbu\!\.})[-1]$.

 Now suppose that there are exact triples ${}''\K^\bu\to{}''\L^\bu
\to{}''\M^\bu$ and ${}'\cQ^\bu\to{}'\cP^\bu\to{}'\N^\bu$ with
$(m-1)$\+coacyclic complexes ${}''\K^\bu$, $''\L^\bu$, ${}'\cQ^\bu$,
and ${}'\cP^\bu$ for certain complexes ${}''\M^\bu=\M^\bu\oplus
{}'\N^\bu$ and ${}'\N^\bu$.
 Let us construct an exact triple $\K^\bu\to\L^\bu\to\M^\bu$ with
$(m-1)$\+coacyclic complexes $\K^\bu$ and $\L^\bu$ (in fact, we
will have $\K^\bu={}''\K^\bu$ and our construction with obvious
modifications will work for the kernel $\M^\bu$ of a surjective
morphism of complexes ${}''\M^\bu\rarrow{}'\N^\bu$).
 Set ${}'''\M^\bu=\M^\bu\oplus{}'\cP^\bu$; then there is a surjective
morphism of complexes ${}'''\M^\bu\rarrow{}''\M^\bu$ with the kernel
${}'\cQ^\bu$.
 Let ${}'''\L^\bu$ be the fibered product of the complexes
${}'''\M^\bu$ and ${}''\L^\bu$ over ${}''\M^\bu$; then there are
exact triples of complexes ${}''\K^\bu\rarrow{}'''\L^\bu\rarrow
{}'''\M^\bu$ and ${}'\cQ^\bu\rarrow{}'''\L^\bu\rarrow{}''\L^\bu$.
 It follows from the latter exact triple that the complex ${}'''\L^\bu$
is $(m-1)$\coacyclic.
 Furthermore, there is an injective morphism of complexes
$\M^\bu\rarrow{}'''\M^\bu$ with the cokernel ${}'\cP^\bu$.
 Let $\L^\bu$ be the fibered product of the complexes $\M^\bu$ and
${}'''\L^\bu$ over ${}'''\M^\bu$; then there are exact triples of
complexes ${}''\K^\bu\rarrow\L^\bu\rarrow\M^\bu$ and $\L^\bu\rarrow
{}'''\L^\bu\rarrow{}'\cP^\bu$.
 It follows from the latter exact triple that the complex $\L^\bu$
is $(m-1)$\coacyclic.

 Part~(a) is proven; the proofs of parts~(b) and~(c) are completely
analogous.
\end{proof}

\begin{rmk}
 It follows from part~(a) of Theorem that (in the same assumptions)
any coacyclic complex of coflat $\C$\comodule s is coacyclic with
respect to the exact category of coflat $\C$\comodule s.
 Indeed, for any complex of $\C$\comodule s $\M^\bu$ coacyclic with
respect to the exact category of $A$\+flat $\C$\comodule s the complex
$\boR_2(\M^\bu)$ is coacyclic with respect to the exact category of
coflat $\C$\comodule s, and for any complex of coflat $\C$\comodule s
$\M^\bu$ the cone of the morphism $\M^\bu\rarrow\boR_2(\M^\bu)$ is
coacyclic with respect to the exact category of coflat $\C$\comodule s
(by Lemma~\ref{absolute-relative-coflat}).
 Analogously, if $\C$ is a flat right $A$\module{} then any coacyclic
complex of $\C/A$\+coflat left $\C$\comodule s is coacyclic with
respect to the exact category of $\C/A$\+coflat left $\C$\comodule s.
 For coprojective $\C$\comodule s, coinjective $\C$\contramodule s,
(quite) $\C/A$\injective{} $\C$\comodule s, and (quite)
$\C/A$\projective{} $\C$\contramodule s even stronger results
are provided by Remark~\ref{co-contra-ctrtor-definition},
Theorem~\ref{comodule-contramodule-subsect},
and Theorem~\ref{co-contra-ctrtor-definition}.
\end{rmk}

\subsection{Derived functors of pull-back and push-forward}
\label{co-contra-pull-push-derived}
 Let $\C\rarrow\D$ be a map of corings compatible with a $k$\+algebra
map $A\rarrow B$.

 Assume that $\C$ is a flat right $A$\module{} and $\D$ is a flat
right $B$\module.
 Then the functor mapping the quotient category of the homotopy
category of complexes of $\D/B$\+coflat left $\D$\comodule s by its
intersection with the thick subcategory of coacyclic complexes to
the coderived category of left $\D$\comodule s is an equivalence of
triangulated categories by Lemma~\ref{semitor-main-theorem}.
 Indeed, for any complex of left $\D$\comodule s $\N^\bu$ there is
a morphism from $\N^\bu$ into a complex of $\D/B$\+coflat
$\D$\comodule s $\boR_2(\N^\bu)$ with a coacyclic cone, which was
constructed in~\ref{cotor-main-theorem}.
 Compose the functor $\N^\bu\mpsto{}_\C\N^\bu$ acting from
the homotopy category of left $\D$\comodule s to the homotopy category
of left $\C$\comodule s with the localization functor
$\Hot(\C\comodl)\rarrow\sD^\co(\C\comodl)$ and restrict it to the full
subcategory of complexes of $\D/B$\+coflat $\D$\comodule s.
 By Theorem~\ref{co-contra-pull-well-defined}(a), this restriction
factorizes through the coderived category of left $\D$\comodule s.
 Let us denote the right derived functor so obtained by
$$
 \N^\bu\mpsto{}^\boR_\C\N^\bu\:\sD^\co(\D\comodl)\lrarrow
 \sD^\co(\C\comodl).
$$
 According to Lemma~\ref{semi-ctrtor-definition}.2, this definition of
a right derived functor does not depend on the choice of a subcategory
of adjusted complexes.

 Assume that $\C$ is a flat left and right $A$\module, $A$ has a finite
weak homological dimension, and $\D$ is a flat right $B$\module.
 Then the functor mapping the quotient category of the homotopy
category of complexes of $A$\+flat $\C$\comodule s by its intersection
with the thick subcategory of coacyclic complexes to the coderived
category of $\C$\comodule s is an equivalence of triangulated
categories by Lemma~\ref{semitor-main-theorem}.
 Indeed, for any complex of $\C$\comodule s $\M^\bu$ there is
a morphism into $\M^\bu$ from a complex of $A$\+flat $\C$\comodule s
$\boL_1(\M^\bu)$ with a coacyclic cone, which was constructed
in~\ref{cotor-main-theorem}.
 Compose the functor $\M^\bu\mpsto{}_B\M^\bu$ acting from the homotopy
category of left $\C$\comodule s to the homotopy category of left
$\D$\comodule s with the localization functor $\Hot(\D\comodl)\rarrow
\sD^\co(\D\comodl)$ and restrict it to the full subcategory of
complexes of $A$\+flat $\C$\comodule s.
 It follows from Theorem~\ref{co-contra-push-well-defined}(a) that
this restriction factorizes through the coderived category of left
$\C$\comodule s.
 Let us denote the left derived functor so obtained by
$$
 \M^\bu\mpsto{}^\boL_B\M^\bu\:\sD^\co(\C\comodl)\lrarrow
 \sD^\co(\D\comodl).
$$
 According to Lemma~\ref{semi-ctrtor-definition}.2, this definition of
a left derived functor does not depend on the choice of a subcategory
of adjusted complexes.

 Analogously, assume that $\C$ is a projective left $A$\module{} and
$\D$ is a projective left $B$\module.
 Then the left derived functor 
$$
 \Q^\bu\mpsto{}_\boL^\C\Q^\bu\:\sD^\ctr(\D\contra)\lrarrow
 \sD^\ctr(\C\contra)
$$
is defined by restricting the functor $\Q^\bu\mpsto{}^\C\Q^\bu$ to
the full subcategory of complexes of $\D/B$\coinjective{} left
$\D$\contramodule s.
 
 Assume that $\C$ is a projective left and a flat right $A$\module,
$A$ has a finite left homological dimension, and $\D$ is a projective
left $B$\module.
 Then the right derived functor
$$
 \P^\bu\mpsto{}_\boR^B\P^\bu\:\sD^\ctr(\C\contra)\lrarrow
 \sD^\ctr(\D\contra)
$$
is defined by restricting the functor $\P^\bu\mpsto{}^B\P^\bu$ to
the full subcategory of complexes of $A$\injective{} left
$\C$\contramodule s.

 Properties of the above-defined derived functors will be studied
(in the greater generality of semimodules and semicontramodules)
in Section~\ref{functoriality-II-section}.
 In particular, the functor $\N^\bu\mpsto{}^\boR_\C\N^\bu$ is right
adjoint to the functor $\M^\bu\mpsto {}^\boL_B\M^\bu$ when the latter
is defined; the functor $\Q^\bu\mpsto{}_\boL^\C\Q^\bu$ is left
adjoint to the functor $\P^\bu\mpsto{}_\boR^B\P^\bu$ when the latter
is defined; the equivalences of categories $\sD^\co(\C\comodl)\simeq
\sD^\ctr(\C\contra)$ and $\sD^\co(\D\comodl)\simeq\sD^\ctr(\D\contra)$,
when they are defined, thansform the functor $\N^\bu\mpsto
{}^\boR_\C\N^\bu$ into the functor $\Q^\bu\mpsto{}_\boL^\C\Q^\bu$;
and there are formulas connecting our derived functors with
the derived functors $\Ctrtor$, $\Cotor$ and $\Coext$.

\subsection{Faithfully flat/projective base ring change}
\label{faithful-base-ring-change}

\subsubsection{}   \label{barr-beck-theorem}
 The main ideas of the following are due to Kontsevich and
Rosenberg~\cite{KR}.

 Let $\C$ be a coring over a $k$\+algebra~$A$ and $A\rarrow B$ be
a $k$\+algebra morphism.
 The coring ${}_B\C_B$ over the $k$\+algebra $B$ is constructed
in the following way.
 As a $B$\+$B$\bimodule, ${}_B\C_B$ is equal to $B\ot_A\C\ot_A B$.
 The comultiplication in ${}_B\C_B$ is defined as the composition
$B\ot_A\C\ot_A B\rarrow B\ot_A\C\ot_A\C\ot_A B\rarrow B\ot_A\C\ot_A B
\ot_A\C\ot_A B = (B\ot_A\C\ot_A B)\ot_B (B\ot_A\C\ot_A B)$ of
the map induced by the comultiplication in~$\C$ and the map induced
by the map $A\rarrow B$.
 The counit in ${}_B\C_B$ is defined as the composition
$B\ot_A\C\ot_A B\rarrow B\ot_A B\rarrow B$ of the map induced by
the counit in~$\C$ and the map induced by the multiplication in~$B$.

 The coring ${}_B\C_B$ is a universal initial object in the category
of corings $\D$ over $B$ endowed with a map $\C\rarrow\D$ compatible
with the map $A\rarrow B$.

 As always, $B$ is called a faithfully flat right $A$\module{} if it
is a flat right $A$\module{} and for any nonzero left $A$\module{} $M$
the tensor product $B\ot_A M$ is nonzero.
 Assuming the former condition, the latter one holds if and only if
the map $M = A\ot_A M\rarrow B\ot_A M$ is injective for any left
$A$\module{} $M$.
 Therefore, $B$ is a faithfully flat right $A$\module{} if and only if
the map $A\rarrow B$ is injective and its cokernel $A/B$ is a flat
right $A$\module.
 Analogously, the ring $B$ is called a faithfully projective left
$A$\module{} if it is a projective generator of the category of left
$A$\module s, i.~e., it is a projective left $A$\module{} and
for any nonzero left $A$\module{} $P$ the module $\Hom_A(B,P)$ is
nonzero.
 Assuming the former condition, the latter one holds if and only if
the map $\Hom_A(B,P)\rarrow \Hom_A(A,P)=P$ is surjective for any
left $A$\module{} $P$.
 Therefore, $B$ is a faithfully projective left $A$\module{} if and
only if the map $A\rarrow B$ is injective and its cokernel $A/B$ is
a projective left $A$\module.

 If the coring $\C$ is a flat right $A$\module{} and the ring $B$
is a faithfully flat right $A$\module, then the functors
$\M\mpsto{}_B\M$ and $\N\mpsto{}_\C\N$ are mutually inverse
equivalences between the abelian categories of left $\C$\comodule s
and left ${}_B\C_B$\comodule s.
 Analogously, if $\C$ is a projective left $A$\module{} and $B$ is
a faithfully projective left $A$\module, then the functors
$\P\mpsto{}^B\P$ and $\Q\mpsto{}^\C\Q$ are mutually inverse
equivalences between the abelian categories of left $\C$\contramodule s
and left ${}_B\C_B$\contramodule s.
 Both assertions follow from the next general Theorem, which is
the particular case of Barr--Beck Theorem~\cite{McL} for abelian
categories and exact functors.

\begin{thm}
 If\/ $\Delta\:\sB\rarrow\sA$ is an exact functor between abelian
categories mapping nonzero objects to nonzero objects and\/
$\Gamma\:\sA\rarrow\sB$ is a functor left (resp., right) adjoint
to\/~$\Delta$, then the natural functor from the category\/ $\sB$ to
the category of modules over the monad\/ $\Delta\Gamma$ (resp.,
comodules over the comonad\/ $\Delta\Gamma$) over the category\/~$\sA$
is an equivalence of abelian categories. \qed
\end{thm}

 To prove the first assertion, it suffices to apply Theorem to
the functor $\Delta\:\allowbreak\C\comodl\rarrow B\modl$ mapping
a $\C$\comodule{} $\M$ to the $B$\module{} $B\ot_A\M$ and the functor
$\Gamma\:B\modl\rarrow \C\comodl$ right adjoint to~$\Delta$ mapping
a $B$\module{} $U$ to the $\C$\comodule{} $\C\ot_A U$.
 To prove the second assertion, apply Theorem to the functor
$\Delta\:\C\contra\rarrow B\modl$ mapping a $\C$\contramodule{} $\P$
to the $B$\module{} $\Hom_A(B,\P)$ and the functor $\Gamma\:B\modl
\rarrow\C\contra$ left adjoint to~$\Delta$ mapping a $B$\module{} $V$
to the $\C$\contramodule{} $\Hom_A(\C,V)$.

\subsubsection{}   \label{cotensor-faithful-change}
  Let $\C$ be a coring over a $k$\+algebra~$A$ and $A\rarrow B$ be
a $k$\+algebra morphism.

 Assume that $\C$ is a flat left and right $A$\module{} and $B$ is
a faithfully flat left and right $A$\module.
 Then it follows from Proposition~\ref{pull-push-cotensor}(a) that
for any right $\C$\comodule{} $\N$ and any left $\C$\comodule{} $\M$
there is a natural map $\N\oc_\C\M\rarrow\N_B\oc_{{}_B\C_B}{}_B\M$,
which is an isomorphism, at least, when one of the $A$\module s
$\N$ and $\M$ is flat or one of the ${}_B\C_B$\comodule s
$\N_B$ and ${}_B\M$ is quasicoflat.
 Analogously, assume that $\C$ is a projective left and a flat right
$A$\module{} and $B$ is a faithfully projective left and a faithfully
flat right $A$\module.
 Then it follows from Proposition~\ref{pull-push-cotensor}(b-c) that
for any left $\C$\comodule{} $\M$ and any left $\C$\contramodule{} $\P$
there is a natural map $\Cohom_{{}_B\C_B}({}_B\M,{}^B\P)\rarrow
\Cohom_\C(\M,\P)$, which is an isomorphism, at least, when
the $A$\module{} $\M$ is projective, the $A$\module{} $\P$
is injective, the ${}_B\C_B$\comodule{} ${}_B\M$ is quasicoprojective,
or the ${}_B\C_B$\contramodule{} ${}^B\P$ is quasicoinjective.

\begin{rmk}
 In general the map $\N\oc_\C\M\rarrow\N_B\oc_{{}_B\C_B}{}_B\M$ is not
an isomorphism, even under the strongest of our assumptions on~$A$,
$B$, and $\C$.
 For example, let $\C=A$ and ${}_B\C_B=B\ot_A B$; then $\N\oc_\C\M=
\N\ot_A\M$, while $\N_B\oc_{{}_B\C_B}{}_B\M$ is the kernel of the pair
of maps $\N\ot_A B\ot_A\M\birarrow\N\ot_A B\ot_A B\ot_A\M$ induced by
the map $A\rarrow B$.
 The sequence $0\rarrow N\ot_A M\rarrow N\ot_A B\ot_A M\rarrow
N\ot_A B\ot_A B\ot_A M$ is exact if one of two $A$\module s $M$
and $N$ is flat or admits a $B$\module{} structure, but in general
the map $N\ot_A M\rarrow N\ot_A B\ot_A M$ is not injective.
 Indeed, let $k$ be a field, $A=k[x]$ be the algebra of polynomials
in one variable, and $B=k[x,\d_x]$ be the algebra of differential
operators in the affine line.
 Let $M=k=N$ be one-dimensional $A$\module s with the trivial
action of~$x$.
 Then the map $N\ot_A M\rarrow N\ot_A B\ot_A M$ is zero, since
$m\ot 1\ot n= m\ot (\d_x x - x\d_x)\ot n = 0$ in $N\ot_A B\ot_A M$.
\end{rmk}

 Assume that $\C$ is a projective left and a flat right $A$\module{}
and $B$ is a faithfully projective left and a faithfully flat right
$A$\module.
 Then the equivalences between the categories $\C\comodl$ and 
${}_B\C_B\comodl$ and between the categories $\C\contra$ and
${}_B\C_B\contra$ transform the functors $\Psi_\C$ and $\Phi_\C$
into the functors $\Psi_{{}_B\C_B}$ and $\Phi_{{}_B\C_B}$.
 Indeed, one has $\Hom_{{}_B\C_B}({}_B\C_B,{}_B\M)=
\Hom_\C(\C_B,\M)=\Hom_A(B,\Hom_\C(\C,\M))$ and
${}_B\C_B\ocn_{{}_B\C_B}{}^B\P={}_B\C\ocn_\C\P=
B\ot_A(\C\ocn_\C\P)$.
 Alternatively, the same isomorphisms can be constructed as
in~\ref{pull-psi-phi} using
Propositions~\ref{cotensor-contratensor-assoc}.1(e)
and~\ref{cotensor-contratensor-assoc}.2(e).

\subsubsection{}   \label{co-contra-faithful-change-derived}
 Let $\C$ be a coring over a $k$\+algebra~$A$ and $A\rarrow B$ be
a $k$\+algebra morphism.

 Obviously, if $\C$ is a flat right $A$\module{} and $B$ is
a faithfully flat right $A$\module, then a complex of left
$\C$\comodule s $\M^\bu$ is coacyclic if any only if the complex of
left ${}_B\C_B$\comodule s ${}_B\M^\bu$ is coacyclic.
 So the functor $\M^\bu\mpsto{}_B\M^\bu$ induces an equivalence of
the coderived categories of left $\C$\comodule s and left
${}_B\C_B$\comodule s.
 If $\C$ is a projective left $A$\module{} and $B$
is a faithfully projective left $A$\module, then a complex of left
$\C$\contramodule s $\P^\bu$ is contraacyclic if and only if
the complex of ${}_B\C_B$\contramodule s ${}^B\P^\bu$ is contraacyclic.
 So the functor $\P^\bu\mpsto{}^B\P^\bu$ induces an equivalence of
the contraderived categories of left $\C$\contramodule s and
left ${}_B\C_B$\contramodule s.

 If $\C$ is a flat left and right $A$\module, $B$ is a faithfully flat
left and right $A$\module, and $A$ and $B$ have finite weak homological
dimensions, then the equivalences of categories $\sD^\co(\comodr\C)
\simeq\sD^\co(\comodr{}_B\C_B)$ and $\sD^\co(\C\comodl)\simeq
\sD^\co({}_B\C_B\comodl)$ transform the derived functor $\Cotor^\C$
into the derived functor $\Cotor^{{}_B\C_B}$.
 If $\C$ is a projective left and a flat right $A$\module,
$B$ is a faithfully projective left and a faithfully flat right
$A$\module, and $A$ and $B$ have finite left homological dimensions,
then the equivalences of categories $\sD^\co(\C\comodl)\simeq
\sD^\co({}_B\C_B\comodl)$ and $\sD^\ctr(\C\contra)\simeq
\sD^\ctr({}_B\C_B\contra)$ transform the derived functor $\Coext_\C$
into the derived functor $\Coext_{{}_B\C_B}$.
 In the same assumptions, the same equivalences of categories
transform the mutually inverse functors $\boR\Psi_\C$ and $\boL\Phi_\C$
into the mutually inverse functors $\boR\Psi_{{}_B\C_B}$ and
$\boL\Phi_{{}_B\C_B}$.
 If $\C$ is a flat right $A$\module, $B$ is a faithfully flat right
$A$\module, and $A$ and $B$ have finite left homological dimensions,
then the above equivalence of categories transforms the functor
$\Ext_\C$ into the functor $\Ext_{{}_B\C_B}$.
 If $\C$ is a projective left $A$\module, $B$ is a faithfully
projective left $A$\module, and $A$ and $B$ have finite left
homological dimensions, then the above equivalences of categories
transform the functors $\Ext^\C$ and $\Ctrtor^\C$ into the functors
$\Ext^{{}_B\C_B}$ and $\Ctrtor^{{}_B\C_B}$.

 These isomorphisms of functors can be deduced from the
uniqueness/universality assertions of Lemmas~\ref{semitor-definition}
and~\ref{semi-ctrtor-definition}.2 or derived from
the preservation/reflection results of the next Remark.
 Besides, they are particular cases of the much more general
isomorphisms constructed in Section~\ref{functoriality-II-section}.

\begin{rmk}
 In the strongest of the above flatness/projectivity and homological
dimension assumptions, almost all the properties of comodules and
contramodules over corings considered in this book are preserved by
the passages from a coring~$\C$ to the coring ${}_B\C_B$ and back.
 This applies to the properties of coflatness, coprojectivity,
coinjectivity, relative coflatness, relative coprojectivity,
relative coinjectivity, injectivity, projectivity, contraflatness,
relative injectivity, relative projectivity, relative contraflatness.
 All of this follows from the facts that an $A$\module{} $M$ is flat
if and only if the $B$\module{} $B\ot_A M$ is flat, an $A$\module{}
$M$ is projective if and only if the $B$\module{} $B\ot_A M$ is
projective, and an $A$\module{} $P$ is injective if and only if
the $B$\module{} $\Hom_A(B,P)$ is injective.
 Indeed, suppose that the left $B$\module{} $B\ot_A M$ is flat.
 Any flat left $B$\module{} is a flat left $A$\module, since the ring
$B$ is a flat left $A$\module.
 Consider the tensor product of complexes $(A\to B)\ot_A\dsb\ot_A
(A\to B)\ot_A M$, where the number of factors $A\rarrow B$ is at least
equal to the weak homological dimension of~$A$.
 This complex is exact everywhere except its rightmost term, since
the map $A\rarrow B$ is injective and $B/A$ is a flat right $A$\module.
 Since all terms of this complex, except possibly the leftmost one,
are flat left $A$\module s, the leftmost term $A$ is also a flat left
$A$\module.
 Alternatively, one can consider the complex $M\rarrow B\ot_A M\rarrow
B\ot_A B\ot_A M\rarrow\dsb$ with the alternating sums of the maps
induced by the map $A\rarrow B$ as the differentials; this complex of
left $A$\module s is acyclic, since the induced complex of left
$B$\module s is contractible.
 Notice that the assumption of finite weak homological dimension of
the ring $A$ is necessary for this argument, since otherwise the ring
$B$ can be absolutely flat while the ring $A$ is not
(see Remark~\ref{morita-change-of-coring}).
 Assuming only that $\C$ is a flat right $A$\module{} and $B$ is
a faithfully flat right $A$\module, the right ${}_B\C_B$\comodule{}
$\N_B$ is coflat if a right $\C$\comodule{} $\N$ is coflat, etc.
 On the other hand, even under the strongest of the above assumptions
there are more quite $\C/A$\injective{} $\C$\comodule s than quite
${}_B\C_B/B$\injective{} ${}_B\C_B$\comodule s and there are more
quite $\C/A$\projective{} $\C$\contramodule s than quite
${}_B\C_B/B$\projective{} ${}_B\C_B$\contramodule s; i.~e., quite
relative injectivity and quite relative projectivity is not preserved
by the equivalences of categories $\M\mpsto{}_B\M$ and
$\P\mpsto{}^B\P$ in general.
 Analogously, there are more quasicoflat $\C$\comodule s than
quasicoflat ${}_B\C_B$\comodule s.
 Indeed, consider the case when $\C=A$ and ${}_B\C_B=B\ot_A B$.
 Then all $\C$\comodule s are coinduced and all $\C$\contramodule s
are induced, while a ${}_B\C_B$\comodule{} is quite
${}_B\C_B/B$\injective, or a ${}_B\C_B$\contramodule{} is quite
${}_B\C_B/B$\projective, if and only if the corresponding
$A$\module{} is a direct summand of an $A$\module{} admitting
a $B$\module{} structure.
 For example, if $A=k[x]$ and $B=k[x,\d_x]$ as in
Remark~\ref{cotensor-faithful-change}, then the one-dimensional
$A$\module{} $M$ with the trivial action of~$x$ is not the direct
summand of any $A$\module{} admitting a $B$\module{} structure,
since the equation $xm=0$ would imply $m=-x\d_xm$.
 At the same time, any projective left $A$\module{} is a direct
summand of a projective left $B$\module{} and any injective left
$A$\module{} is a direct summand of an injective left $B$\module.
 It follows, in particular, that the cokernel of an injective morphism
of quite $\C/A$\injective{} $\C$\comodule s is not always quite
$\C/A$\injective{} and the kernel of a surjective morphism of
quite $\C/A$\projective{} $\C$\contramodule s is not always
quite $\C/A$\projective.
\end{rmk} 

\subsection{Remarks on Morita morphisms}
\label{co-contra-morita-remarks}

\subsubsection{}   
 A \emph{Morita morphism} from a $k$\+algebra $A$ to a $k$\+algebra $B$
is an $A$\+$B$\bimodule{} $E$ such that $E$ is a finitely generated
projective right $B$\module.
 For any Morita morphism $E$ from $A$ to $B$, set $E\dual=
\Hom_{B^\op}(E,B)$; then $E\dual$ is a $B$\+$A$\bimodule{} and
a finitely generated projective left $B$\module.
 To any $k$\+algebra morphism $A\rarrow B$, one can assign a Morita
morphism $E=B=E\dual$ from $A$ to $B$.
 
 Equivalently, a Morita morphism from $A$ to $B$ can be defined as
a pair consisting of an $A$\+$B$\bimodule{} $E$ and
a $B$\+$A$\bimodule{} $E\dual$ endowed with an $A$\+$A$\bimodule{}
morphism $A\rarrow E\ot_B E\dual$ and a $B$\+$B$\bimodule{} morphism
$E\dual\ot_A E\rarrow B$ such that the two compositions
$E\rarrow E\ot_B E\dual\ot_A E\rarrow E$ and $E\dual\rarrow
E\dual\ot_A E\ot_B E\dual\rarrow E\dual$ are equal to the identity
endomorphisms of $E$ and $E\dual$.

 For any Morita morphism $E$ from $A$ to $B$ the functor $N\mpsto{}_A N
=E\ot_B N=\Hom_B(E\dual,N)$ from the category of left $B$\module s
to the category of left $A$\module s has a left adjoint functor
$M\mpsto{}_B M=E\dual\ot_A M$ and a right adjoint functor $P\mpsto
{}^B P=\Hom_A(E,P)$.
 Analogously, the functor $N\mpsto N_A=N\ot_B E\dual=\Hom_{B^\op}(E,N)$
from the category of right $B$\module s to the category of right
$A$\module s has a left adjoint functor $M\mpsto M_B=M\ot_A E$ and
a right adjoint functor $P\mpsto P^B=\Hom_{B^\op}(E\dual,P)$.

 Let $\C$ be a coring over a $k$\+algebra $A$ and $E$ be a Morita
morphism from $A$ to $B$.
 Then there is a coring structure on the $B$\+$B$\bimodule{}
${}_B\C_B=E\dual\ot_A\C\ot_A E$ defined in the following
way~\cite{BKG}.
 The comultiplication in ${}_B\C_B$ is the composition $E\dual\ot_A\C
\ot_A E\rarrow E\dual\ot_A\C\ot_A\C\ot_A E\rarrow E\dual\ot_A\C\ot_A E
\ot_B E\dual\ot_A\C\ot_A E$ of the map induced by the comultiplication
in $\C$ and the map induced by the map $A\rarrow E\ot_B E\dual$.
 The counit in ${}_B\C_B$ is the composition $E\dual\ot_A\C\ot_A E
\rarrow E\dual\ot_A E\rarrow B$, where the first map is induced by
the counit in~$\C$.

 All the results of
\ref{co-contra-compatible-morphisms}--\ref{co-contra-pull-push-derived}
can be generalized to the situation of a Morita morphism $E$ from
a $k$\+algebra $A$ to a $k$\+algebra $B$ and a morphism
${}_B\C_B\rarrow\D$ of corings over~$B$.
 In particular, for any left $\C$\comodule{} $\M$ there is a natural
$\D$\comodule{} structure on the $B$\module{} ${}_B\M=E\dual\ot_A\M$,
and analogously for right comodules and left contramodules.
 For any right $\C$\comodule{} $\M'$ and any left $\C$\comodule{}
$\M''$ there is a natural map $\M'\oc_\C\M''\rarrow\M'_B\oc_\D{}_B\M''$
compatible with the map $\M'\ot_A\M''\rarrow\M'_B\ot_B{}_B\M''$, etc.
 All the results of~\ref{faithful-base-ring-change} can be generalized
to the case of a Morita morphism $E$ from a $k$\+algebra $A$ to
a $k$\+algebra $B$.
 In particular, $E\dual$ is a (faithfully) flat right $A$\module{} if
and only if $E\ot_B E\dual$ is a (faithfully) flat right $A$\module,
etc.

\subsubsection{}   \label{co-contra-morita-morphisms}
 One would like to define a Morita morphism from a coring $\C$ to
a coring $\D$ as a pair consisting of a $\C$\+$\D$\bicomodule{} $\E$
and a $\D$\+$\C$\bicomodule{} $\E\dual$ endowed with maps
$\C\rarrow\E\oc_\D\E\dual$ and $\E\dual\oc_\C\E\rarrow\D$ satisfying
appropriate conditions.
 This works fine for coalgebras over fields, but in the coring
situation it is not clear how to deal with the problems of
nonassociativity of the cotensor product.
 That is why we restrict ourselves to the special case of
coflat/coprojective Morita morphisms.

 Notice that, assuming $\D$ to be a flat right $B$\module,
a $k$\+linear functor $\Lambda\:\C\comodl\rarrow\D\comodl$ is
isomorphic to a functor of the form $\M\mpsto\K\oc_\C\M$ for a certain
$\D$\+$\C$\bicomodule{} $\K$ if and only if it preserves cokernels of
the morphisms coinduced from morphisms of $A$\module s, kernels of
$A$\+split morphisms, and infinite direct sums.
 Analogously, assuming $\D$ to be a projective left $B$\module,
a $k$\+linear functor $\Lambda\:\C\contra\rarrow\D\contra$ is
isomorphic to a functor of the form $\P\mpsto\Cohom_\C(\K,\P)$ for
a certain $\C$\+$\D$\bicomodule{} $\K$ if and only if it preserves
kernels of the morphisms induced from morphisms of $A$\module s,
cokernels of $A$\+split morphisms, and infinite direct products.
 Indeed, let us compose our functor~$\Lambda$ with the induction functor
$A\modl\rarrow\C\contra$ and with the forgetful functor $\D\contra
\rarrow B\modl$; then the functor $A\modl\rarrow B\modl$ so obtained
has the form $U\mpsto\Hom_A(\K,U)$ for an $A$\+$B$\bimodule{} $\K$.
 This follows from a theorem of Watts about representability of
left exact product-preserving covariant functors on the category
of modules over a ring, which is a particular case of the abstract
adjoint functor existence theorem~\cite{McL}.
 The morphism of functors $\Hom_A(\C,\Hom_A(\C,U))\rarrow\Hom_A(\C,U)$
induces a left $\C$\+coaction in $\K$, while the functorial
$\D$\contramodule{} structures on the $B$\module s $\Hom_A(\K,U)$
induce a right $\D$\+coaction in~$\K$.
 Since the functor $\Lambda$ sends the exact sequences
$\Hom_A(\C,\Hom_A(\C,\P))\rarrow\Hom_A(\C,\P) \rarrow\P\rarrow0$
to exact sequences, it is isomorphic to the functor
$\P\mpsto\Cohom_\C(\K,\P)$.

 Let $\C$ be a coring over a $k$\+algebra $A$ and $\D$ be a coring
over a $k$\+algebra $B$.
 Assume that $\C$ is a flat right $A$\module{} and $\D$ is a flat
right $B$\module.
 A \emph{right coflat Morita morphism} from $\C$ to $\D$ is a pair
consisting of a $\D$\+coflat $\C$\+$\D$\bicomodule{} $\E$ and
a $\C$\+coflat $\D$\+$\C$\bicomodule{} $\E\dual$ endowed with
a $\C$\+$\C$\bicomodule{} morphism $\C\rarrow\E\oc_\D\E\dual$ and
a $\D$\+$\D$\bicomodule{} morphism $\E\dual\oc_\C\E\rarrow\D$
such that the two compositions $\E\rarrow\E\oc_\D\E\dual\oc_\C\E
\rarrow\E$ and $\E\dual\rarrow\E\dual\oc_\C\E\oc_\D\E\dual$ are
equal to the identity endomorphisms of $\E$ and $\E\dual$.
 A right coflat Morita morphism $(\E,\E\dual)$ from $\C$ to $\D$
induces an exact functor $\M\mpsto{}_\D\M=\E\dual\oc_\C\M$ from
the category of left $\C$\comodule s to the category of left
$\D$\comodule s and an exact functor $\N\mpsto{}_\C\N=\E\oc_\D\N$
from the category of left $\D$\comodule s to the category of left
$\C$\comodule s; the former functor is left adjoint to the latter one.
 Conversely, any pair of adjoint exact $k$\+linear functors preserving
infinite direct sums between the categories of left $\C$\comodule s
and left $\D$\comodule s is induced by a right coflat Morita morphism.

 Analogously, assume that $\C$ is a projective left $A$\module{} and
$\D$ is a projective left $B$\module.
 A \emph{left coprojective Morita morphism} from $\C$ to $\D$ is
defined as a pair consisting of a $\C$\coprojective{}
$\C$\+$\D$\bicomodule{} $\E$ and a $\D$\+coprojective{}
$\D$\+$\C$\bicomodule{} $\E\dual$ endowed with
a $\C$\+$\C$\bicomodule{} morphism $\C\rarrow\E\oc_\D\E\dual$ and
a $\D$\+$\D$\bicomodule{} morphism $\E\dual\oc_\C\E\rarrow\D$
satisfying the same conditions as above.
 A left coprojective Morita morphism $(\E,\E\dual)$ from $\C$ to $\D$
induces an exact functor $\P\mpsto{}^\D\P=\Cohom_\C(\E,\P)$ from
the category of left $\C$\contramodule s to the category of left
$\D$\contramodule s and an exact functor $\Q\mpsto{}^\C\Q=
\Cohom_\D(\E\dual,\Q)$ from the category of left $\D$\contramodule s
to the category of left $\C$\contramodule s; the former functor
is right adjoint to the latter one.
 Conversely, any pair of adjoint exact $k$\+linear functors preserving
infinite products between the categories of left $\C$\contramodule s
and left $\D$\contramodule s is induced by a left coprojective Morita
morphism.

 All the results of
\ref{co-contra-compatible-morphisms}--\ref{co-contra-pull-push-derived}
can be extended to the situation of a left coprojective and right
coflat Morita morphism from a coring $\C$ to a coring~$\D$.
 In particular, for any right $\C$\comodule{} $\M$ and any left
$\D$\contramodule{} $\Q$ the compositions $(\M\oc_\C\E)\ocn_\D\Q
\rarrow(\M\oc_\C\E)\ocn_\D\Cohom_\C(\E,\Cohom_\D(\E\dual,\Q))\rarrow
\M\ocn_\C\Cohom_\D(\E\dual,\Q)$ and $\M\ocn_\C\Cohom_\D(\E\dual,\Q)
\rarrow(\M\oc_\C\E\oc_\D\E\dual)\ocn_\C\Cohom_\D(\E\dual,\Q)\rarrow
(\M\oc_\C\E)\ocn_\D\Q$ of the maps induced by the morphisms
$\E\dual\oc_\C\E\rarrow\D$ and $\C\rarrow\E\oc_\D\E\dual$ and
the natural ``evaluation'' maps are mutually inverse isomorphisms
between the $k$\module s $\M_\D\ocn_\D\Q$ and $\M\ocn_\C{}^\C\Q$.
 For any left $\D$\contramodule{} $\Q$ there are natural isomorphisms
of left $\C$\comodule s $\Phi_\C({}^\C\Q) =\C\ocn_\C{}^\C\Q\simeq
\C_\D\ocn_\D\Q\simeq\E\ocn_\D\Q\simeq\E\oc_\D(\D\ocn_\D\Q)=
{}_\C(\Phi_\D\Q)$ by
Proposition~\ref{cotensor-contratensor-assoc}.1(e), etc.
 However, one sometimes has to impose the homological dimension
conditions on $A$ and $B$ where they were not previously needed
and strengthen the quasicoflatness (quasicoprojectivity,
quasicoinjectivity) conditions to coflatness (coprojectivity,
coinjectivity) conditions.

\subsubsection{}
 A \emph{right coflat Morita equivalence} between corings $\C$ and $\D$
is a right coflat Morita morphism $(\E,\E\dual)$ from $\C$ to $\D$
such that the bicomodule morphisms $\C\rarrow\E\oc_\D\E\dual$ and
$\E\dual\oc_\C\E\rarrow\D$ are isomorphisms; it can be also considered
as a right coflat Morita morphism $(\E\dual,\E)$ from~$\D$ to~$\C$.
 \emph{Left coflat Morita equivalences} and \emph{left coprojective
Morita equivalences} are defined in the analogous way.
 A right coflat Morita equivalence between corings $\C$ and $\D$
induces an equivalence of the categories of left $\C$\comodule s
and left $\D$\comodule s, and, assuming that $\C$ is a flat right
$A$\module{} and $\D$ is a flat right $B$\module, any equivalence
between these two $k$\+linear categories comes from a right coflat
Morita equivalence.
 Analogously, a left coprojective Morita equivalence between corings
$\C$ and $\D$ induces an equivalence of the categories of left
$\C$\contramodule s and left $\D$\contramodule s, and, assuming
that $\C$ is a projective left $A$\module{} and $\D$ is a projective
left $B$\module, any equivalence between these two $k$\+linear
categories comes from a left coprojective Morita equivalence.

 Let $\C$ be a coring over a $k$\+algebra $A$ and $(E,E\dual)$ be
a Morita morphism from $A$ to~$B$.
 If $\C$ is a flat right $A$\module{} and $E\dual$ is a faithfully flat
right $A$\module, then the pair of bicomodules $\E=\C_B=\C\ot_A E$ and
$\E\dual={}_B\C=E\dual\ot_A\C$ is a right coflat Morita equivalence
between the corings $\C$ and ${}_B\C_B$.
 Analogously, if $\C$ is a projective left $A$\module{} and $E$ is
a faithfully projective left $A$\module, then the same pair of
bicomodules $\E=\C_B$ and $\E\dual={}_B\C$ is a left coprojective
Morita equivalence between the corings $\C$ and ${}_B\C_B$.
 This is a reformulation of the results of~\ref{barr-beck-theorem}
in the case of a Morita morphism of $k$\+algebras.

 All the results of~\ref{co-contra-faithful-change-derived} can be
generalized to the situation of a Morita equivalence, satisfying
appropriate coflatness/coprojectivity conditions, between corings
$\C$ and $\D$.
 The same applies to the results of~\ref{cotensor-faithful-change},
with homological dimension conditions added when necessary and the
quasicoflatness (quasicoprojectivity, quasicoinjectivity) conditions
strengthened to coflatness (coprojectivity, coinjectivity) conditions.

\begin{rmk}
 When the rings $A$ and $B$ are semisimple, one can consider Morita
morphisms from the coring $\C$ to the coring $\D$ without any
coflatness/coprojectivity conditions imposed.
 Moreover, for any Morita morphism $(\E,\E\dual)$ from $\C$ to $\D$
the left $\C$\comodule{} $\E$ is coprojective and the right
$\C$\comodule{} $\E\dual$ is coprojective.
 In particular, any Morita equivalence between $\C$ and $\D$ is left
and right coprojective.
 On the other hand, without such conditions on the rings $A$ and $B$
not every right coflat Morita equivalence between $\C$ and $\D$ is
a left coflat Morita equivalence.
 For example, when $\C$ is a finite-dimensional coalgebra over
a field~$k$, $B$ is the algebra over~$k$ dual to~$\C$, and $\D=B$,
the right coflat Morita equivalence between $\C$ and $\D$ inducing
the equivalence of categories $\C\comodl\simeq B\modl$ is not
left coflat, since this equivalence of categories does not preserve
coflatness of comodules.
\end{rmk}

\Section{Functoriality in the Semialgebra}
\label{functoriality-II-section}

\subsection{Compatible morphisms}  \label{semi-compatible-morphisms}
 Let $\C\rarrow\D$ be a map of corings compatible with a $k$\+algebra
map $A\rarrow B$.
 Let $\S$ be a semialgebra over the coring~$\C$ and $\T$ be
a semialgebra over the coring~$\D$.

\subsubsection{}
 A map $\S\rarrow\T$ is called compatible with the maps $A\rarrow B$
and $\C\rarrow D$ if the biaction maps $A\ot_k\S\ot_k A\rarrow \S$
and $B\ot_k\T\ot_k B\rarrow \T$ form a commutative diagram with
the maps $\S\rarrow\T$ and $A\ot_k\S\ot_k A\rarrow B\ot_k\T\ot_k B$
(that is the map $\S\rarrow\T$ is an $A$\+$A$\+bimodule morphism),
the bicoaction maps $\S\rarrow\C\ot_A\S\ot_A\C$ and $\T\rarrow
\D\ot_B\T\ot_B\D$ form a commutative diagram with the maps
$\S\rarrow\T$ and $\C\ot_A\S\ot_A\C\rarrow\D\ot_B\T\ot_B\D$
(that it the induced map $B\ot_A\S\ot_A B\rarrow\T$ is
a $\D$\+$\D$\+bicomodule morphism), and furthermore,
the semimultiplication maps $\S\oc_\C\S\rarrow\S$ and
$\T\oc_\D\T\rarrow\T$ and the semiunit maps $\C\rarrow\S$ and
$\D\rarrow\T$ form commutative diagrams with the maps
$\C\rarrow\D$, \ $\S\rarrow\T$, and $\S\oc_\C\S\rarrow\T\oc_\D\T$.

 Let $\S\rarrow\T$ be a map of semialgebras compatible with
a map of corings $\C\rarrow\D$ and a $k$\+algebra map $A\rarrow B$.
 Let $\bM$ be a left $\S$\semimodule{} and $\bN$ be a left
$\T$\semimodule.
 A map $\bM\rarrow\bN$ is called compatible with the maps $A\rarrow B$,
\ $\C\rarrow\D$, and $\S\rarrow\T$ if it is compatibe with the maps
$A\rarrow B$ and $\C\rarrow\D$ as a map from a $\C$\comodule{} to
a $\D$\comodule{} and the semiaction maps $\S\oc_\C\bM\rarrow\bM$ and
$\T\oc_\D\bN\rarrow\bN$ form a commutative diagram with the maps
$\bM\rarrow\bN$ and $\S\oc_\C\M\rarrow\T\oc_\D\N$.
 Analogously, let $\bP$ be a left $\S$\semicontramodule{} and $\bQ$
be a left $\T$\semicontramodule.
 A map $\bQ\rarrow\bP$ is called compatible with the maps $A\rarrow B$,
\ $\C\rarrow\D$, and $\S\rarrow\T$ if it is compatibe with the maps
$A\rarrow B$ and $\C\rarrow\D$ as a map from a $\D$\contramodule{} to
a $\C$\contramodule{} and the semicontraaction maps $\bP\rarrow
\Cohom_\C(\S,\bP)$ and $\bQ\rarrow\Cohom_\D(\T,\bQ)$ form a commutative
diagram with the maps $\bQ\rarrow\bP$ and $\Cohom_\D(\T,\bQ)\rarrow
\Cohom_\C(\S,\bP)$.

 Let $\bM'\rarrow\bN'$ be a map from a right $\S$\semimodule{} to
a right $\T$\semimodule{} compatible with the maps $A\rarrow B$, \ 
$\C\rarrow\D$, and $\S\rarrow\T$, and let $\bM''\rarrow\bN''$ be a map
from a left $\S$\semimodule{} to a left $\T$\semimodule{} compatible
with the maps $A\rarrow B$, \ $\C\rarrow\D$, and $\S\rarrow\T$.
 Assume that the triple cotensor products $\bM'\oc_\C\S\oc_\C\bM''$
and $\bN'\oc_\D\T\oc_\D\bN''$ are associative.
 Then there is a natural map of $k$\module s $\bM'\os_\S\bM''\rarrow
\bN'\os_\S\bN''$.
 Analogously, let $\bM\rarrow\bN$ be a map from a left
$\S$\semimodule{} to a left $\T$\semimodule{} compatible with
the maps $A\rarrow B$, \ $\C\rarrow\D$, and $\S\rarrow\T$, and let
$\bQ\rarrow\bP$ be a map from a left $\T$\semicontramodule{} to
a left $\S$\semicontramodule{} compatible with the maps $A\rarrow B$,
\ $\C\rarrow\D$, and $\S\rarrow\T$.
 Assume that the triple cohomomorphisms $\Cohom_\C(\S\oc_\C\bM\;\bP)$
and $\Cohom_\D(\T\oc_\D\bN\;\bQ)$ are associative.
 Then there is a natural map of $k$\module s $\SemiHom_\T(\bN,\bQ)
\rarrow\SemiHom_\S(\bM,\bP)$.

\subsubsection{}
 Let $\S\rarrow\T$ be a map of semialgebras compatible with a map
of corings $\C\rarrow\D$ and a $k$\+algebra map $A\rarrow B$.

 Assume that $\C$ is a flat right $A$\module{} and either $\S$
is a coflat right $\C$\comodule, or $\S$ is a flat right $A$\module{}
and a $\C/A$\+coflat left $\C$\comodule{} and $A$ has a finite weak
homological dimension, or $A$ is absolutely flat.
 Then for any left $\T$\semimodule{} $\bN$ there is a natural
$\S$\semimodule{} structure on the left $\C$\comodule{} ${}_\C\bN$.
 It is constructed as follows: the composition $\S\oc_\C{}_\C\bN
\rarrow\T\oc_\D\bN\rarrow\bN$ of the map induced by the maps
$\S\rarrow\T$ and ${}_\C\bN\rarrow\bN$ with the $\T$\+semiaction
in~$\bN$ is a map from a $\C$\comodule{} to a $\D$\comodule{}
compatible with the maps $A\rarrow B$ and $\C\rarrow\D$, hence
there is a $\C$\comodule{} map $\S\oc_\C{}_\C\bN\rarrow{}_\C\bN$.
 Analogously, assume that $\C$ is a projective left $A$\module{} and
either $\S$ is a coprojective left $\C$\comodule, or $\S$ is
a projective left $A$\module{} and a $\C/A$\+coflat right
$\C$\comodule{} and $A$ has a finite left homological dimension,
or $A$ is semisimple.
 Then for any left $\T$\semicontramodule{} $\bQ$ there is a natural
$\S$\semicontramodule{} structure on the left $\C$\contramodule{}
${}^\C\bQ$.
 Indeed, the composition $\bQ\rarrow\Cohom_\D(\T,\bQ)\rarrow
\Cohom_\C(\S,{}^\C\bQ)$ is a map from a $\D$\contramodule{} to
a $\C$\contramodule{} compatible with the maps $A\rarrow B$ and
$\C\rarrow\D$, hence a $\C$\contramodule{} map ${}^\C\bQ\rarrow
\Cohom_\C(\S,{}^\C\bQ)$.
 Assuming that $\D$ is a flat right $B$\module, $\C$ is a flat right
$A$\module, and $\S$ is a coflat right $\C$\comodule, for any
$\D$\+coflat right $\T$\semimodule{} $\bN$ there is a natural
$\S$\semimodule{} structure on the coflat right $\C$\comodule{}
$\bN_\C$ and for any $\D$\coinjective{} left $\T$\semicontramodule{}
$\bQ$ there is a natural $\S$\semicontramodule{} structure on
the coinjective left $\C$\contramodule{} ${}^\C\bQ$ provided that
$B$ is a flat right $A$\module.

 Assume that $\C$ is a flat right $A$\module, $\S$ is a coflat right
$\C$\comodule, $\D$ is a flat right $B$\module, and $\T$ is a coflat
right $\D$\comodule.
 Then the functor $\bN\mpsto{}_\C\bN$ from the category of left
$\T$\semimodule s to the category of left $\S$\semimodule s has a left
adjoint functor $\bM\mpsto{}_\T\bM$, which is constructed as follows.
 For induced left $\S$\semimodule s, one has ${}_\T(\S\oc_\C\L)=
\T\oc_\D{}_B\L$; to compute the $\T$\semimodule{} ${}_\T\bM$ for
an arbitrary left $\S$\semimodule{} $\bM$, one can represent $\bM$ as
the cokernel of a morphism of induced $\S$\semimodule s.
 Both $k$\module s $\Hom_\S(\bM,{}_\C\bN)$ and $\Hom_\T({}_\T\bM,\bN)$
are isomorphic to the $k$\module{} of all maps of semimodules
$\bM\rarrow\bN$ compatible with the maps $A\rarrow B$, \ $\C\rarrow\D$,
and $\S\rarrow\T$.
 There are also a few situations when the functor $\bM\mpsto{}_\T\bM$
is defined on the full subcategory of induced $\S$\semimodule s.
 Under analogous assumptions, the functor $\bM\mpsto\bM_\T$ left
adjoint to the functor $\bN\mpsto\bN_\C$ acts from the category
of right $\S$\semimodule s to the category of right $\T$\semimodule s.

 Now assume that $\C$ is a flat left and right $A$\module, $\S$ is
a flat left $A$\module{} and a coflat right $\C$\comodule, $A$ has
a finite weak homological dimension, $\D$ is a flat right $B$\module,
and $\T$ is a coflat right $\D$\comodule.
 Then the functor $\bN\mpsto{}_\C\bN$ can be constructed in a different
way: when $\bM$ is a flat left $A$\module, one has ${}_\T\bM=
\T_\C\os_\S\bM$, where $\T_\C=\T\oc_\D{}_B\C$ is
a $\T$\+$\S$\bisemimodule{} with the right $\S$\semimodule{}
structure provided by the above construction.
 To compute the $\T$\semimodule{} ${}_\T\bM$ for an arbitrary left
$\S$\semimodule{} $\bM$, one can represent $\bM$ as the cokernel of
a morphism of $A$\+flat $\S$\semimodule s.
 Assuming only that $\C$ is a flat right $A$\module, $\S$ is a coflat
right $\C$\comodule, $\D$ is a flat right $B$\module, and $\T$ is
a coflat right $\D$\comodule, the functor $\bM\mpsto{}_\T\bM$ can be
defined by the formula ${}_\T\bM=\T_\C\os_\S\bM$ for any $\bM$
whenever $B$ is a flat right $A$\module.
 If $\C$ is a flat left and right $A$\module, $\S$ is a coflat left
and right $\C$\comodule, $\D$ is a flat right $B$\module, and $\T$
is a coflat right $\D$\comodule, the functor $\bM\mpsto{}_\T\bM$ is
given by the formula ${}_\T\bM=\T_\C\os_\S\bM$ on the full subcategory
of $\C$\+coflat $\S$\semimodule s~$\bM$.

 Furthermore, assume that $\C$ is a projective left $A$\module,
$\S$ is a coprojective left $\C$\comodule, $\D$ is a projective left
$B$\module, and $\T$ is a coprojective left $\D$\comodule.
 Then the functor $\bQ\mpsto{}^\C\bQ$ from the category of left
$\T$\semicontramodule s to the category of left $\S$\semicontramodule s
has a right adjoint functor $\bP\mpsto{}^\T\bP$, which is constructed
as follows.
 For coinduced left $\S$\semicontramodule s, one has
${}^\T\!\Cohom_\C(\S,\gR)=\Cohom_\D(\T,{}^B\gR)$; to compute
the $\T$\semicontramodule{} ${}^\T\bP$ for an arbitrary left
$\S$\semicontramodule{} $\bP$, one can represent $\bP$ as the kernel
of a morphism of coinduced $\S$\semicontramodule s.
 Both $k$\module s $\Hom^\S({}^\C\bQ,\bP)$ and $\Hom^\T(\bQ,{}^\T\bP)$
are isomorphic to the $k$\module{} of all maps of semicontramodules
$\bQ\rarrow\bP$ compatible with the maps $A\rarrow B$, \ $\C\rarrow\D$,
and $\S\rarrow\T$.
 There are also a few situations when the functor $\bP\mpsto{}^\T\bP$
is defined on the full subcategory of coinduced
$\S$\semicontramodule s.

 Now assume that $\C$  is a projective left and a flat right
$A$\module, $\S$ is a coprojective left $\C$\comodule{} and a flat
right $A$\module, $A$ has a finite left homological dimension, $\D$
is a projective left $B$\module, and $\T$ is a coprojective left
$\D$\comodule.
 Then the functor $\bP\mpsto{}^\T\bP$ can be constructed in a different
way: when $\bP$ is an injective left $A$\module, $^\T\bP=
\SemiHom_\S({}_\C\T,\bP)$; to compute the $\T$\semicontramodule{}
${}^\T\bP$ for an arbitrary left $\S$\semicontramodule{} $\bP$, one can
represent $\bP$ as the kernel of a morphism of $A$\injective{}
$\S$\semicontramodule s.
 Assuming only that $\C$ is a projective left $A$\module, $\S$ is
a coprojective left $\C$\comodule, $\D$ is a projective left
$B$\module, and $\T$ is a coprojective left $\D$\comodule, the functor
$\bP\mpsto{}^\T\bP$ can be defined by the formula $^\T\bP=
\SemiHom_\S({}_\C\T,\bP)$ for any $\bP$ whenever $B$ is a projective
left $A$\module.
 If $\C$ is a projective left and a flat right $A$\module, $\S$ is
a coprojective left and a coflat right $\C$\comodule, $\D$ is
a projective left $B$\module, and $\T$ is a coprojective left
$\D$\comodule, the functor $\bP\mpsto{}^\T\bP$ is given by the formula
$^\T\bP=\SemiHom_\S({}_\C\T,\bP)$ on the full subcategory of
$\C$\coinjective{} $\S$\semicontramodule s~$\bP$.

 Assume that $\C$ is a projective left $A$\module, $\S$ is
a coprojective left $\C$\comodule, $\D$ is a projective left
$B$\module, and $\T$ is a coprojective left $\D$\comodule.
 Then for any right $\S$\semimodule{} $\bM$ and any left
$\T$\semicontramodule{} $\bQ$ there is a natural isomorphism
$\bM_\T\Ocn_\T\bQ\simeq\bM\Ocn_\S{}^\C\bQ$.
 Moreover, both $k$\module s are isomorphic to the cokernel of the pair
of maps $(\bM\oc_\C\S)_B\ocn_\D\bQ\birarrow\bM_B\ocn_\D\bQ$ one of
which is induced by the $\S$\+semiaction in $\bM$ and the other is
defined in terms of the morphism $(\bM\oc_\C\S)_B\rarrow\bM_B\oc_\D\T$,
the $\T$\+semicontraaction in $\bQ$, and the natural ``evaluation''
map $(\bM_B\oc_\D\T)\ocn_\D\Cohom_\D(\T,\bQ)\rarrow\bM_B\ocn_\D\bQ$.
 This is clear for $\bM\Ocn_\S{}^\C\bQ$, and to construct this
isomorphism for $\bM_\T\Ocn_\T\bQ$ it suffices to represent $\bM$ as
the cokernel of the pair of morphisms of induced $\S$\semimodule s
$\bM\oc_\C\S\oc_\C\S\birarrow\bM\oc_\C\S$.
 In the above situations when $\bM_\T=\bM\os_\S{}_\C\T$, this
isomorphism can be also constructed by representing $\bM_\T$ as
the cokernel of the pair of $\T$\semimodule{} morphisms
$\bM\oc_\C\S\oc_\C{}_\C\T\birarrow\bM\oc_\C{}_\C\T$ and using
the isomorphisms $\bM\oc_\C{}_\C\T\simeq\bM_B\oc_\D\T$.

\subsubsection{}   \label{pull-push-semitensor}
 Let $\S\rarrow\T$ be a map of semialgebras compatible with a map
of corings $\C\rarrow\D$ and a $k$\+algebra map $A\rarrow B$.

\begin{prop}
 \textup{(a)} Let\/ $\bM$ be a left\/ $\S$\semimodule{} and\/ $\bN$ be
a right\/ $\T$\semimodule.
 Then the semitensor product  ${}_\T\bM=\T_\C\os_\S\bM$ can be
endowed with a left $\T$\semimodule{} structure via the construction
of~\ref{semitensor-associative} and the map of semitensor products\/
$\bN_\C\os_\S\bM\rarrow\bN\os_\T{}_\T\bM$ induced by the maps of
semimodules\/ $\bN_\C\rarrow\bN$ and\/ $\bM\rarrow{}_\T \bM$ is
an isomorphism, at least, in the following cases:
\begin{itemize}
 \item $\D$ is a flat right $B$\module, $\T$ is a coflat right\/
       $\D$\comodule, $\bN$ is a coflat right\/ $\D$\comodule,
       $\C$ is a flat left $A$\module, $\S$ is a flat left $A$\module{}
       and a\/ $\C/A$\+coflat right\/ $\C$\comodule, the ring $A$ has
       a finite weak homological dimension, and\/ $\bM$ is a flat left
       $A$\module, or
 \item $\D$ is a flat right $B$\module, $\T$ is a coflat right\/
       $\D$\comodule, $\bN$ is a coflat right\/ $\D$\comodule,
       $\C$ is a flat left $A$\module, $\S$ is a coflat left\/
       $\C$\comodule, and\/ $\bM$ is a coflat left\/ $\C$\comodule, or
 \item $\D$ is a flat right $B$\module, $\T$ is a coflat right\/
       $\D$\comodule, $\bN$ is a coflat right\/ $\D$\comodule,
       $\C$ is a flat right $A$\module, $\S$ is a coflat right\/
       $\C$\comodule, and $B$ is a flat right $A$\module, or
 \item $\D$ is a flat left $B$\module, $\T$ is a flat left $B$\module{}
       and a\/ $\D/B$\+coflat right\/ $\D$\comodule, the ring $B$ has
       a finite weak homological dimension, $\C$ is a flat left
       $A$\module, $\S$ is a coflat left\/ $\C$\comodule, and\/
       $\bM$ is a semiflat left\/ $\S$\semimodule, or
 \item $\D$ is a flat left $B$\module, $\T$ is a coflat left\/
       $\D$\comodule, ${}_B\C$ is a coflat left\/ $\D$\comodule,
       $\C$ is a flat left $A$\module, $\S$ is a coflat left\/
       $\C$\comodule, and\/ $\bM$ is a semiflat left\/ $\S$\semimodule.
\end{itemize}
 When the ring $A$ (resp., $B$) is absolutely flat,
the\/ $\C/A$\+coflatness (resp., $\D/B$\+co\-flat\-ness) assumption
can be dropped. \par
 \textup{(b)} Let\/ $\bP$ be a left\/ $\S$\semicontramodule{} and\/
$\bN$ be a left\/ $\T$\semimodule.
 Then the module of semihomomorphisms\/ ${}^\T\bP=\SemiHom_\S
({}_\C\T,\bP)$ can be endowed with a left\/ $\T$\semicontramodule{}
structure via the construction of~\ref{semihom-associative} and
the map of the semihomomorphism modules\/ $\SemiHom_\T(\bN,{}^\T\bP)
\rarrow\SemiHom_\S({}_\C\bN,\bP)$ induced by the maps of semimodules
and semicontramodules\/ ${}_\C\bN\rarrow\bN$ and\/ ${}^\T\bP\rarrow\bP$
is an isomorphism, at least, in the following cases:
\begin{itemize}
 \item $\D$ is a projective left $B$\module, $\T$ is a coprojective
       left\/ $\D$\comodule, $\bN$ is a coprojective left\/
       $\D$\comodule, $\C$ is a flat right $A$\module, $\S$ is a flat
       right $A$\module{} and a\/ $\C/A$\coprojective{} left\/
       $\C$\comodule, the ring $A$ has a finite left homological
       dimension, and\/ $\bP$ is an injective left $A$\module, or
 \item $\D$ is a projective left $B$\module, $\T$ is a coprojective
       left\/ $\D$\comodule, $\bN$ is a coprojective left\/
       $\D$\comodule, $\C$ is a flat right $A$\module, $\S$ is a coflat
       right $\C$\comodule, and\/ $\bP$ is a coinjective left\/
       $\C$\comodule, or
 \item $\D$ is a projective left $B$\module, $\T$ is a coprojective
       left\/ $\D$\comodule, $\bN$ is a coprojective left\/
       $\D$\comodule, $\C$ is a projective left $A$\module, $\S$ is
       a coprojective left\/ $\C$\comodule, and $B$ is a projective
       left $A$\module, or
 \item $\D$ is a flat right $B$\module, $\T$ is a flat right
       $B$\module{} and a\/ $\D/B$\coprojective{} left\/ $\D$\comodule,
       the ring $B$ has a finite left homological dimension, $\C$ is
       a flat right $A$\module, $\S$ is a coflat right\/ $\C$\comodule,
       and\/ $\bP$ is a semiinjective left\/ $\S$\semicontramodule, or
 \item $\D$ is a flat right $B$\module, $\T$ is a coflat right\/
       $\D$\comodule, $\C_B$ is a coflat right\/ $\D$\comodule,
       $\C$ is a flat right $A$\module, $\S$ is a coflat right\/
       $\C$\comodule, and\/ $\bP$ is a semiinjective left\/
       $\S$\semicontramodule.
\end{itemize}
 When the ring $A$ (resp., $B$) is semisimple,
the\/ $\C/A$\+coprojectivity (resp., $\D/B$\+co\-pro\-jec\-tiv\-ity)
assumption can be dropped. \par
 \textup{(c)} Let\/ $\bM$ be a left $\S$\semimodule{} and\/ $\bQ$ be
a left\/ $\T$\semicontramodule.
 Then the map of semihomomorphism modules\/ $\SemiHom_\T({}_\T\bM,\bQ)
\rarrow\SemiHom_\S(\bM,{}^\C\bQ)$ induced by the map of semimodules\/
$\bM\rarrow{}_\T\bM$ and the map of semicontramodules\/ $\bQ\rarrow
{}^\C\bQ$ is an isomorphism, at least, in the following cases:
\begin{itemize}
 \item $\D$ is a flat right $B$\module, $\T$ is a coflat right\/
       $\D$\comodule, $\bQ$ is a coinjective left\/ $\D$\contramodule,
       $\C$ is a projective left $A$\module, $\S$ is a projective
       left $A$\module{} and a\/ $\C/A$\+coflat right\/ $\C$\comodule,
       the ring $A$ has a finite left homological dimension, and\/
       $\bM$ is a projective left $A$\module, or
 \item $\D$ is a flat right $B$\module, $\T$ is a coflat right\/
       $\D$\comodule, $\bQ$ is a coinjective left\/ $\D$\contramodule,
       $\C$ is a projective left $A$\module, $\S$ is a coprojective
       left\/ $\C$\comodule, and\/ $\bM$ is a coprojective left\/
       $\C$\comodule, or
 \item $\D$ is a flat right $B$\module, $\T$ is a coflat right\/
       $\D$\comodule, $\bQ$ is a coinjective left\/ $\D$\contramodule,
       $\C$ is a flat right $A$\module, $\S$ is a coflat right\/
       $\C$\comodule, and $B$ is a flat right $A$\module, or
 \item $\D$ is a projective left $B$\module, $\T$ is a projective left
       $B$\module{} and a\/ $\D/B$\+coflat right\/ $\D$\comodule,
       the ring $B$ has a finite left homological dimension, $\C$ is
       a projective left $A$\module, $\S$ is a coprojective left\/
       $\C$\comodule, and\/ $\bM$ is a semiprojective left\/
       $\S$\semimodule, or
 \item $\D$ is a projective left $B$\module, $\T$ is a coprojective
       left\/ $\D$\comodule, ${}_B\C$ is a coprojective left\/
       $\D$\comodule, $\C$ is a projective left $A$\module, $\S$ is
       a coprojective left\/ $\C$\comodule, and\/ $\bM$ is
       a semiprojective left\/ $\S$\semimodule.
\end{itemize}
 When the ring $A$ (resp., $B$) is semisimple, the\/ $\C/A$\+coflatness
(resp., $\D/B$\+co\-flat\-ness) assumption can be dropped.
\end{prop}

\begin{proof}
 Part~(a): under our assumptions, there is a natural isomorphism of
right $\S$\semimodule s $\bN_\C\simeq\bN\os_\T\T_\C$.
 For any left $\S$\semimodule{} $\bM$ and right $\T$\semimodule{} $\bN$
for which the iterated semitensor products $(\bN\os_\T\T_\C)\os_\S\bM$
and $\bN\os_\T (\T_\C\os_\S\bM)$ are defined and the triple cotensor
product $\bN\oc_\D\T_\C\oc_\C\bM$ is associative, the map
$(\bN\os_\T\T_\C)\os_\S\bM\rarrow\bN\os_\T(\T_\C\os_\S\bM)$ induced by
the bisemimodule maps $\S\rarrow\T_\C\rarrow\T$ compatible with
the maps $A\rarrow B$, \ $\C\rarrow\D$, and $\S\rarrow\T$ forms
a commutative diagram with the maps $\bN\oc_\D\T_\C\oc_\C\bM\rarrow
(\bN\os_\T\T_\C)\os_\S\bM$ and $\bN\oc_\D\T_\C\oc_\C\bM\rarrow
\bN\os_\T(\T_\C\os_\S\bM)$.
 Indeed, the map $\bN\oc_\D\T_\C\oc_\C\bM\rarrow\bN\oc_\D
(\T_\C\os_\S\bM)$ is equal to the composition of the map
$\bN\oc_\D\T_\C\oc_\C\bM\rarrow\bN\oc_\D\T\oc_\D(\T_\C\os_\S\bM)$
induced by the maps $\T_\C\rarrow\T$ and $\bM\rarrow\T_\C\os_\S\bM$
and the map $\bN\oc_\D\T\oc_\D(\T_\C\os_\S\bM)\rarrow\bN\oc_\D
(\T_\C\os_\S\bM)$ induced by the left $\T$\+semiaction in
$\T_\C\os_\S\bM$.
 To check this, one can notice that the diagram in question is obtained
by taking the cotensor product with $\bN$ of the diagram of maps
$\T_\C\oc_\C\bM\rarrow\T\oc_\D(\T_\C\os_\S\bM)\rarrow\T_\C\os_\S\bM$
and compose the latter diagram with the surjective map
$\T_\C\oc_\C\S\oc_\C\bM\rarrow\T_\C\oc_\C\bM$ induced by
the left $\S$\+semiaction in $\bM$.
 On the other hand, the composition of maps $\bN\oc_\D\T_\C\oc_\C\bM
\rarrow(\bN\os_\T\T_\C)\os_\S\bM\rarrow\bN\oc_\D(\T_\C\os_\S\bM)$ is
equal to the composition of the same map $\bN\oc_\D\T_\C\oc_\C\bM
\rarrow\bN\oc_\D\T\oc_\D(\T_\C\os_\S\bM)$ and the map
$\bN\oc_\D\T\oc_\D(\T_\C\os_\S\bM)\rarrow\bN\oc_\D(\T_\C\os_\S\bM)$
induced by the right $\T$\+semiaction in $\bN$, since both compositions
are equal to the composition of the map $\bN\oc_\D\T_\C\oc_\C\bM\rarrow
\bN\oc_\C\bM$ induced by the composition $\bN\oc_\D\T_\C\rarrow
\bN\os_\T\T_\C\rarrow\bN$ with the map $\bN\oc_\C\bM\rarrow\bN\oc_\D
(\T_\C\os_\S\bM)$ induced by the map $\bM\rarrow\T_\C\os_\S\bM$.
 It remains to apply Proposition~\ref{semitensor-associative}.
 The proofs of parts~(b) and~(c) are completely analogous.
\end{proof}

\subsubsection{}   \label{semi-pull-psi-phi}
 Let $\S\rarrow\T$ be a map of semialgebras compatible with a map
of corings $\C\rarrow\D$ and a $k$\+algebra map $A\rarrow B$.

 Assume that $\C$ is a projective left and a flat right $A$\module,
$\S$ is a coprojective left and a coflat right $\C$\comodule,
$\D$ is a projective left and a flat right $B$\module, and $\T$
is a coprojective left and a coflat right $\C$\comodule.
 Then for any left $\T$\semicontramodule{} $\bQ$ the natural map
of $\C$\comodule s $\Phi_\C({}^\C\bQ)\rarrow{}_\C(\Phi_\D\bQ)$
is an $\S$\semimodule{} morphism $\Phi_\S({}^\C\Q)\rarrow
{}_\C(\Phi_\T\bQ)$.
 Indeed, $\Phi_\S({}^\C\bQ)=\S\Ocn_\S{}^\C\bQ\simeq{}\S_\T\Ocn_\T\bQ
\simeq{}_\C\T\Ocn_\T\bQ$ as a left $\S$\semimodule{} and
${}_\C(\Phi_\T\bQ)={}_\C(\T\Ocn_\T\bQ)$, so there is an
$\S$\semimodule{} morphism $\Phi_\S({}^\C\Q)\rarrow{}_\C(\Phi_\T\bQ)$;
it coincides with the $\C$\comodule{} morphism $\Phi_\C({}^\C\bQ)
\rarrow{}_\C(\Phi_\D\bQ)$ defined in~\ref{pull-psi-phi}.
 Analogously, for any left $\T$\semimodule{} $\bN$ the natural map
of $\C$\contramodule s ${}^\C(\Psi_\D\bN)\rarrow\Psi_\C({}_\C\bN)$
is an $\S$\semicontramodule{} morphism ${}^\C(\Psi_\T\bN)\rarrow
\Psi_\S({}_\C\bN)$.
 Indeed, $\Psi_\S({}_\C\bN)=\Hom_\S(\S,{}_\C\bN)\simeq
\Hom_\T(\S_\T,\bN)\simeq\Hom_\T(\T_\C,\bN)$ as a left
$\S$\semicontramodule{} and
${}^\C(\Psi_\T\bN)={}^\C\!\.\Hom_\T(\T,\bN)$.

 Assume that $\C$ is a projective left $A$\module, $\S$ is
a coprojective left $\C$\comodule, $\D$ is a projective left
$B$\module, $\T$ is a coprojective left $\D$\comodule, and
$B$ is a projective left $A$\module.
 Then the equivalence of categories of $\C$\coprojective{} left
$\S$\semimodule s and $\C$\projective{} left $\S$\semicontramodule s
and the equivalence of categories of $\D$\coprojective{} left
$\T$\semimodule s and $\D$\projective{} left $\S$\semicontramodule s
transform the functor $\bN\mpsto{}_\C\bN$ into the functor
$\bQ\mpsto{}^\C\bQ$.
 Indeed, the above argument shows that for any $\D$\projective{}
left $\T$\semicontramodule{} $\bQ$ the isomorphism $\Phi_\C({}^\C\bQ)
\simeq{}_\C(\Phi_\D\bQ)$ preserves the $\S$\semimodule{} structures.

 Assume that $\C$ is a flat right $A$\module, $\S$ is a coflat right
$\C$\comodule, $\D$ is a flat right $B$\module, $\T$ is a coflat right
$\D$\comodule, and $B$ is a flat right $A$\module.
  Then the equivalence of categories of $\C$\injective{} left
$\S$\semimodule s and $\C$\coinjective{} left $\S$\semicontramodule s
and the equivalence of categories of $\D$\injective{} left
$\T$\semimodule s and $\D$\coinjective{} left $\S$\semicontramodule s
transform the functor $\bN\mpsto{}_\C\bN$ into the functor
$\bQ\mpsto{}^\C\bQ$.
 Indeed, the above argument shows that for any $\D$\injective{}
left $\T$\semimodule{} $\bN$ the isomorphism ${}^\C(\Psi_\D\bN)
\simeq\Psi_\C({}_\C\bN)$ preserves the $\S$\semicontramodule{}
structures.

 Assume that $\C$ is a projective left and a flat right $A$\module,
$\S$ is a coprojective left $\C$\comodule{} and a flat right
$A$\module, $\D$ is a projective left and a flat right $B$\module,
$\T$ is a coprojective left $\D$\comodule{} and a flat right
$B$\module, and the rings $A$ and $B$ have finite left homological
dimensions.
 Then the equivalence of categories of $\C/A$\injective{} left
$\S$\semimodule s and $\C/A$\projective{} left $\S$\semicontramodule s
and the equivalence of categories of $\D/B$\injective{} left
$\T$\semimodule s and $\D/B$\projective{} left $\T$\semicontramodule s
transform the functor $\bN\mpsto{}_\C\bN$ into the functor
$\bQ\mpsto{}^\C\bQ$.
  Indeed, the above argument shows that for any $\D/B$\projective{}
left $\T$\semicontramodule{} $\bQ$ the isomorphism $\Phi_\C({}^\C\bQ)
\simeq{}_\C(\Phi_\D\bQ)$ preserves the $\S$\semimodule{} structures.
 The analogous result holds when $\S$ is a projective left $A$\module{}
and a coflat right $\C$\comodule{} and $\T$ is a projective
left $B$\module{} and a coflat right $\D$\comodule; it can be proven
by applying the above argument to the isomorphism ${}^\C(\Psi_\D\bN)
\simeq\Psi_\C({}_\C\bN)$ for a $\D/B$\injective{} left
$\T$\semimodule{} $\bN$.

 Finally, assume that the rings $A$ and $B$ are semisimple.
 Then the equivalence of categories of $\C$\injective{} left
$\S$\semimodule s and $\C$\projective{} left $\S$\semicontramodule s
and the equivalence of categories of $\D$\injective{} left
$\T$\semimodule s and $\D$\projective{} left $\T$\semicontramodule s
transform the functor $\bN\mpsto{}_\C\bN$ into the functor
$\bQ\mpsto{}^\C\bQ$.
 One can show this using the semialgebra analogues of the assertions
of~\ref{co-contra-pull-push} related to quasicoflat comodules and
quasicoinjective contramodules.

\subsection{Complexes, adjusted to pull-backs and push-forwards}
\label{semi-pull-push-adjusted}
 Let $\S\rarrow\T$ be a map of semialgebras compatible with a map
of corings $\C\rarrow\D$ and a $k$\+algebra map $A\rarrow B$.
 The following result generalizes
Theorem~\ref{semimodule-semicontramodule-subsect}.

\begin{thm1}
 \textup{(a)} Assume that\/ $\D$ is a flat right $B$\module,
$\T$ is a coflat right $\D$\comodule{} and a\/ $\D/B$\+coflat
(\.$\D/B$\coprojective) left\/ $\D$\comodule, and the ring $B$
has a finite weak (left) homological dimension.
 Then the functor mapping the quotient category of the homotopy
category of complexes of\/ $\D/B$\+coflat (\.$\D/B$\coprojective)
left\/ $\T$\semimodule s by its intersection with the thick subcategory
of\/ $\D$\coacyclic{} complexes into the semiderived category of left\/
$\T$\semimodule s is an equivalence of triangulated categories. \par
 \textup{(b)} Assume that\/ $\D$ is a projective left $B$\module,
$\T$ is a coprojective left and a $\D/B$\+coflat right\/ $\D$\comodule,
and the ring $B$ has a finite left homological dimension.
 Then the functor mapping the quotient category of the homotopy
category of complexes of\/ $\D/B$\coinjective{} left\/
$\T$\semicontramodule s by its intersection with the thick subcategory
of\/ $\D$\contraacyclic{} complexes into the semiderived category of
left\/ $\T$\semicontramodule s is an equivalence of triangulated
categories.
\end{thm1}

\begin{proof}
 To prove part~(a) for $\D/B$\+coflat $\T$\semimodule s, use
Lemma~\ref{coflat-semimodule-injection}, the construction of
the morphism of complexes $\bL^\bu\rarrow\boR_2(\bL^\bu)$ from
the proof of Theorem~\ref{semitor-main-theorem}, and
Lemma~\ref{semitor-main-theorem}.
 To prove part~(a) for $\D/B$\coprojective{} $\T$\semimodule s,
use Lemma~\ref{coproj-coinj-semi-mod-contra}(b).
 To prove part~(b), use Lemma~\ref{coproj-coinj-semi-mod-contra}(a)
and the construction of the morphism of complexes $\boL_2(\bgR^\bu)
\rarrow\bgR^\bu$ from the proof of Theorem~\ref{semiext-main-theorem}.
\end{proof}

 A complex of $\S$\semimodule s is called \emph{quite
$\S/\C/A$\semiflat} (\emph{quite $\S/\C/A$\semiprojective}) if it
belongs to the minimal triangulated subcategory of the homotopy
category of complexes of $\S$\semimodule s containing the complexes
induced from complexes of $A$\+flat ($A$\projective) $\C$\comodule s
and closed under infinite direct sums.
 Analogously, a complex of $\S$\semicontramodule s is called
\emph{quite $\S/\C/A$\semiinjective} if it belongs to the minimal
triangulated subcategory of the homotopy category of complexes of
$\S$\semicontramodule s containing the complexes coinduced from
complexes of $A$\injective{} $\C$\contramodule s and closed under
infinite products.
 Under appropriate assumptions on $\S$, $\C$, and $A$, any quite
$\S/\C/A$\semiflat{} complex of $A$\+flat $\S$\semimodule s is
$\S/\C/A$\semiflat{} in the sense of~\ref{relatively-semiflat},
and analogously for birelative semiprojectivity and semiinjectivity
in the sense of~\ref{relatively-semiproj-semiinj}.
 Any quite $\S/\C/A$\semiflat{} complex of right $\S$\semimodule s is
$\S/\C/A$\contraflat, any quite $\S/\C/A$\semiprojective{} complex of
left $\S$\semimodule s is $\S/\C/A$\projective, and any quite
$\S/\C/A$\semiinjective{} complex of left $\S$\semicontramodule s is
$\S/\C/A$\injective{} in the sense of~\ref{birelatively-adjusted}.

\begin{thm2}
 \textup{(a)} Assume that\/ $\C$ is a flat (projective) left and 
a flat right $A$\module, $\C$ is a flat (projective) left $A$\module{}
and a coflat right\/ $\C$\comodule, and the ring $A$ has a finite weak
(left) homological dimension.
 Then the functor mapping the quotient category of the homotopy
category of quite\/ $\S/\C/A$\semiflat{} (quite\/
$\S/\C/A$\semiprojective) complexes of left\/ $\S$\semimodule s by
its minimal triangulated subcategory containing complexes induced from
coacyclic complexes of $A$\+flat ($A$\projective)\/ $\S$\semimodule s
and closed under infinite direct sums into the semiderived category of
left\/ $\S$\semimodule s is an equivalence of triangulated categories.
\par
 \textup{(b)} Assume that\/ $\C$ is a projective left and a flat right
$A$\module, $\C$ is a coprojective left\/ $\C$\comodule{} and a flat
right $A$\module, and the ring $A$ has a finite left homological
dimension.
 Then the functor mapping the quotient category of the homotopy
category of quite\/ $\S/\C/A$\semiinjective{} complexes of left\/
$\S$\semicontramodule s by its minimal triangulated subcategory
containing complexes of coinduced from contraacyclic complexes of
$A$\injective\/ $\C$\contramodule s and closed under infinite products
into the semiderived category of left\/ $\S$\semicontramodule s is
an equivalence of triangulated categories.
\end{thm2}

\begin{proof}
 Proof of part~(a): for any complex of $\S$\semimodule s $\bK^\bu$
there is a natural morphism into $\bK^\bu$ from a quite
$\S/\C/A$\semiflat{} complex of $\S$\semimodule s
$\boL_3\boL_1(\bK^\bu)$ with a $\C$\coacyclic{} cone.
 Hence it follows from Lemma~\ref{semitor-main-theorem} that
the semiderived category of $\S$\semimodule s is equivalent to
the quotient category of the homotopy category of quite
$\S/\C/A$\semiflat{} complexes of $\S$\semimodule s by its
intersection with the thick subcategory of $\C$\coacyclic{} complexes.
 It remains to show that any $\C$\coacyclic{} quite
$\S/\C/A$\semiflat{} complex of $\S$\semimodule s belongs to
the minimal triangulated subcategory containing the complexes induced
from coacyclic complexes of $A$\+flat $\S$\semimodule s and closed
under infinite direct sums.
 Indeed, if a complex of $A$\+flat left $\S$\semimodule s $\bM^\bu$
is $\C$\coacyclic, then the total complex $\boL_3(\bM^\bu)$ of the bar
bicomplex $\dsb\rarrow\S\oc_\C\S\oc_\C\bM^\bu\rarrow\S\oc_\C\bM^\bu$
up to the homotopy equivalence can be obtained from complexes of
$\S$\semimodule s induced from coacyclic complexes of $A$\+flat
$\C$\comodule s using the operations of cone and infinite direct sum.
 So the same applies to a $\C$\coacyclic{} complex of
$\S$\semimodule s $\bM^\bu$ homotopy equivalent to a complex
of $A$\+flat $\S$\semimodule s.
 On the other hand, if a complex of $\S$\semimodule s $\bM^\bu$ is
induced from a complex of $\C$\comodule s, then the cone of
the morphism of complexes $\boL_3(\bM^\bu)\rarrow\bM^\bu$ is
a contractible complex of $\S$\semimodule s, since it is isomorphic
to the cotensor product over~$\C$ of the bar complex $\dsb\rarrow
\S\oc_\C\S\oc_\C\S\rarrow \S\oc_\C\S\rarrow \S$, which is contractible
as a complex of left $\S$\semimodule s with right $\C$\comodule{}
structures, and a certain complex of left $\C$\comodule s.
 So the same applies to any complex of $\S$\semimodule s $\bM^\bu$
that up to the homotopy equivalence can be obtained from complexes
of $\S$\semimodule s induced from complexes of $\C$\comodule s using
the operations of cone and infinite direct sum.
 Part~(a) for quite $\S/\C/A$\semiflat{} complexes is proven;
the proofs of part~(a) for quite $\S/\C/A$\semiprojective{} complexes
and part~(b) are completely analogous.
\end{proof}

\begin{thm3}
 \textup{(a)} Assume that\/ $\C$ is a flat right $A$\module, $\S$ is
a coflat right\/ $\C$\comodule, $\D$ is a flat right $B$\module, and\/
$\T$ is a coflat right\/ $\D$\comodule.
 Then the functor\/ $\bM^\bu\mpsto{}_\T\bM^\bu$ maps quite\/
$\S/\C/A$\semiflat{} (quite\/ $\S/\C/A$\semiprojective) complexes of
left\/ $\S$\semimodule s to quite\/ $\T/\D/B$\semiflat{} (quite\/
$\T/\D/B$\semiprojective) complexes of left\/ $\T$\semimodule s.
 Assume additionally that\/ $\C$ and\/ $\S$ are flat left $A$\module s
and the ring $A$ has a finite weak homological dimension.
 Then the same functor maps\/ $\C$\coacyclic{} quite\/
$\S/\C/A$\semiflat{} complexes of left\/ $\S$\semimodule s to\/
$\D$\coacyclic{} complexes of left\/ $\T$\semimodule s. \par
 \textup{(b)} Assume that\/ $\C$ is a projective left $A$\module,
$\S$ is a coprojective left\/ $\C$\comodule, $\D$ is a projective
left $B$\module, and\/ $\T$ is a coprojective left\/ $\D$\comodule.
 Then the functor\/ $\bP^\bu\mpsto{}^\T\bP^\bu$ maps quite\/
$\S/\C/A$\semiinjective{} complexes of left\/ $\S$\semicontramodule s
to quite\/ $\T/\D/B$\semiinjective{} complexes of left
$\T$\semicontramodule s.
 Assume additionally that\/ $\C$ and\/ $\S$ are flat right $A$\module s
and the ring $A$ has a finite left homological dimension.
 Then the same functor maps\/ $\C$\contraacyclic{} quite\/
$\S/\C/A$\semiinjective{} complexes of left\/ $\S$\semicontramodule s
to\/ $\D$\contraacyclic{} complexes of left\/ $\T$\semicontramodule s.
\par
 \textup{(c)} Assume that\/ $\C$ is a projective left and a flat right
$A$\module, $\S$ is a coprojective left and a coflat right\/
$\C$\comodule, $A$ has a finite left homological dimension, $\D$ is
a projective left and a flat right $B$\module, $\T$ is a coprojective
left and a coflat right\/ $\D$\comodule, and $B$ has a finite left
homological dimension.
 Then the functor\/ $\bM^\bu\mpsto{}_\T\bM^\bu$ maps\/
$\S/\C/A$\projective{} complexes of left\/ $\S$\semimodule s to\/
$\T/\D/B$\projective{} complexes left\/ $\T$\semimodule s and
the functor\/ $\bP^\bu\mpsto{}^\T\bP^\bu$ maps\/ $\S/\C/A$\injective{}
complexes of left\/ $\S$\semicontramodule s to\/ $\T/\D/B$\injective{}
complexes of left\/ $\T$\semicontramodule s.
 The same functors map\/ $\C$\coacyclic\/ $\S/\C/A$\projective{}
complexes of left\/ $\S$\semimodule s to\/ $\D$\coacyclic{} complexes
of left\/ $\T$\semimodule s and\/ $\C$\contraacyclic\/
$\S/\C/A$\injective{} complexes of left\/ $\S$\semicontramodule s
to\/ $\D$\contraacyclic{} complexes of left\/ $\T$\semicontramodule s.
\end{thm3}

\begin{proof}
 Part~(a): the functor $\bM\mpsto{}_\T\bM$ maps the $\S$\semimodule{}
induced from a $\C$\comodule{} $\L$ to the $\T$\semimodule{} induced
from the $\D$\comodule{} ${}_B\L$.
 The first assertion follows immediately; to prove the second one,
use Theorem~\ref{co-contra-push-well-defined}(a) and Theorem~2(a).
 The proof of part~(b) is completely analogous.
 Part~(c): the first assertion follows from the adjointness of functors
$\bM^\bu\mpsto{}_\T\bM^\bu$ and $\bN^\bu\mpsto{}_\C\bN^\bu$,
the adjointness of functors $\bP^\bu\mpsto{}^\T\bP^\bu$ and
$\bQ^\bu\mpsto{}^\C\bQ^\bu$, and the second assertions of
Theorem~\ref{co-contra-pull-well-defined}(a) and~(b).
 The second assertion follows from the first assertions of
Theorem~\ref{co-contra-pull-well-defined}(a-b), because a complex of
left $\S$\semimodule s $\bM^\bu$ is $\S/\C/A$\projective{} and
$\C$\coacyclic{} if and only if the complex $\Hom_\S(\bM^\bu,\bL^\bu)$
is acyclic for all complexes of $\C/A$\injective{} left
$\S$\semimodule s $\bL^\bu$, and a complex of left
$\S$\semicontramodule s $\bP^\bu$ is $\S/\C/A$\injective{} and
$\C$\contraacyclic{} if and only if the complex
$\Hom^\S(\bR^\bu,\bP^\bu)$ is acyclic for all complexes of
$\C/A$\projective{} left $\S$\semicontramodule s $\bR^\bu$
(and analogously for complexes of $\T$\semimodule s and
$\T$\semicontramodule s).
 This follows from Theorem~\ref{semimodule-semicontramodule-subsect}
and the results of~\ref{semi-ctrtor-definition}, since a complex of
$\S$\semimodule s is $\C$\coacyclic{} iff it represents a zero object
of the semiderived category of $\S$\semimodule s, and a complex of
$\S$\semicontramodule s is $\C$\contraacyclic{} iff it represents
a zero object of the semiderived category of $\S$\semicontramodule s.
\end{proof}

\subsection{Derived functors of pull-back and push-forward}
\label{semi-pull-push-derived}
 Let $\S\rarrow\T$ be a map of semialgebras compatible with a map
of corings $\C\rarrow\D$ and a $k$\+algebra map $A\rarrow B$.

 Assume that $\C$ is a flat right $A$\module, $\S$ is a coflat right
$\C$\comodule, $\D$ is a flat right $B$\module, $\T$ is a coflat right
$\D$\comodule{} and a $\D/B$\+coflat left $\D$\comodule, and $B$ has
a finite weak homological dimension.
 The right derived functor
$$
 \bN^\bu\mpsto{}_\C^\boR\bN^\bu\: \sD^\si(\T\simodl)\lrarrow
 \sD^\si(\S\simodl)
$$
is defined by composing the functor $\bN^\bu\mpsto{}_\C\bN^\bu$ acting
from the homotopy category of left $\T$\semimodule s to the homotopy
category of left $\S$\semimodule s with the localization functor
$\Hot(\S\simodl)\rarrow\sD^\si(\S\simodl)$ and restricting it to
the full subcategory of complexes of $\D/B$\+coflat $\T$\semimodule s.
 By Theorems~\ref{semi-pull-push-adjusted}.1(a)
and~\ref{co-contra-pull-well-defined}(a), this restriction factorizes
through the semiderived category of left $\T$\semimodule s.

 Assume that $\C$ is a flat left and right $A$\module, $\S$ is a flat
left $A$\module{} and a coflat right $\C$\comodule, $A$ has a finite
weak homological dimension, $\D$ is a flat right $B$\module, and $\T$
is a coflat right $\D$\comodule.
 The left derived functor
$$
 \bM^\bu\mpsto{}_\T^\boL\bM^\bu\: \sD^\si(\S\simodl)\lrarrow
 \sD^\si(\T\simodl)
$$
is defined by composing the functor $\bM^\bu\mpsto{}_\T\bM^\bu$ acting
from the homotopy category of left $\S$\semimodule s to the homotopy
category of left $\T$\semimodule s with the localization functor
$\Hot(\T\simodl)\rarrow\sD^\si(\T\simodl)$ and restricting it to
the full subcategory of quite $\S/\C/A$\semiflat{} complexes of
$\S$\semimodule s.
 By Theorems~\ref{semi-pull-push-adjusted}.2(a)
and~\ref{semi-pull-push-adjusted}.3(a), this restriction factorizes
through the semiderived category of left $\S$\semimodule s.

 Analogously, assume that $\C$ is a projective left $A$\module,
$\S$ is a coprojective left $\C$\comodule, $\D$ is a projective left
$B$\module, $\T$ is a coprojective left $\D$\comodule{} and
a $\D/B$\+coflat right $\D$\comodule, and $B$ has a finite left
homological dimension.
 The left derived functor
$$
 \bQ^\bu\mpsto{}^\C_\boL\bQ^\bu\: \sD^\si(\T\sicntr)\lrarrow
 \sD^\si(\S\sicntr)
$$
is defined by composing the functor $\bQ^\bu\mpsto{}^\C\bQ^\bu$ with
the localization functor $\Hot(\S\sicntr)\rarrow\sD^\si(\S\sicntr)$
and restricting it to the full subcategory of complexes of
$\D/B$\coinjective{} $\T$\semicontramodule s.
 By Theorems~\ref{semi-pull-push-adjusted}.1(b)
and~\ref{co-contra-pull-well-defined}(b), this restriction factorizes
through the semiderived category of left $\T$\semicontramodule s.
 According to Lemma~\ref{semi-ctrtor-definition}.2, this definition
of a left derived functor does not depend on the choice of
a subcategory of adjusted complexes.

 Assume that $\C$ is a projective left and a flat right $A$\module,
$\S$ is a coprojective left $\C$\comodule{} and a flat right $A$\module,
$A$ has a finite left homological dimension, $\D$ is a projective left
$B$\module, and $\T$ is a coprojective left $\D$\comodule.
$$
 \bP^\bu\mpsto{}_\T^\boR\bP^\bu\: \sD^\si(\S\sicntr)\lrarrow
 \sD^\si(\T\sicntr)
$$
is defined by composing the functor $\bP^\bu\mpsto{}^\T\bP^\bu$ with
the localization functor $\Hot(\T\simodl)\rarrow\sD^\si(\T\simodl)$ and
restricting it to the full subcategory of quite $\S/\C/A$\semiflat
complexes of $\S$\semicontramodule s.
 By Theorems~\ref{semi-pull-push-adjusted}.2(b)
and~\ref{semi-pull-push-adjusted}.3(b), this restriction factorizes
through the semiderived category of left $\S$\semicontramodule s.
 According to Lemma~\ref{semi-ctrtor-definition}.2, this definition
of a right derived functor does not depend on the choice of
a subcategory of adjusted complexes.

 Notice that in the assumptions of
Theorem~\ref{semi-pull-push-adjusted}.3(c) above and Corollary~1(c)
below one can also define the left derived functor
$\bM^\bu\mpsto{}_\T^\boL\bM^\bu$ in terms of $\S/\C/A$\projective{}
complexes of left $\S$\semimodule s and the right derived functor
$\bP^\bu\mpsto{}_\T^\boR\bP^\bu$ in terms of $\S/\C/A$\injective{}
complexes of left $\S$\semicontramodule s.

 The derived functors $\bN^\bu\mpsto{}_\C^\boR\bN^\bu$ and
$\bQ^\bu\mpsto{}^\C_\boL\bQ^\bu$ in the categories of semimodules
and semicontramodules agree with the derived functors
$\N^\bu\mpsto{}_\C^\boR\N^\bu$ and $\Q^\bu\mpsto{}^\C_\boL\Q^\bu$
in the categories of comodules and contramodules, so our notation
is not ambiguous.

\begin{rmk1}
 Under the assumptions that $\C$ is a flat right $A$\module,
$\S$ is a coflat right $\C$\comodule, $\D$ is a flat right $B$\module,
$\T$ is a coflat right $\D$\comodule, and $B$ has a finite left
homological dimension, one can define the derived functor
$\bN^\bu\mpsto{}_\C^\boR\bN^\bu$ in terms of injective complexes of
left $\T$\semimodule s (see Remark~\ref{semi-ctrtor-definition}).
\end{rmk1}

\begin{cor1}
 \textup{(a)} The derived functor\/ $\bM^\bu\mpsto{}_\T^\boL\bM^\bu$ is
left adjoint to the derived functor\/ $\bN^\bu\mpsto{}_\C^\boR\bN^\bu$
whenever both functors are defined by the above construction. \par
 \textup{(b)} The derived functor\/ $\bP^\bu\mpsto{}^\T_\boR\bP^\bu$ is
right adjoint to the derived functor\/ $\bQ^\bu\mpsto{}^\C_\boL\bQ^\bu$
whenever both functors are defined by the above construction. \par
 \textup{(c)} Assume that $\C$ is a projective left and a flat right
$A$\module, $\S$ is a coprojective left and a coflat right
$\C$\comodule, $A$ has a finite left homological dimension, $\D$ is
a projective left and a flat right $B$\module, $\T$ is a coprojective
left and a coflat right $\D$\comodule, and $B$ has a finite left
homological dimension.
 Then for any objects\/ $\bM^\bu$ in\/ $\sD^\si(\simodr\S)$
and\/ $\bQ^\bu$ in\/ $\sD^\si(\T\sicntr)$ there is a natural
isomorphism\/ $\CtrTor^\T(\bM_\T^{\bu\boL},\bQ^\bu)\simeq
\CtrTor^\S(\bM^\bu,{}^\C_\boL\bQ^\bu)$ in the derived category of
$k$\module s.
\end{cor1}

\begin{proof}
 In the assumptions of part~(c), one can prove somewhat stronger
versions of the assertions~(a) and~(b): for any $\bM^\bu$ in
$\sD^\si(\S\simodl)$ and $\bN^\bu$ in $\sD^\si(\T\simodl)$, there is
a natural isomorphism $\Ext_\T({}_\T^\boL\bM,\bN)\simeq
\Ext_\S(\bM,{}_\C^\boR\bN)$ and for any $\bP^\bu$ in
$\sD^\si(\S\sicntr)$ and $\bQ^\bu$ in $\sD^\si(\T\sicntr)$ there is
a natural isomorphism $\Ext^\T(\bQ^\bu,{}^\T_\boR\bP^\bu)\simeq
\Ext^\S({}^\C_\boL\bQ^\bu,\bP^\bu)$ in the derived category of
$k$\module s.
 To obtain the first isomorphism, it suffices to represent the object
$\bM^\bu$ by an $\S/\C/A$\projective{} complex of left
$\S$\semimodule s and the object $\bN^\bu$ by a complex of
$\D/B$\injective{} left $\T$\semimodule s, and use
Lemma~\ref{rel-inj-proj-co-contra-mod}.2(a),
Theorem~\ref{co-contra-pull-well-defined}(a), and
Theorem~\ref{semi-pull-push-adjusted}.3(c).
 In the second case, one can represent the object $\bP^\bu$ by
an $\S/\C/A$\injective{} complex of left $\S$\semicontramodule s
and the object $\bQ^\bu$ by a complex of $\D/B$\projective{} left
$\T$\semicontramodule s, and use
Lemma~\ref{rel-inj-proj-co-contra-mod}.2(b),
Theorem~\ref{co-contra-pull-well-defined}(b), and
Theorem~\ref{semi-pull-push-adjusted}.3(c).
 To verify part~(c), it suffices to represent the object $\bM^\bu$ by
a quite $\S/\C/A$\semiflat{} complex of right $\S$\semimodule s and
the object $\bQ^\bu$ by a complex of $\D/B$\projective{} left
$\S$\semicontramodule s, and use
Lemma~\ref{rel-inj-proj-co-contra-mod}.2(b),
Theorem~\ref{co-contra-pull-well-defined}(b), and
Theorem~\ref{semi-pull-push-adjusted}.3(a).
 Finally, parts~(a) and~(b) in their weaker assumptions follow
from the next Lemma.
\end{proof}

\begin{lem}
 Let\/ $\sH_1$ and $\sH_2$ be categories, $\sS_1$ and\/ $\sS_2$ be
localizing classes of morphisms in\/ $\sH_1$ and\/ $\sH_2$, and\/
$\sF_1$ and\/ $\sF_2$ be full subcategories in\/ $\sH_1$ and\/ $\sH_2$.
 Assume that for any object $X\in\sH_1$ there exists an object
$U\in\sF_1$ together with a morphism $U\rarrow X$ from\/~$\sS_1$
and for any object $Y\in\sH_2$ there exists an object $V\in\sF_2$
together with a morphism $Y\rarrow V$ from\/~$\sS_2$.
 Let\/ $\Sigma\:\sH_1\rarrow\sH_2$ be a functor and\/ $\Pi\:
\sH_2\rarrow\sH_1$ be a functor right adjoint to\/~$\Sigma$.
 Assume that the morphism\/ $\Sigma(t)$ belongs to\/ $\sS_2$ for any
morphism $t\in\sF_1\cap\sS_1$ and the morphism\/ $\Pi(s)$ belongs to\/
$\sS_1$ for any morphism $s\in\sF_2\cap\sS_2$.
 Then the right derived functor\/ $\boR\Pi\:\sH_2[\sS_2^{-1}]\rarrow
\sH_1[\sS_1^{-1}]$ defined by restricting\/ $\Pi$ to\/~$\sF_2$ is right
adjoint to the left derived functor\/ $\boL\Sigma\:\sH_1[\sS_1^{-1}]
\rarrow\sH_2[\sS_2^{-1}]$ defined by restricting\/ $\Sigma$
to\/~$\sF_1$.
\end{lem}

\begin{proof}
 The functors $\sF_i[(\sF_i\cap\sS_i)^{-1}]\rarrow\sH_i[\sS_i^{-1}]$
are equivalences of categories by Lemma~\ref{semitor-main-theorem},
so the derived functors $\boL\Sigma$ and $\boR\Pi$ can be defined.
 For any objects $U\in\sF_1$ and $V\in\sF_2$ we have to construct
a bijection between the sets $\Hom_{\sH_1[\sS_1^{-1}]}(U,\Pi V)$
and $\Hom_{\sH_2[\sS_2^{-1}]}(\Sigma U, V)$, functorial in $U$ and~$V$.
 Any element of the first set can be represented by a fraction
$U\from U'\to \Pi V$ in $\sH_1$ with the morphism $U'\rarrow U$
belonging to~$\sS_1$.
 By assumption, one can choose $U'$ to be an object of~$\sF_1$.
 Assign to this fraction the element of the second set represented
by the fraction $\Sigma U\from \Sigma U' \to V$.
 By assumption, the morphism $\Sigma U'\rarrow \Sigma U$ belongs
to~$\sS_2$.
 Analogously, any element of the second set can be represented by
a fraction $\Sigma U\to V'\from V$ in $\sH_2$ with the morphism
$V\rarrow V'$ belongning to~$\sS_2$, and one can choose $V'$ to be
an object of $\sF_2$.
 Assign to this fraction the element of the first set represented
by the fraction $U\to \Pi V'\from \Pi V$.
 The compositions of these two maps between sets of morphisms are
identities, since the square formed by the morphisms $U'\rarrow U$,
\ $U\rarrow\Pi V'$, \ $U'\rarrow\Pi V$, and $\Pi V\rarrow\Pi V'$
and the square formed by the morphisms $\Sigma U'\rarrow \Sigma U$, \
$\Sigma U\rarrow V'$, \ $\Sigma U'\rarrow V$, and $V\rarrow V'$
are commutative simultaneously.
\end{proof}

 Let $\bR$ be a semialgebra over a coring $\E$ over a $k$\+algebra $F$,
and $\T\rarrow\bR$ be a map of semialgebras compatible with a map of
corings $\D\rarrow\E$ and a $k$\+algebra map $B\rarrow F$.
 Then the composition provides a map of semialgebras $\S\rarrow\bR$
compatible with a map of corings $\C\rarrow\E$ and a $k$\+algebra
map $A\rarrow B$.

\begin{cor2}
 \textup{(a)} There is a natural isomorphism\/
${}_\C^\boR({}_\D^\boR\bL^\bu)\simeq{}_\C^\boR\bL^\bu$ for any object\/
$\bL^\bu$ in\/ $\sD^\si(\bR\simodl)$ whenever both functors
$\bL^\bu\mpsto{}_\D^\boR\bL^\bu$ and\/ $\bN^\bu\mpsto{}_\C^\boR\bN^\bu$
are defined by the above construction. \par
 \textup{(b)} There is a natural isomorphism\/
${}_\bR^\boL({}_\T^\boL\bM^\bu)\simeq{}_\bR^\boL\bM^\bu$ for any
object\/ $\bM^\bu$ in\/ $\sD^\si(\S\simodl)$ whenever both functors
$\bM^\bu\mpsto{}_\T^\boL\bM^\bu$ and\/ $\bN^\bu\mpsto
{}_\bR^\boL\bN^\bu$ are defined by the above construction. \par
 \textup{(c)} There is a natural isomorphism\/
${}^\C_\boL({}^\D_\boL\bgK^\bu)\simeq{}^\C_\boL\bgK^\bu$ for any
object\/ $\bgK^\bu$ in\/ $\sD^\si(\bR\sicntr)$ whenever both functors
$\bgK^\bu\mpsto{}^\D_\boL\bgK^\bu$ and\/ $\bQ^\bu\mpsto
{}^\C_\boL\bQ^\bu$ are defined by the above construction. \par
 \textup{(d)} There is a natural isomorphism\/
${}^\bR_\boR({}^\T_\boR\bP^\bu)\simeq{}^\bR_\boR\bP^\bu$ for any
object\/ $\bP^\bu$ in\/ $\sD^\si(\S\sicntr)$ whenever both functors
$\bP^\bu\mpsto{}^\T_\boR\bP^\bu$ and\/ $\bQ^\bu\mpsto
{}^\bR_\boR\bQ^\bu$ are defined by the above construction.
\end{cor2}

\begin{proof}
 Part~(a) follows from the first assertion of
Theorem~\ref{co-contra-pull-well-defined}(a), part~(b) follows from
the first assertion of Theorem~\ref{semi-pull-push-adjusted}.3(a),
part~(c) follows from the first assertion of
Theorem~\ref{co-contra-pull-well-defined}(b), part~(d) follows from
the first assertion of Theorem~\ref{semi-pull-push-adjusted}.3(b).
\end{proof}

 Recall that a complex of $\C$\+coflat right $\S$\semimodule s
is called quite semiflat if it belongs to the minimal triangulated
subcategory of the homotopy category of right $\S$\semimodule s
containing the complexes of $\S$\semimodule s induced from complexes
of coflat right $\C$\comodule s and closed under infinite direct sums
(see~\ref{remarks-derived-semitensor-bi}).
 This definition presumes that $\C$ is a flat right $A$\module{} and
$\S$ is a coflat right $\C$\comodule.

\begin{cor3}
 \textup{(a)} Assume that\/ $\C$ is a flat left and right $A$\module,
$\S$ is a coflat left and right\/ $\C$\comodule, $A$ has a finite weak
homological dimension, $\D$ is a flat left and right $B$\module,
$\T$ is a coflat left and right\/ $\D$\comodule, and $B$ has a finite
weak homological dimension.
 Then for any objects\/ $\bM^\bu$ in\/ $\sD^\si(\S\simodl)$ and\/
$\bN^\bu$ in\/ $\sD^\si(\simodr\T)$ there is a natural isomorphism\/
$\SemiTor^\T(\bN^\bu,{}_\T^\boL\bM^\bu)\simeq
\SemiTor^\S(\bN_\C^{\bu\boR},\bM^\bu)$ in\/ $\sD(k\modl)$. \par
 \textup{(b)} Under the assumptions of Corollary~1(c),
for any objects\/ $\bP^\bu$ in\/ $\sD^\si(\S\sicntr)$ and\/
$\bN^\bu$ in\/ $\sD^\si(\T\simodl)$ there is a natural isomorphism\/
$\SemiExt_\T(\bN^\bu,{}^\T_\boR\bP^\bu)\simeq
\SemiExt_\S({}_\C^\boR\bN^\bu,\bP^\bu)$ in\/ $\sD(k\modl)$. \par
 \textup{(c)} Under the assumptions of Corollary~1(c),
for any objects\/ $\bM^\bu$ in $\sD^\si(\S\simodl)$ and\/
$\bQ^\bu$ in\/ $\sD^\si(\T\sicntr)$ there is a natural isomorphism\/
$\SemiExt_\T({}_\T^\boL\bM^\bu,\bQ^\bu)\simeq
\SemiExt_\S(\bM^\bu,{}^\C_\boL\bQ^\bu)$ in $\sD(k\modl)$.
\end{cor3}

\begin{proof}
 Part~(a): represent the object $\bM^\bu$ by a quite semiflat complex
of $\S$\semimodule s and the object $\bN^\bu$ by a semiflat complex
of $\D$\+coflat $\T$\semimodule s, and use the second case of
Proposition~\ref{pull-push-semitensor}(a).
 Alternatively, represent $\bM^\bu$ by a quite $\S/\C/A$\semiflat{}
complex of $A$\+flat $\S$\semimodule s and $\bN^\bu$ by a complex of
$\D$\+coflat $\T$\semimodule s, and use
Theorem~\ref{co-contra-pull-well-defined}(a),
Theorem~\ref{semi-pull-push-adjusted}.3(a),
the result of~\ref{relatively-semiflat}, and the first case of
Proposition~\ref{pull-push-semitensor}(a);
or represent $\bM^\bu$ by a quite semiflat complex of semiflat
$\S$\semimodule s and $\bN^\bu$ by a complex of $\D/B$\+coflat
$\T$\semimodule s, and use the same Theorems, the result
of~\ref{relatively-semiflat}, and the fourth case of
Proposition~\ref{pull-push-semitensor}(a).

 Part~(b): represent the object $\bP^\bu$ by a semiinjective complex
of $\C$\coinjective{} $\S$\semicontramodule s (having in mind
Lemma~\ref{birelatively-adjusted}(c) or
Remark~\ref{birelatively-adjusted}) and the object $\bN^\bu$
by a semiprojective complex of $\D$\coprojective{} $\T$\semimodule s,
and use the second case of Proposition~\ref{pull-push-semitensor}(b).
 Alternatively, represent $\bP^\bu$ by a quite
$\S/\C/A$\semiinjective{} complex of $A$\injective{}
$\S$\semicontramodule s and $\bN^\bu$ by a complex of
$\D$\coprojective{} $\T$\semimodule s, and use
Theorem~\ref{co-contra-pull-well-defined}(a),
Theorem~\ref{semi-pull-push-adjusted}.3(b),
the result of~\ref{relatively-semiproj-semiinj}, and the first case of
Proposition~\ref{pull-push-semitensor}(b);
or represent $\bP^\bu$ by a semiinjective complex of semiinjective
$\S$\semicontramodule s and $\bN^\bu$ by a complex of
$\D/B$\coprojective{} $\T$\semimodule s, and use the same Theorems,
the result of~\ref{relatively-semiproj-semiinj}, and the fourth case of
Proposition~\ref{pull-push-semitensor}(b).

 Part~(c): represent the object $\bM^\bu$ by a semiprojective complex
of $\C$\coprojective{} $\S$\semicontramodule s (having in mind
Lemma~\ref{birelatively-adjusted}(b) or
Remark~\ref{birelatively-adjusted}) and the object $\bQ^\bu$ by
a semiinjective complex of $\D$\coinjective{} $\T$\semicontramodule s,
and use the second case of Proposition~\ref{pull-push-semitensor}(c).
 Alternatively, represent $\bM^\bu$ by a quite
$\S/\C/A$\semiprojective{} complex of $A$\projective{}
$\S$\semimodule s and $\bQ^\bu$ by a complex of $\D$\coinjective{}
$\T$\semicontramodule s, and use
Theorem~\ref{co-contra-pull-well-defined}(b),
Theorem~\ref{semi-pull-push-adjusted}.3(a),
the result of~\ref{relatively-semiproj-semiinj}, and the first case of
Proposition~\ref{pull-push-semitensor}(c);
or represent $\bM^\bu$ by a semiprojective complex of semiprojective
$\S$\semimodule s and $\bQ^\bu$ by a complex of $\D/B$\coinjective{}
$\T$\semicontramodule s, and use the same Theorems, the result
of~\ref{relatively-semiproj-semiinj}, and the fourth case of
Proposition~\ref{pull-push-semitensor}(c).
\end{proof}

{\emergencystretch=0.5em\hfuzz=10pt
\begin{rmk2}
 Suppose that two objects ${}'\bM^\bu$ in $\sD^\si(\simodr\S)$ and
${}'\bN^\bu$ in $\sD^\si(\simodr\T)$ are endowed with a morphism
${}'\bM_\T^{\bu\boL} \rarrow{}'\bN^\bu$, or, which is the same,
a morphism ${}'\bM^\bu \rarrow{}'\bN_\C^{\bu\boR}$, and two objects
${}''\bM^\bu$ in $\sD^\si(\S\simodl)$ and ${}''\bN^\bu$ in
$\sD^\si(\T\simodl)$ are endowed with a morphism
${}_{\phantom{\prime\prime}\T}^{\boL\.\prime\prime}\bM^\bu\rarrow
{}''\bN^\bu$, or, which is the same, a morphism ${}''\bM^\bu\rarrow
{}_{\phantom{\.\prime\prime}\C}^{\boR\,\prime\prime}\bN^\bu$.
 Then the two morphisms
$\SemiTor^\S({}'\bM^\bu,{}''\bM^\bu)\rarrow
\SemiTor^\T({}'\bN^\bu,{}''\bN^\bu)$ in $\sD($k$\modl)$ provided by
the compositions $\SemiTor^\S({}'\bM^\bu,{}''\bM^\bu)\rarrow
\SemiTor^\S({}'\bN_\C^{\bu\boR},\allowbreak{}''\bM^\bu)\simeq
\SemiTor^\T({}'\bN^\bu,{}_{\phantom{\prime\prime}\T}^
{\boL\.\prime\prime}\bM^\bu) \rarrow\SemiTor^\T({}'\bN^\bu,{}''\bN^\bu)$
and $\SemiTor^\S({}'\bM^\bu,{}''\bM^\bu)\rarrow \SemiTor^\S
({}'\bM^\bu,{}_{\phantom{\.\prime\prime}\C}^{\boR\,\prime\prime}
\bN^\bu) \simeq\SemiTor^\T({}'\bM_\T^{\bu\boL},{}''\bN^\bu)\rarrow
\SemiTor^\T({}'\bN^\bu,{}''\bN^\bu)$ coincide with each other.
 Indeed, let us represent the objects ${}'\bM^\bu$ and ${}'\bN^\bu$
by complexes of right $\S$\semimodule s and $\T$\semimodule s in
such a way that the adjoint morphisms ${}'\bM_\T^{\bu\boL}\rarrow
{}'\bN^\bu$ and ${}'\bM^\bu \rarrow{}'\bN_\C^{\bu\boR}$ could be
represented by a map of complexes of semimodule s ${}'\bM^\bu\rarrow
{}'\bN^\bu$ compatible with the maps $A\rarrow B$, \ $\C\rarrow\D$,
and $\S\rarrow\T$.
 Applying to the complexes of ${}'\bM^\bu$ and ${}'\bN^\bu$
simultaneously the constructions from the proof of
Theorem~\ref{semitor-main-theorem}, one can construct a map of quite
semiflat complexes of right semimodules $\boL_3\boR_2\boL_1({}'\bM^\bu)
\rarrow\boL_3\boR_2\boL_1({}'\bN^\bu)$ representing the same adjoint
morphisms in the semiderived categories of left semimodules.
 So one can assume ${}'\bM\bu$ and ${}'\bN^\bu$ to be quite semiflat
complexes.
 Analogously, represent the morphisms ${}_{\phantom{\prime\prime}\T}^
{\boL\.\prime\prime}\bM^\bu\rarrow{}''\bN^\bu$ and ${}''\bM^\bu\rarrow
{}_{\phantom{\.\prime\prime}\C}^{\boR\,\prime\prime}\bN^\bu$ in
the semiderived categories of left semimodules by a map of quite
semiflat complexes of left semimodules ${}''\bM^\bu\rarrow{}''\bN^\bu$
compatible with the maps $A\rarrow B$, \ $\C\rarrow\D$, and
$\S\rarrow\T$.
 Then both compositions in question are represented by the same
map of complexes of $k$\module s ${}'\bM^\bu\os_\S{}''\bM^\bu\rarrow
{}'\bN^\bu\os_\T{}''\bN^\bu$.
 Furthermore, suppose that two objects $\bM^\bu$ in
$\sD^\si(\S\simodl)$ and $\bN^\bu$ in $\sD^\si(\T\simodl)$ are endowed
with a morphism ${}_\T^\boL\bM^\bu\rarrow\bN^\bu$, or, which is
the same, a morphism $\bM^\bu \rarrow{}_\C^\boR\bN^\bu$, and two
objects $\bP^\bu$ in $\sD^\si(\S\sicntr)$ and $\bQ^\bu$ in
$\sD^\si(\T\sicntr)$ are endowed with a morphism $\bQ^\bu\rarrow
{}^\T_\boR\bP^\bu$, or, which is the same, a morphism
${}^\C_\boL\bQ^\bu\rarrow\bP^\bu$.
 Then the two morphisms $\SemiExt_\T(\bN^\bu,\bQ^\bu)\rarrow
\SemiExt_\S(\bM^\bu,\bP^\bu)$ in $\sD($k$\modl)$ provided by
the compositions $\SemiExt_\T(\bN^\bu,\bQ^\bu)\rarrow\SemiExt_\T
(\bN^\bu, {}_\T^\boR\bP^\bu)\simeq\SemiExt_\S({}_\C^\boR\bN^\bu,
\bP^\bu) \rarrow\SemiExt_\S(\bM^\bu,\bP^\bu)$ and
$\SemiExt_\T(\bN^\bu,\bQ^\bu)\rarrow\SemiExt_\T({}_\T^\boL\bM^\bu,
\bQ^\bu)\simeq\SemiExt_\S(\bM^\bu,\allowbreak{}^\C_\boL\bQ^\bu)
\rarrow\SemiExt_\S(\bM^\bu,\bP^\bu)$ coincide with each other.
\end{rmk2}}

\begin{cor4}
 Under the assumptions of Corollary~1(c),
the mutually inverse equivalences of categories\/
$\boR\Psi_\S\:\sD^\si(\S\simodl)\rarrow\sD^\si(\S\sicntr)$ and\/
$\boL\Phi_\S\:\sD^\si(\S\sicntr)\rarrow\sD^\si(\S\simodl)$ and
the mutually inverse equivalences of categories\/
$\boR\Psi_\T\:\sD^\si(\T\simodl)\allowbreak\rarrow\sD^\si(\T\sicntr)$
and\/ $\boL\Phi_\T\:\sD^\si(\T\sicntr)\rarrow\sD^\si(\T\simodl)$
transform the derived functor\/ $\bN^\bu\mpsto{}_\C^\boR\bN^\bu$
into the derived functor\/ $\bQ^\bu\mpsto{}^\C_\boL\bQ^\bu$.
\end{cor4}

\begin{proof}
 To construct the isomorphism $\boL\Phi_\S({}^\C_\boL\bQ^\bu)\simeq
{}_\C^\boR(\boL\Phi_\T\bQ^\bu)$, represent the object $\bQ^\bu$ by
a complex of $\D/B$\projective{} $\C$\contramodule s, and use
Lemma~\ref{rel-inj-proj-co-contra-mod}.2,
Theorem~\ref{co-contra-pull-well-defined}(b), and
the results of~\ref{pull-psi-phi} and~\ref{semi-pull-psi-phi}.
 To construct the isomorphism ${}^\C_\boL(\boR\Psi_\T\bN^\bu)\simeq
\boR\Psi_\S({}_\C^\boR\bN^\bu)$, represent the object $\bN^\bu$ by
a complex of $\D/B$\injective{} $\C$\comodule s, and use
Lemma~\ref{rel-inj-proj-co-contra-mod}.2,
Theorem~\ref{co-contra-pull-well-defined}(a), and
the results of~\ref{pull-psi-phi} and~\ref{semi-pull-psi-phi}.
 To show that these isomorphisms agree, it suffices to check that
for any adjoint morphisms $\boL\Phi_\T\bQ^\bu\rarrow\bN^\bu$ and
$\bQ^\bu\rarrow\boR\Psi_\T\bN^\bu$ in the semiderived categories of
$\T$\semimodule s and $\T$\semicontramodule s the compositions
$\boL\Phi_\S({}^\C_\boL\bQ^\bu)\rarrow{}_\C^\boR(\boL\Phi_\T\bQ^\bu)
\rarrow{}_\C^\boR\bN^\bu$ and ${}^\C_\boL\bQ^\bu\rarrow
{}^\C_\boL(\boR\Psi_\T\bN^\bu)\rarrow\boR\Psi_\S({}_\C^\boR\bN^\bu)$
are adjoint morphisms in the semiderived categories of
$\S$\semimodule s and $\S$\semicontramodule s.
 Here one can represent $\bN^\bu$ by a semiprojective complex of
$\D$\coprojective{} left $\T$\semimodule s and $\bQ^\bu$ by
a semiinjective complex of $\D$\coinjective{} left
$\T$\semicontramodule s (having in mind
Lemmas~\ref{cotensor-contratensor-assoc}
and~\ref{birelatively-adjusted}), and use a result
of~\ref{pull-psi-phi}.
\end{proof}

 Thus we have constructed three functors between the semiderived
categories $\sD^\si(\S\simodl)\simeq\sD^\si(\S\sicntr)$ and
$\sD^\si(\T\simodl)\simeq\sD^\si(\T\sicntr)$: the functor described
in Corollary~4, and two functors adjoint to it from the left and
from the right, described in Corollary~1.

\begin{rmk3}
 One can show that the isomorphisms of derived functors from
Corollary~\ref{semiext-and-ext} are compatible with
the change-of-semi\-algebra isomorphisms from Corollaries~1, 3,
and~4 in the following way.
 To check that the compositions of isomorphisms $\SemiExt_\T
({}_\T^\boL\bM^\bu,\boR\Psi_\T(\bN^\bu))\rarrow \Ext_\T
({}_\T^\boL\bM^\bu,\bN^\bu)\rarrow \Ext_\S(\bM^\bu,{}_\C^\boR\bN^\bu)$
and $\SemiExt_\T({}_\T^\boL\bM^\bu,\boR\Psi_\T(\bN^\bu))\rarrow
\SemiExt_\S(\bM^\bu,{}^\C_\boL(\boR\Psi_\T\bN^\bu))\rarrow
\SemiExt_\S(\bM^\bu,\allowbreak\boR\Psi_\S({}_\C^\boR\bN^\bu))\rarrow
\Ext_\S(\bM^\bu,{}_\C^\boR\bN^\bu)$ coincide, represent the object
$\bM^\bu$ by a semiprojective complex of semiprojective left
$\S$\semimodule s and the object $\bN^\bu$ by a complex of
$\D/B$\injective{} left $\T$\semimodule s, and use the result
of~\ref{relatively-semiproj-semiinj}.
 To check that the compositions of isomorphisms $\CtrTor^\S
(\bM^\bu,{}^\C_\boL\bQ^\bu)\rarrow\CtrTor^\T(\bM_\T^{\bu\boL},\bQ^\bu)
\rarrow\SemiTor^\T(\bM_\T^{\bu\boL},\boL\Phi_\T(\bQ^\bu))$
and $\CtrTor^\S(\bM^\bu,{}^\C_\boL\bQ^\bu)\rarrow
\SemiTor^\S(\bM^\bu,\boL\Phi_\S({}^\C_\boL\bQ^\bu))\allowbreak\rarrow
\SemiTor^\S(\bM^\bu,{}_\C^\boR(\boL\Phi_\T\bQ^\bu))\rarrow
\SemiTor^\T(\bM_\T^{\bu\boL},\boL\Phi_\T(\bQ^\bu))$ coincide, represent
the object $\bM^\bu$ by a quite semiflat complex of semiflat
right $\S$\semimodule s and the object $\bQ^\bu$ by a complex of
$\D/B$\projective{} left $\T$\semicontramodule s, and use
the result of~\ref{relatively-semiflat}.
 Commutativity of the respective diagrams on the level of abelian
categories is straightforward to verify under our assumptions on
the terms of the complexes representing the objects $\bM^\bu$.
\end{rmk3}

\subsection{Remarks on Morita morphisms}

\subsubsection{}  \label{morita-change-of-coring-construction}
 Let $\C$ be a coring over a $k$\+algebra $A$ and $\D$ be a coring
over a $k$\+algebra $B$ such that $\C$ is a flat right $A$\module{}
and $\D$ is a flat right $B$\module.
 Let $(\E,\E\dual)$ be a right coflat Morita morphism from $\C$ to
$\D$ and $\T$ be a semialgebra over the coring~$\D$ such that $\T$
is a coflat right $\D$\comodule.
 In this case, the semialgebra ${}_\C\T_\C$ over the coring~$\C$
is constructed in the following way.
 As a $\C$\+$\C$\bicomodule, ${}_\C\T_\C$ is equal to
$\E\oc_\D\T\oc_\D\E\dual$.
 The semimultiplication in ${}_\C\T_\C$ is defined as the composition
$\E\oc_\D\T\oc_\D\E\dual\oc_\C\E\oc_\D\T\oc_\D\E\dual\rarrow
\E\oc_\D\T\oc_\D\T\oc_\D\E\dual\rarrow\E\oc_\D\T\oc_\D\E\dual$
of the morphism induced by the morphism $\E\dual\oc_\C\E\rarrow\D$
and the morphism induced by the semimultiplication in~$\T$.
 The semiunit in ${}_\C\T_\C$ is defined as the composition
$\C\rarrow\E\oc_\D\E\dual\rarrow\E\oc_\D\T\oc_\D\E\dual$ of
the morphism induced by the morphism $\C\rarrow\E\oc_\D\E\dual$
and the morphism induced by the semiunit in~$\T$.

 For example, if $\C\rarrow\D$ is a map of corings compatible with
a $k$\+algebra map $A\rarrow B$ such that $B$ is a flat right
$A$\module{} and $\C_B$ is a coflat right $\D$\comodule, one can
take $\E=\C_B$ and $\E\dual={}_B\C$.
 Then the algebra ${}_\C\T_\C$ is a universal final object
in the category of semialgebras $\S$ over $\C$ endowed with a map
$\S\rarrow\T$ compatible with the maps $A\rarrow B$ and $\C\rarrow\D$.
 The semialgebra ${}_\C\T_\C=\C_B\oc_\D\T\oc_\D{}_B\C$ can be also
defined, e.~g., when $(E,E\dual)$ is a Morita morphism from
a $k$\+algebra $A$ to a $k$\+algebra $B$ and ${}_B\C_B=
E\dual\ot_A\C\ot_A E\rarrow\D$ is a morphism of corings over~$B$
such that $E\dual$ is a flat right $A$\module, ${}_B\C=E\dual\ot_A\C$
is a $\D/B$\+coflat left $\D$\comodule, $\T$ is a flat right
$B$\module{} and a $\D/B$\+coflat left $\D$\comodule, and the rings
$A$ and $B$ have finite weak homological dimensions.

 All the results
of~\ref{semi-compatible-morphisms}--\ref{semi-pull-push-derived}
can be extended to the situation of a left coprojective and right
coflat Morita morphism $(\E,\E\dual)$ from a coring~$\C$ to
a coring~$\D$ and a morphism $\S\rarrow{}_\C\T_\C$ of semialgebras
over~$\C$.
 In particular, when $\C$ is a flat right $A$\module, $\D$ is
a flat right $B$\module, $\S$ is a coflat right $\C$\comodule,
$\T$ is a coflat right $\D$\comodule, and $(\E,\E\dual)$ is a right
coflat Morita morphism, the functor $\bN\mpsto{}_\C\bN=\E\oc_\D\bN$
from the category of left $\T$\semimodule s to the category of left
$\S$\semimodule s has a left adjoint functor
$\bM\mpsto{}_\T\bM=\T_\C\os_\S\bM$.
 Analogously, when $\C$ is a projective left $A$\module, $\D$ is
a projective left $B$\module, $\S$ is a coprojective left
$\C$\comodule, $\T$ is a coprojective left $\D$\comodule, and
$(\E,\E\dual)$ is a left coprojective Morita morphism, the functor
$\bQ\mpsto{}^\C\bQ=\Cohom_\D(\E\dual,\bQ)$ from the category of
left $\T$\semicontramodule s to the category of left
$\S$\semicontramodule s has a right adjoint functor
$\bP\mpsto{}^\T\bP=\SemiHom_\S({}_\C\T,\bP)$, etc.
 However, one sometimes has to impose the homological dimension
conditions on $A$ and $B$ where they were not previously needed.

\subsubsection{}
 Assume that $\C$ is a flat right $A$\module{} and $\D$ is a flat
right $B$\module.
 A right $\D$\comodule{} $\K$ is called \emph{faithfully coflat} if it
is a coflat $\D$\comodule{} and for any nonzero left $\D$\comodule{}
$\M$ the cotensor product $\K\oc_\D\M$ is nonzero.
 A right coflat Morita morphism $(\E,\E\dual)$ from $\C$ to $\D$ is
called \emph{right faithfully coflat} if the right $\D$\comodule{} $\E$
is faithfully coflat.
 A right coflat Morita morphism $(\E,\E\dual)$ is right faithfully
coflat if and only if the right $\D$\comodule{} $\E\dual\oc_\C\E$
is faithfully coflat and if and only if the morphism
$\E\dual\oc_\C\E\rarrow\D$ is surjective and its kernel is a coflat
right $\D$\comodule.
 Indeed, the cotensor product $\E\oc_\D\M$ is nonzero if and only if
the morphism $\E\dual\oc_\C\E\oc_\D\M\rarrow\M$ is nonzero; this holds
for any nonzero left $\D$\comodule{} $\M$ if and only if the morphism
$\E\dual\oc_\C\E\oc_\D\M\rarrow\M$ is surjective for any left
$\D$\comodule{} $\M$, and it remains to use the results of (the proof
of) Lemma~\ref{absolute-relative-coflat}.
 
 Let $(\E,\E\dual)$ be a right faithfully coflat Morita morphism from
$\C$ to $\D$ and $\T$ be a semialgebra over the coring~$\D$ such that
$\T$ is a coflat right $\D$\comodule.
 Then the functor $\bN\mpsto{}_\C\bN$ is an equivalence of the abelian
categories of left $\T$\semimodule s and left ${}_\C\T_\C$\semimodule s.
 This follows from Theorem~\ref{barr-beck-theorem} applied to
the functor $\Delta\:\T\simodl\rarrow\C\comodl$ mapping
a $\T$\semimodule{} $\bN$ to the $\C$\comodule{} ${}_\C\bN$ and
the functor $\Gamma\:\C\comodl\rarrow\T\simodl$ left adjoint
to~$\Delta$ mapping a $\C$\comodule{} $\M$ to the $\T$\semimodule{}
$\T\oc_\D{}_\D\M$.

 Now assume that $\C$ is a projective left $A$\module{} and $\D$ is
a projective left $B$\module.
 A left $\D$\comodule{} $\K$ is called \emph{faithfully coprojective}
if it is a coprojective $\D$\comodule{} and for any nonzero left
$\D$\contramodule{} $\P$ the cohomomorphism module $\Cohom_\D(\K,\P)$
is nonzero.
 A faithfully coprojective $\D$\comodule{} is faithfully coflat.
 A left coprojective Morita morphism $(\E,\E\dual)$ from $\C$ to $\D$
is called \emph{left faithfully coprojective} if the left
$\D$\comodule{} $\E\dual$ is faithfully coprojective.
 A left coprojective Morita morphism $(\E,\E\dual)$ is left faithfully
coprojective if and only if the left $\D$\comodule{} $\E\dual\oc_\C\E$
is faithfully coflat and if and only if the morphism
$\E\dual\oc_\C\E\rarrow\D$ is surjective and its kernel is
a coprojective left $\D$\comodule.

 Let $(\E,\E\dual)$ be a left faithfully coprojective Morita morphism
from $\C$ to $\D$ and $\T$ be a semialgebra over the coring~$\D$
such that $\T$ is a coprojective left $\D$\comodule.
 Then the functor $\bQ\mpsto{}^\C\bQ$ is an equivalence of the abelian
categories of left $\T$\semicontramodule s and left
${}_\C\T_\C$\semicontramodule s.
 This follows from Theorem~\ref{barr-beck-theorem} applied to
the functor $\Delta\:\T\sicntr\rarrow\C\contra$ mapping
a $\T$\semicontramodule{} $\bQ$ to the $\C$\contramodule{} ${}^\C\bQ$
and the functor $\Gamma\:\C\contra\rarrow\T\sicntr$ right adjoint
to~$\Delta$ mapping a $\C$\contramodule{} $\P$ to
the $\T$\semicontramodule{} $\Cohom_\D(\T,{}^\D\P)$.

\subsubsection{}  \label{morita-change-of-coring}
 Assume that $\C$ is a flat right $A$\module{} and $\D$ is a flat right
$B$\module.
 Let $(\E,\E\dual)$ be a right coflat Morita morphism from $\C$ to $\D$
and $\T$ be a semialgebra over the coring~$\D$ such that $\T$ is
a coflat right $\D$\comodule.
 Then the functor $\bN^\bu\mpsto{}_\C\bN^\bu$ maps $\D$\coacyclic{}
complexes of $\T$\semimodule s to $\C$\coacyclic{} complexes of
${}_\C\T_\C$\semimodule s and the semiderived category of left
${}_\C\T_\C$\semimodule s is a localization of the semiderived
category of left $\T$\semimodule s by the kernel of the functor induced
by $\bN^\bu\mpsto{}_\C\bN^\bu$ (as one can check by computing
the functor $\bM^\bu\mpsto{}_\C({}_\T^\boL\bM^\bu)$ on the semiderived
category of left ${}_\C\T_\C$\semimodule s).
 The triangulated categories $\sD^\si(\T\simodl)$ and
$\sD^\si({}_\C\T_\C\simodl)$ are equivalent when $(\E,\E\dual)$ is
a right coflat Morita equivalence, or more generally when the morphism
$\E\dual\oc_\C\E\rarrow\D$ is an isomorphism.

 Analogously, assume that $\C$ is a flat right $A$\module{} and $\D$
is a projective left $B$\module.
 Let $(\E,\E\dual)$ be a left coprojective Morita morphism from $\C$
to $\D$ and $\T$ be a semialgebra over the coring~$\D$ such that $\T$
is a coprojective left $\D$\comodule.
 Then the functor $\bQ^\bu\mpsto{}^\C\bQ^\bu$ maps $\D$\contraacyclic{}
complexes of $\T$\semicontramodule s to $\C$\contraacyclic{} complexes
of ${}_\C\T_\C$\semicontramodule s and the semiderived category of
left ${}_\C\T_\C$\semicontramodule s is a localization of
the semiderived category of left $\T$\semicontramodule s by the kernel
of the functor induced by $\bQ^\bu\mpsto{}^\C\bQ^\bu$.
 The triangulated categories $\sD^\si(\T\sicntr)$ and
$\sD^\si({}_\C\T_\C\sicntr)$ are equivalent when $(\E,\E\dual)$ is
a left coprojective Morita equivalence, or more generally when
the morphism $\E\dual\oc_\C\E\rarrow\D$ is an isomorphism.

\begin{rmk}
 The semiderived categories of left $\T$\semimodule s and left
${}_\C\T_\C$\semimodule s can be different even when $(\E,\E\dual)$ is
a right faithfully coflat Morita morphism and the abelian categories of
left $\T$\semimodule s and left ${}_\C\T_\C$\semimodule s are 
equivalent.
 Indeed, let $A=B=k$ be a field and $F$ be a finite-dimensional algebra
over~$k$.
 Let $\D=F^*$ and $\C=\End(F)^*$ be the coalgebras over~$k$
dual to the finite-dimensional $k$\+algebras $F$ and $\End(F)$.
 Then there is a coalgebra morphism $\C\rarrow\D$ dual to the algebra
embedding $F\rarrow\End(F)$ related to the action of $F$ in itself
by left multiplications.
 Since $\End(F)$ is a free left $F$\module, $\C$ is a cofree right
$\D$\comodule.
 Set $\E=\C=\E\dual$; this is a right faithfully coprojective Morita
morphism from $\C$ to $\D$.
 Now put $\T=\D$; then the semiderived category of left
$\T$\semimodule s coincides with the coderived category of left
$\D$\comodule s.
 At the same time, the coalgebra $\C$ is semisimple and a complex of
$\C$\comodule s is coacyclic if and only if it is acyclic, so
the semiderived category of left ${}_\C\T_\C$\comodule s is equivalent
to the conventional derived category of left $\D$\comodule s.
 When $F$ is a Frobenius algebra, $\End(F)$ is a free left and right
$F$\module, so $(\E,\E\dual)$ is a left and right faithfully
coprojective Morita morphism, but the categories $\sD^\si(\T\simodl)$
and $\sD^\si({}_\C\T_\C\simodl)$ are still not equivalent when
the homological dimension of $F$ is infinite.
 Alternatively, one can consider the right coprojective Morita morphism
from the coalgebra $\C=k$ to the coalgebra $\D=F^*$ with $\E=F^*$
and $\E\dual=F$ and the same semialgebra $\T=\D$ over~$\D$; then
the semialgebra ${}_\C\T_\C$ over~$\C$ is isomorphic to the algebra $F$
over~$k$; the category $\sD^\si(\T\simodl)$ is the coderived category
of $F^*$\comodule s and the category $\sD^\si({}_\C\T_\C\simodl)$ is
the derived category of $F$\module s.
\end{rmk}

 Assume that $\C$ is a flat left and right $A$\module, $\D$ is a flat
left and right $B$\module, the rings $A$ and $B$ have finite weak
homological dimensions, $\T$ is a coflat left and right $\D$\comodule,
and $(\E,\E\dual)$ is a left and right coflat Morita morphism
from $\C$ to~$\D$.
 Then whenever the functor $\bN^\bu\mpsto{}_\C\bN^\bu$ induces
an equivalence of the semiderived categories of left $\T$\semimodule s
and left ${}_\C\T_\C$\semimodule s and the functor $\bN^\bu\mpsto
\bN^\bu_\C$ induces an equivalence of the semiderived categories of
right $\T$\semimodule s and right ${}_\C\T_\C$\semimodule s, these
equivalences of categories transform the functor $\SemiTor^\T$ into
the functor $\SemiTor^{{}_\C\T_\C}$.

 Assume that $\C$ is a projective left and a flat right $A$\module,
$\D$ is a projective left and a flat right $B$\module, the rings $A$
and $B$ have finite left homological dimensions, $\T$ is a coprojective
left and a coflat right $\D$\comodule, and $(\E,\E\dual)$ is
a left coprojective and right coflat Morita morphism from $\C$ to $\D$.
 Then whenever the functor $\bN^\bu\mpsto{}_\C\bN^\bu$ induces
an equivalence of the semiderived categories of left $\T$\semimodule s
and left ${}_\C\T_\C$\semimodule s and the functor $\bQ^\bu\mpsto
{}^\C\bQ^\bu$ induces an equivalence of the semiderived categories of
left $\T$\semicontramodule s and left ${}_\C\T_\C$\semicontramodule s,
these equivalences of categories transform the functor $\SemiExt_\T$ 
into the functor $\SemiExt_{{}_\C\T_\C}$ and the equivalences of
categories $\boR\Psi_\T$ and $\boL\Phi_\T$ into the equivalences of
categories $\boR\Psi_{{}_\C\T_\C}$ and $\boL\Phi_{{}_\C\T_\C}$.
 The same applies to the functors $\Ext_\T$, \ $\Ext^\T$, and
$\CtrTor^\T$, under the relevant assumptions.

\subsubsection{}   \label{relative-homological-dimension}
 Here are some further partial results about equivalence of
the semiderived categories related to $\T$ and ${}_\C\T_\C$.
 The problem is, essentially, to find conditions under which a complex
of left $\D$\comodule s $\N^\bu$ is coacyclic whenever the complex of
$\C$\comodule s ${}_\C\N^\bu$ is coacyclic, or a complex of left
$\D$\contramodule s $\Q^\bu$ is contraacyclic whenever the complex of
$\C$\contramodule s ${}^\C\Q^\bu$ is contraacyclic.

 Consider the following general setting.
 Let $\sA$ and $\sB$ be exact categories with exact functors of
infinite direct sum, $\Delta\:\sB\rarrow\sA$ be an exact functor
preserving infinite direct sums and such that a complex $C^\bu$
over~$\sB$ is acyclic if the complex $\Delta(C^\bu)$ over~$\sA$ is
contractible, and $\Gamma\:\sA\rarrow\sB$ be an exact functor
left adjoint to~$\Delta$.
 Clearly, if a complex $C^\bu$ is coacyclic then the complex
$\Delta(C^\bu)$ is coacyclic; we would like to know when the converse
holds.

 First, if a complex $C^\bu$ is coacyclic whenever the complex
$\Delta(C^\bu)$ is contractible, then a complex $C^\bu$ is coacyclic
if and only if the complex $\Delta(C^\bu)$ is coacyclic.
 Indeed, consider the bar bicomplex $\dsb\rarrow \Gamma\Delta\Gamma
\Delta(C^\bu)\rarrow\Gamma\Delta(C^\bu)\rarrow C^\bu$ whose
differentials are the alternating sums of morphisms induced by
the adjunction morphism $\Gamma\Delta\rarrow\Id$.
 The total complex of this bicomplex contructed by taking infinite
direct sums along the diagonals becomes contractible after applying
the functor $\Delta$; the contracting homotopy is induced by
the adjunction morphism $\Id\rarrow \Delta\Gamma$.
 By assumption, it follows that the total complex itself is coacyclic
over~$\sB$.
 On the other hand, if the complex $\Delta(C^\bu)$ is coacyclic
over~$\sA$, then every complex $(\Gamma\Delta)^n(C^\bu)$ is coacyclic
over $\sB$, since the functors $\Delta$ and $\Gamma$ are exact and
preserve infinite direct sums.
 The total complex of the bicomplex $\dsb\rarrow \Gamma\Delta\Gamma
\Delta(C^\bu)\rarrow\Gamma\Delta(C^\bu)$ is homotopy equivalent to
a complex obtained from the complexes $(\Gamma\Delta)^n(C^\bu)$ using
the operations of shift, cone, and infinite direct sum; hence
the complex $C^\bu$ is coacyclic.

 By the same argument, a complex $C^\bu$ is acyclic if and only if
the complex $\Delta(C^\bu)$ is acyclic, so if the exact category $\sB$
has a finite homological dimension, then a complex $C^\bu$ is
coacyclic if and only if the complex $\Delta(C^\bu)$ is coacyclic.
 This is the trivial case.

 Finally, let us say that an exact functor $\Delta\:\sB\rarrow\sA$ has
a finite relative homological dimension if the category $\sB$ with
the exact category structure formed by the exact triples in~$\sB$
that split after applying~$\Delta$ has a finite homological dimension.
 We claim that when the functor $\Delta$ has a finite relative
homological dimension, a complex $C^\bu$ over~$\sB$ is coacyclic if and
only if the complex $\Delta(C^\bu)$ is coacyclic, in our assumptions.
 Indeed, consider again the bar bicomplex $\dsb\rarrow \Gamma\Delta
\Gamma\Delta(C^\bu)\rarrow \Gamma\Delta(C^\bu)\rarrow C^\bu$.
 One can assume that the category $\sB$ contains images of idempotent
endomorphisms, as passing to the Karoubian closure doesn't change
coacyclicity.
 One can also assume that the complex $C^\bu$ is bounded from above,
as any acyclic complex bounded from below is coacyclic.
 The complex $\dsb\rarrow \Gamma\Delta\Gamma\Delta(X)\rarrow
\Gamma\Delta(X)$ is split exact in high enough (negative) degrees
for any object $X\in\sB$, since it is exact and the complex of
homomorphisms from it to an object $Y\in\sB$ computes $\Ext(X,Y)$
in the relative exact category.
 Let $d$ be an integer not smaller than the relative homological
dimension; denote by $Z(X)$ the image of the morphism
$(\Gamma\Delta)^{d+1}(X)\rarrow(\Gamma\Delta)^d(X)$.
 Then the total complex of the bicomplex $\dsb\rarrow
(\Gamma\Delta)^{d+2}(C^\bu)\rarrow(\Gamma\Delta)^{d+1}(C^\bu)\rarrow
Z(C^\bu)$ is contractible, while the total complex of the bicomplex
$(\Gamma\Delta)^d(C^\bu)/Z(C^\bu)\rarrow(\Gamma\Delta)^{d-1}(C^\bu)
\rarrow\dsb\rarrow \Gamma\Delta(C^\bu)\rarrow C^\bu$ is coacyclic.
 If the complex $\Delta(C^\bu)$ is coacyclic, the total complex of
the bicomplex $\dsb\rarrow \Gamma\Delta\Gamma\Delta(C^\bu)\rarrow
\Gamma\Delta(C^\bu)$ is also coacyclic; thus the complex $C^\bu$
is coacyclic.

\subsubsection{}  \label{semi-morita-morphisms}
 Let $\S$ be a semialgebra over a coring $\C$ and $\T$ be a semialgebra
over a coring~$\D$.
 Assume that $\C$ is a flat right $A$\module, $\D$ is a flat right
$B$\module, $\S$ is a coflat right $\C$\comodule, and $\T$ is a coflat
right $\D$\comodule.
 A \emph{right semiflat Morita morphism} from $\S$ to $\T$ is a pair
consisting of a $\T$\semiflat{} $\S$\+$\T$\bisemimodule{} $\bE$ and
an $\S$\semiflat{} $\T$\+$\S$\bisemimodule{} $\bE\dual$ endowed with
an $\S$\+$\S$\bisemimodule{} morphism $\S\rarrow\bE\os_\T\bE\dual$
and a $\T$\+$\T$\bisemimodule{} morphism $\bE\dual\os_\S\bE\rarrow\T$
such that the two compositions $\bE\rarrow\bE\os_\T\bE\dual\os_\S\bE
\rarrow\bE$ and $\bE\dual\rarrow\bE\dual\os_\S\bE\os_\T\bE\dual
\rarrow\bE\dual$ are equal to the identity endomorphisms of $\bE$
and $\bE\dual$.
 A right semiflat Morita morphism $(\bE,\bE\dual)$ from $\S$ to $\T$
induces an exact functor $\bM\mpsto{}_\T\bM = \bE\dual\os_\S\bM$ from
the category of left $\S$\semimodule s to the category of left
$\T$\semimodule s and an exact functor $\bN\mpsto{}_\S\bN =
\bE\os_\T\bN$ from the category of left $\T$\semimodule s to
the category of left $\S$\semimodule s; the former functor is
left adjoint to the latter one.
 Conversely, any pair of adjoint exact $k$\+linear functors preserving
infinite direct sums between the category of left $\S$\semimodule s and
left $\T$\semimodule s is induced by a right semiflat Morita morphism.
 Indeed, any exact $k$\+linear functor $\S\simodl\rarrow\T\simodl$
preserving infinite direct sums is the functor of semitensor product
with an $\S$\semiflat{} $\T$\+$\S$\bisemimodule; this can be proven as
in~\ref{co-contra-morita-morphisms}.  

 Analogously, assume that $\C$ is a projective left $A$\module,
$\D$ is a projective left $B$\module, $\S$ is a coprojective left
$\C$\comodule, and $\T$ is a coprojective left $\D$\comodule.
 A \emph{left semiprojective Morita morphism} from $\S$ to $\T$ is
defined as a pair consisting of an $\S$\semiprojective{}
$\S$\+$\T$\bisemimodule{} $\bE$ and a $\T$\semiprojective{}
$\T$\+$\S$\bisemimodule{} $\bE\dual$ endowed with
an $\S$\+$\S$\bisemimodule{} morphism $\S\rarrow\bE\os_\T\bE\dual$
and a $\T$\+$\T$\bisemimodule{} morphism $\bE\dual\os_\S\bE\rarrow\T$
satisfying the same conditions as above.
 A left semiprojective Morita morphism $(\bE,\bE\dual)$ from $\S$
to $\T$ induces an exact functor $\bP\mpsto{}^\T\bP =
\SemiHom_\S(\bE,\bP)$ from the category of left $\S$\semicontramodule s
to the category of left $\T$\semicontramodule s and an exact functor
$\bQ\mpsto{}^\S\bQ = \SemiHom_\T(\bE\dual,\bQ)$ from the category of
left $\T$\semicontramodule s to the category of left
$\S$\semicontramodule s; the former functor is right adjoint to
the latter one.
 Conversely, any pair of adjoint exact $k$\+linear functors preserving
infinite products between the category of left $\S$\semicontramodule s
and left $\T$\semicontramodule s is induced by a left semiprojective
Morita morphism.
 Indeed, any exact $k$\+linear functor $\S\sicntr\rarrow\T\sicntr$
preserving infinite products is the functor of semihomomorphisms
from an $\S$\semiprojective{} $\S$\+$\T$\bisemimodule.

 A right semiflat Morita morphism $(\bE,\bE\dual)$ from $\S$ to $\T$
is called a \emph{right semiflat Morita equivalence} if
the bisemimodule morphisms $\S\rarrow\bE\os_\T\bE\dual$ and
$\bE\dual\os_\S\bE\rarrow\T$ are isomorphisms; then it can be also
considered as a right semiflat Morita morphism $(\bE\dual,\bE)$
from $\T$ to~$\S$.
 \emph{Left semiprojective Morita equivalences} are defined in
the analogous way.
 A right semiflat Morita equivalence between semialgebras $\S$ and $\T$
induces an equivalence of the abelian categories of left
$\S$\semimodule s and left $\T$\semimodule s, and in the relevant
above assumptions any equivalence between these two $k$\+linear
categories comes from a right semiflat Morita equivalence.
 Analogously, a left semiprojective Morita equivalence between $\S$
and $\T$ induces an equivalence of the abelian categories of left
$\S$\semicontramodule s and left $\T$\semicontramodule s, and in
the relevant above assumptions any equivalence between these two
$k$\+linear categories comes from a left semiprojective Morita
equivalence.

 Assume that the coring $\C$ is a flat right $A$\module{} and
the coring $\D$ is a flat right $B$\module.
 Let $\T$ be a semialgebra over $\D$ such that $\T$ is a coflat right
$\D$\comodule{} and $(\E,\E\dual)$ be a right faithfully coflat Morita
morphism from $\C$ to $\D$.
 Then the pair of bisemimodules $\bE={}_\C\T$ and $\bE\dual=\T_\C$ is
a right semiflat Morita equivalence between the semialgebras
$\T$ and ${}_\C\T_\C$.
 Analogously, assume that $\C$ is a projective left $A$\module{} and
$\D$ is a projective left $B$\module.
 Let $\T$ be a semialgebra over $\D$ such that $\T$ is a coprojective
left $\D$\comodule{} and $(\E,\E\dual)$ be a left faithfully
coprojective Morita morphism from $\C$ to $\D$.
 Then the same pair of bisemimodules $\bE$ and $\bE\dual$ is a left
semiprojective Morita equivalence between $\T$ and ${}_\C\T_\C$.

 All the results of~\ref{semi-compatible-morphisms} can be extended
to the case of a left semiprojective and right semiflat Morita morphism
$(\bE,\bE\dual)$ from a semialgebra~$\S$ to a semialgebra~$\T$.
 In particular, for any left $\T$\semimodule{} $\bN$ there are
natural isomorphisms of left $\S$\semicontramodule s
$\Psi_\S({}_\S\bN)=\Hom_\S(\S,{}_\S\bN)\simeq\Hom_\T({}_\T\S,\bN)
\simeq\Hom_\T(\bE\dual,\bN)\simeq\SemiHom_\T(\bE\dual,\Hom_\T(\T,\bN))
={}^\S(\Psi_\T\bN)$ by
Proposition~\ref{semitensor-contratensor-assoc}.2(d), etc.
 However, one sometimes has to strengthen the coflatness
(coprojectivity, coinjectivity) conditions to the semiflatness
(semiprojectivity, semiinjectivity) conditions.

 The first assertions of Theorem~\ref{semi-pull-push-adjusted}.3(a),
(b) and~(c) do \emph{not} hold for Morita morphisms of semialgebras,
though. 
 The derived functors $\bM^\bu\mpsto{}_\T^\boL\bM^\bu$ and
$\bP^\bu\mpsto{}^\T_\boR\bP^\bu$ still can be defined in terms
of $\S/\C$\projective{} (\.= quite $\S/\C$\semiflat) complexes of
$\S$\semimodule s and $\S/\C$\injective{} (\.= quite
$\S/\C$\semiinjective) complexes of $\S$\semicontramodule s.
 The right derived functor $\bN^\bu\mpsto{}_\S^\boR\bN^\bu$ can
be defined in terms of injective complexes of $\T$\semimodule s and
the left derived functor $\bQ^\bu\mpsto{}^\S_\boL\bQ^\bu$ can be
defined in terms of projective complexes of $\T$\semicontramodule s
(see Remark~\ref{semi-ctrtor-definition}).

 The results of Corollaries~\ref{semi-pull-push-derived}.2--%
\ref{semi-pull-push-derived}.4 do \emph{not} hold for Morita morphisms
of semialgebras, as one can see in the example of the Morita
equivalence related to a Frobenius algebra from
Remark~\ref{morita-change-of-coring} considered as a Morita morphism
in the inverse direction.
 The mentioned results remain valid for left semiprojective and right
semiflat Morita morphisms from $\S$ to $\T$ when the categories of
$\C$\comodule s and $\C$\contramodule s have finite homological
dimensions, or the Morita morphism of semialgebras is induced by
a left coprojective and right coflat Morita morphism of corings, or
more generally when the functors $\bN^\bu\mpsto{}_\S\bN^\bu$, \ 
$\bN^\bu\mpsto\bN^\bu_\S$, and $\bQ^\bu\mpsto{}^\S\bQ^\bu$ map
$\D$\coacyclic{} and $\D$\contraacyclic{} complexes to $\C$\coacyclic{}
and $\C$\contraacyclic{} complexes.

 Morita equivalences of semialgebras do \emph{not} induce equivalences 
of the semiderived categories of semimodules and semicontramodules,
except in rather special cases.
 A right semiflat Morita equivalence between $\S$ and $\T$ does induce
an equivalence of the semiderived categories of left $\S$\semimodule s
and left $\T$\semimodule s when the categories of left $\C$\comodule s
and left $\D$\comodule s have finite homological dimensions, or when
the Morita equivalence comes from a right faithfully coflat Morita
morphism of corings and one of the conditions
of~\ref{morita-change-of-coring}--\ref{relative-homological-dimension}
is satisfied.

\subsubsection{}
 A short summary: one encounters no problems generalizing
the results of
\ref{co-contra-compatible-morphisms}--\ref{faithful-base-ring-change}
and \ref{semi-compatible-morphisms}--\ref{semi-pull-push-derived}
to the case of a Morita morphism of $k$\+algebras and related maps
of corings and semialgebras.
 The problems are manageable when one considers Morita morphisms of
corings.
 And there are severe problems with Morita morphisms/equivalences
of semialgebras, which do not always respect the essential structure
of ``an object split in two halves'' (see Introduction).

\Section{Closed Model Category Structures}

 By a \emph{closed model category} we mean a model category in
the sense of Hovey~\cite{Hov}.

\subsection{Complexes of comodules and contramodules}
\label{co-contra-model-struct}
 Let $\C$ be a coring over a $k$\+algebra $A$.
 Assume that $\C$ is a projective left and a flat right $A$\module{}
and the ring $A$ has a finite left homological dimension.

\begin{thm}
 \textup{(a)} The category of complexes of left\/ $\C$\comodule s has
a closed model category structure with the following properties.
 A morphism is a weak equivalence if and only if its cone is coacyclic.
 A morphism is a cofibration if and only if it is injective and its
cokernel is a complex of $A$\projective\/ $\C$\comodule s.
 A morphism is a fibration if and only if it is surjective and its
kernel is a complex of\/ $\C/A$\injective\/ $\C$\comodule s.
 An object is simultaneously fibrant and cofibrant if and only if
it is a complex of coprojective left\/ $\C$\comodule s. \par
 \textup{(b)} The category of complexes of left\/ $\C$\contramodule s
has a closed model category structure with the following properties.
 A morphism is a weak equivalence if and only if its cone is
contraacyclic.
 A morphism is a cofibration if and only if it is injective and its
cokernel is a complex of\/ $\C/A$\projective\/ $\C$\contramodule s.
 A morphism is a fibration if and only if it is surjective and its
kernel is a complex of $A$\injective\/ $\C$\contramodule s.
 An object is simultaneously fibrant and cofibrant if and only if
it is a complex of coinjective left\/ $\C$\contramodule s.
\end{thm}

\begin{proof}
 Part~(a): the category of complexes of left $\C$\comodule s has
arbitrary limits and colimits, since it is an abelian category with
infinite direct sums and products.
 The two-out-of-three property of weak equivalences follows from
the octahedron axiom, since coacyclic complexes form a triangulated
subcategory of the homotopy category of left $\C$\comodule s.
 The retraction properties are clear, since the classes of projective
$A$\module s, $\C/A$\injective{} $\C$\comodule s, and coacyclic
complexes of $\C$\comodule s are closed under direct summands.
 It is also clear that a morphism is a trivial cofibration if and
only if it is injective and its cokernel is a coacyclic complex of
$A$\projective{} $\C$\comodule s, and a morphism is a trivial
fibration if and only if it is surjective and its kernel is
a coacyclic complex of $\C/A$\injective{} $\C$\comodule s.
 Now let us verify the lifting properties.

\begin{lem1}
 Let $U$, $V$, $X$, and $Y$ be four objects of an abelian category\/
$\sA$, \ $U\rarrow V$ be an injective morphism with the cokernel $E$,
and $X\rarrow Y$ be a surjective morphism with the kernel $K$.
 Suppose that\/ $\Ext^1_\sA(E,K)=0$.
 Then for any two morphisms $U\rarrow X$ and $V\rarrow Y$ forming
a commutative square with the above two morphisms there exists
a morphism $V\rarrow X$ forming two commutative triangles with
the given four morphisms.
\end{lem1}

\begin{proof}
 Let us first find a morphism $V\rarrow X$ making a commutative
triangle with the morphisms $U\rarrow X$ and $U\rarrow V$.
 The obstruction to extending the morphism $U\rarrow X$ from $U$
to $V$ lies in the group $\Ext^1_\sA(E,X)$.
 Since the composition $U\rarrow X\rarrow Y$ admits an extension
to~$V$, our element of $\Ext^1_\sA(E,X)$ becomes zero in
$\Ext^1_\sA(E,Y)$ and therefore comes from the group
$\Ext^1_\sA(E,K)$.
 Now let us modify the obtained morphism so that the new morphism
$V\rarrow X$ forms also a commutative triangle with the morphisms
$V\rarrow Y$ and $X\rarrow Y$.
 The difference between the given morphism $V\rarrow Y$ and
the composition $V\rarrow X\rarrow Y$ is a morphism $V\rarrow Y$
annihilating $U$, so it comes from a morphism $E\rarrow Y$.
 We need to lift the latter to a morphism $E\rarrow X$.
 The obstruction to this lies in $\Ext^1_\sA(E,K)$.
\end{proof}

 To verify the condition of Lemma~1, consider an extension $\E^\bu
\rarrow\M^\bu\rarrow\K^\bu$ of a complex of $A$\projective{} left
$\C$\comodule s $\K^\bu$ by a complex of $\C/A$\injective{} left
$\C$\comodule s $\E^\bu$.
 By Lemma~\ref{rel-inj-proj-co-contra-mod}.1(a), this extension
is term-wise split, so it comes from a morphism of complexes of
$\C$\comodule s $\K^\bu\rarrow\E^\bu[1]$.
 Now suppose that one of the complexes $\K^\bu$ and $\E^\bu$
is coacyclic.
 Then any morphism $\K^\bu\rarrow\E^\bu[1]$ is homotopic to zero by
a result of~\ref{co-contra-ctrtor-definition}, hence the extension
of complexes $\E^\bu\rarrow\M^\bu\rarrow\K^\bu$ is split.
 The lifing properties are proven.

 It remains to construct the functorial factorizations.
 These constructions use two building blocks: one is
Lemma~\ref{proj-inj-co-contra-module}(a), the other one
is the following Lemma~2.

\begin{lem2}
 \textup{(a)} There exists a (not always additive) functor assigning
to any left\/ $\C$\comodule{} an injective morphism from it into
a\/ $\C/A$\injective{} left\/ $\C$\comodule{} with an $A$\projective{}
cokernel. \par
 \textup{(b)} There exists a (not always additive) functor assigning
to any left\/ $\C$\contramodule{} a surjective morphism onto it from
a\/ $\C/A$\projective{} left\/ $\C$\contramodule{} with an
$A$\injective{} kernel.
\end{lem2}

\begin{proof}
 Part~(a): let $\M$ be a left $\C$\comodule{} and $\cP(\M)\rarrow\M$ be
the surjective morphism onto it from an $A$\projective{} $\C$\comodule{}
$\cP(\M)$ constructed in Lemma~\ref{proj-inj-co-contra-module}(a).
 Let $\K$ be kernel of the map $\cP(\M)\rarrow\M$ and let
$\cP(\M)\rarrow\C\ot_A\cP(\M)$ be the $\C$\+coaction map.
 Set $\cJ(\M)$ to be the cokernel of the composition $\K\rarrow
\cP(\M)\rarrow\C\ot_A\cP(\M)$.
 Then the composition of maps $\cP(\M)\rarrow\C\ot_A\cP(\M)\rarrow
\cJ(\M)$ factorizes through the surjection $\cP(\M)\rarrow\M$, so there
is a natural injective morphism of $\C$\comodule s $\M\rarrow\cJ(\M)$.
 The $\C$\comodule{} $\cJ(\M)$ is $\C/A$\injective{} as the cokernel of
an injective map of $\C/A$\injective{} $\C$\comodule s $\K\rarrow
\C\ot_A\cP(\M)$.
 The cokernel of the map $\M\rarrow\cJ(\M)$ is isomorphic to
the cokernel of the map $\cP(\M)\rarrow\C\ot_A\cP(\M)$ and hence
$A$\projective.
 Part~(a) is proven; the construction of the surjective morphism of
$\C$\contramodule s $\gF(\P)\rarrow\P$ in part~(b) is completely
analogous.
\end{proof}

 Let us first decompose functorially an arbitrary morphism of complexes
of left $\C$\comodule s $\L^\bu\rarrow\M^\bu$ into a cofibration
followed by a fibration.
 This can be done in either of two dual ways.
 Let us start with a surjective morphism $\cP^+(\M^\bu)\rarrow\M^\bu$
onto the complex $\M^\bu$ from a complex of $A$\projective{} left
$\C$\comodule s $\cP^+(\M^\bu)$ constructed as in the proof of
Theorem~\ref{cotor-main-theorem}.
 Let $\K^\bu$ be the kernel of the morphism $\L^\bu\oplus\cP^+(\M^\bu)
\rarrow\M^\bu$ and let $\K^\bu\rarrow\cJ^+(\K^\bu)$ be an injective
morphism from the complex $\K^\bu$ into a complex of $\C/A$\injective{}
left $\C$\comodule s $\cJ^+(\K^\bu)$ constructed in the analogous way
using Lemma~2.
 The cokernel of this morphism is a complex of $A$\projective{}
$\C$\comodule s.
 Let $\E^\bu$ denote the fibered coproduct of $\L^\bu\oplus
\cP^+(\M^\bu)$ and $\cJ^+(\K^\bu)$ over $\K^\bu$.
 There is a natural morphism of complexes $\E^\bu\rarrow\M^\bu$ whose
composition with the morphism $\cJ^+(\K^\bu)\rarrow\E^\bu$ is zero and
composition with the morphism $\L^\bu\oplus\cP^+(\M^\bu)\rarrow \E^\bu$
is equal to our morphism $\L^\bu\oplus\cP^+(\M^\bu)\rarrow\M^\bu$.
 The morphism $\L^\bu\rarrow\M^\bu$ is equal to the composition
$\L^\bu\rarrow\E^\bu\rarrow\M^\bu$.
 The cokernel of the morphism $\L^\bu\rarrow\E^\bu$ is an extension of
the cokernel of the morphism $\K^\bu\rarrow\cJ^+(\K^\bu)$ and
the complex $\cP^+(\M^\bu)$, hence a complex of $A$\projective{}
$\C$\comodule s.
 The kernel of the morphism $\E^\bu\rarrow\M^\bu$ is isomorphic to
$\cJ^+(\K^\bu)$, which is a complex of $\C/A$\injective{}
$\C$\comodule s.
 Another way is to start with an injective morphism $\L^\bu\rarrow
\cJ^+(\L^\bu)$ and consider the cokernel of the morphism $\L^\bu
\rarrow\M^\bu\oplus\cJ^+(\L^\bu)$.

 Now let us construct a factorization of the morphism $\L^\bu\rarrow
\M^\bu$ into a cofibration followed by a trivial fibration.
 Represent the kernel of the morphism $\E^\bu\rarrow\M^\bu$ as
the quotient complex of a complex of $A$\projective{} left
$\C$\comodule s $\E_1^\bu$ by a complex of $\C/A$\injective{}
$\C$\comodule s; represent the latter complex as the quotient complex
of a complex $\E_2^\bu$ with the analogous properties, etc.
 The complexes $\E_i^\bu$ are also complexes of $\C/A$\injective{}
$\C$\comodule s as extensions of complexes of $\C/A$\injective{}
$\C$\comodule s.
 For $d$ large enough, the kernel $\cZ^\bu$ of the morphism
$\E_d^\bu\rarrow\E_{d-1}^\bu$ will be a complex of $A$\projective{}
$\C$\comodule s.
 Actually, $\E_i^\bu$ and $\cZ^\bu$ are complexes of coprojective
$\C$\comodule s, as a $\C/A$\injective{} $A$\projective{}  left
$\C$\comodule{} $\cQ$ is coprojective (since the injection of
$\C$\comodule s $\cQ\rarrow\C\ot_A\cQ$ splits, $\cQ\rarrow\C\ot_A\cQ
\rarrow\C\ot_A\cQ/\cQ$ being an exact triple of $A$\projective{}
$\C$\comodule s).
 Let $\K^\bu$ be the total complex of the bicomplex $\cZ^\bu\rarrow
\E_d^\bu\rarrow\dsb\rarrow\E_1^\bu\rarrow\E^\bu$.
 Then the morphism $\L^\bu\rarrow\M^\bu$ factorizes through $\K^\bu$
in a natural way, the kernel of the morphism $\K^\bu\rarrow\M^\bu$
is a coacyclic complex of $\C/A$\injective{} $\C$\comodule s, and
the cokernel of the morphism $\L^\bu\rarrow\K^\bu$ is a complex of
$A$\projective{} $\C$\comodule s.
 Notice that the complex $\K^\bu$ is the cone of the natural
morphism $\boL_1(\ker(\E^\bu\to\M^\bu))\rarrow\E^\bu$, where
$\boL_1$ denotes the functor from the proof of
Theorem~\ref{coext-main-theorem}.

 Finally, let us construct a factorization of the morphism $\L^\bu
\rarrow\M^\bu$ into a trivial cofibration followed by a fibration.
 Represent the cokernel of the morphism $\L^\bu\rarrow\E^\bu$
as a subcomplex of a complex of $\C/A$\injective{} left
$\C$\comodule s ${}^1\E^\bu$ such that the quotient complex is
a complex of $A$\projective{} $\C$\comodule s; represent this quotient
complex as a subcomplex of a complex ${}^2\E^\bu$ with the analogous
properties, etc.
 The complexes ${}^i\E^\bu$ are also complexes of $A$\projective{}
$\C$\comodule s as extensions of complexes of $A$\projective{}
$\C$\comodule s (so they are complexes of coprojective $\C$\comodule s).
 Let $\K^\bu$ be the total complex of the bicomplex $\E^\bu\rarrow
{}^1\E^\bu\rarrow{}^2\E^\bu\rarrow\dsb$, constructed by taking
infinite direct sums along the diagonals.
  Then the morphism $\L^\bu\rarrow\M^\bu$ factorizes through $\K^\bu$
in a natural way, the cokernel of the morphism $\L^\bu\rarrow\K^\bu$
is a coacyclic complex of $A$\projective{} $\C$\comodule s, and
the kernel of the morphism $\K^\bu\rarrow\M^\bu$ is a complex of
$\C/A$\injective{} $\C$\comodule s.
 The class of $\C/A$\injective{} $\C$\comodule s is closed under
infinite direct sums by Lemma~\ref{rel-inj-proj-co-contra-mod}.2(a).

 Part~(a) is proven; the proof of part~(b) is completely analogous.
\end{proof}

\begin{rmk}
 It follows from the proof of Lemma~2 that any $\C/A$\injective{}
left $\C$\comodule{} can be obtained from coinduced $\C$\comodule s
by taking extensions, cokernels of injective morphisms, and direct
summands.
 Analogously, any $\C/A$\projective{} left $\C$\contramodule{} can be
obtained from induced $\C$\contramodule s by taking extensions,
kernels of surjective morphisms, and direct summands.
\end{rmk}

 Let $\C$ and $\D$ be two corings satisfying the above assumptions and
$\C\rarrow\D$ be a map of corings compatible with a $k$\+algebra map
$A\rarrow B$.
 Then the pair of adjoint functors $\M^\bu\mpsto{}_B\M^\bu$ and 
$\N^\bu\mpsto{}_\C\N^\bu$ is a Quillen adjunction~\cite{Hov} from
the category of complexes of left $\C$\comodule s to the category of
complexes of left $\D$\comodule s; the pair of adjoint functors
$\Q^\bu\mpsto{}^\C\Q^\bu$ and $\P^\bu\mpsto{}^B\P^\bu$ is a Quillen
adjunction from the category of complexes of left $\D$\contramodule s
to the category of complexes of left $\C$\contramodule s.
 The same applies to the case of a Morita morphism $(E,E\dual)$ from
$A$ to $B$ and a morphism ${}_B\C_B\rarrow\D$ of corings over~$B$.

 The pair of adjoint functors $\Phi_\C$ and $\Psi_\C$ applied to
complexes term-wise is \emph{not} a Quillen equivalence, and not even
a Quillen adjunction, between the model category of complexes of
left $\C$\contramodule s and the model category of complexes of
left $\C$\comodule s.
 This pair of functors \emph{is} a Quillen equivalence, however, when
$\C$ is a coring over a semisimple ring~$A$.
 In general, the model categories of complexes of left $\C$\comodule s
and left $\C$\contramodule s can be connected by a chain of three
Quillen equivalences (see Remark~\ref{semi-model-struct}.2).

\subsection{Complexes of semimodules and semicontramodules}
\label{semi-model-struct}
 Let $\S$ be a semialgebra over a coring $\C$ over a $k$\+algebra $A$.
 Assume that $\C$ is a projective left and a flat right $A$\module,
$\S$ is a coprojective left and a coflat right $\C$\comodule, and
the ring $A$ has a finite left homological dimension.

 A left $\S$\semimodule{} $\bL$ is called \emph{$\S/\C/A$\projective}
if the functor of $\S$\semimodule{} homomorphisms from $\bL$ maps
exact triples of $\C/A$\injective{} left $\S$\semimodule s to exact
triples.
 An $A$\projective{} left $\S$\semimodule{} $\bL$ is called
\emph{$\S/\C/A$\semiprojective} if the functor of semihomomorphisms
from $\bL$ over~$\S$ maps exact triples of $\C/A$\coinjective{}
left $\S$\semicontramodule s to exact triples.
 Analogously, a left $\S$\semicontramodule{} $\bQ$ is called
\emph{$\S/\C/A$\injective} if the functor of $\S$\semicontramodule{}
homomorphisms into $\bQ$ maps exact triples of $\C/A$\projective{}
left $\S$\semicontramodule s to exact triples.
 An $A$\injective{} left $\S$\semicontramodule{} $\bQ$ is called
\emph{$\S/\C/A$\semiinjective} if the functor of semihomomorphisms
into $\bQ$ over~$\S$ maps exact triples of $\C/A$\coprojective{}
left $\S$\semimodule s to exact triples.

 As in Lemma~\ref{birelatively-adjusted}, it follows from
Proposition~\ref{semitensor-contratensor-assoc}.2(c) that
an $A$\projective{} left $\S$\semimodule{} is $\S/\C/A$\projective{}
if and only if it is $\S/\C/A$\semiprojective.
 Analogously, it follows from
Proposition~\ref{semitensor-contratensor-assoc}.3(c) that
an $A$\injective{} left $\S$\semicontramodule{} is
$\S/\C/A$\injective{} if and only if it is $\S/\C/A$\semiinjective.
 It will be shown below that any $\S/\C/A$\projective{} left
$\S$\semimodule{} is $A$\projective{} and any $\S/\C/A$\injective{}
left $\S$\semicontramodule{} is $A$\injective.

\begin{thm}
 \textup{(a)} The category of complexes of left\/ $\S$\semimodule s has
a closed model category structure with the following properties.
 A morphism is a weak equivalence if and only if its cone is\/
$\C$\coacyclic.
 A morphism is a cofibration if and only if it is injective and its
cokernel is an\/ $\S/\C/A$\projective{} complex of\/
$\S/\C/A$\projective\/ $\S$\semimodule s.
 A morphism is a fibration if and only if it is surjective and its
kernel is a complex of\/ $\C/A$\injective\/ $\S$\semimodule s.
 An object is simultaneously fibrant and cofibrant if and only if it is
a semiprojective complex of semiprojective left\/ $\S$\semimodule s.
\par
 \textup{(b)} The category of complexes of left\/
$\S$\semicontramodule s has a closed model category structure with
the following properties.
 A morphism is a weak equivalence if and only if its cone is\/
$\C$\contraacyclic.
 A morphism is a cofibration if and only if it is injective and its
cokernel is a complex of\/ $\C/A$\projective\/ $\S$\semicontramodule s.
 A morphism is a fibration if and only if it is injective and its
cokernel is an\/ $\S/\C/A$\injective{} complex of\/ $\S/\C/A$\injective\/
$\S$\semicontramodule s.
 An object is simultaneously fibrant and cofibrant if and only if it is
a semiinjective complex of semiinjective left\/ $\S$\semicontramodule s.
\end{thm}

\begin{proof}
 Part~(a): existence of limits and colimits, the two-out-of-three
property of weak equivalences, and the retraction properties are
verified as in the proof of Theorem~\ref{co-contra-model-struct}.
 It is clear that a morphism is a trivial cofibration if and only if
it is injective and its cokernel is a $\C$\coacyclic{}
$\S/\C/A$\projective{} complex of $\S/\C/A$\projective{}
$\S$\semimodule s, and a morphism is a trivial fibration if and only if
it is surjective and its kernel is a $\C$\coacyclic{} complex of
$\C/A$\injective{} $\S$\semimodule s.
 To check the lifting properties, use
Lemma~\ref{co-contra-model-struct}.1.
 Consider an extension $\bE^\bu\rarrow\bM^\bu\rarrow\bK^\bu$ of
an $\S/\C/A$\projective{} complex of $\S/\C/A$\projective{} left
$\S$\semimodule s $\bK^\bu$ by a complex of $\C/A$\injective{} left
$\S$\semimodule s $\bE^\bu$.
 By the next Lemma~1, this extension is term-wise split, so it comes
from a morphism of complexes of $\S$\semimodule s
$\bK^\bu\rarrow\bE^\bu[1]$.
 Now suppose that one of the complexes $\bK^\bu$ and $\bE^\bu$ is
$\C$\coacyclic.
 Then any morphism $\bK^\bu\rarrow\bE^\bu[1]$ is homotopic to zero
by a result of~\ref{semi-ctrtor-definition}, hence the extension of
complexes $\bE^\bu\rarrow\bM^\bu\rarrow\bK^\bu$ is split.
 So after Lemma~1 is proven it will remain to construct the functorial
factorizations.

\begin{lem1}
 \textup{(a)} A left\/ $\S$\semimodule\/ $\bL$ is\/
$\S/\C/A$\projective{} if and only if for any\/ $\C/A$\injective{}
left\/ $\S$\semimodule\/ $\bM$ the $k$\module s\/ $\Ext_\S^i(\bL,\bM)$
of Yoneda extensions in the abelian category of left\/
$\S$\semimodule s vanish for all $i>0$.
 The functor of\/ $\S$\semimodule{} homomorphisms into
a\/ $\C/A$\injective\/ $\S$\semimodule{} maps exact triples of\/
$\S/\C/A$\projective{} left\/ $\S$\semimodule s to exact triples.
 The classes of\/ $\S/\C/A$\projective{} left\/ $\S$\semimodule s and\/
$\S/\C/A$\projective{} complexes of\/ $\S/\C/A$\projective{} left\/
$\S$\semimodule s are closed under extensions and kernels of
surjective morphisms. \par
 \textup{(b)} A left\/ $\S$\semicontramodule\/ $\bQ$ is\/
$\S/\C/A$\injective{} if and only if for any\/ $\C/A$\projective{}
left\/ $\S$\semicontramodule\/ $\bP$ the $k$\module s\/
$\Ext^{\S,i}(\bP,\bQ)$ of Yoneda extensions in the abelian category
of left\/ $\S$\semicontramodule s vanish for all $i>0$.
 The functor of\/ $\S$\semicontramodule{} homomorphisms from
a\/ $\C/A$\projective\/ $\S$\semicontramodule{} maps exact triples of\/
$\S/\C/A$\injective{} left\/ $\S$\semicontramodule s to exact triples.
 The classes of\/ $\S/\C/A$\injective{} left $\S$\semicontramodule s
and\/ $\S/\C/A$\injective{} complexes of\/ $\S/\C/A$\injective{}
left\/ $\S$\semicontramodule s are closed under extensions and
cokernels of injective morphisms.
\end{lem1}

\begin{proof}
 Part~(a): the forgetful functor $\S\simodl\rarrow\C\comodl$ preserves
injective objects, since it is right adjoint to the exact functor of
induction.
 Let us show that any $\C/A$\injective{} left $\S$\semimodule{} $\bM$
is a subsemimodule of an injective $\S$\semimodule{} (it will follow
that the category of left $\S$\semimodule s has enough injectives).
 The construction of Lemma~\ref{proj-inj-semi-mod-contra}(b) assigns
to a $\C/A$\coinjective{} left $\S$\semicontramodule{} $\bP$
an injective map from it into a semiinjective $\S$\semicontramodule{}
$\bgI(\bP)$ with a $\C/A$\coinjective{} cokernel.
 Indeed, the cokernel of the map $\bP\rarrow\gI(\bP)$ is
a $\C/A$\coinjective{} $\C$\contramodule{} by
Lemma~\ref{proj-inj-co-contra-module}(b), so $\gI(\bP)$ is
a $\C/A$\coinjective{} $\C$\contramodule{} as an extension of two
$\C/A$\coinjective{} $\C$\contramodule s and a coinjective
$\C$\contramodule{} as an $A$\injective{} and $\C/A$\coinjective{}
$\C$\contramodule.
 Hence $\bgI(\bP)=\Cohom_\C(\S,\gI(\bP))$ is a semiinjective
$\S$\semicontramodule.
 The cokernel of the composition of injective morphisms $\bP\rarrow
\Cohom_\C(\S,\bP)\rarrow\Cohom_\C(\S,\gI(\bP))$ is an extension of
the cokernel of the morphism $\Cohom_\C(\S,\bP)\rarrow
\Cohom_\C(\S,\gI(\bP))$ and the cokernel of the morphism
$\bP\rarrow\Cohom_\C(\S,\bP)$; the former is $\C/A$\coinjective{}
since the cokernel of the morphism $\bP\rarrow\bgI(\bP)$ is, and
the latter is $\C/A$\coinjective{} as a $\C$\contramodule{} direct
summand of $\Cohom_\C(\S,\bP)$.
 Hence the cokernel of the morphism $\bP\rarrow\bgI(\bP)$ is
$\C/A$\coinjective.
 Applying these observations to the $\S$\semicontramodule{} $\bP=
\Phi_\S(\bM)$ and using Lemmas~\ref{rel-inj-proj-co-contra-mod}.2(b)
and~\ref{rel-inj-proj-co-contra-mod}.1(c), we conclude that
$\bM\rarrow\Phi_\S\bgI(\Psi_\S\bM)$ is an injective morphism
of $\S$\semimodule s whenever $\bM$ is a $\C/A$\injective{} left
$\S$\semimodule.
 Now the functor $\Phi_\S$ maps semiinjective $\S$\semicontramodule s
to injective $\S$\semimodule s by
Proposition~\ref{semitensor-contratensor-assoc}.2(a).

 So any $\C/A$\injective{} left $\S$\semimodule{} $\bM$ has
an injective right resolution; by the construction or by
Lemma~\ref{rel-inj-proj-co-contra-mod}.1(a), this resolution is
exact with respect to the exact category of $\C/A$\injective{}
$\S$\semimodule s.
 Applying to this resolution the functor of $\S$\semimodule{}
homomorphisms from an $\S/\C/A$\projective{} left $\S$\semimodule{}
$\bL$, we obtain the desired vanishing $\Ext_\S^i(\bL,\bM)=0$ for
all $i>0$.
 The remaining assertions follow (to verify the assertions related to
complexes, notice that the class of acyclic complexes of $k$\module s
is closed under extensions and cokernels of injective morphisms).
 Part~(a) is proven; the proof of part~(b) is completely analogous
and based on the construction of a semicontramodule projective
resolution.
 Alternatively, one can argue as in the proof of 
Lemma~\ref{rel-inj-proj-co-contra-mod}.1(a-b).
\end{proof}

 The analogous results for $\S/\C/A$\semiprojective{} (complexes of)
left $\S$\semimodule s and $\S/\C/A$\semiinjective{} (complexes of)
left $\S$\semicontramodule s can be obtained by considering
the derived functor $\SemiExt_\S^*(\bM,\bP)$, defined as
the cohomology of the object $\SemiExt_\S(\bM,\bP)$ of $\sD(k\modl)$.
 For an $A$\projective{} $\S$\semimodule{} $\bM$ and
a $\C/A$\coinjective{} $\S$\semicontramodule{} $\bP$ or
a $\C/A$\coprojective{} $\S$\semimodule{} $\bM$ and an $A$\injective{}
$\S$\semicontramodule{} $\bP$ it is computed by the cobar-complex
$\Cohom_\C(\bM,\bP)\rarrow\Cohom_\C(\S\oc_\C\bM\;\bP)\rarrow\dsb$,
hence $\SemiExt_\S^i(\bM,\bP)=0$ for $i>0$ and
$\SemiExt_\S^0(\bM,\bP)=\SemiHom_\S(\bM,\bP)$.

\begin{lem2}
 \textup{(a)} There exists a (not always additive) functor assigning
to any left\/ $\S$\semimodule{} an injective morphism from it into
a\/ $\C/A$\injective\/ $\S$\semimodule{} with
an\/ $\S/\C/A$\projective{} cokernel.
 Furthermore, there exists a functor assigning to any complex of
left\/ $\S$\semimodule s{} an injective morphism from it into a complex
of\/ $\C/A$\injective\/ $\S$\semimodule s such that the cokernel is
a\/ $\C$\coacyclic\/ $\S/\C/A$\projective{} complex of\/
$\S/\C/A$\projective\/ $\S$\semimodule s. \par
 \textup{(b)} There exists a (not always additive) functor assigning
to any left\/ $\S$\semicontramodule{} a surjective morphism onto it
from a\/ $\C/A$\projective\/ $\S$\semicontramodule{} with
an\/ $\S/\C/A$\injective{} kernel.
 Furthermore, there exists a functor assigning to any complex of
left\/ $\S$\semicontramodule s{} a surjective morphism onto it
from a complex of\/ $\C/A$\projective\/ $\S$\semicontramodule s
such that the kernel is a\/ $\C$\contraacyclic{} $\S/\C/A$\injective{}
complex of\/ $\S/\C/A$\injective\/ $\S$\semicontramodule s.
\end{lem2}

\begin{proof}
 Part~(a): modify the construction of the second assertion of
Lemma~\ref{coflat-semimodule-injection}, replacing the injective
morphism of $\C$\comodule s $\M\rarrow\cG(\M)=\C\ot_A\M$ with
the injective morphism of $\C$\comodule s $\M\rarrow\cJ(\M)$
constructed in Lemma~\ref{co-contra-model-struct}.2(a).
 In other words, for any left $\S$\semimodule{} $\bM$ denote by
$\bK(\bM)$ the kernel of the morphism $\S\oc_\C\bM\rarrow\bM$ and by
$\bcQ(\bM)$ the cokernel of the composition $\bK(\bM)\rarrow
\S\oc_\C\bM\rarrow\S\oc_\C\cJ(\bM)$.
 The composition of maps $\S\oc_\C\bM\rarrow\S\oc_\C\cJ(\bM)\rarrow
\bcQ(\bM)$ factorizes through the surjection $\S\oc_\C\bM\rarrow\bM$,
so there is a natural injective morphism of $\S$\semimodule s
$\bM\rarrow\bcQ(\bM)$.
 The cokernel of this morphism is isomorphic to
$\S\oc_\C(\cJ(\bM)/\bM)$, which is an $\S/\C/A$\projective{}
$\S$\semimodule{} because $\cJ(\bM)/\bM$ is an $A$\projective{}
$\C$\comodule.
 As in the proof of Lemma~\ref{coflat-semimodule-injection},
the $\S$\semimodule{} morphism $\bM\rarrow\bcQ(\bM)$ can be lifted
to a $\C$\comodule{} morphism $\bM\rarrow\S\oc_\C\cJ(\bM)$.
 Let $\bcJ(\bM)$ denote the inductive limit of the sequence
$\bM\rarrow\S\oc_\C\cJ(\bM)\rarrow\bcQ(\bM)\rarrow
\S\oc_\C\cJ(\bcQ(\bM))\rarrow\bcQ(\bcQ(\bM))\rarrow\dsb$; it is
the desired $\C/A$\injective{} $\S$\semimodule{} into which $\bM$
maps injectively with an $\S/\C/A$\projective{} cokernel.
 Indeed, $\bcJ(\bM)$ is $\C/A$\injective{} by
Sublemma~\ref{coproj-coinj-semi-mod-contra}.B and the cokernel
of the morphism $\bM\rarrow\bcJ(\bM)$ is $\S/\C/A$\projective{}
by the next Sublemma.

\begin{subl}
 \textup{(a)} Let\/ $0=\bcU_0^\bu\rarrow\bcU_1^\bu\rarrow\bcU_2^\bu
\rarrow\dsb$ be an inductive system of complexes of left\/
$\S$\semimodule s such that the successive cokernels\/
$\coker(\bgU_{i-1}^\bu\to\bgU_i^\bu)$ are\/ $\S/\C/A$\projective{}
complexes of\/ $\S/\C/A$\projective\/ $\S$\semimodule s.
 Then the inductive limit\/ $\ilim\bcU_i^\bu$ is
an\/ $\S/\C/A$\projective{} complex of\/ $\S/\C/A$\projective\/
$\S$\semimodule s. \par
 \textup{(b)} Let\/ $0=\bgU_0^\bu\larrow\bgU_1^\bu\larrow\bgU_2^\bu
\larrow\dsb$ be a projective system of complexes of left\/
$\S$\semicontramodule s such that the successive kernels\/
$\ker(\bgU_i^\bu\to\bgU_{i-1}^\bu)$ are\/ $\S/\C/A$\injective{}
complexes of\/ $\S/\C/A$\injective\/ $\S$\semimodule s.
 Then the projective limit\/ $\plim\bgU_i^\bu$ is
an\/ $\S/\C/A$\injective{} complex of\/ $\S/\C/A$\injective\/
$\S$\semimodule s.
\end{subl}

\begin{proof}
 The forgetful functor $\S\simodl\rarrow A\modl$ preserves
inductive limits, since it preserves cokernels and infinite direct
sums, so one has $\Hom_\S(\ilim\bcU_i^\bu\;\bM^\bu)=
\plim\Hom_\S(\bcU_i^\bu,\bM^\bu)$ for any complex of left
$\S$\semimodule s $\bM^\bu$.
 Analogously, the forgetful functor $\S\sicntr\rarrow A\modl$
preserves projective limits, since it preserves kernels and
infinite products, so one has $\Hom^\S(\bP^\bu\;\plim\bgU_i^\bu)
=\plim\Hom^\S(\bP^\bu,\bgU_i^\bu)$ for any complex of left
$\S$\semicontramodule s $\bP^\bu$.
 As the projective limits of sequences of surjective maps
preserve exact triples and acyclic complexes, the assertions
of Sublemma follow from Lemma~1.
\end{proof}

 The first statement of Lemma~2(a) is proven.
 To prove the second one, consider the functor assigning to a complex
of left $\S$\semimodule s $\bM^\bu$ the injective map from it into
the complex $\bcJ^+(\bM^\bu)$, which is constructed in terms of
the functor $\bM\mpsto\bcJ(\bM)$ as in the proof of
Theorem~\ref{cotor-main-theorem}.
 By Sublemma, the cokernel of the morphism $\bM^\bu\rarrow
\bcJ^+(\bM^\bu)$ is an $\S/\C/A$\projective{} complex of
$\S/\C/A$\projective{} $\S$\semimodule s, since a complex of 
$\S$\semimodule s induced from a complex of $A$\projective{}
$\C$\comodule s belongs to this class.
 Set ${}^0\bcJ^\bu=\bcJ^+(\bM^\bu)$, \ ${}^1\bcJ^\bu=\bcJ^+
(\coker(\bM^\bu\to{}^0\bcJ^\bu))$, etc.
 The complexes ${}^i\bcJ^\bu$ are complexes of $\C/A$\injective{}
$\S$\semimodule s and the complexes $\coker(\bM^\bu\to{}^0\bcJ^\bu)$,
$\coker({}^{i-1}\bcJ^\bu\to{}^i\bcJ^\bu)$ are $\S/\C/A$\projective{}
complexes of $\S/\C/A$\projective{} $\S$\semimodule s.
 The complexes ${}^i\bcJ^\bu$ for $i>0$ are also $\S/\C/A$\projective{}
complexes of $\S/\C/A$\projective{} $\S$\semimodule s as extensions
of complexes with these properties.
 Let $\bK^\bu$ be the total complex of the bicomplex ${}^0\bcJ^\bu
\rarrow{}^1\bcJ^\bu\rarrow\dsb$, constructed by taking infinite
direct sums along the diagonals.
 Then $\bK^\bu$ is a complex of $\C/A$\injective{} $\S$\semimodule s 
and the cokernel of the injective morphism $\bM^\bu\rarrow\bK^\bu$
is a $\C$\coacyclic{} (and even $\S$\coacyclic) $\S/\C/A$\projective{}
complex of $\S/\C/A$\projective{} $\S$\semimodule s.
 To check the latter properties, one can apply Sublemma to
the canonical filtration of the complex ${}^0\bcJ^\bu/\bM^\bu\rarrow
{}^1\bcJ^\bu\rarrow{}^2\bcJ^\bu\rarrow\dsb$

 The proof of Lemma~2(b) is completely analogous and based on
the modification of the construction of the second assertion of
Lemma~\ref{coproj-coinj-semi-mod-contra}(a) using the surjective
morphism of $\C$\contramodule s $\gF(\P)\rarrow\P$ from
Lemma~\ref{co-contra-model-struct}.2(b) in place of
the morphism $\gG(\P)=\Hom_A(\C,\P)\rarrow\P$.
\end{proof}

 In the sequel we will denote by $\bM\mpsto\bcJ(\bM)$ the functor
constructed in Lemma~2 rather than its more simplistic version
from Lemmas~\ref{coflat-semimodule-injection}
and~\ref{coproj-coinj-semi-mod-contra}.

\begin{lem3}
 \textup{(a)} There exists a (not always additive) functor assigning
to any left\/ $\S$\semimodule{} a surjective map onto it from
an\/ $\S/\C/A$\projective\/ $\S$\semimodule{} with
a\/ $\C/A$\injective{} kernel.
 Furthermore, there exists a functor assigning to any complex of
left\/ $\S$\semimodule s a surjective map onto it from
an\/ $\S/\C/A$\projective{} complex of\/ $\S/\C/A$\projective\/
$\S$\semimodule s such that the kernel is a\/ $\C$\coacyclic{}
complex of\/ $\C/A$\injective\/ $\S$\semimodule s. \par
 \textup{(b)} There exists a (not always additive) functor assigning
to any left\/ $\S$\semicontramodule{} an injective map from it into
an\/ $\S/\C/A$\injective\/ $\S$\semicontramodule{} with
a\/ $\C/A$\projective{} cokernel.
 Furthermore, there exists a functor assigning to any complex of
left\/ $\S$\semicontramodule s an injective map from it into
an\/ $\S/\C/A$\injective{} complex of\/ $\S/\C/A$\injective\/
$\S$\semicontramodule s such that the cokernel is a\/
$\C$\contraacyclic{} complex of\/ $\C/A$\projective\/
$\S$\semicontramodule s.
\end{lem3}

\begin{proof}
 Part~(a): for any left $\S$\semimodule{} $\bL$, consider the injective
morphism $\bL\rarrow\bcJ(\bL)$ from Lemma~2 and denote by $\bK(\bL)$
its cokernel.
 The functor $\bM\mpsto\bcP(\bM)$ of
Lemmas~\ref{flat-semimodule-surjection}
and~\ref{proj-inj-semi-mod-contra} assigns to a $\C/A$\injective{}
left $\S$\semimodule{} $\bM$ a surjective morphism onto it from
the $\S$\semimodule{} $\bcP(\bM)$ induced from a coprojective
$\C$\comodule{} $\cP(\bM)$ such that the kernel is
a $\C/A$\injective{} $\S$\semimodule{} (see the proof of Lemma~1).
 Denote by $\bcF(\bL)$ the kernel of the composition $\bcP(\bcJ(\bL))
\rarrow\bcJ(\bL)\rarrow\bK(\bL)$.
 The composition of maps $\bcF(\bL)\rarrow\bcP(\bcJ(\bL))\rarrow
\bcJ(\bL)$ factorizes through the injection $\bL\rarrow\bcJ(\bL)$,
so there is a natural surjective morphism of $\S$\semimodule s
$\bcF(\bL)\rarrow\bL$.
 The $\S$\semimodule s $\bcP(\bcJ(\bL))$ and $\bK(\bL)$ are
$\S/\C/A$\projective, hence the $\S$\semimodule{} $\bcF(\bL)$ is also
$\S/\C/A$\projective.
 The kernel of the morphism $\bcF(\bL)\rarrow\bL$ is $\C/A$\injective,
since it is isomorphic to the kernel of the morphism $\bcP(\bcJ(\bL))
\rarrow\bcJ(\bL)$.
 This proves the first statement of part~(a).

 Now consider the functor assigning to any complex of left
$\S$\semimodule s $\bL^\bu$ the surjective map onto it from
the complex $\bcF^+(\bL^\bu)$.
 The complex $\bcF^+(\bL^\bu)$ is $\S/\C/A$\projective{} as the kernel
of a surjective morphism of $\S/\C/A$\projective{} complexes;
it is also a complex of $\S/\C/A$\projective{} and $A$\projective{}
$\S$\semimodule s.
 Set $\bcF_0^\bu=\bcF^+(\bL^\bu)$, \ $\bcF_1^\bu=\bcF(\ker(\bcF_0^\bu
\to\bL^\bu))$, etc.
 The complexes $\ker(\bcF_0^\bu\to\bL^\bu)$, $\ker(\bcF_i^\bu\to
\bcF_{i-1}^\bu)$ are complexes of $\C/A$\injective{} $\S$\semimodule s,
hence the complexes $\bcF_i^\bu$ for $i>0$ are also complexes of
$\C/A$\injective{} $\S$\semimodule s as extensions of complexes of
$\C/A$\injective{} $\S$\semimodule s.
 For $d$~large enough, the kernel $\bcZ^\bu$ of the morphism
$\bcF_{d-1}^\bu\rarrow\bcF_{d-2}^\bu$ will be a complex of
$A$\projective{} $\S$\semimodule s, and consequently a complex of
$\C$\coprojective{} $\S$\semimodule s.
 Let $\bE^\bu$ be the total complex of the bicomplex
$$
 \dsb\lrarrow\S\oc_\C\S\oc_\C\bcZ^\bu\lrarrow\S\oc_\C\bcZ^\bu\lrarrow
 \bcF_{d-1}^\bu\lrarrow\bcF_{d-2}^\bu\lrarrow\dsb\lrarrow\bcF_0^\bu,
$$
constructed by taking infinite direct sums along the diagonals.
 Then the complex $\bE^\bu$ is a complex of $\S/\C/A$\projective{}
$\S$\semimodule s, and it is an $\S/\C/A$\projective{} complex since
it is homotopy equivalent to a complex obtained from such complexes
using the operations of cone and infinite direct sum.
 The kernel of the morphism $\bE^\bu\rarrow\bL^\bu$ is a complex of
$\C/A$\injective{} $\S$\semimodule s, and it is $\C$\coacyclic{}  
since it contains a $\C$\contractible{} subcomplex of $\S$\semimodule s
such that the quotient complex is the total complex of a finite exact
complex of complexes of $\S$\semimodule s.
 Part~(a) is proven; the proof of part~(b) is completely analogous.
\end{proof}

 Let us show that any $\S/\C/A$\projective{} left $\S$\semimodule{}
$\bL$ is $A$\projective.
 Consider the surjective morphism $\bcF(\bL)\rarrow\bL$ from
Lemma~3.
 Since its kernel is $\C/A$\injective, we have an extension of
an $\S/\C/A$\projective{} left $\S$\semimodule{} by
a $\C/A$\injective{} left $\S$\semimodule, which is always trivial
by Lemma~1.
 Therefore, $\bL$ is a direct summand of $\bcF(\bL)$, while
$\bcF(\bL)$ is $A$\projective{} by the construction.
 Analogously, any $\S/\C/A$\injective{} left $\S$\semicontramodule{}
is $A$\injective.

 Let us show that any $\S/\C/A$\projective{} $\C/A$\injective{} left
$\S$\semimodule{} $\bM$ is a direct summand of the $\S$\semimodule{}
induced from the $\C$\comodule{} coinduced from a projective
$A$\module; in particular, a left $\S$\semimodule{} is simultaneously
$\S/\C/A$\projective{} and $\C/A$\injective{} if and only if it is
semiprojective.
 Consider the exact triple $\bK\rarrow\S\oc_\C\bM\rarrow\bM$, where
$\bK=\ker(\S\oc_\C\bM\to\bM)$.
 If an $\S$\semimodule{} $\bM$ is $\C/A$\injective, then so is
the $\S$\semimodule{} $\S\oc_\C\bM$, since $\C/A$\+in\-jec\-tiv\-ity
is equivalent to $\C/A$\+co\-pro\-jec\-tiv\-ity; then
the $\S$\semimodule{} $\bK$ is $\C/A$\injective{} as a $\C$\comodule{}
direct summand of $\S\oc_\C\bM$.
 If the $\S$\semimodule{} $\bM$ is also $\S/\C/A$\projective, then
our exact triple splits over~$\S$ and $\bM$ is a direct summand
of the induced $\S$\semimodule{} $\S\oc_\C\bM$.
 Since the $\C$\comodule{} $\bM$ is $A$\projective{} and
$\C/A$\injective, it is a direct summand of the $\C$\comodule{}
coinduced from a projective $A$\module.
 Analogously, any $\S/\C/A$\injective{} $\C/A$\projective{} left
$\S$\semicontramodule{} $\bP$ is a direct summand of
the $\S$\semicontramodule{} coinduced from the $\C$\contramodule{}
induced from an injective $A$\module; in particular, a left
$\S$\semicontramodule{} is simultaneously $\S/\C/A$\injective{} and
$\C/A$\projective{} if and only if it is semiinjective.
 In other words, $\bM$ is a direct summand of a direct sum of copies of
the $\S$\semimodule{} $\S$ and $\bP$ is a direct summand of a product
of copies of the $\S$\semicontramodule{} $\Hom_k(\S,k\dual)$.

 An $\S/\C/A$\projective{} complex of $\C$\coprojective{} left
$\S$\semimodule s $\bM^\bu$ is homotopy equivalent to a complex
obtained from complexes of $\S$\semimodule s induced from complexes
of $\C$\coprojective{} $\C$\comodule s using the operations of
cone and infinite direct sum.
 In particular, the complex $\bM^\bu$ it is a semiprojective.
 Indeed, the total complex of the bicomplex $\dsb\rarrow
\S\oc_\C\S\oc_\C\bM\rarrow\S\oc_\C\bM\rarrow\bM$ is contractible,
being a $\C$\coacyclic{} $\S/\C/A$\projective{} complex of
$\C/A$\injective{} left $\S$\semimodule s.
 Analogously, an $\S/\C/A$\injective{} complex of $\C$\coinjective{}
left $\S$\semicontramodule s $\bP^\bu$ is homotopy equivalent to
a complex obtained from complexes of $\S$\semicontramodule s coinduced
from complexes of $\C$\coinjective{} $\C$\contramodule s using
the operations of cone and infinite product.
 In particular, the complex $\bP^\bu$ is semiinjective.

\smallskip
 Finally we turn to the construction of functorial factorizations.
 As in the proof of Theorem~\ref{co-contra-model-struct}, let us first
decompose an arbitrary morphism of complexes of left $\S$\semimodule s
$\bL^\bu\rarrow\bM^\bu$ into a cofibration followed by a fibration.
 This can be done in either of two dual ways.
 Let us start with an injective morphism from the complex $\bL^\bu$
into the complex $\bcJ^+(\bL^\bu)$ constructed in Lemma~2.
 Let $\bK^\bu$ be the cokernel of the morphism $\bL^\bu\rarrow\bM^\bu
\oplus\bcJ^+(\bL^\bu)$ and let $\bcF^+(\bK^\bu)\rarrow\bK^\bu$ be
the surjective morphism onto the complex $\bK^\bu$ from
the complex $\bcF^+(\bK^\bu)$ constructed in Lemma~3.
 Let $\bL^\bu\rarrow\bE^\bu\rarrow\bcF^+(\bK^\bu)$ be the pull-back of
the exact triple $\bL^\bu\rarrow\bM^\bu\oplus\bcJ^+(\bL^\bu)\rarrow
\bK^\bu$ with respect to the morphism $\bcF^+(\bK^\bu)\rarrow\bK^\bu$.
 Then the morphism $\bL^\bu\rarrow\bM^\bu$ is equal to the composition
$\bL^\bu\rarrow\bE^\bu\rarrow\bM^\bu$.
 The cokernel $\bcF^+(\bK^\bu)$ of the morphism $\bL^\bu\rarrow\bE^\bu$
is an $\S/\C/A$\projective{} complex of $\S/\C/A$\projective{}
$\S$\semimodule s.
 The kernel of the morphism $\bE^\bu\rarrow\bM^\bu$ is an extension
of the complex $\bcJ^+(\bL^\bu)$ and the kernel of the morphism
$\bcF^+(\bK^\bu)\rarrow\bK^\bu$, hence a complex of $\C/A$\injective{}
$\S$\semimodule s.
 Another way is to start with the surjective morphism $\bcF^+(\bM^\bu)
\rarrow\bM^\bu$ and consider the kernel of the morphism $\bL^\bu\oplus
\bcF^+(\bM^\bu)\rarrow\bM^\bu$.

 Now let us construct a factorization of the morphism $\bL^\bu\rarrow
\bM^\bu$ into a cofibration followed by a trivial fibration.
 Represent the kernel of the morphism $\bE^\bu\rarrow\bM^\bu$ as
the quotient complex of an $\S/\C/A$\projective{} complex of
$\S/\C/A$\projective{} left $\S$\semimodule s $\bcP^\bu$ by
a $\C$\coacyclic{} complex of $\C/A$\injective{} $\S$\semimodule s.
 Then the complex $\bcP^\bu$ is also a complex of $\C/A$\injective{}
$\S$\semimodule s (so it is even a semiprojective complex of
semiprojective $\S$\semimodule s).
 Let $\bK^\bu$ be the cone of the morphism $\bcP^\bu\rarrow\bE^\bu$.
 Then the morphism $\bL^\bu\rarrow\bM^\bu$ factorizes through $\bK^\bu$
in a natural way, the kernel of the morphism $\bK^\bu\rarrow\bM^\bu$
is a $\C$\coacyclic{} complex of $\C/A$\injective{} $\S$\semimodule s,
and the cokernel of the morphism $\bL^\bu\rarrow\bK^\bu$ is
an $\S/\C/A$\projective{} complex of $\S/\C/A$\projective{}
$\S$\semimodule s.

 It remains to construct a factorization of the morphism $\bL^\bu
\rarrow\bM^\bu$ into a trivial cofibration followed by a fibration.
 Represent the cokernel of the morphism $\bL^\bu\rarrow\bE^\bu$ as
a subcomplex of a complex of $\C/A$\injective{} $\S$\semimodule s
$\bcQ^\bu$ such that the quotient complex is a $\C$\coacyclic{}
$\S/\C/A$\projective{} complex of $\S/\C/A$\projective{}
$\S$\semimodule s.
 Then the complex $\bcQ^\bu$ is also an $\S/\C/A$\projective{}
complex of $\S/\C/A$\projective{} $\S$\semimodule s (hence
a semiprojective complex of semiprojective $\S$\semimodule s).
 Let $\bK^\bu$ be the cocone of the morphism $\bE^\bu\rarrow\bcQ^\bu$.
 Then the morphism $\bL^\bu\rarrow\bM^\bu$ factorizes through $\bK^\bu$
in a natural way, the kernel of the morphism $\bK^\bu\rarrow\bM^\bu$
is a complex of $\C/A$\injective{} $\S$\semimodule s, and the cokernel
of the morphism $\bL^\bu\rarrow\bK^\bu$ is a $\C$\coacyclic{}
$\S/\C/A$\projective{} complex of $\S/\C/A$\projective{}
$\S$\semimodule s.

 Part~(a) of Theorem is proven; the proof of part~(b) is completely
analogous.
\end{proof}

{\hbadness=5500
\begin{rmk1}
 One can obtain descriptions of $\S/\C/A$\projective{} complexes
of $\S/\C/A$\projective{} $\S$\semimodule s, $\C$\coacyclic{}
$\S/\C/A$\projective{} complexes of $\S/\C/A$\projective{}
$\S$\semimodule s, etc., from the proof of the above Theorem.
 Namely, let $\bM^\bu$ be an $\S/\C/A$\projective{} complex of
$\S/\C/A$\projective{} left $\S$\semimodule s; decompose the morphism
$0\rarrow\bM^\bu$ into a cofibration $0\rarrow\bK^\bu$ followed by
a trivial fibration $\bK^\bu\rarrow\bM^\bu$ by the above construction
(this can be also obtained directly from Lemma~3).
 Then the complex $\bM^\bu$ is a direct summand of $\bK^\bu$ and
therefore can be obtained from complexes of $\S$\semimodule s
induced from complexes of $A$\projective{} $\C$\comodule s using
the operations of cone, infinitely iterated extension in the sense
of inductive limit, and kernel of surjective morphism.
 Let $\bM^\bu$ be a $\C$\coacyclic{} $\S/\C/A$\projective{} complex
of $\S/\C/A$\projective{} left $\S$\semimodule s; decompose
the morphism $0\rarrow\bM^\bu$ into a trivial cofibration
$0\rarrow\bK^\bu$ followed by a fibration $\bK^\bu\rarrow\bM^\bu$ by
the above construction.
 Then the complex $\bM^\bu$ is a direct summand of $\bK^\bu$ and
therefore up to the homotopy equivalence can be obtained from the total
complexes of exact triples of $\S/\C/A$\projective{} complexes of
$\S/\C/A$\projective{} $\S$\semimodule s using the operations of cone
and infinite direct sum.
 The analogous results hold for complexes of left
$\S$\semicontramodule s.
\end{rmk1}}

 Let $\S$ and $\T$ be two semialgebras satisfying the above assumptions
and $\S\rarrow\T$ be a map of semialgebras compatible with a map of
corings $\C\rarrow\D$ and a $k$\+algebra map $A\rarrow B$.
 Then the pair of adjoint functors $\bM^\bu\mpsto{}_\T\bM^\bu$ and 
$\bN^\bu\mpsto{}_\C\bN^\bu$ is a Quillen adjunction from the category of
complexes of left $\S$\semimodule s to the category of complexes of left
$\T$\semimodule s; the pair of adjoint functors $\bQ^\bu\mpsto{}^\C
\bQ^\bu$ and $\bP^\bu\mpsto{}^\T\bP^\bu$ is a Quillen adjunction from
the category of complexes of left $\T$\semicontramodule s to
the category of complexes of left $\S$\semicontramodule s.
 These results follow from Theorems~\ref{co-contra-pull-well-defined}
and~\ref{semi-pull-push-adjusted}.3(c).
 They also hold in the case of a left coprojective and right coflat
Morita morphism $(\E,\E\dual)$ from $\C$ to $\D$ and a morphism
$\S\rarrow{}_\C\T_\C$ of semialgebras over~$\C$.

 The pair of adjoint functors $\Phi_\S$ and $\Psi_\S$ applied to
complexes term-wise is \emph{not} a Quillen equivalence, and not even
a Quillen adjunction, between the model category of complexes of
left $\S$\semicontramodule s and the model category of complexes of
left $\S$\semimodule s.
 Instead, this pair of adjoint functors between closed model categories
has the following properties.

 The functor $\Phi_\S$ maps trivial cofibrations of complexes of left
$\S$\semicontramodule s to weak equivalences of complexes of left
$\S$\semimodule s.
 The functor $\Psi_\S$ maps trivial fibrations of complexes of left
$\S$\semimodule s to weak equivalences of complexes of left
$\S$\semicontramodule s.
 For a cofibrant complex of left $\S$\semicontramodule s $\bP^\bu$ and
a fibrant complex of left $\S$\semimodule s $\bM^\bu$, a morphism
$\Phi_\S(\bP^\bu)\rarrow\bM^\bu$ is a weak equivalence if and only if
the corresponding morphism $\bP^\bu\rarrow\Psi_\S(\bM^\bu)$ is a weak
equivalence.
 Furthermore, the functor $\Phi_\S$ maps cofibrant complexes to
fibrant ones, while the functor $\Psi_\S$ maps fibrant complexes to
cofibrant ones.
 The restrictions of the functors $\Phi_\S$ and $\Psi_\S$ are mutually
inverse equivalences between the full subcategories formed by
cofibrant complexes of left $\S$\semicontramodule s and fibrant
complexes of right $\S$\semimodule s.
 These restrictions of functors also send weak equivalences to
weak equivalences.

\begin{rmk2}
 One can connect the above model categories of complexes of
left $\S$\semimodule s and left $\S$\semicontramodule s by a chain
of three Quillen adjunctions by considering other model category
structures on these two categories.
 The above model category structures on the category of complexes
of left $\S$\semimodule s can be called the semiprojective model
structure, and the model category structure on the category of
complexes of left $\S$\semicontramodule s can be called
the semiinjective model structure.
 In addition to these, there is also the injective model structure
on the category of complexes of left $\S$\semimodule s, and
the projective model structure on the category of complexes of
left $\S$\semicontramodule s.
 In these alternative model structures, weak equivalences are still
morphisms with $\C$\coacyclic{} or $\C$\contraacyclic{} cones,
respectively.
 A morphism of complexes of semimodules is a cofibration if and only
if it is injective, and a morphism of complexes of semicontramodules
is a fibration if and only if it is surjective.
 A morphism of complexes of semimodules is a fibration if and only if
it is surjective and its kernel is an injective complex of injective
semimodules in the sense of Remark~\ref{semi-ctrtor-definition} and
the proof of Lemma~1 above; a morphism of complexes of semicontramodules
is a cofibration if and only if it is injective and its cokernel is
an projective complex of projective semicontramodules.
 One checks that these a model category structures in the way analogous
to (and much simpler than) the proof of Theorem above.
 The functors $\Phi_\S$ and $\Psi_\S$ are a Quillen equivalence between
the injective and the projective model category structures; the identity
functors are Quillen equivalences between the semiprojective and
the injective model structures, and between the semiinjective and
the projective model structures.
\end{rmk2}

\Section{A Construction of Semialgebras}

\subsection{Construction of comodules and contramodules}

\subsubsection{}
 Let $A$ and $B$ be associative $k$\+algebras.

 For any projective finitely-generated left $A$\module{} $U$
and any left $A$\module{} $V$ there is a natural isomorphism
$\Hom_A(U,A)\ot_A V\simeq \Hom_A(U,V)$ given by the formula
$u^*\ot v\mpsto (u\mapsto\langle u,u^*\rangle v)$.
 In particular, for any $A$\+$B$\bimodule{} $C$ and any projective
finitely-generated left $B$\module{} $D$ there are natural
isomorphisms $\Hom_A(C\ot_B D\;A)\simeq\Hom_B(D,\Hom_A(C,A))
\simeq\Hom_B(D,B)\ot_B\Hom_A(C,A)$.

 It follows that there is a tensor anti-equivalence between the tensor
category of $A$\+$A$\bimodule s that are projective and
finitely-generated as left $A$\module s and the tensor category of
$A$\+$A$\bimodule s that are projective and finitely-generated as
right $A$\module s, given by the mutually-inverse functors
$C\mpsto\Hom_A(C,A)$ and $K\mpsto\Hom_{A^\rop}(K,A)$.
 Therefore, noncommutative ring structures on a right-projective
and finitely $A$\+$A$\bimodule{} $K$ correspond bijectively to
coring structures on the left-projective and finitely-generated
$A$\+$A$\bimodule{} $\Hom_{A^\rop}(K,A)$.
 On the other hand, for any coring $\C$ over $A$ there is a natural
structure of a $k$\+algebra on $\Hom_A(\C,A)$ together with a morphism
of $k$\+algebras $A\rarrow\Hom_A(\C,A)$.

 Furthermore, let $K$ be a $k$\+algebra endowed with a $k$\+algebra
map $A\rarrow K$ such that $K$ is a finitely-generated projective
right $A$\module, and let $\C=\Hom_{A^\rop}(K,A)$ be the corresponding
coring over~$A$.
 Then the natural isomorphism $N\ot_A\C=\Hom_{A^\rop}(K,N)$ for
a right $A$\module{} $N$ provides a bijective correspondence between 
the structures of right $K$\module{} and right $\C$\comodule{} on~$N$.
 Analogously, the natural isomorphism $\Hom_A(\C,P)=K\ot_A P$
for a left $A$\module{} $P$ provides a bijective correspondence between
the structures of left $K$\module{} and left $\C$\contramodule{} on~$P$.
 In other words, there are isomorphisms of abelian categories
$\comodr\C \simeq \modr K$ and $\C\contra \simeq K\modl$.

\subsubsection{}   \label{ring-coring-pairing}
 Here is a generalization of the situation we just described.
 Let $\C$ be a coring over a $k$\+algebra $A$ and $K$ be a $k$\+algebra
endowed with a $k$\+algebra map $A\rarrow K$.
 Suppose that we are given a pairing $\phi\:\C\ot_A K\rarrow A$ which
is an $A$\+$A$\bimodule{} map satisfying the following conditions of
compatibility with the comultiplication in~$\C$ and the multiplication
in~$K$ and with the counit of~$\C$ and the unit of~$K$.
 First, the composition $\C\ot_A K\ot_A K\rarrow\C\ot_A\C\ot_A
K\ot_A K\rarrow\C\ot_A K\rarrow A$ of the map induced by
the comultiplication in~$\C$, the map induced by the pairing~$\phi$,
and the pairing~$\phi$ itself should be equal to the composition
$\C\ot_A K\ot_A K\rarrow\C\ot_A K\rarrow A$ of the map induced by
the multiplication in~$A$ and the pairing~$\phi$.
 Second, the composition $\C=\C\ot_A A\rarrow\C\ot_A K\rarrow A$
of the map coming from the unit of~$K$ with the pairing~$\phi$ should
be equal to the counit of~$\C$.
 Equivalently, the map $K\rarrow\Hom_A(\C,A)$ induced by~$\phi$
should be a morphism of $k$\+algebras.

 Then for any right $\C$\comodule{} $\N$ the composition $\N\ot_A K
\rarrow\N\ot_A\C\ot_A K\rarrow\N$ of the map induced by
the $\C$\+coaction in $\N$ and the map induced by the pairing~$\phi$
defines a structure of right $K$\module{} on~$\N$.
 Analogously, for any left $\C$\contramodule{} $\P$ the composition
$K\ot_A\P\rarrow\Hom_A(\C,\P)\rarrow\P$ of the map given by the formula
$k'\ot p\mpsto (c\mapsto\phi(c,k')p)$ and the $\C$\+contraaction map
defines a structure of left $K$\module{} on~$\P$.
 So the pairing~$\phi$ induces faithful functors $\Delta_\phi\:
\comodr\C\rarrow\modr K$ and $\Delta^\phi\:\C\contra\rarrow K\modl$.

 In particular, a pairing~$\phi$ provides the coring $\C$ with
a structure of left $\C$\comodule{} endowed with a right action of
the $k$\+algebra $K$ by $\C$\comodule{} endomorphisms.
 Moreover, the data of a right action of $K$ by endomorphisms of
the left $\C$\comodule{} $\C$ agreeing with the right action of $A$
in~$\C$ is equivalent to the data of a pairing~$\phi$.

\subsubsection{}
 When $\C$ is a projective left $A$\module, the functor $\Delta^\phi$
has a left adjoint functor $\Gamma^\phi\:K\modl\rarrow\C\contra$.
 This functor sends the induced left $K$\module{} $K\ot_A V$ to
the induced left $\C$\contramodule{} $\Hom_A(\C,V)$ for any left
$A$\module{}~$V$; to compute $\Gamma^\phi(M)$ for an arbitrary left
$K$\module{} $M$, one can respresent $M$ as the cokernel of a morphism
of $K$\module s induced from $A$\module s.
 Analogously, when $\C$ is a flat right $A$\module, the functor 
$\Delta_\phi$ has a right adjoint functor $\Gamma_\phi\:
\modr K\rarrow\comodr\C$.
 This functor sends the coinduced right $K$\module{}
$\Hom_{A^\rop}(K,U)$ to the coinduced right $\C$\comodule{}
$U\ot_A\C$ for any right $A$\module{}~$U$; to compute
$\Gamma_\phi(N)$ for an arbitrary right $K$\module{} $N$, one can
represent $N$ as the kernel of a morphism of $K$\module s coinduced
from $A$\module s.

 Without any conditions on the coring~$\C$, the composition of
functors $\Psi_\C\:\C\comodl\allowbreak\rarrow\C\contra$ and
$\Delta^\phi\:\C\contra\rarrow K\modl$ has a left adjoint functor
$K\modl\rarrow\C\comodl$ mapping a left $K$\module{} $M$ to
the left $\C$\comodule{} $\C\ot_A M$.
 Analogously, the composition of the functors $\Phi_{\C^\rop}\:
\contraR\C \rarrow\comodr\C$ and $\Delta_\phi\:\comodr\C\rarrow
\modr K$ has a right adjoint functor $\modr K\rarrow\contraR\C$
mapping a right $K$\module{} $N$ to the right $\C$\contramodule{}
$\Hom_{K^\rop}(\C,N)$.
 So one can compute the compositions of functors $\Phi_\C\Gamma^\phi$
and $\Psi_{\C^\rop}\Gamma_\phi$ in this way.

\subsubsection{}
 It is easy to see that the functor $\Delta_\phi$ is fully faithful
whenever for any right $A$\module{} $N$ the map $N\ot_A \C\rarrow
\Hom_{A^\rop}(K,N)$ given by the formula $n\ot c\mpsto (k'\mapsto
n\phi(c,k'))$ is injective.
 In particular, when $A$ is a semisimple ring, the functor
$\Delta_\phi$ is fully faithful if the map $\C\rarrow
\Hom_{A^\rop}(K,A)$ induced by the pairing~$\phi$ is injective, i.~e.,
the pairing $\phi$ is nondegenerate in~$\C$.

\subsection{Construction of semialgebras}
\label{construction-of-semialgebras}

\subsubsection{}  \label{semialgebra-constructed}
 Assume that a coring $\C$ over a $k$\+algebra $A$ is a flat left
$A$\module.
 Let $K$ be a $k$\+algebra endowed with a $k$\+algebra map $A\rarrow K$
and a pairing $\phi\:\C\ot_A K\rarrow A$ satisfying the conditions
of~\ref{ring-coring-pairing}, and let $R$ be a $k$\+algebra endowed
with a $k$\+algebra map $f\:K\rarrow R$ such that $R$ 
is a flat left $K$\module.
 Then the tensor product $\C\ot_K R$ is a coflat left $\C$\comodule{}
endowed with a right action of the $k$\+algebra $K$ (and even of
the $k$\+algebra $R$) by left $\C$\comodule{} endomorphisms.

 Suppose that there exists a structure of right $\C$\comodule{}
on $\C\ot_K R$ inducing the existing structure of right $K$\module{}
and such that the following three maps are right $\C$\comodule{}
morphisms: (i)~the left $\C$\+coaction map $\C\ot_K R\rarrow
\C\ot_A(\C\ot_K R)$, (ii)~the semiunit map $\C=\C\ot_K K\rarrow
\C\ot_K R$, and (iii)~the semimultiplication map $(\C\ot_K R)\oc_\C
(\C\ot_K R)\simeq\C\ot_K R\ot_K R\rarrow\C\ot_K R$, where
the isomorphism in~(iii) is the inverse of the natural isomorphism of
Proposition~\ref{tensor-cotensor-assoc}(a) and the map being composed
is induced by the multiplication in~$R$.
 Then the semiunit and semimultiplication maps (ii) and~(iii) define
a semialgebra structure on the $\C$\+$\C$\bicomodule{} $\S=\C\ot_K R$.

 Notice that the maps (i\+iii) always preserve the right $K$\module{}
structures.
 If the functor $\Delta_\phi$ is fully faithful, then a right
$\C$\comodule{} structure inducing a given right $K$\module{} structure
on $\C\ot_K R$ is unique provided that it exists, and the maps (i\+iii)
preserve this structure automatically.
 If the functor $\Delta_\phi$ is an equivalence of categories, then
a unique right $\C$\comodule{} structure with the desired properties
exists on $\C\ot_K R$.

 The associativity of semimultiplication in $\S$ follows from
the associativity of multiplication in~$R$ and the commutativity
of diagrams of associativity isomorphisms of cotensor products.

\subsubsection{}   \label{semi-mod-contra-described}
 By Proposition~\ref{tensor-cotensor-assoc}(a), for any right
$\C$\comodule{} $\N$ there is a natural isomorphism $\N\oc_\C\S
\simeq\N\ot_K R$, hence every right $\S$\semimodule{} has a natural
structure of right $R$\module.
 So there is a faithful exact functor $\Delta_{\phi,f}\:
\simodr\S\allowbreak\rarrow\modr R$ which agrees with the functor
$\Delta_\phi\:\comodr\C\rarrow\modr K$.
 Moreover, the category of right $\S$\semimodule s is isomorphic to
the category of $k$\module s $\N$ endowed with a right
$\C$\comodule{} and right $R$\module{} structures satisfying 
the following compatibility conditions: first, the induced right
$K$\module{} structures should coincide, and second, the action map
$\N\ot_K R\rarrow\N$ should be a morphism of right $\C$\comodule s,
where the right $\C$\comodule{} structure on $\N\ot_K R$ is provided
by the isomorphism $\N\ot_K R=\N\oc_\C\S$.
 When the functor $\Delta_\phi$ is fully faithful, the category
$\simodr\S$ is simply described as the full subcategory of
the category of right $R$\module s consisting of those modules whose
right $K$\module{} structure comes from a right $\C$\comodule{}
structure.

 Analogously, if $\C$ is a projective left $A$\module{} and $R$ is
a projective left $K$\module, then $\S$ is a coprojective left
$\C$\comodule{} and by Proposition~\ref{co-tensor-co-hom-assoc}.2(a) 
for any left $\C$\contramodule{} $\P$ there is a natural isomorphism
$\Cohom_\C(\S,\P)\simeq\Hom_K(R,\P)$, hence any left
$\S$\semicontramodule{} has a natural structure of left $R$\module.
 So there is a faithful exact functor $\Delta^{\phi,f}\:\S\sicntr
\rarrow R\modl$ which agrees with the functor $\Delta^\phi\:\C\comodl
\rarrow K\modl$.
 Moreover, the category of left $\S$\semicontramodule s is isomorphic
to the category of $k$\module s $\P$ endowed with a left
$\C$\contramodule{} and a left $R$\module{} structures satisfying
the following compatibility conditions: first, the induced left
$K$\module{} structures should coincide, and second, the action map
$\P\rarrow\Hom_K(R,\P)$ should be a morphism of $\C$\contramodule s,
where the left $\C$\contramodule{} structure on $\Hom_K(R,\P)$
is provided by the isomorphism $\Hom_K(R,\P)=\Cohom_\C(\S,\P)$.

\subsubsection{}   \label{contratensor-described}
 When $K$ is a projective finitely-generated right $A$\module{} and
the pairing~$\phi$ corresponds to an isomorphism $\C\simeq
\Hom_{A^\rop}(K,A)$, the isomorphisms of categories
$\Delta_\phi\:\comodr\C\simeq \modr K$ and $\Delta^\phi\:\C\contra
\simeq K\modl$ transform the functor of contratensor product over~$\C$
into the functor of tensor product over~$K$: \ $\N\ocn_\C\P \simeq
\N\ot_K\P$.
 Indeed, $\Hom_A(\C,\P)= K\ot_A P$.
 When in addition $R$ is a projective left $K$\module, the isomorphisms
of categories $\Delta_{\phi,f}\:\simodr\S \simeq \modr R$ and
$\Delta^{\phi,f}\:\S\sicntr \simeq R\modl$ transform the functor of
contratensor product over~$\S$ into the functor of tensor product
over~$R$: \ $\bN\Ocn_\S\bP\simeq\bN\ot_R\bP$.
 Indeed, $\bN\oc_\C\S=\bN\ot_K R$.

\subsubsection{}
 The functor $\Delta_{\phi,f}$ has a right adjoint functor
$\Gamma_{\phi,f}\:\modr R\rarrow\simodr\S$, which agrees with
the functor $\Gamma_\phi\:\modr K\rarrow\comodr\C$.
 The functor $\Gamma_{\phi,f}$ is constructed as follows.
 Let $N$ be a right $R$\module; it has an induced right $K$\module{}
structure.
 Consider the composition $\Delta_\phi(\Gamma_\phi(N)\oc_\C\S)=
\Delta_\phi\Gamma_\phi(N)\ot_K R\rarrow N\ot_K R\rarrow N$ of 
the isomorphism of Proposition~\ref{tensor-cotensor-assoc}(a),
the map induced by the adjunction map $\Delta_\phi\Gamma_\phi(N)
\rarrow N$, and the right action map.
 By adjunction, this composition corresponds to a right
$\C$\comodule{} morphism $\Gamma_\phi(N)\oc_\C\S\rarrow\Gamma_\phi(N)$,
which provides a right $\S$\semimodule{} structure on $\Gamma_\phi(N)$.

 Analogously, if $\C$ is a projective left $A$\module{} and $R$ is
a projective left $K$\module, then the functor $\Delta^{\phi,f}$ has
a left adjoint functor $\Gamma^{\phi,f}\: R\modl\rarrow\S\sicntr$,
which agrees with the functor $\Gamma^\phi\: K\modl\rarrow\C\contra$.
 The functor $\Gamma^{\phi,f}$ is constructed as follows.
 Let $P$ be a left $R$\module; it has an induced left $K$\module{}
structure.
 Consider the composition $P\rarrow\Hom_K(R,P)\rarrow
\Hom_K(R\;\Delta^\phi\Gamma^\phi(P))=\Delta^\phi(\Cohom_\C(\S,
\Gamma^\phi(P)))$ of the action map, the map induced by the adjunction
map $P\rarrow\Delta^\phi\Gamma^\phi(P)$, and the isomorphism of
Proposition~\ref{co-tensor-co-hom-assoc}.2(a).
 By adjunction, this composition corresponds to a left
$\C$\contramodule{} morphism $\Gamma^\phi(P)\rarrow
\Cohom_\C(\S,\Gamma^\phi(P))$, which provides a left
$\S$\semicontramodule{} structure on $\Gamma^\phi(P)$.

 Notice that the semialgebra $\S$ has a structure of left
$\S$\semimodule{} endowed with a right action of the $k$\+algebra
$R$ by left $\S$\semimodule{} endomorphisms.
 So when $\C$ is a flat right $A$\module{} and $\S$ turns out to be
a coflat right $\C$\comodule, there is a functor $\S\simodl\rarrow
R\modl$ mapping a left $\S$\semimodule{} $\bM$ to the left
$R$\module{} $\Hom_\S(\S,\bM)$.
 This functor has a left adjoint functor mapping a left $R$\module{}
$M$ to the left $\S$\semimodule{} $\S\ot_R M=\C\ot_K M$.
 In the case when $\C$ is a projective left $A$\module{} and
$R$ is a projective left $K$\module, the former functor is isomorphic
to $\Delta^{\phi,f}\Psi_\S$, and consequently the latter functor is
isomorphic to $\Phi_\S\Gamma^{\phi,f}$.
 Analogously, when $\C$ is a projective right $A$\module{} and $\S$
turns out to be a coprojective right $\C$\comodule, the functor
$\Psi_{\S^\rop}\Gamma_{\phi,f}$ maps a right $R$\module{} $N$
to the right $\S$\semicontramodule{} $\Hom_{R^\op}(\S,N)=
\Hom_{K^\op}(\C,N)$.

 Let us point out that \emph{no explicit description of the category
of left\/ $\S$\semimodule s is in general available}.
 We only described the categories of right $\S$\semimodule s and left
$\S$\semicontramodule s, and constructed certain functors acting to
and from the category of left $\S$\semimodule s.

\subsubsection{}
 The following observations were inspired by~\cite[section~5]{Ar2}.

 Suppose that there is a commutative diagram of $k$\+algebra maps
$A\rarrow K$, \ $K\rarrow R$, \ $A'\rarrow K'$, \ $K'\rarrow R'$, \
$A\rarrow A'$, \ $K\rarrow K'$, \ $R\rarrow R'$.
 Let $\C$ be a coring over~$A$ and $\C'$ be a coring over~$A'$
endowed with a map of corings $\C\rarrow \C'$ compatible with
the $k$\+algebra map $A\rarrow A'$.
 Assume that $\C$ is a flat left $A$\module, $\C'$ is a flat
left $A'$\module, $R$ is a flat left $K$\module, and $R'$ is
a flat left $K'$\module.
 Let $\phi\:\C\ot_A K\rarrow A$ and $\phi'\:\C'\ot_{A'}K'\rarrow
A'$ be two pairings satisfying the conditions
of~\ref{ring-coring-pairing} and forming a commutative diagram
with the maps $\C\ot_A K\rarrow\C'\ot_{A'}K'$ and $A\rarrow A'$.
 Furthermore, suppose that the natural map $K\ot_A A'\rarrow K'$
is an isomorphism.
 Assume that there is a structure of right $\C$\comodule{} on
$\C\ot_K R$ and a structure of right $\C'$\comodule{} on
$\C'\ot_{K'}R'$ satisfying the conditions
of~\ref{semialgebra-constructed}, so that $\C\ot_K R$ is a semialgebra
over~$\C$ and $\C'\ot_{K'}R'$ is a semialgebra over~$\C'$.
 Then the natural map from the right $\C'$\comodule{}
$\C\ot_K R\ot_A A'$ to the right $\C'$\comodule{} $\C'\ot_{K'}R'$
is a morphism of right $K'$\module s.
 If it is also a morphism of right $\C'$\comodule s, then the map
$\C\ot_K R\rarrow\C'\ot_{K'}R'$ is a map of semialgebras compatible
with the map of corings $\C\rarrow \C'$ and the $k$\+algebra map
$A\rarrow A'$.

 Suppose that there is a commutative diagram of $k$\+algebra maps
$A\rarrow K$, \ $K\rarrow R$, \ $A\rarrow K'$, \ $K'\rarrow R'$, \
$K\rarrow K'$, \ $R\rarrow R'$.
 Let $\C$ and $\C'$ be two corings over~$A$ and $\C'\rarrow\C$ be
a morphism of corings over~$A$.
 Assume that $\C$ and $\C'$ are flat left $A$\module s, $R$ is
a flat left $K$\module, and $R'$ is a flat left $K'$\module.
 Let $\phi\:\C\ot_A K\rarrow A$ and $\phi'\:\C'\ot_A K'\rarrow A$ be
two pairings satisfying the conditions of~\ref{ring-coring-pairing}
and forming a commutative diagram with the maps $\C'\ot_A K\rarrow
\C\ot_A K$ and $\C'\ot_A K\rarrow \C'\ot_A K'$.
 Furthermore, suppose that the natural map $K'\ot_K R\rarrow R'$
is an isomorphism.
 Assume that there is a structure of right $\C$\comodule{} on
$\C\ot_K R$ and a structure of right $\C'$\comodule{} on
$\C'\ot_{K'}R'$ satisfying the conditions
of~\ref{semialgebra-constructed}, so that $\C\ot_K R$ is a semialgebra
over~$\C$ and $\C'\ot_{K'}R'$ is a semialgebra over~$\C'$.
 In this case, if the right $K$\module{} map $\C'\ot_{K'} R'
=\C'\ot_{K'}K'\ot_K R\simeq \C'\ot_K R\rarrow \C\ot_K R$ is
a morphism of right $\C$\comodule s, then it is a map of semialgebras
compatible with the morphism $\C'\rarrow\C$ of corings over~$A$.

\subsection{Entwining structures}  \label{entwining-structures}
 An important particular case of the above construction of semialgebras
was considered in~\cite{Brz}.
 Namely, it was noticed that there is a set of data from which one can
construct \emph{both} a coring and a semialgebra.

\subsubsection{}
 Let $\C$ be a coring over a $k$\+algebra $A$ and $A\rarrow B$ be
a morphism of $k$\+algebras.
 A \emph{right entwining structure} for $\C$ and $B$ over $A$ is
an $A$\+$A$\bimodule{} map $\psi\:\C\ot_A B\rarrow B\ot_A\C$ satisfying
the following equations: (i)~the composition $\C\ot_A B\ot_A B\rarrow
B\ot_A\C\ot_A B\rarrow B\ot_A B\ot_A\C\rarrow B\ot_A\C$ of two maps
induced by the map~$\psi$ and the map induced by the multiplication
in~$B$ is equal to the composition $\C\ot_A B\ot_A B\rarrow\C\ot_A B
\rarrow B\ot_A\C$ of the map induced by the multiplication in~$B$
and the map~$\psi$; (ii)~the map~$\psi$ forms a commutative triangle
with the maps $\C\rarrow\C\ot_A B$ and $\C\rarrow B\ot_A\C$ coming
from the unit of~$B$; (iii)~the composition $\C\ot_A B\rarrow
\C\ot_A\C\ot_A B\rarrow\C\ot_A B\ot_A\C\rarrow B\ot_A B\ot_A\C$
of the map induced by the comultiplication in~$\C$ and two maps
induced by the map~$\psi$ is equal to the composition $\C\ot_A B
\rarrow B\ot_A\C\rarrow B\ot_A\C\ot_A\C$ of the map~$\psi$ and
the map induced by the comultiplication in~$\C$; and (iv)~the map
$\psi$ forms a commutative triangle with the maps $\C\ot_A B\rarrow B$
and $B\ot_A\C\rarrow B$ coming from the counit of~$\C$.

 A \emph{left entwining structure} for $\C$ and $B$ over $A$ is
defined as an $A$\+$A$\bimodule{} map $\psi^\#\:B\ot_A\C\rarrow
\C\ot_A B$ satisfying the opposite equations.
 Notice that whenever a map $\psi\:\C\ot_A B\rarrow B\ot_A\C$ is
invertible the map $\psi$ is a right entwining structure if and only
if the map $\psi^\#=\psi^{-1}$ is a left entwining structure.

\subsubsection{}
 A (right) \emph{entwined module} over a right entwining structure
$\psi\: \C\ot_A B\rarrow B\ot_A \C$ is a $k$\module{} $\N$ endowed
with a right $\C$\comodule{} and a right $B$\module{} structures such
that the corresponding right $A$\module{} structures coincide and
the following equation holds: the composition $\N\ot_A B\rarrow
\N\ot_A\C\ot_A B\rarrow\N\ot_A B\ot_A\C\rarrow\N\ot_A\C$ of
the map induced by the $\C$\+coaction in~$\N$, the map induced by
the map~$\psi$, and the map induced by the $B$\+action in~$\N$
is equal to the composition $\N\ot_A B\rarrow\N\rarrow\N\ot_A\C$
of the $B$\+action map and the $\C$\+coaction map.

 A (left) \emph{entwined contramodule} over a right entwining
structure~$\psi$ is a $k$\module{} $\P$ endowed with a left
$\C$\contramodule{} and a left $B$\module{} structures such that
the corresponding left $A$\module{} structures coincide and
the following equation holds: the composition $\Hom_A(\C,\P)\rarrow
\Hom_A(\C,\Hom_A(B,\P))=\Hom_A(B\ot_A\C\;\P)\rarrow
\Hom_A(\C\ot_A B\;\P)=\Hom_A(B,\Hom_A(\C,\P))\rarrow\Hom_A(B,\P)$
of the map induced by the $B$\+action in~$\P$, the map induced by
the map~$\psi$, and the map induced by the $\C$\+contraaction{} in~$\P$
is equal to the composition $\Hom_A(\C,\P)\rarrow\P\rarrow\Hom_A(B,\P)$
of the $\C$\+contraaction map and the $B$\+action map.

 (Left) \emph{entwined modules} and (right) \emph{entwined
contramodules} over a left entwining structure are defined
in the analogous way.

\subsubsection{}
 Let $\psi\:\C\ot_A B\rarrow B\ot_A \C$ be a right entwining structure.
 Define a coring $\D$ over $B$ as the left $B$\module{} $B\ot_A\C$
endowed with the following right action of~$B$, comultiplication, and
counit.
 The right $B$\+action is the composition $(B\ot_A\C)\ot_A B\rarrow
B\ot_A B\ot_A\C\rarrow B\ot_A\C$ of the map induced by the map~$\psi$
and the multiplication in $B$.
 The comultiplication is the map $B\ot_A\C\rarrow B\ot_A\C\ot_A\C
=(B\ot_A\C)\ot_B(B\ot_A\C)$ induced by the comultiplication
in $\C$.
 The counit is the map $B\ot_A\C\rarrow B\ot_A A=B$ coming
from the counit of $\C$.
 One has to use the equation~(i) on the entwining map~$\psi$ to
check that the right action of~$B$ is associative, the equation~(ii)
to check that the right action of~$B$ agrees with the existing
right action of~$A$, and the equations (iii) and~(iv) to check that
the comultiplication and counit are right $B$\module{} maps.

 Analogously, for a left entwining structure $\psi^\#\:B\ot_A\C\rarrow
\C\ot_A B$ one defines a coring $\D^\#=\C\ot_A B$ over $B$.
 When $\psi^\#=\psi^{-1}$ are two inverse maps satysfying the entwining
structure equations, the maps $\psi$ and $\psi^\#$ themselves are
mutually inverse isomorphisms $\D^\#\simeq\D$ between the corresponding
corings over $B$.

\subsubsection{}
 Let $\psi\:\C\ot_A B\rarrow B\ot_A \C$ be a right entwining structure.
 Define a semialgebra $\S$ over $\C$ as the left $\C$\comodule{}
$\C\ot_A B$ endowed with the following right coaction of~$\C$,
semimultiplication, and semiunit.
 The right $\C$\+coaction is the composition $\C\ot_A B\rarrow
\C\ot_A\C\ot_A B\rarrow(\C\ot_A B)\ot_A\C$ of the map induced by
the comultiplication in $\C$ and the map induced by the map~$\psi$.
 The semimultiplication is the map $(\C\ot_A B)\oc_\C(\C\ot_A B)
=\C\ot_A B\ot_A B\rarrow\C\ot_A B$ induced by the multiplication
in $B$.
 The semiunit is the map $\C=\C\ot_A A\rarrow\C\ot_A B$ coming
from the unit of $B$.
 The multiple cotensor products $\N\oc_\C\S\oc_\C\S\oc_\C\dsb
\oc_\C\S$ and the multiple cohomomorphisms $\Cohom_\C(\S\oc_\C\dsb
\oc_\C\S\;\P)$ are associative for any right $\C$\comodule{} $\N$
and any left $\C$\contramodule{} $\P$ by
Propositions~\ref{cotensor-associative}(e)
and~\ref{cohom-associative}(h).

 Analogously, for a left entwining structure $\psi^\#\:B\ot_A\C\rarrow
\C\ot_A B$ one defines a semialgebra $\S^\#=B\ot_A\C$ over $\C$.
 When $\psi^\#=\psi^{-1}$ are two inverse maps satysfying the entwining
structure equations, the maps $\psi$ and $\psi^\#$ themselves are
mutually inverse isomorphisms $\S\simeq\S^\#$ between the corresponding
semialgebras over $\C$.

\subsubsection{}  \label{entwined-coring-semialgebra}
 An entwined module over a right entwining structure~$\psi$ is the same
that a right $\D$\comodule{} and the same that a right $\S$\semimodule;
in other words, the corresponding categories are isomorphic.
 Analogously, an entwined module over a left entwining
structure~$\psi^\#$ is the same that a left $\D^\#$\comodule{} and
the same that a left $\S^\#$\semimodule.
 Similar assertions apply to contramodules: an entwined contramodule
over a right entwining structure~$\psi$ is the same that a left
$\D$\contramodule{} and the same that a left $\S$\semicontramodule;
analogously for a left entwining structure.

 For any entwined module $\N$ over a right entwining structure $\psi$
there is a natural injective morphism $\N\rarrow\N\ot_B\D\simeq
\N\ot_A\C$ from $\N$ into an entwined module which as
a $\C$\comodule{} is coinduced from an $A$\module.
 Analogously, for any left entwined contramodule $\P$ over $\psi$
there is a natural surjective morphism $\Hom_A(\C,\P)\simeq
\Hom_B(\D,\P)\rarrow\P$ onto $\P$ from an entwined contramodule
which as a $\C$\contramodule{} is induced from an $A$\module.
 So we obtain, \emph{in the entwining structure case}, a functorial
injection from an arbitrary $\S$\semimodule{} into a $\C/A$\injective{}
$\S$\semimodule{} and a functorial surjection onto an arbitrary
$\S$\semicontramodule{} from a $\C/A$\projective{}
$\S$\semicontramodule{} constructed in a way much simpler than that
of Lemmas~\ref{coflat-semimodule-injection}
and~\ref{coproj-coinj-semi-mod-contra} (cf.~\cite{Ar1,Ar2}).

 When the ring $A$ is semisimple, there is also a functorial
surjection onto an arbitrary $\D$\comodule{} $\N$ from
a $B$\projective{} $\D$\comodule{} $\N\oc_\C\S\simeq\N\ot_A B$
and a functorial injection from an arbitrary $\D$\contramodule{} $\P$
into a $B$\injective{} $\D$\contramodule{} $\Cohom_\C(\S,\P)\simeq
\Hom_A(B,\P)$; these are much simpler constructions than those of
Lemmas~\ref{flat-comodule-surjection}
and~\ref{proj-inj-co-contra-module}.

 When $B$ is a flat right $A$\module, the construction of
the semialgebra $\S=\C\ot_A B$ corresponding to an entwining
structure~$\psi$ becomes a particular case of the construction
of the semialgebra $\S=\C\ot_K R$ corresponding to a pairing~$\phi$
(take $K=A$, \ $R=B$, and the only possible~$\phi$).

\subsubsection{}
 When $\psi^\#=\psi^{-1}$ are two inverse entwining structures,
there is an explicit description of \emph{both} the categories of
left and right comodules over $\D^\#\simeq\D$ and \emph{both}
the categories of left and right semimodules over $\S\simeq\S^\#$.

 When $\psi$ is invertible, the multiple cotensor products
$\N\oc_\C\S\oc_\C\dsb\oc_\C\S\oc_\C\M$ and the multiple
cohomomorphisms $\Cohom_\C(\S\oc_\C\dsb \oc_\C\S\oc_\C\M\;\P)$
are associative for any right $\C$\comodule{} $\N$, left
$\C$\comodule{} $\M$, and left $\C$\contramodule{} $\P$ by
Propositions~\ref{cotensor-associative}(f)
and~\ref{cohom-associative}(j), so the functors of semitensor
product and semihomomorphism over~$\S$ are everywhere defined.

\subsection{Semiproduct and semimorphisms}
\label{semi-product-morphisms}
 Let $\psi\:\C\ot_A B\rarrow B\ot_A\C$ be a right entwining structure;
suppose that $\psi$ is an invertible map.
 Let $\S=\C\ot_A B$ and $\D=B\ot_A\C$ be the corresponding semialgebra
over~$\C$ and coring over~$B$.

 One defines~\cite{Sev} the \emph{semiproduct} $\N\ot_B^\C \M$ of
a right entwined module $\N$ over~$\psi$ and a left entwined module
$\M$ over~$\psi^{-1}$ as the image of the composition of maps
$\N\oc_\C\M\rarrow \N\ot_A\M\rarrow\N\ot_B\M$.
 Analogously, one defines the $k$\module{} of \emph{semimorphisms}
$\Hom_B^\C(\M,\P)$ from a left entwined module $\M$ over~$\psi^{-1}$
to a left entwined contramodule $\P$ over~$\psi$ as the image of
the composition of maps $\Hom_B(\M,\P)\rarrow\Hom_A(\M,\P)\rarrow
\Cohom_\C(\M,\P)$.

 There is a natural map of semialgebras $\S\rarrow B$ compatible
with the map $\C\rarrow A$ of corings over~$A$.
 Hence for any entwined modules $\N$ over~$\psi$ and $\M$
over~$\psi^{-1}$ there is a natural injective map from the pair of
morphisms $\N\oc_\C\S\oc_\C\M\birarrow\N\oc_\C\M$ to the pair of
morphisms $\N\ot_A B\ot_A\M\birarrow\N\ot_A\M$ .
 Therefore, we have a natural surjective map $\N\os_\S\M\rarrow
\N\ot_B^\C\M$, which is an isomorphism if and only if the map
$\N\os_\S\M\rarrow\N\ot_B\M$ is injective.
 Analogously, for any entwined module $\M$ over~$\psi^{-1}$ and
entwined contramodule $\P$ over~$\psi$ there is a natural surjective
map from the pair of morphisms $\Hom_A(\M,\P)\birarrow
\Hom_A(B\ot_A\M\;\P)$ to the pair of morphisms $\Cohom_\C(\M,\P)
\birarrow\Cohom_\C(\S\oc_\C\M\;\P)$.
 So we get a natural injective map $\Hom_B^\C(\M,\P)\rarrow
\SemiHom_\S(\M,\P)$, which is an isomorphism if and only if the map
$\Hom_B(\M,\P)\rarrow\SemiHom_\S(\M,\P)$ is surjective.

 Consider the natural injective morphism of entwined modules
$\N\rarrow\N\ot_B\D=\N\ot_A\C$.
 Taking the semitensor product of this morphism with $\M$ over~$\S$,
we obtain the map $\N\os_\S\M\rarrow(\N\ot_A\C)\os_\S\M\simeq
\N\ot_B\M$ that we are interested in.
 Hence the natural map $\N\os_\S\M\rarrow\N\ot_B^\C\M$ is
an isomorphism whenever the semitensor product with $\M$ maps
$A$\+split injections of right $\S$\semimodule s to injections
or $\N$ has such property with respect to left $\S$\semimodule s.
 This includes the cases when $\N$ or $\M$ is an $\S$\semimodule{}
induced from a $\C$\comodule.

 Analogously, consider the natural surjective morphism of entwined
contramodules $\Hom_A(\C,\P)=\Hom_B(\D,\P)\rarrow\P$.
 The map $\Hom_B(\M,\P)\rarrow\SemiHom_\S(\M,\P)$ can be obtained
by taking the semihomomorphisms over~$\S$ from $\M$ to the morphism
$\Hom_A(\C,\P)\rarrow\P$ or by taking the semihomomorphisms over~$\S$
from the morphism $\M\rarrow\C\ot_A\M$ to~$\P$.
 Thus the natural map $\Hom_B^\C(\M,\P)\rarrow\SemiHom_\S(\M,\P)$
is an isomorphism whenever the functor of semihomomorphisms from $\M$
maps $A$\+split surjections of left $\S$\semicontramodule s to
surjections or the functor of semihomomorphisms into $\P$ maps
$A$\+split injections of left $\S$\semimodule s to surjections.
 This includes the cases when $\M$ is an $\S$\semimodule{} induced
from a $\C$\comodule{} or $\P$ is an $\S$\semicontramodule{}
coinduced from a $\C$\contramodule. 

 In the same way one constructs a natural injective map $\N\ot_B^\C\M
\rarrow\N\oc_\D\M$ and shows that it is an isomorphism whenever
the cotensor product with $\N$ or $\M$ over $\D$ maps surjections
of $\D$\comodule s to surjections, in particular, when one of
the $\D$\comodule s $\M$ or $\N$ is quasicoflat.
 Analogously, there is a natural surjective map $\Cohom_\D(\M,\P)
\rarrow\Hom_B^\C(\M,\P)$, which is an isomorphism whenever the functor
of cohomomorphisms from $\M$ over $\D$ maps injections of left
$\D$\contramodule s to injections or the functor of cohomomorphisms
into $\P$ over $\D$ maps surjections of left $\D$\comodule s to
injections, in particular, when $\M$ is a quasicoprojective
$\D$\comodule{} or $\P$ is a quasicoinjective $\D$\contramodule.

\Section{Relative Nonhomogeneous Koszul Duality}
\label{nonhom-koszul-section}

\subsection{Graded semialgebras}  \label{graded-semi}

\subsubsection{}
 All the constructions of Sections~1--10 can be carried out with
the category of $k$\module s replaced by the category of graded
$k$\module s.

 So one would consider a graded $k$\+algebra $A$, a coring
object $\C$ in the tensor category of graded $A$\+$A$\bimodule s,
a ring object $\S$ in a tensor category of graded
$\C$\+$\C$\bicomodule s, assume $A$ to have a finite graded
homological dimension, consider graded $\S$\semimodule s
and graded $\S$\semicontramodule s.
 All of our definitions and results can be transfered to the graded
situation without any difficulties.
 All the functors so obtained commute with the shift of grading
in modules.

 Furthermore, there are \emph{two} forgetful functors $\Sigma$ and
$\Pi$ from the category of graded $k$\module s $k\modl^\sgr$ to
the category $k\modl$, the functor $\Sigma$ sending a graded
$k$\module{} to the infinite direct sum of its components and
the functor $\Pi$ sending it to their infinite product.
 For any graded semialgebra $\S$ over a graded coring $\C$ over
a graded $k$\+algebra $A$, there are natural structures of
a $k$\+algebra on $\Sigma A$, of a coring over $\Sigma A$ on
$\Sigma \C$, and of a semialgebra over $\Sigma \C$ on $\Sigma \S$.
 For any graded $\S$\semimodule{} $\bM$ there is a natural structure
of a $\Sigma\S$\semimodule{} on $\Sigma\bM$ and for any graded
$\S$\semicontramodule{} $\bP$ there is a natural structure of
a $\Sigma\S$\semicontramodule{} on $\Pi\bP$.

 The functors of semitensor product and semihomomorphism defined
in the graded setting are related to their ungraded analogues by
the formulas $\Sigma(\bN\os_\S^\gr\bM)\simeq\Sigma\bN\os_{\Sigma\S}
\Sigma\bM$ and $\Pi\SemiHom_\S^\gr(\bM,\bP)\simeq\SemiHom_{\Sigma\S}
(\Sigma\bM,\Pi\bP)$.
 The functors $\bN\mpsto{}_\C\bN$ and $\bM\mpsto{}_\T\bM$ commute
with the forgetful functors $\Sigma$ and the functors
$\bQ\mpsto{}^\C\bQ$ and $\bP\mpsto{}^\T\bP$ commute with
the forgetful functors $\Pi$.
 The corresponding derived functors $\SemiTor$, $\SemiExt$, etc.,
have the analogous properties.
 However, the functors $\Hom_\S$, $\Hom^\S$, $\CtrTor^\S$, $\Psi_\S$,
and $\Phi_\S$ and their derived functors have \emph{no} properties of
compatibility with the functors of forgetting the grading.
 Thus one has to be aware of the distinction between $\Hom_\S$ and
$\Hom_\S^\gr$, \ $\Phi_\S$ and $\Phi_\S^\gr$, etc.

\subsubsection{}
 Assume that $A$ is a nonnegatively graded $k$\+algebra, $\C$ is
a nonnegatively graded coring over~$A$, and $\S$ is a nonnegatively
graded semialgebra over $\C$.
 Let $\S\simodl^\up$ and $\simodrup\S$ denote the categories of
nonnegatively graded $\S$\semimodule s, and $\S\contra^\down$ denote
the category of nonpositively graded $\S$\semicontramodule s.
 
 All the constructions of Sections~1--4 in their graded versions
preserve the categories of comodules and semimodules graded by
nonnegative integers and the categories of contramodules and
semicontramodules graded by nonpositive integers.
 All the definitions and results of these sections can be transfered
to the described situation of bounded grading and no problems occur.
 In particular, one can apply Lemma~\ref{semitor-definition} to define
the functors SemiTor and SemiExt in the bounded grading case.
 Moreover, the functors so obtained agree with the functors
$\SemiTor^\S_\gr$ and $\SemiExt_\S^\gr$ defined in terms of complexes
with unbounded grading.
 This is so because the constructions of resolutions agree.
 For the same reasons, in the assumptions
of~\ref{semimodule-semicontramodule-subsect} the functors
$\sD^\si(\S\simodl^\up)\rarrow \sD^\si(\S\simodl^\sgr)$ and
$\sD^\si(\S\sicntr^\down)\rarrow \sD^\si(\S\sicntr^\sgr)$ are
fully faithful, and the functor CtrTor defined by applying
Lemma~\ref{semi-ctrtor-definition}.2 in the bounded grading case
agrees with the functor $\CtrTor^\S_\gr$. 
 But the functors $\Psi_\S^\gr$ and $\Phi_\S^\gr$ do \emph{not}
preserve the bounded grading.

\subsection{Differential semialgebras}  \label{diff-semialgebras}

\subsubsection{}
 Let $B$ be a graded $k$\+algebra endowed with an odd derivation $d_B$
of degree~$1$ and $\D$ be a graded coring over $B$.
 A homogeneous map $d_\D\:\D\rarrow \D$ of degree~$1$ is called
a \emph{coderivation of\/ $\D$ with respect to\/ $d_B$} if
the biaction map $B\ot_k\D\ot_k B\rarrow \D$ and the comultiplication map
$\D\rarrow \D\ot_B\D$ are morphisms in the category of graded $k$\module s
endowed with endomorphisms of degree~$1$, where the endomorphisms of
the tensor products are defined by the usual super-Leibniz rule
$d(xy)=d(x)y+(-1)^{|x|}xd(y)$ (the degree of a homogeneous element
$x$ being denoted by $|x|$).
 In this case, it follows that the counit map $\D\rarrow B$
satisfies the same condition.
 In the particular case when $B$ is concentrated in the degree~$0$
and $d_B=0$, the condition on the biaction map simply means that
$d_\D$ is a $B$\+$B$\bimodule{} morphism.

 Now assume that $B$ is a DG\+algebra over~$k$, i.~e., $d_B^2=0$.
 A \emph{DG\+coring} over $B$ is a graded coring $\D$ over the graded
ring $B$ endowed with a coderivation $d_\D\:\D\rarrow\D$ with respect
to~$d_B$ such that $d_\D^2=0$.

 Let $\D$ be a DG\+coring over a DG\+algebra $B$.
 Then the cohomology $H(\D)$ is endowed with a natural structure of
a graded coring over the graded algebra $H(B)$ provided that
the natural maps $H(\D)\ot_{H(B)}H(\D)\rarrow H(\D\ot_B\D)$ and
$H(\D)\ot_{H(B)}H(\D)\ot_{H(B)}H(\D)\rarrow H(\D\ot_B\D\ot_B\D)$
are isomorphisms.
 A map of DG\+corings $\C\rarrow\D$ compatible with a morphism of
DG\+algebras $A\rarrow B$ induces a map of graded corings
$H(\C)\rarrow H(\D)$ compatible with the morphism of graded algebras
$H(A)\rarrow H(B)$ whenever both DG\+corings $\C$ and $\D$ satisfy
the above two conditions.
 Here a map $\C\rarrow\D$ from a DG\+coring $\C$ over a DG\+algebra $A$
to a DG\+coring $\D$ over a DG\+algebra $B$ is called compatible with
a morphism of DG\+algebras $A\rarrow B$ over $k$ if the map of graded
corings $\C\rarrow\D$ is compatible with the morphism of graded algebras
$A\rarrow B$ and the map $\C\rarrow\D$ is a morphism of complexes.

\subsubsection{}  \label{quasi-differential-corings}
 Coderivations of a graded coring $\D$ of degree~$-1$ with respect to
coderivations of a graded $k$\+algebra $B$ of degree~$-1$ are defined
in the same way as above.

 Now let $A$ be an ungraded $k$\+algebra.
 A \emph{quasi-differential coring} $\tD$ over $A$ is a graded
coring over~$A$ endowed with a coderivation~$\d$ of degree~$1$
(with respect to the zero derivation of the $k$\+algebra $A$,
which is considered as a graded $k$\+algebra concentrated in
degree~$0$) such that $\d^2=0$ and the cohomology of~$\d$ vanish.
 If $\tD$ is a quasi-differential coring over a $k$\+algebra $A$, then
the cokernel $\tD/\im\d$ of the derivation~$\d$ has a natural structure
of graded coring over~$A$.
 A \emph{quasi-differential structure} on a graded coring $\D$ is 
the data of a quasi-differential coring $\tD$ together with
an isomorphism of graded corings $\tD/\im\d\simeq\D$.
 We will denote the grading of a quasi-differential coring $\tD$ by
lower indices, even though the differential raises the degree.
 This terminology and notation is explained by the following
construction (cf.~\ref{quasi-differential-rings}).

 We will use Sweedler's notation~\cite{Swe} \ $p\mpsto p_{(1)}\ot
p_{(2)}$ for the comultiplication map of a coring $\D$ over $A$;
here $p\in\D$ and $p_{(1)}\ot p_{(2)}\in \D\ot_A\D$.
 A \emph{CDG\+coring} $\D$ over a $k$\+algebra $A$ is a graded coring
over~$A$ endowed with a coderivation~$d$ of degree~$-1$ (with respect
to the zero derivation of~$A$) and an $A$\+$A$\bimodule{} map
$h\:\D_2\rarrow A$ satisfying the equations $d^2(p) = h(p_{(1)})p_{(2)}
- p_{(1)}h(p_{(2)})$ and $h(d(p))=0$ for all $p\in\D$, where
the map~$h$ is considered to be extended by zero to the components
$\D_i$ with $i\ne2$.
 Given a CDG\+coring $\D$ over a $k$\+algebra $A$ and a CDG\+coring
$\E$ over a $k$\+algebra $B$, a morphism of CDG\+corings $\D\rarrow\E$
compatible with a morphism of $k$\+algebras $A\rarrow B$ is a pair
$(g,a)$, where $g\:\D\rarrow\E$ is a map of graded corings compatible
with a morphism of $k$\+algebras $A\rarrow B$ and $a\:\D_1\rarrow B$
is an $A$\+$A$\bimodule{} map satisfying the equations $d(g(p)) =
g(d(p)) + a(p_{(1)})g(p_{(2)}) + (-1)^{|p|}g(p_{(1)})a(p_{(2)})$ and
$h(g(q)) = h(q) + a(d(q)) + a(q_{(1)})a(q_{(2)})$ hold for all $p\in\D$
and $q\in\D_2$ (where the map $a$ is extended by zero to
the components $\D_i$ with $i\ne1$).

 Composition of morphisms of CDG\+corings is defined by the rule
$(g',a')(g'',a'') = (g'g''\;a'g''+a'')$; identity morphisms are
the morphisms~$(\id,0)$.
 So the category of CDG\+corings is defined.
 Notice that two CDG\+corings of the form $(\D,d',h')$ and
$(\D,d'',h'')$ over a $k$\+algebra $A$ with $d''(p) = d'(p) +
a(p_{(1)})p_{(2)} + (-1)^{|p|}p_{(1)}a(p_{(2)})$ and $h''(q) = h'(q) +
a(d'(q)) + a(q_{(1)})a(q_{(2)})$, where $a\:\D_1\rarrow A$ is
an $A$\+$A$\bimodule{} map, are always naturally isomorphic to each
other, the isomorphism being given by $(\id,a)\:(\D,d',h')\rarrow
(\D,d'',h'')$.

 The category of DG\+corings (over ungraded $k$\+algebras considered
as DG\+algebras concentrated in degree zero) has DG\+corings $\D$ over
$k$\+algebras $A$ as objects and maps of DG\+corings $\D\rarrow\E$
compatible with morphisms of $k$\+algebras $A\rarrow B$ as morphisms.
 The category of quasi-differential corings can be defined as
the full subcategory of the category of DG\+corings whose objects
are the DG\+corings with acyclic differentials.
 One can also consider the category of DG\+corings (over ungraded
$k$\+algebras) with coderivations of degree~$-1$.
 There is an obvious faithful, but not fully faithful functor from
the latter category to the category of CDG\+corings, assigning
the CDG\+coring $(\D,d,h)$ with $h=0$ to a DG\+coring $(\D,d)$ and
the morphism of CDG\+corings $(g,0)$ to a map of
DG\+corings~$g\:\D\rarrow\E$ compatible with a morphism of
$k$\+algebras $A\rarrow B$.

 There is a natural fully faithful functor from the category of
CDG\+corings to the category of quasi-differential corings, whose
image consists of the quasi-differential corings $\tD$ over $A$
for which the counit map $\D_0\til\rarrow A$ can be presented as
the composition of the coderivation component $\d_0\:\D_0\til\rarrow
\D_1\til$ and some $A$\+$A$\bimodule{} map $\delta\:\D_1\til\rarrow A$.
 In other words, a quasi-differential coring comes from a CDG\+coring
if and only if the counit map $\D_0\til/\d_{-1}\D_{-1}\til\rarrow A$
can be extended to an $A$\+$A$\bimodule{} map $\D_1\til\rarrow A$,
where the $A$\+$A$\+bimodule{} $\D_0\til/\d_{-1}\D_{-1}\til$ is
embedded into the $A$\+$A$\+bimodule{} $\D_1\til$ by the map~$\d_0$.
 In particular, the categories of quasi-differential corings and
CDG\+corings over a field $A=k$ (\emph{quasi-differential
coalgebras} and \emph{CDG\+coalgebras} over~$k$) are naturally
equivalent.

 Let us first construct the inverse functor.
 Given a quasi-differential coring $(\tD,\d)$ and a map $\delta\:
\D_1\til\rarrow A$ as above, set $\D=\tD/\im\d$ and define $d$ and $h$
by the formulas $d(\overline p)=\delta(p_{(1)})\overline{p_{(2)}\!\.}\.
+ (-1)^{|p|}\overline{p_{(1)}\!\.}\.\delta(p_{(2)})$ and
$h(\overline q) = \delta(q_{(1)})\delta(q_{(2)})$ for $p\in\tD$ and
$q\in\D_2\til$, where the map $\delta$ is extended by zero to
the components $\D_i\til$ with $i\ne 1$ and $\overline r\in\D$ denotes
the image of an element $r\in\tD$.
 To a map of quasi-differential corings $g\:\tD\rarrow\E\til$
endowed with maps $\delta_\D\:\D_1\til\rarrow A$ and $\delta_\E\:
\E_1\til\rarrow B$ with the above property, compatible with a morphism
of $k$\+algebras $f\:A\rarrow B$, one assigns the morphism of
CDG\+corings $(\overline g\;\delta_\E g-f\delta_\D)$, where
$\overline g\:\D\rarrow\E$ denotes the induced morphism on
the cokernels of the coderivations~$\d$.

 Conversely, to a CDG\+coring $(\D,d,h)$ over a $k$\+algebra $A$
one assigns the quasi-differential coring $(\tD,\d)$ over $A$ whose
graded components are the $A$\+$A$\bimodule s $\D_i\til=
\D_i\oplus\D_{i-1}$, the coderivation~$\d$ is given by the formula
$\d(\tau p+\dbar q)=\dbar p$, and the comultiplication is given by
the formula $\tau p+\dbar q\mpsto \tau p_{(1)}\ot \tau p_{(2)} +
(-1)^{|p_{(1)}|}\tau d(p_{(1)})\ot\dbar p_{(2)} + (-1)^{|p_{(2)}|}
h(p_{(1)})\dbar p_{(2)}\ot \dbar p_{(3)} + \dbar q_{(1)}\ot\tau q_{(2)}
+ (-1)^{|q_{(1)}|}\tau q_{(1)}\ot\dbar q_{(2)} + (-1)^{|q_{(1)}|}
\dbar d(q_{(1)})\ot \dbar q_{(2)}$, where $\tau p+\dbar q=(p,q)$ is
a formal notation for an element of $\bigoplus_i(\D_i\oplus\D_{i-1})$.
 To a morphism of CDG\+corings $(g,a)\:\D\rarrow\E$, the morphism
of quasi-differential corings $\bigoplus_i(\D_i\oplus\D_{i-1})
\rarrow\bigoplus_i(\E_i\oplus\E_{i-1})$ given by the formula
$\tau p+\dbar q\mpsto \tau g(p)+a(p_{(1)})\dbar g(p_{(2)})+\dbar g(q)$
is assigned.
 For a quasi-differential coring $(\tD,\d)$ over a $k$\+algebra $A$
endowed with a map $\delta\:\D_1\til\rarrow A$ with the above property
and the corresponding CDG\+coring $(\D,d,h)$, the natural morphism
of quasi-differential corings $\tD\rarrow\bigoplus_i
(\D_i\oplus\D_{i-1})$ over $A$ is given by the formula $p\mpsto
\tau\overline p + \delta(p_{(1)})\dbar\overline{p_{(2)}\!\.}\.$ for
$p\in\tD$.
 This morphism is an isomorphism, since the induced morphism of
the cokernels of the coderivations~$\d$ is an isomorphism.

\subsubsection{}  \label{diff-semialgebras-ii}
 Let $B$ be a graded $k$\+algebra endowed with a derivation $d_B$
of degree~$1$ and $\D$ be a graded coring over $B$ endowed with
a coderivation $\d_\D$ with respect to~$d_B$.
 Let $\T$ be a graded semialgebra over $\D$.
 A homogeneous map $d_\T\:\T\rarrow\T$ of degree~$1$ is called 
a \emph{semiderivation of\/ $\T$ with respect to\/ $d_\D$ and\/ $d_B$}
if the biaction map $B\ot_k\T\ot_k B\rarrow \T$, the bicoaction map
$\T\rarrow\D\ot_B\T\ot_B\D$, and the semimultiplication map
$\T\oc_\D\T\rarrow\T$ are morphisms in the category of graded
$k$\module s endowed with endomorphisms of degree~$1$.
 In this case, it follows that the semiunit map $\D\rarrow\T$
satisfies the same condition.
 In the particular case when $B$ and $\D$ are concentrated in
degree~$0$ and $d_B=0=d_\D$, the conditions on the biaction and
bicoaction map simply mean that $d_\T$ is a $\D$\+$\D$\bicomodule{}
morphism.

 Let $B$ be a DG\+algebra over~$k$ and $\D$ be a DG\+coring over $B$.
 A \emph{DG\+semialgebra} over $\D$ is a graded semialgebra over
the graded coring $\D$ endowed with a semiderivation
$d_\T$ with respect to $d_\D$ and $d_B$ such that $d_\T^2=0$.

 Let $\T$ be a DG\semialgebra{} over a DG\+coring $\D$.
 Then the cohomology $H(\T)$ is endowed with a natural structure
of graded semialgebra over the graded coring $H(\D)$ provided that
(i)~the natural maps from the tensor products of cohomology to
the cohomology of the tensor products are isomorphisms for
the tensor products $\D\ot_B\D$, \ $\D\ot_B\D\ot_B\D$, \ $\D\ot_B\T$,
\ $\T\ot_B\D$, \ $\D\ot_B\D\ot_B\T$, \ $\T\ot_B\D\ot_B\D$, \
$\D\ot_B\T\ot_B\D$, \ $\T\ot_B\T$, \ $\D\ot_B\T\ot_B\T$, \ 
$\T\ot_B\T\ot_B\D$, \ $\T\ot_B\D\ot_B\T$, \ $\D\ot_B\T\ot_B\D\ot_B\T$,
\ $\T\ot_B\D\ot_B\T\ot_B\D$, \ $\T\ot_B\T\ot_B\T$, \
$\T\ot_B\D\ot_B\T\ot_B\T$, \ $\T\ot_B\T\ot_B\D\ot_B\T$;
(ii)~the multiple cotensor products
$H(\T)\oc_{H(\D)}\dsb\oc_{H(\D)}H(\T)$ are associative, where the graded
$H(\D)$\+$H(\D)$\bicomodule{} structure on $H(\T)$ is well-defined in
view of~(i); and
(iii)~the natural maps $H(\T\oc_\D\T)\rarrow H(\T)\oc_{H(\D)} H(\T)$, \
$H(\D\ot_B\T\oc_\D\T)\rarrow H(\D)\ot_{H(B)}H(\T)\oc_{H(\D)}H(\T)$, \
$H(\T\oc_\D\T\ot_B\D)\rarrow H(\T)\oc_{H(\D)}H(\T)\ot_{H(B)}H(\D)$, and
$H(\T\oc_\D\T\oc_\D\T)\rarrow H(\T)\oc_{H(\D)}H(\T)\oc_{H(\D)}H(\T)$,
which are well-defined in view of (i) and~(ii), are isomorphisms.

 A map of DG\semialgebra s $\S\rarrow \T$ compatible with a map
of DG\+corings $\C\rarrow\D$ and a morphism of DG\+algebras $A\rarrow B$
induces a map of graded semialgebras $H(\S)\rarrow H(\T)$ compatible
with the map of graded corings $H(\C)\rarrow H(\D)$ and the morphism of
graded $k$\+algebras $H(A)\rarrow H(B)$ whenever both DG\semialgebra s
$\S$ and $\T$ satisfy the above three conditions.
 Here a map $\S\rarrow\T$ from a DG\semialgebra{} $\S$ over a DG\+coring
$\C$ to a DG\semialgebra{} $\T$ over a DG\+coring $\D$ is called
compatible with a map of DG\+corings $\C\rarrow\D$ and a morphism of
DG\+algebras $A\rarrow B$ if the map of graded semialgebras
$\S\rarrow\T$ is compatible with the map of graded corings $\C\rarrow\D$
and the morphism of graded $k$\+algebras $A\rarrow B$, and the maps
$\S\rarrow\T$ and $\C\rarrow\D$ are morphisms of complexes.

\subsection{One-sided SemiTor}   \label{one-sided-semitor}
 Let $\S$ be a semialgebra over a coring $\C$ over a $k$\+algebra~$A$.
 We will consider two situations separately.

\subsubsection{}
 Assume that $\C$ is a flat right $A$\module{} and $\S$ is a coflat
right $\C$\comodule.

 Consider the functor of semitensor product over $\S$ on the Carthesian
product of the homotopy category of complexes of $\C$\+coflat right
$\S$\semimodule s and the homotopy category of complexes of left
$\S$\semimodule s.
 The semiderived category of $\C$\+coflat right $\S$\semimodule s
is defined as the quotient category of the homotopy category of
$\C$\+coflat right $\S$\semimodule s by the thick subcategory of
complexes of right $\S$\semimodule s that as complexes of
$\C$\comodule s are coacyclic with respect to the exact category
of coflat right $\C$\comodule s.
 A complex of left $\S$\semimodule s $\bM^\bu$ is called \emph{semiflat
relative to\/ $\C$} if the complex $\bN^\bu\os_\S\bM^\bu$ is acyclic for
any $\C$\contractible{} complex of $\C$\+coflat right $\S$\semimodule s
$\bN^\bu$ (cf.~\ref{relatively-semiflat}).

 The left derived functor $\SemiTor^\S$ on the Carthesian product of
the semiderived category of $\C$\+coflat right $\S$\semimodule s and
the semiderived category of left $\S$\semimodule s is defined by
restricting the functor of semitensor product to the Carthesian product
of the homotopy category of $\C$\+coflat right $\S$\semimodule s and
the homotopy category of complexes of left $\S$\semimodule s semiflat
relative to $\C$, or to the Carthesian product of the homotopy category
of semiflat complexes of right $\S$\semimodule s and the homotopy
category of left $\S$\semimodule s.
 This definition of a derived functor is a particular case of
\emph{both} Lemmas~\ref{semitor-definition}
and~\ref{semi-ctrtor-definition}.2.
 If $\bN^\bu$ is a complex of $\C$\+coflat right $\S$\semimodule s and
$\bM^\bu$ is a complex of left $\S$\semimodule s, then the total complex
of the bar bicomplex $\dsb\rarrow\bN^\bu\oc_\C\S\oc_\C\S\oc_\C\bM^\bu
\rarrow\bN^\bu\oc_\C\S\oc_\C\bM^\bu\rarrow\bN^\bu\oc_\C\bM^\bu$,
constructed by taking infinite direct sums along the diagonals,
represents the object $\SemiTor^\S(\bM^\bu,\bN^\bu)$ in $\sD(k\modl)$.
 When the semiunit map $\C\rarrow\S$ is injective and its cokernel is
a flat right $A$\module{} (and hence a coflat right $\C$\comodule{}
by Lemma~\ref{absolute-relative-coflat}), one can also
use the reduced bar bicomplex $\dsb\rarrow
\bN^\bu\oc_\C\S/\C\oc_\C\S/\C\oc_\C\bM^\bu\rarrow
\bN^\bu\oc_\C\S/\C\oc_\C\bM^\bu\rarrow\bN^\bu\oc_\C\bM^\bu$.

 In the case when $\S$ is a graded semialgebra one analogously defines
the derived functor $\SemiTor^\S_\gr$ acting from the Carthesian
product of the semiderived category of $\C$\+coflat graded right
$\S$\semimodule s and the semiderived category of graded left
$\S$\semimodule s to the derived category of graded $k$\module s.

\subsubsection{}
 Assume that $\C$ is a flat right $A$\module, $\S$ is a flat right
$A$\module{} and a $\C/A$\+coflat left $\C$\comodule, and the ring $A$
has a finite weak homological dimension.

 Consider the functor of semitensor product over $\S$ on the Carthesian
product of the homotopy category of complexes of $A$\+flat right
$\S$\semimodule s and the homotopy category of complexes of
$\C/A$\+coflat left $\S$\semimodule s.
 The semiderived category of $A$\+flat right $\S$\semimodule s
($\C/A$\+coflat left $\S$\semimodule s) is defined as the quotient
category of the homotopy category of $A$\+flat right $\S$\semimodule s
($\C/A$\+coflat left $\S$\semimodule s) by the thick subcategory of
complexes of $\S$\semimodule s that as complexes of $\C$\comodule s are
coacyclic with respect to the exact category of $A$\+flat right
$\C$\comodule s ($\C/A$\+coflat left $\C$\comodule s).
 A complex of $\C/A$\+coflat left $\S$\semimodule s $\bM^\bu$ is called
\emph{semiflat relative to $A$} if the complex of $k$\module s
$\bN^\bu\os_\S\bM^\bu$ is acyclic for any complex of right
$\S$\semimodule s $\bN^\bu$ that as a complex of right $\C$\comodule s
is coacyclic with respect to the exact category of $A$\+flat right
$\C$\comodule s.
 A complex of $A$\+flat right $\S$\semimodule s $\bN^\bu$ is called
\emph{$\S/\C/A$\semiflat} if the complex of $k$\module s
$\bN^\bu\os_\S\bM^\bu$ is acyclic for any $\C$\contractible{} complex of
$\C/A$\+coflat left $\S$\semimodule s $\bM^\bu$
(cf.~\ref{relatively-semiflat}).

 The left derived functor $\SemiTor^\S$ on the Carthesian product of
the semiderived category of $A$\+flat right $\S$\semimodule s and
the semiderived category of $\C/A$\+coflat left $\S$\semimodule s is
defined by restricting the functor of semitensor product to
the Carthesian product of the homotopy category of $A$\+flat right
$\S$\semimodule s and the homotopy category of complexes of 
$\C/A$\+coflat left $\S$\semimodule s semiflat relative to $A$, or
to the Carthesian product of the homotopy category of
$\S/\C/A$\semiflat{} complexes of $A$\+flat right $\S$\semimodule s
and the homotopy category of $\C/A$\+coflat left $\S$\semimodule s.
 This definition of a derived functor is a particular case of
\emph{both} Lemmas~\ref{semitor-definition}
and~\ref{semi-ctrtor-definition}.2.
 If $\bN^\bu$ is a complex of $A$\+flat right $\S$\semimodule s and
$\bM^\bu$ is a complex of $\C/A$\+coflat left $\S$\semimodule s, then
the total complex of the bar bicomplex $\dsb\rarrow
\bN^\bu\oc_\C\S\oc_\C\S\oc_\C\bM^\bu\rarrow\bN^\bu\oc_\C\S\oc_\C\bM^\bu
\rarrow\bN^\bu\oc_\C\bM^\bu$, constructed by taking infinite direct sums
along the diagonals, represents the object
$\SemiTor^\S(\bM^\bu,\bN^\bu)$ in $\sD(k\modl)$.
 When the semiunit map $\C\rarrow\S$ is injective and its cokernel is
a flat right $A$\module{} (the cokernel is a $\C/A$\+coflat left
$\C$\comodule{} by Lemma~\ref{absolute-relative-coflat}), one can
also use the reduced bar bicomplex $\dsb\rarrow
\bN^\bu\oc_\C\S/\C\oc_\C\S/\C\oc_\C\bM^\bu\rarrow
\bN^\bu\oc_\C\S/\C\oc_\C\bM^\bu\rarrow\bN^\bu\oc_\C\bM^\bu$.

 In the case when $\S$ is a graded semialgebra one analogously defines
the derived functor $\SemiTor^\S_\gr$ acting from the Carthesian
product of the semiderived category of $A$\+flat graded right
$\S$\semimodule s and the semiderived category of $\C/A$\+coflat graded
left $\S$\semimodule s to the derived category of graded $k$\module s.

\subsection{Koszul semialgebras and corings}  \label{koszul-semi-co}

\subsubsection{}
 Let $\S$ be a semialgebra over a coring $\C$ over a $k$\+algebra $A$.
 Suppose that $\S$ is endowed with an augmentation, i.~e., a morphism
$\S\rarrow\C$ of semialgebras over~$\C$; let $\S_+$ be the kernel of
this map.
 We will denote by $\Br^\bu(\S,\C)$ the reduced bar complex $\dsb\rarrow
\S_+\oc_\C\S_+\oc_C\S_+\rarrow\S_+\oc_\C\S_+\rarrow\S_+\rarrow\C$.
 It can be also defined as the coring $\bigoplus_{n=0}^\infty\S_+^
{\suboc_\C n}$ over the $k$\+algebra $A$ (the ``cotensor coring'' of
the $\C$\+$\C$\bicomodule{} $\S$)  endowed with the unique grading
such that the component $\S_+$ is situated in degree~$-1$ and
the unique coderivation (with respect to the zero derivation of~$A$) of
degree~$1$ whose component mapping $\S_+\oc_\C\S_+$ to $\S_+$ is equal
to the semimultiplication map $\S_+\oc_\C\S_+\rarrow\S_+$.
 So $\Br^\bu(\S,\C)$ is a DG\+coring over the $k$\+algebra $A$
considered as a DG\+algebra concentrated in degree~$0$.

 Now let $\S$ be a graded semialgebra over a coring $\C$ over
a $k$\+algebra $A$, where $A$ and $\C$ are considered as a graded
$k$\+algebra and a graded coring concentrated in degree~$0$; assume
additionally that $\S$ is concentrated in nonnegative degrees, $\C$
is the component of degree~$0$ in $\S$, and the augmentation map
$\S\rarrow\C$ is simply the projection of $\S$ to its component
of degree~$0$.
 In this case there is a graded version $\Br^\bu_\gr(\S,\C)$ of
the above bar complex, which is a bigraded object with the grading
denoted by upper indices coming from the cotensor powers of $\S_+$
and the grading denoted by lower indices coming from the grading
of $\S_+$ itself.
 Notice that the component $\Br^i_n(\S,\C)$ can be only nonzero
when $0\le -i\le n$.

 Let $\C$ and $\D$ be corings over a $k$\+algebra $A$.
 Suppose that we are given two maps $\C\rarrow\D$ and $\D\rarrow\C$
that are morphisms of corings over $A$ such that the composition
$\C\rarrow\D\rarrow\C$ is the identity; let $\D_+$ be the cokernel of
the map $\C\rarrow\D$.
 Assume that the multiple cotensor products $\D\oc_\C\dsb\oc_\C\D$,
where $\D$ is endowed with a $\C$\+$\C$\bicomodule{} structure via
the morphism $\D\rarrow\C$, are associative.
 We will denote by $\Cb^\bu(\D,\C)$ the reduced cobar complex $\C\rarrow
\D_+\rarrow\D_+\oc_\C\D_+\rarrow\D_+\oc_\C\D_+\oc_\C\D_+\rarrow\dsb$\.
 It can be also defined as the semialgebra $\bigoplus_{n=0}^\infty
\D_+^{\suboc_\C n}$ over the coring $\C$ (the ``cotensor semialgebra''
of the $\C$\+$\C$\bicomodule{} $\D$) endowed with the unique grading
such that the component $\D_+$ is situated in degree~$1$ and
the unique semiderivation (with respect to $d_\C=0$ and $d_A=0$) of
degree~$1$ whose component mapping $\D_+$ to $\D_+\oc_\C\D_+$ is equal
to the comultiplication map $\D_+\rarrow\D_+\oc_\C\D_+$.
 So $\Cb^\bu(\D,\C)$ is a DG\+semialgebra over the coring $\C$ over
the $k$\+algebra $A$, where $A$ and $\C$ are considered as a DG\+algebra
and a DG\+coring concentrated in degree~$0$.

 Now let $\D$ be a graded coring over a $k$\+algebra $A$ considered as
a graded $k$\+algebra concentrated in degree~$0$ and $\C$ be a coring
over $A$; assume additionally that $\D$ is concentrated in nonnegative
degrees, $\C$ is the component of degree~$0$ in $\D$, and the maps
$\C\rarrow\D$ and $\D\rarrow\C$ are simply the embedding of and
the projection to the component of degree~$0$.
 In this case there is a graded version $\Cb^\bu_\gr(\D,\C)$ of
the above cobar complex, which is a bigraded object with the grading
denoted by upper indices coming from the cotensor powers of $\D_+$
and the grading denoted by lower indices coming from the grading
of $\D_+$ itself.
 Notice that the component $\Cb^i_n(\D,\C)$ can be only nonzero when
$0\le i\le n$.

\subsubsection{}
 Let $\C$ be a coring over a $k$\+algebra $A$.
 Assume that $\C$ is a flat right $A$\module.

 A graded semialgebra $\S$ over $\C$ is called \emph{right coflat
Koszul} if (i)~$\S$ is nonnegatively graded and the semiunit
homomorphism is an isomorphism $\C\simeq \S_0$; (ii)~the components
$\S_i$ are flat right $A$\module s; (iii)~the cohomology
$H^i_n\Br^\bu_\gr(\S,\C)$ are only nonzero on the diagonal $-i=n$;
and (iv)~whenever the component $\Br_n^\bu(\S,\C)$ is a complex of
$A$\+flat right $\C$\comodule s, so the diagonal cohomology
$H_n^{-n}\Br_\gr^\bu(\S,\C)$ can be endowed with a right
$\C$\comodule{} structure as the kernel of a morphism in the category
of right $\C$\comodule s, it is a coflat right $\C$\comodule.

 When the ring $A$ has a finite weak homological dimension, there is
an analogous definition of a \emph{right flat and left relatively
coflat Koszul} semialgebra $\S$ over~$\C$.
 One imposes the same conditions (i-iii) and replaces (iv) with 
the condition (iv$'$)~the diagonal cohomology $H_n^{-n}\Br_\gr^\bu
(\S,\C)$ is a $\C/A$\+coflat left $\C$\comodule{} for all~$n$.

 A graded coring $\D$ over the $k$\+algebra $A$ endowed with
a morphism $\D\rarrow\C$ of corings over~$A$ is called a \emph{right
coflat Koszul coring over~$\C$} if (i)~$\D$ is nonnegatively graded
and the morphism $\D\rarrow\C$ vanishes on the components of positive
degree in $\D$ and induces an isomorphism $\D_0\simeq\C$; (ii)~whenever
a component $\D_n$ is a flat right $A$\module, it is a coflat right
$\C$\comodule; (iii)~whenever all the multiple cotensor products
entering into the construction of the component $\Cb_n^\bu(\D,\C)$
are associative, so this component is well-defined, the cohomology
$H^i\Cb^\bu_n(\D,\C)$ is only nonzero on the diagonal $i=n$; and
(iv)~in the assumptions of~(iii), the diagonal cohomology
$H^n\Cb^\bu_n(\D,\C)$ is a flat right $A$\module.

 When the ring $A$ has a finite weak homological dimension, there is
an analogous definition of a \emph{right flat and left relatively
coflat Koszul coring} $\D$ over~$\C$.
 One imposes the same conditions (i-ii), (iv), and replaces (iii) with
the condition (iii$'$)~the component $\D_n$ is a $\C/A$\+coflat left
$\C$\comodule{} for all~$n$.

\subsubsection{}  \label{homogeneous-quadratic-duality}
 The objects of the category of right coflat Koszul semialgebras are
right coflat Koszul semialgebras $\S$ over corings $\C$ over
$k$\+algebras $A$ such that $\C$ is a flat right $A$\module.
 Morphisms are maps of graded semialgebras $\S\rarrow\S'$ compatible
with maps of corings $\C\rarrow\C'$ and morphisms of $k$\+algebras
$A\rarrow A'$.
 Imposing the additional assumption that $A$ has a finite weak
homological dimension, one analogously defines the category of right
flat and left relatively coflat Koszul semialgebras.

 The objects of the category of right coflat Koszul corings are right
coflat Koszul corings $\D$ over corings $\C$ over $k$\+algebras $A$
such that $\C$ is a flat right $A$\module.
 Morphisms are maps of graded corings $\D\rarrow\D'$ compatible with
morphisms of $k$\+algebras $A\rarrow A'$.
 Imposing the additional assumption that $A$ has a finite weak
homological dimension, one analogously defines the category of
right flat and left relatively Koszul corings.

\begin{thm}
 The category of right coflat Koszul semialgebras is equivalent to
the category of right coflat Koszul corings.
 Analogously, the category of right flat and left relatively coflat
Koszul semialgebras is equivalent to the category of right flat and
left relatively coflat Koszul corings.
 In both cases, the mutually inverse equivalences are provided by
the functor assigning to a Koszul semialgebra $\S$ the coring of
cohomology of the graded DG\+coring $\Br^\bu_\gr(\S,\C)$ and the functor
assigning to a Koszul coring $\D$ the semialgebra of cohomology of
the graded DG\semialgebra{} $\Cb^\bu_\gr(\D,\C)$.
\end{thm}

\begin{proof}
 The assertions of Theorem follow from Propositions~1 and~2 below.
 To check the conditions of~\ref{diff-semialgebras} needed for
the coring of cohomology and the semialgebra of cohomology to be
defined, use Lemma~\ref{absolute-relative-coflat} and
Proposition~\ref{cotensor-associative}.
\end{proof}

 Let $\C$ be a coring over a $k$\+algebra $A$.

\begin{prop1}
 \textup{(a)} Assume that\/ $\C$ is a flat right $A$\module.
 Then a graded semialgebra\/ $\S$ over\/ $\C$ is right coflat Koszul
if and only if (i)\/~$\S$ is nonnegatively graded and the semiunit map
is an isomorphism\/ $\C\simeq\S_0$; (ii)~for any $n\ge1$ the natural
map from the quotient $k$\module{} of the cotensor power\/
$\S_1^{\suboc_\C n}$ by the sum of the kernels of its maps to cotensor
products\/ $\S_1^{\suboc_\C i-1}\oc_\C\S_2\oc_\C\S_1^{\suboc_\C n-i-1}$,
\ $i=1,\dsc,n-1$ to the component\/ $\S_n$ is an isomorphism;
(iii)~the lattice of submodules of the $k$\module\/ $\S_i^{\suboc_\C n}$
generated by these $n-1$ kernels is distributive; (iv)~all the quotient
modules of embedded submodules belonging to the mentioned lattice are
flat right $A$\module s in their natural right $A$\module{} structures;
and (v)~all the quotient modules of embedded submodules belonging to
this lattice are coflat right\/ $\C$\comodule s in their right\/
$\C$\comodule{} structures that are well-defined in view of~(iv). \par
 \textup{(b)} Assume that\/ $\C$ is a flat right $A$\module{} and
$A$ has a finite weak homological dimension.
 Then a graded semialgebra\/ $\S$ over\/ $\C$ is right flat and left
relatively coflat Koszul if and only if it satisfies the conditions
(i-iv) of~(a) and the condition (v\/$'$)~all the quotient modules of
embedded submodules belonging to the lattice under consideration are\/
$\C/A$\+coflat left\/ $\C$\comodule s in their natural left\/
$\C$\comodule{} structures.
\end{prop1}

\begin{prop2}
 \textup{(a)} Assume that\/ $\C$ is a flat right $A$\module.
 Then a graded coring\/ $\D$ endowed with a morphism\/ $\D\rarrow\C$ of
corings over $A$ is a right coflat Koszul coring over\/ $\C$ if and
only if (i)\/~$\D$ is nonpositively graded and the morphism\/
$\D\rarrow\C$ vanishes on the components of positive degrees in\/ $\D$
and induces an isomorphism\/ $\D_0\simeq\C$; (ii)~for any $n\ge1$
the natural map from the component\/ $\D_n$ to the intersection of
images of the maps from cotensor products\/
$\D_1^{\suboc_\C i-1}\oc_\C\D_2\oc_\C\D_1^{\suboc_\C n-i-1}$, \
$i=1,\dsc,n-1$ to the cotensor power\/ $\D_1^{\suboc_\C n}$ is
an isomorphism;
(iii)~the lattice of submodules of the $k$\module{} $\D_1^{\suboc_\C n}$
generated by these $n-1$ images is distributive; (iv) all the quotient
modules of the embedded submodules belonging to the mentioned lattice
are flat right $A$\module s in their natural right $A$\module{}
structures; and (v)~all the quotient modules of embedded submodules
belonging to this lattice are coflat right $\C$\comodule s in their
right $\C$\comodule{} structures that are well-defined in view of~(iv).
\par
 \textup{(b)} Assume that\/ $\C$ is a flat right $A$\module{} and
$A$ has a finite weak homological dimension.
 Then a graded coring\/ $\D$ endowed with a morphism\/ $\D\rarrow\C$ of
corings over $A$ is a right flat and left relatively coflat Koszul
coring over\/ $\C$ if and only if it satisfies the conditions (i-iv)
of~(a) and the condition (v\/$'$)~all the quotient modules of embedded
submodules belonging to the lattice under consideration are\/
$\C/A$\+coflat left\/ $\C$\comodule s in their natural left\/
$\C$\comodule{} structures.
\end{prop2}

\begin{proof}[Proof of Propositions 1 and 2]
 Both Propositions follow by induction in the internal degree~$n$ from
Lemma~\ref{absolute-relative-coflat},
Proposition~\ref{cotensor-associative}, and the next Lemma~1 
(parts (a)$\Longleftrightarrow$(c), (a)$\Longleftrightarrow$(c*)),
and the final assertion) and Lemma~2.
\end{proof}

\begin{lem1}
 Let $W$ be a $k$\module{} and $X_1,\dsc,X_{n-1}\subset W$ be
a collection of submodules such that any proper subset $X_1,\dsc,\hX_k,
\dsc,X_{n-1}$ generates a distributive lattice of submodules in $W$.
 Then the following conditions are equivalent:
  \begin{itemize}
    \item[(a)] the collection of submodules $X_1,\dsc,X_{n-1}$ generates
               a distributive lattice of submodules in $W$;
    \item[(b)] the following complex of $k$\module s
               $K_\bu(W;X_1,\dsc,X_{n-1})$ is exact
      \begin{multline*}
      0\rarrow X_1\cap\dsb\cap X_{n-1}\rarrow X_2\cap\dsb\cap X_{n-1}
      \rarrow X_3\cap\dsb\cap X_{n-1}/X_1\rarrow \\
      \dsb\rarrow \bigcap_{s=i+1}^{n-1} X_s \big/
      \sum\nolimits_{t=1}^{i-1} X_t\rarrow\dsb\rarrow\\
      X_{n-1}/(X_1+\dsb+X_{n-3})\rarrow W/(X_1+\dsb+X_{n-2})
      \rarrow W/(X_1+\dsb+X_{n-1})\rarrow 0,
    \end{multline*}
               where we denote $Y/Z=Y/Y\cap Z$;
    \item[(c)] the following complex of $k$\module s
               $B_\bu(W;X_1,\dsc,X_{n-1})$
      \begin{multline*}
       W\rarrow \bigoplus\nolimits_t W/X_t\rarrow\dsb\rarrow
       \bigoplus\nolimits_{t_1<\dsb<t_{n-i}}
       W/{\textstyle\sum_{s=1}^{n-i} X_{t_s}}\rarrow\\
       \dsb\rarrow W/{\textstyle\sum_s X_s}\rarrow 0
     \end{multline*}
               is exact everywhere except for the leftmost term;
   \item[(c*)] the following complex of $k$\module s
               $B^\bu(W;X_1,\dsc,X_{n-1})$
       $$
        0\rarrow {\textstyle\bigcap_s X_s}
        \rarrow\dsb\rarrow\bigoplus\nolimits_{t_1<\dsb<t_{n-i}}
        {\textstyle\bigcap_{s=1}^{n-i} X_{t_s}}\rarrow\dsb\rarrow
        \bigoplus\nolimits_t X_t\rarrow W
       $$
               is exact everywhere except for the rightmost term.
    \end{itemize}
 Besides, the complex in~(c) is always exact at its two rightmost
nontrivial terms, and the complex in~(c*) is always exact at its
two leftmost nontrivial terms.
\end{lem1}

\begin{proof}
 See the proof of~\cite[Proposition~7.2 of Chapter~1]{PP}.
\end{proof}

 Assume that the coring $\C$ is a flat right $A$\module.

\begin{lem2}
 Let $W$ be a $k$\module{} and $X_1,\dsc,X_{n-1}\subset W$ be
a collection of submodules generating a distributive lattice of
submodules in~$W$. \par
 \textup{(a)} Suppose that $W$ is a right $A$\module{} and $X_s$ are
its $A$\+submodules.
 Then all the subquotient modules in the lattice of submodules
generated by $X_s$ are flat right $A$\module s if and only if for any
$1\le t_1<\dsb<t_{m-1}\le n-1$ the quotient module $W/(X_{t_1}+\dsb
+X_{t_{m-1}})$ is a flat right $A$\module. \par
 \textup{(b)} Suppose that $W$ is a left\/ $\C$\comodule{} and $X_s$
are its\/ $\C$\subcomodule s.
 Then all the subquotient modules in the lattice of submodules
generated by $X_s$ are\/ $\C/A$\+coflat left $\C$\comodule s if and
only if for any $1\le t_1<\dsb<t_{m-1}\le n-1$ the submodule
$X_{t_1}\cap\dsb\cap X_{t_{m-1}}$ is a\/ $\C/A$\+coflat left\/
$\C$\comodule. \par
 \textup{(c)} Suppose that $W$ is a right $\C$\comodule{} and $X_s$ are
its\/ $\C$\subcomodule s such that all the subquotient modules in
the lattice of submodules generated by $X_s$ are flat right
$A$\module s.
 Then all these subquotient modules are coflat right\/
$\C$\comodule s if and only if for any $1\le t_1<\dsb<t_{m-1}\le n-1$
the submodule $X_{t_1}\cap\dsb\cap X_{t_{m-1}}$ is a coflat right\/
$\C$\comodule.
\end{lem2}

\begin{proof}
 Part~(a): proceed by induction in $n$.
 Since the lattice is distributive, any subquotient module can be
presented as an iterated extenion of subquotient modules of the form
$\bigcap_{j\in J} X_j/\bigcap_{j\in J} X_j\cap \sum_{i\notin J} X_i$,
where $J\subset\{1,\dsc,n-1\}$.
 Whenever the inclusion $J\subset\{1,\dsc,n-1\}$ is proper, this
subquotient module can be presented as an element of a smaller
lattice generated by the submodules $X_j/X_j\cap\sum_{i\notin J}X_i$
in the quotient module $W/\sum_{i\notin J} X_i$.
 It follows from the induction hypothesis that all the submodules
belonging to this smaller lattice are flat right $A$\module s.
 It remains to show that the submodule $X_1\cap\dsb\cap X_{n-1}$
is a flat right $A$\module.
 But the latter submodule is the only nonzero cohomology module
at the leftmost term of the complex of flat right $A$\module s
$B_\bu(W;X_1,\dsc,X_{n-1})$ from Lemma~1(c).
 The proofs of parts~(b) and~(c) are completely analogous, except
for the use of Lemma~\ref{absolute-relative-coflat}.
\end{proof}

 A right coflat (right flat and left relatively coflat) Koszul
semialgebra and a right coflat (right flat and left relatively
coflat) Koszul coring corresponding to each other under
the equivalence of categories from the above Theorem are called
\emph{quadratic dual} to each other.

\subsubsection{}
 Let $\S$ be a right coflat (right flat and left relatively coflat)
Koszul semialgebra over a coring $\C$ and $\D$ be the right coflat
(right flat and left relatively coflat) Koszul coring over $\C$
quadratic dual to $\S$.
 Then on the cotensor products $\S\oc_\C\D$ and $\D\oc_\C\S$ there are
structures of graded complexes whose differentials are the compositions
$\S_i\oc_\C\D_j\to \S_i\oc_\C\D_1\oc_\C\D_{j-1}\simeq\S_i\oc_\C\S_1
\oc_\C\D_{j-1}\to \S_{i+1}\oc_\C\D_{j-1}$ of the maps induced by
the comultiplication in $\D$ and the maps induced by
the semimultiplication in $\S$ (and analogously for $\D\oc_\C\S$).
 These complexes are called the \emph{Koszul complexes} of
the semialgebra $\S$ and the coring $\D$.
 All the grading components of the Koszul complexes with respect
to the grading $i+j$, except the component of degree $i+j=0$,
are acyclic.
 This follows from Lemma~\ref{homogeneous-quadratic-duality}.1
((a)$\Longleftrightarrow$(b)).

\subsection{Central element theorem}  \label{central-element-theorem}
 Let $\C$ be a coring over a $k$\+algebra $A$.
 Assume that $\C$ is a flat right $A$\module.

 A \emph{right coflat increasing filtration} $F$ on a semialgebra $\tS$
over a coring $\C$ is a family of $\C$\+$\C$\bicomodule s $F_n\tS$
endowed with injective morphisms of $\C$\+$\C$\bicomodule s
$F_{n-1}\tS\rarrow F_n\tS$ and an isomorphism of $\C$\+$\C$\bicomodule s
$\ilim F_n\tS\simeq\tS$ such that (i)~$F_i\tS=0$ for $i<0$, \
$F_0\tS=\C$, and the map $F_0\tS\rarrow\tS$ is the semiunit map;
(ii)~the compositions $F_i\tS\oc_\C F_j\tS\rarrow\S\oc_\C\S\rarrow\S$
of the maps induced by the injections $F_n\tS\rarrow\S$ and
the semimultiplication map factorize through $F_{i+j}\tS$;
(iii)~the successive quotients $F_n\tS/F_{n-1}\tS$ are flat right
$A$\module s; and (iv)~the filtration components $F_n\tS$ are
coflat right $\C$\comodule s (then the successive quotients are also
coflat right $\C$\comodule s).
 Assuming that $A$ has a finite weak homological dimension, one
analogously defines \emph{right flat and left relatively coflat
increasing filtrations} by replacing the condition~(iv) with
the condition (iv$'$)~the filtration components $F_n\tS$ are
$\C/A$\+coflat left $\C$\comodule s (then the successive quotients
are also $\C/A$\+coflat left $\C$\comodule s).

\begin{thm}
 Let $\tS$ be a semialgebra over a coring $\C$ endowed with a right
coflat (right flat and left relatively coflat) increasing
filtration~$F$.
 Then the graded semialgebra $\T=\bigoplus_n F_n\tS$ over
the coring~$\C$ is right coflat (right flat and left relatively coflat)
Koszul if and only if the graded semialgebra $\S=\bigoplus_n F_n\tS/
F_{n-1}\tS$ overthe coring~$\C$ is right coflat (right flat and left
relatively coflat) Koszul.
\end{thm}

\begin{proof}
 Consider the reduced bar resolution $\dsb\rarrow\T_+\oc_\C\T_+\oc_\C\T
\rarrow\T_+\oc_\C\T\rarrow\T$ of the right $\T$\semimodule{} $\C$
and denote by $\bcX^\bu$ its semitensor product $\dsb\rarrow
\T_+\oc_\C\T_+\oc_\C\S\rarrow\T_+\oc_\C\S\rarrow\S$ with
the left $\T$\semimodule{} $\S$.
 Denote by $\bcY^\bu$ the two-term complex of graded right
$\S$\semimodule s $\T_1\rarrow \S_0\oplus\S_1$, where $\T_1$ is endowed
with a right $\S$\semimodule{} structure via the augmentation of~$\S$
and $\S_0\oplus\S_1$ is the quotient semimodule $\S/\bigoplus_{n\ge2}
\S_n$; the components of the differential in this complex are
the zero map $\T_1\rarrow\S_0$ and the projection $\T_1\rarrow\S_1$.
 There is a natural morphism of complexes of graded right
$\S$\semimodule s $\bcX^\bu\rarrow \bcY^\bu$ whose components are
the projections $\T_+\oc_\C\S\rarrow\T_1\oc_\C\S_0\simeq\T_1$ and 
$\S\rarrow\S_0\oplus\S_1$.

 All the three complexes $\bcX^\bu$, \ $\bcY^\bu$, and
$\ker(\bcX^\bu\to\bcY^\bu)$ are complexes of $\C$\+coflat right
$\S$\semimodule s ($A$\+flat right $\S$\semimodule s).
 Let us show that the complex $\ker(\bcX^\bu\to\bcY^\bu)$ is
is coacyclic with respect to the exact category of coflat graded
right $\C$\comodule s ($A$\+flat right $\C$\comodule s).
 Indeed, denote by $\bcZ^\bu$ the kernel of the map from the reduced
bar resolution of the right $\T$\semimodule{} $\C$ (written down 
above) to $\C$ itself.
 The complex of graded $\T$\semimodule s $\bcZ^\bu$ has a natural
endomorphism $z$ of internal degree~$1$ and cohomological degree~$0$
induced by the endomorphism of the reduced bar resolution acting by
the identity on the cotensor factors $\T_+$ and by the natural
injections $\T_{n-1}\to\T_n$ on the cotensor factors $\T$.
 Since $\bcZ^\bu$ is a contractible complex of coflat graded right
$\C$\comodule s ($A$\+flat right $\C$\comodule s), the endomorphism $z$
is injective, and its cokernel is a complex of coflat right
$\C$\comodule s ($A$\+flat right $\C$\comodule s), this cokernel is
coacyclic with respect to the exact category of  coflat graded right
$\C$\comodule s ($A$\+flat right $\C$\comodule s).
 Now the kernel $\ker(\bcX^\bu\to\bcY^\bu)$ is isomorphic as
a complex of right $\C$\comodule s to the kernel of a surjective
morphism from $\coker(z)$ to the contractible two-term complex of
coflat right $\C$\comodule s ($A$\+flat right $\C$\comodule s)
$\T_1\rarrow\T_1$.

 Since the semitensor product $\bcX^\bu\os_\S\C$ is isomorphic to
$\Br^\bu_\gr(\T,\C)$, it represents the object $\SemiTor^\S_\gr(\C,\C)$
in the derived category of graded $k$\module s
(see~\ref{one-sided-semitor}).
 On the other hand, since $\bcX^\bu$ is a bounded from above complex
whose terms considered as one-term complexes are semiflat complexes of
graded right $\S$\semimodule s ($\S/\C/A$\semiflat{} complexes of
graded right $\S$\semimodule s), $\bcX^\bu$ is a semiflat complex of
graded right $\S$\semimodule s ($\S/\C/A$\semiflat{} complex of
graded right $\S$\semimodule s).
 The cone of the morphism $\bcX^\bu\rarrow\bcY^\bu$ is coacyclic with
respect to the exact category of coflat graded right $\C$\comodule s
($A$\+flat graded right $\C$\comodule s), so the semitensor product
$\bcX^\bu\os_\S\C$ represents also the object
$\SemiTor^\T_\gr(\bcY^\bu,\C)$ in the derived category of graded
$k$\module s.

 In the semiderived category of graded $\C$\+coflat ($A$\+flat) right
$\S$\semimodule s there is a distinguished triangle $\C(-1)[1]\rarrow
\bcY^\bu\rarrow\C\rarrow\C(-1)[2]$ (where the number in round brackets
denotes the shift of internal grading $M(1)_n=M_{n+1}$).
 It follows from the induced long exact sequence of cohomology of
the objects $\SemiTor^\S_\gr({-},\C)$ by induction in the internal
degree that $\Br^\bu_\gr(\S,\C)$ has nonzero cohomology on the diagonal
$-i=n$ only if and only if $\Br^\bu_\gr(\T,\C)$ has nonzero cohomology
on the diagonal $-i=n$ only.
 Assume that this is so; then there are short exact sequences
$0\rarrow H_{n-1}^{-n+1}\Br^\bu_\gr(\S,\C)\rarrow
H_n^{-n}\Br^\bu_\gr(\T,\C)\rarrow H_n^{-n}\Br^\bu_\gr(\S,\C)\rarrow0$.
 Furthermore, the diagonal cohomology $H_n^{-n}\Br^\bu_\gr(\T,\C)$ and
$H_n^{-n}\Br^\bu_\gr(\S,\C)$ are flat right $A$\module s by
Lemma~\ref{homogeneous-quadratic-duality}.2(a), and so are endowed
with $\C$\+$\C$\bicomodule{} structures.
 The maps $H_n^{-n}\Br^\bu_\gr(\T,\C)\rarrow H_n^{-n}\Br^\bu_\gr(\S,\C)$
in the short exact sequences above are induced by the morphism of
semialgebras $\T\rarrow\S$, hence they are morphisms of
$\C$\+$\C$\bicomodule s.

 Let us describe the compositions $H_n^{-n}\Br^\bu_\gr(\T,\C)\rarrow 
H_n^{-n}\Br^\bu_\gr(\S,\C)\rarrow H_{n+1}^{-n-1}\Br^\bu_\gr(\T,\C)$, which
will be denoted by $\d_n$.
 Let $t\:\C\rarrow\T_1$ be the natural injection.
 Consider the endomorphism $\d_\bcX$ of internal degree~$1$ and
cohomological degree~$-1$ of the complex of graded right
$\S$\semimodule s $\bcX^\bu$ that is defined by the following formulas:
the component $\S$ maps to $\T_+\oc_\C\S$ by $t\oc\id$, the component
$\T_+\oc_\C\S$ maps to $\T_+\oc_\C\T_+\oc_\C\S$ by $t\oc\id\oc\id-
\id\oc t\oc\id$, etc.
 Consider also the endomorphism $\d_\bcY$ of internal degree~$1$ and
cohomological degree~$-1$ of the complex of graded right
$\S$\semimodule s $\bcY$ mapping $\S_0\oplus \S_1$ to $\T_1$ by
the composition of the projection $\S_0\oplus\S_1\to \C$ and
the embedding~$t$.
 Then the endomorphisms $\d_\bcX$ and $\d_\bcY$ form a commutative
diagram with the morphism $\bcX^\bu\rarrow\bcY^\bu$.

 Since the endomorphism $\d_\bcY$ represents in the semiderived category
of $\C$\+coflat ($A$\+flat) graded $\S$\semimodule s the composition of
morphisms $\bcY^\bu\rarrow\C\rarrow\bcY^\bu(1)[-1]$ from
the distinguished triangle above, the desired maps $\d_n$ are induced
by the endomorphism $\d_{\Br}$ of the bar complex $\Br_\gr^\bu(\T,\C)=
\bcX^\bu\os_\S\C$ that is induced by the endomorphism $\d_\bcX$ of
the complex $\bcX^\bu$.
 The endomorphism $\d_{\Br}$ maps the component $\C$ to $\T_+$ by~$t$,
the component $\T_+$ to $\T_+\oc_\C\T_+$ by $t\oc_\C\id-\id\oc_\C t$,
etc.
 Since $\d_{\Br}$ is an endomorphism of complexes of
$\C$\+$\C$\bicomodule s, $\d_n$ are also endomorphisms of
$\C$\+$\C$\bicomodule s.
 Hence the maps $H_{n-1}^{-n+1}\Br^\bu_\gr(\S,\C)\rarrow
H_n^{-n}\Br^\bu_\gr(\T,\C)$ in the short exact sequences above are 
morphisms of $\C$\+$\C$\bicomodule s.
 Now it follows easily by induction using
Lemma~\ref{absolute-relative-coflat} that all
$H_n^{-n}\Br^\bu_\gr(\T,\C)$ are coflat right $\C$\comodule s
($\C/A$\+coflat left $\C$\comodule s) if and only if all
$H_n^{-n}\Br^\bu_\gr(\S,\C)$ are coflat right $\C$\comodule s
($\C/A$\+coflat left $\C$\comodule s).
\end{proof}

 A semialgebra $\tS$ over a coring $\C$ endowed with a right coflat
(right flat and left relatively coflat) increasing filtration~$F$
is called a \emph{right coflat} (\emph{right flat and left relatively
coflat}) \emph{nonhomogeneous Koszul semialgebra} over $\C$ if
the equivalent conditions of Theorem are satisfied for it, i.~e.,
the graded semialgebras $\bigoplus_n F_n\tS$ and
$\bigoplus_n F_n\tS/F_{n-1}\tS$ are right coflat (right flat and
left relatively coflat) Koszul semialgebras over~$\C$.

\subsection{Poincare--Birkhoff--Witt theorem}  \label{pbw-theorem}
 Let $\C$ be a coring over a $k$\+algebra~$A$; assume that $\C$ is
a flat right $A$\module.
 A quasi-differential coring $\tD$ over $A$ concentrated in
the nonnegative degrees and endowed with an isomorphism
$\C\simeq\D_0\til$ is called \emph{right coflat} (\emph{right flat
and left relatively coflat}) \emph{Koszul} over $\C$ if the graded
coring $\tD/\im\d$ is right coflat (right flat and left relatively
coflat) Koszul over~$\C$.

\begin{lem}
 Let\/ $\T$ be a right coflat (right flat and left relatively coflat)
Koszul semialgebra over\/ $\C$ and\/ $\E$ be the quadratic dual right
coflat (right flat and left relatively coflat) Koszul coring over\/ $\C$.
 Then a\/ $\C$\+$\C$\bicomodule{} morphism\/ $\C\rarrow\T_1\simeq\E_1$
can be extended to a graded\/ $\T$\+$\T$\bisemimodule{} morphism\/
$\T\rarrow\T$ of degree~$1$ (i.~e., represents a ``central element''
of\/~$\T$) if and only if it can be extended to a coderivation\/
$\E\rarrow\E$ of degree~$1$ of the coring\/ $\E$ (with respect to
the zero coderivation of $A$).
 Both the\/ $\T$\+$\T$\bisemimodule{} morphism and the coderivation
of\/ $\E$ with the given component\/ $\C\rarrow\T_1\simeq\E_1$ are
unique if they exist; the coderivation always has a zero square.
\end{lem}

\begin{proof}
 Both conditions hold if and only if the difference of the two
maps $\T_1\simeq\C\oc_\C\T_1\rarrow\T_1\oc_\C\T_1$ and $\T_1\simeq
\T_1\oc_\C\C\rarrow\T_1\oc_\C\T_1$ induced by our map $\C\to\T_1$
factorizes through the injection $\E_2\rarrow\E_1\oc_\C\E_1\simeq
\T_1\oc_\C\T_1$.
\end{proof}

 The objects of the category of right coflat nonhomogeneous Koszul
semialgebras are right coflat nonhomogeneous Koszul semialgebras
$(\tS,F)$ over corings $\C$ over $k$\+algebras $A$ such that $\C$
is a flat right $A$\module.
 Morphisms are maps of semialgebras $\tS\rarrow\tS'$ compatible with
maps of corings $\C\rarrow\C'$ and morphisms of $k$\+algebras 
$A\rarrow A'$ which map the filtration components $F_n\tS$ into
the filtration components $F'_n\tS'$.
 Imposing the additional assumption that $A$ has a finite weak
homological dimension, one analogously defines the category of right
flat and left relatively coflat nonhomogeneous Koszul semialgebras.

 The objects of the category of right coflat Koszul quasi-differential
corings are right coflat Koszul quasi-differential corings $\tD$ over
corings $\C$ over $k$\+algebras $A$ such that $\C$ is a flat right
$A$\module.
 Morphisms are maps of graded corings $\tD\rarrow\tD'$ compatible with
morphisms of $k$\+algebras $A\rarrow A'$ and making a commutative
diagram with the coderivations $\d$ and~$\d'$.
 Imposing the additional assumption that $A$ has finite weak
homological dimension, one analogously defines the category of right
flat and left relatively coflat Koszul quasi-differential corings.

\begin{thm}
 The category of right coflat (right flat and left relatively coflat)
nonhomogeneous Koszul semialgebras is equivalent to the category of
right coflat (right flat and left relatively coflat) Koszul
quasi-differential corings.
 If a filtered semialgebra\/ $\tS$ over a coring\/ $\C$ and
a quasi-differential coring\/ $\tD$ correspond to each other under
this duality, then the graded semialgebra\/ $\T=\bigoplus_nF_n\tS$
and the graded coring\/ $\tD$ are quadratic dual right coflat (right
flat and left relatively coflat) Koszul semialgebra and coring over\/
$\C$; the graded semialgebra\/ $\S=\bigoplus_n F_n\tS/F_{n-1}\tS$
and the graded coring\/ $\D=\tD/\im\d$ are quadratic dual right coflat
(right flat and left relatively coflat) Koszul semialgebra and
coring over\/ $\C$; the related isomorphisms $F_1\tS\simeq\D_1\til$
and $F_1\tS/F_0\tS\simeq\D_1\til/\d_0\D_0\til$ are compatible with
each other; and the injection $F_0\tS\rarrow F_1\tS$ corresponds to
the coderivation component $\d_0\:\D_0\til\rarrow\D_1\til$ under the
isomorphisms $F_0\tS\simeq\C\simeq\D_0\til$ and $F_1\tS\simeq\D_1\til$.
\end{thm}

\begin{proof}
 It follows from Lemma that the category of right coflat (right flat
and left relatively coflat) Koszul semialgebras $\T$ endowed with
a $\T$\+$\T$\bisemimodule{} morphism $\T\rarrow\T$ of degree~$1$ is
equivalent to the category of right coflat (right flat and left
relatively coflat) Koszul corings $\E$ endowed with a coderivation
of degree~$1$.
 It remains to prove that semialgebras $\T$ with maps $\T\rarrow\T$
of degree~$1$ coming from right coflat (right flat and left relatively
coflat) nonhomogeneous Koszul semialgebras $\tS$ correspond under
this equivalence to right coflat (right flat and left relatively coflat)
Koszul quasi-differential corings $\tD=\E$ and vice versa.
 Besides, we will have to show that whenever for a quasi-differential
coring $\tD$ the graded coring $\tD/\im\d$ is a right coflat (left
relatively coflat) Koszul coring over a coring $\C$, the graded coring
$\tD$ is also a right coflat (left flat and right relatively coflat)
Koszul coring over~$\C$.

 According to the proof of Theorem~\ref{central-element-theorem}, for
any right coflat (right flat and left relatively coflat) nonhomogeneous 
Koszul semialgebra $\tS$ there is a right coflat (right flat and left
relatively coflat) Koszul quasi-differential coring $\tD$.
 Indeed, set $\tD=\bigoplus_n H_n^{-n}\Br^\bu_\gr(\T,\C)$, where
$\T=\bigoplus_n F_n\tS$; then the endomorphism $\d$ of
the $\C$\+$\C$\bicomodule{} $\tD$ induced by the endomorphism $\d_{\Br}$
of the reduced bar construction $\Br^\bu_\gr(\T,\C)$ is a coderivation of
degree~$1$ (with respect to the zero coderivation of~$A$) and its
restriction to $\D_0\til$ coincides with the injection
$\D_0\til\simeq\T_0\rarrow\T_1\simeq\D_1\til$.
 It also follows from this proof that the right coflat (right flat and
left relatively coflat) Koszul semialgebra $\S=\bigoplus
F_n\tS/F_{n-1}\tS$ is quadratic dual to the coring $\D=\tD/\im\d$,
which is therefore right coflat (right flat and left relatively coflat)
Koszul over~$\C$.

 Let us now construct the nonhomogeneous Koszul semialgebra
corresponding to a right coflat (right coflat and left relatively coflat)
Koszul quasi-differential coring $\tD$ over a coring~$\C$.
 Set $\D=\tD/\im\d$.
 Consider the bigraded coring $\K$ over the $k$\+algebra $A$
(which is considered as a bigraded $k$\+algebra concentrated in
the bidegree~$(0,0)$) with the components $\K^{p,q}=\D_{q-p}\til$ for
$p\le0$, $q\le0$ and $\K^{p,q}=0$ otherwise.
 The coring $\K$ considered as a graded coring in the total grading~$p+q$
has a coderivation $\d_\K$ (with respect to the zero coderivation of~$A$)
mapping the component $\K^{p,q}$ to $\K^{p,q+1}$ by $\d_{q-p}$; one has
$\d_\K^2=0$. 
 There is a morphism of bigraded corings $\K\rarrow\D$ inducing
an isomorphism of the corings of cohomology, where the coring $\D$
is placed in the bigrading $\D^{p,0}=\D_{-p}$ and endowed with the zero
differential.

 Denote by $\K_+$ the cokernel of the injection $\C\simeq\K^{0,0}
\rarrow\K$.
 Let $\bR=\bigoplus_{r=0}^\infty \K_+^{\suboc_\C r}$ be the ``cotensor
semialgebra'' of the bigraded $\C$\+$\C$\bicomodule{} $\K_+$.
 By the definition, $\bR$ is a trigraded semialgebra over the coring~$\C$
(which is considered as a trigraded coring concentrated in
the tridegree~$(0,0,0)$) with the gradings $p$ and $q$ inherited from
the bigrading of~$\K_+$ and the additional grading $r$ by the number of
cotensor factors.
 We will consider $\bR$ as a graded semialgebra in the total grading
$p+q+r$.
 The semialgebra $\bR$ is endowed with three semiderivations (with
respect to the zero derivation of the coring~$\C$) of total degree~$1$,
which we will now introduce.

 Let $\d_\bR$ be the only semiderivation of $\bR$ which preserves
$\K_+\subset\bR$ (embedded as the part of degree $r=1$) and whose
restriction to $\K_+$ is equal to~$-\d_\K$.
 Let $d_\bR$ be the only semiderivation of $\bR$ which maps $\K_+$ to
$\K_+\oc_\C\K_+$ by the composition of the comultiplication map
$\K_+\rarrow\K_+\oc_\C\K_+$ and the sign automorphism of $\K_+\oc_\C\K_+$
acting on the component $\K^{p',q'}\oc_\C\K^{p'',q''}$ as $(-1)^{p'+q'}$.
 Finally, let $\delta_\bR$ be the only semiderivation of $\bR$ whose
restriction to $\K_+$ is the identity map of the component
$\K^{-1,-1}\simeq\C$ to the semiunit component $\bR^{0,0,0}=\C$ and zero
on all the remaining components of $\K_+$.
 All the three differentials are constructed so that they satisfy
the super-Leibniz rule in the parity $p+q+r$.
 The semiderivations $\d_\bR$, \ $d_\bR$, and $\delta_\bR$ have
tridegrees $(0,1,0)$, \ $(0,0,1)$, and $(1,1,-1)$, respectively,
in the trigrading $(p,q,r)$.
 All the three semiderivations have zero squares, and they
pairwise anti-commute.

 There is a right coflat (right flat and left relatively coflat)
increasing filtration $F$ on the graded semialgebra $\bR$ whose
component $F_n\bR$ is the direct sum of all trigrading components
$\bR^{p,q,r}$ with $-p\le n$.
 This filtration is compatible with the differentials $\d_\bR$, \
$d_\bR$, and $\delta_\bR$; the semialgebra $\bigoplus_nF_n\bR/F_{n-1}\bR$
with the differential induced by $\d_\bR+d_\bR+\delta_\bR$ is 
naturally isomorphic to the semialgebra $\bR$ with the differential
$\d_\bR+d_\bR$.

 Consider the following sign-modified version of cobar construction
${}'\!\Cb(\D,\C)$.
 Define ${}'\!\Cb(\D,\C)$ as the ``tensor semialgebra'' $\bigoplus_r
\D_+^{\suboc_\C r}$ of the $\C$\+$\C$\bicomodule{} $\D_+$ and endow it
with the grading $p$ coming from the grading $\D^p=\D_{-p}$ of $\D_+$
and the grading $r$ by the number of cotensor factors.
 We will consider ${}'\!\Cb(\D,\C)$ as a graded semialgebra over $\C$
in the total grading $p+r$.
 Let $d_{\Cb}'$ be the only coderivation of ${}'\!\Cb(\D,\C)$ which
maps $\D_+\subset{}'\!\Cb(\D,\C)$ to $\D_+\oc_\C\D_+$ by
the composition of the comultiplication map $\D_+\rarrow\D_+\oc_\C\D_+$
and the sign automorphism of $\D_+\oc_\C\D_+$ acting on the component
$\D^{p'}\oc_\C\D^{p''}$ as $(-1)^{p'}$.
 Then one has $d_{\Cb}^{\prime2}=0$.
 Notice that the differential $d_{\Cb}'$ satisfies the super-Leibniz
rule in the parity $p+r$, while the differential $d_{\Cb}$ of
the cobar construction $\Cb^\bu_\gr(\D,\C)$ satisfies the super-Leibniz
rule in the parity~$r$.
 The automorphism of $\bigoplus_r\D_+^{\suboc_\C r}$ acting on
the component $\D^{p_1}\oc_\C\dsb\oc_\C\D^{p_r}$ by minus one to
the power $\sum_{s=1}^r p_s(p_s+1)/2+\sum_{1\le s<t\le r} p_s(p_t+1)$
transforms $d_{\Cb}$ to $d_{\Cb}'$, so the semialgebras of cohomology
of the DG\semialgebra s ${}'\!\Cb(\D,\C)$ and $\Cb^\bu_\gr(\D,\C)$
are naturally isomorphic in the Koszul case.

 Consider the morphism of DG\semialgebra s $(\bR\;\d_\bR+d_\bR)\rarrow
({}'\!\Cb(\D,\C)\;d'_{\Cb})$ induced by the morphism of corings
$\K\rarrow\D$.
 This morphism of DG\semialgebra s induces an isomorphism of
the semialgebras of cohomology. 
 Indeed, the components of fixed grading~$p$ of the DG\semialgebra{}
$(\bR\;\d_\bR+d_\bR)$ are the total components of finite bicomplexes
whose components of fixed grading~$r$ are multiple cotensor products
of components of fixed grading~$p$ of the DG\+coring $\K$, and
the natural maps from these multiple cotensor products to
the corresponding multiple cotensor products of components of
the coring $\D$ are quasi-isomorphisms.
 Hence $H^0_{\d_\bR+d_\bR}(\bR)\simeq\S$ and $H^i_{\d_\bR+d_\bR}(\bR)=0$
for $i\ne 0$, where $\S$ denotes the right coflat (right flat and left
relatively coflat) Koszul semialgebra quadratic dual to $\D$.
 Analogously, the morphism of DG\semialgebra s $(\bR\;\d_\bR+d_\bR)
\rarrow ({}'\!\Cb(\D,\C)\;d'_{\Cb})$ induces quasi-isomorphisms of
the tensor and cotensor products related to these DG\semialgebra s
that were listed in (i) and (iii) of~\ref{diff-semialgebras-ii}.

 The associated graded quotient complexes to the tensor and cotensor
product of the DG\semialgebra{} $(\bR\;\d_\bR+d_\bR+\delta_\bR)$ listed
in (i) and (iii) of~\ref{diff-semialgebras-ii} with respect to
the filtrations induced by the filtration~$F$ are naturally isomorphic
to the corresponding tensor and cotensor products of the DG\semialgebra{}
$(\bR\;\d_\bR+d_\bR)$.
 Therefore, the associated graded modules of the cohomology of these
tensor and cotensor products of the DG\semialgebra{}
$(\bR\;\d_\bR+d_\bR+\delta_\bR)$ are isomorphic to the cohomology of
the corresponding tensor and cotensor products of the DG\semialgebra{}
$(\bR\;\d_\bR+d_\bR)$.
 In particular, set $\tS=H_{\d_\bR+d_\bR+\delta_\bR}^0(\bR)$; then $\tS$
is endowed with an increasing filtration~$F$ such that
$\bigoplus_n F_n\tS/F_{n-1}\tS\simeq \S$, while
$H_{\d_\bR+d_\bR+\delta_\bR}^i(\bR)=0$ for $i\ne0$.
 Since $\S$ is a coflat right $\C$\comodule{} (a flat right
$A$\module{} and a $\C/A$\+coflat left $\C$\comodule), the associated
graded quotient modules to the tensor and cotensor products under
consideration of the cohomology module $\tS$ are isomorphic to
the corresponding tensor and cotensor products of~$\S$.
 Thus $\tS$ is a semialgebra over~$\C$ and $F$ is its right coflat
(right flat and left relatively coflat) increasing filtration.

 Since the semialgebra $\S$ is right coflat (right flat and left
relatively coflat) Koszul, so is the semialgebra $\T=\bigoplus_nF_n\tS$.
 Let $\tD'$ be the right coflat (right flat and left relatively
coflat) coring quadratic dual to~$\T$; then $\tD'$ is endowed
with a coderivation $\d'$, making it a right coflat (right flat
and left relatively coflat) Koszul quasi-differential coring, as
we have already proven.
 Moreover, the cokernel $\D'$ of the coderivation~$\d'$ is
quadratic dual to $\S$, hence there is a natural isomorphism
of graded corings $\D\simeq\D'$.
 Furthermore, the embedding of the component $\D_1\til=\bR_{-1,0,1}
\rarrow\bR$ induces an isomorphism $\D_1\til\simeq F_1\tS$.
 The composition $\D_2\til\rarrow\D_1\til\oc_\C\D_1\til\simeq
F_1\tS\oc_\C F_1\tS\rarrow F_2\tS$ of the comultiplication and
semimultiplication maps vanishes, so there is a natural morphism
of graded corings $\tD\rarrow\tD'$.
 Since the embedding $F_0\tS\rarrow F_1\tS$ corresponds to
the map $\d_0\:\D_0\til\rarrow\D_1\til$ under the isomorphisms
$F_0\tS\simeq\C\simeq\D_0\til$ and $F_1\tS\simeq\D_1\til$,
the morphism $\tD\rarrow\tD'$ forms a commutative diagram with
the differentials $\d$ in $\tD$ and $\d'$ in $\tD'$.
 The induced morphism $\tD/\im\d\rarrow\tD'/\im\d'$ coincides
with the natural isomorphism $\D\rarrow\D'$ on the components
of degree~$1$, and consequently on the other components as
well.
 Hence the morphism of corings $\tD\rarrow\tD'$ is also
an isomorphism.
 Thus the coring $\tD$ is right coflat (right flat and left relatively
coflat) Koszul over~$\C$, and the semialgebra $\T$ quadratic dual to
it together with its $\T$\+$\T$\bicomodule{} endomorphism of
degree~$1$ comes from the right coflat (right flat and left
relatively coflat) nonhomogeneous Koszul semialgebra~$\tS$.
\end{proof}

 A right coflat (right flat and left relatively coflat) nonhomogeneous
Koszul semialgebra and a right coflat (right flat and left relatively
coflat) Koszul quasi-differential coring corresponding to each other
under the equivalence of categories from the above Theorem are called
\emph{nonhomogeneous quadratic dual} to each other.

 All the definitions and results
of~\ref{koszul-semi-co}--\ref{pbw-theorem} have their obvious
analogues for the coflatness conditions replaced with
coprojectivity ones.
 So, when $\C$ is a projective left $A$\module{} one can speak
of left coprojective Koszul semialgebras and corings.
 When $\C$ is a flat right $A$\module{} and $A$ has a finite left
homological dimension, one can define right flat and left
relatively projective Koszul semialgebras and corings.
 When $\C$ is a projective left $A$\module{} and $A$ has a finite
left homological dimension, one can consider left projective and
right relatively coflat Koszul semialgebras and corings.
 Of course, when $\C$ is a flat left $A$\module, one can define
left coflat Koszul semialgebras and corings, etc.

\begin{rmk}
 All the results of~\ref{koszul-semi-co}--\ref{pbw-theorem} have
their analogues for semimodules and semicontramodules over
semialgebras, comodules and contramodules over corings.
 In particular, for any right coflat Koszul semialgebra $\S$ and
the right coflat Koszul coring $\D$ quadratic dual to $\S$ there
is a natural equivalence between the categories of Koszul left
$\S$\semimodule s and Koszul left $\D$\comodule s given by
the functors of cohomology of the reduced bar and cobar constructions
with coefficients in the semimodules and comodules.
 No (co)flatness conditions need to be imposed on the semimodules
and comodules in this setting.
 For any left coprojective Koszul semialgebra $\S$ and the left
coprojective Koszul coring $\D$ quadratic dual to $\S$ there is
an equivalence between the categories of Koszul left
$\S$\semicontramodule s and Koszul left $\D$\contramodule s
(which are nonpositively graded).
 For any right flat and left relatively coflat Koszul semialgebra
$\S$ and the right flat and left relatively coflat Koszul coring
$\D$ quadratic dual to $\S$ there is an equivalence between
the categories of $A$\+flat Koszul right $\S$\semimodule s and
$A$\+flat Koszul right $\D$\comodule s, etc.
 Furthermore, for a right coflat nonhomogeneous Koszul semialgebra
$\tS$ and a left semimodule $\bM\til$ over $\tS$ endowed with
an increasing filtration $F$ compatible with the filtration of
$\tS$, the semimodule $\bigoplus_n F_n\bM\til$ is Koszul
over the semialgebra $\bigoplus_n F_n\tS$ if and only if
the semimodule $\bigoplus_n F_n\bM\til/F_{n-1}\bM\til$ is Koszul
over the semialgebra $\bigoplus_n F_n\tS/F_{n-1}\tS$.
 A filtered semimodule $\M\til$ satisfying these conditions can be
called a nonhomogeneous Koszul semimodule over~$\tS$.
 The Koszul $\D$\comodule{} quadratic dual to the second of these
graded semimodules is naturally isomorphic to the Koszul
$\tD$\comodule{} quadratic dual to the first semimodule with
the induced $\D$\comodule{} structure.
 A Koszul quasi-differential left $\tD$\comodule{} is a graded left
$\tD$\comodule{} that is Koszul as a $\D$\comodule; then it is
also Koszul as a $\tD$\comodule.
 There is a natural equivalence between the categories of
nonhomogeneous Koszul left semimodules over $\tS$ and Koszul
quasi-differential comodules over $\tD$.
 When $\tS$ is a left coprojective nonhomogeneous Koszul semialgebra,
there is an analogous equivalence of categories of nonhomogeneous
Koszul left semicontramodules over $\tS$ and Koszul quasi-differential
contramodules over $\tD$ (where nonhomogeneous Koszul semicontramodules
are semicontramodules endowed with complete decreasing filtrations).
\end{rmk}

\subsection{Quasi-differential comodules and contramodules}

\subsubsection{}  \label{quasi-differential-comodules}
 Let $(\tD,\d)$ be a quasi-differential coring over a $k$\+algebra $A$;
set $\D=\tD/\im\d$.
 Assume that $\D$ is a flat graded right $A$\module.

 A \emph{quasi-differential left comodule} over $\tD$ is just a graded
left $\tD$\comodule{} (without any differential).
 The DG\category{} of quasi-differential left $\tD$\comodule s
$\DG(\tD\qcomodl)$ is defined as follows.
 The objects of $\DG(\tD\qcomodl)$ are quasi-differential left
$\tD$\comodule s.
 Let us construct the complex of morphisms in the category
$\DG(\tD\qcomodl)$ between quasi-differential left $\tD$\comodule s
$\L$ and $\M$, denoted by $\Hom_\D^\bu(\L,\M)$.
 The component $\Hom_\D^n(\L,\M)$ of this complex is the $k$\module{}
of all homogeneous maps $\L\rarrow\M$ of degree~$-n$ supercommuting
with the coaction of $\D$ in $\L$ and $\M$.
 This means that an element $f\in\Hom_\D^n(\L,\M)$ should satisfy
the equation $\overline{f(x)_{(-1)}\!\.}\.\ot f(x)_{(0)} =
(-1)^{n|x_{(-1)}|} \overline{x_{(-1)}\!\.}\.\ot f(x_{(0)})$ in
Sweedler's notation~\cite{Swe}, where $z\mpsto z_{(-1)}\ot z_{(0)}$
denotes the left coaction maps, $\overline{p}$ is the image of
an element $p\in\tD$ in $\D$, and $|p|$ is the degree of
a homogeneous element~$p$.
 To define the differential $d(f)$ of an element $f$, consider
the supercommutator of the coaction maps $\L\rarrow\tD\oc_\D\L$ and
$\M\rarrow\tD\oc_\D\M$ with $f$, that is the map $\delta_f\:
\L\rarrow\tD\oc_\D\M$ given by the formula $x\mpsto f(x)_{(-1)}
\oc f(x)_{(0)} - (-1)^{n|x_{(-1)}|}x_{(-1)}\oc f(x_{(0)})$.
 For any $f\in\Hom_\D^n(\L,\M)$, the map $\delta_f$ factorizes
through the injection $\M\simeq\D\oc_\D\M\rarrow\tD\oc_\D\M$ induced
by the homogeneous morphism $\dbar\:\D\rarrow\tD$ of degree~$1$
given by $\dbar(\overline{p})=\d(p)$, hence the desired map
$d(f)\:\L\rarrow\M$ of degree $-n-1$.

 Since the map $\delta_f$ and the morphism $\dbar$ supercommute
with the left coactions of~$\D$, so does the map $d(f)$.
 Let us check that $d^2(f)=0$, in other words, that $d(f)$ 
supercommutes with the left coactions of $\tD$ in $\L$ and $\M$.
 Consider the two homogeneous maps $\L\rarrow\tD\oc_\D\M$ given
by the formulas $x\mpsto (df)(x)_{(-1)}\oc(df)(x)_{(0)}$ and
$x\mpsto (-1)^{(n+1)|x_{(-1)}|}x_{(-1)}\oc (df)(x_{(0)})$; we
have to check that these two maps coincide.
 Consider the image of the former map under the map
$\Hom^{n+1}_\D(\L\;\tD\oc_\D\M)\rarrow \Hom^{n+1}_\D(\L\;
\tD\oc_\D\tD\oc_\D\M)$ given by the formula $g\mpsto
(x\mapsto (-1)^{(n+1)|x_{(-1)}|}\d(x_{(-1)})\oc g(x_{(0)}))$
and the image of the second map under the map
$\Hom^{n+1}_\D(\L\;\tD\oc_\D\M)\rarrow \Hom^{n+1}_\D(\L\;
\tD\oc_\D\tD\oc_\D\M)$ given by the formula $g\mpsto
(x\mapsto (-1)^{|g(x)_1|}g(x)_1\oc\d(g(x)_{2(-1)})\oc g(x)_{2(0)})$,
where $y=y_1\oc y_2$ is a notation for an element $y\in\tD\oc_\D\M$.
 The sum of these two elements of $\Hom^{n+1}_\D(\L\;\tD\oc_\D\tD
\oc_\D\M)$ is equal to the image of the element $\delta_f$ under
the map $\Hom^{n+1}_\D(\L\;\tD\oc_\D\M)\rarrow \Hom^{n+1}_\D(\L\;
\tD\oc_\D\tD\oc_\D\M)$ induced by the comultiplication map
$\tD\rarrow\tD\oc_\D\tD$.
 There is a commutative square formed by the diagonal embedding
$\tD\rarrow\tD\oplus\tD$, the morphism $\tD\oplus\tD\rarrow
\tD\oc_\D\tD$ given by the formula $(x,y)\mpsto \d(x_{(1)})
\oc x_{(2)} + (-1)^{|y_{(1)}|}y_{(1)}\oc\d(y_{(2)})$, the morphism
$\d\:\tD\rarrow\tD$, and the comultiplication morphism $\tD\rarrow
\tD\oc_\D\tD$.
 Considering the filtrations originating from the two-term
filtration $\d(\tD)\subset\tD$, one can check that this square
is Carthesian.
 It remains Carthesian after applying the functors ${-}\oc_\D\M$
and $\Hom_\D^{n+1}(\L,{-})$, so we are done.

 Let $\M$ be a quasi-differential left $\tD$\comodule{} and
$q\:\M\rarrow\M$ be an element of $\Hom_\D^1(\M,\M)$ satisfying
the Maurer--Cartan equation $d(q)+q^2=0$.
 The quasi-differential left $\tD$\comodule{} structure on $\M$ twisted
with~$q$ is constructed as follows.
 First of all, the structure of a graded $\D$\comodule{} on $\M$ does
not change under twisting.
 Next, the twisted coaction map $\M\rarrow\tD\oc_\D\M$ is the sum
of the original coaction map and the composition $\M\rarrow\M
\rarrow\tD\oc_\D\M\rarrow\tD\oc_\D\M$ of the map $q\:\M\rarrow \M$,
the coaction map $\M\rarrow\tD\oc_\D\M$, and the map $\tD\oc_\D\M
\rarrow\tD\oc_\D\M$ induced by the morphism~$\d\:\tD\rarrow\tD$.
 Denote the quasi-differential $\tD$\comodule{} we have constructed
by $\M(q)$.
 For any quasi-differential left $\tD$\comodule{} $\L$ the differential
in the complex $\Hom_\D^\bu(\L,\M(q))$ differs from the differential
in the complex $\Hom_\D^\bu(\L,\M)$ according to the formula
$d_q(f)=d(f)+q\circ f$.

 Since infinite direct sums and shifts of objects clearly exist
in the DG\category{} $\DG(\tD\qcomodl)$, it follows from the above
construction, in particular, that cones exist in it.
 Therefore, the homotopy category $\Hot(\tD\qcomodl)$, whose objects are
the objects of $\DG(\tD\qcomodl)$ and morphisms are the zero cohomology
of the complexes of morphisms in $\DG(\tD\qcomodl)$, is naturally
triangulated.
 Furthermore, for any complex of quasi-differential left $\tD$\comodule s
one can define the total graded quasi-differential left $\tD$\comodule{}
such that the corresponding graded left $\D$\comodule{} will coincide 
with the infinite direct sum of the shifts of the terms of the complex
considered as graded left $\D$\comodule s.
 In particular, one can speak about the total quasi-differential left
$\tD$\comodule s of exact triples of quasi-differential left
$\tD$\comodule s, which allows us to define the coderived category
of quasi-differential left $\tD$\comodule s $\sD^\co(\tD\qcomodl)$ as
the quotient category of the homotopy category $\Hot(\tD\qcomodl)$
by its minimal triangulated subcategory containing such objects
associated to exact triples and closed under infinite direct sums.

\subsubsection{}
 Let $(\tD,\d)$ be a quasi-differential coring over $A$; 
assume that $\D=\tD/\im\d$ is a flat graded left $A$\module.
 
 A \emph{quasi-differential right comodule} over $\tD$ is just a graded
right $\tD$\comodule.
 Let us define the DG\category{} of quasi-differential right
comodules $\DG(\qcomodr\tD)$ over $\tD$.
 The objects of $\DG(\qcomodr\tD)$ are quasi-differential right
$\tD$\comodule s.
 The complex of morhisms $\Hom_\tD^\bu(\R,\N)$ in the category
$\DG(\qcomodr\tD)$ between quasi-differential right $\tD$\comodule s
$\R$ and $\N$ is constructed as follows.
 The component $\Hom_\tD^n(\R,\N)$ of this complex is the $k$\module{}
of all homogeneous maps $\R\rarrow\N$ of degree~$-n$ commuting with
the $\D$\comodule{} structures (without any signs).
 To define the differential of an element $f\in\Hom^\bu_\D(\R,\N)$,
consider the map $\delta_f\:\R\rarrow\N\oc_\D\tD$ given by
the formula $x\mpsto f(x)_{(0)}\oc f(x)_{(1)} - f(x_{(0)})\oc x_{(1)}$,
where $z\mpsto z_{(0)}\ot z_{(1)}$ denotes the right coaction maps.
 The map~$\delta_f$ factorizes through the injection $\N\rarrow\N
\oc_\D\tD$ given by the formula $y\mpsto (-1)^{|y_{(0)}|}y_{(0)}\oc
\d(y_{(1)})$, hence the desired map $d(f)\:\R\rarrow\N$.

 Let $\N$ be a quasi-differential right $\tD$\comodule{} and
$q\in\Hom_\D^1(\N,\N)$ be an element satisfying the equation
$d(q)+q^2=0$.
 To define quasi-differential right $\tD$\comodule{} structure on $\N$
twisted with~$q$, set the new coaction map $\N\rarrow\N\oc_\D\tD_1$
to be the sum of the original coaction map and the composition 
$\N\rarrow\N\rarrow\tD\oc_\D\N$ of the map $q\:\N\rarrow\N$ and
the map $\N\rarrow\N\oc_\D\tD$ given by the formula $y\mpsto
(-1)^{|y_{(0)}|}y_{(0)}\oc \d(y_{(1)})$.
 Denote the quasi-differential $\tD$\comodule{} so constructed
by $\N(q)$; for any quasi-differential right $\tD$\comodule{} $\R$
the differential in the complex $\Hom_\D^\bu(\R,\N(q))$ differs from
the differential in the complex $\Hom_\D^\bu(\R,\N)$ by the rule
$d_q(f)=d(f)+q\circ f$.

 The definitions of the homotopy category of quasi-differential right
$\tD$\comodule s $\Hot(\qcomodr\tD)$ and the coderived category of
quasi-differential right $\tD$\comodule s $\sD^\co(\qcomodr\tD)$ are
the same as in the left quasi-differential comodule case.

\subsubsection{}  \label{quasi-differential-contramodules}
 Let $(\tD,\d)$ be a quasi-differential coring over $A$; assume 
that $\D=\tD/\im\d$ is a projective graded left $A$\module.

 A \emph{quasi-differential left contramodule} over $\tD$ is just
a graded left $\tD$\contramodule.
 Let us define the DG\category{} of quasi-differential left
$\tD$\contramodule s $\DG(\tD\qcontra)$.
 The objects of $\DG(\tD\qcontra)$ are quasi-differential left
$\tD$\contramodule s.
 The complex of morphisms $\Hom^{\D,\bu}(\P,\Q)$ in the category
$\DG(\tD\qcontra)$ between quasi-differential left $\tD$\contramodule s
$\P$ and $\Q$ is constructed as follows.
 The component $\Hom^{\D,n}(\P,\Q)$ of this complex is the $k$\module{}
of all homogeneous maps $\P\rarrow\Q$ of degree~$-n$ supercommuting
with the $\D$\contramodule{} structures.
 This means that an element $f\in\Hom^{\D,n}(\P,\Q)$ should satisfy
the equation $\pi_\P(f\circ x)=(-1)^{mn}f(\pi_\P(x))$ for any
$x\in\Hom_A(\D_m,\P)$, where $\pi_\P$ denotes the contraaction map.
 To define the differential of an element $f\in\Hom^{\D,n}(\P,\Q)$,
consider the map $\delta_f\:\Cohom_\D(\tD,\P)\rarrow\Q$ given by
the formula $\overline{x}\mpsto \pi_\P(f\circ x)-(-1)^{mn}f(\pi_\P(x))$
for $x\in\Hom_A(\D_m\til,\P)$, where $\overline{x}$ denotes the class
of~$x$ in $\Cohom_\D(\tD,\P)$.
 The map $\delta_f$ factorizes through the surjection
$\Cohom_\D(\tD,\P)\rarrow\Cohom_\D(\D,\P)\simeq\P$ induced by
the morphism $\dbar\:\D\rarrow\tD$, hence the desired map
$d(f)\:\P\rarrow\Q$.

 Let $\P$ be a quasi-differential left $\tD$\contramodule{} and
$q\in\Hom^{\D,1}(\P,\P)$ be an element satisfying the equation
$d(q)+q^2=0$.
 To define the quasi-differential left $\tD$\contramodule{} structure
on $\P$ twisted with~$q$, set the new contraaction map
$\Cohom_\D(\tD,\P)\rarrow\P$ to be the sum of the original
contraaction map and the composition of the map $\Cohom_\D(\tD,\P)
\rarrow\P$ induced by $\dbar\:\D\rarrow\tD$ and the map
$q\:\P\rarrow\P$.
 Denote the quasi-differential left $\tD$\contramodule{} so
constructed by $\P(q)$; for any quasi-differential left 
$\tD$\contramodule{} $\Q$ the differential in the complex
$\Hom^{\D,\bu}(\Q,\P(q))$ differs from the differential in the complex
$\Hom^{\D,\bu}(\Q,\P)$ by the rule $d_q(f)=d(f)+q\circ f$.

 The definitions of the homotopy category of quasi-differential left
$\tD$\contramodule s $\Hot(\tD\qcontra)$ and the contraderived
category of quasi-differential left $\tD$\contramodule s
$\sD^\ctr(\tD\qcontra)$ are completely analogous to the corresponding
definitions in the comodule case; the only difference is that one
considers infinite products instead of infinite direct sums.

\begin{rmk}
 One can define the DG\+categories of CDG\comodule s and
CDG\contramodule s over a CDG\+coring
(see~\ref{quasi-differential-corings}
and~\ref{quasi-differential-rings}) and identify them with
the DG\+categories of quasi-differential comodules and contramodules
in the case when a quasi-differential coring corresponds to
a CDG\+coring.
 More generally, let $(\tD,\d)$ be a quasi-differential coring over
a $k$\+algebra $A$ such that the components $\D_i$ of the coring
$\D=\tD/\im\d$ are projective left $A$\module s.
 Then there is a natural structure of quasi-differential $k$\+algebra
on the graded $A$\+$A$\bimodule{} with the components
$R^n{}\til=\Hom_A(\D_n\til,A)$.
 Let $(R,d,h)$ be a CDG\+algebra over~$k$ corresponding to $R\til$.
 In this situation the DG\category{} of quasi-differential right
$\tD$\comodule s is isomorphic to a full subcategory of the DG\category{}
of right CDG\module s over $(R,d,h)$, and there is a forgetful functor 
from the DG\category{} of quasi-differential left $\tD$\contramodule s
to the DG\category{} of left CDG\module s over $(R,d,h)$.
\end{rmk}

\subsection{Koszul duality}   \label{koszul-duality-theorem}
 Let $\C$ be a coring over a $k$\+algebra $A$.

\begin{thm}
 \textup{(a)} Assume that\/ $\C$ is a flat right $A$\module.
 Let\/ $\tS$ be a right coflat nonhomogeneous Koszul semialgebra
over the coring\/ $\C$ and\/ $\tD$ be the right coflat Koszul
quasi-differential coring over\/ $\C$ nonhomogeneous quadratic dual
to\/~$\tS$.
 Then the semiderived category of left\/ $\tS$\semimodule s is
naturally equivalent to the coderived category of quasi-differential
left\/ $\tD$\comodule s. \par
 \textup{(b)} Assume that\/ $\C$ is a flat left $A$\module.
 Let\/ $\tS$ be a left coflat nonhomogeneous Koszul semialgebra
over the coring\/ $\C$ and\/ $\tD$ be the left coflat Koszul
quasi-differential coring over\/ $\C$ nonhomogeneous quadratic dual
to\/~$\tS$.
 Then the semiderived category of right\/ $\tS$\semimodule s is
naturally equivalent to the coderived category of quasi-differential
right\/ $\tD$\comodule s. \par
 \textup{(c)} Assume that\/ $\C$ is a projective left $A$\module.
 Let\/ $\tS$ be a left coprojective nonhomogeneous Koszul semialgebra
over the coring\/ $\C$ and\/ $\tD$ be the left coprojective Koszul
quasi-differential coring over\/ $\C$ nonhomogeneous quadratic dual
to\/~$\tS$.
 Then the semiderived category of left\/ $\tS$\semicontramodule s
is naturally equivalent to the contraderived category of
quasi-differential left\/ $\tD$\contramodule s.
\end{thm}

\begin{proof}
 Part~(a): let us construct a pair of adjoint functors between
the DG\category{} of complexes of left $\tS$\semimodule s and
the DG\category{} of quasi-differential left $\tD$\comodule s.
 The functor $\Xi$ assigns to a left $\S$\semimodule{} $\bM$
the graded left $\D$\comodule{} $\D\oc_\C\bM$ endowed with the following
left $\tD$\comodule{} structure.
 Consider the map $\D_n\til\oc_\C\bM\rarrow\D_n\til\oc_\C\bM$ equal
to the sum of the identity map and $(-1)^n$ times the composition
$\D_n\til\oc_\C\bM\rarrow\D_{n-1}\til\oc_\C\D_1\til\oc_\C\bM\rarrow
\D_{n-1}\til\oc_\C\bM\rarrow\D_n\til\oc_\C\bM$ of the map induced by
the comultiplication morphism $\D_n\til\rarrow\D_{n-1}\til\oc_\C
\D_1\til$, the map induced by the semiaction morphism $F_1\tS\oc_\C\bM
\rarrow\bM$, and the map induced by the morphism $\d_{n-1}\:\D_{n-1}\til
\rarrow\D_n\til$.
 This map factorizes through the surjection $\D_n\til\oc_\C\bM\rarrow
\D_n\oc_\C\bM$, since its composition with the map $\D_{n-1}\til\oc_\C\bM
\rarrow\D_n\til\oc_\C\bM$ induced by the morphism $\d_{n-1}$ vanishes.
 So we obtain a natural map $\D_n\oc_\C\bM\rarrow\D_n\til\oc_\C\bM$.
 Now the compositions $\D_{i+j}\oc_\C\bM\rarrow\D_{i+j}\til\oc_\C\bM
\rarrow\D_i\til\oc_\C\D_j\til\oc_\C\bM\rarrow\D_i\til\oc_\C\D_j\oc_\C\bM$
of the maps $\D_{i+j}\oc_\C\bM\rarrow\D_{i+j}\til\oc_\C\bM$ we have
constructed with the maps induced by the comultiplication maps
$\D_{i+j}\til\rarrow\D_i\til\oc_\C\D_j\til$ and the maps induced by
the natural surjections $\D_j\til\rarrow\D_j$ define the desired
graded $\tD$\comodule{} structure on $\D\oc_\C\bM$.
 To a complex of left $\tS$\semimodule s $\bM^\bu$ the functor $\Xi$
assigns the total quasi-differential $\tD$\comodule{} of the complex
of quasi-differential $\tD$\comodule s $\D\oc_\C\bM^\bu$.

 The functor $\Ups$ assigns to a quasi-differential left $\tD$\comodule{}
$\L$ the complex of left $\tS$\semimodule s $\Ups^\bu(\L)=\tS\oc_\C\L$
with the terms $\Ups^i(\L)=\tS\oc_\C\L_{-i}$ and the differential defined
as the composition $\tS\oc_\C\L\rarrow\tS\oc_\C\D_1\til\oc_\C\L\rarrow
\tS\oc_\C\L$ of the map induced by the coaction morphism $\L\rarrow
\D_1\til\oc_\C\L$ and the map induced by the semimultiplication morphism
$\tS\oc_\C F_1\tS\rarrow\tS$.
 The functor $\Xi$ is right adjoint to the functor $\Ups$, since both
complexes $\Hom_\D^\bu(\L,\Xi(\bM^\bu))$ and $\Hom_\tS(\Ups^\bu(\L),
\bM^\bu)$ are naturally isomorphic to the total complex of the bicomplex
$\Hom_\C(\L_i,\bM^j)$, one of whose differentials is induced by
the differential in $\bM^\bu$ and the other one assigns to
a $\C$\comodule{} morphism $f\:\L_i\rarrow\bM^j$ the composition
$\L_{i+1}\rarrow\D_1\til\oc_\C\L_i\rarrow\D_1\til\oc_\C\bM^j\simeq
F_1\tS\oc_\C\bM^j\rarrow\bM^j$ of the coaction map, the map induced
by the morphism~$f$, and the semiaction map.

{\hbadness=2000
 Let us show that the functors $\Xi$ and $\Ups$ induce mutually inverse
equivalences between the semiderived category $\sD^\si(\tS\simodl)$ and
the coderived category $\sD^\co(\tD\qcomodl)$.
 Firstly, the functor $\Xi$ sends $\C$\coacyclic{} complexes of 
$\S$\semimodule s to coacyclic quasi-differential $\tD$\comodule s.
 Indeed, for any complex of left $\S$\semimodule s $\bM^\bu$
the quasi-differential $\tD$\comodule{} $\Xi(\bM^\bu)=\D\oc_\C\bM^\bu$
has an increasing filtration by quasi-differential $\tD$\subcomodule s
defined by the formula $F_n(\D\oc_\C\bM^\bu)=\bigoplus_{i\le n}
\D_i\oc_\C\bM^\bu$.
 The associated graded quotient quasi-differential comodule to this
filtration is described as follows.
 There is a functor from the DG\category{} of complexes of
$\C$\comodule s to the DG\category{} of quasi-differential
$\tD$\comodule s assigning to a complex of $\C$\comodule s the total
quasi-differential $\tD$\comodule{} of the complex of quasi-differential
$\tD$\comodule s whose terms are the terms of the original complex of
$\C$\comodule s endowed with the graded $\tD$\comodule{} structure
via the embedding $\C\simeq\D_0\til\rarrow\tD$.
 (This functor can be also described in terms of the quasi-differential
subcoring in $\tD$ whose components are $\D_0\til$ and $\d_0(\D_0\til)
\subset\D_1\til$.)
 Clearly, this functor sends coacyclic complexes of $\C$\comodule s
to coacyclic quasi-differential $\tD$\comodule s.
 Now the quasi-differential $\tD$\comodule s $F_n\Xi(\bM^\bu)/
\allowbreak F_{n-1}\Xi(\bM^\bu)$ are isomorphic to the images of
the $\C$\comodule s $\D_n\oc_\C\bM^\bu$ under this functor, and
are, therefore, coacyclic whenever $\bM^\bu$ is $\C$\coacyclic. \par}

 Secondly, the functor $\Ups$ sends coacyclic quasi-differential
$\tD$\comodule s to complexes of $\tS$\semimodule s that are
coacyclic not only over $\C$, but even over~$\tS$. 

 Thirdly, let us check that for any complex of left $\S$\semimodule s
$\bM^\bu$ the cone of the natural morphism of complexes of 
$\S$\semimodule s $\Ups^\bu\Xi(\bM^\bu)\rarrow\bM^\bu$ is coacyclic as
a complex of left $\C$\comodule s.
 The complex of $\C$\comodule s $\Ups^\bu\Xi(\bM^\bu)=
\tS\oc_\C\D\oc_\C\bM^\bu$ has an increasing filtration given by 
the formula $F_n\Ups^\bu\Xi(\bM^\bu)=\sum_{i+j\le n}
F_i\tS\oc_\C\D_j\oc_\C\bM^\bu\subset \tS\oc_\C\D\oc_\C\bM^\bu$.
 The cone of the morphism $\Ups^\bu\Xi(\bM^\bu)\rarrow\bM^\bu$ has
an induced filtration $F$ whose components are the cones of
the morphisms $F_n\Ups^\bu\Xi(\bM^\bu)\rarrow\bM^\bu$.
 The quotient complex $\cone(F_n\Ups^\bu\Xi(\bM^\bu)\to\bM^\bu)/
\allowbreak \cone(F_{n-1}\Ups^\bu\Xi(\bM^\bu)\to\bM^\bu)$
is isomorphic to the cone of the identity endomorphism of $\bM^\bu$
for $n=0$ and to the cotensor product of a positive-degree component
of the Koszul complex $\S\oc_\C\D$ and the complex $\bM^\bu$ for $n>0$
(where, as always, $\S=\bigoplus_n F_n\tS/F_{n-1}\tS$).
 Thus in both cases the quotient complex is coacyclic.

 Fourthly, it remains to check that for any quasi-differential left 
$\tD$\comodule{} $\L$ the cone of the natural morphism of
quasi-differential $\tD$\comodule s $\L\rarrow\Xi\Ups^\bu(\L)$ is
coacyclic.
 First let us show that it suffices to consider the case when $\L$ 
is a graded $\C$\comodule{} endowed with a graded $\tD$\comodule{}
structure via the embedding of corings $\D_0\til\rarrow\tD$.
 To this end, consider the increasing filtration of $\L$ by
quasi-differential $\tD$\subcomodule s $G_n(\L)=\nu_L^{-1}
(\bigoplus_{i\le n}\tD_i\oc_\C\L)$, where $\nu_L\:\L\rarrow\tD\oc_\C\L$
denotes the coaction map.
 The quotient quasi-differential comodules $G_n(\L)/G_{n-1}(\L)$ are
are graded $\tD$\comodule s originating from graded $\C$\comodule s;
the filtration~$G$ induces a filtration on the cone of the morphism
$\L\rarrow\Xi\Ups^\bu(\L)$ whose components are the cones of
the morphisms $G_n(\L)\rarrow\Xi\Ups^\bu(G_n(\L))$; and the associated
quotient quasi-differential comodules of the latter filtration are
the cones of the morphisms $G_n(\L)/G_{n-1}(\L)\rarrow\Xi\Ups^\bu
(G_n(\L)/G_{n-1}(\L))$.

 Now assume that $\L$ is a graded $\C$\comodule{} with the induced
graded $\tD$\comodule{} structure, or even a complex of $\C$\comodule s
with the induced quasi-differential $\tD$\comodule{} structure.
 In this case, the quasi-differential $\tD$\comodule{} $\Xi\Ups^\bu(\L)
=\D\oc_\C\tS\oc_\C\L$ has an increasing filtration by
quasi-differential subcomodules given by the formula
$F_n\Xi\Ups^\bu(\L)=\sum_{i+j\le n} \D_j\oc_\C F_i\tS\oc_\C\L\subset
\D\oc_\C\tS\oc_\C\L$.
 The cone of the morphism $\L\rarrow\Xi\Ups^\bu(\L)$ has the induced
filtration $F$ whose components are the cones of the morphisms
$\L\rarrow F_n\Xi\Ups^\bu(\L)$.
 The associated quotient comodules of the latter filtration are
coacyclic quasi-differential $\tD$\comodule s.
 Indeed, the component $F_0\cone(\L\rarrow\Xi\Ups^\bu(\L))$ is
isomorphic to the cone of the identity endomorphism of $\L$, while
the quotient quasi-differential comodules with $n>0$ are isomorphic
to cotensor products of positive-degree components of the Koszul
complex $\D\oc_\C\S$ and the $\C$\comodule{} $\L$, endowed with
the quasi-differential $\tD$\comodule{} structures originating from
their structures of complexes of $\C$\comodule s.
 Thus all these quotient quasi-differential comodules are coacyclic.

 Part~(b): we will only construct a pair of adjoint functors
between the DG\category{} of complexes of right $\tS$\semimodule s
and the DG\category{} of quasi-differential right $\tD$\comodule s;
the rest of the proof is identical to that of part~(a).
 The functor $\Xi$ assigns to a right $\S$\semimodule{} $\bN$
the graded right $\D$\comodule{} $\bN\oc_\C\D$ endowed with
a right $\tD$\comodule{} structure in terms of the following
maps $\bN\oc_\C\D_n\rarrow\bN\oc_\C\D_n\til$.
 Consider the map $\bN\oc_\C\D_n\rarrow\bN\oc_\C\D_n$ equal to
the difference of the identity map and the composition
$\N\oc_\C\D_n\til\rarrow\N\oc_\C\D_1\til\oc_\C\D_{n-1}\til\rarrow
\N\oc_\C\D_{n-1}\til\rarrow\N\oc_\C\D_n\til$ of the map induced
by the comultiplication morphism, the map induced by the semiaction
morphism, and the map induced by the morphism~$\d_{n-1}$.
 This difference factorizes through the surjection
$\bN\oc_\C\D_n\til\rarrow\bN\oc_\C\D_n$, hence the desired map.
 The functor $\Ups$ assigns to a quasi-differential right
$\tD$\comodule{} $\R$ the complex of right $\tS$\semimodule s
$\Ups^\bu(\R)=\R\oc_\C\tS$ with the terms $\Ups^i(\R)=\R_{-i}\oc_\C\tS$
and the differentials $d^i$ defines as $(-1)^i$ times the composition
$\R_{-i}\oc_\C\tS\rarrow\R_{-i-1}\oc_\C\tD_1\oc_\C\tS\rarrow
\R_{-i-1}\oc_\C\tS$ of the map induced by the coaction morphism and
the map induced by the semimultiplication morphism.

 Part~(c): let us construct a pair of adjoint functors between
the DG\category{} of complexes of left $\tS$\semicontramodule s and
the DG\category{} of quasi-differential left $\tD$\contramodule s.
 The functor $\Xi$ assigns to a left $\S$\semicontramodule{} $\bP$
the graded left $\D$\contramodule{} $\Cohom_\C(\D,\bP)$ endowed with
a left $\tD$\contramodule{} structure in terms of the following
maps $\Cohom_\C(\D_n\til,\bP)\rarrow\Cohom_\C(\D_n,\bP)$.
 Consider the map $\Cohom_\C(\D_n\til,\bP)\rarrow\Cohom_\C(\D_n\til,\bP)$
equal to the difference of the identity map and the composition
$\Cohom_\C(\D_n\til,\bP)\rarrow\Cohom_\C(\D_{n-1}\til,\bP)\rarrow
\Cohom_\C(\D_{n-1}\til,\Cohom_\C(\D_1\til,\bP))\rarrow
\Cohom_\C(\D_n\til,\bP)$ of the map induced by the morphism~$\d_{n-1}$,
the map induced by the semicontraaction morphism $\bP\rarrow
\Cohom_\C(F_1\tS,\bP)$, and the map induced by the comultiplication
morphism $\D_n\til\rarrow\D_1\til\oc_\C\D_{n-1}\til$.
 This difference factorizes through the injection $\Cohom_\C(\D_n,\bP)
\rarrow\Cohom_\C(\D_n\til,\bP)$, hence the desired map.
 The functor $\Ups$ right adjoint to $\Xi$ assigns to
a quasi-differential left $\tD$\contramodule{} $\Q$ the complex
of left $\tS$\semicontramodule s $\Ups^\bu(\Q)=\Cohom_\C(\tS,\Q)$
with the terms $\Ups^i(\Q)=\Cohom_\C(\tS,\Q_{-i})$ and
the differential defined as the composition $\Cohom_\C(\tS,\Q)
\rarrow\Cohom_\C(F_1\tS\oc_\C\tS\;\Q)\rarrow\Cohom_\C(\tS,Q)$ of
the map induced by the semimultiplication morphism
$F_1\tS\oc_\C\tS\rarrow\tS$ and the map induced by the contraaction
morphism $\Cohom_\C(\D_1\til,\Q)\rarrow\Q$.

 The rest of the proof is analogous to that of part~(a), with
the exception of the argument related to the filtration~$G$
(the first step of the fourth part of the proof).
 The problem here is that the decreasing filtration $G$ of a graded
$\tD$\contramodule{} $\Q$ whose components are the images $G^n\Q$ of
the contraaction maps $\Cohom_\C(\tD/\bigoplus_{i\le n}\D_i\til\;\Q)
\rarrow\Q$ is not in general separated, i.~e., the intersection of
$G^n\Q$ may be nonzero (see Appendix~\ref{contramodules-appendix}).
 What one should do is replace an arbitrary quasi-differential left
$\tD$\contramodule{} $\Q$ with the total quasi-differential 
contramodule $\gR$ of its bar resolution $\dsb\rarrow
\Cohom_\C(\tD\oc_\C\tD\oc_\C\tD\;\Q)\rarrow\Cohom_\C(\tD\oc_\C\tD\;\Q)
\rarrow\Cohom_\C(\tD,\Q)$.
 Since the cone of the natural morphism of quasi-differential
$\tD$\contramodule s $\gR\rarrow\Q$ is contraacyclic, one can
consider the quasi-differential $\tD$\contramodule{} $\gR$
instead of~$\Q$.

 In addition to the filtration $G$ introduced above, consider also
the decreasing filtration ${}'\!G$ of a graded $\tD$\contramodule{}
$\Q$ whose components are the images ${}'\!G^n\Q$ of
the contraaction maps $\Cohom_\C(\D/\bigoplus_{i\le n}\D_i\;\Q)
\rarrow\Q$.
 It is clear that $\gR\simeq\plim_n\gR/{}'\!G^n\gR$.
 Next one can either show that ${}'\!G$ is a filtration by
graded $\tD$\+sub\-con\-tra\-mod\-ules and use the filtration ${}'\!G$
of $\gR$, or show that the filtrations $G$ and ${}'\!G$ are
commensurable, ${}'\!G^n\gR\subset G^n\gR\subset {}'\!G^{n-1}\gR$,
and use the filtration $G$ of $\gR$.
 (The quotient quasi-differential $\tD$\contramodule s
$G^n\gR/G^{n+1}\gR$ originate from graded $\C$\contramodule s, while
the quotient quasi-differential $\tD$\contramodule s
${}'\!G^n\gR/{}'\!G^{n+1}\gR$ originate from complexes of
$\C$\contramodule s, which is also sufficient.)
 Both assertions for an arbitrary graded $\tD$\contramodule{} follow
from the fact that the composition $\D_i\til\rarrow
\D_1\til\oc_\C\D_{i-1}\til\rarrow\D_1\til\oc_\C\D_{i-1}$ of
the comultiplication map and the map induced by the natural surjection
$\D_{i-1}\til\rarrow\D_{i-1}$ is injective and its cokernel, being
isomorphic to the cokernel of the comultiplication map $\D_i\rarrow
\D_1\oc_\C\D_{i-1}$, is a coprojective left $\C$\comodule.
 To check the latter, consider the composition of the map
$\D_i\til\rarrow\D_1\til\oc_\C\D_{i-1}$ in question with
the map $\d_{i-1}$.

 Alternatively, one can replace an arbitrary quasi-differential
$\tD$\contramodule{} $\Q$ with the cone of the morphism
$\ker(\Cohom_\C(\tD,\Q)\to\Q)\rarrow\Cohom_\C(\tD,\Q)$ and use
the appropriate generalization of Lemma~\ref{contra-nakayama-lemma}.3.
\end{proof}

\begin{rmk}
 Notice that no homological dimension condition on the $k$\+algebra
$A$ is assumed in the above Theorem.
 In particular, when $\C=A$, so $\tS$ is just a filtered $k$\+algebra,
Theorem provides a description of certain semiderived categories of
$\tS$\module s relative to $F_0\tS=A$.
 A description of the conventional derived category can also be
obtained.
 Namely, in the assumptions of part~(a) of Theorem the conventional
derived category of left $\tS$\semimodule s is equivalent to
the quotient category of the coderived category of quasi-differential
left $\tD$\comodule s by its minimal triangulated subcategory
containing all the quasi-differential $\tD$\comodule s originating
from acyclic complexes of left $\C$\comodule s and closed under
infinite direct sums.
 This is so because for any acyclic complex of $\tS$\semimodule s
$\bM^\bu$ the quasi-differential $\tD$\comodule s $F_n\Xi(\bM^\bu)/
F_{n-1}\Xi(\bM^\bu)$ originate from acyclic complexes of
$\C$\comodule s, and conversely, for any quasi-differential
$\tD$\comodule{} $\L$ originating from an acyclic complex of
$\C$\comodule s the complex of $\tS$\semimodule s $\Ups^\bu(\L)$
is acyclic.
 The analogous result holds for right $\tS$\semimodule s in
the assumptions of part~(b); and in the assumptions of part~(c)
the conventional derived category of left $\tS$\semicontramodule s
is equivalent to the quotient category of the contraderived category
of quasi-differential left $\tD$\contramodule s by its minimal
triangulated subcategory containing all the quasi-differential
contramodules originating from acyclic complexes of left
$\C$\contramodule s and closed under infinite products.
\end{rmk}

\subsection{SemiTor and Cotor, SemiExt and Coext}

\subsubsection{}
 Let $(\tD,\d)$ be a quasi-differential coring over a $k$\+algebra $A$;
assume that $\D=\tD/\im\d$ is a flat left and right $A$\module.

 Let $\N$ be a quasi-differential right $\tD$\comodule{} and $\M$
be a quasi-differential left $\tD$\comodule.
 Assume that one of the graded $A$\module s $\N$ and $\M$ is flat.
 Then on the cotensor product $\N\oc_\D\M$ of the graded comodules
$\N$ and $\M$ over the graded coring $\D$ there is a natural
differential with zero square, which is defined as follows.

 Consider the map $\delta\:\N\oc_\D\M\rarrow\N\oc_\D\tD\oc_\D\M$
given by the formula $x\oc y\mpsto -x_{(0)}\oc x_{(1)}\oc y +
x\oc y_{(-1)}\oc y_{(0)}$.
 This map factorizes through the injection $\N\oc_\D\M\rarrow \N\oc_\D
\tD\oc_\D\M$ given by the formula $x\oc y\mpsto (-1)^{|x_{(0)}|}
x_{(0)}\oc \d(x_{(1)})\oc y = (-1)^{|x|}x\oc \d(y_{(-1)})\oc y_{(0)}$,
hence the desired map $d\:\N\oc_\D\M\rarrow\N\oc_\D\M$.
 Let us check that $d^2=0$, that is the image of $d$ is contained in
$\N\oc_\tD\M$.
 Set $d(x\oc y)=x'\oc y'$.
 Consider the two elements $x'\oc y'_{(-1)}\oc y'_{(0)}$ and
$x'_{(0)}\oc x'_{(1)}\oc y$ of the cotensor product $\N\oc_\D\tD
\oc_\D\M$; we have to check that these two elements coincide.
 Consider the image of the former element under the map
$\N\oc_\D\tD\oc_\D\M\rarrow\N\oc_\D\tD\oc_\D\tD\oc_\D\M$ given by
the formula $u\oc b\oc v\mpsto (-1)^{|u_{(0)}|}u_{(0)}\oc\d(u_{(1)})
\oc b\oc v$ and the image of the latter element under the map
$\N\oc_\D\tD\oc_\D\M\rarrow\N\oc_\D\tD\oc_\D\tD\oc_\D\M$ given by
the formula $u\oc b\oc v\mpsto (-1)^{|u|+|b|}u\oc b\oc \d(v_{(-1)})
\oc v_{(0)}$.
 The sum of these two elements of $\N\oc_\D\tD\oc_\D\tD\oc_\D\M$ is
equal to the image of the element $\delta(x\oc y)$ under the map
$\N\oc_\D\tD\oc_\D\M\rarrow\N\oc_\D\tD\oc_\D\tD\oc_\D\M$ induced by
the comultiplication map $\tD\rarrow \tD\oc_\D\tD$.
 It remains to notice that the Carthesian square formed by the maps
$\tD\rarrow \tD\oplus\tD$, \ $\tD\oplus \tD\rarrow \tD\oc_\D\tD$, \
$\tD\rarrow\tD$, and $\tD\rarrow \tD\oc_\D\tD$ constructed
in~\ref{quasi-differential-comodules} remains Carthesian after
taking the cotensor product with $\N$ and~$\M$.
 We will denote the complex we have constructed by $\N\oc_\D^\bu\M$;
its terms are $\N\oc_\D^n\M=(\N\oc_\D\M)_{-n}$.

 Now assume that the ring $A$ has a finite weak homological dimension. 
 In order to define the double-sided derived functor of the functor
$\oc_\D^\bu$, we will show that the coderived category of quasi-differential
$\tD$\comodule s is equivalent to the quotient category of the homotopy
category of $\D$\+coflat quasi-differential $\tD$\comodule s by its
intersection with the thick subcategory of coacyclic quasi-differential
$\tD$\comodule s.

 The argument is analogous to that of either
Theorem~\ref{cotor-main-theorem} or Theorem~\ref{semitor-main-theorem}.
 First let us construct for any quasi-differential left $\tD$\comodule{}
$\K$ a morphism into it from an $A$\+flat quasi-differential left
$\tD$\comodule{} $\boL_1(\K)$ with a coacyclic cone.
 Use the graded version of Lemma~\ref{flat-comodule-surjection} to
obrain a finite resolution $0\rarrow \cZ\rarrow\cP_{d-1}(\K)\rarrow
\dsb\rarrow\cP_0(\K)\rarrow\K$ of a graded $\tD$\comodule{} $\K$
consisting of $A$\+flat graded $\tD$\comodule s.
 The total quasi-differential $\tD$\comodule{} of the complex of
quasi-differential $\tD$\comodule s $\cZ\rarrow\cP_{d-1}(\K)\rarrow
\dsb\rarrow\cP_0(\K)$ is an $A$\+flat quasi-differential
$\tD$\comodule{} whose morphism into $\K$ has a coacyclic cone.
 Indeed, the total quasi-differential $\tD$\comodule{} of any
acyclic complex of quasi-differential $\tD$\comodule s bounded
from below is coacyclic, since it has an increasing filtration
by quasi-differential $\tD$\subcomodule s such that the associated
quotient $\tD$\comodule s are isomorphic to cones of identity
endomorphisms of certain quasi-differential $\tD$\comodule s.

 Now let us construct for any $A$\+flat quasi-differential left
$\tD$\comodule{} $\L$ a morphism from it into a $\D$\+coflat
quasi-differential left $\tD$\comodule{} $\boR_2(\L)$ with
a coacyclic cone.
 Consider the cobar construction $\tD\ot_A\L\rarrow\tD\ot_A\tD\ot_A\L
\rarrow\dsb$\.
 Notice that $\tD$ is a coflat graded left $\D$\comodule, since
there is an exact triple of left $\D$\comodule s $\D(-1)\rarrow\tD
\rarrow\D$ (where $\D(-1)_i=\D_{i-1}$).
 Hence the total quasi-differential $\tD$\comodule{} of this cobar
complex of quasi-differential $\tD$\comodule s is a $\D$\+coflat
quasi-differential $\tD$\comodule{} such that the map into it from
the quasi-differential $\tD$\comodule{} $\L$ has a coacyclic cone.

 It is easy to see that the cotensor product of a quasi-differential
right $\tD$\comodule{} and a quasi-differential left $\tD$\comodule{}
is an acyclic complex whenever one of the two quasi-differential
$\tD$\comodule s is coacyclic and the other one is $\D$\+coflat.
 The derived functor $\Cotor^{\tD}_\qq$ on the Carthesian product
of coderived categories of right and left quasi-differential
$\tD$\comodule s is defined by restricting the functor $\oc_\D^\bu$ to
the Carthesian product of the homotopy category of quasi-differential
right $\tD$\comodule s and the homotopy category of $\D$\+coflat
quasi-differential left $\tD$\comodule s or to the Carthesian product
of the homotopy category of $\D$\+coflat quasi-differential right
$\tD$\comodule s and the homotopy category of quasi-differential
left $\tD$\comodule s, and composing it with the localization
functor $\Hot(k\modl)\rarrow\sD(k\modl)$.

\subsubsection{}
 Let $(\tD,\d)$ be a quasi-differential coring over a $k$\+algebra $A$;
assume that $\D=\tD/\im\d$ is a projective left and a flat right
$A$\module.

 Let $\M$ be a quasi-differential left $\tD$\comodule{} and $\P$ be
a quasi-differential left $\tD$\contramodule.
 Assume that either the graded $A$\module{} $\M$ is projective, or
the graded $A$\module{} $\P$ is injective.
 Then on the graded $k$\module{} of cohomomorphisms $\Cohom_\D(\M,\P)$
from the graded comodule $\M$ to the graded contramodule $\P$ over
the graded coring $\D$ there is a natural differential with zero square,
which is defined as follows.
 Consider the map $\delta\:\Cohom_\D(\M,\Cohom_\D(\tD,\P))\simeq
\Cohom_\D(\tD\oc_\D\M\;\P)\rarrow\Cohom_\D(\M,\P)$ defined by
the formula $f\mpsto \pi_\P\circ f-f\circ\nu_\M$ (where $\pi_\P$ and
$\nu_\M$ denote the contraaction and coaction morphisms).
 This map factorizes through the surjection
$\Cohom_\D(\tD\oc_\D\M\;\P)\rarrow\Cohom_\D(\M\;\P)$ induced
by the morphism $\dbar\:\D\rarrow\tD$, hence the map
$d\:\Cohom_\D(\M,\P)\rarrow\Cohom_\D(\M,\P)$.
 We will denote the complex we have constructed by
$\Cohom_\D^\bu(\M,\P)$; its terms are $\Cohom_\D^n(\M,\P)=
\Cohom_\D(\M,\P)_{-n}$.

 Assume that the ring $A$ has a finite left homological dimension.
 Then the coderived category of quasi-differential left $\tD$\comodule s
is equivalent to the quotient category of the homotopy category of
$\D$\coprojective{} quasi-differential left $\tD$\comodule s by its
intersection with the thick subcategory of coacyclic quasi-differential
$\tD$\comodule s.
 Analogously, the contraderived category of quasi-differential left
$\tD$\contramodule s is equivalent to the quotient category of
the homotopy category of $\D$\coinjective{} quasi-differential left
$\tD$\contramodule s by its intersection with the thick subcategory of
contraacyclic quasi-differential $\tD$\contramodule s.
 The double-sided derived functor $\Coext_{\tD}^\qq$ on the Carthesian
product of the coderived category of quasi-differential left
$\tD$\comodule s and the contraderived category of quasi-differential
left $\tD$\contramodule s is defined by restricting the functor
$\Cohom_\D^\bu$ to the Carthesian product of the homotopy category of
$\D$\coprojective{} quasi-differential left $\tD$\comodule s and
the homotopy category of quasi-differential left $\tD$\contramodule s
or to the Carthesian product of the homotopy category of
quasi-differential left $\tD$\comodule s and the homotopy category
of $\D$\coinjective{} quasi-differential left $\tD$\contramodule s,
and composing in with the localization functor $\Hot(k\modl)\rarrow
\sD(k\modl)$.

\subsubsection{}
 Let $\C$ be a coring over a $k$\+algebra $A$.
 Assume that $\C$ is a flat left and right $A$\module{} and $A$ has
a finite weak homological dimension.
 Let $\tS$ be a left and right coflat nonhomogeneous Koszul semialgebra
over~$\C$, and $\tD$ be the left and right coflat Koszul
quasi-differential coring nonhomogeneous quadratic dual to~$\tS$.

\begin{cor}
 \textup{(a)} The equivalences of categories\/ $\sD^\si(\simodr\tS)\simeq
\sD^\co(\qcomodr\tD)$ and\/ $\sD^\si(\tS\simodl)\simeq\sD^\co(\tD\qcomodl)$
transform the derived functor\/ $\SemiTor^{\tS}$ into the derived
functor\/ $\Cotor^{\tD}_\qq$. \par
 \textup{(b)} Assume additionally that\/ $\C$ is a projective left
$A$\module, $A$ has a finite left homological dimension, and\/ $\tS$ is
a left coprojective nohnomogeneous Koszul semialgebra.
 Then the equivalences of categories\/ $\sD^\si(\tS\simodl)\simeq
\sD^\co(\tD\qcomodl)$ and\/ $\sD^\si(\tS\sicntr)\simeq
\sD^\ctr(\tD\qcontra)$ transform the derived functor\/ $\SemiExt_{\tS}$
into the derived functor\/ $\Coext_{\tD}^\qq$.
\end{cor}

\begin{proof}
 Part~(a): for any complex of right $\tS$\semimodule s $\bN^\bu$ and any
quasi-differential left $\tD$\comodule{} $\L$ there is a natural
isomorphism of complexes of $k$\module s $\Xi(\bN^\bu)\oc_\D^\bu\L
\simeq \bN^\bu\os_{\tS}\Ups^\bu(\L)$.
 Indeed, both complexes are isomorphic to the total complex of
the bicomplex $\bN^i\oc_\C\L_j$, one of whose differentials is induced
by the differential in $\bN^\bu$ and the other is equal to
the composition $\bN^i\oc_\C\L_j\rarrow\bN^i\oc_\C\D_1\til\oc_\C\L_{j-1}
\rarrow\bN^i\oc_\C\L_{j-1}$ of the map induced by the $\tD$\+coaction
in $\L$ and the map induced by the $\tS$\+semiaction in $\bN^i$.
 Now let $\bN^\bu$ be a semiflat complex of $\C$\+coflat right
$\tS$\semimodule s and $\bM^\bu$ be a complex of left $\tS$\semimodule s.
 Then there is an isomorphism $\Xi(\bN^\bu)\oc_\D^\bu\Xi(\bM^\bu)\simeq
\bN^\bu\os_{\tS}\Ups^\bu\Xi(\bM^\bu)$ and a quasi-isomorphism
$\bN^\bu\os_{\tS}\Ups^\bu\Xi(\bM^\bu)\rarrow\bN^\bu\os_\tS\bM^\bu$.
 Analogously, for a complex of right $\tS$\semimodule s $\bN^\bu$ and
a semiflat complex of $\C$\+coflat left $\tS$\semimodule s $\bM^\bu$
there is an isomorphism $\Xi(\bN^\bu)\oc_\D^\bu\Xi(\bM^\bu)\simeq
\Ups^\bu\Xi(\bN^\bu)\os_{\tS}\bM^\bu$ and a quasi-isomorphism
$\Ups^\bu\Xi(\bN^\bu)\os_{\tS}\bM^\bu\rarrow\bN^\bu\os_{\tS}\bM^\bu$.
 It is easy to check that the square diagram formed by these maps
is commutative.
 The proof of part~(b) is completely analogous.
\end{proof}

\begin{qst}
 Can one construct a comodule-contramodule correspondence
(equivalence between the coderived and contraderived categories)
for quasi-differential comodules and contramodules?
 Also, is there a natural closed model category structure on
the category of quasi-differential comodules (contramodules)?
\end{qst}

\appendix
\Section{Contramodules over Coalgebras over Fields}
\label{contramodules-appendix}

 Let $\C$ be a coassociative coalgebra with counit over a field~$k$.
 It is well-known~\cite{Swe} that $\C$ is the union of its
finite-dimensional subcoalgebras and any $\C$\comodule{} is
a union of finite-dimensional comodules over finite-dimensional
subcoalgebras of~$\C$.
 The dual assertion for $\C$\contramodule s is \emph{not} true:
for the most common of nonsemisimple infinite-dimensional coalgebras
$\C$ there exist $\C$\contramodule s $\P$ such that the intersection of
the images of $\Hom_\C(\C/\cU,\P)$ in $\P$ over all finite-dimensional
subcoalgebras $\cU\subset\C$ is nonzero.
 A weaker statement holds, however: if the map $\Hom_\C(\C/\cU,\P)
\rarrow\P$ is surjective for any finite-dimensional subcoalgebra $\cU$
of~$\C$, then $\P=0$.
 Besides, even though adic filtrations of contramodules are not in
general separated, they are always \emph{complete}.
 Using the related techniques we show that any contraflat
$\C$\contramodule{} is projective, generalizing the well-known result
that any flat module over a finite-dimensional algebra is
projective~\cite{Bas}.

\subsection{Counterexamples}

\subsubsection{} Let $\C$ be the coalgebra for which the dual algebra
$\C^*$ is isomorphic to the algebra of formal power series $k[[x]]$.
 Then a $\C$\contramodule{} $\P$ can be equivalently defined as
a $k$\+vector space endowed with the following operation of summation
of sequences of vectors with formal coefficients~$x^n$: for any
elements $p_0$, $p_1$,~\dots\ in~$\P$, an element of~$\P$ denoted by
$\sum_{n=0}^\infty x^np_n$ is defined.
 This operation should satisfy the following equations:
$\sum_{n=0}^\infty x^n(ap_n+bq_n)=a\sum_{n=0}^\infty x^np_n+
b\sum_{n=0}^\infty x^nq_n$ for $a$,~$b\in k$, \ $p_n$,~$q_n\in\P$
(linearity); $\sum_{n=0}^\infty x^np_n=p_0$ when $p_1=p_2=\dsb=0$
(counity); and $\sum_{i=0}^\infty x^i \bigl(\sum_{j=0}^\infty x^jp_{ij}
\bigr) = \sum_{n=0}^\infty x^n \bigl(\sum_{i+j=n}p_{ij}\bigl)$
for any $p_{ij}\in\P$, $i,j=0,1,\dsc$ (contraassociativity).
 Here the interior summation sign in the right hand side denotes
the conventional finite sum of elements of a vector space, while
the three other summation signs refer to the contramodule infinite
summation operation.

 The following examples of $\C$\contramodule s are revealing.
 Let $\gE$ denote the free $\C$\contramodule{} generated by
the sequence of symbols $e_0$, $e_1$,~\dots; its elements can be
represented as formal sums $\sum_{i=0}^\infty a_i(x)e_i$, where $a_i(x)$
are formal power series in~$x$ such that the sequence of their orders
of zero $\ord_x a_i(x)$ at $x=0$ tends to infinity as $i$ increases.
 Let $\gF$ denote the free $\C$\contramodule{} generated by
the sequence of symbols $f_1$, $f_2$,~\dots; then $\C$\contramodule{}
homomorphisms from $\gF$ to $\gE$ correspond bijectively to
sequences of elements of $\gE$ that are images of the elements~$f_i$.
 We are interested in the map $g\:\gF\rarrow\gE$ sending $f_i$ to
$x^ie_i-e_0$; in other words, an element $\sum_{i=1}^\infty b_i(x)f_i$
of $\gF$ is mapped to the element $\sum_{i=1}^\infty x^ib_i(x)e_i -
\bigl(\sum_{i=1}^\infty b_i(x)\bigl)e_0$.
 It is clear from this formula that the element $e_0\in\gE$ does not
belong to the image of~$g$.
 Let $\P$ denote the cokernel of the morphism~$g$ and $p_i$ denote
the images of the elements~$e_i$ in~$\P$.
 Then one has $p_0=x^np_n$ in~$\P$; in other words, the element $p_0$
belongs to the image of $\Hom_\C(\C/\cU,\P)$ under the contraaction
map $\Hom_\C(\C,\P)\rarrow\P$ for any finite-dimensional subcoalgebra
$\cU = (k[[x]]/x^n)^*$ of~$\C$.

 Now let $\gE'$ be the free $\C$\contramodule{} generated by
the symbols $e_1$, $e_2$,~\dots, \ $\P'$ denote the cokernel of the map
$g'\:\gF\rarrow\gE'$ sending $f_i$ to~$x_ie_i$, and $p'_i$ denote
the images of~$e'_i$ in $\P'$.
 Then the result of the contramodule infinite summation
$\sum_{n=1}^\infty x^np'_n$ is nonzero in~$\P'$, even though every
element $x^np'_n$ is equal to zero.
 Therefore, \emph{the contramodule summation operation cannot be
understood as any kind of limit of finite partial sums}.
 Actually, the $\C$\contramodule{} $\P'$ is just the direct sum of
the contramodules $k[[x]]/x^nk[[x]]$ over $n=1$, $2$,~\dots\
in the category of $\C$\contramodule s.
 Notice that the element $\sum_{n=1}^\infty x^np'_n$ also belongs to
the image of $\Hom_\C(\C/\cU,\P')$ in $\P'$ for any 
finite-dimensional subcoalgebra $\cU\subset\C$.

\begin{rmk}
 In the above notation, a $\C$\contramodule{} structure on
a $k$\+vector space $\P$ is uniquely determined by the underlying
structure of a module over the algebra of polynomials $k[x]$;
the natural functor $\C\contra\rarrow k[x]\modl$ is fully faithful.
 Indeed, for any $p_0$, $p_1$,~\dots\ in~$\P$ the sequence
$q_m=\sum_{n=0}^\infty x^np_{m+n}\in\P$, \ $m=0$, $1$,~\dots\ is
the unique solution of the system of equations $q_m = p_m + xq_{m+1}$.
 The image of this functor is a full abelian subcategory closed under
kernels, cokernels, extensions, and infinite products; it consists
of all $k[x]$\module s $P$ such that $\Ext^i_{k[x]}(k[x,x^{-1}],P)=0$
for $i=0$,~$1$.
 It follows that if $\D$ is a coalgebra for which the dual algebra
$\D^*$ is isomorphic to a quotient algebra of the algebra of formal
power series $k[[x_1,\dsc,x_m]]$ in a finite number of (commuting)
variables by a closed ideal, then the natural functor $\D\contra
\rarrow k[x_1,\dsc,x_m]\modl$ is fully faithful.
\end{rmk}

\subsubsection{}
 Now let us give an example of \emph{finite-dimensional} (namely,
two-dimensional) contramodule $\P$ over a coalgebra $\C$ such that
the intersection of the images of $\Hom_\C(\C/\cU,\P)$ in $\P$ is
nonzero.
 Notice that for any coalgebra $\C$ there are natural left
$\C^*$\module{} structures on any left $\C$\comodule{} and any left
$\C$\contramodule; that is there are natural faithful functors
$\C\comodl\rarrow\C^*\modl$ and $\C\contra\rarrow\C^*\modl$
(where $\C^*$ is considered as an abstract algebra without any
topology).
 The functor $\C\comodl\rarrow\C^*\modl$ is fully faithful, while
the functor $\C\contra\rarrow\C^*\modl$ is fully faithful on
finite-dimensional contramodules.

 Let $V$ be a vector space and $\C$ be the coalgebra such that
the dual algebra $\C^*$ has the form $ki_2\oplus i_2V^*i_1\oplus ki_1$,
where $i_1$ and $i_2$ are idempotent elements with $i_1i_2=i_2i_1=0$
and $i_1+i_2=1$.
 Then left $\C^*$\module s are essentially pairs of $k$\+vector spaces
$M_1$, $M_2$ endowed with a map $V^*\ot_k M_1\rarrow M_2$, left
$\C$\comodule s are pairs of vector spaces $\M_1$, $\M_2$ endowed with
a map $\M_1\rarrow V\ot_k\M_2$, and left $\C$\contramodule s are pairs
of vector spaces $\P_1$, $\P_2$ endowed with a map $\Hom_k(V,\P_1)
\rarrow\P_2$.
 In particular, the functor $\C\contra\rarrow\C^*\modl$ is not
surjective on morphisms of infinite-dimensional objects, while
the functor $\C\comodl\rarrow\C^*\modl$ is not surjective on
the isomorphism classes of finite-dimensional objects.
 (Neither is in general the functor $\C\contra\rarrow\C^*\modl$, as one
can see in the example of an analogous coalgebra with three idempotent
linear functions instead of two and three vector spaces instead of one;
when $k$ is a finite field and $\C$ is the countable direct sum of
copies of the coalgebra~$k$, there even exists a one-dimensional
$\C^*$\module{} which comes from no $\C$\comodule{} or
$\C$\contramodule.)

 Let $\P$ be the $\C$\contramodule{} with $\P_1=k=\P_2$ corresponding
to a linear function $V^*\rarrow k$ coming from no element of~$V$.
 Then the intersection of the images of $\Hom_\C(\C/\cU,\P)$ in $\P$
over all finite-dimensional subcoalgebras $\cU\subset\C$ is equal
to~$\P_2$.

 More generally, for any coalgebra $\C$ any finite-dimensional left
$\C$\comodule{} $\M$ has a natural left $\C$\contramodule{} structure
given by the composition $\Hom_k(\C,\M)\simeq\C^*\ot_k\M \rarrow
\C^*\ot_k\C\ot_k\M\rarrow\M$ of the map induced by the $\C$\+coaction
in $\M$ and the map induced by the pairing $\C^*\ot_k\C\rarrow k$.
 The category of finite-dimensional left $\C$\comodule s is isomorphic
to a full subcategory of the category of finite-dimensional left
$\C$\contramodule s; a finite-dimensional $\C$\contramodule{} comes
from a $\C$\comodule{} if and only if it comes from a contramodule
over a finite-dimensional subcoalgebra of~$\C$.
 We will see below that every \emph{irreducible} $\C$\contramodule{}
is a finite-dimensional contramodule over a finite-dimensional
subcoalgebra of~$\C$; it follows that the above functor provides
a bijective correspondence between irreducible left $\C$\comodule s
and irreducible left $\C$\contramodule s.

 Compairing the cobar complex for comodules with the bar complex for
contramodules, one discovers that for any finite-dimensional left
$\C$\comodule s $\L$ and $\M$ there is a natural isomorphism
$\Ext^{\C,i}(\L,\M)\simeq\Ext_\C^i(\L,\M)^{{*}{*}}$.
 In other words, the $\Ext$ spaces between finite-dimensional
$\C$\comodule s in the category of arbitrary $\C$\contramodule s
are the completions of the $\Ext$ spaces in the category of
finite-dimensional $\C$\comodule s with respect to
the profinite-dimensional topology.

\subsection{Nakayama's Lemma}  \label{contra-nakayama-lemma}
 A coalgebra is called \emph{cosimple} if it has no nontrivial proper
subcoalgebras.
 A coalgebra $\C$ is called \emph{cosemisimple} if it is a union of
finite-dimensional coalgebras dual to semisimple $k$\+algebras, or
equivalently, if abelian the category of (left or right)
$\C$\comodule s is semisimple.
 Any cosemisimple coalgebra can be decomposed into an (infinite)
direct sum of cosimple coalgebras in a unique way.
 For any coalgebra~$\C$, let $\C^\ss$ denote its maximal cosemisimple
subcoalgebra; it contains all other cosemisimple subcoalgebras of~$\C$.

\begin{lem1}
 Let\/ $\C$ be a coalgebra over a field~$k$ and\/ $\P$ be a nonzero
left\/ $\C$\contramodule.
 Then the image of the space\/ $\Hom_k(\C/\C^\ss,\P)$ under
the contraaction map\/ $\Hom_k(\C,\P)\rarrow\P$ is not equal to\/~$\P$.
\end{lem1}

\begin{proof}
 Notice that the coalgebra without counit $\D=\C/\C^\ss$ is
\emph{conilpotent}, that is any element of $\D$ is annihilated by
the iterated comultiplication map $\D\rarrow\D^{\ot n}$ with
a large enough~$n$.
 We will show that for any contramodule $\P$ over a conilpotent
coalgebra $\D$ surjectivity of the map $\Hom_k(\D,\P)\rarrow\P$ implies
vanishing of~$\P$.
 Indeed, assume that the contraaction map $\pi_\P$ of
a $\D$\contramodule{} $\P$ is surjective.
 Let $p$ be an element of~$\P$; it is equal to $\pi_\P(f_1)$ for
a certain map $f_1\:\D\rarrow\P$.
 Since the map $\pi_\P$ is surjective, the map $f_1$ can be lifted
to a certain map $\D\rarrow\Hom_k(\D,\P)$, which supplies a map
$f_2\:\D\ot_k\D\rarrow\P$.
 So one constructs a sequence of maps $f_i\:\D^{\ot i}\rarrow\P$
such that $f_{i-1}=\pi_{\P,1}(f_i)$, where $\pi_{\P,1}$ signifies
the application of~$\pi_\P$ at the first tensor factor of $\D^{\ot i}$.
 Set $g_i=\mu_{\D,2..i}(f_i)=f_i\circ\mu_{\D,2..i}$, where
$\mu_{\D,2..i}$ denotes the tensor product of the identity map
$\D\rarrow\D$ with the iterated comultiplication map
$\D\rarrow\D^{\ot i-1}$.
 Then $g_i$ is a map $\D\otimes\D\rarrow\P$ defined for each $i\ge2$.
 We have $\pi_{\P,1}(g_i)=\mu_{\D,1..i-1}(f_{i-1})$ and
$\mu_\D(g_i)=\mu_{\D,1..i}(f_i)$.
 Notice that by conilpotency of the coalgebra $\D$ the series
$\sum_{i=2}^\infty g_i$ converges in the sense of point-wise limit
of functions $\D\ot_k\D\rarrow\P$, and even of functions
$\D\rarrow\Hom_k(\D,\P)$.
 (As always, we presume the identification $\Hom_k(U\ot_k V\; W)
=\Hom_k(V,\Hom_k(U,W))$ when we consider left contramodules.)
 Therefore, $\pi_{\P,1}\bigl(\sum_{i=2}^\infty g_i\bigr)=
\sum_{i=2}^\infty \mu_{\D,1..i-1}(f_{i-1})$ and
$\mu_\D\bigl(\sum_{i=2}^\infty g_i\bigr)=
\sum_{i=2}^\infty\mu_{\D,1..i}(f_i)$, hence
$\pi_{\P,1}\bigl(\sum_{i=2}^\infty g_i\bigr)-
\mu_\D\bigl(\sum_{i=2}^\infty g_i\bigr) = f_1$.
 By the contraassociativity equation, it follows that
$p=\pi_\P(f_1)=0$.
\end{proof}

\begin{lem2}
 Let coalgebra\/ $\C$ be the direct sum of a family of coalgebras\/
$\C_\alpha$.
 Then any left contramodule\/ $\P$ over\/ $\C$ is the product of
a uniquely defined family of left contramodules\/ $\P_\alpha$
over\/~$\C_\alpha$.
\end{lem2}

\begin{proof}
 Uniqueness and functoriality is clear, since the component $\P_\alpha$
can be recovered as the image of the projector corresponding to
the linear function on~$\C$ that is equal to the counit on $\C_\alpha$
and vanishes on $\C_\beta$ for all $\beta\ne\alpha$.
 Existence is obvious for a free $\C$\contramodule.
 Now suppose that a $\C$\contramodule{} $\Q$ is the product of
$\C_\alpha$\contramodule s $\Q_\alpha$; let us show that any
subcontramodule $\gR\subset\Q$ is the product of its images
$\gR_\alpha$ under the projections $\Q\rarrow\Q_\alpha$.
 Let $r_\alpha$ be a family of elements of $\gR$.
 Consider the linear map $f\:\C\rarrow\gR$ whose restriction to
$\C_\alpha$ is equal to the composition $\C_\alpha\rarrow k\rarrow\gR$
of the counit of~$\C_\alpha$ and the map sending $1\in k$ to $r_\alpha$.
 Set $r=\pi_\gR(f)$.
 Then the image of the element~$r$ under the projection $\gR\rarrow
\gR_\alpha$ is equal to the image of~$r_\alpha$ under this projection.
 Thus $\gR$ is identified with the product of $\gR_\alpha$.
 It remains to notice that any $\C$\contramodule{} is isomorphic to
the quotient contramodule of a free contramodule by one of
its subcontramodules.
\end{proof}

\begin{cor}
 For any coalgebra\/ $\C$ and any nonzero contramodule\/ $\P$
over\/~$\C$ there exists a finite-dimensional (and even cosimple)
subcoalgebra\/ $\cU\subset\C$ such that the image of\/
$\Hom_k(\C/\cU,\P)$ under the contraaction map\/ $\Hom_k(\C,\P)
\rarrow\P$ is not equal to\/~$\P$.
\end{cor}

\begin{proof}
 By Lemma~1, the image of the map $\Hom_k(\C/\C^\ss,\P)\rarrow\P$
is not equal to $\P$.
 Denote this image by~$\Q$; it is a subcontramodule of $\P$ and
the quotient contramodule $\P/\Q$ is a contramodule over $\C^\ss$.
 By Lemma~2, there exists a cosimple subcoalgebra $\C_\alpha$
of $\C^\ss$ such that $\P/\Q$ has a nonzero quotient which
is a contramodule over $\C_\alpha$.
\end{proof}

\begin{lem3}
 Let\/ $\C_0\subset\C_1\subset\C_2\subset\dsb\subset\C$ be a coalgebra
with an increasing sequence of subcoalgebras.
 For a\/ $\C$\contramodule\/ $\P$, denote by $G^i\P$ the image of
the contraaction map\/ $\Hom_k(\C/\C_i,\P)\rarrow\P$.
 Then for any\/ $\C$\contramodule\/ $\P$ the natural map
$\P\rarrow\plim_i\P/G^i\P$ is surjective.
\end{lem3}

\begin{proof}
 The assertion is obvious for a free $\C$\contramodule.
 Represent an arbitrary $\C$\contramodule{} $\P$ as the quotient
contramodule $\Q/\gK$ of a free $\C$\contramodule{} $\Q$.
 Since the maps $G^i\Q\rarrow G^i\P$ are surjective, there are
short exact sequences $0\rarrow\gK/\gK\cap G^i\Q\rarrow\Q/G^i\Q
\rarrow\P/G^i\P\rarrow0$.
 Passing to the projective limits, we see that the map
$\plim_i\Q/G^i\Q\rarrow\plim_i\P/G^i\P$ is surjective.
\end{proof}

 When $\C=\bigcup_i\C_i$, it follows, in particular, that
$\gK\simeq\plim_i\gK/G^i\gK$ for any $\C$\contramodule{} $\gK$
which is a subcontramodule of a projective $\C$\contramodule.

\subsection{Contraflat contramodules}  \label{contraflat-contramodules}

\begin{lem}
 Let\/ $\C$ be a coalgebra over a field $k$.  Then a left\/
$\C$\contramodule{} is contraflat if and only if it is projective.
\end{lem}

\begin{proof}
 For any $\C$\contramodule{} $\Q$ and any subcoalgebra $\cV\subset\C$
denote by ${}^\cV\Q=\coker(\Hom_\C(\C/\cV,\Q)\to\Q)\simeq
\Cohom_\C(\cV,\Q)$ the maximal quotient contramodule of~$\Q$ that is
a contramodule over~$\cV$.
 The key step is to construct for any $\C^\ss$\contramodule{} $\gR$
a projective $\C$\contramodule{} $\Q$ such that ${}^{\C^\ss}\Q
\simeq\gR$.
 By Lemma~\ref{contra-nakayama-lemma}.2, $\gR$ is a product of
contramodules over cosimple components $\C_\alpha$ of $\C^\ss$.
 Any contramodule over $\C_\alpha$ is, in turn, a direct sum of
copies of the unique irreducible $\C_\alpha$\contramodule.
 Hence it suffices to consider the case of an irreducible
$\C_\alpha$\contramodule{} $\gR$.
 Let $e_\alpha$ be an idempotent element of the algebra $\C_\alpha^*$
such that $\gR$ is isomorphic to $\C_\alpha^*e_\alpha$.
 Consider the idempotent linear function $e_\ss$ on $\C^\ss$ equal
to $e_\alpha$ on $\C_\alpha$ and zero on $\C_\beta$ for all
$\beta\ne\alpha$.
 It is well-known that for any surjective map of rings $A\rarrow B$
whose kernel is a nil ideal in~$A$ any idempotent element of~$B$
can be lifted to an idempotent element of~$A$.
 Using this fact for finite-dimensional algebras and Zorn's Lemma,
one can show that any idempotent linear function on $\C^\ss$ can be
extended to an idempotent linear function on~$\C$.
 Let $e$ be an idempotent linear function on~$\C$ extending $e_\ss$;
set $\Q=\C^*e$.
 Then one has ${}^{\C^\ss}\Q\simeq({}^{\C^\ss}\C^*)e\simeq
\C^{\ss*}e_\ss \simeq\gR$ as desired.
 Now let $\P$ be a contraflat left $\C$\contramodule.
 Consider a projective left $\C$\contramodule{} $\Q$ such that
${}^{\C^\ss}\Q\simeq{}^{\C^\ss}\P$.
 Since $\Q$ is projective, the map $\Q\rarrow{}^{\C^\ss}\P$ can be
lifted to a $\C$\contramodule{} morphism $f\:\Q\rarrow\P$.
 Since ${}^{\C^\ss}(\coker f)=\coker({}^{\C^\ss}\!\.f)=0$, it
follows from Lemma~\ref{contra-nakayama-lemma}.1 that the morphism
$f$ is surjective; it remains to show that $f$ is injective.
 For any right comodule $\N$ over a subcoalgebra $\cU\subset\C$
there is a natural isomorphism $\N\ocn_\C\P\simeq\N\ocn_\cU{}^\cU\P$,
hence the $\cU$\contramodule{} ${}^\cU\P$ is contraflat.
 Now let $\cU$ be a finite-dimensional subcoalgebra; then
${}^\cU\P$ is a flat left $\cU^*$\module.
 Denote by $K$ the kernel of the map ${}^\cU\!\.f\:{}^\cU\Q\rarrow
{}^\cU\P$.
 For any right $\cU^*$\module{} $N$ we have a short exact sequence
$0\rarrow N\ot_{\cU^*}K\rarrow N\ot_{\cU^*}{}^\cU\Q\rarrow N\ot_{\cU^*}
{}^\cU\P\rarrow0$.
 Since for any cosimple subcoalgebra $\cU_\alpha\subset\cU$ the map
$\cU_\alpha^*\ot_{\cU^*}{}^\cU\!\.f = {}^{\cU_\alpha}\!f$ is
an isomorphism, we can conclude that the module $\cU_\alpha^*
\ot_{\cU^*}K={}^{\cU_\alpha^*}K$ is zero.
 It follows that $K=0$ and the map ${}^\cU\!\.f$ is an isomorphism.
 Finally, let $\gK$ be the kernel of the map $\Q\rarrow\P$.
 Since ${}^\cU\!\.f$ is an isomorphism, the subcontramodule
$\gK\subset\Q$ is contained in the image of $\Hom_k(\C/\cU,\Q)$
in~$\Q$ for any finite-dimensional subcoalgebra $\cU\subset\C$.
 But the intersection of such images is zero, because the $\Q$ is
a projective $\C$\contramodule.
\end{proof}

\begin{rmk}
 Much more generally, one can define left contramodules over
an arbitrary complete and separated topological ring $R$ where open
right ideals form a base of neighborhoods of zero (cf.~\cite{Beil}).
 Namely, for any set $X$ let $R[X]$ denote the set of all formal
linear combinations of elements of $X$ with coefficients coverging
to zero in~$R$, i.~e., the set of all formal sums $\sum_{x\in X} r_xx$,
with $r_x\in R$ such that for any neighborhood of zero $U\subset R$
one has $r_x\in U$ for all but a finite number of elements $x\in X$.
 Then for any set $X$ there is a natural map of ``opening
the parentheses'' $R[R[X]]\rarrow R[X]$ assigning a formal linear
combination to a formal linear combination of formal linear
combinations; it is well-defined in view of our condition on~$R$.
 There is also a map $X\rarrow R[X]$ defined in terms of the unit
element of $R$; taken together, these two natural maps make
the functor $X\mpsto R[X]$ a monad on the category of sets.
 Left contramodules over~$R$ are, by the definition, modules over
this monad.
 One can see that the category of left $R$\contramodule s is abelian
and there is a forgetful functor from it to the category of
$R$\module s; this functor is exact and preserves infinite products.
 For example, when $R=\boZ_l$ is the topological ring of $l$\+adic 
integers, the category of $R$\contramodule s is isomorphic to
the category of weakly $l$\+complete abelian groups in the sense
of Jannsen~\cite{Jan}, i.~e., abelian groups $P$ such that
$\Ext_\boZ^i(\boZ[l^{-1}],P)=0$ for $i=0$,~$1$.
 When $R$ is a topological algebra over a field, the above definition
of an $R$\contramodule{} is equivalent to the definition given
in~\ref{assoc-contramodules}.
 Now if $T$ is a topological ring without unit satisfying the above
condition, and $T$ is pronilpotent, that is for any neiborhood of zero
$U\subset T$ there exists $n$ such that $T^n\subset U$, then any
left $T$\contramodule{} $P$ such that the contraaction map
$T[P]\rarrow P$ is surjective vanishes.
 Besides, any left contramodule over the topological product $R$ of
a family of rings (with units) $R_\alpha$ satisfying the above condition
is naturally the product of a family of left $R_\alpha$\contramodule s.
 Finally, let $R$ be a topological ring satisfying the above condition
and endowed with a decreasing filtration $R\supset G^1R\supset G^2R
\supset\dsb$ by closed left ideals such that any neighborhood of zero
in $R$ contains $G^iR$ for large~$i$.
 For any left $R$\contramodule{} $P$, denote by $G^iP$ the image of
the contraaction map $G^iR[P]\rarrow P$; then the natural map
$P\rarrow\plim_i P/G^iP$ is surjective.
 The proofs of these assertions are analogous to those of
Lemmas~\ref{contra-nakayama-lemma}.1--\ref{contra-nakayama-lemma}.3.
 When open two-sided ideals form a base of neighborhoods of zero in $R$
and the discrete quotient rings of $R$ are right Artinian, a left
$R$\contramodule{} $P$ is projective if and only if for any open
two-sided ideal $J\subset R$ the cokernel of the contraaction map
$J[P]\rarrow P$ is a projective left $R/J$\+module.
 The proof of this result is analogous to that of
Lemma~\ref{contraflat-contramodules}.
 It follows that the class of projective left $R$\contramodule s is
closed under infinite products under these assumptions.
 For a profinite ring $R$ (defined equivalently as a projective limit
of finite rings endowed with the topology of projective limit or as
a profinite abelian group endowed with a continuous associative
multiplication with unit) one can even obtain
the comodule-contramodule correspondence.
 Namely, the coderived category of discrete left $R$\module s is
equivalent to the contraderived category of left $R$\contramodule s;
this equivalence is constructed in the way analogous
to~\ref{prelim-co-contra-correspondence} with the role of 
a coalgebra~$\C$ played by the discrete $R$\+$R$\bimodule{} of
continuous abelian group homomorphisms $R\rarrow\boQ/\boZ$.
 More generally, let $R$ be a topological ring where open two-sided
ideals form a base of neighborhods of zero and the discrete
quotient rings are right Artinian.
 A pseudo-compact right $R$\module{} is a topological right
$R$\module{} whose open submodules form a base of neighborhoods
of zero and discrete quotient modules have finite length.
 The category of pseudo-compact right $R$\module s is an abelian
category with exact functors of infinite products and enough
projectives; the projective pseudo-compact right $R$\module s
are the direct summands of infinite products of copies of
the pseudo-compact right $R$\module~$R$.
 There is a natural anti-equivalence between the contraderived
categories of pseudo-compact right $R$\module s and left
$R$\contramodule s provided by the derived functors of
pseudo-compact module and contramodule homomorphisms into~$R$.
\end{rmk}

\Section{Comparison with Arkhipov's $\Ext^{\infty/2+*}$ \protect\\
and Sevostyanov's $\Tor_{\infty/2+*}$}

 Semi-infinite cohomology of associative algebras was introduced
by S.~Arkhi\-pov~\cite{Ar1,Ar2}; later A.~Sevostyanov studied it
in~\cite{Sev}.
 The constructions of derived functors $\SemiTor$ and $\SemiExt$
in this book are based on three key ideas which were not known
in the '90s: namely, (i)~the notion of a semialgebra and the functors
of semitensor product and semihomomorphisms; (ii)~the constructions
of adjusted objects from Lemmas~\ref{coflat-semimodule-injection}
and~\ref{coproj-coinj-semi-mod-contra}; and (iii)~the definitions of
semiderived categories.
 We have discussed already Sevostyanov's substitute for~(i) 
in~\ref{semi-product-morphisms} and mentioned Arkhipov's substitute
for~(ii) in~\ref{entwined-coring-semialgebra}.
 Here we consider Arkhipov's substitute for~(i) and suggest
an Arkhipov and Sevostyanov-style substitute for~(iii).
 Combining these together, we obtain comparison results
relating out $\SemiExt$ to Arkhipov's $\Ext^{\infty/2+*}$ and
our $\SemiTor$ to Sevostyanov's $\Tor_{\infty/2+*}$.

 Throughout this appendix we will freely use the notation and remarks
of~\ref{graded-semi}.

\subsection{Algebras $R$ and $R^\#$}

\subsubsection{}  \label{complementary-graded-subalgebras}
 Let $R$ be a graded associative algebra over a field~$k$ endowed with
a pair of subalgebras $K$ and $B\subset R$.
 Assume that all the components $K_i$ are finite dimensional,
$K_i=0$ for $i$ large negative, and $B_i=0$ for $i$ large positive.
 Set $\C_i=K_{-i}^*$ and $\C=\bigoplus_i \C_i$; this is the coalgebra
graded dual to the algebra $K$.
 The coalgebra structure on $\C$ exists due to the conditions imposed
on the grading of $K$.
 There is a natural pairing $\phi\:\C\ot_k K\rarrow k$ satisfying
the conditions of~\ref{ring-coring-pairing}.

 Notice that a structure of graded (left or right) $\C$\comodule{} on
a graded $k$\+vector space $M$ with $M_i=0$ for $i\gg0$ is the same
that a structure of graded (left or right) $K$\module{} on $M$.
 Analogously, a structure of graded (left or right) $\C$\contramodule{}
on a graded $k$\+vector space $P$ with $P_i=0$ for $i\ll0$ is the same
that a structure of graded (left or right) $K$\module{} on $P$.
 Indeed, one has $\Hom_k^\gr(\C,P)\simeq K\ot_k P$.

 Furthermore, assume that the multiplication map $K\ot_k B\rarrow R$
is an isomorphism of graded vector spaces.
 The algebra $R$ is uniquely determined by the algebras $K$ and $B$
and the ``permutation'' map $B\ot_k K\rarrow K\ot_k B$ obtained
by restricting the multiplication map $R\ot_k R\rarrow R\simeq
K\ot_k B$ to the subspace $B\ot_k K\subset R\ot_k R$.
 Transfering the tensor factors $K$ to the other sides of this arrow,
one obtains a map $\C\ot_k B\rarrow\Hom_k^\gr(K,B)$ given by
the formula $c\ot b\mpsto (k'\mapsto (\phi\ot\id_B)(c\ot bk'))$,
where the graded $\Hom$ space in the right hand side is defined, as
always, as direct sum of the spaces of homogeneous maps of various
degrees.
 By the conditions imposed on the gradings of $K$ and $B$, we have
$B\ot_k\C\simeq\Hom_k^\gr(K,B)$, so we get a homogeneous map
$\psi\:\C\ot_k B \rarrow B\ot_k \C$.
 One can check that the map~$\psi$ is a right entwining structure for
the graded coalgebra $\C$ and the graded algebra $B$ over~$k$.

 Conversely, if the map~$\psi$ corresponding to a ``permutation''
map $B\ot_k K\rarrow K\ot_k B$ satisfies the entwining structure
equations, then the latter map can be extended to an associative
algebra structure on $R=K\ot_k B$ with subalgebras $K$ and $B\subset R$.
 However, \emph{not every homogeneous map\/ $\C\ot_k B\rarrow
B\ot_k \C$ comes from a homogeneous map $B\ot_k K\rarrow K\ot_k B$}.

 In the described situation the constructions
of~\ref{construction-of-semialgebras} and~\ref{entwining-structures}
produce the same graded semialgebra $\C\ot_K R = \S\simeq \C\ot_k B$.
 The pairing $\phi\:\C\ot_k K\rarrow k$ is nondegenerate in~$\C$,
so the functor $\Delta_\phi$ is fully faithful and in order to show
that the construction of~\ref{construction-of-semialgebras} works one
only has to check that there exists a right $\C$\comodule{} structure
on $\C\ot_K R$ inducing the given right $K$\module{} structure.
 This is so because $\S_i=0$ for $i\gg0$ according to the conditions
imposed on the gradings of $K$ and $B$.

\subsubsection{}  \label{r-r-sharp}
 Now suppose that we are given two graded algebras $R$ and $R^\#$
with the same two graded subalgebras $K$, $B\subset R$ and $K$,
$B\subset R^\#$ such that the multiplication maps $K\ot_k B\rarrow R$
and $B\ot_k K\rarrow R^\#$ are isomorphisms of vector spaces.
 Assume that $\dim_k K_i<\infty$ for all~$i$, \ 
$K_i=0$ for $i\ll0$, and $B_i=0$ for $i\gg0$.
 Furthermore, assume that the right entwining structure
$\psi\:\C\ot_k B \rarrow B\ot_k \C$ coming from the ``permutation'' map
in $R$ and the left entwining structure $\psi^\#\:B\ot_k\C\rarrow
\C\ot_k B$ coming from the ``permutation'' map in $R^\#$ are inverse
to each other.

 Then there are isomorphisms of graded semialgebras $\S=\C\ot_K R
\simeq\C\ot_k B\simeq B\ot_k\C\simeq R^\#\ot_K\C=\S^\#$, which allow
one to describe left and right $\S$\semimodule s and
$\S$\semicontramodule s in terms of left and right $R$\module s
and $R^\#$\module s.
 In particular, $\S$ has a natural structure of graded
$R^\#$\+$R$\bimodule.

 By the graded version of the result of~\ref{semi-mod-contra-described},
a structure of graded right $\S$\semimodule{} on a graded $k$\+vector
space $N$ with $N_i=0$ for $i\gg0$ is the same that a structure of
graded right $R$\module{} on $N$.
 A structure of graded left $\S$\semimodule{} on a graded $k$\+vector
space $M$ with $M_i=0$ for $i\gg0$ is the same that a structure of
graded left $R^\#$\module{} on $M$.
 A structure of graded left $\S$\semicontramodule{} on a graded
$k$\+vector space $P$ with $P_i=0$ for $i\ll0$ is the same that
a structure of graded left $R$\module{} on $P$.
 In other words, there are isomorphisms of the corresponding categories
of graded modules and homogeneous morphisms between them.

 Besides, for any graded right $R$\module{} $N$ with $N_i=0$ for
$i\gg0$ and any graded left $R$\module{} $P$ with $P_i=0$ for $i\ll0$
there is a natural isomorphism $N\Ocn_\S^\gr P\simeq N\ot_R P$.
 Indeed, one has $N\ocn^\gr_\C P \simeq N\ot_K P$ and $(N\oc_\C^\gr\S)
\ocn_\C^\gr P \simeq N\ot_K R\ot_K P$.

\subsubsection{}  \label{graded-semihom-described}
 A graded $K$\module{} $M$ with $M_i=0$ for $i\gg0$ is injective as
a graded $\C$\comodule{} if and only if it is injective as a graded
$K$\module{} and if and only if it is injective in the category of
graded $K$\module s with the same restriction on the grading.
 Analogously, a graded $K$\module{} $P$ with $P_i=0$ for $i\ll0$ is
projective as a graded $\C$\contramodule{} if and only if it is
projective as a graded $K$\module{} and if and only if it is
projective in the category of graded $K$\module s with the same
restriction on the grading.

 By the graded version of
Proposition~\ref{semitensor-contratensor-assoc}.1(a), for any graded
right $R$\module{} $N$ with $N_i=0$ for $i\gg0$ and any $K$\injective{}
graded left $R^\#$\module{} $M$ with $M_i=0$ for $i\gg0$ there are
natural isomorphisms 
 $$
  N\os_\S^\gr M\simeq N\Ocn_\S^\gr\Psi_\S^\gr(M)\simeq
  N\Ocn_\S^\gr \Hom_{R^\#}^\gr(\S,M).
 $$
 Analogously, for any $K$\injective{} graded right $R$\module{} $N$
with $N_i=0$ for $i\gg0$ and any graded left $R^\#$\module{} $M$
with $M_i=0$ for $i\gg0$ there is a natural isomorphism
 $$
  N\os_\S^\gr M\simeq M\Ocn_{\S^\op}^\gr\Hom_{R^\rop}^\gr(\S,N)
 $$
 The contratensor products in the right hand sides of these formulas
\emph{cannot} be in general replaced by the tensor product over $R$
and $R^\#$, as the graded $\S$\semicontramodule{} $\Psi_\S^\gr(M)$
does not have zero components in large negative degrees.
 In this situation the contratensor product is a certain quotient
space of the tensor product.

 By the graded version of
Proposition~\ref{semitensor-contratensor-assoc}.3(a), for any
$K$\injective{} graded left $R^\#$\module{} $M$ with $M_i=0$ for
$i\gg0$ and any graded left $R$\module{} $P$ with $P_i=0$ for $i\ll0$
there are natural isomorphisms
 $$
  \SemiHom_\S^\gr(M,P)\simeq\Hom^\S_\gr(\Psi_\S^\gr(M),P)\simeq
  \Hom^\S_\gr(\Hom_{R^\#}^\gr(\S,M),P).
 $$
 Here the homomorphisms of graded $\S$\semicontramodule s again
\emph{cannot} be replaced by homomorphisms of graded left $R$\module s.
 The former homomorphism spaces are certain subspaces of the latter
ones.

 By the graded version of
Proposition~\ref{semitensor-contratensor-assoc}.2(a), for any
graded left $R^\#$\module{} $M$ with $M_i=0$ for $i\gg0$ and any
$K$\projective{} graded left $R$\module{} $P$ with $P_i=0$ for $i\ll0$
there are natural isomorphisms
 $$
  \SemiHom_\S^\gr(M,P)\simeq\Hom_\S^\gr(M,\Phi_\S^\gr(P))\simeq
  \Hom_{R^\#}^\gr(M\;\S\ot_R^\gr P).
 $$
 Here the homomorphisms of graded left $\S$\semimodule s \emph{can} be
replaced by the homomorphisms of graded left $R^\#$\module s, since
the functor $\Delta_\phi$ is fully faithful, and consequently
so is the functor $\Delta_{\phi,f}$.

 All of these formulas except the last one have ungraded versions:
\begin{gather*}
  N\os_\S M\simeq N\Ocn_\S \Hom_{R^\#}(\S,M), \qquad
  N\os_\S M\simeq M\Ocn_{\S^\op}\Hom_{R^\rop}(\S,N), \\
  \SemiHom_\S(M,P)\simeq\Hom^\S(\Hom_{R^\#}(\S,M),P)
\end{gather*}
under the appropriate $K$\+injectivity conditions.

\subsection{Finite-dimensional case}
 When the subalgebra $K\subset R$ is finite-dimensional, the algebra
$R^\#$ can be constructed without any reference to the grading or
the complementary subalgebra $B$.

\subsubsection{}
 Let $K$ be a finite-dimensional $k$\+algebra and $\C=K^*$ be 
the coalgebra dual to $K$.
 Then the categories of left $\C$\comodule s and left
$\C$\contramodule s are isomorphic to the category of left
$K$\module s and the category of right $\C$\comodule s is
isomorphic to the category of right $K$\module s.

 The adjoint functors $\Phi_\C$ and $\Psi_\C$ can be therefore
considered as adjoint endofunctors on the category of left
$K$\module s defined by the formulas $P\mpsto \C\ot_K P$ and
$M\mpsto K\oc_\C M\simeq\Hom_K(\C,M)$.
 The restrictions of these functors define an equivalence between
the categories of projective and injective left $K$\module s.

 By Proposition~\ref{tensor-cotensor-assoc}(a-b), the mutually
inverse equivalences $P\mpsto \C\ot_K P$ and $M\mpsto K\oc_\C M$
between the category of $K$\+$K$\bicomodule s that are projective as
left $K$\module s and the category of $K$\+$K$\bicomodule s that are
injective as left $K$\module s transforms the functor of tensor
product over~$K$ in the former category into the functor of
cotensor product over~$\C$ in the latter one.
 In other words, these two tensor categories are equivalent, and
therefore there is a correspondence between ring objects in
the former and the latter tensor category.

\subsubsection{}   \label{finite-dimensional-r-r-sharp}
 Let $K$ be a finite-dimensional $k$\+algebra and $K\rarrow R$ be
a morphism of $k$\+algebras.
 By the above argument, if $R$ is a projective left $K$\module,
then the tensor product $\S=\C\ot_K R$ has a natural structure of
semialgebra over~$\C$.
 Furthermore, if $\S$ is an injective right $K$\module, then
the cotensor product $R^\#=\S\oc_\C K$ has a natural structure of
$k$\+algebra endowed with a $k$\+algebra morphism $K\rarrow R^\#$.
 In this case the semialgebra $\S$ can be also obtained as
the tensor product $R^\#\ot_K\C$.

 By the result of~\ref{semi-mod-contra-described}, a structure of
right $\S$\semimodule{} on a $k$\+vector space $N$ is the same
that a structure of right $R$\module{} on $N$.
 A structure of left $\S$\semimodule{} on a $k$\+vector space $M$
is the same that a structure of left $R^\#$\module{} on $M$.
 A structure of left $\S$\semicontramodule{} on a $k$\+vector
space $P$ is the same that a structure of left $R$\module{} on $P$.
 In other words, the corresponding categories are isomorphic.
 Besides, for any right $R$\module{} $N$ and any left $R$\module{}
$P$ there is a natural isomorphism $N\Ocn_\S P\simeq N\ot_R P$
(see~\ref{contratensor-described}).

\begin{rmk}
 The case of a Frobenius algebra $K$ is of special interest.
 In this case the $k$\+algebra $R^\#$ is isomorphic to
the $k$\+algebra $R$, but the $k$\+algebra morphisms $K\rarrow R$
and $K\rarrow R^\#$ differ by the Frobenius automorphism of $K$.
\end{rmk}

\subsubsection{}  \label{finite-dim-semitensor-semihom-described}
 By Proposition~\ref{semitensor-contratensor-assoc}.1(a), for any
right $R$\module{} $N$ and any $K$\injective{} left $R$\module{} $M$
there are natural isomorphisms
 $$
  N\os_\S M\simeq N\Ocn_\S\Psi_\S(M)\simeq
  N\ot_R \Hom_{R^\#}(\S,M).
 $$
 Analogously, for any $K$\injective{} right $R$\module{} $N$ and any
left $R^\#$\module{} $M$ there is a natural isomorphism
 $$
  N\os_\S M\simeq \Hom_{R^\rop}(\S,N)\ot_{R^\#}M.
 $$
 By Proposition~\ref{semitensor-contratensor-assoc}.3(a), for any
$K$\injective{} left $R^\#$\module{} $M$ and any left $R$\module{} $P$
there are natural isomorphisms
 $$
  \SemiHom_\S(M,P)\simeq\Hom^\S(\Psi_\S(M),P)\simeq
  \Hom_R(\Hom_{R^\#}(\S,M),P).
 $$
 By Proposition~\ref{semitensor-contratensor-assoc}.2(a), for any
left $R^\#$\module{} $M$ and any $K$\projective{} left $R$\module{} $P$
there are natural isomorphisms
 $$
  \SemiHom_\S(M,P)\simeq\Hom_\S(M,\Phi_\S(P))\simeq
  \Hom_{R^\#}(M\;\S\ot_R P).
 $$
 All of these formulas have obvious graded versions.

\subsection{Semijective complexes}  \label{semijective-complexes}
 Let $\S$ be a graded semialgebra over a graded coalgebra~$\C$
over a field~$k$.
 Suppose that $\S_i=0=\C_i$ for $i>0$ and $\C_0=k$.
 Assume also that $\S$ is an injective left and right graded
$\C$\comodule.
 Let $\C\comodl^\down$ and $\comodrdown\C$ denote the categories
of $\C$\comodule s graded by nonpositive integers, $\C\contra^\up$
denote the category of left $\C$\comodule s graded by nonnegative
integers, $\S\simodl^\down$, $\simodrdown\S$, and $\S\sicntr^\up$
denote the categories of graded $\S$\semimodule s and
$\S$\semicontramodule s with analogously bounded grading.

 Any acyclic complex over $\C\comodl^\down$ is coacyclic with
respect to $\C\comodl^\down$.
 Analogously, any acyclic complex over $\C\contra^\up$ is
contraacyclic with respect to $\C\contra^\up$.
 Indeed, let $\K^\bu$ be an acyclic complex of nonpositively
graded $\C$\comodule s.
 As before, we denote by upper indices the homological grading
and by lower indices the internal grading.
 Introduce an increasing filtration on $\K^\bu$ by the complexes of
graded subcomodules $F_n \K^j=\bigoplus_{i\ge -n}\K_i^j$.
 Then the acyclic complexes of trivial $\C$\comodule s $F_n\K^\bu/
F_{n-1}\K^\bu$ are clearly coacyclic.

 So we have $\sD^\si(\S\simodl^\down)=\sD(\S\simodl^\down)$ and
$\sD^\si(\S\sicntr^\up)=\sD(\S\sicntr^\up)$.

 A complex $\M^\bu$ over $\C\comodl^\down$ is called \emph{injective}
if the complex of homogeneous homomorphisms into $\M^\bu$ from any
acyclic complex over $\C\comodl^\down$ is acyclic.
 In this case the complex of homogeneous homomorphisms into $\M^\bu$
from any acyclic complex over $\C\comodl^\sgr$ is also acyclic.
 Analogously, a complex $\P^\bu$ over $\C\contra^\up$ is called
\emph{projective} if the complex of homogeneous homomorphisms from
$\P^\bu$ into any acyclic complex over $\C\contra^\up$ is acyclic.
 In this case the complex of homogeneous homomorphisms from $\P^\bu$
into any acyclic complex over $\C\contra^\sgr$ is also acyclic.

 By Lemma~\ref{comodule-contramodule-subsect}, any complex of injective
objects in $\C\comodl^\down$ is injective and any complex of projective
objects in $\C\contra^\up$ is projective.

 A complex $\bM^\bu$ over $\S\simodl^\down$ is called \emph{quite
$\S/\C$\projective} if the complex of homogeneous homomorphisms
from $\bM^\bu$ into any $\C$\contractible{} complex over
$\S\simodl^\down$ is acyclic.
 Equivalently, $\bM^\bu$ should belong to the minimal triangulated
subcategory of $\Hot(\S\simodl^\down)$ containing the complexes
of graded semimodules induced from complexes over $\C\comodl^\down$
and closed under infinite direct sums.
 Indeed, any quite $\S/\C$\projective{} complex of graded
$\S$\semimodule s is homotopy equivalent to the total complex of
its bar resolution.

 Analogously, a complex $\bP^\bu$ over $\S\sicntr^\up$ is called
\emph{quite $\S/\C$\injective} if the complex of homogeneous
homomorphisms into $\bP^\bu$ from any $\C$\contractible{} complex
over $\S\sicntr^\up$ is acyclic.
 Equivalently, $\bP^\bu$ should belong to the minimal triangulated
subcategory of $\Hot(\S\sicntr^\up)$ containing the complexes of
graded semicontramodules coinduced from complexes over
$\C\contra^\up$ and closed under infinite products.

 A complex $\bM^\bu$ over $\S\simodl^\down$ is called
\emph{semijective} if it is $\C$\injective{} and quite
$\S/\C$\projective.
 Analogously, a complex $\bP^\bu$ over $\S\sicntr^\up$ is called
\emph{semijective} if it is $\C$\projective{} and quite
$\S/\C$\injective.
 Clearly, any acyclic semijective complex of semimodules or
semicontramodules is contractible.

 By the graded version of Remark~\ref{birelatively-adjusted}, any
semiprojective complex of nonpositively graded $\C$\injective{}
$\S$\semimodule s is semijective and any semiinjective complex of
nonnegatively graded $\C$\projective{} $\S$\semicontramodule s is
semijective.
 Hence the homotopy category of semijective complexes over
$\S\simodl^\down$ or $\S\sicntr^\up$ is equivalent to the derived
category $\sD(\S\simodl^\down)$ or $\sD(\S\sicntr^\up)$, any
semijective complex is semiprojective or semiinjective, and one can
use semijective complexes to compute the derived functors
$\SemiTor^\S$ and $\SemiExt_\S$.

\subsection{Explicit resolutions}
 Let us return to the situation of~\ref{r-r-sharp}, but make
the stronger assumptions that $\dim_k K_i<\infty$ for all~$i$, \ 
$K_i=0$ for $i<0$, \ $K_0=k$, and $B_i=0$ for $i>0$.
 Set $\C_i=K_{-i}^*$ and $\S=\C\ot_K R\simeq R^\#\ot_K\C=\S^\#$.

\subsubsection{}  \label{arkhipov-compare}
 For any complex of nonnegatively graded left $R$\module s $P^\bu$
denote by $\boL_2(P^\bu)$ the total complex of the reduced relative
bar complex
$$
 \dsb\lrarrow R\ot_B R/B\ot_B R/B\ot_B P^\bu\lrarrow
 R\ot_B R/B\ot_B P^\bu\lrarrow R\ot_B P^\bu.
$$
 It does not matter whether to construct this total complex by
taking infinite direct sums or infinite products in the category
of graded $R$\module s, as the two total complexes coincide.
 The complex $\boL_2(P^\bu)$ is a complex of $K$\projective{}
nonnegatively graded left $R$\module s quasi-isomorphic to
the complex $P^\bu$.

 For any complex of nonpositively graded left $R^\#$\module s $M^\bu$
denote by $\boL_3(M^\bu)$ the total complex of the reduced relative
bar complex
$$
 \dsb\lrarrow R^\#\ot_K R^\#/K\ot_K R^\#/K\ot_K M^\bu\lrarrow
 R^\#\ot_K R^\#/K\ot_K M^\bu\lrarrow R^\#\ot_K M^\bu,
$$
constructed by taking infinite direct sums along the diagonals.
 The complex $\boL_3(M^\bu)$ is a quite $\S/\C$\projective{} complex
of nonpositively graded left $\S$\semimodule s quasi-isomorphic
to the complex $M^\bu$.

 By~\ref{relatively-semiproj-semiinj}
and~\ref{graded-semihom-described}, the complex
$\Hom_{R^\#}^\gr(\boL_3(M^\bu)\;\S\ot_R^\gr \boL_2(P^\bu))$ represents
the object $\SemiExt_\S^\gr(M^\bu,P^\bu)$ in $\sD(k\vect^\sgr)$.
 We have reproduced Arkhipov's explicit complex~\cite{Ar1,Ar2}
computing $\Ext^{\infty/2+*}_R(M^\bu,P^\bu)$.

\subsubsection{}   \label{sevostyanov-compare}
 For any complex of nonpositively graded left $R^\#$\module s $M^\bu$
denote by $\boR_2(M^\bu)$ the total complex of the reduced relative
cobar complex
 \begin{multline*}
  \Hom_B^\gr(R^\#,M^\bu)\lrarrow\Hom_B^\gr(R^\#/B\ot_B R\; M^\bu)\\
  \lrarrow\Hom_B^\gr(R^\#/B\ot_B R^\#/B\ot_B R\; M^\bu)\lrarrow\dsb
 \end{multline*}
 It does not matter whether to construct this total complex by
taking infinite direct sums or infinite products in the category
of graded $R^\#$\module s, as the two total complexes coincide.
 The complex $\boR_2(M^\bu)$ is a complex of $K$\injective{}
nonpositively graded left $R^\#$\module s quasi-isomorphic to
the complex $M^\bu$.

 For any complex of $K$\injective{} nonpositively graded left
$R^\#$\module s $M^\bu$ the complex $\boL_3(M^\bu)$ defined
in~\ref{arkhipov-compare} is a semiprojective complex of
$\C$\injective{} left $\S$\semimodule s, since it is isomorphic
to the total complex of the reduced bar complex
$$
 \dsb\lrarrow\S\oc_\C\S/\C\oc_\C\S/\C \oc_\C M^\bu\lrarrow
 \S\oc_\C\S/\C\oc_\C M^\bu\lrarrow\S\oc_\C M^\bu
$$
and the left $\C$\comodule{} $\S/\C$ is injective in our assumptions.

 Let $N^\bu$ be a complex of nonpositively graded right $R$\module s
and $M^\bu$ be a complex of nonpositively graded left $R^\#$\module s.
 By~\ref{semi-product-morphisms}, the complex
$N^\bu\ot_B^\C\boL_3\boR_2(M^\bu)$ represents the object
$\SemiTor^\S_\gr(N^\bu,M^\bu)$ in $\sD(k\vect^\sgr)$.
 We have reproduced Sevostyanov's explicit complex~\cite{Sev}
computing $\Tor^R_{\infty/2+*}(N^\bu,M^\bu)$.

\subsection{Explicit resolutions for a finite-dimensional subalgebra}
 Let us consider the situation of an associative algebra $R$
endowed with a pair of subalgebras $K$ and $B\subset R$ such that
the multiplication map $K\ot_k B\rarrow R$ is an isomorphism of
vector spaces and $K$ is a finite-dimensional algebra.
 Let $\C=K^*$ be the coalgebra dual to~$K$.
 Then the construction
of~\ref{complementary-graded-subalgebras}--\ref{r-r-sharp}
is applicable, e.~g., with $R$ endowed by the trivial grading,
and whenever the entwining map $\psi\:B\ot_k\C\rarrow\C\ot_k B$
turns out to be invertible, this construction produces an algebra
$R^\#$ with subalgebras $K$ and $B$ and isomorphisms of semialgebras
$\S=\C\ot_K R \simeq\C\ot_k B\simeq B\ot_k\C\simeq R^\#\ot_K\C=\S^\#$.

\subsubsection{}  \label{finite-dim-entwining-resolutions}
 For any complex of right $R$\module s $N^\bu$ denote by
$\boR_2(N^\bu)$ the total complex of the reduced relative cobar
complex
 \begin{multline*}
  \Hom_{B^\op}(R,N^\bu)\lrarrow\Hom_{B^\op}(R\ot_B R/B\; N^\bu)\\
  \lrarrow\Hom_{B^\op}(R\ot_B R/B\ot_B R/B\; N^\bu)\lrarrow\dsb,
 \end{multline*}
constructed by taking infinite direct sums along the diagonals.
 The complex $\boR_2(N^\bu)$ is a complex of $K$\injective{} right
$R$\module s and the cone of the morphism $N^\bu\rarrow\boR_2(N^\bu)$
is $K$\coacyclic{} (and even $R$\coacyclic).
 For any complex of left $R^\#$\module s $M^\bu$ the complex
$\boR_2(M^\bu)$ is constructed in the analogous way.

 For any complex of left $R$\module s $P^\bu$ denote by $\boL_2(P^\bu)$
the total complex of the reduced relative bar complex
$$
 \dsb\lrarrow R\ot_B R/B\ot_B R/B\ot_B P^\bu\lrarrow
 R\ot_B R/B\ot_B P^\bu\lrarrow R\ot_B P^\bu,
$$
constructed by taking infinite products along the diagonals.
 The complex $\boL_2(P^\bu)$ is a complex of $K$\projective{} left
$R$\module s and the cone of the morphism $\boL_2(P^\bu)\rarrow P^\bu$
is $K$\contraacyclic{} (and even $R$\contraacyclic).

 For any complex of right $R$\module s $N^\bu$ denote by
$\boL_3(N^\bu)$ the total complex of the reduced relative bar
complex
$$
 \dsb\lrarrow N^\bu\ot_K R/K\ot_K R/K\ot_K R\lrarrow
 N^\bu\ot_K R/K\ot_K R\lrarrow N^\bu\ot_K R,
$$
constructed by taking infinite direct sums along the diagonals.
 The complex $\boL_3(N^\bu)$ is a quite $\S/\C$\projective{} complex
of right $\S$\semimodule s and the cone of the morphism $\boL_3(N^\bu)
\rarrow N^\bu$ is $\C$\contractible.
 Whenever $N^\bu$ is a complex of $\C$\injective{} right
$\S$\semimodule s, $\boL_3(N^\bu)$ is a semiprojective complex
of $\C$\injective{} right $\S$\semimodule s, as it was explained
in~\ref{sevostyanov-compare}.
 For a complex of left $R^\#$\module s $M^\bu$ the complex
$\boL_3(M^\bu)$ is constructed in the analogous way.

 For any complex of left $R$\module s $P^\bu$ denote by $\boR_3(P^\bu)$
the total complex of the reduced relative cobar complex
 \begin{multline*}
  \Hom_K(R,P^\bu)\lrarrow\Hom_K(R/K\ot_K R\; P^\bu)\\
  \lrarrow\Hom_K(R/K\ot_K R/K\ot_KR \; M^\bu)\lrarrow\dsb,
 \end{multline*}
constructed by taking infinite products along the diagonals.
 The complex $\boR_3(P^\bu)$ is a quite $\S/\C$\injective{} complex
of left $\S$\semicontramodule s and the cone of the morphism
$P^\bu\rarrow\boR_3(P^\bu)$ is $\C$\contractible.
 Whenever $P^\bu$ is a complex of $\C$\projective{} left
$\S$\semicontramodule s, $\boR_3(P^\bu)$ is a semiinjective complex
of $\C$\projective{} right $\S$\semimodule s, since it is
isomorphic to the total complex of the reduced cobar complex
 \begin{multline*}
  \Cohom_\C(\S,P^\bu)\lrarrow\Cohom_\C(\S/\C\oc_\C\S\;P^\bu) \\
  \lrarrow\Cohom_\C(\S/\C\oc_\C\S/\C\oc_\C\S\;P^\bu) \lrarrow\dsb
 \end{multline*}
and the right $\C$\comodule{} $\S/\C$ is injective in our assumptions.

\subsubsection{}
 One can use these resolutions in various ways to compute the derived
functors $\SemiTor^\S$, $\SemiExt_\S$, $\Psi_\S$, $\Phi_\S$,
$\Ext_\S$, $\Ext^\S$, and $\CtrTor^\S$.

 Specifically, for any complex of right $R$\module s $N^\bu$ and
any complex of left $R^\#$\module s $M^\bu$ the object
$\SemiTor^\S(N^\bu,M^\bu)$ in $\sD(k\vect)$ is represented by
either of the four complexes
\begin{gather*}
  N^\bu\ot_R \Hom_{R^\#}(\S,\boL_3\boR_2(M^\bu)), 
  \quad \boL_3(N^\bu)\ot_R \Hom_{R^\#}(\S,\boR_2(M^\bu)), \\
  \Hom_{R^\rop}(\S,\boL_3\boR_2(N^\bu))\ot_{R^\#}M^\bu,
  \quad \Hom_{R^\op}(\S,\boR_2(N^\bu))\ot_{R^\#}\boL_3(M^\bu)
\end{gather*}
according to the formulas
of~\ref{finite-dim-semitensor-semihom-described}
and the results of~\ref{relatively-semiflat}.
 For any complex of left $R^\#$\module s $M^\bu$ and any complex of
left $R$\module s $P^\bu$ the object $\SemiExt_\S(M^\bu,P^\bu)$ in
$\sD(k\vect)$ is represented by either of the four complexes
\begin{gather*}
  \Hom_R(\Hom_{R^\#}(\S,\boL_3\boR_2(M^\bu)),P^\bu),
  \quad \Hom_R(\Hom_{R^\#}(\S,\boR_2(M^\bu)),\boR_3(P^\bu)), \\
  \Hom_{R^\#}(M^\bu\;\S\ot_R\boR_3\boL_2(P^\bu)),
  \quad \Hom_{R^\#}(\boL_3(M^\bu)\;\S\ot_R\boL_2(P^\bu))
\end{gather*}
according to the formulas
of~\ref{finite-dim-semitensor-semihom-described}
and the results of~\ref{relatively-semiproj-semiinj}.

 One can also use the constructions of~\ref{semi-product-morphisms}
instead of the formulas~\ref{finite-dim-semitensor-semihom-described}.

 For any complex of left $R^\#$\module s $M^\bu$ the object
$\Psi_\S(M^\bu)$ in $\sD^\si(\S\sicntr)$ is represented by the complex
of left $R$\module s $\Hom_{R^\#}(\S,\boR_2(M^\bu))$.
 For any complex of left $R$\module s $P^\bu$ the object
$\Phi_\S(P^\bu)$ in $\sD^\si(\S\simodl)$ is represented by
the complex of left $R^\#$\module s $\S\ot_R\boL_2(P^\bu)$.

 For any complexes of left $R^\#$\module s $L^\bu$ and $M^\bu$
the object $\Ext_\S(L^\bu,M^\bu)$ in $\sD(k\vect)$ is represented
by the complex
 $
  \Hom_{R^\#}(\boL_3(L^\bu),\boR_2(M^\bu)).
 $
 For any complexes of left $R$\module s $P^\bu$ and $Q^\bu$
the object $\Ext^\S(P^\bu,Q^\bu)$ in $\sD(k\vect)$ is represented
by the complex
 $
  \Hom_R(\boL_2(P^\bu),\boR_3(Q^\bu)).
 $
 For any complex of right $R$\module s $N^\bu$ and any complex of
left $R$\module s $P^\bu$ the object $\CtrTor^\S(N^\bu,P^\bu)$
in $\sD(k\vect)$ is represented by the complex
 $
  \boL_3(N^\bu)\ot_R\boL_2(P^\bu).
 $
 These assertions follow from the results
of~\ref{semi-ctrtor-definition}.

 In the situation of~\ref{finite-dimensional-r-r-sharp} (with no
complementary subalgebra $B$) one has to use the constructions of
resolutions $\boR_2(N)$, $\boR_2(M)$, and $\boL_2(P)$ from
the proofs of Theorems~\ref{semitor-main-theorem}
and~\ref{semiext-main-theorem} instead of the constructions
of~\ref{finite-dim-entwining-resolutions}.

\Section{Semialgebras Associated to Harish-Chandra Pairs}
\label{hopf-appendix}

\centerline{by Leonid Positselski and Dmitriy Rumynin}
\bigskip\medskip

 Recall that in~\ref{construction-of-semialgebras} we described
the categories of right semimodules and left semicontramodules over
a semialgebra of the form $\S=\C\ot_K R$, but no satisfactory
description of the category of left semimodules over~$\S$ was found.
 Here we consider the situation when $\C$ and $K$ are Hopf algebras
over a field~$k$, and, under certain assumptions, construct a Morita
equivalence between the semialgebras $\C\ot_K R$ and $R\ot_K\C$.
 This includes the case of an algebraic Harish-Chandra pair $(\g,H)$
with a smooth affine algebraic group~$H$.

\subsection{Two semialgebras}

\subsubsection{}   \label{two-hopf-algebras}
 Let $K$ and $\C$ be two Hopf algebras over a field~$k$ with invertible
antipodes~$s$.
 Using Sweedler's notation~\cite{Swe}, we will denote the multiplications
in $K$ and $\C$ by $x\ot y\mpsto xy$ and the comultiplications by
$x\mpsto x_{(1)}\ot x_{(2)}$; the units will be denoted by~$e$ and
the counits by~$\eps$, so that one has $s(x_{(1)})x_{(2)}=\eps(x)e=
x_{(1)}s(x_{(2)})$ for $x\in K$ or~$\C$.
 Let $\lan\,\.,\,\ran\:K\ot_k\C\rarrow k$ be a pairing between $K$
and~$\C$ such that $\lan xy,c\ran=\lan x,c_{(1)}\ran \lan y,c_{(2)}\ran$
and $\lan x,cd\ran=\lan x_{(1)},c\ran \lan x_{(2)},d\ran$ for
$x$, $y\in K$, \ $c$, $d\in\C$; besides, one should have
$\lan x,e\ran=\eps(x)$ and $\lan e,c\ran=\eps(c)$.
 The pairing $\lan\,\.,\.\ran$ will be also presumed to be compatible
with the antipodes, $\lan s(x),c\ran=\lan x,s(c)\ran$.
 
 Finally, suppose that we are given an ``adjoint'' right coaction
of $\C$ in $K$ denoted by $x\mapsto x_{[0]}\ot x_{[1]}$, satisfying
the equations $x_{(1)}ys(x_{(2)})=\lan x,y_{[1]}\ran y_{[0]}$ and
$(xy)_{[0]}\ot (xy)_{[1]}=x_{[0]}y_{[0]}\ot x_{[1]}y_{[1]}$; besides, 
assume that $e_{[0]}\ot e_{[1]}= e\ot e$.
 This coaction should also satisfy the equations of compatibility with
the squares of antipodes $(s^2x)_{[0]}\ot (s^2x)_{[1]}=s^2(x_{[0]})\ot
s^{-2}(x_{[1]})$ and compatibility with the pairing
$\lan s^{-1}(x_{[0]}),c_{(2)}\ran s(c_{(1)})c_{(3)}x_{[1]}=
\lan s^{-1}(x),c\ran e$.
 When the pairing $\lan\,\.,\,\ran$ is nondegenerate in~$\C$, the latter
four equations follow from the first one and the previous assumptions.
 Indeed, one has $\lan y,s^2(x)_{[1]}\ran (s^2x)_{[0]}=y_{(1)}s^2(x)
s(y_{(2)})=s^2((s^{-2}y)_{(1)}xs((s^{-2}y)_{(2)}))=
\lan s^{-2}(y),x_{[1]}\ran s^2(x_{[0]})=\lan y,s^{-2}(x_{[1]})\ran
s^2(x_{[0]})$ and $\lan s^{-1}(x_{[0]}),c_{(2)}\ran \allowbreak
\lan y,s(c_{(1)})c_{(3)}x_{[1]}\ran =
\lan s^{-1}(x_{[0]}),c_{(2)}\ran \allowbreak\lan y_{(3)},x_{[1]}\ran
\allowbreak\lan y_{(1)},s(c_{(1)})\ran\allowbreak\lan y_{(2)},c_{(3)}
\ran = \lan s^{-1}(y_{(3)}xs(y_{(4)})),c_{(2)}\ran\allowbreak
\lan s(y_{(1)}),c_{(1)}\ran\allowbreak\lan y_{(2)},c_{(3)}\ran =
\lan s(y_{(1)})y_{(4)}s^{-1}(x)s^{-1}(y_{(3)})y_{(2)},c\ran =
\lan s^{-1}(x),c\ran\eps(y)$; analogously for the second and the third
equations.

\subsubsection{}   \label{right-hopf-semialgebra-constructed}
 Let $R$ be an associative algebra over~$k$ endowed with a morphism of
algebras $f\:K\rarrow R$ and a right coaction of the coalgebra $\C$,
which we will denote by $u\mpsto u_{[0]}\ot u_{[1]}$, \ $u\in R$.
 Assume that $f$ is a morphism of right $\C$\comodule s and the right
coaction of $\C$ in $R$ satisfies the equations $f(x)uf(s(x))=
\lan x,u_{[1]}\ran u_{[0]}$ and $(uv)_{[0]}\ot (uv)_{[1]}=u_{[0]}v_{[0]}
\ot u_{[1]}v_{[1]}$ for $u$, $v\in R$, \ $x\in K$.

 Define a pairing $\phi_r\:\C\ot_kK\rarrow k$ by the formula
$\phi_r(c,x)=\lan s^{-1}(x),c\ran$.
 The pairing $\phi_r$ satisfies the conditions
of~\ref{ring-coring-pairing}; in particular, it induces a right
action of $K$ in $\C$ given by the formula $c\ract x=
\lan s^{-1}(x),c_{(2)}\ran c_{(1)}$.
 Assume that the morphism of $k$\+algebras~$f$ makes $R$ a flat
left $K$\module.
 We will now apply the construction of~\ref{semialgebra-constructed}
in order to obtain a semialgebra structure on the tensor product
$\S^r=\C\ot_KR$.

 Define a right $\C$\comodule{} structure on $\C\ot_K R$ by
the formula $c\ot_K u\mpsto c_{(1)}\ot_K u_{[0]}\ot c_{(2)}u_{[1]}$.
 First let us check that this coaction is well-defined.
 We have $(c\ract x)\ot u = \lan s^{-1}x, c_{(2)}\ran c_{(1)}\ot u
\mpsto c_{(1)}\ot_K u_{[0]}\ot \lan s^{(-1)}x,c_{(3)}\ran c_{(2)}
u_{[1]}$ and $c\ot f(x)u\mpsto c_{(1)}\ot_K f(x_{[0]})u_{[0]}\ot
c_{(2)}x_{[1]}u_{[1]} = (c_{(1)}\ract x_{[0]})\ot_K u_{[0]}\ot
c_{(2)}x_{[1]}u_{[1]} = c_{(1)}\ot_K u_{[0]}\ot
\lan c_{(2)},s^{-1}(x_{[0]})\ran c_{(3)}x_{[1]}u_{[1]}$; now
$\lan s^{-1}(x),d_{(2)}\ran d_{(1)}=\lan s^{-1}(x_{[0]}),d_{(3)}
\ran d_{(1)}s(d_{(2)})d_{(4)}x_{[1]}=\lan s^{-1}(x_{[0]},d_{(1)}\ran
d_{(2)}x_{[1]}$ for $d\in\C$, \ $x\in K$.
 Furthermore, this right $\C$\comodule{} structure agrees with
the right $K$\module{} structure on $\C\ot_K R$, since
$\lan s^{-1}(x),c_{(2)}u_{[1]}\ran c_{(1)}\ot_K u_{[0]}=
\lan s^{-1}(x_{(2)}),c_{(2)}\ran c_{(1)}\ot_K 
\lan s^{-1}(x_{(1)}),u_{[1]}\ran u_{[0]}=(c\ract x_{(3)})\ot_K
f(s^{-1}(x_{(2)}))uf(x_{(1)})=c\ot_K uf(x)$.
 It is easy to see that this right coaction of~$\C$ commutes with
the left coaction of~$\C$ and that the semiunit map
$\C\rarrow\C\ot_K R$ is a morphism of right $\C$\comodule s.
 Finally, to show that the semimultiplication map
$(\C\ot_K R)\oc_\C(\C\ot_K R)\rarrow\C\ot_K R$ is a morphism of
right $\C$\comodule s, one defines a right $\C$\comodule{} structure
on $\C\ot_K R\ot_K R$ by the formula $c\ot_K u\ot_K v\mpsto
c_{(1)}\ot_K u_{[0]}\ot_K v_{[0]}\ot c_{(2)}u_{[1]}v_{[1]}$ and checks
that both the isomorphism $\C\ot_K R\ot_K R\simeq
(\C\ot_K R)\oc_\C(\C\ot_K R)$ and the map $\C\ot_K R\ot_K R\rarrow
\C\ot_K R$ are morphisms of right $\C$\comodule s.

\subsubsection{}  \label{left-hopf-semialgebra-constructed}
 Define a pairing $\phi_l\:K\ot_k\C\rarrow k$ by the formula
$\phi_l(x,c)=\lan s(x),c\ran$.
 The pairing $\phi_l$ induces a left action of $K$ in $\C$ given
by the formula $x\lact c = \lan s(x),c_{(1)}\ran c_{(2)}$.
 Assume that the morphism of $k$\+algebras $f$ makes $R$ a flat
right $K$\module.
 We will apply the opposite version of the construction
of~\ref{semialgebra-constructed} in order to obtain a semialgebra
structure on $\S^l=R\ot_K\C$.

 Define a left $\C$\comodule{} structure on $R\ot_K\C$ by the formula
$u\ot_K c\mpsto c_{(1)}s^{-1}(u_{[1]})\ot u_{[0]}\ot_K c_{(2)}$.
 This coaction is well-defined, since one has $u\ot(x\lact c)
=u\ot\lan s(x),c_{(1)}\ran c_{(2)}\mpsto \lan s(x),c_{(1)}\ran
c_{(2)}s^{-1}(u_{[1]})\ot u_{[0]}\ot_K c_{(3)}$ and $uf(x)\ot c\mpsto
c_{(1)}s^{-1}(u_{[1]}x_{[1]})\ot u_{[0]}f(x_{[0]})\ot_K c_{(2)} = 
c_{(1)}s^{-1}(x_{[1]})s^{-1}(u_{[1]})\ot u_{[0]}\ot_K
(x_{[0]}\lact c_{(2)}) = \lan s(x_{[0]}),c_{(2)}\ran
c_{(1)}s^{-1}(x_{[1]})s^{-1}(u_{[1]})\ot u_{[0]}\ot_K c_{(3)}$; now
the identity $\lan s(x),c_{(1)}\ran c_{(2)} = \lan s(x_{[0]}),c_{(2)}
\ran c_{(1)}s^{-1}(x_{[1]})$ follows from the identities
$\lan s^{-1}(x),d_{(2)}\ran d_{(1)}=\lan s^{-1}(x_{[0]}),d_{(1)}\ran
\allowbreak d_{(2)}x_{[1]}$ and $(s^2x)_{[0]}\ot (s^2x)_{[1]} =
s^2(x_{[0]})\ot s^{-2}(x_{[1]})$.
 This left $\C$\comodule{} structure agrees with the left $K$\module{}
structure on $R\ot_K\C$, since $\lan s(x),c_{(1)}s^{-1}(u_{[1]}))\ran
u_{[0]}\ot_K c_{(2)} = \lan s(x_{(1)}),s^{-1}(u_{[1]})\ran \ot_K
\lan s(x_{(2)}),c_{(1)}\ran c_{(2)} = \lan x_{(1)},u_{[1]}\ran u_{[0]}
\ot_K (x_{(2)}\lact c)=x_{(1)}us(x_{(2)})x_{(3)}\ot_K c = xu\ot_K c$.
 The rest is analogous to~\ref{right-hopf-semialgebra-constructed};
the left $\C$\comodule{} structure on $R\ot_K R\ot_K \C$ is defined
by the formula $u\ot_K v\ot_K c\mpsto c_{(1)}s^{-1}(v_{[1]})s^{-1}
(u_{[1]})\ot u_{[0]}\ot_K v_{[0]}\ot_K c_{(2)}$.

\subsubsection{}  \label{hopf-semi-contra-modules-described}
 According to~\ref{semi-mod-contra-described}, the category of right
$\S^r$\semimodule s is isomorphic to the category of $k$\+vector spaces
$\bN$ endowed with right $\C$\comodule{} and right $R$\module{}
structures such that $\lan s^{-1}(x),n_{(1)}\ran n_{(0)}=nf(x)$
and $(nr)_{(0)}\ot (nr)_{(1)}=n_{(0)}r_{[0]}\ot n_{(1)}r_{[1]}$
for $n\in\bN$, \ $x\in K$, \ $r\in R$, where $n\mpsto n_{(0)}\ot
n_{(1)}$ denotes the right $\C$\+coaction map and $n\ot r\mpsto nr$
denotes the right $R$\+action map.

 Assuming that $R$ is a projective left $K$\module, the category of left
$\S^r$\semicontramodule s is isomorphic to the category of $k$\+vector
spaces $\bP$ endowed with left $\C$\contramodule{} and left $R$\module{}
structures such that $\pi_\bP(c\mapsto\lan s^{-1}(x),c\ran p)=
f(x)p$ and $\pi_\bP(c\mapsto r_{[0]}g(cr_{[1]}))=r\pi_\bP(g)$ for
$p\in\bP$, \ $x\in K$, \ $c\in \C$, \ $r\in R$, \ $g\in\Hom_k(\C,\bP)$,
where $\pi_\bP$ denotes the contraaction map and $r\ot p\rarrow rp$
denotes the left action map.

 The category of left $\S^l$\semimodule s is isomorphic to the category
of $k$\+vector spaces $\bM$ endowed with left $\C$\comodule{} and 
left $R$\module{} structures such that $\lan s(x),m_{(-1)}\ran m_{(0)}
= f(x)m$ and $(rm)_{(-1)}\ot (rm)_{(0)} = m_{(-1)}s^{-1}(r_{[1]})\ot
r_{[0]}m_{[0]}$ for $m\in \bM$, \ $x\in K$, \ $r\in R$, where
$m\mpsto m_{(-1)}\ot m_{(0)}$ denotes the left $\C$\+coaction map
and $r\ot m\mpsto rm$ denotes the left $R$\+action map.

\subsection{Morita equivalence}

\subsubsection{}   \label{right-hopf-bicomodule-structure}
 Let $\E$ be a $k$\+vector space endowed with a $\C$\+$\C$\bicomodule{}
structure and a right $\C$\module{} structure satisfying the equation
$(jc)_{(-1)}\ot (jc)_{(0)}\ot (jc)_{(1)}=j_{(-1)}c_{(1)}\ot
j_{(0)}c_{(2)}\ot j_{(1)}c_{(3)}$ for $j\in\E$, \ $c\in\C$, where
$j\mpsto j_{(-1)}\ot j_{(0)}\ot j_{(1)}$ denotes the bicoaction map and
$j\ot c\mpsto jc$ denotes the right action map.

 In particular, $\E$ is a right Hopf module~\cite{Swe} over~$\C$,
hence $\E$ is isomorphic to the tensor product $E\ot_k\C$ as a right
$\C$\comodule{} and a right $\C$\module, where $E$ is a $k$\+vector
space which can be defined as the subspace in $\E$ consisting of all
$i\in\E$ such that $i_{(0)}\ot i_{(1)}=i\ot e$.
 One can see that $E$ is a left $\C$\subcomodule{} in $\E$, so
$\E$ can be identified with the tensor product $E\ot_k\C$ endowed
with the bicoaction $(i\ot c)_{(-1)}\ot (i\ot c)_{(0)}\ot (i\ot c)_{(1)}
=i_{(-1)}c_{(1)}\ot (i_{(0)}\ot c_{(2)})\ot c_{(3)}$ and the right
action $(i\ot c)d = i\ot cd$.
 Besides, the isomorphism $E\ot_k\C\simeq \C\ot_k E$ given
by the formulas $i\ot c\mpsto i_{(-1)}c\ot i_{(0)}$ and $c\ot i\mpsto
i_{(0)}\ot s^{-1}(i_{(-1)})c$ identifies $\E$ with the tensor product
$\C\ot_k E$ endowed with the bicoaction $(c\ot i)_{(-1)}\ot
(c\ot i)_{(0)}\ot (c\ot i)_{(1)} = c_{(1)}\ot (c_{(2)}\ot i_{(0)})\ot
s^{-1}(i_{(-1)})c_{(3)}$ and the right action $(c\ot i)d=cd\ot i$.

\subsubsection{}   \label{hopf-morita-main-isomorphism}
 The pairings $\phi_l$ and $\phi_r$ induce left and right actions
of $K$ in $\E$ given by the formulas $x\lact j =
\lan s(x),j_{(-1)}\ran j_{(0)}$ and $j\ract x = 
\lan s^{-1}(x),j_{(1)}\ran j_{(0)}$ for $j\in \E$, \ $x\in K$.
 Assume that these two actions satisfy the equation
$$
 x_{[0]}\lact (jx_{[1]}) = j\ract x, \quad \text{or equivalently,}
 \quad x\lact j = (js^{-1}(x_{[1]}))\ract x_{[0]}.
$$
 Let us construct an isomorphism $\E\ot_K R\simeq R\ot_K\E$.
 Set the map $\E\ot_K R\rarrow R\ot_K\E$ to be given by the formula
$j\ot_K u\mpsto u_{[0]}\ot_K ju_{[1]}$ and the map $R\ot_K\E\rarrow
\E\ot_KR$ to be given by the formula $u\ot_K j\mpsto
js^{-1}(u_{[1]})\ot_K u_{[0]}$ for $j\in\E$, \ $u\in R$.
 We have to check that these maps are well-defined.

 One has $x\lact (jc)=\lan s(x), j_{(-1)}c_{(1)}\ran j_{(0)}c_{(1)}
=\lan s(x_{(2)}),j_{(-1)}\ran \lan s(x_{(1)})c_{(1)}\ran j_{(0)}c_{(1)}
=(x_{(2)}\lact j)(x_{(1)}\lact c)$.
 Therefore, $\lan x_{(1)},d_{(1)}\ran x_{(2)}\lact(jd_{(2)})=
(x_{(3)}\lact j)\allowbreak(x_{(2)}\lact(s^{-1}(x_{(1)})\lact d))=
(x\lact j)d$.
 Now we have $(j\ract x)\ot u\mpsto u_{[0]}\ot_K(j\ract x)u_{[1]}$
and $j\ot f(x)u\mpsto f(x_{[0]})u_{[0]}\ot_K jx_{[1]}u_{[1]}=
f(x_{[0](1)})u_{[0]}f(s(x_{[0](2)}))f(x_{[0](3)})\ot_K jx_{[1]}u_{[1]}=
\lan x_{[0](1)},u_{[1]}\ran u_{[0]}f(x_{[0](2)})\ot_K jx_{[1]}u_{[2]}=
u_{[0]}\ot_K \lan x_{[0](1)},u_{[1]}\ran x_{[0](2)}\lact(jx_{[1]}
u_{(2)})=u_{[0]}\ot_K(x_{[0]}\lact(jx_{[1]}))u_{[1]}$.

 Analogously, one has $(jc)\ract x=\lan s^{-1}(x),j_{(1)}c_{(2)}\ran
j_{(0)}c_{(1)}=\lan s^{-1}(x_{(2)}),j_{(1)}\ran\allowbreak
\lan s^{-1}(x_{(1)}),c_{(2)}\ran j_{(0)}c_{(1)}=
(j\ract x_{(2)})(c\ract x_{(1)})$.
 Therefore, $\lan s^{-1}(x_{(1)}),d_{(1)}\ran (js^{-1}(d_{(2)}))
\ract x_{(2)} = \lan x_{(1)},s^{-1}(d_{(1)})\ran (js^{-1}(d_{(2)}))
\ract x_{(2)} = (j\ract x_{(3)})((s^{-1}(d)\ract s(x_{(1)}))
\ract x_{(2)}) = (j\ract x)s^{-1}(d)$.
 Now we have $u\ot (x\lact j)\mpsto (x\lact j)s^{-1}(u_{[1]})
\ot_K u_{[0]}$ and $uf(x)\ot j\mpsto js^{-1}(x_{[1]})s^{-1}(u_{[1]})
\ot_K u_{[0]}f(x_{[0]}) = js^{-1}(x_{[1]})s^{-1}(u_{[1]})\ot_K 
f(x_{[0](3)})f(s^{-1}(x_{[0](2)}))u_{[0]}\allowbreak f(x_{[0](1)}) =
js^{-1}(x_{[1]})s^{-1}(u_{[2]})\ot_K\lan s^{-1}(x_{[0](1)}),u_{[1]}\ran
f(x_{[0](2)})u_{[0]} = \lan s^{-1}(x_{[0](1)}),u_{[1]}\ran \allowbreak
(js^{-1}(x_{[1]})s^{-1}(u_{[2]}))\ract x_{[0](2)}\ot_K u_{[0]} =
(js^{-1}(x_{[1]}))\ract x_{[0]}\ot_K u_{[0]}$.

 Checking that these two maps are mutually inverse is easy.

\subsubsection{}
 Assume that $R$ is a flat left and right $K$\module.
 Then the vector space $\E\ot_KR\simeq \E\oc_\C(\C\ot_KR)$ is endowed
with the structures of left $\C$\comodule, right $\S^r$\semimodule,
and right $R$\module.
 The vector space $R\ot_K\E\simeq (R\ot_K\C)\oc_\C\E$ is endowed with
the structures of right $\C$\comodule, left $\S^l$\semimodule, and
left $R$\module.

 Let us check that the isomorphism $\E\ot_K R\simeq R\ot_K\E$
preserves the $\C$\+$\C$\bicomodule{} structures.
 Indeed, one has $(j\ot_K u)_{(-1)}\ot (j\ot_K u)_{(0)}\ot
(j\ot_K u)_{(1)} = j_{(-1)}\ot (j_{(0)}\ot_K u_{[0]})\ot j_{(1)}u_{[1]}$
and $(u\ot_K j)_{(-1)}\ot (u\ot_K j)_{(0)}\ot (u\ot_K j)_{(1)} =
j_{(-1)}s^{-1}(u_{[1]})\ot (u_{[0]}\ot_K j_{(0)})\ot j_{(1)}$, hence
$(u_{[0]}\ot_K ju_{[1]})_{(-1)}\ot(u_{[0]}\ot_K ju_{[1]})_{(-0)}\ot
(u_{[0]}\ot_K ju_{[1]})_{(1)} = j_{(-1)}u_{[2]}s^{-1}(u_{[1]})\ot
(u_{[0]}\ot_K j_{(0)}u_{[3]})\ot j_{(1)}u_{[4]} = 
j_{(-1)}\ot (u_{[0]}\ot_K j_{(0)}u_{[1]})\ot j_{[1]}u_{[2]}$.

 Set $\E\ot_KR=\bE\simeq R\ot_K\E$.
 The left and right actions of $R$ in $\bE$ commute; indeed, the left
and the induced right actions of $R$ in $R\ot_K \E$ are given by
the formulas $w(u\ot_K j)= wu\ot_K j$ and $(u\ot_K j)v =
uv_{[0]}\ot_K jv_{[1]}$.

 It follows easily that $\E\oc_\C\S^r\simeq\bE\simeq\S^l\oc_\C\E$
is an $\S^l$\+$\S^r$\bisemimodule.

\subsubsection{}  \label{morita-auto-equivalence}
 Now assume that the $\C$\+$\C$\bicomodule{} $\E$ can be included
into a Morita autoequivalence $(\E,\E\dual)$ of~$\C$.
 This means that a $\C$\+$\C$\bicomodule{} $\E\dual$ is given
together with isomorphisms of $\C$\+$\C$\bicomodule s
$\E\oc_\C\E\dual\simeq\C\simeq\E\dual\oc_\C\E$ such that the two
induced isomorphisms $\E\oc_\C\E\dual\oc_\C\E\birarrow\E$ coincide
and the two induced isomorphisms $\E\dual\oc_\C\E\oc_\C\E\dual
\birarrow\E\dual$ coincide (see~\ref{co-contra-morita-remarks}).
 The Morita equivalence $(\E,\E\dual)$ is unique if it exists,
and it exists if and only if the left $\C$\comodule{} $E$ is
one-dimensional.
 In the latter case, the bicomodule $\E\dual$ is constructed
as follows.

 The left $\C$\+coaction in $E$ has the form $i_{(-1)}\ot i_{(0)}=
c_E\ot i$ for a certain element $c_E\in\C$ such that
$c_{E(1)}\ot c_{E(2)}=c_E\ot c_E$ and $\eps(c_E)=1$.
 Set $E\dual=\Hom_k(E,k)$ and define a left coaction of $\C$
in $E\dual$ by the formula $\check\imath_{(-1)}\ot
\check\imath_{(0)}=s(c_E)\ot\check\imath$.
 Take $\E\dual=E\dual\ot_k\C$ and define the $\C$\+$\C$\bicomodule{}
structure on $\E\dual$ by the formula $(\check{\imath}\ot c)_{(-1)}\ot
(\check{\imath}\ot c)_{(0)}\ot(\check{\imath}\ot c)_{(1)}=
\check{\imath}_{(-1)}c_{(1)}\ot(\check{\imath}_{(0)}\ot c_{(2)})\ot
c_{(3)}$.
 Then one has $\E\oc_\C\E\dual\simeq E\ot_kE\dual\ot_k\C\simeq\C$ and
$\E\dual\oc_\C\E\simeq E\dual\ot_kE\ot_k\C\simeq\C$.
 There is also an isomorphism $E\dual\ot_k\C\simeq\C\ot_kE\dual$ given
by the formulas analogous to~\ref{right-hopf-bicomodule-structure}.

\subsubsection{}
 Taking the cotensor product of the isomorphism $\E\oc_\C\S^r\simeq
\S^l\oc_\C\E$ with the $\C$\+$\C$\bicomodule{} $\E\dual$ on the left,
we obtain an isomorphism $\S^r\simeq\E\dual\oc_\C\S^l\oc_\C\E$.
 Define a semialgebra structure on $\E\dual\oc_\C\S^l\oc_\C\E$ in
terms of the semialgebra structure on $\S^l$ and the isomorphisms
$\E\oc_\C\E\dual\simeq\C\simeq\E\dual\oc_\C\E$
(see~\ref{morita-change-of-coring-construction}).
 Using the facts that $\E\oc_\C\S^r\simeq\S^l\oc_\C\E$ is an
$\S^l$\+$\S^r$\bisemimodule{} and the isomorphism $\E\oc_\C\S^r\simeq
\S^l\oc_\C\E$ forms a commutative diagram with the maps
$\E\rarrow\E\oc_\C\S^r$ and $\E\rarrow\S^l\oc_\C\E$ induced by
the semiunit morphisms of $\S^r$ and $\S^l$, one can check that
$\S^r\simeq\E\dual\oc_\C\S^l\oc_\C\E$ is an isomorphism of
semialgebras over~$\C$.

 It follows that $\S^l$ and $\S^r$ are left and right coflat
$\C$\comodule s.
 Set $\S^r\oc_\C\E\dual=\bE\dual\simeq\E\dual\oc_\C\S^l$; then
$\bE\dual$ is an $\S^r$\+$\S^l$\bicomodule{} and the pair
$(\bE,\bE\dual)$ is a left and right coflat Morita equivalence
between $\S^l$ and $\S^r$ (see~\ref{semi-morita-morphisms}).

 The category of left $\S^r$\semimodule s can be now described.
 Namely, it is equivalent to the category of left $\S^l$\semimodule s;
this equivalence assigns to a left $\S^r$\semimodule{} $\bL$
the left $\S^l$\semimodule{} $\E\oc_\C\bL$ and to a left
$\S^l$\semimodule{} $\bM$ the left $\S^r$\semimodule{}
$\E\dual\oc_\C\bM$.
 On the level of $\C$\comodule s, one has $\E\oc_\C\L\simeq
E\ot_k\L$ and $\E\dual\oc_\C\M\simeq E\dual\ot_k\M$; the left
$\C$\+coaction in $E\ot_k\L$ and $E\dual\ot_k\M$ is given
by the formulas $(i\ot l)_{(-1)}\ot (i\ot l)_{(0)} =
c_El_{(-1)}\ot (i\ot l_{(0)})$ and $(\check\imath\ot m)_{(-1)}
\ot (\check\imath\ot m)_{(0)} = s(c_E)m_{(-1)}\ot
(\check\imath\ot m_{(0)})$.

\subsubsection{}   \label{hopf-unimodular-case}
 In particular, when the left and right actions of $K$ in $\C$
satisfy the equation
$$
 x_{[0]}\lact (cx_{[1]}) = c\ract x, \quad \text{or equivalently,}
 \quad x\lact c = (cs^{-1}(x_{[1]}))\ract x_{[0]},
$$
one can take $\E=\C=\E\dual$.
 Thus the semialgebras $\S^r$ and $\S^l$ are isomorphic in this
case, the isomorphism being given by the formulas $c\ot_K u\mpsto
u_{[0]}\ot_K cu_{[1]}$ and $u\ot_K c\mpsto cs^{-1}(u_{[1]})\ot_K
u_{[0]}$ for $c\in\C$ and $u\in R$.

\begin{rmk}
 One can construct an isomorphism between versions of the semialgebras
$\S^r$ and $\S^l$ in slightly larger generality.
 Namely, let $\chi_r$, $\chi_l\:K\rarrow k$ be $k$\+algebra
homomorphisms satisfying the equations $\chi(x_{[0]})x_{[1]}=\chi(x)e$
and $\chi(x_{[0](2)})x_{[0](1)}\ot x_{[1]}=\chi(x_{(2)})
x_{(1)[0]}\ot x_{(1)[1]}$.
 These equations hold automatically when the pairing $\lan\,\.,\,\ran$
is nondegenerate in~$\C$ (apply $\id\ot\chi$ to the identity
$\lan y, x_{[1]}\ran x_{[0](1)}\ot x_{[0](2)} = (\lan y,x_{[1]}\ran
x_{[0]})_{(1)}\ot \lan y,x_{[1]}\ran x_{[0]})_{(2)} = 
(y_{(1)}xs(y_{(2)}))_{(1)}\ot(y_{(1)}xs(y_{(2)}))_{(2)} =
y_{(1)}x_{(1)}s(y_{(4)})\ot y_{(2)}x_{(2)}s(y_{(3)})$).
 Define pairings $\phi_r\:\C\ot_kK\rarrow k$ and $\phi_l\:K\ot_k\C
\rarrow k$ by the formulas $\phi_r(c,x)=\chi_r(x_{(2)})
\lan s^{-1}(x_{(1)}),c\ran$ and $\phi_l(x,c)=\chi_l(x_{(2)})
\lan s(x_{(1)}),c\ran$, and modify the definitions of the right
and left actions of $K$ in $\C$ accordingly, $c\ract x=
\phi_r(c_{(2)},x)c_{(1)}$ and $x\lact c=\phi_l(x,c_{(1)})c_{(2)}$.
 Then the tensor products $\C\ot_KR$ and $R\ot_K\C$ are semialgebras
over $\C$ with the $\C$\+$\C$\bicomodule{} structures given by
the same formulas as in~\ref{right-hopf-semialgebra-constructed}--%
\ref{left-hopf-semialgebra-constructed}.
 Assuming that the modified left and right actions of $K$ in $\C$
satisfy the above equation, the maps $c\ot_K u\mpsto u_{[0]}\ot_K
cu_{[1]}$and $u\ot_K c\mpsto cs^{-1}(u_{[1]})\ot_K u_{[0]}$ are
mutually inverse isomorphisms between these two semialgebras.
\end{rmk}

\subsubsection{}   \label{hopf-morita-anti-involution}
 The \emph{opposite coring} $\D^\rop$ to a coring $\D$ over
a $k$\+algebra $B$ and the \emph{opposite semialgebra} $\T^\rop$ to
a semialgebra $\T$ over $\D$ are defined in the obvious way;
$\D^\rop$ is a coring over $B^\rop$ and $\T^\rop$ is a semialgebra
over $\D^\rop$.

 In the above assumptions, notice the identity $s(x\lact c)=
s(c)\ract s(x)$ for $x\in K$, \ $c\in\C$.
 Suppose that the $k$\+algebra $R$ is endowed with an anti-endomorphism
$s$ satisfying the equations $f(s(x))=s(f(x))$ and
$c_{(1)}\ot_K (s(u))_{[0]}\ot c_{(2)}(s(u))_{[1]} = c_{(1)}\ot_K
s(u_{[0]})\ot u_{[1]}c_{(2)}$; the second equation follows from
the first one if the pairing $\lan\,\.,\,\ran$ is nondegenerate in~$\C$.
 Then there is a map of semialgebras $s\: (R\ot_K\C)^\rop\rarrow
\C\ot_KR$ compatible with the isomorphism of coalgebras $s\: \C^\rop
\simeq\C$; it is defined by the formula $s(u\ot_K c) = s(c)\ot_K s(u)$.

 Suppose that we are given a map $s\:\E\rarrow\E$ satisfying
the equation $(s(j))_{(0)}\ot (s(j))_{(1)} = s(j_{(0)})\ot s(j_{(-1)})$.
 Then one has $s(x\lact j)=s(j)\ract s(x)$.
 The induced map $s\:R\ot_K\E\rarrow\E\ot_K R$ given by the formula
$s(u\ot_K j)=s(j)\ot_K s(u)$ is a map of right semimodules compatible
with the isomorphism of coalgebras $s\:\C^\rop\simeq\C$ and the map of
semialgebras $s\:\S^{l\mskip.75\thinmuskip\rop}\rarrow\S^r$, where
the right $\S^{l\mskip.75\thinmuskip\rop}$\semimodule{} structure on
$R\ot_K\E$ corresponds to its left $\S^l$\semimodule{} structure.
 
 Now assume that $\C$ is commutative, $K$ is cocommutative, and our
data satisfy the equations $s^2(u)=u$, \ $s^2(j)=j$, \
$s(jc)=s(j)s(c)$, and $(s(u))_{[0]}\ot (s(u))_{[1]} = s(u_{[0]})\ot
u_{[1]}$; the latter equation holds automatically when the pairing
$\lan\,\.,\,\ran$ is nondegenerate in~$\C$.
 Then the composition of the isomorphism of
$\S^l$\+$\S^r$\bisemimodule s $\E\ot_KR\simeq R\ot_K\E$ and the map
$s\:R\ot_K\E\rarrow\E\ot_KR$ in an involution of the bisemimodule
$\bE$ transforming its left $\S^l$\semimodule{} and right
$\S^r$\semimodule{} structures into each other in a way compatible
with the isomorphism of coalgebras $s\:\C^\rop\rarrow\C$ and
the isomorphism of semialgebras $s\:\S^{l\mskip.75\thinmuskip\rop}
\simeq\S^r$.
 In particular, in the situation of~\ref{hopf-unimodular-case}
the map $s\:R\ot_K\C\rarrow\C\ot_KR$ becomes an involutive
anti-automorphism of the semialgebra $\S^l\simeq\S^r$ compatible
with the anti-automorphism~$s$ of the coalgebra~$\C$.

\subsection{Semitensor product and semihomomorphisms,
SemiTor and SemiExt}
 Let us return to the assumptions of~\ref{two-hopf-algebras}--%
\ref{morita-auto-equivalence}.

\subsubsection{}
 Let $\N$ be a right $\C$\comodule{} and $\M$ be a left $\C$\comodule.
 Then one can easily check that the two injections
$\N\oc_\C\E\dual\oc_\C\M\simeq\N\oc_\C(E\dual\ot_k\C)\oc_\C\M\simeq
\N\oc_\C(E\dual\ot_k\M)\rarrow\N\ot_kE\dual\ot_k\M$ and
$\N\oc_\C\E\dual\oc_\C\M\simeq\N\oc_\C(\C\ot_kE\dual)\oc_\C\M\simeq
(\N\ot_kE\dual)\oc_\C\M\rarrow\N\ot_kE\dual\ot_k\M$ coincide.

 Let $\bN$ be a right $\S^r$\semimodule{} and $\bM$ be a left
$\S^l$\semimodule{} (see~\ref{hopf-semi-contra-modules-described}).
 Then the isomorphism $(\bN\ot_K R)\oc_\C (E\dual\ot_k\bM)\simeq
\bN\oc_\C(\C\ot_KR)\oc_\C\E\dual\oc_\C\bM\simeq
\bN\oc_\C\E\dual\oc_\C(R\ot_K\C)\oc_\C\bM\simeq
(\bN\ot_kE\dual)\oc_\C(R\ot_K\bM)$ induced by the isomorphism
$(\C\ot_KR)\oc_\C\E\dual\simeq\E\dual\oc_\C(R\ot_K\C)$ can be
computed as follows.

 There is an isomorphism $(\C\ot_kR)\oc_\C\E\dual\simeq
\E\dual\oc_\C(R\ot_k\C)$ defined by the same formulas that
the isomorphism $\S^r\oc_\C\E\dual\simeq\E\dual\oc_\C\S^l$
($\ot_K$ being replaced with $\ot$).
 Hence the induced isomorphism $(\bN\ot_k R)\oc_\C(E\dual\ot_k\bM)
\simeq (\bN\ot_kE\dual)\oc_\C(R\ot_k\bM)$, which is given by
the simple formula $n\ot r\ot \check\imath\ot m\mpsto
n\ot\check\imath\ot r\ot m$.
 The isomorphisms $(\bN\ot_K R)\oc_\C (E\dual\ot_k\bM)\simeq
(\bN\ot_kE\dual)\oc_\C(R\ot_K\bM)$ and $(\bN\ot_k R)\oc_\C
(E\dual\ot_k\bM)\simeq (\bN\ot_kE\dual)\oc_\C(R\ot_k\bM)$ form
a commutative diagram with the natural maps $(\bN\ot_k R)\oc_\C
(E\dual\ot_k\bM)\rarrow(\bN\ot_K R)\oc_\C(E\dual\ot_k\bM)$ and
$(\bN\ot_kE\dual)\oc_\C(R\ot_k\bM)\rarrow (\bN\ot_kE\dual)\oc_\C
(R\ot_K\bM)$.
 This provides a description of the first isomorphism whenever
the latter two maps are surjective---in particular, when either
$\bN$ or $\bM$ is a coflat $\C$\comodule.
 To compute the desired isomorphism in the general case, it
suffices to represent either $\bN$ or $\bM$ as the kernel of
a morphism of $\C$\+coflat semimodules (both sides of this
isomorphism preserve kernels, since $R$ is a flat left and
right $K$\module).

\subsubsection{}
 Let $\M$ be a left $\C$\comodule{} and $\P$ be a left
$\C$\contramodule.
 Then one can check that the two surjections
$\Hom_k(E\dual\ot_k\M\;\P)\rarrow \Cohom_\C(E\dual\ot_k\M\;\P)\simeq
\Cohom_\C(\E\dual\oc_\C\M\;\P)$ and $\Hom_k(\M,\Hom_k(E\dual,\P))
\rarrow\Cohom_\C(\M,\Hom_k(E\dual,\P))\allowbreak\simeq
\Cohom_\C(\M,\Cohom_\C(\E\dual,\P))$ coincide.

 Let $\bM$ be a left $\S^l$\semimodule{} and $\bP$ be a left
$\S^r$\semicontramodule.
 Assuming that $R$ is a projective left $K$\module, the isomorphism
$\Cohom_\C(E\dual\ot_k\bM\;\allowbreak\Hom_K(R,\bP))\simeq
\Cohom_\C(\E\dual\oc_\C\bM\;\Cohom_\C(\C\ot_KR\;\bP))\simeq
\Cohom_\C((R\ot_K\C)\oc_\C\bM\;\Cohom_\C(\E\dual\,\bP))\simeq
\Cohom_\C(R\ot_K\bM\;\Hom_k(E\dual,\bP))$ induced by the isomorphism
$(\C\ot_KR)\oc_\C\E\dual\simeq\E\dual\oc_\C(R\ot_K\C)$ can be
computed as follows.

{\hbadness=1500
 The isomorphism $\Cohom_\C(E\dual\ot_k\bM\;\Hom_k(R,\bP))\simeq
\Cohom_\C(R\ot_k\bM\;\allowbreak\Hom_k(E\dual,\bP))$ induced by
the isomorphism $(\C\ot_kR)\oc_\C\E\dual\simeq\E\dual\oc_\C(R\ot_k\C)$
is given by the simple formula $g\mpsto h$, \ $h(r\ot m)(\check\imath)
= g(\check\imath\ot m)(r)$.
 The isomorphisms $\Cohom_\C(E\dual\ot_k\bM\;\Hom_K(R,\bP))\simeq
\Cohom_\C(R\ot_K\bM\;\Hom_k(E\dual,\bP))$ and $\Cohom_\C(E\dual\ot_k\bM
\;\Hom_k(R,\bP))\simeq\Cohom_\C(R\ot_k\bM\;\Hom_k(E\dual,\bP))$ form
a commutative diagram with the natural maps $\Cohom_\C(E\dual\ot_k\bM\;
\Hom_K(R,\bP))\rarrow\Cohom_\C(E\dual\ot_k\bM\;\Hom_k(R,\bP))$ and
$\Cohom_\C(R\ot_K\bM\;\Hom_k(E\dual,\bP))\rarrow\Cohom_\C(R\ot_k\bM\;
\Hom_k(E\dual,\bP))$.
 This provides a description of the first isomorphism in the case
when the latter two maps are injective---in particular, when
either $\bM$ is a coprojective $\C$\comodule, or $\bP$ is
a coinjective $\C$\contramodule.
 To compute the desired isomorphism in the general case, it suffices
to either represent $\bM$ as the kernel of a morphism of
$\C$\coprojective{} semimodules, or represent $\bP$ as the cokernel
of a morphism of $\C$\coinjective{} semicontramodules. \par}

\subsubsection{}    \label{hopf-semitensor-product-determined-by}
 Assume that the $k$\+algebra $R$ is endowed with a Hopf algebra
structure $u\mpsto u_{(1)}\ot u_{(2)}$, \ $u\mpsto\eps(u)$
with invertible antipode~$s$ such that $f\:K\rarrow R$
is a Hopf algebra morphism.
 Let $\bN$ be a right $\S^r$\semimodule{} and $\bM$ be a left
$\S^l$\semimodule; assume that either $\bN$ or $\bM$ is a coflat
$\C$\comodule.
 Define right $R$\module{} and right $\C$\comodule{} structures
on the tensor product $\bN\ot_kE\dual\ot_k\bM$ by the formulas
$(n\ot \check\imath\ot m)r = nr_{(2)}\ot \check\imath\ot
s^{-1}(r_{(1)})m$ and $(n\ot \check\imath\ot m)_{(0)}\ot
(n\ot\check\imath\ot m)_{(1)} = (n_{(0)}\ot\check\imath\ot m_{(0)})
\ot s^{-1}(m_{(-1)})c_En_{(1)}$.
 Then the semitensor product $\bN\os_{\S^r}\bE\dual\os_{\S^l}\bM$
(which is easily seen to be associative) is uniquely determined
by these right $R$\module{} and right $\C$\comodule{} structures
on $\bN\ot_kE\dual\ot_k\bM$.

 Indeed, the subspace $\bN\oc_\C(E\dual\ot_k\bM)\simeq
\bN\oc_\C\E\dual\oc_\C\bM\simeq(\bM\ot_kE\dual)\oc_\C\bM$ of
the space $\bN\ot_kE\dual\ot_k\bM$ can be defined by
the equation $n_{(0)}\ot\check\imath\ot m_{(0)}\ot
s^{-1}(m_{(-1)})c_En_{(0)}=n\ot\check\imath\ot m\ot e$.
 The isomorphism $\bN\ot_kR\ot_kE\dual\ot_k\bM\simeq \bN\ot_kE\dual
\ot_k\bM\ot_kR$ given by the formulas $n\ot r\ot\check\imath\ot m
\mpsto n\ot\check\imath\ot r_{(1)}m\ot r_{(2)}$ and
$n\ot \check\imath\ot m\ot r\mpsto n\ot r_{(2)}\ot \check\imath\ot
s^{-1}(r_{(1)})m$ transforms the pair of maps $\bN\ot_kR\ot_kE\dual
\ot_k\bM\birarrow\bN\ot_kE\dual\ot_k\bM$ given by the formulas
$n\ot r\ot \check\imath\ot m\mpsto nr\ot\check\imath\ot m$, \
$n\ot\check\imath\ot rm$ into the pair of maps $\bN\ot_kE\dual
\ot_k\bM\ot_k R\birarrow\bN\ot_kE\dual\ot_k\bM$ given by the formulas
$n\ot \check\imath\ot m\ot r\mpsto nr_{(2)}\ot\check\imath
\ot s^{-1}(r_{(1)})m$, \ $\eps(r)n\ot\check\imath\ot m$.
 This isomorphism also transforms the subspace $(\bN\ot_kR)\oc_\C
(E\dual\ot_k\bM)$ of $\bN\ot_kR\ot_kE\dual\ot_k\bM$, which can be
defined by the equation $n_{(0)}\ot r_{[0]}\ot\check\imath\ot m_{(0)}
\ot s^{-1}(m_{(-1)})c_En_{(1)}r_{[1]} = n\ot r\ot\check\imath\ot
m\ot e$, into the subspace of $\bN\ot_kE\dual\ot_k\bM\ot_kR$ defined by
the equation $n_{(0)}\ot \check\imath\ot r_{(2)[0](1)}
(s^{-1}(r_{(1)}))_{[0]}m_{(0)}\ot r_{(2)[0](2)}\ot
s^{-2}((s^{-1}(r_{(1)}))_{[1]})s^{-1}(m_{(-1)})c_En_{(1)}r_{(2)[1]} =
n\ot\check\imath\ot m\ot r\ot e$.
 Finally, the same isomorphism transforms the quotient space
$\bN\ot_KR\ot_kE\dual\ot_k\bM$ of the space $\bN\ot_kR\ot_kE\dual
\ot_k\bM$ into the quotient space $(\bN\ot_kE\dual\ot_k\bM)\ot_KR$
of the space $\bN\ot_kE\dual\ot_k\bM\ot_kR$, as one can check using
the isomorphism $\bN\ot_kK\ot_kR\ot_kE\dual\ot_k\bM\simeq
\bN\ot_kE\dual\ot_k\bM\ot_kK\ot_kR$ given by the formulas
$n\ot x\ot r\ot\check\imath\ot m\mpsto n\ot\check\imath\ot
x_{(1)}r_{(1)}m\ot x_{(2)}\ot r_{(2)}$ and $n\ot\check\imath\ot
m\ot x\ot r\mpsto n\ot x_{(2)}\ot r_{(2)}\ot \check\imath\ot
s^{-1}(r_{(1)})s^{-1}(x_{(1)})m$.

\subsubsection{}
 Let $\bM$ be a left $\S^l$\semimodule{} and $\bP$ be a left
$\S^r$\semicontramodule; assume that either $\bM$ is a coprojective
$\C$\comodule, or $\bP$ is a coinjective $\C$\contramodule.
 Define left $R$\module{} and left $\C$\contramodule{} structures
on the space $\Hom_k(E\dual\ot_k\bM\;\bP)$ by the formulas
$rg(\check\imath\ot m)=r_{(2)}g(\check\imath\ot s^{-1}(r_{(1)})m)$ and
$\pi(h)(\check\imath\ot m)=\pi_\bP(c\mapsto h(s^{-1}(m_{(-1)})c_Ec)
(\check\imath\ot m_{(0)}))$ for $g\in \Hom_k(E\dual\ot_k\bM\;\bP)$
and $h\in\Hom_k(\C,\Hom_k(E\dual\ot_k\bM\;\bP))$, where $\pi_\bP$
denotes the $\C$\+contraaction in~$\bP$.
 Then the semihomomorphism space $\SemiHom_{\S^r}(\bE\os_{\S^l}\bM\;
\bP)$ is uniquely determined by these $R$\module{} and
$\C$\contramodule{} structures on $\Hom_k(E\dual\ot_k\bM\;\bP)$.

 This is established in the way analogous
to~\ref{hopf-semitensor-product-determined-by} using the isomorphism
$\Hom_k(R\ot_kE\dual\ot_k\bM\;\bP)\simeq\Hom_k(R,\Hom_k(E\dual\ot_k\bM\;
\bP))$ given by the formulas $g\mpsto h$, \ $g(r\ot\check\imath\ot m) =
h(r_{(2)})(\check\imath\ot r_{(1)})$, \ $h(r)(\check\imath\ot m) =
g(r_{(2)}\ot\check\imath\ot s^{-1}(r_{(1)})m)$.

\subsubsection{}  \label{hopf-semitensor-product-factorizes-through}
 Now assume that $\C$ is commutative, $K$ is cocommutative, and
the equations $(s(u))_{[0]}\ot (s(u))_{[1]} = s(u_{[0]})\ot u_{[1]}$, \
$\eps(u_{[0]})u_{[1]}=\eps(u)e$, and $u_{(1)[0]}\ot u_{(2)[0]}\ot
u_{(1)[1]}u_{(2)[1]} = u_{[0](1)}\ot u_{[0](2)}\ot u_{[1]}$ are
satisfied for $u\in R$; when the pairing $\lan\,\.,\,\ran$ is
nondegenerate in~$\C$, these equations hold automatically.

 Let $\bN$ be a right $\S^r$\semimodule{} and $\bM$ be a left
$\S^l$\semimodule.
 Define right $R$\module{} and right $\C$\comodule{} structures
on the tensor product $\bN\ot_k\bM$ by the formulas
$(n\ot m)r = nr_{(2)}\ot s^{-1}(r_{(1)})m$ and $(n\ot m)_{(0)}\ot
(n\ot m)_{(1)} = (n_{(0)}\ot m_{(0)})\ot s^{-1}(m_{(-1)})n_{(1)}$.
 These right action and right coaction satisfy the equations 
of~\ref{hopf-semi-contra-modules-described}, so they define
a right $\S^r$\semimodule{} structure on $\bN\ot_k\bM$.
 The ground field~$k$, endowed with the trivial left $R$\module{}
and left $\C$\comodule{} structures $ra=\eps(r)a$ and
$a_{(-1)}\ot a_{(0)}=e\ot a$ for $a\in k$, becomes a left
$\S^l$\semimodule.
 Then there is a natural isomorphism $\bN\os_{\S^r}\bE\dual\os_{\S^l}\bM
\simeq (\bN\ot_k\bM)\os_{\S^r}\bE\dual\os_{\S^l}k$.

 Indeed, let us first assume that either $\bN$ or $\bM$ is a coflat
$\C$\comodule; notice that $\bN\ot_k\bM$ is then a coflat $\C$\comodule,
too.
 The isomorphism $\bN\oc_\C(E\dual\ot_k\bM)\simeq(\bN\ot_k\bM)\oc_\C
E\dual$ given by the formula $n\ot\check\imath\ot m\mpsto
n\ot m\ot\check\imath$ and the isomorphism $(\bN\ot_KR)\oc_\C
(E\dual\ot_k\bM)\simeq((\bN\ot_k\bM)\ot_KR)\oc_\C E\dual$ constructed
in~\ref{hopf-semitensor-product-determined-by} transform the pair
of maps whose cokernel is $\bN\os_{\S^r}\bE\dual\os_{\S^l}\bM$ into
the pair of maps whose cokernel is $(\bN\ot_k\bM)\os_{\S^r}\bE\dual
\os_{\S^l}k$.
 In the general case, represent $\bN$ or $\bM$ as the kernel of
a morphism of $\C$\+coflat semimodules; then the pair of maps whose
cokernel is $\bN\os_{\S^r}\bE\dual\os_{\S^l}\bM$ and the pair of maps
whose cokernel is $(\bN\ot_k\bM)\os_{\S^r}\bE\dual\os_{\S^l}k$
become the kernels of isomorphic morphisms of pairs of maps.

\subsubsection{}  \label{hopf-semihom-factorizes-through}
 Let $\bM$ be a left $\S^l$\semimodule{} and $\bP$ be a left 
$\S^r$\semicontramodule.
 Define left $R$\module{} and left $\C$\contramodule{} structures
on the space $\Hom_k(\bM,\bP)$ by the formulas $rg(m) =
r_{(2)}g(s^{-1}(r_{(1)})m)$ and $\pi(h)(m)=
\pi_\bP(c\mapsto h(s^{-1}(m_{(-1)})c)(m_{(0)}))$ for
$g\in \Hom_k(\bM,\bP)$ and $h\in\Hom_k(\C,\Hom_k(\bM,\bP))$.
 These left action and left contraaction satisfy the equations
of~\ref{hopf-semi-contra-modules-described}, so they define
a left $\S^r$\semicontramodule{} structure on $\Hom_k(\bM,\bP)$.
 Then there is a natural isomorphism $\SemiHom_{\S^r}(\bE\dual\os_{\S^l}
\bM\;\bP)\simeq\SemiHom_{\S^r}(\bE\dual\os_{\S^l}k\;\Hom_k(\bM,\bP))$.

\subsubsection{}
 Let $\bN^\bu$ be a complex of right $\S^r$\semimodule s and $\bM^\bu$
be a complex of left $\S^l$\semimodule s.
 Then there are natural isomorphisms $\SemiTor^{\S^r}(\bN^\bu\;
\bE\dual\os_{\S^l}\bM^\bu)\simeq\SemiTor^{\S^l}(\bN^\bu\os_{\S^r}
\bE\dual\;\bM^\bu)\simeq\SemiTor^{\S^r}(\bN^\bu\ot_k\bM^\bu\;
\bE\dual\os_{\S^l}k)$ in the derived category of $k$\+vector spaces.
 The isomorphism between the first two objects is provided by
the results of~\ref{morita-change-of-coring}, and the isomorphism
between either of the first two objects and the third object
follows from~\ref{hopf-semitensor-product-factorizes-through}.

 Indeed, assume that the complex $\bN^\bu$ is semiflat.
 Then the complex $\bN^\bu\ot_k\bM^\bu$ is also semiflat, since
$(\bN^\bu\ot_k\bM^\bu)\os_{\S^r}\bL^\bu\simeq(\bN^\bu\ot_k\bM^\bu)
\os_{\S^r}\bE\dual\os_{\S^l}\bE\os_{\S^r}\bL^\bu\simeq(\bN^\bu\ot_k
\bM^\bu\ot_k(\bE\os_{\S^r}\bL^\bu))\os_{\S^r}\bE\dual\os_{\S^l}k\simeq
\bN^\bu\os_{\S^r}\bE\dual\os_{\S^l}(\bM^\bu\ot_k(\bE\os_{\S^r}\bL^\bu))$
for any complex of left $\S^r$\semimodule s $\bL^\bu$.

\subsubsection{}
 Let $\bM^\bu$ be a complex of left $\S^l$\semimodule s and $\bP^\bu$
be a complex of left $\S^r$\semicontramodule s.
 Then there are natural isomorphisms
$\SemiExt_{\S^r}(\bE\dual\os_{\S^l}\bM^\bu\;\allowbreak\bP^\bu)\simeq
\SemiExt_{\S^l}(\bM^\bu,\SemiHom_{\S^r}(\bE\dual,\bP^\bu))\simeq
\SemiExt_{\S^r}(\bE\dual\os_{\S^l}k\;\Hom_k(\bM^\bu,\bP^\bu))$.

\subsection{Harish-Chandra pairs}

\subsubsection{}
 Let $(\g,H)$ be an algebraic Harish-Chandra pair over a field~$k$,
that is $H$ is an algebraic group, which we will assume to be affine,
$\g$ is a Lie algebra into which the Lie algebra $\h$ of the algebraic
group $H$ is embedded, and an action of $H$ by  Lie algebra
automorphisms of $\g$ is given.
 Two conditions should be satisfied: $\h$ is an $H$\+submodule in $\g$
where $H$ acts by the adjoint action of $H$ in~$\h$, and the action
of $\h$ in $\g$ obtained by differentiating the action of $H$ in $\g$
coincides with the adjoint action of $\h$ in~$\g$.
 Notice that the dimension of $\h$ is presumed to be finite, though
the dimension of $\g$ may be infinite.

 Set $K=U(\h)$ and $R=U(\g)$ to be the universal enveloping algebras
of $\h$ and~$\g$, and $f\:K\rarrow R$ to be the morphism induced by
the embedding $\h\rarrow\g$.
 Let $\C=\C(H)$ be the coalgebra of functions on~$H$.
 Then $\C$, \ $K$, and $R$ are Hopf algebras; the adjoint action of
$H$ in $\h$ and the given action of $H$ in $\g$ provide us with
right coactions $x\mpsto x_{[0]}\ot x_{[1]}$ and $u\mpsto u_{[0]}\ot
u_{[1]}$ of $\C$ in $K$ and $R$; and there is a natural pairing
$\lan\,\.,\,\ran\:K\ot_k\C\rarrow k$ such that the equations
of~\ref{two-hopf-algebras}--\ref{right-hopf-semialgebra-constructed}
and~\ref{hopf-semitensor-product-factorizes-through} are satisfied.

\subsubsection{}  \label{hopf-harish-chandra-semi-contra-modules}
 So we obtain two opposite semialgebras $\S^l=\S^l(\g,H)$ and
$\S^r=\S^r(\g,H)$ such that the categories of left $\S^l$\semimodule s
and right $\S^r$\semimodule s are isomorphic to the category of
Harish-Chandra modules over~$(\g,H)$.
 Recall that a Harish-Chandra module $\bN$ over $(\g,H)$ is
a $k$\+vector space endowed with $\g$\module{} and $H$\module{}
structures such that the two induced $\h$\module{} structures coincide
and the action map $\g\ot_k \bN\rarrow \bN$ is a morphism of
$H$\module s.
 The assertion follows from~\ref{hopf-semi-contra-modules-described}; 
indeed, it suffices to notice that the equations
of~\ref{hopf-semi-contra-modules-described} hold whenever they hold
for $x$ and $r$ belonging to some sets of generators of the algebras
$K$ and $R$.

 Analogously, the category of left $\S^r(\g,H)$\semicontramodule s is
isomorphic to the category of $k$\+vector spaces $\bP$ endowed with
$\g$\module{} and $\C(H)$\contramodule{} structures such that
the two induced $\h$\module{} structures coincide and the action
map $\bP\rarrow\Hom_k(\g,\bP)$ is a morphism of $\C(H)$\contramodule s.
 Here a left $\C$\contramodule{} structure induces an $\h$\module{}
structure by the formula $xp=-\pi_\bP(c\mapsto\lan x,c\ran p)$ for
$p\in\bP$, \ $x\in\h$, \ $c\in\C$; for a left $\C$\comodule{} $\M$
and a left $\C$\contramodule{} $\P$, the left $\C$\contramodule{}
structure on $\Hom_k(\M,\P)$ is defined by the formula 
$\pi(g)(m)=\pi_\P(c\mapsto g(s^{-1}(m_{(-1)})c)(m_{(0)}))$ for
$m\in\M$, \ $g\in\Hom_k(\C,\Hom_k(\M,\P))$.

\subsubsection{}
 Now assume that the algebraic group $H$ is smooth (i.~e., reduced).
 Let $\E$ be the $\C$\+$\C$\bicomodule{} and right $\C$\module{} of
differential top forms on~$H$, with the bicomodule structure coming from
the action of $H$ on itself by left and right shifts and the module
structure given by the multiplication of top forms with functions.
 Let $\E\dual$ be the $\C$\+$\C$\bicomodule{} of top polyvector fields
on~$H$.
 Then the equation of~\ref{right-hopf-bicomodule-structure} is clearly
satisfied, and one can see that $(\E,\E\dual)$ is a Morita
autoequivalence of~$\C$.
 The left $\C$\comodule s $E$ and $E\dual$ can be identified with 
the top exterior powers of the vector spaces $\Hom_k(\h,k)$ and~$\h$,
respectively; $c_E$ is the modular character of~$H$.
 The inverse element anti-automorphism of $H$ induces a map $s\:\E
\rarrow\E$ satisfying the equations
of~\ref{hopf-morita-anti-involution}.

 Let us show that the equation of~\ref{hopf-morita-main-isomorphism}
holds for~$\E$.
 First let us check that it suffices to prove the desired equation for
all $x$ belonging to a set of generators of the algebra~$K$.
 Indeed, one has $(xy)_{[0]}\lact (j(xy)_{[1]}) = x_{[0]}y_{[0]}\lact
(jx_{[1]}y_{[1]}) = x_{[0](1)}y_{[0]}s(x_{[0](2)})x_{[0](3)}\lact
(jx_{[1]}y_{[1]}) = \lan x_{[0](1)},y_{[1]}\ran y_{[0]}x_{[0](2)}\lact
(jx_{[1]}y_{[2]}) = y_{[0]}\lact(x_{[0](2)}\lact
(\lan x_{[0](1)},y_{[1]}\ran jx_{[1]}y_{[2]})) =
y_{[0]}\lact((x_{[0]}\lact (jx_{[1]}))y_{[1]})$ for $j\in\E$, \
$x$, $y\in K$, since $x_{(2)}\lact (\lan x_{(1)},c_{(1)}\ran jc_{(2)})
= x_{(2)}\lact (j(s^{-1}(x_{(1)})\lact c)) = (x_{(3)}\lact j)
(x_{(2)}\lact (s^{-1}(x_{(1)})\lact c)) = (x\lact j)c$ for $j\in\E$, \
$x\in K$, \ $c\in\C$.
 So it remains to check that the equation holds for $x\in\h\subset K$.

 For $h\in\h$, let $r_h$ and $l_h$ denote the left- and right-invariant
vector fields on $H$ corresponding to~$h$.
 Then one has $r_h=h_{[1]}l_{h_{[0]}}$, hence $\omega\ract h = 
\Lie_{r_h}\omega = \Lie_{l_{h_{[0]}}}(\omega h_{[1]}) = 
h_{[0]}\lact(\omega h_{[1]})$ for $\omega\in\E$, where $\Lie_v\omega$
denotes the Lie derivative of a top form $\omega$ along a vector
field~$v$.

\subsubsection{}
 Thus there is a left and right coflat Morita equivalence 
$(\bE(\g,H),\bE\dual(\g,H))$ between the semialgebras $\S^r(\g,H)$
and $\S^l(\g,H)$.
 The functors of semitensor product with $\bE(\g,H)$ and
$\bE\dual(\g,H)$ provide mutually inverse equivalences between
the category of left $\S^r$\semimodule s and the category of left
$\S^l$\semimodule s.
 The $\S^l$\+$\S^r$\bisemimodule{} $\bE(\g,H)$ is endowed with
an involutive automorphism transforming its two semimodule structures
into each other in a way compatible with the antipode isomorphisms
$\C^\rop\simeq\C$ and $\S^{l\mskip.75\thinmuskip\rop}\simeq\S^r$.
 When the algebraic group $H$ is unimodular, the semialgebras
$\S^l(\g,H)$ and $\S^r(\g,H)$ are naturally isomorphic and endowed with
an involutive anti-automorphism.

\begin{rmk}
 Let us assume for simplicity that the field~$k$ has characteristic~$0$
and the Harish-Chandra pair $(\g,H)$ originates from an embedding of
affine algebraic groups $H\subset G$.
 Then the bisemimodule $\bE(\g,H)$ can be interpreted geometrically as
the Harish-Chandra bimodule of distributions on $G$, supported on $H$
and regular along~$H$.
 Technically, the desired vector space of distributions can be defined
as the direct image of the right $\Diff_H$\module{} of top forms on $H$
under the closed embedding $H\rarrow G$, where $\Diff$ denotes
the rings of differential operators~\cite{Bern}.
 This vector space has two commuting structures of a Harish-Chandra
module over $(\g,H)$, one given by the action of $H$ by left shifts
and the action of $\g$ by right invariant vector fields, the other
in the opposite way; so it can be considered as
an $\S^l$\+$\S^r$\bisemimodule.
 The desired map from the vector space $\bE\simeq \E\ot_KR$ to the space
of distributions can be defined as the unique map forming a commutative
diagram with the embeddings of the space of top forms $\E$ into both
vector spaces and preserving the right $R$\module{} structures.
 To prove that this map is an isomorphism, it suffices to consider
the filtration of $\bE$ induced by the natural filtration of
the universal enveloping algebra $R$ and the filtration of the space
of distributions induced by the filtration of $\Diff_G$ by the order
of differential operators.
 When the algebraic group $H$ is unimodular, one can identify
the semialgebra $\S^l\simeq\S=\S^r$ itself with the above vector space
of distributions by choosing a nonzero biinvariant top form~$\omega$
on~$H$.
 The semiunit and semimultiplication in $\S$ are then described as
follows.
 Given a function on~$\C$, one has to multiply it with $\omega$ and take
the push-forward with respect to the closed embedding $H\rarrow G$ to
obtain the corresponding distribution under the semiunit map.
 To describe the semimultiplication, denote by $G\times_HG$ the quotient
variety of the Carthesian product $G\times G$ by the equivalence
relation $(g'h,g'')\sim(g',hg'')$.
 Then the pull-back of distributions with respect to the smooth map
$G\times G\rarrow G\times_HG$ using the relative top form $\omega$
identifies $\S\oc_\C\S$ with the space of distributions on $G\times_HG$
supported in $H\subset G\times_HG$ and regular along~$H$.
 The push-forward of distributions with respect to the multiplication
map $G\times_HG\rarrow G$ provides the semimultiplication in~$\S$.
 (Cf.\ Appendix~\ref{groupoid-appendix}.)
\end{rmk}

\subsection{Semiinvariants and semicontrainvariants}

\subsubsection{}   \label{finite-dim-semiinvariants-defining-map}
 Let $\h\subset \g$ be a Lie algebra with a finite-dimensional
subalgebra; let $N$ be a $\g$\module.
 Then there is a natural map $(\det(\h)\ot_k\g/\h\ot_k N)^\h\rarrow
(\det(\h)\ot_kN)^\h$, where $\det(V)$ is the top exterior power
of a finite-dimensional vector space $V$ and the superindex~$\h$
denotes the $\h$\invariant s.
 This natural map is constructed as follows.

 Tensoring the action map $\g\ot_kN\rarrow N$ with $\det(\h)$ and
passing to the $\h$\invariant s, we obtain a map
$(\det(\h)\ot_k\g\ot_k N)^\h\rarrow(\det(\h)\ot_kN)^\h$.
 Let us check that the composition $(\det(\h)\ot_k\h\ot_k N)^\h\rarrow
(\det(\h)\ot_k\g\ot_k N)^\h\rarrow(\det(\h)\ot_kN)^\h$ vanishes.
 Notice that this compostion only depends on the $\h$\module{}
structure on~$N$.
 Let $n$ be an $\h$\invariant{} element of $\det(\h)\ot_k\h\ot_kN$;
it can be also considered as an $\h$\module{} map $\check n\:
(\det(\h)\ot_k\h)\dual\rarrow N$, where $V\dual$ denotes the dual
vector space $\Hom_k(V,k)$.
 Let $t$ denote the trace element of the tensor product
$(\det(\h)\ot_k\h)\ot_k(\det(\h)\ot_k\h)\dual$; then $t$ is
an $\h$\invariant{} element and one has $(\id\ot\check n)(t)=n$.
 So it remains to check that the image of~$t$ under the action map
$(\det(\h)\ot_k\h)\ot_k(\det(\h)\ot_k\h)\dual\rarrow\det(\h)\ot_k
(\det(\h)\ot_k\h)\dual\simeq\h\dual$ vanishes; this is straightforward.

 We have constructed a map $(\det(\h)\ot_k\g\ot_k N)^\h/(\det(\h)\ot_k
\h\ot_k N)^\h\rarrow(\det(\h)\ot_kN)^\h$.
 When $N$ is an injective $U(\h)$\module, this provides the desired
map $(\det(\h)\ot_k\g/\h\ot_k N)^\h\rarrow (\det(\h)\ot_kN)^\h$.
 To construct the latter map in the general case, it suffices to 
represent $N$ as the kernel of a morphism of $U(\h)$\injective{}
$U(\g)$\module s (notice that any injective $U(\g)$\module{} is 
an injective $U(\h)$\module).
 Indeed, both the left and the right hand sides of the desired map
preserve kernels.

 The vector space of \emph{$(\g,\h)$\semiinvariant s} $N_{\g,\h}$ of
a $\g$\module{} $N$ is defined as the cokernel of the map
$(\det(\h)\ot_k\g/\h\ot_k N)^\h\rarrow(\det(\h)\ot_kN)^\h$ that we
have obtained.
 The $(\g,\h)$\semiinvariant s are a mixture of invariants along~$\h$
and coinvariants in the direction of $\g$ relative to~$\h$.

\subsubsection{}   \label{finite-dim-semicontrainvs-defining-map}
 Let $P$ be another $\g$\module.
 Then there is a natural map $\Hom_k(\det(\h),P)_\h\allowbreak\rarrow
\Hom_k(\det(\h)\ot_k\g/\h\;P)_\h$, where the subindex $\h$ denotes
the $\h$\coinvariant s.
 This map is constructed as follows.

 Tensoring the action map $P\rarrow\Hom_k(\g,P)$ with $\det(\h)\dual$
and passing to the $\h$\coinvariant s, we obtain a map
$\Hom_k(\det(\h),P)_\h\rarrow\Hom_k(\det(\h)\ot_k\g\;P)_\h$.
 Let us check that the composition $\Hom_k(\det(\h),P)_\h\rarrow
\Hom_k(\det(\h)\ot_k\g\;P)_\h\rarrow\Hom_k(\det(\h)\ot_k\h\;P)_\h$
vanishes.
 Notice that this composition only depends on the $\h$\module{}
structure on~$P$.
 Let $p$ be an $\h$\invariant{} map $\Hom_k(\det(\h)\ot_k\h\;P)
\rarrow k$; it can be also considered as an $\h$\module{} map
$\check p\: P\rarrow \det(\h)\ot\h$.
 Then the map $p$ factorizes through the map $\Hom_k(\det(\h)\ot_k\h\;P)
\rarrow\Hom_k(\det(\h)\ot_k\h\;\det(\h)\ot\h)$ induced by $\check p$.
 So it suffices to consider the case of a finite-dimensional
$\h$\module{} $P=\det(\h)\ot\h$, when the assertion follows by duality
from the result of~\ref{finite-dim-semiinvariants-defining-map}.

 We have constructed a map from $\Hom_k(\det(\h),P)_\h$ to
the kernel of the map $\Hom_k(\det(\h)\ot_k\g\;P)_\h\rarrow
\Hom_k(\det(\h)\ot_k\h\;P)_\h$.
 When $P$ is a projective $U(h)$\module, this provides the desired
map $\Hom_k(\det(\h),P)_\h\rarrow\Hom_k(\det(\h)\ot_k\g/\h\;P)_\h$.
 In the general case, represent $P$ as the cokernel of a morphism of
$U(\h)$\projective{} $U(\g)$\module s and notice that both the left
and the right hand side of the desired map preserve cokernels.

 The vector space of \emph{$(\g,\h)$\semicontrainvariant s} $P^{\g,\h}$
of a $\g$\module{} $P$ is defined as the kernel of the map
$(P\ot_k\det(\h)\dual)_\h\rarrow\Hom_k(\det(\h)\ot_k\g/\h,P)_\h$
that we have obtained.
 The $(\g,\h)$\semicontrainvariant s are a mixture of coinvariants
along~$\h$ and invariants in the direction of $\g$ relative to~$\h$.

\subsubsection{}
 Now let $(\g,H)$ be an algebraic Harish-Chandra pair.
 Let $\bN$ be a right $\S^r(\g,H)$\semimodule, that is a Harish-Chandra
module over~$(\g,H)$.
 Then the action map $\g\ot_k\bN\rarrow\bN$ is a morphism of
Harish-Chandra modules.
 Tensoring it with $\det(\h)$ and passing to the $H$\invariant s, we
obtain a map $(\det(\h)\ot_k\g\ot_k\bN)^H\rarrow(\det(\h)\ot_k\bN)^H$.
 By the result of~\ref{finite-dim-semiinvariants-defining-map},
the composition $(\det(\h)\ot_k\h\ot_k\bN)^H\rarrow
(\det(\h)\ot_k\g\ot_k\bN)^H\rarrow(\det(\h)\ot_k\bN)^H$ vanishes.
 When $\bN$ is a coflat $\C(H)$\comodule, this provides a natural
map $(\det(\h)\ot_k\g/\h\ot_k\bN)^H\rarrow(\det(\h)\ot_k\bN)^H$;
to define this map in the general case, it suffices to represent $\bN$
as the kernel of a morphism of $\C$\+coflat{} $\S^r$\semimodule s
(see Lemma~\ref{coflat-semimodule-injection}).

 The vector space of \emph{$(\g,H)$\semiinvariant s} $\bN_{\g,H}$ is
defined as the cokernel of the map $(\det(\h)\ot_k\g/\h\ot_k\bN)^H
\rarrow(\det(\h)\ot_k\bN)^H$ that we have constructed.

\subsubsection{}
 Let $\bP$ be a left $\S^r(\g,H)$\semicontramodule{}
(see~\ref{hopf-harish-chandra-semi-contra-modules}).
 Then the action map $\bP\rarrow\Hom_k(\g,\bP)$ is a morphism of
$\S^r$\semicontramodule s.
 Applying to it the functor $\Hom_k(\det(\h),{-})$, we get a morphism
of $\C(H)$\contramodule s.
 Passing to the $H$\coinvariant s, i.~e., the maximal quotient
$\C$\contramodule s with the trivial contraaction, we obtain a map
$\Hom_k(\det(\h),\bP)_H\rarrow\Hom_k(\det(\h)\ot_k\g\;\bP)_H$.
 By the result of~\ref{finite-dim-semicontrainvs-defining-map},
the composition $\Hom_k(\det(\h),\bP)_H\rarrow\Hom_k(\det(\h)\ot_k\g\;
\bP)_H\rarrow\Hom_k(\det(\h)\ot_k\h\;\bP)_H$ vanishes.
 When $\bP$ is a coinjective $\C$\contramodule, this provides a natural
map $\Hom_k(\det(\h),\bP)_H\rarrow\Hom_k(\det(\h)\ot_k\g/\h\;\bP)_H$;
to define this map in the general case, it suffices to represent $\bP$
as the cokernel of a morphism of $\C$\coinjective{}
$\S^r$\semicontramodule s (see
Lemma~\ref{coproj-coinj-semi-mod-contra}).

 The vector space of \emph{$(\g,H)$\semicontrainvariant s} $\bP^{\g,H}$
is defined as the kernel of the map $\Hom_k(\det(\h),\bP)_H\rarrow
\Hom_k(\det(\h)\ot_k\g/\h\;\bP)_H$ that we have constructed.

\subsubsection{}  \label{finite-dim-semitensor-product-semiinvariants}
 Let $\bN$ be a right $\S^r(\g,H)$\semimodule{} and $\bM$ be a left
$\S^l(\g,H)$\semimodule; assume that either $\bN$ or $\bM$ is
a coflat $\C(H)$\comodule.
 Then there is a natural isomorphism $\bN\os_{\S^r}\bE\dual\os_{\S^r}\bM
 \simeq(\bN\ot_k\bM)_{\g,H}$, where $\bN\ot_k\bM$ is considered as
the tensor product of Harish-Chandra modules $\bN$ and $\bM$.

 Indeed, introduce an increasing filtration $F$ of the $k$\+algebra
$R=U(\g)$ whose component $F_tR$, \ $t=0$, $1$,~\dots is the linear span
of all products of elements of~$\g$ where at most $t$ factors do not
belong to~$\h$.
 In particular, we have $F_0R\simeq K=U(\h)$.
 Set $F_t\S^r=\C\ot_KF_tR$; then we have $F_0\S^r\simeq\C$, the natural
maps $F_{t-1}\S^r\rarrow F_t\S^r$ are injective, their cokernels are
coflat left and right $\C$\comodule s, $\S^r\simeq\ilim F_t\S^r$, and
the semimultiplication map $F_p\S^r\oc_\C F_q\S^r \rarrow\S\oc_\C\S
\rarrow\S$ factorizes through $F_{p+q}\S^r$.
 Moreover, the maps $F_p\S^r\oc_\C F_q\S^r\rarrow F_{p+q}\S^r$ are
surjective and their kernels are coflat left and right
$\C$\comodule s.  (Cf.~\ref{central-element-theorem}.)

 Let $\bN$ be a right $\S^r$\semimodule{} and $\bL$ be a left
$\S^r$\semimodule{} such that either $\bN$ or $\bL$ is a coflat
$\C$\comodule.
 Denote by $\eta_t\:\bN\oc_\C F_t\S^r\oc_\C\bL\rarrow\bN\oc_\C\bL$
the map equal to the difference of the map induced by the semiaction
map $F_t\S^r\oc_\C\bL\rarrow\bL$ and the map induced by the semiaction
map $\bN\oc_\C F_t\S^r\rarrow\bN$.
 Let us show that the images of $\eta_t$ coincide for $t\ge 1$.
 Let $p$, $q\ge 1$; then the map $\bN\oc_\C F_p\S^r\oc_\C F_q\S^r\oc_\C
\bL\rarrow\bN\oc_\C F_{p+q}\S^r\oc_\C\bL$ is surjective in view of
our assumption on $\bN$ and~$\bL$.
 The composition of the map $\bN\oc_\C F_p\S^r\oc_\C F_q\S^r\oc_\C\bL
\rarrow\bN\oc_\C F_{p+q}\S^r\oc_\C\bL$ with the map $\eta_{p+q}$
is equal to the sum of the composition of the map $\bN\oc_\C F_p\S^r
\oc_\C F_q\S^r\oc_\C\bL\rarrow\bN\oc_\C F_p\tS\oc_\C\bL$ induced by
the semiaction map $F_q\S^r\oc_\C\bL\rarrow\bL$ and the map $\eta_p$,
and the composition of the map $\bN\oc_\C F_p\S^r\oc_\C F_q\S^r\oc_\C\bL
\rarrow\bN\oc_\C F_q\S^r\oc_\C\bL$ induced by the semiaction map
$\bN\oc_\C F_p\S^r\rarrow\bN$ and the map~$\eta_q$.
 So the assertion follows by induction.
 Therefore, the semitensor product $\bN\os_\S\bL$ is isomorphic to
the cokernel of the map~$\eta_1$.

 On the other hand, the map $\eta_0$ vanishes.
 In view of our assumption on $\bN$ and $\bL$, the quotient space
$(\bN\oc_\C F_1\S^r\oc_\C\bL)/(\bN\oc_\C F_0\S^r\oc_\C\bL)$ is
isomorphic to $\bN\oc_\C F_1\S^r/F_0\S^r\oc_\C\bL$.
 Hence the semitensor product $\bN\os_\S\bL$ is isomorphic to
the cokernel of the induced map $\bar\eta_1\:\bN\oc_\C
(F_1\S^r/F_0\S^r)\oc_\C\bL\rarrow \bN\oc_\C\bL$.

 Now when $\bL=\bE\dual\os_\S\bM$ for a left $\S^l$\semimodule{} $\bM$,
the natural isomorphisms $\bN\oc_\C\E\dual\oc_\C\bM\simeq (\bN\ot_k
E\dual\ot_k\bM)^H\simeq(E\dual\ot_k\bN\ot_k\bM)^H$ and
$\bN\oc_\C (F_1\S^r/F_0\S^r)\oc_\C\E\dual\oc_\C\bM\simeq
\bN\oc_\C(\C\ot_k\g/\h)\oc_\C\E\dual\oc_\C\bM\simeq (\bN\ot_k\g/\h
\ot_kE\dual\ot_k\bM)^H\simeq(E\dual\ot_k\g/\h\ot_k\bN\ot_k\bM)^H$
given by the formulas $\check\imath\ot n\ot m\mpsto n\ot \check\imath
\ot m$ and $n\ot\bar z\ot \check\imath\ot m\mpsto \check\imath\ot
\bar z\ot n\ot m$ identify the map $\bar\eta_1$ with the map whose
cokernel is, by the definition, the space of semiinvariants
$(\bN\ot_k\bM)_{\g,H}$.

\subsubsection{}
 Let $\bM$ be a left $\S^l(\g,H)$\semimodule{} and $\bP$ be a left
$\S^r(\g,H)$\semicontramodule; assume that either $\bM$ is
an coprojective $\C(H)$\comodule, or $\bN$ is a coinjective
$\C(H)$\contramodule.
 Then there is a natural isomorphism $\SemiHom_{\S^r}(\bE\dual
\os_{\S^l}\bM\;\bP)\simeq\Hom_k(\bM,\bP)^{\g,H}$, where the structure
of left $\S^r$\semicontramodule{} on $\Hom_k(\bM,\bP)$ was introduced
in~\ref{hopf-semihom-factorizes-through}.
 The proof is analogous to that
of~\ref{finite-dim-semitensor-product-semiinvariants}.

\Section{Tate Harish-Chandra Pairs and Tate Lie Algebras}
\label{tate-appendix}

\centerline{by Sergey Arkhipov and Leonid Positselski}
\bigskip\medskip

 In order to formulate the comparison theorem relating the functors
$\SemiTor$ and $\SemiExt$ to the semi-infinite (co)homology of Tate
Lie algebras, one has to consider Harish-Chandra pairs $(\g,H)$ with
a Tate Lie algebra $\g$ and a proalgebraic group~$H$ corresponding
to a compact open subalgebra $\h\subset\g$.
 In such a situation, the construction of a Morita equivalence from
Appendix~\ref{hopf-appendix} no longer works; instead, there is
an isomorphism of ``left'' and ``right'' semialgebras correspoding
to different central charges.
 The proof of this isomorphism is based on the nonhomogeneous
quadratic duality theory developed in
Section~\ref{nonhom-koszul-section} (see also~\ref{koszul-over-ring}).
 Once the isomorphism of semialgebras is constructed and the standard
semi-infinite (co)homological complexes are introduced, the proof
of the comparison theorem becomes pretty straightforward.
 The equivalence between the semiderived categories of Harish-Chandra
modules and Harish-Chandra contramodules with complementary (or
rather, shifted) central charges follows immediately from
the isomorphism of semialgebras.

\subsection{Continuous coactions}

\subsubsection{}
 Let $k$ be a fixed ground field.
 A \emph{linear topology} on a vector space over~$k$ is a topology
compatible with the vector space structure for which open vector
subspaces form a base of neighborhoods of zero.
 In the sequel, by a topological vector space we will mean a $k$\+vector
space endowed with a complete and separated linear topology.
 Equivalently, a topological vector space is a filtered projective
limit of discrete vector spaces with its projective limit topology.
 Accordingly, the (separated) completion of a vector space endowed
with a linear topology is just the projective limit of its quotient
spaces by open vector subspaces.

 The category of topological vector spaces and continuous linear maps
between them has an exact category structure in which a triple of
topological vector spaces $V'\rarrow V\rarrow V''$ is exact if it is
an exact triple of vector spaces strongly compatible with
the topologies, i.~e., the map $V'\rarrow V$ is closed and the map
$V\rarrow V''$ is open.
 Any open surjective map of topological vector spaces is an admissible
epimorphism.
 Any closed injective map from a topological vector space
admitting a countable base of neighborhoods of zero is a split
admissible monomorphism.

 A topological vector space is called (\emph{linearly}) \emph{compact}
if it has a base of neighborhoods of zero consisting of vector
subspaces of finite codimension.
 Equivalently, a topological vector space is compact if it is
a projective limit of finite-dimensional discrete vector spaces.
 A \emph{Tate vector space} is a topological vector space admitting
a compact open subspace.
 Equivalently, a topological vector space is a Tate vector space if it
is topologically isomorphic to the direct sum of a compact vector space
and a discrete vector space.
 The \emph{dual Tate vector space} $V\dual$ to a Tate vector space $V$
is defined as the space of continuous linear functions $V\rarrow k$
endowed with the topology where annihilators of compact open subspaces
of $V$ form a base of neighborhoods of zero.
 In particular, the dual Tate vector spaces to compact vector spaces
are discrete and vice versa; for any Tate vector space $V$,
the natural map $V\rarrow (V\dual)\dual$ is a topological isomorphism.

\subsubsection{}
 The projective limit of a projective system of topological vector
spaces endowed with the topology of projective limit is a topological
vector space.
 This is called the \emph{topological projective limit}.

 The inductive limit of an inductive system of topological vector
spaces can be endowed with the topology of inductive limit of vector
spaces with linear topologies; we will call the inductive limit
endowed with this topology the \emph{uncompleted inductive limit}.
 The \emph{completed inductive limit} is the (separated) completion
of the uncompleted inductive limit.
 For any countable filtered inductive system formed by closed
embeddings of topological vector spaces the uncompleted and completed
inductive limits coincide.
 Morover, let $V_\alpha$ be a filtered inductive system of
topological vector spaces satisfying the following condition.
 For any increasing sequence of indices $\alpha_1\le\alpha_2\le\dsb$
the uncompleted inductive limit of $V_{\alpha_i}$ is a direct
summand of the uncompleted inductive limit of $V_\alpha$ considered
as an object of the category of vector spaces endowed with
noncomplete linear topologies.
 Then the uncompleted and completed inductive limits of~$V_\alpha$
coincide.

\subsubsection{}   \label{three-tensor-products}
 We will consider three operations of tensor product of topological
vector spaces~\cite{Beil}.
 For any two topological vector spaces $V$ and $W$, denote by $V\ot^!W$
the completion of the tensor product $V\ot_k W$ with respect to 
the topology with a base of neighborhoods of zero consisting of
the vector subspaces $V'\ot W + V\ot W'$, where $V'\subset V$ and
$W'\subset W$ are open vector subspaces in $V$ and~$W$.
 Furthermore, denote by $V\ot^*W$ the completion of $V\ot_k W$ with
respect to the topology formed by the subspaces of $V\ot W$
satisfying the following conditions: a vector subspace $T\subset V\ot W$
is open if (i)~there exist open subspaces $V'\subset V$, \ $W'\subset W$
such that $V'\ot W'\subset T$, (ii)~for any vector $v\in V$ there exists
a subspace $W''\subset W$ such that $v\ot W''\subset T$, and (iii)~for
any vector $w\in W$ there exists a subspace $V''\subset V$ such that
$V''\ot w\subset T$.
 Finally, denote by $V\vot W$ the completion of $V\ot_kW$ with respect
to the topology formed by the subspaces satisfying the following
conditions: a vector subspace $T\subset V\ot_k W$ is open if (i)~there
exists an open subspace $W'\subset W$ such that $V\ot_kW'\subset T$,
and (ii)~for any vector $w\in W$ there exists an open subspace
$V''\subset V$ such that $V''\ot w\subset T$.
 Set $V\wot W = W\vot V$.
 
 The topological tensor products $\ot^!$ and $\ot^*$ define two
structures of associative and commutative tensor category on 
the category of topological vector spaces.
 The topological tensor product $\vot$ defines a structure of
associative, but not commutative tensor category on the category of
topological vector spaces.
 For any topological vector spaces $V_1$,~\dots, $V_n$ and $W$
the vector space of continuous polylinear maps $V_1\times\dsb\times V_n
\rarrow W$ is naturally isomorphic to the vector space of continuous
linear maps $V_1\ot^*\dsb\ot^*V_n\rarrow W$.
 When both topological vector spaces $V$ and $W$ are compact (discrete),
the topological tensor product $V\ot^*W \simeq V\vot W\simeq V\wot W
\simeq V\ot^!W$ is also compact (discrete).
 The functor $\ot^!$ preserves topological projective limits.
 The functor $\ot^*$ preserves (uncompleted or completed) inductive
limits of filtered inductive systems of open injections.
 The topological tensor product $V\vot W$ is the topological
projective limit of $\vot$\+products of $V$ with discrete quotient
spaces of~$W$.
 The functor $(V,W)\mpsto V\vot W$ preserves completed inductive
limits in its second argument~$W$.
 The underlying vector space of the topological tensor product
$V\vot W$ is determined by (the topological vector space $W$ and)
the underlying vector space of the topological vector space~$V$.

 For Tate vector spaces $V_1$,~\dots, $V_n$ and a topological vector
space $U$, consider the vector space of continuous polylinear maps
$\prod_i V_i\rarrow U$ endowed with the topology with a base of
neighborhoods of zero formed by the subspaces of all polylinear maps
mapping the Carthesian product of a collection of compact subspaces
$V'_i\subset V_i$ into an open subspace $U'\subset U$
(the ``compact-open'' topology).
 This vector space is naturally topologically isomorphic to
the topological tensor product
$V_1\dual\ot^!\dsb\ot^!V_n\dual\ot^!W$~\cite{BD1}.
 For any topological vector spaces $U$, \ $W$ and Tate vector space $V$,
the vector space of continuous linear maps $V\ot^*W\rarrow U$ is
naturally isomorphic to the vector space of continuous linear maps
$W\rarrow V\dual\ot^!U$.

\subsubsection{}
 Let $\C$ be a coalgebra over the field~$k$ and $V$ be a topological
vector space.
 A \emph{continuous right coaction} of $\C$ in $V$ is a continuous
linear map $V\rarrow V\ot^!\C$, where $\C$ is considered as
a discrete vector space, satisfying the coassociativity and counity
equations.
 Namely, the map $V\rarrow V\ot^!\C$ should have equal compositions
with the two maps $V\ot^!\C\birarrow V\ot^!\C\ot^!\C$ induced by
the map $V\rarrow V\ot^!\C$ and the comultiplication in $\C$, and
the composition of the map $V\rarrow V\ot^!\C$ with the map
$V\ot^!\C\rarrow V$ induced by the counit of $\C$ should be equal
to the identity map.
 Equivalently, a continuous right coaction of $\C$ in $V$ can be defined
as a continuous linear map $V\ot^*\C\dual\rarrow V$, where $\C\dual$
is considered as a compact vector space, satisfying the associativity 
and unity equations.
 \emph{Continuous left coactions} are defined in the analogous way.

 A closed subspace $W\subset V$ of a topological vector space $V$
endowed with a continuous right coaction of a coalgebra $\C$ is
said to be invariant with respect to the continuous coaction
(or $\C$\invariant) if the image of $W$ under the continuous coaction
map $V\rarrow V\ot^!\C$ is contained in the closed subspace
$W\ot^!\C\subset V\ot^!\C$.
 It follows from the next Lemma that any topological vector space with
a continuous coaction of a coalgebra $\C$ is a filtered projective
limit of discrete vector spaces endowed with $\C$\comodule{} structures.

\begin{lem}
 For any topological vector space $V$ endowed with a continuous
coaction $V\rarrow V\ot^!\C$ of a coalgebra $\C$, open subspaces of~$V$
invariant under the continuous coaction form a base of neighborhoods of
zero in~$V$.
\end{lem}

\begin{proof}
 Let $U\subset V$ be an open subspace; then the full preimage $U'$ 
of the open subspace $U\ot^!\C\subset V\ot^!\C$ under the continuous
coaction map $V\rarrow V\ot^!\C$ is an invariant open subspace in $V$
contained in~$U$.
 To check that $U'$ is $\C$\invariant, use the fact the functor of
$\ot^!$\+product preserves kernels in the category of topological
vector spaces, and in particular, the $\ot^!$\+product with $\C$
preserves the kernel of the composition $V\rarrow V\ot^!\C
\rarrow V/U\ot^!\C$.
 To check that $U'$ is contained in~$U$, use the counity equation for
the continuous coaction.
\end{proof}

 The category of topological vector spaces endowed with a continuous
coaction of a coalgebra~$\C$ has an exact category structure such that
a triple of topological vector spaces with continuous coactions of~$\C$
is exact if and only if it is exact as a triple of topological
vector spaces.

 If $V$ is a Tate vector space with a continuous right coaction of $\C$,
then the dual Tate vector space $V\dual$ is endowed with a continuous
left coaction of~$\C$.

 Let $V$ be a topological vector space with a continuous right coaction
of a coalgebra $\C$ and $W$ be a topological vector space with
a continuous coaction of a coalgebra $\D$.
 Then all the three topological tensor products $V\ot^!W$, \ $V\ot^*W$,
and $V\vot W$ are endowed with continuous right coactions of
the coalgebra $\C\ot_k\D$.
 To construct the continuous coaction on $V\ot^!W$, one uses the natural
isomorphism $(V\ot^!\C)\ot^!(W\ot^!\D)\simeq (V\ot^!W)\ot^!(\C\ot_k\D)$.
 The continuous coaction on $V\ot^*W$ is defined in terms of
the natural continuous map $(V\ot^!\C)\ot^*(W\ot^!\D)\rarrow(V\ot^*W)
\ot^! (\C\ot_k\D)$, which exists for any topological vector spaces
$V$, $W$ and any discrete vector spaces $\C$, $\D$.
 The continuous coaction on $V\vot W$ is defined in terms of the natural
continuous map $(V\ot^!\C)\vot(W\ot^!\D)\rarrow
(V\vot W)\ot^!(\C\ot_k\D)$.

 It follows that for a commutative Hopf algebra $\C$ the topological
tensor products $V\ot^!W$, \ $V\ot^*W$, and $V\vot W$ of topological
vector spaces with continuous right coactions of $\C$ are also endowed
with continuous right coactions of~$\C$.
 Besides, one can transform a continuous left coaction of $\C$ in~$V$
into a continuous right coaction using the antipode.

 Now let $W$, \ $U$ be topological vector spaces and $V$ be a Tate
vector space; suppose that $W$, \ $U$, and $V$ are endowed with
continuous coactions of a commutative Hopf algebra~$\C$.
 Let $f\:V\ot^* W\rarrow U$ and $g\:W\rarrow V\dual\ot^!U$ be continuous
linear maps corresponding to each other under the isomorphism
from~\ref{three-tensor-products}; then $f$ preserves the continuous
coactions of $\C$ if and only if $g$~does.

\subsubsection{}
 A topological Lie algebra $\g$ is topological vector space
endowed with a Lie algebra structure such that the bracket is
a continuous bilinear map $\g\times\g\rarrow\g$.
 Topological associative algebras are defined in the analogous way.
 For example, let $V$ be a Tate vector space.
 Denote by $\End(V)$ the associative algebra of continuous endomorphisms
of~$V$ endowed with the compact-open topology and by $\gl(V)$ the Lie
algebra corresponding to $\End(V)$.
 Then $\End(V)$ is a topological associative algebra and $\gl(V)$ is
a topological Lie algebra.

 Let $U$, $V$, $W$ be topological vector spaces endowed with
continuous coactions of a commutative Hopf algebra~$\C$.
 Then a continuous bilinear map $V\times W\rarrow U$ is called
compatible with the continuous coactions of~$\C$ if the corresponding
linear map $V\ot^*W\rarrow U$ preserves the continuous coactions
of~$\C$.
 So one can speak about compatibility of continuous pairings, Lie or
associative algebra structures, Lie or associative actions, etc.,
with continuous coactions of a commutative Hopf algebra.
 
 Explicitly, a bilinear map $V\times W\rarrow U$ is continuous and
compatible with the continuous coactions of $\C$ if and only if
the following condition holds.
 For any $\C$\invariant{} open subspace $U'\subset U$ and any
finite-dimensional subspaces $E\subset V$, \ $F\subset W$ there
should exist invariant open subspaces $V'\subset V''\subset V$, \
$W'\subset W''\subset W$ such that $E\subset V''$, \ $F\subset W''$,
the map $V''\ot_k W''\rarrow U/U'$ factorizes through
$V''/V'\ot_kW''/W'$, and the induced map $V''/V'\ot_k W''/W'\rarrow
U/U'$ is a morphism of $\C$\comodule s.

\subsubsection{}
 For any Tate vector spaces $V$ and $W$, there is a split exact triple
of topological vector spaces $V\ot^*W\rarrow V\vot W\oplus W\vot V 
\rarrow V\ot^!W$, where the first map is the sum of the natural maps
$V\ot^*W\rarrow V\vot W$, \ $V\ot^*W\rarrow W\vot V$, while the second
map is the difference of the natural maps $V\vot W\rarrow V\ot^!W$, \
$W\vot V\rarrow V\ot^!W$.
 Let us take $W=V\dual$.
 Then $V\ot^!V\dual$ is naturally isomorphic to $\gl(V)$; the spaces
$V\vot V\dual$ and $V\dual\vot V$ can be identified with the subspaces
in $\gl(V)$ formed by the linear operators with open kernel and
compact closure of image, respectively; and $V\ot^* V\dual$ is
the intersection of $V\vot V\dual$ and $V\dual\vot V$ in $\gl(V)$.

 Taking the push-forward of the exact triple $V\ot^*V\dual\rarrow
V\vot V\dual\oplus V\dual\vot V\rarrow V\ot^!V\dual$ with respect to
the natural trace map $\tr\:V\ot^*V\dual\rarrow k$ corresponding to
the pairing $V\times V\dual\rarrow k$, one obtains an exact triple
of topological vector spaces $k\rarrow \gl(V)\til\rarrow\gl(V)$.
 This is also an exact triple of $\gl(V)$\module s, which allows to
define a Lie algebra structure on $\gl(V)\til$ making it a central
extension of the Lie algebra~$\gl(V)$.
 The anti-commutativity and the Jacobi identity follow from the fact
that the commutator of an operator with open kernel and an operator
with compact closure of image has zero trace.

 Now assume that a Tate vector space $V$ is endowed with a continuous
coaction of a commutative Hopf algebra $\C$.
 Then $V\ot^*V\dual\rarrow V\vot V\dual\oplus V\dual\vot V\rarrow
V\ot^!V\dual$ is an exact triple of topological vector spaces endowed
with continuous coactions of~$\C$; the trace map also preserves
the continuous coactions.
 Thus the topological vector space $\gl(V)\til$ acquires a continuous
coaction of~$\C$.

\subsubsection{}   \label{clifford-central-extension-definition}
 Here is another construction of the Lie algebra $\gl(V)\til$
(see~\cite{BD1}).
 Consider the quotient space of the vector space $V\ot_kV\dual\oplus
V\dual\ot_kV\oplus k$ by the relation $v\ot g+g\ot v=\lan g,v\ran$,
where $\lan\,\.,\,\ran$ denotes the pairing of $V\dual$ with~$V$.
 This vector space is a Lie subalgebra of the Clifford algebra
$\Cl(V\oplus V\dual)$ of the vector space $V\oplus V\dual$ with
the symmetric bilinear form given by the pairing $\lan\,\.,\,\ran$;
the Lie bracket on this subalgebra is given by the formulas
$[v_1\ot g_1\;v_2\ot g_2]= \lan g_1,v_2\ran v_1\ot g_2 - 
\lan g_2,v_1\ran v_2\ot g_1$, \ $[v\ot g\; 1]=0$.
 This Lie algebra acts in the vector space $V$ by the formulas
$(v\ot g)(v')=\lan g,v'\ran v$, \ $1(v)=0$.
 There is a separated topology on this Lie algebra with a base of
neighborhoods of zero formed by the Lie subalgebras $V\ot W'+
V'\ot V\dual$, where $V'\subset V$ and $W'\subset V\dual$ are open
subspaces such that $\lan W',V'\ran=0$.
 The completion of this Lie algebra with respect to this topology can
be easily identified with the Lie algebra $\gl(V)\til$ defined above.

 Hence the Lie bracket on $\gl(V)\til$ is continuous.
 In addition, we need to check that when $V$ is endowed with
a continuous coaction of a commutative Hopf algebra~$\C$,
the Lie bracket is compatible with the continuous coaction of $\C$
in $\gl(V)\til$.
 The latter follows from the existence of a well-defined commutator map
$$
 \Hom(X_4,X_3,X_1;\medspace X/X_1,X_4/X_1,X_2/X_1)^{\ot 2}\rarrow
 \gl(X)\til/(X\ot X_3^\perp+X_2\ot X\dual)
$$
for any flag of finite-dimensional vector spaces $X_1\subset X_2\subset
X_3 \subset X_4\subset X$, where the $\Hom$ space in the left hand side
consists of all maps $X_4\rarrow X/X_1$ sending $X_3$ to $X_4/X_1$ and
$X_1$ to~$X_2/X_1$, and $Y^\perp\subset X\dual$ denotes the orthogonal
complement to a vector subspace $Y\subset X$.

\subsubsection{}
 A \emph{Tate Lie algebra} is a Tate vector space endowed with
a topological Lie algebra structure.
 Let $\g$ be a Tate Lie algebra endowed with a continuous coaction
of a commutative Hopf algebra $\C$ such that the Lie algebra structure
is compatible with the continuous coaction.
 Then $\C$\invariant{} compact open subalgebras form a base of
neighborhoods of zero in~$\g$.

 Indeed, choose a $\C$\invariant{} compact open subspace $U\subset\g$;
let $\h$ be the normalizer of $U$ in $\g$, i.~e., the subspace
of all $x\in\g$ such that $[x,U]\subset U$.
 Then $\h$ is a $\C$\invariant{} open subalgebra in $\g$, since it
is the kernel of the adjoint action map $\g\rarrow\Hom_k(U,\g/U)$.
 Therefore, the intersection $\h\cap U$ is an $\C$\invariant{} compact
open subalgebra in~$\g$ contained in~$U$.

 The canonical central extension $\g\til$ of a Tate Lie algebra~$\g$
is defined as the fibered product of $\g$ and $\gl(\g)\til$ over
$\gl(\g)$, where $\g$ maps to $\gl(\g)$ by the adjoint representation.
 The vector space $\g\til$ is endowed with the topology of fibered
product; this makes $\g\til$ a Tate Lie algebra.
 The central extension $\g\til\rarrow\g$ splits canonically and
continuously over any compact open Lie subalgebra $\h\subset\g$.
 Indeed, the image of $\h$ in $\gl(\g)$ is contained in the open Lie
subalgebra $\gl(\h,\g)\subset\gl(\g)$ of endomorphisms preserving~$\h$, 
and the map from the open Lie subalgebra of $\gl(\g)\til$ constructed
as the completion of $\g\ot_k\h^\perp+\g\dual\ot_k\h$ onto $\gl(\h,\g)$
is a topological isomorphism.

 The natural continuous coaction of $\C$ in $\g\til$ is constructed
as the fibered product of the coactions in $\g$ and $\gl(\g)\til$;
it is clear that the Lie algebra structure on $\g\til$ is compatible
with the continuous coaction.
 If $\h\subset\g$ is a $\C$\invariant{} compact open subalgebra,
then the canonical splitting $\h\rarrow\g\til$ preserves 
the continuous coactions.

 When a Tate vector space $V$ is decomposed into a direct sum
$V\simeq E\oplus F$ of a compact vector space $E$ and a discrete vector
space $F$, there is a natural section $\gl(V)\rarrow\gl(V)\til$
of the central extension $\gl(V)\til\rarrow\gl(V)$; the image of
this section is the completion of $V\ot F\dual+F\ot V\dual+V\dual\ot E
+E\dual\ot V$.
 Consequently, when a Tate Lie algebra $\g$ is decomposed into a direct
sum $\g\simeq\h\oplus \b$ of a compact open Lie subalgebra $\h$ and
a discrete vector subspace $\b$, there is a natural section
$\g\rarrow\g\til$ of the central extension $\g\til\rarrow\g$; this
section agrees with the natural splitting $\h\rarrow\g\til$.
 
\subsection{Construction of semialgebra}

\subsubsection{}   \label{tate-harish-chandra-pair-defined}
 We will sometimes use Sweedler's notation~\cite{Swe} \
$c\mpsto c_{(1)}\ot c_{(2)}$ for the comultiplication map in
a coassociative coalgebra~$\C$.
 The analogous notation for coactions of $\C$ in a right
$\C$\comodule{} $\N$ and a left $\C$\comodule{} $\M$ is
$n\mpsto n_{(0)}\ot n_{(1)}$ and $m\mpsto m_{(-1)}\ot m_{(0)}$,
where $n$, $n_{(0)}\in\N$, \ $m$, $m_{(0)}\in\M$, and
$n_{(1)}$, $m_{(-1)}\in\C$.

 A \emph{Lie coalgebra} $\L$ is a $k$\+vector space endowed with
a $k$\+linear map $\L\rarrow \bigwedge^2_k\L$ from $\L$ to the second
exterior power of $\L$ denoted by $l\mpsto l_{\{1\}}\wedge l_{\{2\}}$,
which should satisfy the dual version of Jacobi identity
$l_{\{1\}\{1\}}\wedge l_{\{1\}\{2\}}\wedge l_{\{2\}} = l_{\{1\}}\wedge
l_{\{2\}\{1\}}\wedge l_{\{2\}\{2\}}$, where $l'\wedge l''\wedge l'''$
denotes an element of $\bigwedge^3_k\L$.
 A \emph{comodule} $\M$ over a Lie coalgebra $\L$ is a $k$\+vector 
space endowed with a $k$\+linear map $\M\rarrow \L\otimes\M$ denoted
by $m\mpsto m_{\{-1\}}\ot m_{\{0\}}$ satisfying the equation
$m_{\{-1\}}\wedge m_{\{0\}\{-1\}}\ot m_{\{0\}\{0\}} = m_{\{-1\}\{1\}}
\wedge m_{\{-1\}\{2\}}\ot m_{\{0\}}$, where $l'\wedge l''\ot m$ denotes
an element of $\bigwedge^2_k\L\ot_k\M$.

 A \emph{Tate Harish-Chandra pair} $(\g,\C)$ is a set of data
consisting of a Tate Lie algebra~$\g$, a commutative Hopf algebra $\C$,
a continuous coaction of $\C$ in $\g$ such that the Lie algebra
structure on~$\g$ is compatible with the continuous coaction,
a $\C$\invariant{} compact open subalgebra $\h\subset\g$, and
a continuous pairing $\psi\:\C\times
\h\rarrow k$, where $\C$ is considered with the discrete topology.
 This data should satisfy the following conditions (cf.~\cite{BFM}):
\begin{itemize}
 \item[(i)] The pairing~$\psi$ is compatible with the multiplication and
comultiplication in $\C$, i.~e., the map $\check\psi\:\C\rarrow\h\dual$
corresponding to~$\psi$ is a morphism of Lie coalgebras such that
$\check\psi(c'c'')=\eps(c')\check\psi(c'')+\eps(c'')\check\psi(c')$ for
$c'$, $c''\in\C$.
 Here the Lie coalgebra structure on $\C$ is defined by the formula
$c\mpsto c_{(1)}\wedge c_{(2)}$ and the Lie coalgebra structure on
$\h\dual$ is given by the formula $\lan x^*, [x',x'']\ran =
\lan x^*_{\{1\}},x''\ran \lan x_{\{2\}},x'\ran - \lan x^*_{\{1\}},x'\ran 
\lan x^*_{\{2\}},x''\ran$ for $x^*\in\h\dual$, \ $x'$, $x''\in\h$.
 By $\eps$ we denote the counit of~$\C$.
 \item[(ii)] The pairing~$\psi$ is compatible with the continuous
coaction of~$\C$ in $\h$ obtained by restricting the coaction in~$\g$
and the adjoint coaction of $\C$ in itself.
 The latter is defined by the formula $c\mpsto c_{[0]}\ot c_{[1]} =
c_{(2)}\ot s(c_{(1)})c_{(3)}$, where $s$ denotes the antipode map of
the Hopf algebra $\C$ (the square brackets are used to
avoid ambiguity of notation).
 The compatibility means that the continuous linear map
$\C\ot^*\h\rarrow k$ corresponding to~$\psi$ preserves the continuous
coactions, or equivalently, the map $\check\psi$ is a morphism of
$\C$\comodule s.
 \item[(iii)] The action of $\h$ in $\g$ induced by the continuous
coaction of $\C$ in $\g$ and the pairing~$\psi$ coincides with
the adjoint action of $\h$ in~$\g$.
 Here the former action is constructed as the projective limit of
the actions of $\h$ in quotient spaces of $\g$ by $\C$\invariant{}
open subspaces; for a right $\C$\comodule{} $\N$, the $\h$\module{}
structure on $\N$ induced by the pairing~$\psi$ is defined by
the formula $xn=-\psi(n_{(1)},x)n_{(0)}$ for $x\in\h$, \ $n\in\N$.
\end{itemize}

 Given a Tate Harish-Chandra pair $(\g,\C)$, one can construct
a Tate Harish-Chandra pair $(\g\til,\C)$ with the same Lie
subalgebra~$\h$, where $\g\til$ is the canonical central extension
of a Tate Lie algebra~$\g$.
 A continuous coaction of $\C$ in $\g\til$ and a canonical embedding
of $\h$ into $\g$ preserving the continuous coactions of $\C$ were
constructed above; it remains to check the condition~(iii).
 Here it suffices to notice that the adjoint action of $\gl(\g)$
in $\gl(\g)\til$ coincides with the action of $\gl(\g)$ in $\gl(\g)\til$
induced by the action of $\gl(\g)$ in $\g$, hence the adjoint action of
$\h$ in $\gl(\g)\til$ coincides with the action of $\h$ in $\gl(\g)\til$
induced by the coaction of $\C$ in $\gl(\g)\til$ and the pairing~$\psi$.

\subsubsection{}  \label{tate-harish-chandra-central-extension}
 Let $(\g',\C)$ be a Tate Harish-Chandra pair such that the Tate Lie
algebra $\g'$ is a central extension of a Tate Lie algebra $\g$ with
the kernel identified with~$k$; assume that $\C$ coacts trivially on
$k\subset\g'$ and the Lie subalgebra $\h\subset \g'$ that is a part of
the Tate Harish-Chandra pair structure does not contain~$k$.
 Then $(\g,\C)$ is naturally also a Tate Harish-Chandra pair with
the induced continuous coaction of $\C$ in~$\g$ and the Lie subalgebra
$\h\subset\g$ defined as the image of $\h$ in~$\g$.
 In this case, we will say that $(\g',\C)\rarrow(\g,\C)$ is a central
extension of Tate Harish-Chandra pairs with the kernel~$k$.
 One example of a central extension of Tate Harish-Chandra pairs is
the canonical central extension $(\g\til,\C)\rarrow(\g,\C)$.

 Let $\kap\:(\g',\C)\rarrow(\g,\C)$ be a central extension of Tate
Harish-Chandra pairs with the kernel~$k$.
 Consider the tensor product $\S^r_\kap(\g,\C)=\C\ot_{U(\h)}U_\kap(\g)$,
where $U(\h)$ and $U(\g')$ denote the universal enveloping algebras of
the Lie algebras $\h$ and $\g'$ considered as Lie algebras without any
topologies, $U_\kap(\g)=U(\g')/(1_{U(\g')}-1_{\g'})$ is the modification
of the universal enveloping algebra of $\g$ corresponding to the central
extension $k\rarrow\g'\rarrow g$, and $1_{U(\g')}$ and $1_{\g'}$ denote
the unit elements of the algebra $U(\g')$ and the vector subspace
$k\subset\g'$, respectively.
 The structure of right $U(\h)$\module{} on $\C$ comes from the pairing
$\phi\:\C\ot_kU(\h)\rarrow k$ corresponding to the algebra morphism
$U(\h)\rarrow\C\dual$ induced by the Lie algebra morphism
$\Check{\Check\psi}\:\h\rarrow\C\dual$, where the multiplication on
$\C^*$ is defined by the formula $\lan c'{}^*c''{}^*,c\ran =
\lan c'{}^*,c_{(2)}\ran \lan c''{}^*,c_{(1)}\ran$ for $c'{}^*$,
$c''{}^*\in\C^*$, \ $c\in\C$ and the Lie bracket is given by
the formula $[c'{}^*,c''{}^*]=c'{}^*c''{}^*-c''{}^*c'{}^*$.

 We claim that the vector space $\S^r_\kap(\g,\C)$ has a natural
structure of semialgebra over the coalgebra~$\C$ provided by
the general construction of~\ref{semialgebra-constructed}.
 The construction of this semialgebra structure becomes a little
simpler if one assumes that
\begin{itemize}
\item[(iv)] the pairing $\phi\:\C\ot_kU(\h)\rarrow k$ is nondegenerate
in~$\C$,
\end{itemize}
but this is not necessary.

\subsubsection{}
 To construct a right $\C$\comodule{} structure on $\S^r_\kap(\g,\C)$,
we will have to approximate this vector space by finite-dimensional
spaces.
 Let $V_1$,~\dots, $V_t$ be a sequence of $\C$\invariant{} compact
open subspaces of $\g'$ containing $\h$ and~$k$ such that $V_i+
[V_i,V_i]\subset V_{i-1}$.
 Let $\N$ be a finite-dimensional right $\C$\comodule.
 Choose a $\C$\invariant{} compact open subspace $W_1\subset\h$ such
that the $\C$\comodule{} $\N$ is annihilated by the action of $W_1$
obtained by restricting the action of~$\h$ induced by the pairing~$\psi$.
 For each $i=2$,~\dots, $t$ choose a $\C$\invariant{} compact open
subspace $W_i\subset \h$ such that $W_i+[V_i,W_i]\subset W_{i-1}$.
 Denote by $\S^r_\kap(V_1,\dsc,V_t;\N)$ the quotient space of the vector
space $\N\ot_k(k\oplus V_1/W_1\oplus\dsb\oplus(V_t/W_t)^{\ot t})$ by
the obvious relations imitating the relations in the enveloping algebra
$U_\kap(\g)$ and its tensor product with $\N$ over $U(\h)$.
 It is easy to see that this quotient space does not depend on
the choice of the subspaces~$W_i$.
 In other words, denote by $R(V_1,\dsc,V_t)$ the subspace
$U(\h)(k+V_1+\dsb+V_t^t)\subset U_\kap(\g)$; it is
an $U(\h)$\+$U(\h)$\+subbimodule of $U_\kap(\g)$ and a free left
$U(\h)$\module.
 The tensor product $\N\ot_{U(\h)}R(V_1,\dsc,V_t)$ is naturally
isomorphic to $\S^r_\kap(V_1,\dsc,V_t;\N)$.
 This is an isomorphism of right $U(\h)$\module s; when $\N=\D$ is
a finite-dimensional subcoalgebra of~$\C$, this is also an isomorphism
of left $\C$\comodule s.
 Clearly, the inductive limit of $\S^r_\kap(V_1,\dsc,V_t;\D)$ over
increasing $t$, \ $V_i$, and finite-dimensional subcoalgebras
$\D\subset\C$ is naturally isomorphic to $\S^r_\kap(\g,\C)$.
 
 Now the vector space $\S^r_\kap(V_1,\dsc,V_t;\N)$ has a right
$\C$\comodule{} structure induced by the right $\C$\comodule{}
structure on $\N\ot_k(k\oplus V_1/W_1\oplus\dsb\oplus(V_t/W_t)^{\ot t})$
obtained by taking the tensor product of the $\C$\comodule{} structures
on $V_i/W_i$ and the right $\C$\comodule{} structure on~$\D$.
 The inductive limit of these $\C$\comodule{} structures for $\N=\D$
provides the desired right $\C$\comodule{} structure on
$\S^r_\kap(\g,\C)$.
 It commutes with the left $\C$\comodule{} structure on
$\S^r_\kap(\g,\C)$ and agrees with the right $U(\h)$\module{}
structure, since such commutativity and agreement hold on the level
of the spaces $\S^r_\kap(V_1,\dsc,V_t;\D)$.
 Furthermore, by the (classical) Poincare--Birkhoff--Witt theorem
$U_\kap(\g)$ is a free left $U(\h)$\module.
 If the condition~(iv) holds, the construction of the semialgebra
$\S^r_\kap(\g,\C)$ is finished; otherwise, we still have to check
that the semiunit map $\C\rarrow\S^r_\kap(\g,\C)$ and
the semimultiplication map $\S^r_\kap(\g,\C)\oc_\C\S^r_\kap(\g,\C)
\rarrow\S^r_\kap(\g,\C)$ are morphisms of right $\C$\comodule s.

 The former is clear, and the latter can be proven in the following way.
 Any finite-dimensional $\C$\comodule{} $\N$ is a comodule over
a finite-dimensional subcoalgebra $\E\subset\C$.
 There is a natural isomorphism $\N\ot_{U(\h)}R(V_1,\dsc,V_t)\simeq
\N\oc_\C(\E\ot_{U(\h)}R(V_1,\dsc,V_t))$.
 The corresponding isomorphism $\S^r_\kap(V_1,\dsc,V_t;\N)\simeq
\N\oc_\C\S^r_\kap(V_1,\dsc,V_t;\E)$, which is induced by
the isomorphism $\N\simeq\N\oc_\C\E$, preserves the right
$\C$\comodule{} structures.
 All of this is applicable to the case of $\N=
\S^r_\kap(V'_1,\dsc,V'_t;\D)$, where $V'_1$,~\dots, $V'_t$ is another
sequence of subspaces of $\g'$ satisfying the above conditions.
 Now let $V''_1$,~\dots, $V''_{2t}\subset\g'$ be a sequence of 
subspaces satisfying the above conditions and such that
$V'_i$, $V_i\subset V''_{t+i}$.
 The map $\C\ot_{U(\h)}U_\kap(\g)\ot_{U(\h)}U_\kap(\g)\rarrow 
\C\ot_{U(\h)}U_\kap(\g)$ induced by the multiplication map
$U_\kap(\g)\ot_{U(\h)}U_\kap(\g)\rarrow U_\kap(\g)$ is the inductive
limit of the maps $\D\ot_{U(\h)}R(V'_1,\dsc,V'_t)\ot_{U(\h)}
R(V_1,\dsc,V_t)\rarrow \D\ot_{U(\h)}R(V''_1,\dsc,V''_{2t})$ over
increasing $t$, \ $V_i$, $V'_i$, $V''_i$, and~$\D$.
 The corresponding map $\S^r_\kap(V_1,\dsc,V_t;
\S^r_\kap(V'_1,\dsc,V'_t;\D))\rarrow\S^r_\kap(V''_1,\dsc,V''_{2t};\D)$
is induced by the map $\D\ot_k(k\oplus V'_1/W'_1\oplus\dsb\oplus
(V'_t/W'_t)^{\ot t})\ot_k(k\oplus V_1/W_1\oplus\dsb\oplus
(V_t/W_t)^{\ot t})\rarrow\D\ot_k(k\oplus V''_1/W''_1\oplus\dsb\oplus
(V''_{2t}/W''_{2t})^{\ot 2t})$, where the sequences of subspaces
$W'_i$, $W_i$, $W''_i$ satisfy the above conditions with respect to
the sequences of subspaces $V'_i$, $V_i$, $V''_i$, and the right
$\C$\comodule s $\D$, \ $\D\ot_k(k\oplus V'_1/W'_1\oplus\dsb\oplus
(V'_t/W'_t)^{\ot t})$, \ $\D$, respectively, and the additional
condition that $W'_i$, $W_i\subset W''_{t+i}$.
 One can easily see that the latter map is a morphism of right
$\C$\comodule s.

\subsubsection{}
 The semialgebra $\S^r_\kap(\g,\C)$ over the coalgebra~$\C$ is
constructed.
 Analogously one defines a semialgebra structure on the tensor product
$\S^l_\kap(\g,\C)=U_\kap(\g)\ot_{U(\h)}\C$.
 The semialgebras $\S^r_\kap=\S^r_\kap(\g,\C)$ and $\S^l_\kap=
\S^l_\kap(\g,\C)$ are essentially opposite to each other
(see~\ref{hopf-morita-anti-involution}). 
 More precisely, the antipode anti-automorphisms of $U(\g')$ and $\C$
induce a natural isomorphism of semialgebras $\S^r_\kap\simeq
\S^{l\mskip.75\thinmuskip\rop}_{-\kap}$ compatible with the isomorphism
of coalgebras  $\C^\rop\simeq\C$, where we denote by $-\kap$
the central extension of Tate Harish-Chandra pairs with the kernel~$k$
that is obtained from the central extension $\kap$ by multiplying
the embedding $k\rarrow\g'$ with~$-1$.

\subsubsection{}   \label{tate-harish-chandra-modules}
 A \emph{discrete module} $M$ over a topological Lie algebra $\g$ is 
$\g$\module{} such that the action map $\g\times M\rarrow M$ is
continuous with respect to the discrete topology of $M$.
 Equivalently, a $\g$\module{} $M$ is discrete if the annihilator of
any element of $M$ is an open Lie subalgebra in~$\g$.
 In particular, if $\psi\:\C\times\h\rarrow k$ is a continuous pairing
between a compact Lie algebra $\h$ and a coalgebra $\C$ such that
the map $\check\psi\:\C\rarrow\h\dual$ is a morphism of Lie coalgebras,
then the $\h$\module{} structure induced by a $\C$\comodule{} structure
by the formula of~\ref{tate-harish-chandra-pair-defined}\.(iii) is 
always discrete.

 Let $\kap\:(\g',\C)\rarrow(\g,\C)$ be a central extension of Tate
Harish-Chandra pairs with the kernel~$k$.
 Then the category of left semimodules over $\S^l_\kap(\g,\C)$ is
isomorphic to the category of $k$\+vector spaces $\bM$ endowed with
$\C$\comodule{} and discrete $\g'$\module{} structures such that
the induced discrete $\h$\module{} structures coincide, the action map
$\g/U\ot_k\L \rarrow \bM$ is a morphism of $\C$\comodule s for any
finite-dimensional $\C$\subcomodule{} $\L\subset\bM$ and any
$\C$\invariant{} compact open subspace $U\subset\g$ annihilating~$\L$,
and the unit element of $k\subset\g'$ acts by the identity in~$\bM$.
 The second of these three conditions can be reformulated as follows:
for any $\C$\+invariant compact subspace $V\subset\g'$, the natural Lie
coaction map $\bM\rarrow V\dual\ot_k\bM$ is a morphism of
$\C$\comodule s.
 When the assumption~(iv) of~\ref{tate-harish-chandra-central-extension}
is satisfied, the second condition is redundant.

 Abusing terminology, we will call vector spaces $\bM$ endowed with
such a structure \emph{Harish-Chandra modules over\/ $(\g,\C)$ with
the central charge\/~$\kap$}.
 Analogously, the category of right semimodules over $\S^r_\kap(\g,\C)$
is isomorphic to the category of Harish-Chandra modules over $(\g,\C)$
with central charge~$-\kap$.

\subsubsection{}  \label{lie-contramodules}
 For a topological vector space $V$ and a vector space $P$, denote by
$V\ot\comp P$ the tensor product $V\ot^!P = V\wot P$ \emph{considered
as a vector space without any topology}, where $P$ is endowed with
the discrete topology for the purpose of making the topological tensor
product.
 In other words, one has $V\ot\comp P=\plim_U V/U\ot_k P$, where
the projective limit is taken over all open subspaces $U\subset V$.

 For a topological vector space~$V$, denote by $\bigwedge^{*,2}(V)$
the completion of $\bigwedge^2_k(V)$ with respect to the topology
with the base of neighborhoods of zero formed by all the subspaces
$T\subset\bigwedge^2_k(V)$ such that there exists an open subspace
$V'\subset V$ for which $\bigwedge^2_k(V')\subset T$ and for any
vector $v\in V$ there exists an open subspace $V''\subset V$ for
which $v\wedge V''\subset T$.
 For any topological vector spaces $V$ and~$W$, the vector space of
continuous skew-symmetric bilinear maps $V\times V\rarrow W$ is
naturally isomorphic to the vector space of continuous linear maps
$\bigwedge^{2,*}(V)\rarrow W$.
 The space $\bigwedge^{*,2}(V)$ is a closed subspace of the space
$V\ot^*V$; the skew-symmetrization map $V\ot^*V\rarrow V\ot^*V$
factorizes through $\bigwedge^{*,2}(V)$.

 Let $\g$ be a topological Lie algebra.
 A \emph{contramodule} over $\g$ is a vector space $P$ endowed with
a linear map $\g\ot\comp P\rarrow P$ satisfying the following version
of Jacobi equation.
 Consider the vector space $\bigwedge^{*,2}(\g)\ot\comp P$.
 There is a natural map $\bigwedge^{*,2}(\g)\ot\comp P\rarrow\g\ot\comp
P$ induced by the bracket map $\bigwedge^{*,2}(\g)\rarrow\g$.
 Furthermore, there is a natural map $(\g\ot^*\g)\ot\comp P\rarrow
\g\ot\comp(\g\ot\comp P)$, which is constructed as follows.
 For any open subspace $U\subset\g$ there is a natural surjection
$(\g\ot^*\g)\ot\comp P\rarrow (\g/U\ot^*\g)\ot\comp P$ and for any
discrete vector space $F$ there is a natural isomorphism
$(F\ot^*\g)\ot\comp P\simeq F\ot_k(\g\ot\comp P)$, so the desired map
is obtained as the projective limit over~$U$.
 Composing the map $\bigwedge^{*,2}(\g)\ot\comp P\rarrow(\g\ot^*\g)
\ot\comp P$ induced by the embedding $\bigwedge^{*,2}(\g)\rarrow
\g\ot^*\g$ with the map $(\g\ot^*\g)\ot\comp P\rarrow\g\ot\comp
(\g\ot\comp P)$ that we have constructed and with the map
$\g\ot\comp(\g\ot\comp P)\rarrow\g\ot\comp P$ induced by
the contraaction map $\g\ot\comp P\rarrow P$, we obtain a second map
$\bigwedge^{*,2}(\g)\ot\comp P\rarrow P$.
 Now the contramodule Jacobi equation claims that the two maps
$\bigwedge^{*,2}\ot\comp P\rarrow\g\ot\comp P$ should have equal
compositions with the contraaction map $\g\ot\comp P\rarrow P$.

 Alternatively, the map $(\g\ot^*\g)\ot\comp P\rarrow
\g\ot\comp(\g\ot\comp P)$ can be constructed as the composition
$(\g\ot^*\g)\ot\comp P\rarrow(\g\wot\g)\ot\comp P\simeq
\g\wot\g\wot P\simeq\g\ot\comp(\g\ot\comp P)$ of the map induced by
the natural continuous map $\g\ot^*\g\rarrow\g\wot\g$ and the natural
isomorphisms whose existence follows from the fact that the topological
tensor product $W\wot V$ considered as a vector space without any
topology does not depend on the topology of~$V$.
 The following comparison between the definitions of a discrete
$\g$\+module and a $\g$\+contramodule can be made: a discrete
$\g$\+module structure on a vector space $M$ is given by a continuous
linear map $\g\ot^*M\simeq\g\vot M\rarrow M$, while
a $\g$\+contramodule structure on a vector space $P$ is given by
a discontinuous linear map $\g\ot^!P\simeq\g\wot P\rarrow P$,
where $M$ and $P$ are endowed with discrete topologies.
 For any topological Lie algebra $\g$, the category of
$\g$\contramodule s is abelian (cf.~\ref{assoc-contramodules}).
 There is a natural exact forgetful functor from the category of
$\g$\contramodule s to the category of modules over the Lie
algebra~$\g$ considered without any topology.

 For any discrete $\g$\module{} $M$ and any vector space $E$ there is
a natural structure of $\g$\contramodule{} on the space of linear maps
$\Hom_k(M,E)$.
 The contraaction map $\g\ot\comp\Hom_k(M,E)\rarrow \Hom_k(M,E)$
is constructed as the projective limit over all open subspaces
$U\subset \g$ of the maps $\g/U\ot_k\Hom_k(M,E)\rarrow\Hom_k(M^U,E)$
given by the formula $\bar z \ot g\mpsto (m\mpsto -g(\bar zm))$
for $\bar z \in\g/U$, \ $g\in\Hom_k(M,E)$, and $m\in M^U$, where
$M^U\subset M$ denotes the subspace of all elements of~$M$ annihilated
by~$U$.

 More generally, for any discrete module $M$ over a topological Lie
algebra $\g_1$ and any contramodule{} $P$ over a topological Lie
algebra $\g_2$ there is a natural $(\g_1\oplus\g_2)$\contramodule{}
structure on $\Hom_k(M,P)$ with the contraaction map
$(\g_1\oplus\g_2)\ot\comp\Hom_k(M,P)\rarrow \Hom_k(M,P)$ defined
as the sum of two commuting contraactions of $\g_1$ and $\g_2$ in 
$\Hom_k(M,P)$, one of which is introduced above and the other one is
given by the composition $\g_2\ot\comp\Hom_k(M,P)\rarrow
\Hom_k(M\;\g_2\ot\comp P)\rarrow \Hom_k(M,P)$ of the natural map
$\g_2\ot\comp\Hom_k(M,P)\rarrow\Hom_k(M\;\g_2\ot\comp P)$
and the map $\Hom_k(M\;\g_2\ot\comp P)\rarrow \Hom_k(M,P)$ induced by
the $\g_2$\+contraaction in~$P$.
 Hence for any discrete $\g$\module{} $M$ and any $\g$\contramodule{}
$P$ there is a natural $\g$\contramodule{} structure on $\Hom_k(M,P)$
induced by the diagonal embedding of Lie algebras $\g\rarrow\g\oplus\g$.

\subsubsection{}   \label{tate-lie-contramodules}
 When $\g$ is a Tate Lie algebra, a $\g$\contramodule{} $P$ can be also
defined as a $k$\+vector space endowed with a linear map $\Hom_k
(V\dual,P)\rarrow P$ for every compact open subspace $V\subset \g$.
 These linear maps should satisfy the following two conditions:
when $U\subset V\subset \g$ are compact open subspaces, the maps
$\Hom_k(U\dual,P)\rarrow P$ and $\Hom_k(V\dual,P)\rarrow P$ should
form a commutative diagram with the map $\Hom_k(U\dual,P)\rarrow
\Hom_k(V\dual,P)$ induced by the natural surjection $V\dual\rarrow
U\dual$, and for any compact open subspaces $V'$, $V''$, $W\subset \g$
such that $[V',V'']\subset W$ the composition $\Hom_k(V''{}\dual\ot_k
V'{}\dual\;P)\rarrow\Hom_k(W\dual,P)\rarrow P$ of the map induced by
the cobracket map $W\dual\rarrow V''{}\dual\ot_kV'{}\dual$ and
the contraaction map $\Hom_k(W\dual,P)\rarrow P$ should be equal to
the difference of the iterated contraaction map $\Hom_k(V''{}\dual\ot_k
V'{}\dual\;P)\simeq\Hom_k(V'{}\dual,\Hom_k(V''{}\dual,P))\rarrow
\Hom_k(V'{}\dual,P)\rarrow P$ and the composition of the isomorphism
$\Hom_k(V''{}\dual\ot_k V'{}\dual\;P)\simeq \Hom_k(V'{}\dual\ot_k
V''{}\dual\;P)$ induced by the isomorphism $V'{}\dual\ot_kV''{}\dual
\simeq V''{}\dual\ot_kV'{}\dual$ with the iterated contraaction map
$\Hom_k(V'{}\dual\ot_k V''{}\dual\;P)\simeq\Hom_k(V''{}\dual,
\Hom_k(V'{}\dual,P))\rarrow P$.

 For a Tate Lie algebra $\g$, a discrete $\g$\module{} $M$, and
an $\g$\contramodule{} $P$, the structure of $\g$\contramodule{}
on $\Hom_k(M,P)$ defined above is given by the formula
$\pi(g)(m)=\pi_P(x^*\mapsto g(x^*)(m)) - g(m_{\{-1\}})(m_{\{0\}})$
for a compact open subspace $V\subset\g$, a linear map
$g\in\Hom_k(V\dual,\Hom_k(M,P))$, and elements $x^*\in V\dual$, \
$m\in M$, where $m\mpsto m_{\{-1\}}\ot m_{\{0\}}$ denotes the map
$M\rarrow V\dual\ot_kM$ corresponding to the $\g$\+action map
$V\times M\rarrow M$ and $\pi_P$ denotes the $\g$\+contraaction map
$\Hom_k(V\dual,P)\rarrow P$.

 If $\psi\:\C\times\h\rarrow k$ is a continuous pairing between
a coalgebra $\C$ and a compact Lie algebra $\h$ such that the map
$\check\psi\:\C\rarrow\h\dual$ is a morphism of Lie coalgebras, then
for any left $\C$\contramodule{} $\P$ the induced contraaction of
$\h$ in $\P$ is defined as the composition $\Hom_k(\h\dual,\P)\rarrow
\Hom_k(\C,\P)\rarrow\P$ of the map induced by the map $\check\psi$
and the $\C$\+contraaction map.

\subsubsection{}  \label{tate-harish-chandra-contramodules}
 Let $\kap\:(\g',\C)\rarrow(\g,\C)$ be a central extension of Tate
Harish-Chandra pairs with the kernel~$k$.
 Then the category of left semicontramodules over the semialgebra
$\S^r_\kap(\g,\C)$ is isomorphic to the category of $k$\+vector spaces
$\bP$ endowed with a left $\C$\contramodule{} and a $\g'$\contramodule{}
structures such that the induced $\h$\contramodule{} structures
coincide, for any $\C$\invariant{} compact open subspace
$V\subset\g$ the $\g$\+contraaction map $\Hom_k(V\dual,\bP)\rarrow\bP$
is a morphism of $\C$\contramodule s, and the unit element of
$k\subset\g'$ acts by the identity in~$\bP$.
 Here the left $\C$\contramodule{} structure on the vector space
$\Hom_k(\M,\P)$ for a left $\C$\comodule{} $\M$ and a left
$\C$\contramodule{} $\P$ is defined by the formula
$\pi(g)(m)=\pi_\P(c\mapsto g(s(m_{(-1)})c)(m_{(0)}))$ for
$m\in\M$, \ $g\in\Hom_k(\C,\Hom_k(\M,\P))$.

 Indeed, according to~\ref{semi-mod-contra-described}, a left
$\S^r_\kap$\semicontramodule{} structure on $\bP$ is the same that
a left $\C$\contramodule{} and a left $U_\kap(\g)$\+module structures
such that induced $U(\h)$\module{} structures on $\bP$ coincide and
the (semicontra)action map $\bP\rarrow\Hom_{U(\h)}(U_\kap(\g),\bP)\simeq
\Cohom_\C(\S^r_\kap,\bP)$ is a morphism of $\C$\contramodule s.
 The latter condition is equivalent to the map $\bP\rarrow
\Hom_{U(\h)}(U(\h)\cdot V\;\bP)\simeq\Cohom_\C(\C\ot_{U(\h)}U(\h)\cdot
V\;\bP)$ being a morphism of $\C$\contramodule s for any compact
$\C$\invariant{} subspace $\h\oplus k\subset V\subset\g'$, where
$U(\h)\cdot V\subset U_\kap(\g)$.
 Given this data, one can use the short exact sequences
$\h\ot_k\bP\rarrow\h\ot\comp\bP\oplus V\ot_k\bP\rarrow V\ot\comp\bP$
to construct the Lie contraaction maps $V\ot\comp\bP\rarrow\bP$.
 Then the map $\bP\rarrow\Hom_{U(\h)}(U(\h)\cdot V\;\bP)$ is a morphism
of $\C$\contramodule s if and only if the map $\Hom_k(V\dual,\bP)
\rarrow \bP$ is a morphism of $\C$\contramodule s.
 To check this, one can express the first condition in terms of
the equality of two appropriate maps $\Hom_k(\C,\bP)\birarrow
\Hom_k(V,\bP)$ and the second condition in terms of the equality of
two maps $\Hom_k(V\dual\ot_k\C\;\bP)\birarrow\bP$.
 These two pairs of maps correspond to each other under a natural
isomorphism $V\dual\ot_k\C\simeq\C\ot_kV\dual$.
 In particular, our maps $\Hom_k(V\dual,\bP)\rarrow\bP$ are morphisms
of $\h$\contramodule s, and it follows that they define
a $\g'$\contramodule{} structure.

 We will call vector spaces $\bP$ endowed with such a structure 
\emph{Harish-Chandra contramodules over $(g,\C)$ with the central
charge~$\kap$}.
 If $\kap_1\:(\g',\C)\rarrow(\g,\C)$ and $\kap_2\:(\g'',\C)\rarrow
(\g,\C)$ are two central extensions of Tate Harish-Chandra pairs
with the kernels~$k$, and $\bM$ and $\bP$ are a Harish-Chandra
module and a Harish-Chandra contramodule over $(\g,\C)$ with
the central charges $\kap_1$ and $\kap_2$, respectively, then
the vector space $\Hom_k(\bM,\bP)$ has a natural structure of
Harish-Chandra contramodule with the central charge $\kap_2-\kap_1$.
 Here $\kap_2-\kap_1:(\g''',\C)\rarrow(\g,\C)$ denotes the Baer
difference of the central extensions $\kap_2$ and~$\kap_1$.
 This Harish-Chandra contramodule structure consists of the
$\g'''$\contramodule{} and $\C$\contramodule{} structures on
$\Hom_k(\bM,\bP)$ defined by the above rules.

\subsection{Isomorphism of semialgebras}
\label{tate-isomorphism-of-semialgebras}

\subsubsection{}   \label{semialgebra-isomorphism-characterized}
 For any two central extensions of Tate Harish-Chandra pairs
$\kap'\:(\g',\C)\rarrow(\g,\C)$ and $\kap''\:(\g'',\C)\rarrow(\g,\C)$
with the kernels identified with~$k$ we denote by $\kap'+\kap''$
their Baer sum, i.~e., the central extension of Tate Harish-Chandra
pairs $(\g''',\C)\rarrow(\g,\C)$ with $\g'''=\ker(\g'\oplus \g''
\to \g)/\im k$, where the map $\g'\oplus \g''\rarrow \g$ is
the difference of the maps $\g'\rarrow\g$ and $\g''\rarrow\g$, and
the map $k\rarrow \g'\oplus \g''$ is the difference of the maps
$k\rarrow\g'$ and $k\rarrow\g''$.
 The canonical central extension $(\g\til,\C)\rarrow(\g,\C)$ will be
denoted by~$\kap_0$.

 We claim that for any central extension of Tate Harish-Chandra pairs
$\kap\:(g',\C)\rarrow(\g,\C)$ with the kernel~$k$ such that 
the condition~(iv) of~\ref{tate-harish-chandra-central-extension} is
satisfied there is a natural isomorphism $\S^r_{\kap+\kap_0}(\g,\C)
\simeq\S^l_\kap(\g,\C)$ of semialgebras over the coalgebra~$\C$.
 This isomorphism is characterized by the following three properties.

\begin{itemize}
 \item[(a)] Consider the increasing filtration $F$ of the $k$\+algebra
$U_\kap(\g)$ with the components $F_iU_\kap(\g) =(k+\g'+\dsb+\g'{}^i)
U(\h) = U(\h)(k+\g'+\dsb+\g'{}^i)$ and the induced filtration
$F_i\S^l_\kap=F_iU_\kap(\g)\ot_{U(\h)}\C$ of the semialgebra
$\S^l_\kap=\S^l_\kap(\g,\C)$.
 Then we have $F_0\S^l_\kap\simeq\C$, \
$\S^l_\kap\simeq\ilim F_i\S^l_\kap$, and the semimultiplication maps
$F_i\S^l_\kap\oc_\C F_j\S^l_\kap\rarrow \S^l_\kap\oc_\C\S^l_\kap
\rarrow\S^l_\kap$ factorize through $F_{i+i}\S^l_\kap$.
 There is an analogous filtration $F_i\S^r_{\kap+\kap_0}=\C\ot_{U(\h)}
F_iU_{\kap+\kap_0}(\g)$ of the semialgebra $\S^r_{\kap+\kap_0}=
\S^r_{\kap+\kap_0}(\g,\C)$.
 The desired isomorphism $\S^r_{\kap+\kap_0}\simeq\S^l_\kap$ preserves
the filtrations~$F$.
 \item[(b)] The natural maps $F_{i-1}\S^l_\kap\rarrow F_i\S^l_\kap$ are
injective and their cokernels are coflat left and right
$\C$\comodule s, so the associated graded quotient semialgebra
$\gr_F\S^l_\kap=\bigoplus_i F_i\S^l_\kap/F_{i-1}\S^l_\kap$ is defined
(cf.~\ref{central-element-theorem}).
 The semialgebra $\gr_F\S^l_\kap$ is naturally isomorphic to the tensor
product $\Sym_k(\g/\h)\ot_k\C$ of the symmetric algebra $\Sym_k(\g/\h)$
of the $k$\+vector space $\g/\h$ and the coalgebra~$\C$, endowed
with the semialgebra structure corresponding to the left entwining
structure $\Sym_k(\g/\h)\ot_k\C\rarrow\C\ot_k\Sym_k(\g/\h)$ for
the coalgebra~$\C$ and the algebra $\Sym_k(\g/\h)$
(see~\ref{entwining-structures}).
 Here the entwining map is given by the formula $u\ot c\mpsto
cu_{(-1)}\ot u_{(0)}$, where $u\mpsto u_{(-1)}\ot u_{(0)}$ denotes
the $\C$\+coaction in $\Sym_k(\g/\h)$ induced by the $\C$\+coaction
in~$\g/\h$.
 Analogously, the semialgebra $\gr_F\S^r_{\kap+\kap_0}$ is naturally
isomorphic to the tensor product $\C\ot_k\Sym_k(\g/\h)$ endowed with
the semialgebra structure corresponding to the right entwining
structure $\C\ot_k\Sym_k(\g/\h)\rarrow\Sym_k(\g/\h)\ot_k\C$.
 Here the entwining map is given by the formula $c\ot u\mpsto
u_{(0)}\ot cu_{(1)}$, where the right coaction
$u\mpsto u_{(0)}\ot u_{(1)}$ is obtained from the above left coaction
$u\mpsto u_{(-1)}\ot u_{(0)}$ by applying the antipode.
 These left and right entwining maps are inverse to each other,
hence there is a natural isomorphism of semialgebras
$\gr_F\S^l_\kap\simeq\gr_F\S^r_{\kap+\kap_0}$.
 This isomorphism can be obtained by passing to the associated graded
quotient semialgebras in the desired isomorphism
$\S^l_\kap\simeq\S^r_{\kap+\kap_0}$.
 \item[(c)] Choose a section $b'\:\g/\h\rarrow\g'$ of the natural
surjection $\g'\rarrow\g'/(\h\oplus k)\simeq\g/\h$.
 Composing $b'$ with the surjection $\g'\rarrow\g$, we obtain
a section~$b$ of the natural surjection $\g\rarrow\g/\h$, hence
a direct sum decomposition $\g\simeq\h\oplus b(\g/\h)$.
 So there is the corresponding section $\g\rarrow\g\til$ of
the canonical central extension $\g\til\rarrow\g$; denote
by $\tilde b$ the composition $\g/\h\rarrow\g\rarrow\g\til$ of
the section~$b$ and the section $\g\rarrow\g\til$.
 The Baer sum of the sections $b'$ and $\tilde b$ provides a section
$b''\:\g/\h\rarrow \g''$, where $(\g'',\C)\rarrow(\g,\C)$ denotes
the central extension $\kap+\kap_0$.
 Now the composition $\g/\h\ot_k\C\rarrow \g'\ot_k\C\simeq F_1U_\kap(\g)
\ot_k\C\rarrow F_1\S^l_\kap$ of the map induced by the map~$b'$,
the isomorphism induced by the natural isomorphism
$\g'\simeq F_1U_\kap(\g)$, and the surjection $U_\kap(\g)\ot_k\C\rarrow
U_\kap(\g)\ot_{U(\h)}\C$ provides a section of the natural surjection
$F_1\S^l_\kap\rarrow F_1\S^l_\kap/F_0\S^l_\kap\simeq\g/\h\ot_k\C$.
 This section is a morphism of right $\C$\comodule s.
 Hence the corresponding retraction $F_1\S^l_\kap\rarrow F_0\S^l_\kap
\simeq\C$ is also a morphism of right $\C$\comodule s.
 Analogously, the composition $\C\ot_k\g/\h\rarrow \C\ot_k\g''\simeq
\C\ot_kF_1U_{\kap+\kap_0}(\g)\rarrow F_1\S^r_{\kap+\kap_0}$, where
the first morphism is induced by the map~$b''$, is a section of
the natural surjection $F_1\S^r_{\kap+\kap_0}\rarrow
F_1\S^r_{\kap+\kap_0}/F_0\S^r_{\kap+\kap_0}\simeq\C\ot_k\g/\h$;
this section is a morphism of left $\C$\comodule s.
 Hence so is the corresponding retraction $F_1\S^r_{\kap+\kap_0}\rarrow
F_0\S^r_{\kap+\kap_0}\simeq\C$.
 The desired isomorphism $F_1\S^l_\kap\simeq F_1\S^r_{\kap+\kap_0}$
identifies the compositions $F_1\S^l_\kap\rarrow\C\rarrow k$ and
$F_1\S^r_{\kap+\kap_0}\rarrow\C\rarrow k$ of the retractions
$F_1\S^l_\kap\rarrow \C$ and $F_1\S^r_{\kap+\kap_0}\rarrow\C$ with
the counit map $\C\rarrow k$.
 This condition holds for all sections~$b'$.
\end{itemize}

\begin{thm}
 There exists a unique isomorphism of semialgebras\/
$\S^r_{\kap+\kap_0}(\g,\C)\simeq\S^l_\kap(\g,\C)$ over~$\C$
satisfying the above properties~(a-c).
\end{thm}

\begin{proof}
 Uniqueness is clear, since a morphism from a $\C$\+$\C$\bicomodule{}
to the bicomodule $\C$ is determined by its composition with
the counit map $\C\rarrow k$.
 The proof of existence occupies
subsections~\ref{co-pbw-decreasing-filtration}--%
\ref{tate-change-of-section}.

\medskip
 The next result is obtained by specializing
the semimodule-semicontramodule correspondence theorem to the case
of Harish-Chandra modules and contramodules.

\begin{cor}
 There is a natural equivalence $\boR\Psi_{\S^l_\kap} =
\boL\Phi_{\S^r_{\kap+\kap_0}}^{-1}$ between the semiderived category
of Harish-Chandra modules with the central charge~$\kap$ over $(\g,\C)$
and the semiderived category of Harish-Chandra contramodules with
the central charge~$\kap+\kap_0$ over~$(\g,\C)$.
 Here the semiderived category of Harish-Chandra modules is defined
as the quotient category of the homotopy category of complexes of
Harish-Chandra modules by the thick subcategory of\/ $\C$\coacyclic{}
complexes; the semiderived category of Harish-Chandra contramodules
is analogously defined as the quotient category by the thick
subcategory of\/ $\C$\contraacyclic{} complexes.
\end{cor}

\begin{proof}
 This follows from the results of~\ref{tate-harish-chandra-modules}
and \ref{tate-harish-chandra-contramodules}, the above Theorem, and
Corollary~\ref{semimodule-semicontramodule-subsect}.
\end{proof}

\begin{rmk}
 The main property of the equivalence of semiderived categories
provided by the above Corollary is that it transforms
the Harish-Chandra modules $\bM$ that, considered as $\C$\comodule s,
are the cofree comodules $\C\ot_k E$ cogenerated by a vector space
$E$, into the Harish-Chandra contramodules $\bP$ that, considered as
$\C$\contramodule s, are the free contramodules $\Hom_k(\C,E)$
generated by the same vector space $E$, and vice versa.
 The similar assertion holds for any complexes of $\C$\+cofree
Harish-Chandra modules and $\C$\+free Harish-Chandra contramodules.
 The above Corollary is a way to formulate the classical duality
between Harish-Chandra modules with the complementary central
charges $\kap$ and $-\kap-\kap_0$ \cite{FeFu,RCW}.
 Of course, there is no hope of establishing an \emph{anti-equivalence}
between any kinds of exotic derived categories of \emph{arbitrary}
Harish-Chandra modules over $(\g,\C)$ with the complementary central
charges, as the derived category of vector spaces is not
anti-equivalent to itself.
 At the very least, one would have to impose some finiteness conditions
on the Harish-Chandra modules.
 The introduction of contramodules allows to resolve this problem.
 Still one can use the functor $\Phi_\S$ to construct
a \emph{contravariant functor} between the semiderived categories of
Harish-Chandra modules with the complementary central charges.
 Choose a vector space $U$; for example, $U=k$.
 Consider the functor $\bN\mpsto \Hom_k(\bN,U)$ acting from
the semiderived category of Harish-Chandra semimodules over $(\g,\C)$
with the central charge $-\kap-\kap_0$ to the semiderived category
of Harish-Chandra semicontramodules over $(\g,\C)$ with the central
charge $\kap+\kap_0$.
 Composing this functor $\Hom_k({-},U)$ with the functor
$\boL\Phi_{\S^r_{\kap+\kap_0}}$, one obtains a contravariant functor
$\sD^\si(\simodr\S^r_{\kap+\kap_0})\rarrow\sD^\si(\S^l_\kap\simodl)$.
 The latter functor transforms the Harish-Chandra modules that
as $\C$\comodule s are cofreely cogenerated by a vector space $E$
into the Harish-Chandra modules that as $\C$\comodule s are 
cofreely cogenerated by the vector space $\Hom_k(E,U)$, and
similarly for complexes of $\C$\+cofree Harish-Chandra modules.
 One cannot avoid using the exotic derived categories in this
construction, because the functor $\boL\Phi_\S$ does not preserve
acyclicity, in general (see~\ref{prelim-co-contra-acyclicity}).
\end{rmk}

\subsubsection{}  \label{co-pbw-decreasing-filtration}
 The semialgebras $\S^l_\kap$ and $\S^r_{\kap+\kap_0}$ endowed with
the increasing filtrations~$F$ are left and right coflat nonhomogeneous
Koszul semialgebras over the coalgebra~$\C$
(see~\ref{central-element-theorem}).
 Indeed, there are natural isomorphisms of complexes of
$\C$\+$\C$\bicomodule s $\Br^\bu_\gr(\gr_F\S^l_\kap,\C)\simeq
\Br^\bu_\gr(\Sym_k(\g/\h),k)\ot_k\C$ and
$\Br^\bu_\gr(\gr_F\S^r_{\kap+\kap_0},\C)\simeq
\C\ot_k\Br^\bu_\gr(\Sym_k(\g/\h,k)$, and
the $k$\+algebra $\Sym_k(\g/\h)$ is Koszul.
 Here the left $\C$\+coaction in $\Br^\bu_\gr(\Sym_k(\g/\h),k)\ot_k\C$
is the tensor product of the $\C$\+coaction in
$\Br^\bu_\gr(\Sym_k(\g/\h),k)$ induced by the $\C$\+coaction in $\g/\h$
and the left $\C$\+coaction in $\C$, while the right $\C$\+coaction
in $\Br^\bu_\gr(\Sym_k(\g/\h),k)\ot_k\C$ is induced by the right
$\C$\+coaction in~$\C$. 
 The $\C$\+$\C$\bicomodule{} structure on
$\C\ot_k\Br^\bu_\gr(\Sym_k(\g/\h,k)$ is defined in the analogous way
(with the left and right sides switched).

 The left and right coflat Koszul coalgebras $\D^l$ and $\D^r$ over~$\C$
quadratic dual to the left and right coflat Koszul semialgebras
$\gr_F\S^l_\kap$ and $\gr_F\S^r_{\kap+\kap_0}$ are described as follows.
 One has $\D^l\simeq\bigwedge_k(\g/\h)\ot_k\C$, where
$\bigwedge_k(\g/\h)$ denotes the exterior coalgebra of the $k$\+vector
space $\g/\h$, i.~e., the coalgebra quadratic dual to the symmetric
algebra $\Sym_k(\g/\h)$.
 The counit of $\bigwedge_k(\g/\h)\ot_k\C$ is the tensor product of
the counits of $\bigwedge_k(\g/\h)$ and $\C$, while
the comultiplication in $\bigwedge_k(\g/\h)\ot_k\C$ is constructed as
the composition $\bigwedge_k(\g/\h)\ot_k\C\rarrow\bigwedge_k(\g/\h)
\ot_k\bigwedge_k(\g/\h)\ot_k\C\ot_k\C\rarrow\bigwedge_k(\g/\h)\ot_k\C
\ot_k\bigwedge_k(\g/\h)\ot_k\C$ of the map induced by
the comultiplications in $\bigwedge_k(\g/\h)$ and $\C$ and the map
induced by the ``permutation'' map $\bigwedge_k(\g/\h)\ot_k\C\rarrow
\C\ot_k\bigwedge_k(\g/\h)$.
 The latter map is given by the formula $u\ot c\mpsto
cu_{[-1]}\ot u_{[0]}$ for $u\in\bigwedge_k(\g/\h)$ and $c\in\C$,
where $u\mpsto u_{[-1]}\ot u_{[0]}$ denotes the $\C$\+coaction in
$\bigwedge_k(\g/\h)$ induced by the $\C$\+coaction in~$\g/\h$.

 Analogously, one has $\D^r\simeq\C\ot_k\bigwedge_k(\g/\h)$, where
the counit of $\C\ot_k\bigwedge_k(\g/\h)$ is the tensor product of
the counits of $\bigwedge_k(\g/\h)$ and $\C$, while
the comultiplication in $\C\ot_k\bigwedge_k(\g/\h)$ is defined in
terms of the ``permutation'' map $\C\ot_k\bigwedge_k(\g/\h)\rarrow
\bigwedge_k(\g/\h)\ot_k\C$.
 The latter map is given by the formula $c\ot u\mpsto u_{[0]}\ot
cu_{[1]}$, where the right coaction $u\mpsto u_{[0]}\ot u_{[1]}$
is obtained from the left coaction $u\mpsto u_{[-1]}\ot u_{[0]}$
by applying the antipode.
 Both coalgebras $\bigwedge_k(\g/\h)\ot_k\C$ and $\C\ot_k
\bigwedge_k(\g/\h)$ have gradings induced by the grading of
$\bigwedge_k(\g/\h)$.
 The two ``permutation'' maps are inverse to each other, and they
provide an isomorphism of graded coalgebras $\D^l\simeq\D^r$.

 Now recall that we have assumed the condition~(iv)
of~\ref{tate-harish-chandra-central-extension}.
 Denote by $\dsb\subset V^2\C\subset V^1\C\subset V^0\C=\C$
the decreasing filtration of $\C$ orthogonal to the natural increasing
filtration of the universal enveloping algebra $U(\h)$, that is
$V^i\C\subset\C$ consists of all $c\in\C$ such that $\phi(c,x)=0$ for
all $x\in k+\h+\dsb+\h^{i-1}\subset U(\h)$.
 Notice that the decreasing filtration $V$ is compatible with both
the coalgebra and algebra structures on~$\C$; in particular,
it is a filtration by ideals with respect to the multiplication.
 The subspace $V^1\C$ is the kernel of the counit map $\C\rarrow k$;
the subspace $V^2\C$ is the kernel of the map $\C\rarrow\h\dual\oplus k$
which is the sum of the map $\check\psi$ and the counit map.

 Define decreasing filtrations $V$ on the coalgebras $\D^l$ and $\D^r$
by the formulas $V^i\D^l\simeq \bigwedge_k(\g/\h)\ot_k V^i\C$ and
$V^i\D^r\simeq V^i\C\ot_k\bigwedge_k(\g/\h)$; these filtrations are
compatible with the coalgebra structures on $\D^l$ and $\D^r$, and
correspond to each other under the isomorphism $\D^l\simeq\D^r$.
 Set $\D^l=\D\simeq\D^r$.
 The coalgebra $\D$ is cogenerated by the maps $\D\rarrow\D_0/V^2\D_0$
and $\D\rarrow\D_1/V^1\D_1$, i.~e., the iterated comultiplication map
from $\D$ to the direct product of all tensor powers
$\D_0/V^2\D_0\oplus\D_1/V^1\D_1$ is injective.
 Moreover, the decreasing filtration $V$ on~$\D$ is cogenerated by
the filtrations on $\D_0/V^2\D_0$ and $\D_1/V^1\D_1$, i.~e.,
the subspaces $V^i\D$ are the full preimages of the subspaces of
the induced filtration on the product of all tensor powers of
$\D_0/V^2\D_0\oplus\D_1/V^1\D_1$ under the iterated
comultiplication map.

\subsubsection{}
 Composing the equivalences of categories
from~\ref{quasi-differential-corings} and Theorem~\ref{pbw-theorem},
we obtain an equivalence between the category of left (right) coflat
nonhomogeneous Koszul semialgebras over~$\C$ and the category of
left (right) coflat Koszul CDG\+coalgebras over~$\C$.
 Here a CDG\+coagebra $(\D,d,h)$ is called Koszul over~$\C$ if
the underlying graded coalgebra~$\D$ is Koszul over~$\C$.
 Recall that for a left (right) coflat nonhomogeneous Koszul
semialgebra $\tS$ and the corresponding quasi-differential coalgebra
$\tD$ one has $F_1\tS\simeq\D_1\til$, so to construct a specific
CDG\+coalgebra $(\D,d,h)$ corresponding to a given filtered semialgebra
$\tS$ one has to choose a linear map $\delta\:F_1\tS\rarrow k$
such that the composition of the injection $\C\simeq F_0\tS\rarrow
F_1\tS$ with~$\delta$ coincides with the counit map of~$\C$.

 Choose a section $b'\:\g/\h\rarrow\g'$ and construct the related
section~$b''\:\g/\h\rarrow\g''$; denote by $\delta^l_{b'}\:F_1\S^l_\kap
\rarrow k$ and $\delta^r_{b''}\:F_1\S^r_{\kap+\kap_0}\rarrow k$
the corresponding linear functions constructed in~(c)
of~\ref{semialgebra-isomorphism-characterized}.
 In order to obtain an isomorphism of left and right coflat
nonhomogeneous Koszul semialgebras $\S^l_\kap\simeq\S^r_{\kap+\kap_0}$,
we will construct an isomorphism between the CDG\+coalgebras
$(\D^l,d_{b'}^l,h_{b'}^l)$ and $(\D^r,d_{b''}^r,h_{b''}^r)$
corresponding to the filtered semialgebras $\S^l_\kap$ and
$\S^r_{\kap+\kap_0}$ endowed with the linear functions
$\delta^l_{b'}$ and $\delta^r_{b''}$.
 The isomorphism of coalgebras $\D^l\simeq\D^r$ is already defined;
all we have to do is to check that it identifies $d_{b'}^l$ with
$d_{b''}^r$ and $h_{b'}^l$ with $h_{b''}^r$.
 Besides, we need to show that the isomorphism $\S^l_\kap\simeq
\S^r_{\kap+\kap_0}$ so obtained does not depend on the choice of~$b'$.
 Here it suffices to check that changing the section~$b'$ to $b'_1$
leads to isomorphisms of CDG\+coalgebras $(\id,a^l)\:
(\D^l,d_{b'}^r,h_{b'}^r)\rarrow(\D^l,d_{b'_1}^r,h_{b'_1}^r)$ and
$(\id,a^r)\:(\D^r,d_{b''}^r,h_{b''}^r)\rarrow
(\D^r,d_{b''_1}^r,h_{b''_1}^r)$ with the linear functions
$a^l$ and $a^r$ being identified by the isomorphism $\D^l\simeq\D^r$.

 Since the coalgebra $\D^l=\D\simeq\D^r$ is cogenerated by the maps
$\D\rarrow \D_0/V^2\D_0$ and $\D\rarrow\D_1/V_1\D_1$, it suffices to
check that the compositions of $d_{b'}^l$ and $d_{b''}^r$ with
these two maps coincide in order to show that $d_{b'}^l=d_{b''}^r$.
 We will also see that these compositions factorize through
$\D_1/V^2\D_1$ and $\D_2/V^1\D_2$, respectively, and the induced map
$\D_1/V^2\D_1\rarrow\D_0/V^2\D_0$ preserves the images of~$V^1$
(actually, even maps the whole of $\D_1/V^2\D_1$ into
$V^1\D_0/V^2\D_1$), hence it will follow that the differential
$d_{b'}^l=d_{b''}^r$ preserves the decreasing filtration~$V$.
 Besides, we will see that the linear function $h_{b'}^l=h_{b''}^r$
annihilates the subspace $V^2\D_2$ and the linear function $a^l=a^r$
corresponding to a change of section $b'$ annihilates the subspace
$V^2\D_1$.

\subsubsection{}   \label{commutator-components-notation}
 Let us introduce notation for the components of the commutator map
with respect to the direct sum decomposition $\g'\simeq \h\oplus
b'(\g/\h)\oplus k$.
 As above, the Lie coalgebra structure on $\h\dual$ is denoted by
$x^*\mpsto x^*_{\{1\}}\wedge x^*_{\{2\}}$.
 Denote the Lie coaction of $\h\dual$ in $\g/\h$, i.~e., the map
$\g/\h\rarrow\h\dual\ot_k\g/\h$ corresponding to the commutator map
$\h\times \g/\h\rarrow\g/\h$, by $u\mpsto u_{\{-1\}}\ot u_{\{0\}}$.
 These two maps do not depend on the choice of the section~$b'$;
the rest of them do.

 Denote by $u\ot x^*\mpsto u(x^*)$ the map $\g/\h\ot_k\h\dual\rarrow
\h\dual$ corresponding to the projection of the commutator map
$b'(\g/\h)\times\h\rarrow \g'\rarrow\h$.
 Denote by $u\wedge v\mpsto \{u,v\}$ the map $\bigwedge^2_k(\g/\h)
\rarrow\g/\h$ corresponding to the commutator map
$\bigwedge^2_kb'(\g/\h) \rarrow \g/\h$.
 Denote by $u\wedge v\ot x^*\mpsto (u,v)_{x^*}$ the map $\bigwedge^2_k
(\g/\h)\ot_k\h\dual\rarrow k$ corresponding to the projection of
the commutator map $\bigwedge^2_kb'(\g/\h)\rarrow\g'\rarrow\h$.

 The above five maps only depend on the Lie algebra $\g$ with
the subalgebra $\h$ and the section $b\:\g/\h\rarrow \g$, but
the following two will depend essentially on~$\g'$ and~$b'$.
 Denote by $\rho'\:\g/\h\rarrow\h\dual$ the map corresponding to
the projection of the commutator map $b'(\g/\h)\times \h\rarrow
\g'\rarrow k$.
 Denote by $\sigma'\:\bigwedge^2_k(\g/\h)\rarrow k$ the map
corresponding to the projection of the commutator map
$\bigwedge^2_kb'(\g/\h)\rarrow \g'\rarrow k$.
 Denote by $\tilde\rho$, $\tilde\sigma$ and $\rho''$, $\sigma''$
the analogous maps corresponding to the central extensions
$\g\til\rarrow \g$ and $\g''\rarrow\g$ with the sections $\tilde b$
and $b''$.
 Clearly, we have $\rho''=\rho'+\tilde\rho$ and $\sigma''=
\sigma'+\tilde\sigma$.

 Set $\b=b(\g/\h)\subset\g$.
 The composition $\theta$ of the commutator map in $\gl(\g)\til$ with
the projection $\gl(\g)\til\rarrow k$ corresponding to the section
$\gl(\g)\rarrow\gl(\g)\til$ coming from the direct sum decomposition
$\g\simeq\h\oplus\b$ is written explicitly as follows.
 For any continuous linear operator $A\:\g\rarrow\g$ denote
by $A_{\h\to\b}\:\h\rarrow\b$, \ $A_{\b\to\h}\:\b\rarrow\h$, etc., its
components with respect to our direct sum decomposition.
 Then the cocycle $\theta$ is given by the formula $\theta(A\wedge B)
=\tr(A_{\b\to\h}B_{\h\to\b})-\tr(B_{\b\to\h}A_{\h\to\b})$, where $\tr$
denotes the trace of a linear operator $\h\rarrow\h$ with an open
kernel.

 Using this formula, one can find that, in the above notation, 
$\tilde\rho(u)=-u_{\{0\}}(u_{\{-1\}})$ and $\tilde\sigma(u\wedge v)=
-(u,v_{\{0\}})_{v_{\{-1\}}}+(v,u_{\{0\}})_{u_{\{-1\}}}$.

\subsubsection{}
 We have $\D^l_0\simeq \C$, \ $\D^l_1\simeq\g/\h\ot_k\C$, and
$\D^l_2\simeq \bigwedge^2_k(\g/\h)\ot_k\C$.
 The composition of the map $d^l_{b'}\:\g/\h\ot_k\C\rarrow\C$ with
the counit map $\eps\:\C\rarrow k$ vanishes, since $d^l_{b'}$ is
a coderivation.
 Let us start with computing the composition of the map $d^l_{b'}$
with the map $\check\psi\:\C\rarrow\h\dual$.

 The class of an element $u\ot c\in \g/\h\ot_k\C$ can be represented
by the element $b'(u)\ot_{U(\h)} c\in F_1U_\kap(\g)\ot_{U(\h)}\C\simeq
\D^l_1{}\til$ in the quasi-differential coalgebra $\D^l{}\til$
corresponding to the filtered semialgebra $\S^l_\kap$.
 Denote the image of $b'(u)\ot_{U(\h)}c$ under the comultiplication map
$\D^l_1{}\til\rarrow\C\ot_k\D^l_1{}\til$ by $c_1\ot (z\ot_{U(\h)} c_2)$,
where $z\in F_1U_\kap(\g)$.
 The total comultiplication of $b'(u)\ot_{U(\h)}c$ is then equal to
$c_1\ot (z\ot_{U(\h)} c_2)+ (b'(u)\ot_{U(\h)}c_{(1)})\ot c_{(2)}$.
 We have $d^l_{b'}(u\ot c)=\delta_{b'}^l(b'(u)\ot_{U(\h)}c_{(1)})c_{(2)}
- \delta_{b'}^l(z\ot_{U(\h)}c_2)c_1 = 
- \delta_{b'}^l(z\ot_{U(\h)}c_2)c_1$.
 Furthermore, $\psi(d^l_{b'}(u\ot c),x) = -\psi(c_1,x)
\delta_{b'}^l(z\ot_{U(\h)}c_2)=-\delta_{b'}^l(xb'(u)\ot_{U(\h)}c)=
-\delta_{b'}^l([x,b'(u)]\ot_{U(\h)}c)-\delta_{b'}^l(b'(u)\ot_{U(\h)}xc)
=\delta_{b'}^l([b'(u),x]\ot_{U(\h)}c)=\lan x, u(\check\psi(c))\ran
+\lan x,\rho'(u)\ran \eps(c)$ for $x\in\h$, since
$\psi(x,c_1)z\ot_{U(\h)}c_2=xb'(u)\ot_{U(\h)}c$.

 So the composition of the map  $d^l_{b'}\:\g/\h\ot_k\C\rarrow\C$ with
the map $\check\psi\:\C\rarrow\h\dual$ is equal to the composition of
the map $\id\ot(\check\psi,\eps)\:\g/\h\ot_k\C\rarrow\g/\h\ot_k
(\h\dual\oplus k)$ with the map $\g/\h\ot_k(\h\dual\oplus k)\rarrow
\h\dual$ given by the formula $u\ot x^*+v\mpsto u(x^*)+\rho'(v)$.

 Now let us compute the composition of the map $d_{b'}^l\:\bigwedge^2_k
(\g/\h)\ot_k\C\rarrow\g/\h\ot_k\C$ with the map $\id\ot\eps\:\g/\h
\ot_k\C\rarrow\g/\h$.
 The vector space $\D_2^l{}\til$ is the kernel of the semimultiplication
map $F_1\S^l_\kap\oc_\C F_1\S^l_\kap\rarrow F_2\S^l_\kap$, which can be
identified with the kernel of the map $F_1U_\kap(\g)\ot_{U(\h)}
F_1U_\kap(\g)\ot_{U(\h)}\C\rarrow F_2U_\kap(\g)\ot_{U(\h)}\C$ induced by
the multiplication map $F_1U_\kap(\g)\ot_{U(\h)}F_1U_\kap(\g)\rarrow
F_2U_\kap(\g)$.
 The class of an element $u\wedge v\ot c\in\bigwedge^2_k(\g/\h)\ot_k\C$
can be represented by the element $b'(u)\ot_{U(\h)}b'(v)\ot_{U(\h)}c -
b'(v)\ot_{U(\h)}b'(u)\ot_{U(\h)}c - [b'(u),b'(v)]\ot_{U(\h)} 1
\ot_{U(\h)} c$ in the latter kernel.

 Denote the image of $b'(v)\ot_{U(\h)}c$ under the comultiplication map
$\D^l_1{}\til\rarrow\C\ot_k\D^l_1{}\til$ by $c_1\ot (z\ot_{U(\h)} c_2)$,
where $z\in F_1U_\kap(\g)$; then the image of $b'(u)\ot_{U(\h)}b'(v)
\ot_{U(\h)}c$ under the comultiplication map $\D^l_2{}\til\rarrow
\D^l_1{}\til\ot_k\D^l_1{}\til$ is equal to $(b'(u)\ot_{U(\h)}c_1)\ot
(z\ot_{U(\h)}c_2)$.
 The image of $[b'(u),b'(v)]\ot_{U(\h)}1\ot_{U(\h)}c$ under
the same map $\D^l_2{}\til\rarrow\D^l_1{}\til\ot_k
\D^l_1{}\til$ is equal to $([b'(u),b'(v)]\ot_{U(\h)}c_{(1)})\ot
(1\ot_{U(\h)}c_{(2)})$.
 We have $\delta_{b'}^l(b'(u)\ot_{U(\h)}c_1)\overline{z\ot_{U(\h)}c_2}
=0$ and $\delta_{b'}^l(z\ot_{U(\h)}c_2)(\id\ot\eps)
\overline{b'(u)\ot_{U(\h)}c_1}=\eps(c_1)\delta_{b'}^l(z\ot_{U(\h)}c_2)u
=\delta_{b'}^l(b'(v)\ot_{U(\h)}c)u=0$, where $\overline p\in\D_1^l$
denotes the image of an element $p\in\D_1^l{}\til$.
 Furthermore, $\delta_{b'}^l([b'(u),b'(v)]\ot_{U(\h)}c_{(1)})\overline
{1\ot_{U(\h)}c_{(2)}}=0$.
 Hence $(\id\ot\eps)d_{b'}^l(u\wedge v\ot c) = -\delta_{b'}^l
(1\ot_{U(\h)}c_{(2)})(\id\ot\eps)\overline{[b'(u),b'(v)]\ot_{U(\h)}
c_{(1)}} = -(\id\ot\eps)\overline{[b'(u),b'(v)]\ot_{U(\h)}c} =
-\{u,v\}\eps(c)$.

 So the composition of the map $d_{b'}^l\:\bigwedge^2_k(\g/\h)\ot_k\C
\rarrow\g/\h\ot_k\C$ with the map $\id\ot\eps\:\g/\h\ot_k\C\rarrow\g/\h$
is equal to the composition of the map $\id\ot\eps\:\bigwedge^2_k(\g/\h)
\ot_k\C\rarrow\bigwedge^2_k(\g/\h)$ with the map $\bigwedge^2_k(\g/\h)
\rarrow\g/\h$ given by the formula $u\wedge v\mpsto -\{u,v\}$.

 Finally, let us compute the linear function $h_{b'}^l\:
\bigwedge^2_k(\g/\h)\ot_k\C\rarrow k$.
 We have $\delta_{b'}^l(b'(u)\ot_{U(\h)}c_1)\delta_{b'}^l
(z\ot_{U(\h)}c_2)=0$, hence $h(u\wedge v\ot c) = -\delta_{b'}^l
([b'(u),b'(v)]\ot_{U(\h)}c_{(1)})\delta_{b'}^l(1\ot_{U(\h)}c_{(2)}) =
-\delta_{b'}^l([b'(u),b'(v)]\ot_{U(\h)}c) = -(u,v)_{\check\psi(c)} 
-\sigma'(u\wedge v)\eps(c)$.
 So the linear function $h_{b'}^l$ is equal to the composition of
the map $\id\ot (\check\psi,\eps)\:\bigwedge^2_k(\g/\h)\ot_k\C\rarrow
\bigwedge^2_k(\g/\h)\ot_k(\h\dual\oplus k)$ and the linear function
$\bigwedge^2_k(\g/\h)\ot_k(\h\dual\oplus k)\rarrow k$ given by
the formula $u_1\wedge v_1\ot x^* + u\wedge v\mpsto
-(u_1,v_1)_{x^*}-\sigma'(u\wedge v)$.

\subsubsection{}
 Analogously, we have $\D^r_0\simeq\C$, \ $\D^r_1\simeq\C\ot_k\g/\h$,
and $\D^r_2\simeq\C\ot_k\bigwedge^2_k(\g/\h)$.
 The composition of the map $d_{b''}^r\:\C\ot_k\g/\h\rarrow\C$ with
the map $\eps\:\C\rarrow k$ vanishes.
 The composition of the map $d_{b''}^r\:\C\ot_k\g/\h\rarrow\C$ with
the map $\check\psi\:\C\rarrow\h\dual$ is equal to the composition
of the map $(\check\psi,\eps)\ot\id\:\C\ot_k\g/\h\rarrow(\h\dual
\oplus k)\ot_k\g/\h$ and the map $(\h\dual\oplus k)\ot_k\g/\h\rarrow
\h\dual$ given by the formula $x^*\ot u+v\mpsto u(x^*)+\rho''(v)$.
 The composition of the map $d_{b''}^r\:\C\ot_k\bigwedge^2_k(\g/\h)
\rarrow\C\ot_k\g/\h$ with the map $\eps\ot\id\:\C\ot_k\g/\h\rarrow
\g/\h$ is equal to the composition of the map $\eps\ot\id\:
\C\ot_k\bigwedge^2_k(\g/\h)\rarrow\bigwedge^2_k(\g/\h)$ and the map
$\bigwedge^2_k(\g/\h)\rarrow\g/\h$ given by the formula $u\wedge v
\mpsto-\{u,v\}$.
 The linear function $h_{b''}^r\:\C\ot_k\bigwedge^2_k(\g/\h)\rarrow k$
is equal to the composition of the map $(\check\psi,\eps)\ot\id\:
\C\ot_k\bigwedge^2_k(\g/\h)\rarrow (\h\dual\oplus k)\ot_k
\bigwedge^2_k(\g/\h)$ and the linear function $(\h\dual\oplus k)\ot_k
\bigwedge^2_k(\g/\h)\rarrow k$ given by the formula $x^*\ot u_1\wedge
v_1 + u\wedge v\mpsto -(u_1,v_1)_{x^*} - \sigma''(u\wedge v)$.

 The isomorphism $\g/\h\ot_k\C\simeq\C\ot_k\g/\h$ forms a commutative
diagram with the map $\g/\h\ot_k\C\rarrow\g/\h\ot_k(\h\dual\oplus k)$,
the map $\C\ot_k\g/\h\rarrow(\h\dual\oplus k)\ot_k\g/\h$, and
the isomorphism $\g/\h\ot_k(\h\dual\oplus k)\simeq (\h\dual\oplus k)
\ot_k\g/\h$ given by the formula $u\ot x^* + v \mpsto x^*\ot u +
v_{\{-1\}}\ot v_{\{0\}} + v$.
 Analogously, the isomorphism $\bigwedge^2_k(\g/\h)\ot_k\C\simeq
\C\ot_k\bigwedge^2_k(\g/\h)$ forms a commutative diagram with the map
$\bigwedge^2_k(\g/\h)\ot_k\C\rarrow\bigwedge^2_k(\g/\h)\ot_k
(\h\dual\oplus k)$, the map $\C\ot_k\bigwedge^2_k(\g/\h)\rarrow
(\h\dual\oplus k)\ot_k\bigwedge^2_k(\g/\h)$, and the isomorphism
$\bigwedge^2_k(\g/\h)\ot_k(\h\dual\oplus k)\simeq (\h\dual\oplus k)
\ot_k\bigwedge^2_k(\g/\h)$ given by the formula $u_1\wedge v_1\ot x^* +
u\wedge v\mpsto x^*\ot u_1\wedge v_1 + u_{\{-1\}}\ot u_{\{0\}}\wedge v
+ v_{\{-1\}}\ot u\wedge v_{\{0\}} + u\wedge v$.

 Now it is straightforward to check that the isomorphism
$\D^l\simeq\D^r$ identifies $d^l_{b'}$ with $d^r_{b''}$ modulo
$V^2\D_0\oplus V^1\D_1 \oplus\D_2\oplus\D_3\oplus\dsb$ and
$h^l_{b'}$ with $h^r_{b''}$.
 Indeed, one has $u(x^*)+v_{\{0\}}(v_{\{-1\}})+\rho''(v) =
u(x^*)+\rho'(v)$ and $-(u_1,v_1)_{x^*} - (u_{\{0\}},v)_{u_{\{-1\}}} -
(u,v_{\{0\}})_{v_{\{-1\}}} - \sigma''(u\wedge v) = -(u_1,v_1)_{x^*} -
\sigma'(u\wedge v)$.

\subsubsection{}    \label{tate-change-of-section}
 Finally, let $b_1'\:\g/\h\rarrow\g'$ be another section of
the surjection $\g'\rarrow\g'/(\h\oplus k)\simeq\g/\h$.
 Then we can write $b_1'=b+t+t'$ with $t\:\g/\h\rarrow\h$ and
$t'\:\g/\h\rarrow k$.
 Analogously, the sections $\tilde b_1\:\g/\h\rarrow \g\til$ and
$b_1''\:\g/\h\rarrow \g''$ corresponding to $b_1'$ have the form
$\tilde b_1 = \tilde b+t+\tilde t$ and $b_1''=b''+t+t''$ with
$t''=t'+\tilde t$.

 Denote by $\tau\;\tau_1\:\gl(\g)\rarrow\gl(\g)\til$ the sections
corresponding to direct sum decompositions $\g=\h\oplus b(\g/\h)$
and $\g=\h\oplus b_1(\g/\h)$ with $b_1=b+t$, \ $t\:\g/\h\rarrow\h$.
 Then one has $\tau_1(A)-\tau(A)=\tr(tA_{\h\to\g/\h})$ for any
$A\in\gl(\g)$, where $A_{\h\to\g/\h}$ denotes the composition
$\h\rarrow\g\rarrow\g\rarrow\g/\h$ of the endomorphism $A$ with
the injection $\h\rarrow\g$ and the surjection $\g\rarrow\g/\h$.

 Using this formula, one can find that $\tilde t(u)=
-\lan u_{\{-1\}},t(u_{\{0\}})\ran$, where $\lan\,\.,\,\ran$
denotes the natural pairing $\h\dual\times\h\rarrow k$.

 The natural isomorphism $(\id,a^l)\:(\D^l,d^l_{b'},h^l_{b'})\rarrow
(\D^l,d^l_{b'_1},h^l_{b'_1})$ between the CDG\+coalgebras corresponding
to the sections $b'$ and $b'_1$ can be computed easily; the linear
function $a^l\:\D^l_1\rarrow k$ is the composition of the map
$\id\ot(\check\psi,\eps)\:\g/\h\ot_k\C\rarrow\g/\h\ot_k
(\h\dual\oplus k)$ and the linear function $\g/\h\ot_k(\h\dual\oplus k)
\rarrow k$ given by the formula $u\ot x^* + v \mpsto -\lan x^*,t(u)\ran
- t'(v)$.
 Analogously, the linear function $a^r$ in the natural isomorphism
$(\id,a^r)\:(\D^r,d^r_{b''},h^r_{b''})\rarrow
(\D^r,d^r_{b''_1},h^r_{b''_1})$ between the CDG\+coalgebras
corresponding to the sections $b''$ and $b''_1$ is the composition of
the map $(\check\psi,\eps)\ot\id\:\C\ot_k\g/\h\rarrow(\h\dual\oplus k)
\ot_k\g/\h$ and the linear function $(\h\dual\oplus k)\ot_k\g/\h
\rarrow k$ given by the formula $x^*\ot u + v \mpsto -\lan x^*,t(u)\ran
- t''(v)$.

 Now it is straightforward to check that the isomorphism
$\D^l\simeq\D^r$ identifies $a^l$ with~$a^r$.
 Indeed, $-\lan x^*,t(u)\ran - \lan v_{\{-1\}},t(v_{\{0\}})\ran -
t''(v) = -\lan x^*,t(u)\ran - t'(v)$.

\medskip
 Theorem~\ref{semialgebra-isomorphism-characterized} is proven.
\end{proof}

\subsection{Semiinvariants and semicontrainvariants}

\subsubsection{}
 Let $\g$ be a Tate Lie algebra, $\g\til\rarrow\g$ be the canonical
central extension, and $\h\subset\g$ be a compact open subalgebra;
recall that the central extension $\g\til\rarrow\g$ splits
canonically over~$\h$.
 Let $N$ be a discrete $\g\til$\module{} where the unit element of
$k\subset\g\til$ acts by minus the identity.
 We would like to construct a natural map $(\g/\h\ot_kN)^\h\rarrow
N^\h$, where the superindex $\h$ denotes the $\h$\invariant s and
the action of $\h$ in $N$ is defined in terms of the canonical
splitting $\h\rarrow\g\til$.

 Choose a section $b\:\g/\h\rarrow\g$ of the surjection $\g\rarrow
\g/\h$.
 The direct sum decomposition $\g\simeq\h\oplus b(\g/\h)$ leads to
a section of the central extension $\gl(\g)\til\rarrow\gl(\g)$, and
consequently to a section of the central extension $\g\til\rarrow\g$.
 Composing the section~$b$ with the latter section, we get
a section $\tilde b\:\g/\h\rarrow\g\til$ of the surjection
$\g\til\rarrow\g\til/(\h\oplus k)\simeq\g/\h$.

 Consider the composition $(\g/\h\ot_kN)^\h\rarrow \g/\h\ot_kN\rarrow
\g\til\ot_kN\rarrow N$ of the natural injection $(\g/\h\ot_kN)^\h
\rarrow\h$, the map $\g/\h\ot_k N\rarrow\g\til\ot_kN$ induced by
the section $\tilde b\:\g/\h\rarrow\g\til$, and the $\g\til$\+action
map $\g\til\ot_kN\rarrow N$.
 Let us check that this composition does not depend on the choice of~$b$
and its image lies in the subspace of invariants $N^\h\subset N$,
so it provides the desired natural map $(\g/\h\ot_kN)^\h\rarrow N^\h$.

 Let $u\ot n$ be a formal notation for an element of $\g/\h\ot_kN$.
 Denote by $n\mpsto n_{\{-1\}}\ot n_{\{0\}}$ the map
$N\rarrow\h\dual\ot_kN$ corresponding to the $\h$\+action map
$\h\times N\rarrow N$.
 Rewriting the identity $x\tilde b(u)n=\tilde b(u)xn+[x,\tilde b(u)]n$
for $x\in\h$ in the notation of~\ref{commutator-components-notation},
we obtain the identity $(\tilde b(u)n)_{\{-1\}}\ot
(\tilde b(u)n)_{\{0\}} = n_{\{-1\}}\ot \tilde b(u)n_{\{0\}} +
u_{\{-1\}}\ot \tilde b(u_{\{0\}})n - u(n_{\{-1\}})\ot n_{\{0\}} -
u_{\{0\}}(u_{\{-1\}})\ot n$.
 Now whenever $u\ot n$ is an $\h$\invariant{} element of $\g/\h\ot_kN$
one has $n_{\{-1\}}\ot u\ot n_{\{0\}} + u_{\{-1\}}\ot u_{\{0\}}\ot n
= 0$, hence $(\tilde b(u)n)_{\{-1\}}\ot (\tilde b(u)n)_{\{0\}} = 0$
and $\tilde b(u)n$ is an $\h$\invariant{} element of~$N$.

 Let $b_1\:\g/\h\rarrow\h$ be another section of the surjection
$\g\rarrow\g/\h$ and $\tilde b_1\:\g/\h\rarrow\g\til$ be
the corresponding section of the surjection $\g\til\rarrow\g/\h$.
 According to~\ref{tate-change-of-section}, we have $\tilde b_1=
\tilde b + t + \tilde t$ with a map $t=b_1-b\:\g/\h\rarrow \h$ and
the linear function $\tilde t\:\g/\h\rarrow k$ given by the formula
$\tilde t(u)=-\lan u_{\{-1\}},t(u_{\{0\}})\ran$.
 Let $u\ot n$ be an $\h$\invariant{} element of $\g/\h\ot N$; then
the equation $n_{\{-1\}}\ot u\ot n_{\{0\}} + u_{\{-1\}}\ot u_{\{0\}}
\ot n = 0$ implies $\lan n_{\{-1\}},t(u)\ran n_{\{0\}} +
\lan u_{\{-1\}},t(u_{\{0\}})\ran n = 0$ and $t(u)n-\tilde t(u)n=0$.

 The cokernel $N_{\g,\h}$ of the natural map $(\g/\h\ot_kN)^\h
\rarrow N^\h$ that we have constructed is called the space of
\emph{$(\g,\h)$\semiinvariant s} of a discrete $\g$\module{} $N$.
 The $(\g,\h)$\semiinvariant s are a mixture of $\h$\invariant s
and ``coinvariants along $\g/\h$''.

\subsubsection{}
 For a topological Lie algebra $\h$ and an $\h$\contramodule{} $P$
the space  of $\h$\coinvariant s $P_\h$ is defined as the maximal
quotient contramodule of $P$ where $\h$ acts by zero, i.~e.,
the cokernel of the contraaction map $\h\ot\comp P\rarrow P$.

 Let $\g$ be a Tate Lie algebra with a compact open subalgebra~$\h$. 
 Let $P$ be a $\g\til$\contramodule{} where the unit element of
$k\subset\g\til$ acts by the identity.
 We would like to construct a natural map $P_\h\rarrow
\Hom_k(\g/\h,P)_\h$, where the $\h$\+contraaction in $\Hom_k(\g/\h,P)$
is induced by the discrete action of~$\h$ in $\g/\h$ and
the $\h$\+contraaction in~$P$ as explained
in~\ref{lie-contramodules}--\ref{tate-lie-contramodules}.

 As above, choose a section $b\:\g/\h\rarrow\g$ and construct
the corresponding section $\tilde b\:\g/\h\rarrow\g\til$.
 Consider the composition $P\rarrow\Hom_k(\g\til,P)\rarrow
\Hom_k(\g/\h,P)\rarrow\Hom_k(\g/\h,P)_\h$ of the map $P\rarrow
\Hom_k(\g\til,P)$ corresponding to the action of $\g\til$ in $P$
induced by the contraaction of $\g\til$ in $P$, the map
$\Hom_k(\g\til,P)\rarrow\Hom_k(\g/\h,P)$ induced by
the section~$\tilde b$, and the natural surjection $\Hom_k(\g/\h,P)
\rarrow\Hom_k(\g/\h,P)_\h$.
 Let us check that this composition factorizes through the natural
surjection $P\rarrow P_\h$ and does not depend on the choice of~$b$,
so it defines the desired map $P_\h\rarrow\Hom_k(\g/\h,P)_\h$.

 Let $f$ a linear function $\h\dual\rarrow P$ and $\pi_P(f)\in P$
be its image under the contraaaction map.
 The image of $\pi_P(f)$ under the composition $P\rarrow
\Hom_k(\g\til,P)\rarrow\Hom_k(\g/\h,P)$ is given by the formula
$u\mpsto \tilde b(u)\pi_P(f) = \pi_P(x^*\mapsto \tilde b(u)f(x^*))
- \tilde b(u_{\{0\}})f(u_{\{-1\}}) + \pi_P(x^*\mapsto f(u(x^*))) -
f(u_{\{0\}}(u_{\{-1\}}))$ in the notation
of~\ref{commutator-components-notation}.
 This element of $\Hom_k(\g/\h,P)$ is the image of the element
$g\in\Hom_k(\h\dual,\Hom_k(\g/\h,P))$ given by the formula $g(x^*)(u) =
\tilde b(u)f(x^*) + f(u(x^*))$ under the contraaction map.

 If $b_1\:\g/\h\rarrow\g$ is a different section, then $\tilde b_1(u)
= b_1(u) + t(u) - \lan u_{\{-1\}},t(u_{\{0\}})\ran$ and for any
$p\in P$ the element of $\Hom_k(\g/\h,P)$ given by the formula $u\mpsto
t(u)p - \lan u_{\{-1\}},t(u_{\{0\}})\ran p$ is the image of the element
$g\in\Hom_k(\h\dual,\Hom_k(\g/\h,P))$ given by the formula $g(x^*)(u)
= \lan x^*, t(u)\ran p$ under the contraaction map.

 The kernel $P^{\g,\h}$ of the natural map $P_\h\rarrow
\Hom_k(\g/\h,P)_\h$ is called the space of
\emph{$(\g,\h)$\semicontrainvariant s} of a $\g$\contramodule{} $P$.
  The $(\g,\h)$\semicontrainvariant s are a mixture of 
$\h$\coinvariant s and ``invariants along $\g/\h$''.

\begin{rmk}
 The above definitions of $(\g,\h)$\semiinvariant s and 
$(\g,\h)$\semicontrainvariant s agree with the definitions
from~\ref{finite-dim-semiinvariants-defining-map}--%
\ref{finite-dim-semicontrainvs-defining-map} up to twists with
a one-dimensional vector space $\det(\h)$, essentially for
the following reason.
 When $\g$ is a discrete Lie algebra, the central extension
$\g\til\rarrow\g$ has a canonical splitting induced by the canonical
splitting of the central extension $\gl(\g)\til\rarrow\gl(\g)$.
 When $\g$ is a Tate Lie algebra and $\h\subset\g$ is a compact
open Lie subalgebra, the central extension $\g\til\rarrow\g$
has a canonical splitting over~$\h$.
 When $\g$ is a discrete Lie algebra and $\h\subset\g$ is 
a finite-dimensional Lie subalgebra, these two splittings do
\emph{not} agree over~$\h$; instead, they differ by the modular
character of the Lie algebra~$\h$.
\end{rmk}

\subsubsection{}
 Let $\C$ be a coalgebra endowed with a coaugmentation (morphism of
coalgebras) $e\: k\rarrow\C$, \ $\h$ be a compact Lie algebra, and
$\psi\:\C\times\h\rarrow k$ be a pairing such that the map
$\check\psi\:\h\dual\rarrow\C$ is a morphism of coalgebras
and $\psi$ annihilates~$e(k)$.
 For any right $\C$\comodule{} $\N$, the maximal subcomodule of
$\N$ where the coaction of~$\C$ is trivial can be described as
the cotensor product $\N\oc_\C k$.
 Here a coaction of~$\C$ is called trivial if it is induced by~$e$;
the vector space~$k$ is endowed with the trivial coaction.
 There is a natural injective map $\N\oc_\C k\rarrow \N^\h$, which is
an isomorphism provided that the assumption~(iv)
of~\ref{tate-harish-chandra-central-extension} holds.

 Analogously, for any left $\C$\contramodule{} $\P$ the maximal
quotient $\C$\contramodule{} of $\P$ with the trivial contraaction
can be described as the space of cohomomorphisms $\Cohom_\C(k,\P)$.
 There is a natural surjective map $\P_\h\rarrow\Cohom_\C(k,\P)$,
which is also an isomorphism provided that the condition~(iv) holds.
 Indeed, it suffices to consider the case when $\P=\Hom_k(\C,E)$ is
an induced $\C$\contramodule; in this case one only has to check that
the kernel of the composition $\C\rarrow\C\ot_k\C\rarrow\C\ot_k\h\dual$
of the comultiplication map and the map induced by the map~$\check\psi$
coincides with~$e(k)$.

 Let $\C$ be a commutative Hopf algebra.
 Then for any right $\C$\comodule{} $\N$ and left $\C$\comodule{} $\M$
there is a natural isomorphism $\N\oc_\C\M\simeq(\N\ot_k\M)\oc_\C k$,
where the right $\C$\comodule{} structure on $\N\oc_\C\M$ is defined
using the antipode and multiplication in~$\C$.
 Analogously, for any left $\C$\comodule{} $\M$ and left
$\C$\contramodule{} $\P$ there is a natural isomorphism
$\Cohom_\C(\M,\P)\simeq\Cohom_\C(k,\Hom_k(\M,\P))$.

\subsubsection{}    \label{tate-semitensor-product-semiinvariants}
 Now let $\kap\:(\g',\C)\rarrow(\g,\C)$ be a central extension of
Tate Harish-Chandra pairs with the kernel~$k$ satisfying
the assumption~(iv) of~\ref{tate-harish-chandra-central-extension}
and $\S^l_\kap=\S\simeq\S^r_{\kap+\kap_0}$ be the corresponding
semialgebra over~$\C$.

\begin{lem}
\textup{(a)} Let\/ $\bN$ be a Harish-Chandra module with the central
charge $-\kap-\kap_0$ and\/ $\bM$ be a Harish-Chandra module with
the central charge~$\kap$ over\/~$(\g,\C)$; in other words, $\bN$ is
a right\/ $\S^r_{\kap+\kap_0}$\semimodule{} and\/ $\bM$ is a left\/
$\S^l_\kap$\semimodule.
 Assume that either\/ $\bN$ or\/ $\bM$ is a coflat\/ $\C$\comodule.
 Then there is a natural isomorphism\/ $\bN\os_\S\bM\simeq
(\bN\ot_k\bM)_{\g,\h}$, where the tensor product\/ $\bN\ot_k\bM$ is
a Harish-Chandra module with the central charge~$-\kap_0$. \par
\textup{(b)} Let\/ $\bM$ be a Harish-Chandra module with the central
charge~$\kap$ and\/ $\bP$ be a Harish-Chandra contramodule with
the central charge $\kap+\kap_0$ over\/~$(\g,\C)$; in other words,
$\bM$ is a left\/ $\S^l_\kap$\semimodule{} and\/ $\bP$ is a left\/
$\S^r_{\kap+\kap_0}$\semicontramodule.
 Assume that either\/ $\bM$ is a coprojective $\C$\comodule, or\/ $\bP$
is a coinjective $\C$\contramodule.
 Then there is a natural isomorphism\/ $\SemiHom_\S(\bM,\bP)\simeq
\Hom_k(\bM,\bP)^{\g,\h}$, where the space\/ $\Hom_k(\bM,\bP)$ is 
a Harish-Chandra contramodule with the central charge~$\kap_0$.
\end{lem}

\begin{proof}
 Part~(a): denote by $\eta_i\:\bN\oc_\C F_i\S\oc_\C\bM\rarrow
\bN\oc_\C\bM$ the map equal to the difference of the map induced by
the semiaction map $F_i\S\oc_\C\bM\rarrow\bM$ and the map induced by
the semiaction map $\bN\oc_\C F_1\S\rarrow\bN$.
 The map $\eta_0$ vanishes and the quotient space $(\bN\oc_\C F_1\S
\oc_\C\bM)/(\bN\oc_\C F_0\S\oc_\C\bM$ is isomorphic to $\bN\oc_\C
(F_1\S/F_0\S)\oc_\C\bM$, hence the induced map $\bar\eta_1\:\bN\oc_\C
(F_1\S/F_0\S)\oc_\C\bM\rarrow\bN\oc_\C\bM$.
 The cokernel of the map $\bar\eta_1$ coincides with the semitensor
product $\bN\os_\S\bM$ for the reasons explained
in~\ref{finite-dim-semitensor-product-semiinvariants}.
 The cotensor product $\bN\oc_\C(F_1\S/F_0\S)\oc_\C\bM\simeq
\bN\oc_\C(\g/\h\ot_k\C)\oc_\C\bM\simeq\bN\oc_\C(\g/\h\ot_k\bM)$ is
isomorphic the space of invariants $(\g/\h\ot_k\bM\ot_k\bN)^\h$ in
view of the assumption~(iv); this isomorphism coincides with
the isomorphism $\bN\oc_\C(F_1\S/F_0\S)\oc_\C\bM\simeq
(\g/\h\ot_k\bM\ot_k\bN)^\h$ induced by the isomorphism
$F_1\S/F_0\S\simeq\C\ot_k\g/\h$.
 Let us check that this isomorphism identifies the map~$\bar\eta_1$
with the map whose cokernel is, by the definition, the space of
semiinvariants $(\bN\ot_k\bM)_{\g,\h}$.

 Choose a section $b'\:\g/\h\rarrow\g'$ and consider the corresponding
section $b''\:\g/\h\rarrow\g''$.
 There is an isomorphism of right $\C$\comodule s $\C\oplus\g/\h\ot_k\C
\simeq F_1\S^l_\kap$ given by the formula $c_1'+u'\ot c'\mpsto
1\ot_{U(\h)}c_1'+b'(u')\ot_{U(\h)}c'$ and an analogous isomorphism of
left $\C$\comodule s $\C\oplus\C\ot_k\g/\h\simeq F_1\S^r_{\kap+\kap_0}$
given by the formula $c_1''+c''\ot u''\mpsto c_1''\ot_{U(\h)}1 +
c''\ot_{U(\h)}b''(u'')$.
 The induced isomorphism $\bM\oplus \g/\h\ot_k\bM\simeq F_1\S^l_\kap
\oc_\C\bM\simeq F_1U_\kap(\g)\ot_{U(\h)}\bM$ is given by the formula
$m_1 + u\ot m\mpsto 1\ot_{U(\h)}m_1 + b'(u)\ot_{U(\h)}m$.
 Now let $z = n\ot (c_1'+u'\ot_{U(\h)} c')\ot m = n\ot (c_1''+
c''\ot_{U(\h)} u'')\ot m$ be an element of $\bN\oc_\C F_1\S\oc_\C\bM$.
 Then the corresponding element of $\bN\oc_\C F_1U_\kap(\g)\ot_{U(\h)}
\bM$ can be written as $\eps(c_1')n\ot 1\ot_{U(\h)}m + 
\eps(c')n\ot b'(u')\ot_{U(\h)}m$, hence the image of $z$ under
the map $\bN\oc_\C F_1\S\oc_\C\bM\rarrow\bN\oc_\C\bM$ induced by
the semiaction map $F_1\S\oc_\C\bM\rarrow\bM$ is equal to
$\eps(c_1')n\ot m+\eps(c')n\ot b'(u')m$.
 Analogously, the image of $z$ under the map $\bN\oc_\C F_1\S\oc_\C\bM
\rarrow\bN\oc_\C\bM$ induced by the semiaction map $\bN\oc_\C F_1\S
\rarrow\bN$ is equal to $\eps(c_1'')n\ot m - \eps(c'')b''(u'')n\ot m$.
 One has $\eps(c_1')=\eps(c_1'')$ by the condition~(c)
of~\ref{semialgebra-isomorphism-characterized}.
 Thus $\eta_1(z) = n\ot b'(u)m + b''(u)n\ot m = \tilde b(u)(n\ot m)$,
where $u=\eps(c')u'=\eps(c'')u''$.
 Part~(a) is proven; the proof of part~(b) is completely analogous.
\end{proof}

\subsection{Semi-infinite homology and cohomology}

\subsubsection{}
 A \emph{discrete right module} $N$ over a topological associative
algebra $R$ is a right $R$\module{} such that the action map
$N\times R\rarrow N$ is continuous with respect to the discrete
topology of~$N$.
 Equivalently, a right $R$\module{} $N$ is discrete if the annihilator
of any element of $N$ is an open right ideal in~$R$.

 Let $A$ and $B$ be topological associative algebras in which open
right ideals form bases of neighborhoods of zero.
 Then the topological tensor product $A\ot^!B$ has a natural structure
of topological associative algebra with the same property.
 The tensor product of a discrete right $A$\module{} and a discrete
right $B$\module{} is naturally a discrete right $A\ot^!B$\module.

 Let $\kap\:\g'\rarrow\g$ be a central extension of topological Lie
algebras with the kernel~$k$.
 Then the modified enveloping algebra $U_\kap(\g)=
U(\g')/(1_{U(\g')}-1_{\g'})$ can be endowed with the topology where
right ideals generated by open subspaces of $\g'$ form a base of
neighborhoods of zero~\cite{BD1,Beil}.
 Denote the completion of $U_\kap(\g)$ with respect to this topology
by $U_\kap\comp(\g)$; this is a topological associative algebra.
 The category of discrete $\g'$\module s where the unit element of
$k\subset\g'$ acts by minus the identity is isomorphic to the category
of discrete right $U_\kap\comp(\g)$\module s.

\subsubsection{}  \label{assoc-contramodules}
 Let $R$ be a topological associative algebra where open right ideals
form a base of neighborhoods of zero.
 Then for any $k$\+vector space $P$ there is a natural map $R\ot\comp
(R\ot\comp P)\rarrow R\ot\comp P$ induced by the multiplication in~$R$;
it is constructed as the projective limit over all open right ideals
$U\subset R$ of the maps $R/U\ot_k(R\ot\comp P)\rarrow R/U\ot_k P$
induced by the discrete right action of $R$ in $R/U$.
 A \emph{left contramodule} over $R$ is a vector space $P$ endowed with
a linear map $R\ot\comp P\rarrow P$ satisfying the following
contraassociativity and unity equations.
 First, the two maps $R\ot\comp(R\ot\comp P)\rarrow R\ot\comp P$,
one induced by the multiplication in $R$ and the other induced by
the contraaction map $R\ot\comp P\rarrow P$, should have equal
compositions with the contraaction map.
 Second, the composition $P\rarrow R\ot\comp P\rarrow P$ of the map
induced by the unit of~$R$ and the contraaction map should be
equal to the identity endomorphism of~$P$.

 The category of left $R$\contramodule s is abelian and there is
a natural exact forgetful functor from it to the category of left
modules over the algebra $R$ considered without any topology
(cf.\ Remark~\ref{contraflat-contramodules}).
 Notice also the isomorphisms $R\ot\comp(R\ot\comp P)\simeq
R\wot R\wot P\simeq (R\wot R)\ot\comp P$, demonstrating the similarity
of the above definition with the definition of a contramodule over
a Lie algebra given in~\ref{lie-contramodules}.
 The above natural map $R\ot\comp(R\ot\comp P)\rarrow R\ot\comp P$ is
induced by the continuous multiplication map $R\wot R\rarrow R$, which
exists for any topological associative algebra $R$ where open right
ideals form a base of neighborhoods of zero.
 Just as for Lie algebras, a structure of a discrete right $R$\module{}
on a vector space $N$ is given by a continuous linear map $N\ot^*R
\simeq N\wot R\rarrow N$, while a structure of a left
$R$\contramodule{} on a vector space $P$ is given by a discontinuous
linear map $R\ot^!P\simeq R\wot P\rarrow P$, where $N$ and $P$
are endowed with discrete topologies.

 For any discrete right $R$\module{} $N$ and any $k$\+vector space $E$,
the vector space $\Hom_k(N,E)$ has a natural structure of left
contramodule over~$R$.
 The contraaction map $R\ot\comp\Hom_k(N,E)\rarrow \Hom_k(N,E)$ is
constructed as the projective limit over all open right ideals
$U\subset R$ of the maps $R/U\ot_k\Hom_k(N,E)\rarrow\Hom_k(N^U,E)$
given by the formulas $\bar r \ot_k g\mpsto (n\mpsto g(n\bar r))$
for $\bar r \in R/U$, \ $g\in\Hom_k(N,E)$, and $n\in N^U$, where
$N^U\subset N$ denotes the subspace of all elements of~$N$
annihilated by~$U$.

 More generally, let $A$ and $B$ be topological associative algebras
where open right ideals form bases of neighborhoods of zero, $N$ be
a discrete right $B$\module, and $P$ be a left $A$\contramodule.
 Then the vector space $\Hom_k(N,P)$ has a natural structure of
contramodule over $A\ot^!B$.
 The contraaction map $(A\ot^!B)\ot\comp\Hom_k(N,E)\rarrow \Hom_k(N,E)$
is constructed as the projective limit over all open right ideals
$U\subset B$ of the compositions $A\ot\comp(B/U\ot_k\Hom_k(N,P))\rarrow
A\ot\comp\Hom_k(N^U,P)\rarrow\Hom_k(N^U\;A\ot\comp P)\rarrow
\Hom_k(N^U,P)$, where the first map is induced by the right $B$\+action
in $N$ and the third map is induced by the $A$\+contraaction in~$P$.

\subsubsection{}  \label{contra-enveloping-action}
 Let $\kap\:\g'\rarrow\g$ be a central extension of topological
Lie algebras with the kernel~$k$.

\begin{thm}
 Assume that the topological Lie algebra $\g'$ has a countable base
of neighborhoods of zero consisting of open Lie subalgebras.
 Then the category of $\g'$\contramodule s where the unit element of
$k\subset\g'$ acts by the identity is isomorphic to the category of
left contramodules over the topological algebra $U_\kap\comp(\g)$.
\end{thm}

\begin{proof}
 It is easy to see that the composition $\g'\ot\comp P\rarrow
U_\kap\comp(\g)\ot\comp P\rarrow P$ defines a $\g'$\contramodule{}
structure on any left $U_\kap\comp(\g)$\contramodule{} $P$ (so, in
particular, $U_\kap\comp(\g)$ itself is a $\g'$\contramodule).
 Let us construct the functor in the opposite direction.

 The standard homological Chevalley complex $\dsb\rarrow
\bigwedge^2_k\g'\ot_kU_\kap(\g)\rarrow\g'\ot_kU_\kap(\g)\rarrow
U_\kap(\g)\rarrow 0$ is acyclic.
 For any open Lie subalgebra $\h\subset\g'$ not containing
$k\subset\g'$, the complex $\dsb\rarrow\bigwedge^2_k\h\ot_kU_\kap(\g)
\rarrow\h\ot_k U_\kap(\g)\rarrow \h U_\kap(\g)\rarrow0$ is an acyclic
subcomplex of the previous complex.
 Taking the quotient complex and passing to the projective limit
over~$\h$, we obtain a split exact complex of topological vector spaces
$\dsb\rarrow\bigwedge^{\s,2}\g'\ot^!U_\kap(\g)\rarrow\g'\ot^!U_\kap(\g)
\rarrow U_\kap\comp(\g)\rarrow0$, where we denote by
$\bigwedge^{\s,i}\g'$ the completion of $\bigwedge^i_k\g'$ with respect
to the topology with a base of neighborhoods of zero formed by
the subspaces $\bigwedge^i_k\h$ and the enveloping algebra $U_\kap(\g)$
is considered as a discrete topological vector space.
 Applying the functor $\ot\comp P$, we obtain an exact sequence of
vector spaces $\bigwedge^{\s,2}\g'\ot\comp(U_\kap(\g)\ot_k P)\rarrow
\g'\ot\comp(U_\kap(\g)\ot_kP)\rarrow U_\kap\comp(\g)\ot\comp P\rarrow0$
for any $k$\+vector space~$P$.

 Now let $P$ be a $\g'$\contramodule{} where the unit element of
$k\subset\g'$ acts by the identity; then, in particular, $P$ is
a $\g'$\module{} and a $U_\kap(\g)$\module.
 It is clear from the above exact sequence that the composition
$\g'\ot\comp(U_\kap(\g)\ot_k P)\rarrow\g'\ot\comp P\rarrow P$ of the map
induced by the $U_\kap(\g)$\+action map and the $\g'$\+contraaction map
factorizes through $U_\kap\comp(\g)\ot\comp P$, providing
the desired contraaction map $U_\kap\comp(\g)\ot\comp P\rarrow P$.

 Let us check that this contraaction map satisfies
the contraassociativity equation.
 Any element $z$ of $U_\kap\comp(\g)\ot\comp P$ can be presented in
the form $z=\sum_{i=0}^\infty u_i\ot p_i$ with $u_i\in U_\kap(\g)$ and
$p_i\in P$, where $u_i\to 0$ in $U_\kap\comp(\g)$ as $i\to\infty$ and
the infinite sum is understood as the limit in the topology of
$U_\kap\comp(\g)\ot^! P$.
 Let us denote the image of the element $\sum_i u_i\ot p_i$ under
the contraaction map $U_\kap\comp(\g)\ot\comp P\rarrow P$ by
$\sum_i u_ip_i\in P$.
 In this notation, the $U_\kap\comp(\g)$\+contraaction map is defined
by the formula $\sum_i (x_{i_1}x_{i_2}\dsb x_{i_{k_i}})p_i = 
\sum_i x_{i_1}(x_{i_2}\dsb x_{i_{k_i}}p_i)$ for any $x_{i_t}\in \g'$
and $p_i\in P$ such that $x_{i_1}\to 0$ in $\g'$ as $i\to\infty$.
 We have to show that $\sum_i u_i\sum_j v_{ij}p_{ij} = \sum_{i,j}
(u_iv_{ij})p_{ij}$ for any $u_i$, $v_{ij}\in U_\kap(\g)$ and
$p_{ij}\in P$ such that $u_i\to 0$ as $i\to\infty$ and $v_{ij}\to 0$
as $j\to\infty$ for any~$i$.

 Let us first check that $\sum_i x_i\sum_jy_{ij}p_{ij} = \sum_{i,j}
(x_iy_{ij})p_{ij}$ for any $x_i$, $y_{ij}\in\g'$ and $p_{ij}\in P$
such that $x_i\to 0$ in $\g'$ as $i\to\infty$ and $y_{ij}\to 0$
in $\g'$ as $j\to\infty$ for any~$i$.
 Choose an integer $j_i$ for each~$i$ such that $\{y_{ij}\mid
j>j_i\}$ converges to zero in $\g'$ as $i+j\to\infty$.
 Then we have $\sum_{i,j}(x_iy_{ij})p_{ij} = \sum_{j\le j_i}
x_i(y_{ij}p_{ij}) + \sum_{j>j_i}y_{ij}(x_ip_{ij}) + \sum_{j>j_i}
[x_i,y_{ij}]p_{ij}$.
 To check that $\sum_i x_i\sum_{j>j_i}y_{ij}p_{ij} =
\sum_{j>j_i}y_{ij}(x_ip_{ij}) + \sum_{j>j_i}[x_i,y_{ij}]p_{ij}$,
apply the equation on the contraaction map of a contramodule
over a topological Lie algebra to the element $\sum_{j>j_i}
x_i\wedge y_{ij}\ot p_{ij}$ of the vector space
$\bigwedge^{*,2}(\g')\ot\comp P$.

 It follows that $\sum_i x_i\sum_j v_{ij}p_{ij} = \sum_{i,j}
(x_iv_{ij})p_{ij}$ for any $x_i\in\g'$, \ $v_{ij}\in U_\kap(\g)$, and
$p_{ij}\in P$ such that $x_i\to 0$ in $\g'$ as $i\to\infty$ and
$v_{ij}\to 0$ in $U_\kap\comp(\g)$ as $j\to\infty$ for any~$i$.
 Indeed, assuming that $v_{ij}=y_{ij1}y_{ij2}\dsb y_{ijk_{ij}}$, where
$y_{ijt}\in\g'$ and $y_{ij1}\to 0$ in $\g'$ as $j\to\infty$, we have
$\sum_i x_i\sum_j (y_{ij1}y_{ij2}\dsb y_{ijk_{ij}})p_{ij} =
\sum_i x_i\sum_j y_{ij1}(y_{ij2}\dsb y_{ijk_{ij}}p_{ij}) =
\sum_{i,j}(x_iy_{ij1})(y_{ij2}\dsb y_{ijk_{ij}}p_{ij}) =
\sum_{i,j}(x_iy_{ij1}y_{ij2}\dsb y_{ijk_{ij}})p_{ij}$.
 Furthermore, it follows that $x_1\dsb x_s\sum_j v_jp_j=\sum_j
(x_1\dsb x_sv_j)p_j$ for any $x_t\in\g'$, \ $v_j\in U_\kap(\g)$,
and $p_j\in P$ such that $v_j\to 0$ in $U_\kap\comp(\g)$ as
$j\to\infty$.
 Now to check that $\sum_i u_i\sum_j v_{ij}p_{ij} = \sum_{i,j}
(u_iv_{ij})p_{ij}$, we can assume that $u_i=x_{i_1}x_{i_2}\dsb
x_{i_{k_i}}$, where $x_{i_t}\in\g'$ and $x_{i_1}\to 0$ in $\g'$
as $i\to\infty$.
 Then we have $\sum_i(x_{i_1}x_{i_2}\dsb x_{i_{k_i}})\sum_j v_{ij}p_{ij}
= \sum_i x_{i_1}(x_{i_2}\dsb x_{i_{k_i}}\sum_j v_{ij}p_{ij}) =
\sum_i x_{i_1}\sum_j (x_{i_2}\dsb x_{i_{k_i}}v_{ij}p_{ij}) = 
\sum_{i,j} (x_{i_1}x_{i_2}\dsb x_{i_{k_i}}v_{ij})p_{ij}$.
\end{proof}

\begin{qst}
 Can one construct an isomorphism between the categories of
``$\g$\contramodule s with central charge~$\kap$'' and left
$U_\kap\comp(\g)$\contramodule s without the countability
assumption on the topology of~$\g'$?
\end{qst}

\subsubsection{}  \label{noncountable-contra-identity}
 The following weaker version of Theorem~\ref{contra-enveloping-action}
holds without the countability assumption.
 Let $\kap\:\g'\rarrow\g$ be a central extension of topological
Lie algebras with the kernel~$k$; assume that open subalgebras from
a base of neighborhoods of zero in~$\g'$.
 Let $B$ be a topological associative algebra where open right ideals
form a base of neighborhoods of zero, $N$ be a discrete right
$B$\module, and $P$ be a $\g'$\contramodule{} where the unit element
of $k\subset\g'$ acts by the identity.
 Then one can define the contraaction map $(\g'\ot^!B)\ot\comp
\Hom_k(N,P)\rarrow\Hom_k(N,P)$ as in~\ref{assoc-contramodules}.

 Consider the iterated contraaction map $((\g'\ot^!B)\wot(\g'\ot^!B))
\ot\comp\Hom_k(N,P)\simeq(\g'\ot^!B)\ot\comp((\g'\ot^!B)\ot\comp
\Hom_k(N,P))\rarrow\Hom_k(N,P)$.
 It was noticed in~\cite{Beil} that a topological associative algebra
$A$ has the property that open right ideals form a base of neighborhoods
of zero if and only if the multiplication map $A\ot^*A\rarrow A$
factorizes through $A\wot A$.
 Let $K$ denote the kernel of the multiplication map $(\g'\ot^!B)
\wot(\g'\ot^!B)\rarrow U_\kap\comp(g)\ot^!B$.
 We claim that the composition of the injection $K\ot\comp\Hom_k(N,P)
\rarrow ((\g'\ot^!B)\wot(\g'\ot^!B))\ot\comp\Hom_k(N,P)$ and
the iterated contraaction map $((\g'\ot^!B)\wot(\g'\ot^!B))
\ot\comp\Hom_k(N,P)\rarrow\Hom_k(N,P)$ vanishes.

 For any topological vector spaces $U$, $V$, $X$, $Y$ there is
a natural map $(U\ot^!X)\wot(V\ot^!Y)\rarrow(U\wot V)\ot^!(X\wot Y)$.
 The composition $(\g'\ot^!B)\wot(\g'\ot^!B)\rarrow(\g'\wot\g')
\ot^!(B\wot B)\rarrow (\g'\wot\g')\ot^!B$ induces the map
$((\g'\ot^!B)\wot(\g'\ot^!B))\ot\comp\Hom_k(N,P)\rarrow
((\g'\wot\g')\ot^!B)\ot\comp\Hom_k(N,P)$.
 A contraaction map $((\g'\wot\g')\ot^!B)\ot\comp\Hom_k(N,P)\rarrow
\Hom_k(N,P)$ can be defined in terms of the discrete right action
of $B$ in~$N$ and the iterated contraaction map $(\g'\wot\g')\ot\comp P
\rarrow P$.
 The iterated contraaction map $((\g'\ot^!B)\wot(\g'\ot^!B))
\ot\comp\Hom_k(N,P)\rarrow\Hom_k(N,P)$ is equal to the composition
$((\g'\ot^!B)\wot(\g'\ot^!B))\ot\comp\Hom_k(N,P)\rarrow
((\g'\wot\g')\ot^!B)\ot\comp\Hom_k(N,P)\rarrow\Hom_k(N,P)$ of
the above induced map and contraaction map.
 Let $Q$ denote the kernel of the multiplication map $\g'\wot\g'
\rarrow U_\kap\comp(\g)$.
 The image of $K$ under the map $(\g'\ot^!B)\wot(\g'\ot^!B)\rarrow
(\g'\wot\g')\ot^!B$ is contained in $Q\ot^!B$.
 So it suffices to check that the composition of the injection
$Q\ot\comp P\rarrow(\g'\wot\g')\ot\comp P$ and the iterated
contraaction map $(\g'\wot\g')\ot\comp P\rarrow P$ vanishes.

 The topological vector space $Q$ is the topological projective limit
of the kernels of multiplication maps $\g'/\h\ot^*\g'\rarrow
U_\kap(\g)/\h U_\kap(\g)$ over all open subalgebras $\h\subset\g'$
not containing $k\subset\g'$.
 Since the intersection of $\h U_\kap(\g)$ and $\g'{}^2$ inside
$U_\kap(\g')$ is equal to $\h\g'$, the kernel of the (nontopological)
multiplication map $\g'\ot_k\g'\rarrow U_\kap(\g')$ maps surjectively
onto the kernels we are interested in.
 This nontopological kernel is the image of the map $\bigwedge_k^2(\g')
\rarrow\g'\ot_k\g'$ given by the formula $x\wedge y\mpsto x\ot y -
y\ot x - 1\ot[x,y]$.
 The kernel of the composition $\bigwedge^2_k(\g')\rarrow\g'\ot_k\g'
\rarrow\g'/\h\ot_k\g'$ is the subspace $\bigwedge^2_k(\h)\subset
\bigwedge^2_k(\g')$.
 Hence the kernel of the map $\g'/\h\ot^*\g'\rarrow
U_\kap(\g)/\h U_\kap(\g)$ is the subspace $\bigwedge^2_k(\g')/
\bigwedge^2_k(\h)\subset\g'/\h\ot^*\g'$, embedded by the above formula,
endowed with the induced topology of a closed subspace.
 One can easily check that this topology on $\bigwedge^2_k(\g')/
\bigwedge^2_k(\h)$ is the topology of the quotient space
$\bigwedge^{*,2}(\g')/\bigwedge^{*,2}(\h)$.
 Thus the topological vector space $K$ is isomorphic to
$\bigwedge^{*,2}(\g')$. 

\subsubsection{}
 The following constructions are due to Beilinson and
Drinfeld~\cite{BD1}.

 Let $V$ be a Tate vector space and $E\subset V$ be a compact open
subspace.
 The graded vector space of \emph{semi-infinite forms}
$\bigwedge_E^{\infty/2}(V)=\bigoplus_i\bigwedge_E^{\infty/2+i}(V)$ is
defined as the inductive limit of the spaces $\bigwedge_k(V/U)\ot_k
\det(E/U)\dual$ over all compact open subspaces $U\subset E$.
 Here $\det(X)$ denotes the top exterior power of a finite-dimensional
vector space $X$ and $\bigwedge_k(W)$ denotes the direct sum of all
exterior powers of a vector space~$W$; the grading on $\bigwedge_k(V/U)
\ot_k\det(E/U)\dual$ is defined so that $\bigwedge^j_k(V/U)\ot_k
\det(E/U)\dual$ is the component of degree $j-\dim(E/U)$.
 The limit is taken over the maps induced by the natural maps
$\bigwedge^j_k(V/U')\ot_k\det(U'/U'')\rarrow\bigwedge^{j+m}_k(V/U'')$,
where $U''\subset U'$ and $m=\dim(U'/U'')$.
 The spaces of semi-infinite forms corresponding to different compact
open subspaces $E\subset V$, only differ by a dimensional shift and
a determinantal twist: if $F\subset V$ is another compact open subspace,
then there are natural isomorphisms $\bigwedge_F^{\infty/2+i}(V)\simeq
\bigwedge_E^{\infty/2+i+\dim(E,F)}(V)\ot_k\det(E,F)$, where $\dim(E,F)=
\dim(E/E\cap F)-\dim(F/E\cap F)$ and $\det(E,F)=\det(E/E\cap F)
\ot_k\det(F/E\cap F)\dual$.

 Denote by $\oCl(V)$ the algebra of endomorphisms of the vector space
$\bigwedge_E^{\infty/2}(V)$ endowed with the topology where annihilators
of finite-dimensional subspaces of $\bigwedge_E^{\infty/2}(V)$ form
a base of neighborhoods of zero.
 Clearly, the topological associative algebra $\oCl(V)$ does not
depend on the choice of a compact open subspace $E\subset V$; open
left ideals form a base of neighborhoods of zero in $\oCl(V)$.
 Denote by $\oCl^i(V)$ the closed subspace of homogeneous
endomorphisms of degree~$i$ in $\oCl(V)$.

 The Clifford algebra $\Cl(V\oplus V\dual)$ acts naturally in
$\bigwedge_E^{\infty/2}(V)$, so there is a morphism of associative
algebras $e\:\Cl(V\oplus V\dual)\rarrow\oCl(V)$; in particular,
the map~$e$ sends $V$ to $\oCl^1(V)$ and $V\dual$ to $\oCl^{-1}(V)$.
 Let $\bigwedge^{!,i}(V\dual)$ denote the completion of
$\bigwedge^i_k(V\dual)$ with respect to the topology with the base
of neighborhoods of zero formed by the subspaces
$U\wedge \bigwedge^{i-1}_k(V\dual)\subset\bigwedge^i_k(V\dual)$,
where $U\subset V\dual$ is an open subspace.
 The composition $\bigwedge^i_k(V\dual)\rarrow\Cl(V\oplus V\dual)
\rarrow\oCl(V)$ can be extended by continuity to a map
$\bigwedge^{!,i}_k(V\dual)\rarrow\oCl(V)$,
which we will denote also by~$e$.
 The construction of~\ref{clifford-central-extension-definition}
provides a morphism of topological Lie algebras $\gl(V)\til
\rarrow\oCl^0(V)$.

 Let $\g$ be a Tate Lie algebra and $\kap_0\:\g\til\rarrow\g$ be its
canonical central extension.
 Consider the topological tensor product $\g\til\ot^!\oCl(\g)^\rop$,
where $\oCl(\g)^\rop$ denotes the topological algebra opposite to
$\oCl(\g)$; this topological tensor product is a bimodule
over $\oCl(\g)^\rop$.
 The unit elements of $\oCl(\g)^\rop$ and $k\subset\g\til$ induce
embeddings of $\g\til$ and $\oCl(\g)^\rop$ into
$\g\til\ot^!\oCl(\g)^\rop$.
 Consider the difference of the composition $\g\til\rarrow
\gl(\g)\til\rarrow\oCl(\g)\simeq\oCl(\g)^\rop\rarrow\g\til\ot^!
\oCl(\g)^\rop$ and the embedding $\g\til\rarrow\g\til\ot^!\oCl(\g)^
\rop$; this difference maps $k\subset\g\til$ to zero and so induces
a natural map $l\:\g\rarrow\g\til\ot^!\oCl^0(\g)^\rop$.
 The composition of the map~$l$ with the embedding $\g\til\ot^!
\oCl^{-1}(\g)^\rop\rarrow U_{\kap_0}\comp(\g)\ot^!\oCl^{-1}(\g)^\rop$
is an anti-homomorphism of Lie algebras, i.~e., it transforms
the commutators to minus the commutators.

 Denote by $\delta\:\g\dual\rarrow\bigwedge^{!,2}(\g\dual)$
the continuous linear map given by the formula $\delta(x^*)=
x^*_{\{1\}}\wedge x^*_{\{2\}}$, where $\lan x^*, [x',x'']\ran =
\lan x^*_{\{1\}},x''\ran \lan x_{\{2\}},x'\ran - \lan x^*_{\{1\}},x'\ran
\lan x^*_{\{2\}},x''\ran$ for $x^*\in\g\dual$, \ $x'$, $x''\in\g$.
 Define the map $\chi\:\g\ot\g\dual\rarrow\g\til\ot^!\oCl^{-1}(\g)^\rop$
by the formula $\chi(x\ot x^*) = l(x)e(x^*)^\rop-e(x)^\rop
e(\delta(x^*))^\rop=e(x^*)^\rop l(x)-e(\delta(x^*))^\rop e(x)^\rop$,
where $a^\rop$ denotes the element of $\Cl(\g)^\rop$ corresponding to
an element $a\in\Cl(\g)$, and extend this map by continuity to
a map $\g\ot^!\g\dual\rarrow\g\til\ot^!\oCl^{-1}(\g)^\rop$.
 Identify $\g\ot^!\g\dual$ with $\End(\g)$ and set $\dd=\chi(\id_\g)\in
\g\til\ot^!\oCl^{-1}(\g)^\rop$. 

 Denote the image of $\dd$ under the embedding $\g\til\ot^!
\oCl^{-1}(\g)^\rop\rarrow U_{\kap_0}\comp(\g)\ot^!\oCl^{-1}(\g)^\rop$
also by~$\dd$.
 Using the equality $[l(x),e(y)^\rop]=-e([x,y])^\rop$, one can check
that $[\dd,e(x)^\rop]=l(x)$ and $[[\dd^2,e(x)^\rop],e(y)^\rop]=0$ for
all $x$, $y\in\g$, where $[\,\.,\,]$ denotes the supercommutator with
respect to the grading in which $U_{\kap_0}\comp(\g)\ot^!\oCl^i(\g)^
\rop$ lies in the degree~$i$.
 It is easy to see that any element of $\oCl^i(\g)$ supercommuting
with $e(x)$ for all $x\in\g$ is zero when $i<0$; hence the same applies
to elements of $U_{\kap_0}\comp(\g)\ot^!\oCl^i(\g)^\rop$ with $i<0$.
 It follows that $\dd$ is the unique element of $U_{\kap_0}\comp(\g)
\ot^!\oCl^{-1}(\g)^\rop$ satisfying the equation $[\dd,e(x)^\rop]=l(x)$
and that $\dd^2=0$.

\subsubsection{}
 Let $\g$ be a Tate Lie algebra and $E\subset\g$ be a compact open
vector subspace.
 Let $N$ be a discrete $\g\til$\module{} where the unit element of
$k\subset\g\til$ acts by minus the identity.
 Then $N$ can be considered as a discrete right
$U_{\kap_0}\comp(\g)$\module{} and $\bigwedge^{\infty/2}_E(\g)$ is
a discrete right $\oCl(\g)^\rop$\module, so the tensor product
$\bigwedge^{\infty/2}_E(\g)\ot_k N$ is a discrete right module over
$U_{\kap_0}\comp(\g)\ot^!\oCl(\g)^\rop$.
 The action of the element $\dd\in U_{\kap_0}\comp(\g)\ot^!
\oCl^{-1}(\g)^\rop$ defines a differential $d_{\infty/2}$ of
degree~$-1$ on the graded vector space $C_{\infty/2+\bu}^E(\g,N)$
with the components $C_{\infty/2+i}^E(\g,N) = 
\bigwedge_E^{\infty/2+i}(\g)\ot_k N$.
 One has $d_{\infty/2}^2=0$, since $\dd^2=0$; so
$C_{\infty/2+\bu}^E(\g,N)$ becomes a complex.
 This complex is called the \emph{semi-infinite homological complex}
and its homology is called the \emph{semi-infinite homology} of
a Tate Lie algebra $\g$ with coefficients in a discrete
$\g\til$\module~$N$.
 
 Let $P$ be a $\g\til$\contramodule{} where the unit element of
$k\subset\g\til$ acts by the identity.
 First assume that $\g$ has a countable base of neighborhoods of zero.
 Then $P$ can be considered as a left
$U_{\kap_0}\comp(\g)$\contramodule{} and $\bigwedge^{\infty/2}_E(\g)$
is a discrete right $\oCl(\g)^\rop$\module, so the space of
homomorphisms $\Hom_k(\bigwedge^{\infty/2}_E(\g),P)$ is a left
contramodule over $U_{\kap_0}\comp(\g)\ot^!\oCl(\g)^\rop$.
 The action of the element $\dd\in U_{\kap_0}\comp(\g)\ot^!
\oCl^{-1}(\g)^\rop$ defines a differential $d^{\infty/2}$ of degree~$1$
on the graded vector space $C^{\infty/2+\bu}_E(\g,P)$ with
the components $C^{\infty/2+i}_E(\g,P) =
\Hom_k(\bigwedge_E^{\infty/2+i}(\g),P)$.
 One has $(d^{\infty/2})^2=0$, since $\dd^2=0$.
 Without the countability assumption, the element $\dd\in\g\til\ot^!
\oCl^{-1}(\g)^\rop$ still acts on the graded vector space
$C^{\infty/2+\bu}_E(\g,P)$ by an operator $d^{\infty/2}$ of degree~$1$.
 By the result of~\ref{noncountable-contra-identity}, the identity
$\dd^2=0$ in $U_{\kap_0}\comp(\g)\ot^!\oCl(\g)^\rop$ implies
the equation $(d^{\infty/2})^2=0$; so $C^{\infty/2+\bu}_E(\g,P)$
becomes a complex.
 This complex is called the \emph{semi-infinite cohomological complex}
and its cohomology is called the \emph{semi-infinite cohomology} of
a Tate Lie algebra $\g$ with coefficients in
a $\g\til$\contramodule{}~$P$.

\subsection{Comparison theorem}

\subsubsection{}
 The correspondence $\h\mpsto \L=\h\dual$ provides an anti-equivalence
between the categories of compact Lie algebras and Lie coalgebras.
 The correspondence between the action maps $\h\times M\rarrow M$
and the coaction maps $M\rarrow\h\dual\ot_k M$ defines an equivalence
between the categories of discrete $\h$\module s and $\L$\comodule s.

 A Lie coalgebra $\L$ is called \emph{conilpotent} if it is a filtered
inductive limit of finite-dimensional Lie coalgebras dual to 
finite-dimensional nilpotent Lie algebras; in other words, $\L$ is 
conilpotent if the dual compact Lie algebra $\h$ is \emph{pronilpotent}.
 A comodule $\M$ over a Lie coalgebra $\L$ is called \emph{conilpotent}
if it is an inductive limit of finite-dimensional comodules which can
be represented as iterated extensions of trivial comodules (that is
comodules with a zero coaction map); analogously one defines
\emph{nilpotent} discrete modules over topological Lie algebras.
 A coassociative coalgebra $\C$ endowed with a coaugmentation map
$k\rarrow\C$ is called \emph{conilpotent} if for every element $c$
of the coalgebra without counit $\C/\im k$ there exists a positive
integer~$i$ such that the iterated comultiplication map
$\C/\im k\rarrow (\C/\im k)^{\ot i}$ annihilates~$c$.

 For a conilpotent Lie coalgebra~$\L$, the \emph{conilpotent
coenveloping coalgebra} $\C(\L)$ is constructed as follows.
 Consider the category of finite-dimensional conilpotent
$\L$\comodule s together with the forgetful functor from it to
the category of finite-dimensional vector spaces; by a result
of~\cite{DM}, this category is equivalent (actually, isomorphic) to
the category of finite-dimensional left comodules over a certain
uniquely defined coalgebra $\C(\L)$ together with the forgetful functor
from this category to the category of  finite-dimensional vector spaces.
 Clearly, the category of (arbitrary) left $\C(\L)$\comodule s is
isomorphic to the category of conilpotent $\L$\comodule s.
 The trivial $\L$\comodule{} $k$ defines a coaugmentation
$k\rarrow\C(\L)$; since this is the only irreducible left
$\C(\L)$\comodule, the coalgebra $\C(\L)$ is conilpotent.

 The coalgebra $\C(\L)$ is the universal final object in the category
of conilpotent coalgebras $\C$ endowed with a Lie coalgebra morphism
$\C\rarrow\L$ such that the composition $k\rarrow\C\rarrow\L$ vanishes.
 Indeed, there is a morphism of Lie coalgebras $\C(\L)\rarrow\L$,
since there is a natural $\L$\comodule{} structure on every left
$\C(\L)$\comodule, and in particular, on the left comodule $\C(\L)$.
 Conversely, a morphism $\C\rarrow\L$ with the above properties defines
a functor assigning to a left $\C$\comodule{} $\M$ a conilpotent
$\L$\comodule{} structure on the same vector space $\M$, hence a left
$\C(\L)$\comodule{} structure on $\M$; this induces a coalgebra
morphism $\C\rarrow\C(\L)$.
 Since the category of finite-dimensional conilpotent $\L$\comodule s
is a tensor category with duality, the coalgebra $\C(\L)$ acquires
a Hopf algebra structure.

 Let $\h$ be the compact Lie algebra dual to~$\L$; then the pairing
$\phi\:\C(\L)\times U(\h)\rarrow k$ is nondegenerate in~$\C$, since
the morphism $\C(\L)\rarrow\L$ factorizes through the quotient
coalgebra of $\C(\L)$ by the kernel of~$\phi$, so a nonzero kernel
would be contradict the universality property.

 Let $\M$ be a conilpotent $\C$\comodule; set $\C=\C(\L)$ and
$\C_+=\C/\im k$.
 Then the natural surjective morphism from the reduced cobar complex
$\M\rarrow\C_+\ot_k\M\rarrow \C_+\ot_k\C_+\ot_k\M\rarrow\dsb$
computing $\Cotor^\C(k,\M)\simeq\Ext_\C(k,\M)$ onto the cohomological
Chevalley complex $\M\rarrow\L\ot_k\M\rarrow\bigwedge^2_k\L\ot_k\M
\rarrow\dsb$ is a quasi-isomorphism.
 It suffices to check this for a finite-dimensional Lie coalgebra
$\L$ dual to a finite-dimensional nilpotent Lie algebra~$\h$;
essentially, one has show that the fully faithful functor from
the category of nilpotent $\h$\module s to the category of arbitrary
$\h$\module s induces isomorphisms on the $\Ext$ spaces.
 This well-known fact can be proven by induction on the dimension
of~$\h$ using the Serre--Hochschild spectral sequences for both
types of cohomology under consideration.
 The key step is to check that for a Lie subcoalgebra $\E\subset\L$
the $\C(\L/\E)$\comodule{} $\C(\L)$ is injective and
the $\E$\comodule{} $\C(\E)$ is the comodule of $\L/\E$\invariant s
in the $\L$\comodule{} $\C(\L)$; it suffices to consider the case
when $\L/\E$ is one-dimensional.

\subsubsection{}  \label{tate-comparison-theorem-proved}
 Let $\g$ be a Tate Lie algebra and $\h\subset\g$ be a compact open Lie
subalgebra.
 Assume that $\h$ is pronilpotent and the discrete $\h$\module{} $\g/\h$
is nilpotent.
 Then the conilpotent coalgebra $\C=\C(\h\dual)$ coacts continuously
in~$\g$, making $(\g,\C)$ a Tate Harish-Chandra pair.
 Let $k\rarrow\g'\rarrow\g$ be a central extension of Tate Lie algebras
endowed with a splitting over~$\h$; then there are a Tate Harish-Chandra
pair $(\g',\C)$ and a central extension of Tate Harish-Chandra pairs
$\kap\:(\g',\C)\rarrow(\g,\C)$ with the kernel~$k$.
 Denote by $\kap_0\:(\g\til,\C)\rarrow(\g,\C)$ the canonical central
extension.
 Set $\S^l_\kap=\S^l_\kap(\g,\C)\simeq\S^r_{\kap+\kap_0}(\g,\C)=
\S^r_{\kap+\kap_0}$.

\begin{thm} \textup{(a)}
 Let\/ $\bN^\bu$ be a complex of right\/
$\S^r_{\kap+\kap_0}$\semimodule s and\/ $\bM^\bu$ be a complex of
left\/ $\S^l_\kap$\semimodule s; in other words, $\bN^\bu$ is a complex
of Harish-Chandra modules with the central charge $-\kap-\kap_0$ and\/
$\bM^\bu$ is a complex of Harish-Chandra modules with the central
charge~$\kap$ over\/~$(\g,\C)$.
 Then the total complex of the semi-infinite homological bicomplex
$C_{\infty/2+\bu}^\h(\g\;\bN^\bu\ot_k\bM^\bu)$ constructed by taking
infinite direct sums along the diagonals represents the object\/
$\SemiTor^{\S^l_\kap}(\bN^\bu,\bM^\bu)$ in the derived category of
$k$\+vector spaces.
 Here the tensor product\/ $\bN^\bu\ot_k\bM^\bu$ is a complex of
Harish-Chandra modules with the central charge~$-\kap_0$. \par
\textup{(b)}
 Let\/ $\bM^\bu$ be a complex of left\/ $\S^l_{\kap}$\semimodule s
and\/ $\bP^\bu$ be a complex of left\/
$\S^r_{\kap+\kap_0}$\semicontramodule s; in other words, $\bM^\bu$ is
a complex of Harish-Chandra modules with the central charge~$\kap$
and\/ $\bP^\bu$ is a complex of Harish-Chandra contramodules with
the central charge $\kap+\kap_0$ over\/~$(\g,\C)$.
 Then the total complex of the semi-infinite cohomological bicomplex
$C_\h^{\infty/2+\bu}(\g\;\Hom_k(\bM^\bu,\bP^\bu))$ constructed by
taking infinite products along the diagonals represents the object\/
$\SemiExt_{\S^l_\kap}(\bM^\bu,\bP^\bu)$ in the derived category of
$k$\+vector spaces.
 Here $\Hom_k(\bM^\bu,\bP^\bu)$ is a complex of Harish-Chandra
contramodules with the central charge~$\kap_0$.
\end{thm}

\begin{proof}
 Part~(a): set $\S^l_{-\kap_0}\simeq\S=\S^r_0$.
 Consider the semi-infinite homological complex
$C_{\infty/2+\bu}^\h(\g,\S)$ of the $\g\til$\module{} $\S$ with
the discrete $\g\til$\module{} structure originating from the left
$\S^l_{-\kap_0}$\semimodule{} structure.
 The complex $C_{\infty/2+\bu}^\h(\g,\S)$ is a complex of right
$\S^r_0$\semimodule s.
 Let us check that it is a semiflat complex naturally isomorphic to
the semimodule $k$ in the semiderived category of right
$\S^r_0$\semimodule s.

 Let $F_i\S$ denote the increasing filtration of the semialgebra $\S$
introduced in~\ref{semialgebra-isomorphism-characterized}.
 Set $F_i\bigl(\bigwedge_\h^{\infty/2}(\g)\bigr)=\bigwedge^i_k(\g)
\wedge\bigwedge_\h^{\infty/2}(\h)$.
 Denote by $F$ the induced filtration of the tensor product
$C_{\infty/2+\bu}^\h(\g,\S)=\bigwedge_\h^{\infty/2}(\g)\ot_k\S$; this is
an increasing filtration of the complex of right $\S^r_0$\semimodule s
$C_{\infty/2+\bu}^\h(\g,\S)$ by subcomplexes of right $\C$\comodule s.
 The complex $\gr_FC_{\infty+\bu}^\h(\g,\S)$ can be identified with
the total complex of the cohomological Chevalley bicomplex
$$
 \bigwedge\nolimits_k(\h\dual)\ot_k\bigwedge
 \nolimits_k(\g/\h)\ot_k\Sym_k(\g/\h)\ot_k\C
$$
of the complex of $\h\dual$\comodule s $\bigwedge_k(\g/\h)\ot_k
\Sym_k(\g/\h)\ot_k\C$ obtained as the tensor product of the Koszul
complex $\bigwedge_k(\g/\h)\ot_k\Sym_k(\g/\h)$ with the coaction
of $\h\dual$ induced by the coaction in~$\g/\h$ and the left
$\C$\comodule{} $\C$ with the induced $\h\dual$\comodule{} structure.
 It follows that the cone of the injection
$F_0C_{\infty/2+\bu}^\h(\g,\S)\rarrow C_{\infty/2+\bu}^\h(\g,\S)$
is a coacyclic complex of right $\C$\comodule s.

 The complex $F_0C_{\infty/2+\bu}^\h(\g,\S)$ is naturally isomorphic to
the cohomological Chevalley complex $\bigwedge_k(\h\dual)\ot_k\C$;
it is a complex of right $\C$\comodule s bounded from below and
endowed with a quasi-isomorphism of complexes of right
$\C$\comodule s $k\rarrow F_0C_{\infty/2+\bu}^\h(\g,\S)$.
 Let us check that the right $\S^r_0$\semimodule{} structure on
$H_0C_{\infty/2+\bu}^\h(\g,\S)\simeq k$ corresponds to the trivial
$\g$\module{} structure.
 The unit element of this homology group can be represented by
the cycle $\lambda\ot 1\in C_{\infty/2+0}^\h(\g,\S)$, where $\lambda$
denotes the unit element of $k\simeq\bigwedge_\h^{\infty/2+0}(\h)
\subset\bigwedge_\h^{\infty/2+0}(\g)$ and $1\in\C\subset\S$ is
the unit (coaugmentation) element of~$\C$.
 Then for any $z\in\g$ one has $(\lambda\ot 1)z = \lambda\ot
(1\ot_{U(\h)}z) = \lambda \ot (\tilde{b}(\overline{z}_{(0)})
\ot_{U(\h)} s(\overline{z}_{(-1)})) = d_{\infty/2}((\overline{z}_{(0)}
\wedge\lambda)\ot s(\overline{z}_{(-1)}))$, where $\overline{z}$
denotes the image of $z$ in $\g/\h$ and
$\tilde{b}\:\g/\h\rarrow\g\til$ is the section corresponding
to any section $b\:\g/\h\rarrow\g$.
 The second equation holds, since the elements $1\ot_{U(\h)}z$
and $\tilde{b}(\overline{z}_{(0)})\ot_{U(\h)} s(\overline{z}_{(-1)})$
of $F_1\S$ have the same images in $F_1\S/F_0\S$ and are both
annihilated by the left action of $\h$ and the map
$\delta^l_{\tilde{b}}=\delta^r_b$.
 To check the third equation, one can use the supercommutation
relation $[\dd,e(y)^\rop]=l(y)$ for $y\in\g$.

 Now let $\tau_{\ge0}C_{\infty/2+\bu}(\g,\S)$ denote the quotient
complex of canonical truncation of the complex of right
$\S^r_0$\semimodule s $C_{\infty/2+\bu}^\h(\g,\S)$ concentrated in
the nonnegative cohomological (nonpositive homological) degrees;
then there are natural morphisms of complexes of right
$\S^r_0$\semimodule s
$$
 k\lrarrow\tau_{\ge0}C_{\infty/2+\bu}^\h(\g,\S)\llarrow
 C_{\infty/2+\bu}^\h(\g,\S)
$$
with $\C$\coacyclic{} cones.
 Indeed, recall that any acyclic complex bounded from below is
coacyclic.
 The embedding $F_0C_{\infty/2+\bu}^\h(\g,\S)\rarrow
C_{\infty/2+\bu}^\h(\g,\S)$ has a $\C$\coacyclic{} cone, as
has the composition $F_0C_{\infty/2+\bu}^\h(\g,\S)\rarrow
C_{\infty/2+\bu}^\h(\g,\S)\rarrow\tau_{\ge0}C_{\infty/2+\bu}^\h
(\g,\S)$, so the cone of the map $C_{\infty/2+\bu}^\h(\g,\S)
\rarrow\tau_{\ge0}C_{\infty/2+\bu}^\h(\g,\S)$ is also $\C$\coacyclic.

 For any complex of left $\S^l_{-\kap_0}$\semimodule s $\bK^\bu$,
the semitensor product $C_{\infty/2+\bu}^\h(\g,\S)\allowbreak
\os_\S\bK^\bu$ is naturally isomorphic to the total complex of
the bicomplex $C_{\infty/2+\bu}^\h(\g,\bK^\bu)$, constructed by
taking infinite direct sums along the diagonals.
 Consider the increasing filtration of the total complex of
$C_{\infty/2+\bu}^\h(\g,\bK^\bu)=\bigwedge^{\infty/2}_\h(\g)\ot_k
\bK^\bu$ induced by the above filtration $F$ of
$\bigwedge^{\infty/2}_\h(\g)$.
 The associated graded quotient complex of this filtration can
be identified with the total complex of the cohomological Chevalley
bicomplex $\bigwedge_k(\h\dual)\ot_k\bigwedge_k(\g/\h)\ot_k\bK^\bu$
of the tensor product of the graded $\h\dual$\comodule{}
$\bigwedge_k(\g/\h)$ and the complex $\bK^\bu$ with
the $\h\dual$\comodule{} structure induced by the left
$\C$\comodule{} structure.
 It follows that the complex $C_{\infty/2+\bu}^\h(\g,\S)\os_\S\bK^\bu$
is acyclic whenever a complex of left $\S^l_{-\kap_0}$\semimodule s
$\bK^\bu$ is $\C$\coacyclic, so the complex of right 
$\S^r_0$\semimodule s $C_{\infty/2+\bu}^\h(\g,\S)$ is semiflat.

 The tensor product $C_{\infty/2+\bu}^\h(\g,\S)\ot_k\bN^\bu$ of
the complex of Harish-Chandra modules $C_{\infty/2+\bu}^\h(\g,\S)$
with central charge~$0$ and the complex of Harish-Chandra modules
$\bN^\bu$ with central charge $-\kap-\kap_0$ is a complex of
Harish-Chandra modules with the central charge $-\kap-\kap_0$.
 This complex of right $\S^r_{\kap+\kap_0}$\semimodule s is semiflat
and naturally isomorphic to $\bN^\bu$ in the semiderived category
of right $\S^r_{\kap+\kap_0}$\semimodule s.
 The latter is clear, and to check the former, notice the isomorphisms
of Lemma~\ref{tate-semitensor-product-semiinvariants}(a)
\begin{multline*}
 (C_{\infty/2+\bu}^\h(\g,\S)\ot_k\bN^\bu)\os_{\S^l_\kap}\bL^\bu\simeq
 (C_{\infty/2+\bu}^\h(\g,\S)\ot_k\bN^\bu\ot_k\bL^\bu)_{\g,\h} \\
 \simeq C_{\infty/2+\bu}^\h(\g,\S)\os_{\S^r_0}(\bN^\bu\ot_k\bL^\bu)
\end{multline*}
for any complex of left $\S^l_\kap$\semimodule s $\bL^\bu$.
 Now the object $\SemiTor^{\S^l_\kap}(\bN^\bu,\bM^\bu)$ is represented
by the complex 
$$(C_{\infty/2+\bu}^\h(\g,\S)\ot_k\bN^\bu)\os_{\S^l_\kap}
\bM^\bu\simeq C_{\infty/2+\bu}^\h(\g,\S)\os_{\S^r_0}(\bN^\bu\ot_k
\bM^\bu)\simeq C_{\infty/2+\bu}^\h(\g\;\bN^\bu\ot_k\bM^\bu)$$
in the derived category of $k$\+vector spaces.

 Another way to identify $\SemiTor^{\S^l_\kap}(\bN^\bu,\bM^\bu)$
with $C_{\infty/2+\bu}^\h(\g\;\bN^\bu\ot_k\bM^\bu)$ is to consider
the semiflat complex of left $\S^l_\kap$\semimodule s
$C_{\infty/2+\bu}^\h(\g,\S) \ot_k\bM^\bu$ naturally isomorphic to
$\bM^\bu$ in the semiderived category of left $\S^l_\kap$\semimodule s.
 To check that these two identifications coincide, represent the images
of $\bN^\bu$ and $\bM^\bu$ in the semiderived categories of semimodules
by arbitrary semiflat complexes.
 The proof of part~(b) is completely analogous.
\end{proof}

\begin{qst}
 Can one obtain the semi-infinite homology of arbitrary discrete modules
over a Tate Lie algebra with a fixed compact open subalgebra (rather
than only Harish-Chandra modules under the nilpotency conditions) as
a kind of double-sided derived functor of the functor of semiinvariants
on an appropriate exotic derived category of discrete modules?
 Notice that the cohomology of the Chevalley complex $\M\rarrow
\L\ot_k\M\rarrow\bigwedge^2_k\L\ot_k\M\rarrow\dsb$ for a comodule $\M$
over a Lie coalgebra $\L$ is indeed the right derived functor of
the functor of $\L$\invariant s on the abelian category of
$\L$\comodule s, since the category of $\L$\comodule s has enough
injectives and the cohomology of the Chevalley complex is
an effaceable cohomological functor.
 The former holds since the category of discrete modules over a compact
Lie algebra~$\h=\L\dual$ has exact functors of filtered inductive
limits preserved by the forgetful functor to the category of $k$\+vector
spaces, and the discrete $\h$\module s $U(\h)\ot_{U(\aa)}k$ induced from
trivial modules over open subalgebras $\aa\subset\h$ form a set of
generators, so the forgetful functor even has a right adjoint.
 To check the latter, one can represent cocycles in the cohomological
Chevalley complex by discrete $\h$\module{} morphisms into $\M$
from the relative homological complexes $\dsb\rarrow\bigwedge^2_k
(\h/\aa)\ot_{U(\aa)}U{\h}\rarrow\h/\aa\ot_{U(\aa)}U(\h)\rarrow U(\h)$,
which are quotient complexes of the Chevalley homological complex of
the $\h$\module{} $U(\h)$ and finite discrete $\h$\module{}
resolutions of the trivial $\h$\module~$k$.
 Furthermore, the semiderived category of discrete modules over a Tate
Lie algebra~$\g$, defined as the quotient category of the homotopy
category of discrete $\g$\module s by the thick subcategory of complexes
coacyclic as complexes of discrete $\h$\module s, does not depend on
the choice of an open compact subalgebra $\h\subset\g$.
 This can be demonstrated by considering the tensor product over~$k$
of the above relative homological complex with a complex of discrete
$\h$\module s coacyclic over~$\aa$.
 Notice that in the above proof we have essentially shown that
the semi-infinite homology is a functor on the semiderived category
of discrete $\g$\module s.
\end{qst}

\subsubsection{}
 We keep the assumptions and notation
of~\ref{tate-comparison-theorem-proved}, and also use
the notation of Corollary~\ref{semialgebra-isomorphism-characterized}.
 The following result makes use of the semimodule-semicontramodule
correspondence in order to express the semi-infinite homology and
cohomology in terms of compositions of one-sided derived functors.

\begin{cor}
 \textup{(a)} Let\/ $\bM^\bu$ be a complex of Harish-Chandra modules
with the central charge~$\kap$ and\/ $\bP^\bu$ be a complex of
Harish-Chandra contramodules with the central charge~$\kap+\kap_0$
over $(\g,\C)$.
 Then the semi-infinite cohomological complex $C^{\infty/2+\bu}_\h
(\g,\Hom_k(\bM^\bu,\bP^\bu))$ represents the object\/
$\Ext_{\S^l_\kap}(\bM^\bu,\boL\Phi_{\S^r_{\kap+\kap_0}}(\bP^\bu))\simeq
\Ext^{\S^r_{\kap+\kap_0}}(\boR\Psi_{\S^l_\kap}(\bM^\bu),\bP^\bu)$ in
the derived category of $k$\+vector spaces. \par
 \textup{(b)} Let\/ $\bM^\bu$ be a complex of Harish-Chandra modules
with the central charge~$\kap$ and\/ $\bN^\bu$ be a complex of
Harish-Chandra modules with the central charge~$-\kap-\kap_0$ over
$(\g,\C)$.
 Then the semi-infinite homological complex $C_{\infty/2+\bu}^\h
(\g\;\bN^\bu\ot_k\bM^\bu)$ represents the object\/
$\Ctrtor^{\S^r_{\kap+\kap_0}}(\bN^\bu,\boR\Psi_{\S^l_\kap}(\bM^\bu))
\simeq\Ctrtor^{\S^r_{-\kap}}(\bM^\bu,\boR\Psi_{\S^l_{-\kap-\kap_0}}
(\bN^\bu))$ in the derived category of $k$\+vector spaces.
\end{cor}

\begin{proof}
 This follows from Theorem~\ref{tate-comparison-theorem-proved} and
Corollary~\ref{semiext-and-ext}.
\end{proof}

\begin{rmk}
 Set $\S^l_0=\S\simeq\S^r_{\kap_0}$ and consider the complex of left
$\S$\semimodule s $\bR^\bu=\S\ot_k\bigwedge^\h_{\infty/2+\bu}(\g)$.
 This is a semiprojective complex of semiprojective left
$\S$\semimodule s isomorphic to the trivial $\S$\semimodule~$k$
in $\sD^\si(\S\simodl)$.
 Assume that the pronilpotent Lie algebra $\h$ is
infinite-dimensional (cf.~\ref{prelim-co-contra-acyclicity}).
 Then the complex of left $\S$\semicontramodule s $\Psi_\S(\bR^\bu)$
is acyclic.
 Indeed, it suffices to check that the complex of left
$\C$\contramodule s obtained by applying the functor $\Psi_\C$
to the cohomological Chevalley complex $\C\ot_k\bigwedge_k(\h\dual)$
is acyclic; one can reduce this problem to the case of an abelian Lie
algebra by considering the decreasing filtration $\h\supset[\h,\h]
\supset[\h,[\h,\h]]\supset\dsb$ on $\h$ and the induced increasing
filtrations on $\h\dual$ and~$\C$.
 The complex $\Psi_\S(\bR^\bu)$ is also a projective complex of
projective left $\S$\semicontramodule s (see
Remark~\ref{semi-ctrtor-definition} and
subsection~\ref{semi-model-struct}); it can be thought of as
the ``projective $\S$\semicontramodule{} resolution of
a (nonexistent) one-dimensional left $\S$\semicontramodule{}
placed in the degree~$+\infty$''.
 For any complex of right $\S$\semimodule s $\bN^\bu$,
the contratensor product complex $\bN^\bu\Ocn_\S\Psi_\S(\bR^\bu)$
computes the semi-infinite homology of $\g$ with coefficients
in $\bN^\bu$.
 For any complex of left $\S$\semicontramodule s $\bP^\bu$,
the complex of semicontramodule homomorphisms
$\Hom^\S(\Psi_\S(\bR^\bu),\bP^\bu)$ computes the semi-infinite
cohomology of $\g$ with coefficients in $\bP^\bu$.
 (Cf.~\cite[subsection~3.11.4]{Vor}.)
\end{rmk}

\Section{Groups with Open Profinite Subgroups}

 To a locally compact totally disconnected topological group $G$
and a commutative ring $k$ one associates a family of left and right
semiprojective Morita equivalent semialgebras $\S=\S_k(G,H)$ numbered
by open profinite subgroups $H\subset G$.
 As explained in~\ref{semi-morita-morphisms}, Morita equivalences of
semialgebras do not have to preserve the semiderived categories or 
the derived functors $\SemiTor$ and $\SemiExt$, and this is indeed
the case here: $\SemiTor^\S$ and $\SemiExt_\S$ depend very essentially
on~$H$.
 For a complex of smooth $G$\module s $N^\bu$ over~$k$ and a complex
of $k$\+flat smooth $G$\module s $M^\bu$ over~$k$, we show that
$\SemiTor^{\S_k(G,H)}(N^\bu,M^\bu)$ only depends on the complex of
smooth $G$\module s $N^\bu\ot_kM^\bu$, and analogously for
$\SemiExt_{\S_k(G,H)}$.

 When $k$ is a field of zero characteristic, one can climb one step
higher and assign to a ``good enough'' group object $\boG$ in
the category of ind\+pro\+topological spaces with a subgroup object
$\boH$ belonging to the category of pro\+topological groups a coring
object in the tensor category of representations of $\boH\times\boH$
in pro-vector spaces over~$k$.
 So there is a functor of cotensor product~\cite{GK2}
on certain categories of representations of central extensions
of~$\boG$; it has a double-sided derived functor $\ProCotor$.
 
\subsection{Morita equivalent semialgebras}

\subsubsection{}
 In the sequel, all topological spaces and topological groups are
presumed to be locally compact and totally disconnected.

 For a topological space $X$ and an abelian group~$A$, denote by
$A(X)$ the abelian group of locally constant compactly supported
$A$\+valued functions on~$X$.
 For any proper map of topological spaces $X\rarrow Y$, the pull-back
map $A(Y)\rarrow A(X)$ is defined.
 For any \'etale map (local homeomorphism) of topological spaces
$X\rarrow Y$, the push-forward map $A(X)\rarrow A(Y)$ is defined.

 For any topological spaces $X$ and $Y$ and an abelian group $A$,
there is a natural isomorphism $A(X\times Y)\simeq A(X)(Y)$.
 For any topological space $X$, a commutative ring~$k$, and 
a $k$\module{} $A$, there is a natural isomorphism $A(X)\simeq
A\ot_k k(X)$.
 
 For a topological space $X$ and an abelian group $A$, denote by
$A[[X]]$ the abelian group of finitely-additive compactly supported
$A$\+valued measures defined on the open-closed subsets of~$X$.
 For any map of topological spaces $X\rarrow Y$, the push-forward
map $A[[X]]\rarrow A[[Y]]$ is defined.

 For any compact topological space $X$, a commutative ring~$k$,
and a $k$\module{} $A$, there is a natural isomorphism
$A[[X]]\simeq\Hom_k(k(X),A)$.

\subsubsection{}
 Let $H$ be a profinite group and $k$ be a commutative ring.
 Then the module of locally constant functions $k(H)$ has a natural
structure of coring over~$k$ where the left and right actions of~$k$
coincide.
 This coring structure, which we denote by $\C=\C_k(H)$, is defined
as follows.
 The counit map $k(H)\rarrow k$ is the evaluation at the unit element
$e\in H$.
 The comultiplication map $k(H)\rarrow k(H)\ot_k k(H)$ is provided
by the pull-back map $k(H)\rarrow k(H\times H)$ induced by
the multiplication map $H\times H\rarrow H$ together with
the identification $k(H\times H)\simeq k(H)\ot_kk(H)$.

 Let $G$ be a topological group and $H\subset G$ be an open profinite
subgroup.
 Then the module of locally constant compactly supported functions
$k(G)$ has a natural structure of semialgebra over the coring
$\C_k(H)$.
 This semialgebra structure, which we denote by $\S=\S_k(G,H)$,
is defined as follows.
 The bicoaction map $k(G)\rarrow k(H)\ot_kk(G)\ot_kk(H)\simeq
k(H\times G\times H)$ is the pull-back map induced by
the multiplication map $H\times G\times H\rarrow G$.
 The semiunit map $k(H)\rarrow k(G)$ is the push-forward map induced
by the injection $H\rarrow G$.

 Denote by $G\times_H G$ the quotient space of the Carthesian square
$G\times G$ by the equivalence relation $(g'h,g'')\sim(g',hg'')$
for $g'$, $g''\in G$ and $h\in H$.
 The pull-back map $k(G\times_H G)\rarrow k(G\times G)$ induced by
the natural surjection $G\times G\rarrow G\times_H G$ identifies
$k(G\times_HG)$ with the cotensor product $k(G)\oc_\C k(G)\subset
k(G)\ot_kk(G)\simeq k(G\times G)$.
 The semimultiplication map $k(G)\oc_\C k(G)\rarrow k(G)$ is
the push-forward map induced by the multiplication map $G\times_HG
\rarrow G$.

 The involutions $k(H)\rarrow k(H)$ and $k(G)\rarrow k(G)$ induced
by the inverse element maps $H\rarrow H$ and $G\rarrow G$ provide
the isomorphism of semialgebras $\S_k(G,H)^\rop\simeq\S_k(G,H)$
compatible with the isomorphism of corings $\C_k(H)^\rop\simeq\C_k(H)$
over~$k$ (see~\ref{hopf-morita-anti-involution}).

 Now let $H_1$, $H_2\subset G$ be two open profinite subgroups of~$G$.
 Then the $k$\module{} $k(G)=\S_k(G,H_1)$ has a natural left
$\S_k(G,H_1)$\semimodule{} structure and at the same time
$k(G)=S_k(G,H_2)$ has a natural right $\S_k(G,H_2)$\semimodule{}
structure.
 Obviously, these two semimodule structures commute; we denote this
bisemimodule structure on $k(G)$ by $\S_k(G,H_1,H_2)$.
 For any three open profinite subgroups $H_1$, $H_2$, $H_3\subset G$,
there is a natural isomorphism $\S_k(G,H_1,H_2)\os_{S_k(G,H_2)}
\S_k(G,H_2,H_3)\simeq\S_k(G,H_2)\os_{S_k(G,H_2)}\S_k(G,H_2)\simeq
\S_k(G,H_2)\simeq\S_k(G,H_1,H_3)$; this is an isomorphism of
$\S_k(G,H_1)$\+$\S_k(G,H_3)$\bisemimodule s.
 
 One can check that $k(H)$ is a projective $k$\module.
 Clearly, $\S_k(G,H)$ is a coprojective left and right
$\C_k(H)$\comodule.
 So the pair $(\S_k(G,H_1,H_2),\S_k(G,H_2,H_1))$ is a left and right
semiprojective Morita equivalence between the semialgebras
$\S_k(G,H_1)$ and $\S_k(G,H_2)$.

\subsubsection{}
 The semialgebra $\S_k(G,H)$ can be also obtained by the construction
of~\ref{semialgebra-constructed}.

 Denote by $k[H]$ and $k[G]$ the group $k$\+algebras of the groups
$H$ and $G$ considered as groups without any topology.
 There is a pairing $\phi\:\C_k(H)\ot_kk[H]\rarrow k$ satisfying
the conditions of~\ref{ring-coring-pairing}, given by the formula
$(c,h)\mpsto c(h^{-1})$ for any $c\in k(H)$ and $h\in H$
 The induced functor $\Delta_\phi\:\comodr\C_k(H)\rarrow\modr k[H]$
is fully faithful; its image is described as follows.

 A module $M$ over a topological group~$G$ is called \emph{smooth}
(discrete), if the action map $G\times M \rarrow M$ is continuous
with respect to the discrete topology of $M$; equivalently, $M$
is smooth if the stabilizer of every its element is an open
subgroup in~$G$.
 The functor $\Delta_\phi$ identifies the category of right
$\C_k(H)$\comodule s with the category of smooth $H$\module s
over~$k$.

 The tensor product $\C_k(H)\ot_{k[H]}k[G]$ is a smooth $G$\module{}
with respect to the action of $G$ by right multiplications, so it
becomes a semialgebra over $\C_k(H)$.
 This semialgebra can be identified with $\S_k(G,H)$ by the formula
by the formula $(c\ot g)\mpsto (g'\mpsto c(g'g^{-1}))$, where 
a locally constant function $c\:H\rarrow k$ is presumed to be
extended to~$G$ by zero.
 By the result of~\ref{semi-mod-contra-described}, the category of
right $\S_k(G,H)$\semimodule s is isomorphic to the category of smooth
$G$\module s over~$k$.

 One also obtains the following description of the category of
left $\S_k(G,H)$\semicontramodule s.
 For a topological group~$G$ and a commutative ring~$k$,
a \emph{$G$\+contramodule{} over~$k$} is a $k$\module{} $P$ endowed
with a $k$\+linear map $P[[G]]\rarrow P$ satisfying the following
conditions.
 First, the point measure supported in the unit element $e\in G$
and corresponding to an element $p\in P$ should map to the element~$p$.
 Second, the composition $P[[G\times G]]\rarrow P[[G]][[G]]\rarrow P$
of the natural map $P[[G\times G]]\rarrow P[[G]][[G]]$ and the iterated
contraaction map $P[[G]][[G]]\rarrow P[[G]]\rarrow P$ should be equal
to the composition $P[[G\times G]]\rarrow P[[G]]\rarrow P$ of
the push-forward map $P[[G\times G]]\rarrow P[[G]]$ induced by
the multiplication map $G\times G\rarrow G$ and the contraaction map
$P[[G]]\rarrow P$.

 The images of the point measures under the contraaction map define
the forgetful functor from the category of $G$\contramodule s over~$k$
to the category of (nontopological) $G$\module s over~$k$.
 The category of left $\S_k(G,H)$\semicontramodule s is isomorphic to
the category of $G$\contramodule s over~$k$.

\subsubsection{} {\hbadness=3000
 For any smooth $G$\module{} $M$ over $k$ and any $k$\module{} $E$
there is a natural $G$\contramodule{} structure on the space of
$k$\+linear maps $\Hom_k(M,E)$.
 The contraaction map $\Hom_k(M,E)[[G]]\rarrow\Hom_k(M,E)$ is
constructed as the projective limit over all open subgroups
$U\subset G$ of the compositions $\Hom_k(M,E)[[G]]\rarrow
\Hom_k(M,E)[G/U]\rarrow\Hom_k(M^U,E)$ of the maps $\Hom_k(M,E)[[G]]
\rarrow\Hom_k(M,E)[G/U]$ induced by the surjections $G\rarrow G/U$
and the maps $\Hom_k(M,E)[G/U]\rarrow\Hom_k(M^U,E)$ induced by 
the action maps $G/U\times M^U\rarrow M$, where $M^U$ denotes
the $k$\+submodule of $U$\invariant s in~$M$ and $G/U$ is the set
of all left cosets of $G$ modulo~$U$. \par}

 More generally, let $G_1$ and $G_2$ be topological groups.
 Then for any smooth $G_1$\module{} $M$ over~$k$ and
$G_2$\contramodule{} $P$ over~$k$ there is a natural
$G_1\times G_2$\contramodule{} structure on $\Hom_k(M,P)$ with
the contraaction map $\Hom_k(M,P)[[G_1\times G_2]]\rarrow\Hom_k(M,P)$
defined as either of the compositions $\Hom_k(M,P)[[G_1\times G_2]]
\rarrow\Hom_k(M,P)[[G_1]][[G_2]]\rarrow\Hom_k(M,P)[[G_2]]\rarrow
\Hom_k(M,P)$ or $\Hom_k(M,P)[[G_1\times G_2]]\rarrow
\Hom_k(M,P)[[G_2]][[G_1]]\rarrow\Hom_k(M,P)[[G_1]]\rarrow
\Hom_k(M,P)$, where the $G_1$\+contraaction map $\Hom_k(M,P)[[G_1]]
\rarrow \Hom_k(M,P)$ is defined above and the $G_2$\+contraaction
map $\Hom_k(M,P)[[G_2]]\rarrow\Hom_k(M,P)$ is constructed as
the composition $\Hom_k(M,P)[[G_2]]\rarrow\Hom_k(M,P[[G_2]])
\rarrow\Hom_k(M,P)$.
 Hence for any smooth $G$\module{} $M$ over~$k$ and any
$G$\contramodule{} $P$ over~$k$ there is a natural $G$\contramodule{}
structure on $\Hom_k(M,P)$ induced by the diagonal map of topological
groups $G\rarrow G\times G$.

\subsection{Semiinvariants and semicontrainvariants}

\subsubsection{}
 Let $G$ be a topological group and $H\subset G$ be an open profinite
subgroup.
 For a smooth $H$\module{} $M$ over~$k$, let $\Ind_H^GM$ denote
the induced $G$\module{} $k[G]\ot_{k[H]}M$.
 For any smooth $G$\module{} $N$ over~$k$ we will construct a pair of
maps $(\Ind_H^GN)^H\birarrow N^H$, where the superindex $H$ denotes
the $k$\+submodule of $H$\invariant s.
 Namely, the first map $(\Ind_H^GN)^H\rarrow N^H$ is obtained by
applying the functor of $H$\invariant s to the map $(\Ind_H^GN)
\rarrow N$ given by the formula $g\ot n\mpsto g(n)$.
 The second map $(\Ind_H^GN)^H\rarrow N^H$ only depends on
the $H$\module{} structure on~$N$.

 To define this second map, identify the induced representation
$\Ind_H^GN$ with the $k$\module{} of all compactly supported functions
$G\rarrow N$ transforming the right action of $H$ in $G$ into
the action of $H$ in~$N$; this identification assigns to a function
$g\mpsto n_g$ the formal linear combination $\sum_{g\in G/H} g\ot n_g$,
where $G/H$ denotes the set of all left cosets of $G$ modulo~$H$.
 An $H$\invariant{} element of $\Ind_H^GN$ is then represented by
a compactly supported function $G\rarrow N$ denoted by $g\mpsto n_g$
and satisfying the equations $n_{gh}=h^{-1}(n_g)$ and $n_{hg}=n_g$ for
$g\in G$, \ $h\in H$.
 The second map $(\Ind_H^GN)^H\rarrow N^H$ sends a function $g\mpsto n_g$
to the element $\sum_{g\in H\backslash G} n_g\in N$, where $H\backslash G$
denotes the set of all right cosets of $G$ modulo~$H$.

 The cokernel $N_{G,H}$ of this pair of maps $(\Ind_H^GN)^H\birarrow N^H$
is called the module of \emph{$(G,H)$\semiinvariant s} of a smooth
$G$\module~$N$.
 The $(G,H)$\semiinvariant s are a mixture of $H$\invariant s and
coinvariants along $G$ relative to~$H$.

\subsubsection{}
 For an $H$\contramodule{} $Q$ over~$k$, let $\Coind_H^G(Q)$ denote
the coinduced $G$\contramodule{} $\Hom_{k[H]}(k[G],Q)\simeq
\Cohom_{\C_k(H)}(\S_k(G,H),Q)$.
 For any $G$\contramodule{} $P$ over~$k$ we will construct a pair 
of maps $P_H\birarrow\Coind_H^G(P)_H$, where the subindex $H$ denotes
the $k$\module{} of $H$\coinvariant s, i.~e., the maximal quotient
$H$\contramodule{} with the trivial contraaction.
 Namely, the first map $P_H\rarrow\Coind_H^G(P)_H$ is obtained by
applying the functor of $H$\coinvariant s to the semicontraaction
map $P\rarrow\Coind_H^G(P)$, which is given by the formula
$p\mpsto (g\mpsto g(p))$.

 The second map $P_H\rarrow\Coind_H^G(P)_H$ only depends on
the $H$\contramodule{} structure on~$P$.
 It is given by the formula $p\mpsto \sum_{g\in G/H}g*p$ for $p\in P$,
where $g*p\:k[G]\rarrow P$ is the $k[H]$\+linear map defined by
the rules $hg\mpsto h(p)$ for $h\in H$ and $g'\mpsto 0$ for all
$g'\in G$ not belonging to the right coset $Hg$.
 Clearly, the infinite sum over $G/H$ converges element-wise on
$k[G]$ for any choice of representatives of the left cosets; its
image in the module of coinvariants does not depend on this choice
and is determined by the image of~$p$ in~$P_H$.

 The kernel $P^{G,H}$ of this pair of maps $P_H\birarrow\Coind_H^G(P)_H$
is called the module of \emph{$(G,H)$\semicontrainvariant s} of
a $G$\contramodule{} $P$.
 The $(G,H)$\semicontrainvariant s are a mixture of $H$\coinvariant s
and invariants along $G$ relative to~$H$.

\subsubsection{}
 Denote by $\chi\:G\rarrow\boQ^*$ the modular character of~$G$,
i.~e., the character with which $G$ acts by left shifts on
the one-dimensional $\boQ$\+vector space of $\boQ_\chi$ of right
invariant $\boQ$\+valued measures defined on the open compact
subsets of~$G$.
 Equivalently, $\chi(g)$ is equal to the ratio of the number of left
cosets contained in the double coset $HgH$ to the number of right
cosets.
 Whenever the commutative ring~$k$ contains $\boQ$, there is a natural
isomorphism $N_{G,H}\simeq (N\ot_\boQ\boQ_\chi)_G$ for any smooth
$G$\module{} $N$ over~$k$, where the subindex~$G$ denotes
the $k$\module{} of $G$\coinvariant s.

 Indeed, the composition $M^H\rarrow M\rarrow M_H$ is an isomorphism
for any smooth $H$\module{} $M$ over~$k$, so in particular there are
isomorphisms $N^H\simeq N_H$ and $(\Ind_H^GN)^H\rarrow(\Ind_H^G)_H$.
 These isomorphisms transform the above pair of maps 
$(\Ind_H^GN)^H\birarrow N^H$ into the pair of maps
$(\Ind_H^GN)_H\birarrow N_H$ given by the formulas
$g\ot n\mpsto g(n)$ and $g\ot n\mpsto \chi^{-1}(g)n$.

\subsubsection{}
 For any $G$\contramodule{} $P$ over~$k$ denote by $P^G$
the $k$\module{} of $G$\invariant s in~$P$ defined as the submodule
of all $p\in P$ such that for any measure $m\in k[[G]]$ the image
of the measure $pm\in P[[G]]$ under the contraaction map 
$P[[G]]\rarrow P$ is equal to the value $m(G)$ of~$m$ at~$G$.

 Assuming that the commutative ring $k$ contains~$\boQ$,
the composition $Q^H\rarrow Q\rarrow Q_H$ is an isomorphism for any
$H$\contramodule{} $Q$ over~$k$, as one can show using the action of
the Haar measure of the profinite group $H$ in the contramodule~$Q$.
 One can also use the Haar measure to check that $Q^H$ is
the maximal subcontramodule of~$Q$ with the trivial contraaction
of~$H$, and it follows that $P^G$ is the maximal subcontramodule
of~$P$ with the trivial contraaction of~$G$, under our assumption.
 Finally, when $k\supset\boQ$ there is a natural isomorphism
$P^{G,H}\simeq\Hom_\boQ(\boQ_\chi,P)^G$.

\subsubsection{}
 Let $N$ be a left $\S_k(G,H)$\semimodule{} and $M$ be a right
$\S_k(G,H)$\semimodule; in other words, $N$ and $M$ are smooth
$G$\module s over~$k$.
 Then there is a natural isomorphism $N\os_{\S_k(G,H)}M\simeq
(N\ot_kM)_{G,H}$, where $N\ot_kM$ is considered as a smooth
$G$\module{} over~$k$.
 Here the semitensor product is well-defined by
Proposition~\ref{cotensor-associative}(f).

 Indeed, there is an obvious isomorphism $N\oc_{\C_k(H)}M\simeq
(N\ot_kM)^H$.
 The $k$\module{} $N\oc_{\C_k(H)}\S_k(G,H)\oc_{\C_k(H)}M$ can be
identified with the module of locally constant compactly
supported functions $f\:G\rarrow N\ot_kM$ satisfying
the equations $f(hg)=hf(g)$ and $f(gh)=f(g)h$ for $g\in G$, \
$h\in H$, where $(g,a)\mpsto ga$ and $(a,g)\mpsto ag$ denote
the actions of $G$ in $N\ot_kM$ induced by the actions
in $N$ and $M$, respectively.
 At the same time, the $k$\module{} $(\Ind_H^G(N\ot_kM))^H$
can be identified with the module of locally constant compactly
supported functions $f'\:G\rarrow N\ot_kM$ satisfying
the equations $f(hg)=f(g)$ and $f(gh)=h^{-1}f(g)h$.
 The formula $f'(g)=g^{-1}f(g)$ defines an isomorphism between
these two $k$\module s transforming the pair of maps whose cokernel
is  $N\os_{\S_k(G,H)}M$ into the pair of maps whose cokernel is
$(N\ot_kM)_{G,H}$.
 
\subsubsection{}
 Let $M$ be a left $\S_k(G,H)$\semimodule{} and $P$ be a left
$\S_k(G,H)$\semicontramodule.
 Then there is a natural isomorphism $\SemiHom_{\S_k(G,H)}(M,P)
\simeq\Hom_k(M,P)^{G,H}$, where $\Hom_k(M,P)$ is considered
as a $G$\contramodule{} over~$k$ and the semihomomorphism module is
well-defined by Proposition~\ref{cohom-associative}(j).

{\hbadness=1500
 Indeed, the quotient modules $\Cohom_{\C_k(H)}(M,P)$ and $\Hom_k(M,P)_H$
of the $k$\module{} $\Hom_k(M,P)$ coincide.
 There are two commuting contraactions of $G$ in $\Hom_k(M,P)$ induced
by the smooth action in $M$ and the contraaction in~$P$; denote 
these contraaction maps by $\pi_M$ and~$\pi_P$, and the corresponding
actions of $G$ in $\Hom_k(M,P)$ by $(g,x)\mpsto g_M(x)$ and 
$(g,x)\mpsto g_P(x)$.
\par}

 The $k$\module{} $\Cohom_{\C_k(H)}(\S_k(G,H)\oc_{\C_k(H)}M\;P)$ can be
identified with the quotient module of the module of all
finitely-additive measures defined on compact open subsets of $G$ and
taking values in $\Hom_k(M,P)$ by the submodule generated by measures
of the form $U\mpsto \pi_M(W\mpsto \mu(U\times W^{-1})) - 
\mu(\{(g,h)\mid gh\in U\})$ and $U\mpsto\pi_P(W\mpsto\nu(U\times W)) - 
\nu(\{(h,g)\mid hg\in U\})$, where $\mu$ and $\nu$ are finitely-additive
measures defined on compact open subsets of $G\times H$ and $H\times G$,
respectively, and taking values in $\Hom_k(M,P)$.
 Here $U$ denotes a compact open subset of $G$ and $W$ is an open-closed
subset of $H$; by $W^{-1}$ we denote the (pre)image of $W$ under
the inverse element map.

 At the same time, the $k$\module{} $\Coind_H^G(\Hom_k(M,P))_H$ can be
identified with the quotient module of the same module of measures
on~$G$ by the submodule generated by measures of the form
$U\mpsto \mu'(U\times H) - \mu'(\{(g,h)\mid gh\in U\})$ and
$U\mpsto \pi_P(W_2\mpsto \pi_M(W_1\mpsto\nu'((W_1\cap W_2)\times U)))
- \nu'(\{(h,g)\mid hg\in U\})$, where $\mu'$ and $\nu'$ are 
finitely-additive measures defined on compact open subsets of
$G\times H$ and $H\times G$, respectively, and taking values
in $\Hom_k(M,P)$.
 Here $W_1$ and $W_2$ denote open-closed subsets of~$H$.
 The formulas $m'(U)=\pi_M(V\mpsto m(U\cap V^{-1}))$ and $m(U)=
\pi_M(V\mpsto m'(U\cap V))$, where $V$ denotes an open-closed subset
of~$G$, define as isomorphism between these two quotient spaces of
measures.

 The pair of maps $\Cohom_{\C_k(H)}(M,P)\birarrow\Cohom_{\C_k(H)}(\S_k(G,H)
\oc_{\C_k(H)}M\;P)$ whose kernel is $\SemiHom_{\S_k(G,H)}(M,P)$ is given
by the rules $x\mpsto \sum_{g\in G/H}g_M^{-1}(x)\delta_g$ and
$x\mpsto \sum_{g\in H\backslash G}g_P(x)\delta_g$, where $y\delta_g$
denotes the $\Hom_k(M,P)$\+valued measure defined on compact open
subsets of~$G$ that is supported in the point~$g\in G$ and 
corresponds to an element $y\in\Hom_k(M,P)$.
 The pair of maps $\Hom_k(M,P)_H\birarrow\Coind_H^G(\Hom_k(M,P))_H$
whose kernel is $\Hom_k(M,P)^{G,H}$ is given by the rules
$x\mpsto\sum_{g\in G/H}x\delta_g$ and $x\mpsto \sum_{g\in H\backslash G}
g_Pg_M(x)\delta_g$.
 The above isomorphism between two quotient spaces of measures
transforms one of these two pairs of maps into the other.

\subsubsection{}
 Let $H_1\subset H_2$ be two open profinite subgroups of a topological
group $G$ and $N$ be a smooth $G$\module{} over~$k$. 
 Then, at least, in the following two situations there is a natural
isomorphism $N_{G,H_1}\simeq N_{G,H_2}$: when $N$ as a $G$\module{} is
induced from a smooth $H_1$\module{} over~$k$, and when $N$ as
an $H_2$\module{} over~$k$ is coinduced from a module over the trivial
subgroup $\{e\}$ (i.~e., $N$ is a coinduced $\C_k(H_2)$\comodule).

 These isomorphisms are constructed as follows.
 In the first case, one shows that the triple semitensor product
$N\os_{\S_k(G,H_1)}\S_k(G,H_1,H_2)\os_{\S_k(G,H_2)}k$ is associative
in the sense that the conclusion of
Proposition~\ref{semitensor-associative} applies to it.
 In the second case, one shows that the triple semitensor product
$N\os_{\S_k(G,H_2)}\S_k(G,H_2,H_1)\os_{\S_k(G,H_1)}k$ is associative
in the similar sense.
 In both cases, the argument is analogous to that of
Proposition~\ref{semitensor-associative}.

\subsubsection{}
 Let $H_1\subset H_2$ be two open profinite subgroups of~$G$ and
$P$ be a $G$\contramodule{} over~$k$.
 Then, at least, in the following two situations there is a natural
isomorphism $P^{G,H_1}\simeq P^{G,H_2}$: when $P$ is coinduced from
an $H_1$\contramodule{} over~$k$, and when $P$ as an $H_2$\contramodule{}
over~$k$ is induced from a contramodule over the trival subgroup
$\{e\}$ (i.~e., $P$ is an induced $\C_k(H_2)$\contramodule).

\subsection{SemiTor and SemiExt}

\subsubsection{}
 Assume that the ring~$k$ has a finite weak homological dimension.
 Then for any complexes of smooth $G$\module s $N^\bu$ and $M^\bu$
over~$k$ the object $\SemiTor^{\S_k(G,H)}\allowbreak(N^\bu,M^\bu)$ in
the derived category of $k$\module s is defined.
 Furthermore, whenever either $N^\bu$ or $M^\bu$ is a complex
of $k$\+flat smooth $G$\module s over~$k$ there is a natural
isomorphism $\SemiTor^{\S_k(G,H)}(N^\bu,M^\bu)\simeq
\SemiTor^{\S_k(G,H)}(N^\bu\ot_kM^\bu\;k)$ in $\sD(k\modl)$.

 Indeed, assume that $M^\bu$ is a complex of flat $k$\module s.
 If $N^\bu$ is a semiflat complex of right $\S_k(G,H)$\semimodule s,
then so is the tensor product $N^\bu\ot_kM^\bu$, since $(N^\bu\ot_kM^\bu)
\os_{\S_k(G,H)}L^\bu\simeq(N^\bu\ot_kM^\bu\ot_kL^\bu)_{G,H}\simeq
N^\bu\os_{\S_k(G,H)}(M^\bu\ot_kL^\bu)$ for any complex of smooth
$G$\module s $L^\bu$ over~$k$, and the complex of
left $\C_k(H)$\comodule s $M^\bu\ot_kL^\bu$ is coacyclic whenever
the complex of left $\C_k(H)$\comodule s $L^\bu$~is.
 It remains to use the natural isomorphism $N^\bu\os_{\S_k(G,H)}M^\bu
\simeq(N^\bu\ot_kM^\bu)_{G,H}$.

\subsubsection{}
 Let $k\rarrow k'$ be a morphism of commutative rings of finite
weak homological dimension.
 Then for any complex of smooth $G$\module s $N^\bu$ over~$k'$
the image of the object $\SemiTor^{\S_{k'}(G,H)}(N^\bu,k')$ under
the restriction of scalars functor $\sD(k'\modl)\rarrow\sD(k\modl)$
is naturally isomorphic to the object $\SemiTor^{\S_k(G,H)}(N^\bu,k)$.
 This follows from Corollary~\ref{semi-pull-push-derived}.3(a).

\begin{rmk}
 When $G$ is a discrete group, $H=\{e\}$ is the trivial subgroup,
and $N$ is a $G$\module{} over~$k$, the cohomology of the object
$\SemiTor^{\S_k(G,\{e\})}(k,N)$ coincides with the discrete group
homology $H_*(G,N)$ and is concentrated in the nonpositive
cohomological degrees.
 When $G$ is a profinite group, $H=G$ is the whole group, and
$N$ is a smooth (discrete) $G$\module{} over~$k$, the cohomology
of $\SemiTor^{\S_k(G,G)}(k,N)$ coincides with the profinite group
cohomology $H^*(G,N)$ and is concentrated in the nonnegative
cohomological degrees.
 More generally, when $G=G_1\times G_2$ is the product of a discrete
group $G_1$ and a profinite group $G_2$, \ $H=G_2\subset G$, \ $k$
is a field, and $N=N_1\ot_k N_2$ is the tensor product of
a $G_1$\module{} and a smooth $G_2$\module, the cohomology of
$\SemiTor^{\S_k(G,H)}(k,N)$ is isomorphic to the tensor product
$H_*(G_1,N_1)\ot_k H^*(G_2,N_2)$.
 Applying these observations to the case of a finite group~$G$, one
can see that the semi-infinite homology of topological groups
$\SemiTor^{\S_k(G,H)}(k,N^\bu)$ does depend essentially on the choice of
an open profinite subgroup $H\subset G$.
\end{rmk}

\subsubsection{}
 Assume that the ring~$k$ has a finite homological dimension.
 Then for any complex of smooth $G$\module s $M^\bu$ over~$k$
and any complex of $G$\contramodule s $P^\bu$ over~$k$ the object
$\SemiExt_{\S_k(G,H)}(M^\bu,P^\bu)$ in the derived category of
$k$\module s is defined.
 Furthermore, whenever either $M^\bu$ is a complex of projective
$k$\module s or $P^\bu$ is a complex of injective $k$\module s
there is a natural isomorphism $\SemiExt_{\S_k(G,H)}(M^\bu,P^\bu)
\simeq\SemiExt_{\S_k(G,H)}(k,\Hom_k(M^\bu,P^\bu))$ in $\sD(k\modl)$.

\subsubsection{}
 Let $k\rarrow k'$ be a morphism of commutative rings of finite
homological dimension.
 Then for any complex of $G$\contramodule s $P^\bu$ over~$k'$
the image of the object $\SemiExt_{\S_{k'}(G,H)}(k',P^\bu)$ under
the restriction of scalars functor $\sD(k'\modl)\rarrow\sD(k\modl)$
is naturally isomorphic to the object $\SemiExt_{\S_k(G,H)}(k,P^\bu)$.

\subsection{Remarks on the Gaitsgory--Kazhdan construction}
 This subsection contains some comments on
the papers~\cite{GK1,GK2}.

\subsubsection{}
 Let $G$ be a topological group and $k$ be a field of 
characteristic~$0$.
 Then the category of discrete $G\times G$\module s over~$k$ has
a tensor category structure with the tensor product of two modules
$K'\ot_G K''$ defined as the module of coinvariants of the action
of $G$ in $K'\ot_k K''$ induced by the action of the second copy
of~$G$ in $K'$ and the action of the first copy of~$G$ in $K''$.
 The category of discrete $G$\module s over~$k$ has structures of
a left and a right module category over this tensor category;
and the functor $(N,M)\mpsto N\ot_G M=(N\ot_k M)_G$ defines
a pairing between these two module categories taking values
in the category of $k$\+vector spaces.
 A finitely-additive $k$\+valued measure defined on compact open
subsets of $G$ is called \emph{smooth} if it is equal to the product
of a locally constant $k$\+valued function on $G$ and a (left or
right invariant) Haar measure on $G$.
 The $G\times G$\module{} of compactly supported smooth $k$\+valued
measures on $G$ is the unit object of the above tensor category.

 Notice that the above $G\times G$\module{} has a natural
$k$\+algebra structure, given by the convolution of measures.
 However, this algebra has no unit and the category of (left
or right) modules over it contains the category of smooth
$G$\module s over~$k$ as a proper full subcategory.

\subsubsection{}
 For any category $\sC$, denote by $\Pro\sC$ and $\Inj\sC$
the categories of pro-objects and ind-objects in~$\sC$.
 Let $\Set_\fin$ denote the category of finite sets.
 We will identify the category of compact topological spaces with
the category $\Pro\Set_\fin$, the category of (discrete) sets
with the category $\Inj\Set_\fin$, and the category of topological
spaces with a full subcategory of the category $\Inj\Pro\Set_\fin$
formed by the inductive systems of profinite sets and their open
embeddings.
 In addition, we will consider the category $\Inj\Pro\Set_\fin$
as a full subcategory in $\Pro\Inj\Pro\Set_\fin$ and the latter
category as a full subcategory in $\Inj\Pro\Inj\Pro\Set_\fin$.

 Let $\boH$ be a group object in $\Pro\Inj\Pro\Set_\fin$ such
that $\boH$ can be represented by a projective system of topological
groups in $\Inj\Pro\Set_\fin$ and open surjective morphisms between
them.
 A representation of $\boH$ in $k\vect$ is just a smooth 
$G$\module{} over $k$, where $G$ is a quotient group object of $\boH$
that is a topological group in $\Inj\Pro\Set_\fin$.
 The category $\Rep_k(\boH)$ of representations of $\boH$ in
$\Pro(k\vect)$, defined in~\cite{GK1}, is equivalent to the category of
pro-objects in the category of representations of $\boH$ in $k\vect$.
 So the category $\Rep_k(\boH\times\boH)$ has a natural tensor category
structure with the unit object given by the projective system formed by
the modules of smooth compactly supported measures on topological
groups that are quotient groups of $\boH$ and the push-forward maps
between the spaces of measures.
 The category $\Rep_k(\boH)$ is a left and a right module category
over this tensor category, and there is a pairing between these left
and right module categories taking values in $\Pro(k\vect$).

\subsubsection{}
 Let $\boG$ be a group object in $\Inj\Pro\Inj\Pro\Set_\fin$ and
let $\boH\subset\boG$ be a subgroup object which belongs to
$\Pro\Inj\Pro\Set_\fin$.
 Assume that the object $\boG$ is given by an inductive system of
objects $\boG_\alpha\in\Pro\Inj\Pro\Set_\fin$ and the group object
$\boH$ is given by a projective system of its quotient group objects
$\boH/\boH^i\in\Inj\Pro\Set_\fin$ satisfying the following conditions.
 It is convenient to assume that $\boG_0=\boH=\boH^0$ is the initial
object in the inductive system of $\boG_\alpha$ and the final object
in the projective system of $\boH^i$.
 As the notation suggests, $\boH^i$ are normal subgroup objects of
$\boH$.
 The group object $\boH$ acts in $\boG_\alpha$ in a way compatible
with the action of $\boH$ in $\boG$ by right multiplications.
 The quotient objects $\boG_\alpha/\boH^i$ are topological spaces
in $\Inj\Pro\Set_\fin$.
 The morphisms $\boG_\alpha/\boH^i\rarrow\boG_\beta/\boH^i$ are
closed embeddings of topological spaces; the morphisms 
$\boG_\alpha/\boH^i\rarrow\boG_\alpha/\boH^j$ are principal
$\boH^i/\boH^j$\+bundles.
 Finally, the quotient objects $\boG_\alpha/\boH$ are compact
topological spaces, i.~e., belong to $\Pro\Set_\fin\subset
\Inj\Pro\Set_\fin$.

 Let $\boG'$ be a group object in $\Inj\Pro\Inj\Pro\Set_\fin$
endowed with a central subgroup object identified with
the multiplicative group $k^*$, which is considered as a discrete
topological group, that is a group object in $\Inj\Set_\fin
\subset\Pro\Inj\Pro\Set_\fin$.
 Suppose that the quotient group object $\boG'/k^*$ is identified
with $\boG$ and the central extension $\boG'\rarrow\boG$ is
endowed with a splitting over $\boH$, i.~e., $\boH$ is a subgroup
object in $\boG'$. 
 Moreover, denote by $\boG_\alpha'$ of the preimages of $\boG_\alpha$
in $\boG'$ and assume that the morphisms $\boG_\alpha'/\boH^i\rarrow
\boG_\alpha/\boH^i$ are principal $k^*$\+bundles of topological
spaces in $\Inj\Pro\Set_\fin$.

\subsubsection{}
 One example of such a group $\boG'$ is provided by the canonical
central extension $\boG\til$ of the group $\boG$ with the kernel
$k^*$, which is constructed as follows.
 For each $\alpha$ let us choose $i$ such that
$\Ad_{\boG_\alpha^{-1}}(\boH^i)\subset\boH$ and $j$ such that
$\Ad_{\boG_\alpha}(\boH^j)\subset\boH^i$.
 For any topological group $\boG$ denote by $\mu(G)$
the one-dimensional vector space of left invariant finitely-additive
$k$\+valued Haar measures defined on compact open subsets of~$\boG$.
 An element of the topological space $\boG\til_\alpha/\boH$ is a pair
consisting of an element of $g\in\boG_\alpha/\boH$ and an isomorphism
$\mu(\boH^i/\Ad_g(\boH^j))\simeq\mu(\boH^i/\boH^j)$.
 The topology on $\boG\til_\alpha/\boH$ is defined by the condition
that the following set of sections of the $k^*$\+torsor
$\boG\til_\alpha/\boH\rarrow\boG_\alpha/\boH$ consists of continuous
maps.
 Choose $m$ such that $\Ad_{\boG_\alpha^{-1}}(\boH^m)\subset\boH^j$,
a compact open subset $U$ in the quotient group $\boH^i/\boH^m$, and
an element $a\in k^*$; for each $g\in\boG_\alpha/\boH$ define
the isomorphism $\mu(\boH^i/\Ad_g(\boH^j))\simeq\mu(\boH^i/\boH^j)$
so that the left-invariant measure on $\boH^i/\Ad_g(\boH^j)$ for which
the measure of the image of $U$ is equal to~$1$ corresponds to
the left-invariant measure on $\boH^i/\boH^j$ for which the measure
of the image of $U$ is equal to~$a$.
 The ratio of any two such sections is a locally constant 
function.
 Now the object $\boG\til_\alpha$ of $\Pro\Inj\Pro\Set_\fin$ is
the fibered product of $\boG\til_\alpha/\boH$ and $\boG_\alpha$
over $\boG_\alpha/\boH$; it is easy to define the group structure
on $\boG\til$ (one should first check that the construction of
$\boG\til_\alpha/\boH$ does not depend on the choice of $\boH^i$
and $\boH^j$).

\subsubsection{}
 Let $c'\:\boG'\rarrow\boG$ be a central extension satisfying
the above conditions.
 Denote by $\Rep_{c'}(\boG)$ the category of representations of
$\boG'$ in $\Pro(k\vect)$ in which the central subgroup 
$k^*\subset\boG'$ acts tautologically by automorphisms proportional
to the identity, as defined in~\cite{GK1}.
 Then the forgetful functor $\Rep_{c'}(\boG)\rarrow\Rep(\boH)$
admits a right adjoint functor, which can be described as
the functor of tensor product over $\boH$ with a certain
representation of $\boG'\times\boH$ in $\Pro(k\vect)$.
 The underlying pro-vector space of this representation, denoted by
$\bC_{c'}(\boG,\boH)$, is the space of ``pro-semimeasures on $\boG$
relative to $\boH$ on the level~$c'$''; it is given by the projective
system formed by the vector spaces $k_{c'}(\boG_\alpha/\boH^i)\ot_k
\mu(\boH/\boH^i)$, where the first factor is the space of locally
constant compactly supported functions on $\boG'_\alpha/\boH^i$ which
transform by the tautological character under the action of~$k^*$.
 The morphism in this projective system corresponding to
a change of~$\alpha$ is the pull-back map with respect to a closed
embedding, while the morphism corresponding to a change of~$i$
is the map of integration along the fibers of a principal bundle.
 Hence the representation $\bC_{c'}(\boG,\boH)$ considered as
an object of the tensor category $\Rep(\boH\times\boH)$ is endowed
with a structure of coring (with counit), and it follows from
Theorem~\ref{barr-beck-theorem} that the category $\Rep_{c'}(\boG)$
is equivalent to the category of left comodules over
$\bC_{c'}(\boG,\boH)$ in $\Rep(\boH)$.

 Now let $c''\:\boG''\rarrow\boG$ be the central extension satisfying
the same conditions that is complementary to $c'$, i.~e., the Baer
sum $c'+c''$ is identified with minus the canonical central
extension $c_0\:\boG\til\rarrow\boG$.
 Gaitsgory and Kazhdan noticed that one can extend the right action
of $\boH$ in $\bC_{c'}(\boG,\boH)$ to an action of $\boG''$ commuting
with the left action of $\boG'$, with the central subgroup $k^*$ of
$\boG''$ acting tautologically.
 Moreover, in~\cite{GK2} there is a construction of a natural
anti-isomorphism of corings $\bC_{c'}(\boG,\boH)\simeq
\bC_{c''}(\boG,\boH)$ permuting the left and right actions
of~$\boG'$ and $\boG''$.
 Thus the category $\Rep_{c''}(\boG)$ is equivalent to the category
of right comodules over $\bC_{c'}(\boG,\boH)$ in $\Rep(\boH)$.

\subsubsection{}
 So there is the functor of cotensor product
$$
 \oc_{\bC_{c'}(\boG,\boH)}\:\Rep_{c''}(\boG)\times\Rep_{c'}(\boG)
 \lrarrow\Pro(k\vect),
$$
which is called ``semi-invariants'' in~\cite{GK2}.
 This functor is neither left, nor right exact in general.
 One can construct its double-sided derived functor in the way
analogous to that of Remark~\ref{semitor-definition}, at least
when the set of indices~$i$ is countable.

 The semi-derived category $\sD^\si(\Rep_{c'}(\boG))$ is defined as
the quotient category of the homotopy category $\Hot(\Rep_{c'}(\boG))$
by the thick subcategory of complexes that are contraacyclic as
complexes over the abelian category $\Rep(\boH)$.
 Then Lemma~\ref{semitor-definition} allows to define the double-sided
derived functor
$$
 \ProCotor_{\bC_{c'}(\boG,\boH)}\:\sD^\si(\Rep_{c''}(\boG))\times
 \sD^\si(\Rep_{c'}(\boG))\lrarrow\sD(\Pro(k\vect))
$$
in terms of coflat complexes in $\Hot(\Rep_{c''}(\boG))$ and
$\Hot(\Rep_{c'}(\boG))$.
 The key step is to construct for any object of $\Rep_{c'}(\boG)$
a surjective map onto it from an object of $\Rep_{c'}(\boG)$
that is flat as a representation of~$\boH$.
 This construction is dual to that of
Lemma~\ref{coflat-semimodule-injection}; it is based on the fact
that any module over a topological group~$G$ induced from
the trivial module~$k$ over a compact open subgroup $H\subset G$
is flat with respect to the tensor product of discrete $G$\module s
over~$k$.

\begin{qst}
 Can the cotensor product $\bN\oc_{\bC_{c'}(\boG,\boH)}\bM$
of an object $\bN\in\Rep_{c''}(\boG)$ and an object
$\bM\in\Rep_{c'}(\boG)$ be recovered from the tensor product
$\bN\ot_k^{\mathrm p}\bM$ in the category of pro-vector spaces,
considered as a representation of $\boG\til$ with the diagonal action?
\end{qst}

\Section{Algebraic Groupoids with Closed Subgroupoids}
\label{groupoid-appendix}

 To any smooth affine groupoid $(M,H)$ one can associate a coring
$\C(H)$ over a ring $A(M)$ and a natural left and right coflat
Morita autoequivalence $(\E,\E\dual)$ of $\C(H)$; it has the form
$\C(H)\ot_{A(M)}E\simeq\E\simeq E\ot_{A(M)}\C(H)$ and $E\dual\ot_{A(M)}
\C(H)\simeq \E\dual\simeq\C(H)\ot_{A(M)}E\dual$, where $(E,E\dual)$
is a certain pair of mutually inverse invertible $A(M)$\module s.
 To any groupoid $(M,G)$ containing $(M,H)$ as a closed subgroupoid,
one can assign two opposite semialgebras $\S^l(G,H)$ and $\S^r(G,H)$
over $\C(H)$ together with a natural left and right semiflat
Morita equivalence $(\bE,\bE\dual)$ between them formed by
the bisemimodules $\S^l(G,H)\oc_{\C(H)}\E\simeq\bE\simeq\E
\oc_{\C(H)}\S^r(G,H)$ and $\E\dual\oc_{\C(H)}\S^l(G,H)\simeq\bE\dual
\simeq\S^r(G,H)\oc_{\C(H)}\E\dual$.
 To obtain these results, we will have to assume the existence of
a quotient variety~$G/H$.

\smallskip
 In this appendix, by a variety we mean a smooth algebraic variety
(smooth separated scheme) over a fixed ground field~$k$ of zero
characteristic.
 The structure sheaf of a variety $X$ is denoted by $O=O_X$ and
the sheaf of differential top forms by $\Omega=\Omega_X$; for
any invertible sheaf $L$ over $X$, its tensor powers are denoted
by~$L^n$ for $n\in\boZ$.

\subsection{Coring associated to affine groupoid} \label{groupoid-coring}
 A \emph{(smooth) groupoid} $(M,G)$ is a set of data consisting
of two varieties $M$ and $G$, two smooth morphisms $s$, $t\:G\birarrow M$
of \emph{source} and \emph{target}, a \emph{unit} morphism
$e\:M\rarrow G$, a \emph{multiplication} morphism $m\: G\times_M G
\rarrow G$ (where the first factor $G$ in the fibered product
$G\times_MG$ maps to $M$ by the morphism~$s$ and the second factor $G$
maps to $M$ by the morphism~$t$), and an \emph{inverse element}
morphism $i\:G\rarrow G$.
 The following equations should be satisfied: first, $se=\id_M=te$
and $m(e\times \id_G)=\id_G=m(\id_G\times e)$ (unity); second,
$tm=tp_1$, \ $sm=sp_2$, and $m(m\times\id_G) = m(\id_G\times m)$,
where $p_1$ and $p_2$ denote the canonical projections of the fibered
product $G\times_MG$ to the first and the second factors, respectively
(associativity); third, $si=t$, \ $ti=s$, \ $m(i\times\id_G)\Delta_t=es$,
and $m(\id_G\times i)\Delta_s=et$, where $\Delta_s$ and $\Delta_t$
denote the diagonal embeddings of $G$ into the fibered squares of $G$
over $M$ with respect to the morphisms $s$ and $t$, respectively
(inverseness).
 It follows from these equations that $i^2=\id_G$.

 A groupoid $(M,H)$ is said to be \emph{affine} if $M$ and $H$ are
affine varieties.
 In the sequel $(M,H)$ denotes an affine groupoid; its structure
morphisms are denoted by the same letters $s$, $t$, $e$, $m$, $i$.

 Let $A=A(M)=O(M)$ and $\C=\C(H)=O(H)$ be the rings of functions
on $M$ and $H$, respectively.
 The maps of source and target $s$, $t\:H\birarrow M$ induce two maps
of rings $A\birarrow\C$, which endow $\C$ with two structures of
$A$\module; we will consider the $A$\module{} structure on $\C$
coming from the morphism~$t$ as a left module structure and
the $A$\module{} structure on $\C$ coming from the morphism~$s$ as
a right module structure.
 Then there is a natural isomorphism $O(H\times_M H)\simeq \C\ot_A\C$,
hence the multiplication morphism $m\:H\times_M H\rarrow H$ induces
a comultiplication map $\C\rarrow\C\ot_A\C$.
 Besides, the unit map $e\:M\rarrow H$ induces a counit map
$\C\rarrow A$.
 It follows from the associativity and unity equations of
the groupoid $(M,H)$ that these comultiplication and counit maps
are morphisms of $A$\+$A$\bimodule s satisfying the coassociativity
and counity equations; so $\C$ is a coring over~$A$.
 Clearly, $\C$ is a coflat left and right $A$\module.

\subsection{Canonical Morita autoequivalence}
\label{canonical-morita-auto}
 Denote by $V=V_H$ the invertible sheaf $\Omega_H\ot s^*(\Omega_M^{-1})
\ot t^*(\Omega_M^{-1})$ on~$H$.
 Let $q_1$ and $q_2$ denote the canonical projections of the fibered
product $H\times_MH$ to the first and the second factors, respectively.
 Then there are natural isomorphisms $q_1^*(V)\simeq m^*(V)\simeq
q_2^*(V)$ of invertible sheaves on $H\times_M H$.
 Indeed, one has $q_1^*(\Omega_H)\simeq\Omega_{H\times_MH}\ot q_2^*
(\Omega_H^{-1})\ot q_2^*t^*(\Omega_M)$ and $m^*(\Omega_H)\simeq
\Omega_{H\times_MH}\ot q_2^*(\Omega_H^{-1})\ot q_2^*s^*(\Omega_M)$.
 Now denote by $U$ the invertible sheaf $e^*\Omega_H\ot \Omega_M^{-2}$
on~$M$.
 Applying the functors of inverse image with respect to the morphisms
$e\times\id_H$, $\id_H\times e\:H\birarrow H\times_MH$ to the above
isomorphisms, one obtains natural isomorphisms of invertible sheaves
$t^*U\simeq V\simeq s^*U$.

 Set $\E=V(H)$ and $\E\dual=V^{-1}(H)$.
 Then $\E$ and $\E\dual$ are $\C$\module s, and consequently 
$A$\+$A$\bimodule s.
 The pull-back map $V(H)\rarrow m^*(V)(H\times_MH)$ with respect to
the multiplication morphism~$m$ together with the isomorphisms
$m^*(V)(H\times_MH)\simeq q_1^*(V)(H\times_MH)\simeq\E\ot_A\C$ and
$m^*(V)(H\times_MH)\simeq q_2^*(V)(H\times_MH)\simeq\C\ot_A\E$
defines the right and left coactions of $\C$ in~$\E$. 
 It follows from the associativity equation of~$H$ that these
coactions commute; so $\E$ is a $\C$\+$\C$\bicomodule.
 Analogously one defines a $\C$\+$\C$\bimodule{} structure
on~$\E\dual$.
 Set $E=U(M)$ and $E\dual=U^{-1}(M)$; then there are natural
isomorphisms of $\C$\comodule s $\C\ot_AE\simeq\E\simeq E\ot_A\C$
and $E\dual\ot_A\C\simeq\E\dual\simeq\C\ot_AE\dual$.
 These isomorphisms have the property that two maps $\E\simeq
\C\ot_AE\rarrow E$ and $\E\simeq E\ot_A\C\rarrow E$ induced by
the counit map $\C\rarrow A$ coincide, and analogously
for~$\E\dual$.
 Besides, there are obvious isomorphisms $E\ot_AE\dual\simeq A
\simeq E\dual\ot_AE$.
 
 It follows that $\E\oc_\C\E\dual\simeq (E\ot_A\C)\oc_\C
(\C\ot_A E\dual)\simeq E\ot_A\C\ot_AE\dual\simeq\C$, and
analogously $\E\dual\oc_\C\E\simeq\C$.
 So the pair $(\E,\E\dual)$ is a left and right coflat Morita
equivalence (see~\ref{co-contra-morita-remarks}) between $\C$
and itself.
 Since the bicomodules $\E$ and $\E\dual$ can be expressed in
the above form in terms of $A$\module s $E$ and $E\dual$, it
follows that there are natural isomorphisms of corings
$E\ot_A\C\ot_AE\dual\simeq\C\simeq E\dual\ot_A\C\ot_AE$. 

\subsection{Distributions and generalized sections}
 Let $X\supset Z$ be a variety with a (smooth) closed subvariety
and $L$ be a locally constant sheaf on~$X$.
 The sheaf $L^Z$ of generalized sections of $L$, supported in $Z$
and regular along $Z$, can be defined as the image of~$L$ with
respect to the $d$\+th right derived functor of the functor
assigning to any quasi-coherent sheaf on $X$ its maximal subsheaf
supported set-theoretically in $Z$, where $d=\dim X - \dim Z$.
 The sheaf $L^Z$ is a quasi-coherent sheaf on $X$ supported
set-theoretically in~$Z$.

 There is a natural isomorphism $L^Z\simeq L\ot_{O_X} O_X^Z$.
 The sheaf $\Omega_X^Z$ can be alternatively defined as the direct
image of the constant right module $\Omega_Z$ over the sheaf of
differential operators $\Diff_Z$ under the closed embedding
$Z\rarrow X$ (see~\cite{Bern}); this makes $\Omega_X^Z$ not only
an $O_X$\module, but even a $\Diff_X$\module.
 The sheaf $\Omega_X^Z$ is called the sheaf of distributions on $X$,
supported in $Z$ and regular along~$Z$.

 Let $g\:Y\rarrow X$ be a morphism of varieties and $Z\subset X$ be
a closed subvariety.
 Assume that the fibered product $Z\times_XY$ is smooth and
$\dim Y - \dim Z\times_XY = \dim X - \dim Z$ if $Z\times_XY$
is nonempty.
 Then there is a natural isomorphism $g^*(L^Z)\simeq (g^*L)^{Z\times_XY}$
of quasi-coherent sheaves on~$Y$.
 In particular, there is a natural pull-back map of the modules of
global generalized sections $g^+\:L^Z(X)\rarrow(g^*L)^{Z\times_XY}(Y)$.
 
 Let $h\:W\rarrow X$ be a morphism of varieties and $Z\subset W$
be a closed subvariety such that the composition $Z\rarrow W\rarrow X$
is also a closed embedding.
 Then there is a natural push-forward map $h_*(\Omega^Z_W)\rarrow
\Omega^Z_X$ of quasi-coherent sheaves on~$X$.
 Consequently, for any locally constant sheaves $L'$ on $X$ and
$L''$ on $W$ endowed with an isomorphism $L''\ot\Omega_W^{-1}
\simeq h^*(L'\ot\Omega_X^{-1})$ there is a push-forward map
$h_*(L''{}^Z)\rarrow L'{}^Z$.
 In particular, there is a natural push-forward map of the modules
of generalized sections $h_+\:L''{}^Z(W)\rarrow L'{}^Z(X)$.

 Let $g\:Y\rarrow X$ and $h\:W\rarrow X$ be morphisms of varieties 
satisfying the above conditions with respect to a closed subvariety
$Z\subset W$.
 Assume also that the fibered product $W\times_XY$ is smooth and 
$\dim W\times_XY + \dim X = \dim W + \dim Y$ if $W\times_XY$ is
nonempty.
 Set $\tilde g =\id_W\times g\:W\times_XY\rarrow W$ and
$\tilde h = h\times\id_Y\: W\times_XY\rarrow Y$.
 Then there is a natural isomorphism of invertible sheaves
$\tilde g^*\Omega_W\ot\Omega_{W\times_X Y}^{-1}\simeq \tilde h^*
(g^*\Omega_X\ot\Omega_Y^{-1})$ on $W\times_XY$, as one can see
by decomposing $g$ and $h$ into closed embeddings followed by
smooth morphisms.
 For any quasi-coherent sheaf $F$ on $W$ supported set-theoretically
in~$Z$ there is a natural isomorphism $\tilde h_*\tilde g^* F\simeq
g^*h_*F$ of quasi-coherent sheaves on $W\times_XY$.
 The push-forward maps of the sheaves of distributions with
respect to the morphisms $h$ and $\tilde h$ are compatible with
the pull-back isomorphisms with respect to $g$ and~$\tilde g$
in the obvious sense.

 The sheaf of generalized sections $L^Z$ of a locally constant sheaf
$L$ on a variety $X\supset Z$ is endowed with a natural increasing
filtration by coherent subsheaves $F_nL^Z$ of generalized sections
of order no greater than~$n$.
 This filtration is preserved by all the above natural isomorphisms
and maps.
 The associated graded sheaf $\gr_F\Omega_X^Z$ is the direct image
under our closed embedding $\iota\:Z\rarrow X$ of a sheaf of
$O_Z$\module s naturally isomorphic to the tensor product
$\Omega_Z\ot_{O_Z}\Sym_{O_Z}\!N_{Z,X}$, where $N_{Z,X}$ is the normal
bundle to $Z$ in~$X$ and $\Sym$ denotes the symmetric algebra.

 In particular, there is a natural isomorphism $\lambda_0\:
\iota_*\Omega_Z\rarrow F_0\Omega^Z_X$.
 Furthermore, there is a natural map of sheaves of $k$\+vector
spaces $\lambda_1\:\iota_*\Omega_Z\ot_{O_X}T_X\rarrow F_1\Omega^Z_X$
which induces the isomorphism $\iota_*(\Omega_Z\ot_{O_Z}N_{Z,X})
\simeq F_1\Omega^Z_X/F_0\Omega^Z_X$, where $T=T_X$ denotes
the tangent bundle of~$X$.
 The map $\lambda_1$ satisfies the equation $\lambda_1(f\omega\ot v)
= \lambda_1(\omega\ot fv) = f\lambda_1(\omega\ot v) - 
\lambda_0(v(f)\omega)$ for local sections $f\in O_X$, \
$\omega\in\Omega_Z$, and $v\in T_X$, where $(v,f)\mpsto v(f)$
denotes the action of vector fields in functions.

\subsection{Lie algebroid of a groupoid}
 A \emph{Lie algebroid} $\g$ over a commutative ring $A$ is
an $A$\module{} endowed with a Lie algebra structure and a Lie
action of $\g$ by derivations of~$A$ satisfying the equations
$[x,ay]=a[x,y]+x(a)y$ and $(ax)(b)=a(x(b))$ for $a$, $b\in A$,
\ $x$, $y\in\g$.
 The \emph{enveloping algebra} $U_A(\g)$ of a Lie algebroid $\g$
over~$A$ is generated by $A$ and $\g$ with the relations
$a\cdot b=ab$, \ $a\cdot x= ax$, \ $x\cdot a = ax + x(a)$,
and $x\cdot y - y\cdot x = [x,y]$, where $(u,v)\mpsto u\cdot v$
denotes the multiplication in $U_A(\g)$.
 The algebra $U_A(\g)$ is endowed with a natural increasing
filtration $F_nU_A(g)$ defined by the rules $F_0U_A(g) = \im A$,
\ $F_1U_A(\g) = \im A+\im\g$, and $F_nU_A(\g)=F_1U_A(\g)^n$ for $n\ge1$.
 When $\g$ is a flat $A$\module, the associated graded algebra
$\gr_FU_A(\g)$ is isomorphic to the symmetric algebra $\Sym_A(\g)$
of the $A$\module~$\g$.

 Let $(M,G)$ be a groupoid with an affine base variety~$M$; set
$A=A(M)=O(M)$.
 Then the $A$\module{} $\g=N_{e(M),G}(M)$ has a natural Lie
algebroid structure.
 To define the action of $\g$ in $A$, consider the $A$\module{}
$(e^*T_G)(M)$ of vector fields on~$e(M)$ tangent to~$G$.
 There are natural push-forward morphisms $s_+$, $t_+\:
(e^*T_G)(M)\birarrow T(M)$.
 Identify $\g$ with the kernel of the morphism $t_+$; then the action
of $\g$ in $A$ is defined in terms of the map $s_+\:\g\rarrow T(M)$.
 To define the Lie algebra structure on~$\g$, we will embed $\g$
into a certain module of generalized sections on~$G$, supported
in $e(M)$ and regular along~$e(M)$.

 Set $K^l(G)=(t^*\Omega_M^{-1}\ot\Omega_G)^{e(M)}(G)$; this module of
generalized sections in endowed with a natural filtration~$F$.
 The $O(G)$\module{} structure on $K^l(G)$ induces
an $A$\+$A$\bimodule{} structure; as in~\ref{groupoid-coring},
we consider the $A$\module{} structure coming from the morphism~$t$
as a left $A$\module{} structure and the $A$\module{} structure
coming from the morphism~$s$ as a right $A$\module{} structure.
 There is a natural isomorphism of $A$\+$A$\bimodule s
$A\simeq F_0K^l(G)$.
 Define a $k$\+linear map of sheaves $e_*N_{e(M),G}\rarrow
(t^*\Omega_M^{-1}\ot\Omega_G)^{e(M)}$ locally by the formula
$v\mpsto t^*\omega^{-1}\ot\lambda_1(\omega\ot v)$, where $v$ is
a local vector field on $e(M)$ tangent to $G$ such that $t_*(v)=0$
and $\omega$ is a local nonvanishing top form on~$M$; it is easy
this expression does not depend on the choice
of~$\omega$.
 Passing to the global sections, we obtain an injective map
$\g\rarrow F_1K^l(G)$ inducing an isomorphism $\g\simeq
F_1K^l(G)/F_0K^l(G)$.
 This injective map and the $A$\+$A$\bimodule{} structure satisfy
the compatibility equations $a\cdot x=ax$ and $x\cdot a=x(a)+ax$
for $x\in\g\subset F_1K^l(G)$ and $a\in A$, where $(a,x)\mpsto ax$
denotes the action of $A$ in~$\g$, while $(a,u)\mpsto a\cdot u$
and $(u,a)\mpsto u\cdot a$ denote the left and right actions
of $A$ in $F_1K^l(G)$.

 Let us define a $k$\+algebra structure on $K^l(G)$.
 There is a natural isomorphism $p_1^*(t^*\Omega_M^{-1}\ot\Omega_G)
\ot p_2^*(t^*\Omega_M^{-1}\ot\Omega_G)\simeq p_1^*t^*\Omega_M^{-1}
\ot \Omega_{G\times_MG}$ of invertible sheaves on $G\times_MG$.
 The pull-back with respect to the closed embedding $G\times_MG
\rarrow G\times_{\Spec k}G$ provides an isomorphism
$(p_1^*t^*\Omega_M^{-1}\ot \Omega_{G\times_MG})^{(e\times e)(M)}
(G\times_MG)\simeq K^l(G)\ot_A K^l(G)$, and the push-forward
with respect to the multiplication map $G\times_MG\rarrow G$ defines
an associative multiplication $K^l(G)\ot_A K^l(G)\rarrow K^l(G)$.
 The associated graded algebra $\gr_FK^l(G)$ is naturally isomorphic
to $\Sym_AF_1K^l(G)/F_0K^l(G)$.
 The formula $(u,a)\mpsto t_+(s^+(a)u)$ defines a left action of
$K^l(G)$ in $A$; the subspace $\g$ in $F_1K^l(G)$ is characterized
as the annihilator of the unit element of~$A$ under this action.
 Hence $\g$ is a Lie subalgebra of $K^l(G)$; this makes it a Lie algebra
and a Lie algebroid over~$A$.
 It follows that there is a natural isomorphism $U_A(\g)\simeq K^l(G)$.

 Analogously one defines an algebra structure on
the $A$\+$A$\bimodule{} of generalized sections $K^r(G)=
(s^*\Omega_M^{-1}\ot\Omega_G)^{e(M)}(G)$; then there is a natural
isomorphism of $k$\+algebras $U_A(\g)^\rop\simeq K^r(G)$.

\subsection{Two Morita equivalent semialgebras}
 Let $(M,H)\rarrow(M,G)$ be a morphism of smooth groupoids with
the same base variety $M$ such that the groupoid $H$ is affine
and the morphism of varieties $H\rarrow G$ is a closed embedding.
 Denote by $\h=N_{e(M),H}(M)$ the Lie algebroid of the groupoid~$H$;
then there is a natural injective morphism $\h\rarrow\g$ of Lie
algebroids over~$A(M)$.

 The formula $\phi_l(u,c) = t_+(i^+(c)u)$ defines a pairing
$\phi_l\:K^l(H)\ot_A\C\rarrow A$ between the algebra $K^l(H)$ and
the coring~$\C$ satisfying the conditions of~\ref{ring-coring-pairing}
with the left and right sides switched.
 The push-forward with respect to the closed embedding $H\rarrow G$
defines an injective morphism of $k$\+algebras $K^l(H)\rarrow K^l(G)$;
this is the enveloping algebra morphism induced by the morphism
of Lie algebroids $\h\rarrow\g$.
 Since $\h$ and $\g/\h$ are projective $A$\module s, $K^l(G)$ is 
a projective left and right $K^l(H)$\module.
 Set $\S^l=\S^l(G,H)=K^l(G)\ot_{K^l(H)}\C$; we will use the construction
of~\ref{semialgebra-constructed} to endow $\S^l$ with a structure of
semialgebra over~$\C$.

 Consider the cotensor product $\S^l\oc_\C\E\simeq K^l(G)\ot_{K^l(H)}\E$.
 Denote by $p_1$ and $q_2$ the projections of the fibered product
$G\times_MH$ to the first and the second factors.
 There is a natural isomorphism $p_1^*(t^*\Omega_M^{-1}\ot\Omega_G)
\ot q_2^*(V_H)\simeq p_1^*t^*\Omega_M^{-1}\ot\Omega_{G\times_MH}\ot
q_2^*s^*\Omega_M^{-1}$ of invertible sheaves on $G\times_MH$, where
$V_H=\Omega_H\ot s^*(\Omega_M^{-1})\ot t^*(\Omega_M^{-1})$ is 
the invertible sheaf on~$H$ defined in~\ref{canonical-morita-auto}.
 The pull-back with respect to the closed embedding $G\times_MH\rarrow
G\times_{\Spec k}H$ identifies the tensor product $K^l(G)\ot_A\E$ with
the module of generalized sections $(p_1^*t^*\Omega_M^{-1}\ot
\Omega_{G\times_MH}\ot q_2^*s^*\Omega_M^{-1})^{e(M)\times_MH}(G\times_MH)$.
 The push-forward with respect to the multiplication morphism
$G\times_MH\rarrow G$ defines a natural map $K^l(G)\ot_A\E\rarrow\bE$,
where $\bE = V_G^H(G)$ is the space of generalized sections of
the invertible sheaf $V_G=\Omega_G\ot s^*(\Omega_M^{-1})\ot
t^*(\Omega_M^{-1})$ on~$G$.
 It follows from the associativity equation for the two iterated
multiplication maps $G\times_MH\times_MH\birarrow G$ that this 
map factorizes through $K^l(G)\ot_{K^l(H)}\E$.
 The induced map $K^l(G)\ot_{K^l(H)}\E\rarrow\bE$ is an isomorphism,
since the associated graded map with respect to the filtrations~$F$ is.

 Denote by $q_1$ and $p_2$ the projections of the fibered product
$H\times_MG$ to the first and second factors.
 One constructs a natural isomorphism $m^*(V_G)\simeq p_2^*(V_G)$
of invertible sheaves on $H\times_MG$ in the same way as
in~\ref{canonical-morita-auto}.
 The pull-back with respect to the multiplication map $m\:H\times G
\rarrow G$ together with this isomorphism provide a left coaction
of $\C$ in~$\bE$.
 Analogously one defines a right coaction of $\C$ in~$\bE$; it follows
from the associativity equation for $H\times G\times H\birarrow G$
that these two coactions commute.
 We have $\S^l\oc_\C\E\simeq\bE$, thus the isomorphism
$\S^l\simeq\bE\oc_\C\E\dual$ provides a left coaction of $\C$ in $\S^l$
commuting with the natural right coaction of $\C$ in $\S^l$.

 Since the natural map $\E\rarrow\bE$ provided by the push-forward
with respect to the closed embedding $H\rarrow G$ is a morphism of
$\C$\+$\C$\bicomodule s, so is the semiunit map $\C\rarrow\S^l$.
 It remains to show that the semimulplication map $\S^l\oc_\C\S^l
\rarrow\S^l$ is a morphism of left $\C$\comodule s; here it suffices
to check that the map $\S^l\oc_\C\S^l\oc_\C\E\rarrow\S^l\oc_\C\E$
is a morphism of left $\C$\comodule s.
 After we have done with this verification, the latter map will define
a left $\S^l$\semimodule{} structure on~$\bE$.

 Analogously, define a pairing $\phi_r\:\C\ot_AK^r(H)\rarrow A$ by
the formula $\phi_r(c,u)=s_+(i^+(c)u)$ and set $\S^r=\S^r(G,H)=
\C\ot_{K^r(H)}K^r(G)$.
 The same construction makes $\S^r$ a semialgebra over $\C$ and
$\bE$ a right $\S^r$\semimodule.
 We will have to check that the left $\S^l$\semimodule{} and
the right $\S^r$\semimodule{} structures on $\bE$ commute.

 After this is done, we get an $\S^l$\+$\S^r$\bisemimodule{} 
$\S^l\oc_\C\E\simeq\bE\simeq\E\oc_\C\S^r$, where both the maps
$\E\birarrow\bE$ induced by the semiunit maps $\C\rarrow\S^l$
and $\C\rarrow\S^r$ coincide with the push-forward map
$V(H)\rarrow V_G^H(G)$ under the closed embedding $H\rarrow G$.
 The isomorphisms $\E\dual\oc_\C\S^l\simeq\E\dual\oc_\C\S^l\oc_\C\E
\oc_\C\E\dual \simeq\E\dual\oc_\C\E\oc_\C\S^r\oc_\C\E\dual\simeq
\S^r\oc_\C\E\dual$ define an $\S^r$\+$\S^l$\bisemimodule{}
$\E\dual\oc_\C\S^l=\bE\dual\simeq\S^r\oc_\C\E\dual$ endowed with
bisemimodule isomorphisms $\bE\os_{S^r}\bE\dual\simeq
(\E\oc_\C\S^r)\os_{\S^r}(\S^r\oc_\C\E\dual)\simeq\E\oc_\C\S^r
\oc_\C\E\dual\simeq\S^l$ and $\bE\dual\os_{\S^l}\bE\simeq
(\E\dual\oc_\C\S^l)\os_{\S^l}(\S^l\oc_\C\E)\simeq\E\dual\oc_\C\S^l
\oc_\C\E\simeq\S^r$.
 This provides a left and right semiflat Morita equivalence
$(\bE,\bE\dual)$ between the semialgebras $\S^l$ and $\S^r$,
and isomorphisms of semialgebras $\S^r\simeq\E\dual\oc_\C\S^l\oc_\C\E$
and $\S^l\simeq\E\oc_\C\S^r\oc_\C\E\dual$.
 (See~\ref{semi-morita-morphisms}
and~\ref{morita-change-of-coring-construction} for the relevant
definition and construction.)

\subsection{Compatibility verifications}
 In order to check that the map $\S^l\oc_\C\bE\simeq\S^l\oc_\C\S^l
\oc_\C\E\rarrow\S^l\oc_\C\E\simeq\bE$ is a morphism of left
$\C$\comodule s, we will identify this map with a certain
push-forward map of appropriate modules of generalized sections.

 Here we will need to assume the existence of a variety of left
cosets $G/H$ such that $G\rarrow G/H$ is a smooth surjective morphism
and the fibered square $G\times_{G/H}G$ can be identified with
$G\times_MH$ so that the canonical projection maps $G\times_{G/H}G
\birarrow G$ correspond to the projection and multiplication maps
$p_1$, $m\:G\times_MH\birarrow G$.
 Actually, we are interested in the quotient variety $G\times_HG$
of $G\times_MG$ by the equivalence relation $(g'h,g'')\sim
(g',hg'')$ for $g'$, $g''\in G$ and $h\in H$; it can be constructed
as either of the fibered products $H\backslash G\times_MG\simeq
G\times_HG = G\times_MG/H$, where $H\backslash G$ denotes the variety
of right cosets, $H\backslash G\simeq G/H$.
 Analogously one can construct the quotient variety
$G\times_HG\times_HG$ of the triple fibered product
$G\times_MG\times_MG$ by the equivalence relation $(g'h_1,g''h_2,
g''')\sim(g',h_1g'',h_2g''')$ for $g'$, $g''$, $g'''\in G$ and
$h_1$, $h_2\in H$.

 We have $\S^l\oc_\C\bE\simeq\bE\oc_\C\E\dual\oc_\C\bE$.
 Consider the natural map $r\:G\times_MH\times_MG\rarrow
G\times_HG$ given by the formula $(g',h,g'')\mpsto (g'h,g'')
=(g',hg'')$.
 Let $p_1$, $q_2$, $p_3$ denote the projections of the triple
fibered product $G\times_MH\times_MG$ to the three factors
and $n\:G\times_HG\rarrow G$ denote the multiplication morphism.
 There are natural isomorphisms $p_1^*t^*(\Omega_M^{-1})\ot
\Omega_{G\times_MH\times_MG}\ot p_3^*s^*(\Omega_M^{-1})\simeq
p_1^*(V_G)\ot q_2^*(\Omega_H)\ot p_3^*(V_G)$ and therefore
$r^*(n^*t^*\Omega_M^{-1}\ot\Omega_{G\times_HG}\ot n^*s^*\Omega_M^{-1})
\simeq p_1^*(V_G)\ot q_2^*(V_H^{-1})\ot p_3^*(V_G)$ of invertible
sheaves on $G\times_MH\times_MG$.

 The pull-back with respect to the closed embedding 
$G\times_MH\times_MG\rarrow G\times_{\Spec k}H\times_{\Spec k}G$
provides an isomorphism $(p_1^*V_G\ot q_2^*V_H^{-1}\ot p_3^*V_G)^
{H\times_MH\times_MH}(G\times_MH\times_MG)\simeq\bE\ot_A\E\dual\ot_A\bE$.
  The pull-back with respect to the smooth morphism~$r$ identifies
the module of generalized sections $(n^*t^*\Omega_M^{-1}\ot
\Omega_{G\times_HG}\ot n^*s^*\Omega_M^{-1})^{H\times_HH}(G\times_HG)$ with
the submodule $\bE\oc_\C\E\dual\oc_\C\bE\subset\bE\ot_A\E\dual
\ot_A\bE$, as one can see by identifying the tensor product $\bE\ot_A
\C\ot_A\E\dual\ot_A\C\ot_A\bE$ with a module of generalized sections
on the fibered square of $G\times_MH\times_MG$ over $G\times_HG$.
 Now our map $\bE\oc_\C\E\dual\oc_\C\bE\rarrow\bE$ is identified
with the push-forward map with respect to the multiplication
morphism~$n$; to check this, one can first identify the map
$K^l(G)\ot_A\bE\rarrow K^l(G)\ot_{K^l(H)}\bE\simeq\bE\oc_\C\E\dual
\oc_\C\bE$ with a push-forward map with respect to the morphism
$G\times_MG\rarrow G\times_HG$.
 The desired compatibility with the left $\C$\comodule{} structures
now follows from the commutation of the pull-back and push-forward
maps of generalized sections.

 To check that the left $\S^l$\semimodule{} and the right
$\S^r$\semimodule{} structures on $\bE$ commute, one can 
identify $\S^l\oc_\C\bE\oc_\C\S^r$ with a module of generalized
sections on $G\times_HG\times_HG$ and use the associativity
equation for the iterated multiplication maps $G\times_HG\times_HG
\birarrow G\times_HG\rarrow G$.

\end{document}